\documentclass[10pt,twoside,leqno,makeidx]{book}
\usepackage{amsfonts,amsmath}
\newtheorem{teo}{Theorem}[section]
\newtheorem{defi}{Definition}[section]
\newtheorem{prop}{Proposition}[section]
\newtheorem{rem}{Remark}[section]
\newtheorem{lem}{Lemma}[section]
\newtheorem{cor}{Corollary}[section]
\makeindex

\begin{document}

\setlength{\textwidth}{12.5cm}

\setlength{\textheight}{19.5cm}

\title{\textbf{THE GEOMETRY OF HIGHER-ORDER\\
HAMILTON SPACES\\} \textbf{\ Applications to Hamiltonian
Mechanics}}
\author{\textbf{} \\
\\
\textbf{Radu MIRON}\\
\textbf{'Al.I.Cuza' University of Iassy, Romania}}
\date{\bf 2002}
%\date{\vspace*{6.5cm}
%\textbf{Kluwer Academic Publishers}\\
%\textbf{\ Dordrecht. The Netherland}\\

\pagenumbering{roman}

\maketitle

\tableofcontents

\pagestyle{myheadings}

\markboth{\it{THE GEOMETRY OF HIGHER-ORDER HAMILTON SPACES\ \ \ \ \ }}{\it{CONTENTS}}

\chapter*{Preface}

\markboth{\it{THE GEOMETRY OF HIGHER-ORDER HAMILTON SPACES\ \ \ \ \ }}{\it{Preface}}

As is known, the Lagrange and Hamilton geometries have appeared relatively recently [76, 86].
Since 1980 these geometries have been intensively studied by mathematicians and physicists from Romania, Canada,
Germany, Japan, Russia, Hungary, U.S.A. etc.

Scientific prestigious manifestations devoted to Lagrange and
Hamilton geometries and their applications have been organized in
the above mentioned countries and have been published a number of
books and monographs by specialists in the field: R. Miron [94,
95], R. Miron and M. Anastasiei [99, 100], R. Miron, D. Hrimiuc,
H. Shimada and S. Sab\u au [115], P.L. Antonelli, R. Ingarden and
M. Matsumoto [7]. Finsler spaces, which form a subclass of the
class of Lagrange spaces, have been generated some excellent
books, due to M. Matsumoto [76], M. Abate and G. Patrizio [1], D.
Bao, S.S. Chern and Z. Shen [17] and A. Bejancu and H.R.Farran
[20]. Also, we remark on the monographs of M. Crampin [34], O.
Krupkova [72] and D.Opri\c s, I. Butulescu [125], D.Saunders [144]
- which contain pertinent applications in Analytical Mechanics and
in the Theory of Partial Differential Equations. But the direct
applications in Mechanics, Cosmology, Theoretical Physics and
Biology can be found in the well known books of P.L. Antonelli and
T. Zawstaniak [11], G.S. Asanov [14], S. Ikeda [59], M. de Leone
and P.Rodrigues [73].

The importance of Lagrange and Hamilton geometries consists of the fact that the variational problems for
important Lagrangians or Hamiltonians have numerous applications in various fields, as: Mathematics, Theory of
Dynamical Systems, Optimal Control, Biology, Economy etc.

In this respect, P.L. Antonelli's remark is interesting:

'There is now strong evidence that the symplectic geometry of
Hamiltonian dynamical  systems is deeply connected to Cartan
geometry, the dual of Finsler geometry', (see V.I.Arnold,
I.M.Gelfand and V.S.Retach [13]).

But, all of the above mentioned applications have imposed also the introduction of the notions of higher order
Lagrange spaces and, of course, of higher order Hamilton spaces. The base manifolds of these spaces are the bundles
of accelerations of superior order. The methods used in the construction of these geometry are the natural extensions
of the classical methods used in the edification of Lagrange and Hamilton geometries. These methods allows us to solve
an old problem of differential geometry formulated by Bianchi and Bompiani [94], more than 100 years ago. Namely, the
problem of prolongation of a Riemann structure $g$ defined on the base manifold $M,$ to the tangent bundle $T^kM,\ k>1.$
By means of this solution of the previous problem, we can construct, for the first time, good examples of regular
Lagrangians and Hamiltonians of higher order.

While the higher order Lagrange geometry has been sufficiently
developed we can not say the same thing about the higher order
Hamilton geometry. However the beginning was done only for the case
$k=2,$ in the book [115]. The reason comes from the fact that, in
the year 2001 we hadn't solved the variational problem for a
Hamiltonian which depends on the higher order accelerations and
momenta. This problem was solved this year, [98]. Another reason
was due to the absence of a consistent theory of subspaces in the
Hamilton space of order $k,\ k\geq 1,$ such kind of theory being
indispensable in applications.

In the present book, we give the general geometrical theory of the
Hamilton spaces of order $k\geq 1$. This is not a simple
generalization of the theory expounded for the case $k=2$ in the
monograph 'The Geometry of Hamilton and Lagrange Spaces' Kluwer
Acad. Publ. FTPH, no. 118, written by the present author together
with D. Hrimiuc, H. Shimada and V. S. Sab\u au, but it is a global
picture of this new geometry, extremely useful in applications
from Hamiltonian Mechanics, Quantum Physics, Optimal Control and
Biology.

Consequently, this book must be considered a direct continuation
of the monographs [94], [95], [99] and [115 ]. It contains new
developments of the subjects: variational principles for higher
order Hamiltonians; higher order energies; laws of conservations;
N\"other theorems; the Hamilton subspaces of order $k$ and their
fundamental equations. Also, the Cartan spaces of order $k$ are
investigated in details as dual of Finsler spaces of the same
order.

In this respect, a more explicit argumentation is as follows.

The geometry of Lagrange space of order $k\geq 1$ is based on the
geometrical edifice of the $k$-accelerations bundle $(T^kM,\pi ^k,M)$.

In Analytical Mechanics the manifold $M\,$ is the space of configurations of
a physical system. A point $x=(x^i)$, ($i=1,...,n=\dim M$) in $M$ is called
a configuration. A mapping $c:t\in I\rightarrow (x^i(t))\in U\subset M$ is a
law of moving (evolution), $t$ is time, a pair $(t,x)$ is an event and the $%
k $-uple $\left( \displaystyle\frac{dx^i}{dt},\cdots ,\displaystyle\frac
1{k!}\displaystyle\frac{d^kx^i}{dt^k}\right) $ gives the velocity and
generalized accelerations of order $1$, ..., $k-1$. The factors $%
\displaystyle\frac 1{h!}$ ($h=1,...,k$) are introduced here for
the simplicity of calculus. In this book we omit the word
'generalized' and say shortly, the accelerations of order $h$, for
$\displaystyle\frac 1{h!}\displaystyle\frac{ d^hx^i}{dt^h}$. A law
of moving $c:t\in I\rightarrow c\left( t\right) \in U$ will be
called a curve parametrized by time $t$.

A Lagrangian of order $k\geq 1$ is a real scalar function $%
L(x,y^{(1)},...,y^{(k)})$ on $T^kM$, where $y^{(h)i}=\displaystyle\frac
1{h!} \displaystyle\frac{d^hx^i}{dt^h}$. This definition is for autonomous
Lagrangians. A similar definition can be formulated for nonautonomous
Lagrangians of order $k$, by

\[
L:(t,x,y^{(1)},...,y^{(k)})\in \mathbf{R}\times T^kM\rightarrow
L(t,x,y^{(1)},...,y^{(k)})\in \mathbf{R,}
\]
$L$ being the scalar functions on the manifold $\mathbf{R}\times T^kM$.

The previous considerations can be done for the Hamiltonians of
order $k$. Let $p_i=\displaystyle\frac
12\displaystyle\frac{\partial L}{\partial y^{(k)i}}$ be the
'momenta' determined by the Lagrangian $L$ of order $k$. Then a
scalar function
\[
H:(x,y^{(1)},...,y^{(k-1)},p)\in T^{*k}M\rightarrow
H(x,y^{(1)},...,y^{(k-1)},p)\in \mathbf{R,}
\]
is an autonomous Hamiltonian of order $k$. It is a function of the
configurations $x$, accelerations $y^{(1)}$, ..., $y^{(k-1)}$ of order $1$,
..., $k-1$ and momenta $p$.

A similar definition can be formulated for nonautonomous Hamiltonian of order $k$.

For us it is preferable to study the autonomous Lagrangians and Hamiltonians,
because the notions of Lagrange space of order $k$ or Hamilton space of
order $k$ are geometrical concepts and one can construct these geometries
over the differentiable manifolds $T^kM$ and $T^{*k}M$, respectively. Of
course, the geometries of nonautonomous Lagrangians $%
L(t,x,y^{(1)},...,y^{(k)})$ and nonautonomous Hamiltonians $%
H(t,x,y^{(1)},...,y^{(k-1)},p)$ can be constructed by means of the same
methods. One obtains the rheonomic Lagrange spaces of order $k$ and
rheonomic Hamiltonian spaces of order $k.$

Now we know the usefulness of the geometry of Higher order
Lagrange space (see the book [94]). But why do we need a geometry
of Hamilton spaces of order $k \geq 1$? Clearly, for the same
reason as for the Lagrange case: For the determination of the
adequate geometrical models for the Hamiltonian Mechanics of order
k. This must be a natural extension of the classical Hamiltonian
Mechanics, expounded by V.I. Arnold in the book [12] or R.M.
Santilli in the book [139].

The problem is why did we use the manifold $T^{*k}M$ as the
background for the construction of the Hamilton geometry of order
$k$. The answer is as follows. We need a 'dual' of the $k$ -
acceleration bundle $(T^kM,\pi ^k,M)$ denoted by $(T^{*k}M,\pi
^{*k},M)$ which must have the following properties:

\qquad 1$^{\circ }.$ $T^{*1}M=T^{*}M$, (($T^{*}M,\pi ^{*},M)$ is the cotangent
bundle).

\qquad 2$^{\circ }.$ dim$T^{*k}M=\dim T^kM=(k+1)n$.

\qquad 3$^{\circ }.$ The manifold $T^{*k}M$ carries a natural presympletic
structure.

\qquad 4$^{\circ }.$ $T^{*k}M$ carries a natural Poisson structure.

\qquad 5$^{\circ }.$ $T^{*k}M$ is local diffeomorphic to $T^kM.$

We solved this problem by considering the differentiable bundle \\ $(T^{*k}M,\pi
^{*k},M)$ as the fibred bundle $(T^{k-1}M\times _MT^{*}M,\pi ^{k-1}\times _M\pi
^{*},M)$. So we have \[
T^{*k}M=T^{k-1}M\times _MT^{*}M,\quad \pi ^{*k}=\pi ^{k-1}\times _M\pi ^{*}.
\]

A point $u\in T^{*k}M$ is of the form $u=(x,y^{(1)},...,y^{(k-1)},p).$ It is
determined by a configuration $x=(x^i),$ the accelerations \newline $y^{(1)i\text{ }%
}= \displaystyle\frac{dx^i}{dt},..., y^{(k-1)i}=\displaystyle\frac 1{(k-1)!} %
\displaystyle\frac{d^{k-1}x^i}{dt^{k-1}}$ and the momenta $p=(p_i).$

All previous conditions 1$^{\circ }$-5$^{\circ }$ are satisfied.
These considerations imply the fact that the geometries of higher
order Lagrange space and higher order Hamilton space are dual.

The duality is obtained via a Legendre transformation.

For a good understanding of the important concept of duality we had to make a short introduction to the geometrical theory of Lagrange
and Finsler spaces of order $k$ and then continue with the main subject of
the book, the geometry of Hamilton and Cartan spaces of order $k$.

The Lagrange spaces of order $k$ are defined as the pairs $L^{(k)n}=(M,L)$,
where $L$ is a regular Lagrangian of order $k$. By means of variational calculus the integral of action
$ I(c)=\int\limits_0^1L\left( x(t),\displaystyle\frac{dx}{dt}(t),...,
\displaystyle\frac 1{k!}\displaystyle\frac{d^kt}{dt^k}\right) dt$ gives the
Euler-Lagrange equations and the Craig-Synge equations. The last equations
determine a canonical $k$-semispray $S$. The geometry of the space $L^{(k)n}$
can be developed by means of the fundamental function $L$, of the
fundamental tensor $g_{ij}=\displaystyle\frac 12\displaystyle\frac{\partial
^2L}{\partial y^{(k)i}\partial y^{(k)j}}$ and of the canonical $k$-semispray $S$.
The lifting of the previous geometrical edifice to the total space $T^kM$
will give us a metrical almost contact structure, canonically related to the
Lagrange space of order $k$, $L^{(k)n}$. Of course, this structure involves
the geometry of the space $L^{(k)n}$.

An important problem was to find some remarkable examples of spaces $%
L^{(k)n} $, for $k>1$. By solving the problem of prolongations to $T^kM$ of a
Riemannian structures $g$ given on the base manifold $M$, we found interesting examples of Lagrange spaces of order $k$.

For the applications, one studies the notions of energy of order
$1,2,...,k$ and one proves the law of conservation for the energy
of order $k$ and a N\"other type theorem.

The spaces $L^{(k)n}$ have two important particular cases.

The Finsler spaces of order $k$, $F^{(k)n}$, obtained when the fundamental
function $L$ is homogeneous with respect to accelerations $y^{(1)}$, ..., $%
y^{(k)}$, and the Riemann spaces of order $k$, $\mathcal{R}^{(k)n}$ are the spaces $L^{(k)n}$ for which the fundamental tensor $g_{ij}$ does not
depend on the accelerations $y^{(1)}$, ..., $y^{(k)}$.

Therefore we have the following sequence of inclusions, [94]:

$$
\left\{ \mathcal{R}^{(k)n}\right\} \subset \left\{ F^{(k)n}\right\} \subset
\left\{ L^{(k)n}\right\} \subset \left\{ GL^{(k)n}\right\}.
\leqno{(*)}
$$

In the case $k=1$ this sequence admits a {\it dual}, which is
obtained via Legendre transformation. In the book [115] we have
introduced the 'dual' of the sequence (*) for the case $k=2$.

Now, the main problem for us is to define and study a {\it dual} sequence of the
inclusions (*) in the {\it dual} Hamilton space of order $k$.

A Hamilton space of order $k$ is a pair $H^{(k)n}=(M,H)$, where

$H:(x,y^{(1)},...,y^{(k-1)},p)\in T^{*k}M\rightarrow
H(x,y^{(1)},...,y^{(k-1)},p)\in \mathbf{R}$ is a regular Hamiltonian. Here, the
regularity means: the Hessian of $H$, with respect to the momenta $p_i$, is
not singular. The elements of the Hessian matrix are $g^{ij}=\displaystyle
\frac 12\displaystyle\frac{\partial ^2H}{\partial p_i\partial p_j}$. Thus $H$
is called the fundamental function and $g^{ij}$ the fundamental tensor of
the space $H^{(k)n}$. The geometry of the space $H^{(k)n}$ can be based on these two geometrical object fields: $H$ and $g^{ij}.$
In the case when $ g^{ij}=g^{ij}(x)$ we have a particular class of Hamilton spaces $\mathcal{R}
^{*(k)n}$ called Riemannian. If the fundamental function $H$ is $2k$
-homogeneous on the fibres of the bundle $T^{*k}M$, the spaces $H^{(k)n}$
are called Cartan spaces of order $k$ and denoted $\mathcal{C}^{(k)n}$.

Finally, a pair $GH^{(k)n}=(M,g^{ij})$, where $g^{ij}(x,y^{\left( 1\right)
},...,y^{\left( k-1\right) },p)$ is a symmetric, nonsingular, distinguished
tensor field which is called a generalized Hamilton space of order $k$.

Consequently, we obtain the sequence of inclusions
$$
\left\{ \mathcal{R}^{*(k)n}\right\} \subset \left\{ \mathcal{C}
^{*(k)n}\right\} \subset \left\{ H^{(k)n}\right\} \subset \left\{
GH^{(k)n}\right\}.
\leqno{(**)}
$$
This is the {\it dual} sequence of the sequence (*) via the
Legendre transformation.

The main goal of this book is to study the classes of spaces from the
sequence (**). The chapters 4-11 of the book are devoted to this subject.

Therefore we begin with the geometry of the total space $T^{*k}M$ of the
dual bundle $(T^{*k}M,\pi ^{*k},M)$ of the $k$-tangent bundle $(T^kM,\pi
^k,M)$ underline: vertical distributions; Liouville vector fields; Liouville
$1$-form $\omega =p_idx^i$; the closed $2$-form $\theta =d\omega $ which
defines a natural presymplectic structure on $T^{*k}M.$ In the chapter 5 a
new theory of variational problem for the Hamiltonian $H$ of order $k$ is
developed starting from the integral of action of $H$ defined by

\[
I(c)=\int_0^1[p_i\displaystyle\frac{dx^i}{dt}-\displaystyle\frac 12H(x, %
\displaystyle\frac{dx}{dt},...,\displaystyle\frac 1{(k-1)!}\displaystyle
\frac{d^{k-1}x}{dt^{k-1}},p)]dt
\]

One proves that: the extremal curves are the solutions of the following
Hamilton-Jacobi equations:

\[
\begin{array}{l}
\displaystyle\frac{dx^i}{dt}=\displaystyle\frac 12\displaystyle\frac{
\partial H}{\partial p_i}, \\
\\
\displaystyle\frac{dp_i}{dt}=-\displaystyle\frac 12[\displaystyle\frac{
\partial H}{\partial x^i}-\displaystyle\frac d{dt}\displaystyle\frac{
\partial H}{\partial y^{\left( 1\right) i}}+...+\left( -1\right) ^{k-1} %
\displaystyle\frac 1{\left( k-1\right) !}\displaystyle\frac{d^{k-1}}{%
dt^{k-1} }\displaystyle\frac{\partial H}{\partial y^{\left( k-1\right) i}}].
\end{array}
\]

These equations are fundamental in the whole construction of the
geometry of Hamiltonians of order $k$. They allow the introduction of the notion of
energy of order $k-1,...,1,$ $\mathcal{E}^{k-1}(H),...,\mathcal{E}^1(H)$ and
prove a law of conservation for $\mathcal{E}^{k-1}(H)$ along extremal
curves. Now we can introduce in a natural way the Jacobi-Ostrogradski momenta and the
Hamilton-Jacobi-Ostrogradski equations.

A theory of symmetries of the Hamiltonians $H$ and the N\"other type
theorems are investigated, too; a specific theory of tangent structure $J
$ and its adjoint $J^{*}$; canonical Poisson structure; the notion of dual
semispray, which can be defined only by $k\geq 2$; nonlinear connection $N$;
the dual coefficients of $N$; the almost product structure $\Bbb{P}$,
almost contact structure $\Bbb{F}$ and Riemannian structure $\Bbb{G}$, are studied.

We pay special attention to the theory of $N$-linear connections; curvatures
and torsions; parallelism and structures equations.

Chapter 8 is devoted to the main subject from the book: Hamilton
spaces of order $k,$ $H^{(k)n}=(M,H(x,y^{(1)},..,y^{(k-1)},p)).$
To begin with, we prove the existence of these spaces and the
existence of a natural presymplectic structure, as well as of a
natural Poisson structure. Using the Legendre mapping from a
Lagrange space of order $k,$ $L^{(k)n}=(M,L)$ to the Hamilton
space of order $k,$ $H^{(k)n}=(M,H)$ one proves that there is a
local diffeomorphism between these spaces. A direct consequence of
previous results one can determine some important geometric object
fields on the Hamilton spaces $H^{(k)n},$ namely: the canonical
nonlinear connection, the N-linear metrical connection given by
generalized Christoffel symbols. Evidently, the structure
equations, curvatures and torsions of above mentioned connections
are pointed out. The Hamilton-Jacobi equations and an example from
Electrodynamics end this chapter.

A theory of subspaces
$\stackrel{\vee}{H}^{(k)m}=(\!\!\stackrel{\vee}{M},
\stackrel{\vee}{H})$ in the Hamilton spaces $H^{(k)n}=(M,H)$
appears for the first time in this book, ch. 9. Of course, it is
absolutely necessary, especially for applications. But
$\!\!\stackrel{\vee}{M}$ being a submanifold in the manifold $M$,
the immersion $i:$ $\!\!\stackrel{\vee}{M}\rightarrow M$ does not
implies automatically an immersion of
$T^{*k}\!\!\stackrel{\vee}{M}$ into the dual manifold $T^{*k}M$.

So, by means of an immersion of the cotangent bundle $T^{*}\!\!\stackrel{\vee}{M}$
into $T^{*}M$ we construct $T^{*k}\!\!\stackrel{\vee}{M}$ as an immersed
submanifold of the manifold $T^{*k}M$.

The Hamilton space $H^{(k)n}=(M,H)$ induces Hamilton subspaces $\stackrel{
\vee}{H}^{(k)m}=(\!\!\stackrel{\vee}{M},\stackrel{\vee}{H})$. So, we study
the intrinsic geometrical object fields on $\stackrel{\vee}{H}^{(k)m}$ and
the induced geometrical object fields, as well as the relations between
them. These problems are studied using the method of moving frame -
suggested by the theory of subspaces in Lagrange spaces of order $k$.
The Gauss-Weingarten formulae and the Gauss-Codazzi equations are important results.

In chapter 10 we investigate the notion of Cartan space of order
$k\geq 1$ as dual of that of Finsler space of same order. We point
out the canonical linear connection, N-metrical connection,
structure equations, the fundamental equations of Hamilton Jacobi
and the Riemannian almost contact model of these spaces.

The last chapter, Ch. 11, is devoted to the Generalized Hamilton spaces $GH^{(k)n}$,
Generalized Cartan spaces $GC^{(k)n}$ and applications in the Hamiltonian Relativistic Optics.

Now, some remarks. The book can be divided in three parts: the
Lagrange geometry of order $k$, presented in the first three
chapters, the geometrical theory of the dual manifolds $T^{*k}M$ -
chapters 4-7 and the geometry of Hamilton spaces of order $k$ and
their subspaces, contained in the last four chapters. They are
studied directly and as 'dual' geometry, via Legendre
transformation. More details for Lagrange geometry of order $k$
can be found in the book [94]. Also, the particular case, $k=2$,
of the geometry of Hamilton spaces $H^{(k)n}$ can be found in the
book [115].

For these reasons, the book is accessible for readers from graduate
students to researchers in Mathematics, Mechanics, Physics, Biology,
Informatics, etc.

\textbf{Acknowledgments.} I would like to express my gratitude to P.L. Antonelli, M. Anastasiei, A. Bejancu, M. Matsumoto, R.M. Santilli, P.S. Morey Jr. and Izu Vaisman for their continuous moral support and numerous valuable suggestions, as well to my collaborators: H. Shimada, D. Hrimiuc and V.S. Sab\u au, for the realization of the joint book [115].

Special thanks to I. Bucataru, M. Roman (University of Ia\c si), to Ph. D. Assistants L. Popescu, F. Munteanu and to Professor P. Stavre
(University of Craiova) who gave the manuscript a meticulous reading and typeset the manuscript into its final excellent form.

Finally, I should like to thanks to Kluwer Academic Publishers for co-operation.

\chapter{Geometry of the $k$-Tangent Bundle $T^{k}M$}

\markboth{\it{THE GEOMETRY OF HIGHER-ORDER HAMILTON SPACES\ \ \ \ \ }}{\it{Geometry of the $k$-Tangent Bundle} $T^{k}M$}

\pagenumbering{arabic}

The notion of $k$- tangent bundle (or $k$-accelerations bundle or $k$-
osculator bundle), ($T^kM,\pi ^k,M)$ is sufficiently known. It was presented
in the book [115].

The manifold $T^kM$ carries some geometrical object fields as the vertical
distributions $V_1,...,V_k,$ the Liouville independent vector fields $%
\stackrel{1}{\Gamma },...,\stackrel{k}{\Gamma },$ with the properties $%
\stackrel{k}{\Gamma }$ belongs to $V_1$, $\stackrel{k-1}{\Gamma }$ belongs
to $V_2$,... and $\stackrel{1}{\Gamma }$ belongs to the distribution $V_k.$
On $T^kM$ is defined a $k$- tangent structure $J$ which maps $\stackrel{1}{
\Gamma }$ on $\stackrel{2}{\Gamma }$, $\stackrel{2}{\Gamma }$ on $\stackrel{%
3 }{\Gamma }$ , $\stackrel{k-1}{\Gamma }$ on $\stackrel{k}{\Gamma
}$ and $J\stackrel{k}{\Gamma }=0.$

Besides these fundamental notion on $T^kM$ we can introduce new concepts
as the $k$- semisprays $S$, nonlinear connections $N$ and the $N$- linear
connections $D.$ But for $D$ we can get the curvatures, torsions, structure
equations, geodesics, etc. The $k$- semispray $S$ is defined by the
conditions $JS=\stackrel{k}{\Gamma }.$ It is important to remark that $S$ is
used for introducing those notions as nonlinear connection, or $N$-
linear connection.

Concluding the geometry of $k$-accelerations bundle is basic for a
geometrical theory of higher order Lagrange spaces or higher order Finsler
spaces. In this book we need it for a theory of duality between higher order
Lagrange spaces and higher order Hamilton spaces.

\section{The Category of $k$-Accelerations Bundles}

Let $M$ be a real $n$-dimensional manifolds of $C^\infty $ class  and $ (T^kM,\pi ^k,M$ $)$ its bundle of accelerations of order $k$. It can be identified with the $k$-osculator bundles [94] or with the tangent bundle of order $k$. In the case $k=1$, $(T^1M,\pi ^1,M)$ is the tangent bundle of the manifold $M$.

A point $u\in T^kM$ will be written as $u=(x,y^{(1)},...,y^{(k)})$ and $\pi^k(u)=x,$ $x\in M.$ The canonical coordinates of $u$ are $(x^i,y^{\left(1\right) i},...,y^{\left( k\right) i})$, $i=\overline{1,n}$, $n=\dim M.$

These coordinates have a geometrical meaning. If $c:I\rightarrow M\,$ is a
differentiable curve, $c(0)=x_0\in M,$ and $Im$ $c\subset U,$ U being a local
chart of the base manifold $M,$ and the mapping $c:I\rightarrow M$ is
represented by $x^i=x^i(t),$ $t\in I,$ then the osculating space of the
curve $c,$ in the point $x_0^i=x^i(0)$ is characterized by the set of numbers:
\begin{equation}
x_0^i=x^i(0),\ y_0^{(1)i}=\displaystyle\frac{dx^i}{dt}\left( 0\right) ,...,\
y^{\left( k\right) i}=\displaystyle\frac 1{k!}\displaystyle\frac{d^kx^i}{%
dt^k }\left( 0\right)  \tag{1.1.1}
\end{equation}

Thus, the formulas (1.1.1) give us the canonical coordinates of a point $%
u_0=(x_0,y_0^{(1)},...,y_0^{(k)})$ of the domain of the local chart $(\pi
^k)^{-1}(U)\subset T^kM.$

Starting from (1.1.1) it is not difficult to see which are the changing rules of the
local coordinates on $T^kM:(x^i,y^{\left( 1\right) i},...,y^{\left( k\right)
i})\rightarrow (\widetilde{x}^i,\widetilde{y}^{\left( 1\right) i},...,
\widetilde{y}^{\left( k\right) i}).$

We deduce:
\begin{equation}
\left\{
\begin{array}{l}
\widetilde{x}^i=\widetilde{x}^i(x^1,...,x^n),\ rank\left\| \displaystyle
\frac{\partial \widetilde{x}^i}{\partial x^j}\right\| =n \\
\\
\widetilde{y}^{(1)i}=\displaystyle\frac{\partial \widetilde{x}^i}{\partial
x^j}y^{(1)j} \\
\\
2\widetilde{y}^{(2)i}=\displaystyle\frac{\partial \widetilde{y}^{(1)i}}{
\partial x^j}y^{(1)j}+2\displaystyle\frac{\partial \widetilde{y}^{\left(
1\right) i}}{\partial y^{\left( 1\right) j}}y^{\left( 2\right) j} \\
................................................. \\
k\widetilde{y}^{(k)i}=\displaystyle\frac{\partial \widetilde{y}^{(k-1)i}}{
\partial x^j}y^{(1)j}+2\displaystyle\frac{\partial \widetilde{y}^{\left(
k-1\right) i}}{\partial y^{\left( 1\right) j}}y^{\left( 2\right) j}+...+k %
\displaystyle\frac{\partial \widetilde{y}^{\left( k-1\right) i}}{\partial
y^{\left( k-1\right) j}}y^{\left( k\right) j}.
\end{array}
\right.  \tag{1.1.2}
\end{equation}

But we must remark the following identities:
\begin{equation}
\displaystyle\frac{\partial \widetilde{y}^{(\alpha )i}}{\partial x^j}= %
\displaystyle\frac{\partial \widetilde{y}^{(\alpha +1)i}}{\partial y^{\left(
1\right) j}}=...=\displaystyle\frac{\partial \widetilde{y}^{(k)i}}{\partial
y^{\left( k-\alpha \right) j}},\ (\alpha =0,...,k-1;\ y^{\left( 0\right)
}=x).  \tag{1.1.3}
\end{equation}

In the following $T^0M$ is canonically identified to $M.$ Sometimes we
employ the notations $y^{(0)i}=x^i.$ The projections:
\[
\pi _l^k:(x,y^{\left( 1\right) },...,y^{\left( k\right) })\in
T^kM\rightarrow (x,y^{\left( 1\right) },...,y^{\left( l\right) })\in T^lM,\
(0\leq l<k)
\]
are submersions. Clearly $\pi _0^k=\pi ^k.$

A section $S:M\rightarrow T^kM$ of the projection $\pi ^k$ is a
differentiable mapping with the property $\pi ^k\circ S=1_M.$ It is a local
section if $\pi ^k\circ S\left| _U\right. =1_U.$ Of course a section $%
S:M\rightarrow T^kM$ along a curve $c:I\rightarrow M$ has the property $\pi
^k(S(c))=c.$

If $c:I\rightarrow M$ is locally represent on $U\subset M$ by $x^i=x^i(t)$
then the mapping $\widetilde{c}:I\rightarrow T^kM$ given by:
\begin{equation}
x^i=x^i(t),\ \ y^{\left( 1\right) i}=\displaystyle\frac 1{1!}\displaystyle
\frac{dx^i}{dt}(t),...,\ \ y^{\left( k\right) i}=\displaystyle\frac 1{k!} %
\displaystyle\frac{dx^{\left( k\right) i}}{dt^k}(t),\ t\in I  \tag{1.1.4}
\end{equation}
is the extension of order $k$ of $c$. Of course $\pi ^k\circ \widetilde{c}
=c. $ So $\widetilde{c}$ is a section of $\pi ^k$ along $c.$

More general, if $V$ is a vector field on the domain of a chart $U$ and $%
c:I\rightarrow U$ is a curve, then the mapping

$S_V:c\rightarrow (\pi ^k)^{-1}(U)\subset T^kM,$ defined by
\begin{equation}
S_V:x^i=x^i(t),\ y^{\left( 1\right) i}=V^i(x(t)),...,\ y^{(k)i}=%
\displaystyle \frac 1{k!}\displaystyle\frac{d^{k-1}V^i(x\left( t\right) )}{%
dt^{k-1}},\ t\in I  \tag{1.1.5}
\end{equation}
is a section of the projection $\pi ^k$ along curve $c.$

Of course the notion of the section of $\pi _l^k$ along $T^lM$ can be
defined, as in the previous case.

The following property hold:

\begin{teo}
If the differentiable manifold $M$ is paracompact, then $T^kM$ is a
paracompact manifold.
\end{teo}

We can see, that
\[
T^k:Man\rightarrow Man,
\]
where $Man$ is the category of differentiable manifolds, is a covariant
functor.

Indeed we define:

$T^k:M\in ObMan\rightarrow T^kM\in ObMan$ and

$T^k:\{f:M\rightarrow M^{\prime }\}\rightarrow \{T^kf:T^kM\rightarrow T^kM^{\prime }\}$ as follows: \\
if $f(x)$ in the local coordinate of $M$ is given by $x^{i^{\prime
}}=x^{i^{\prime }}(x^1,...,x^n),$ \break
$\ i^{^{\prime }},j^{\prime }=1,...,m=\dim
M^{\prime }$, then the morpfism $T^kf:T^kM\rightarrow T^kM^{\prime }$ is
defined by:
\begin{equation}
\left\{
\begin{array}{l}
x^{i^{\prime }}=x^{i^{\prime }}(x^1,...,x^n), \\
\\
y^{(1)i^{\prime }}=\displaystyle\frac{\partial x^{i^{\prime }}}{\partial x^j}
y^{(1)j}, \\
\\
2y^{(2)i^{\prime }}=\displaystyle\frac{\partial y^{(1)i^{\prime }}}{\partial
x^j}y^{(1)j}+2\displaystyle\frac{\partial y^{\left( 1\right) i^{\prime }}}{
\partial y^{\left( 1\right) j}}y^{\left( 2\right) j}, \\

................................................. \\
ky^{(k)i^{\prime }}=\displaystyle\frac{\partial y^{(k-1)i^{\prime }}}{
\partial x^j}y^{(1)j}+2\displaystyle\frac{\partial y^{\left( k-1\right)
i^{\prime }}}{\partial y^{\left( 1\right) j}}y^{\left( 2\right) j}+...+k %
\displaystyle\frac{\partial y^{\left( k-1\right) i^{\prime }}}{\partial
y^{\left( k-1\right) j}}y^{\left( k\right) j}.
\end{array}
\right.  \tag{1.1.6}
\end{equation}
Remarking that
\begin{equation}
\displaystyle\frac{\partial y^{(\alpha )i^{\prime }}}{\partial x^j}= %
\displaystyle\frac{\partial y^{(\alpha +1)i^{\prime }}}{\partial y^{\left(
1\right) j}}=...=\displaystyle\frac{\partial y^{(k)i^{\prime }}}{\partial
y^{\left( k-\alpha \right) j}},\ (\alpha =0,...,k-1;\ y^{\left( 0\right)
}=x),  \tag{1.1.7}
\end{equation}
we can prove without difficulties that $T^k$ is a covariant functor.

\section{Liouville Vector Fields. $k$-Semisprays}

A local coordinate changing (1.1.2) transforms the natural basis of the
tangent space $T_u(T^kM)$ by the following rule:
\begin{equation}
\left\{
\begin{array}{l}
\displaystyle\frac \partial {\partial x^i}=\displaystyle\frac{\partial
\widetilde{x}^j}{\partial x^i}\displaystyle\frac \partial {\partial
\widetilde{x}^j}+\displaystyle\frac{\partial \widetilde{y}^{(1)j}}{\partial
x^i}\displaystyle\frac \partial {\partial \widetilde{y}^{(1)j}}+...+ %
\displaystyle\frac{\partial \widetilde{y}^{(k)j}}{\partial x^i}\displaystyle %
\frac \partial {\partial \widetilde{y}^{(k)j}} \\
\\
\displaystyle\frac \partial {\partial y^{(1)i}}=\qquad \qquad \displaystyle
\frac{\partial \widetilde{y}^{(1)j}}{\partial y^{(1)i}}\displaystyle\frac
\partial {\partial \widetilde{y}^{(1)j}}+...+\displaystyle\frac{\partial
\widetilde{y}^{(k)j}}{\partial y^{(1)i}}\displaystyle\frac \partial
{\partial \widetilde{y}^{(k)j}} \\
............................................................... \\
\displaystyle\frac \partial {\partial y^{(k)i}}=\qquad \qquad \qquad \qquad
\qquad \qquad \displaystyle\frac{\partial \widetilde{y}^{(k)j}}{\partial
y^{(k)i}}\displaystyle\frac \partial {\partial \widetilde{y}^{(k)j}},
\end{array}
\right. \tag{1.2.1}
\end{equation}
calculated at the point $u\in T^kM.$

These formulas imply the transformation of the natural cobasis at the point $u \in T^kM$ by the rule:
\begin{equation}
\left\{
\begin{array}{l}
d\widetilde{x}^i=\displaystyle\frac{\partial \widetilde{x}^i}{\partial x^j}
dx^j, \\
\\
d\widetilde{y}^{\left( 1\right) i}=\displaystyle\frac{\partial \widetilde{y}
^{(1)i}}{\partial x^j}dx^j+\displaystyle\frac{\partial \widetilde{y}^{(1)i}}{
\partial y^{(1)j}}dy^{\left( 1\right) j}, \\
............................................................... \\
d\widetilde{y}^{\left( k\right) i}=\displaystyle\frac{\partial \widetilde{y}
^{(k)i}}{\partial x^j}dx^j+\displaystyle\frac{\partial \widetilde{y}^{(k)i}}{
\partial y^{(1)j}}dy^{\left( 1\right) i}+...+\displaystyle\frac{\partial
\widetilde{y}^{(k)i}}{\partial y^{(k)j}}dy^{\left( k\right) j}.
\end{array}
\right.  \tag{1.2.1'}
\end{equation}

The matrix of coefficients of second member of (1.2.1) is the Jacobian matrix
of the changing of coordinates (1.1.2). Since
\[
\displaystyle\frac{\partial \widetilde{x}^i}{\partial x^j}=\displaystyle
\frac{\partial \widetilde{y}^{(1)i}}{\partial y^{\left( 1\right) j}}=...= %
\displaystyle\frac{\partial \widetilde{y}^{(k)i}}{\partial y^{\left(
k\right) j}},\
\]
it follows.

\begin{teo}
If the number $k$ is odd, then the manifold $T^kM$ is orientable.
\end{teo}

Also the formulae (1.2.1), (1.2.1') allow to determine some important geometric
object fields on the total space of accelerations bundle $T^kM.$

The distribution $V_1:u\in T^kM\rightarrow V_{1,u}\subset T_u(T^kM)$
generated by the tangent vectors $\{\displaystyle\frac \partial {\partial
y^{\left( 1\right) i}},...,\displaystyle\frac \partial {\partial y^{\left(
k\right) i}}\}_u,$ $\forall u \in T^kM $ is the vertical distribution on the bundle $T^kM.$ Its
dimension is $kn.$

$V_2:u\in T^kM\rightarrow V_{2,u}\subset T_u(T^kM)$ generated by $\{ %
\displaystyle\frac \partial {\partial y^{\left( 2\right) i}},..., %
\displaystyle\frac \partial {\partial y^{\left( k\right) i}}\}_u,$ $\forall u \in T^kM $ is a
subdistribution of $V_1$ of local dimension $(k-1)n,$ ...,

$V_k:u\in T^kM\rightarrow V_{k,u}\subset T_u(T^kM)$ generated by $\{ %
\displaystyle\frac \partial {\partial y^{\left( k\right) i}}\}_u,$ $\forall u \in T^kM $ is a subdistribution of dimension $n$ of the distribution $V_{k-1}$. All these subdistributions are integrable and the following sequence holds:
\[
V_1\supset V_2\supset ...\supset V_k
\]
Using again (1.2.1) we deduce:

\begin{teo}
The following operators in the algebra of functions $\mathcal{F}(T^kM):$
\begin{equation}
\begin{array}{l}
\stackrel{1}{\Gamma }=y^{\left( 1\right) i}\displaystyle\frac \partial
{\partial y^{\left( k\right) i}}, \\
\stackrel{2}{\Gamma }=y^{\left( 1\right) i}\displaystyle\frac \partial
{\partial y^{\left( k-1\right) i}}+2y^{\left( 2\right) i}\displaystyle\frac
\partial {\partial y^{\left( k\right) i}}, \\
.............................................. \\
\stackrel{k}{\Gamma }=y^{\left( 1\right) i}\displaystyle\frac \partial
{\partial y^{\left( 1\right) i}}+2y^{\left( 2\right) i}\displaystyle\frac
\partial {\partial y^{\left( 2\right) i}}+...+ky^{\left( k\right) i}%
\displaystyle\frac \partial {\partial y^{\left( k\right) i}}
\end{array}
\tag{1.2.2}
\end{equation}
are vector fields globally defined on $T^kM$ and independent on the manifold
$\widetilde{T^kM}=T^kM\setminus \{0\},$ $\stackrel{1}{\Gamma }$ belongs to
distribution $V_k,$ $\stackrel{2}{\Gamma }$ belongs to distribution $V_{k-1},...,\stackrel{k}{\Gamma }$ belongs to
distribution $V_1.$
\end{teo}
Taking into account (1.2.1') it is not hard to prove:

\begin{teo}
For any differentiable function $L(x,y^{(1)},...,y^{\left( k\right) })$ on
the \break manifold $\widetilde{T^kM}$, the following entries $d_0L,...,d_kL$ are $1
$-form fields on $\widetilde{T^kM}:$
\begin{equation}
\begin{array}{l}
d_0L=\displaystyle\frac{\partial L}{\partial y^{\left( k\right) i}}dx^i, \\
\\
d_1L=\displaystyle\frac{\partial L}{\partial y^{\left( k-1\right) i}}dx^i+%
\displaystyle\frac{\partial L}{\partial y^{\left( k\right) i}}dy^{\left(
1\right) i}, \\
.............................................. \\
d_kL=\displaystyle\frac{\partial L}{\partial x^i}dx^i+\displaystyle\frac{%
\partial L}{\partial y^{\left( 1\right) i}}dy^{\left( 1\right) i}+...+%
\displaystyle\frac{\partial L}{\partial y^{\left( k\right) i}}dy^{\left(
k\right) i}.
\end{array}
\tag{1.2.2'}
\end{equation}
\end{teo}

Evidently, $d_kL=dL$ is the differential of the function $L.$ But:

$d_0L$ vanish on the distribution $V_1$,

$d_1L$ vanish on the distribution $V_2$,

.................................................

$d_{k-1}L$ vanish on the distribution $V_k.$\\
In applications we shall use also the following nonlinear operator
\begin{equation}
\Gamma =y^{\left( 1\right) i}\displaystyle\frac \partial {\partial
x^i}+2y^{\left( 2\right) i}\displaystyle\frac \partial {\partial y^{\left(
1\right) i}}+...+ky^{\left( k\right) i}\displaystyle\frac \partial {\partial
y^{(k-1)i}}.  \tag{1.2.3}
\end{equation}

A $k$-tangent structure $J$ on $T^kM$ is defined as usual by the following \newline
$\mathcal{F}(T^kM)$-linear mapping $J:\mathcal{X}(T^kM)\rightarrow \mathcal{%
X }(T^kM):$
\begin{equation}
\begin{array}{l}
J(\displaystyle\frac \partial {\partial x^i})=\displaystyle\frac \partial
{\partial y^{\left( 1\right) i}};\quad J(\displaystyle\frac \partial
{\partial y^{\left( 1\right) i}})=\displaystyle\frac \partial {\partial
y^{\left( 2\right) i}};...; \\
\\
J(\displaystyle\frac \partial {\partial y^{\left( k-1\right) i}})= %
\displaystyle\frac \partial {\partial y^{\left( k\right) i}};\quad J( %
\displaystyle\frac \partial {\partial y^{\left( k\right) i}})=0.
\end{array}
\quad  \tag{1.2.4}
\end{equation}
$J$ is globally defined on $T^kM.$ It is a tensor fields of type (1,1) on $%
T^kM,$ locally expressed by
\begin{equation}
J=\displaystyle\frac \partial {\partial y^{\left( 1\right) i}}\otimes dx^i+ %
\displaystyle\frac \partial {\partial y^{\left( 2\right) i}}\otimes
dy^{\left( 1\right) i}+...+\displaystyle\frac \partial {\partial y^{\left(
k\right) i}}\otimes dy^{\left( k-1\right) i}.  \tag{1.2.4'}
\end{equation}
The last form of $J$ implies that $J$ is an integrable structure. The $k$
-tangent structure $J$ has the following properties:

1$^{\circ }$. $ImJ=V_1,\quad KerJ=V_k.$

2$^{\circ }$. $rankJ=kn.$

3$^{\circ }$. $J\stackrel{k}{\Gamma }=\stackrel{k-1}{\Gamma },...,$ $J
\stackrel{2}{\Gamma }=\stackrel{1}{\Gamma },$ $J\stackrel{1}{\Gamma }=0.$

4$^{\circ }$. $J\circ ...\circ J=0$, $(k+1$ factors).

A $k$-semispray on the manifold $T^kM$ is a vector field $S$ on $T^kM$ with
the property
\begin{equation}
JS=\stackrel{k}{\Gamma }.  \tag{1.2.5}
\end{equation}
The notion of local $k$-semispray is obvious.

Any $k$-semispray $S$ can be uniquely written in the form
\begin{equation}
S=y^{\left( 1\right) i}\displaystyle\frac \partial {\partial
x^i}+...+ky^{\left( k\right) i}\displaystyle\frac \partial {\partial
y^{(k-1)i}}-(k+1)G^i(x,y^{\left( 1\right) },...,y^{\left( k\right) }) %
\displaystyle\frac \partial {\partial y^{\left( k\right) i}},  \tag{1.2.6}
\end{equation}
or shortly
\begin{equation}
S=\Gamma -(k+1)G^i(x,y^{\left( 1\right) },...,y^{\left( k\right) }) %
\displaystyle\frac \partial {\partial y^{\left( k\right) i}}.  \tag{1.2.6'}
\end{equation}

The set of functions $G^i$ is the set of \textit{coefficient}s of $S.$ With
respect to (1.1.2), $G^i$ are transformed as following:
\begin{equation}
(k+1)\widetilde{G}^i=(k+1)G^j\displaystyle\frac{\partial \widetilde{x}^i}{
\partial x^j}-\Gamma \widetilde{y}^{\left( k\right) i}.  \tag{1.2.6''}
\end{equation}

A curve $c:I\rightarrow M$ is called a $k$\textit{-path} on $M$ with respect
to a $k$-semispray $S$ if its extension $\widetilde{c}$ to $T^kM$ is an integral curve of $S.$
If $S$ is given in the form (1.2.6), then its $k$-paths are characterized by the system of differential equations:
\begin{equation}
\displaystyle\frac{d^{k+1}x^i}{dt^{k+1}}+(k+1)!G^i(x,\displaystyle\frac{dx}{
dt},...,\displaystyle\frac{d^kx}{dt^k})=0.  \tag{1.2.7}
\end{equation}
Indeed, the solutions curves of $S$ are given by the system of ordinary differential equations
\[
\displaystyle\frac{dx^i}{dt}=y^{\left( 1\right) i},\quad \displaystyle\frac{
dy^{\left( 1\right) i}}{dt}=2y^{\left( 2\right) i},...,\quad \displaystyle
\frac{dy^{\left( k-1\right) i}}{dt}=ky^{\left( k\right) i},
\]
\[
\displaystyle\frac{dy^{\left( k\right) i}}{dt}=-(k+1)G^i(x,y^{\left(
1\right) },...,y^{\left( k\right) }).
\]
If we eliminate $y^{\left( 1\right) i},...,y^{\left( k\right) i}$ we obtain
(1.2.7).

The previous theory will be used in the geometry of higher order Lagrange
spaces, which is based on the regular Lagrangians of higher order.

Now, let us considered the \textit{adjunct} $k$-tangent structure $J^{*}.$
It is the endomorphism of the module $\mathcal{X}^{*}(T^kM),$ defined by:
\begin{equation}
\left.
\begin{array}{l}
J^{*}(dy^{\left( k\right) i})=dy^{\left( k-1\right) i}=dy^{\left( k\right)
i}\circ J, \\
\\
J^{*}(dy^{\left( k-1\right) i})=dy^{\left( k-2\right) i}=dy^{\left(
k-1\right) i}\circ J, \\
................................................... \\
J^{*}(dy^{\left( 2\right) i})=dy^{\left( 1\right) i}=dy^{\left( 2\right)
i}\circ J, \\
\\
J^{*}(dy^{\left( 1\right) i})=dx^i=dy^{\left( 1\right) i}\circ J, \\
\\
J^{*}(dx^i)=0.
\end{array}
\right.  \tag{1.2.8}
\end{equation}
By using the formula (1.2.1') is not difficult to prove:

\begin{teo}
$J^{*}$ is globally defined on $T^kM.$
\end{teo}
If $\omega \in X^{*}(T^kM)$ is given by
\[
\omega =\stackrel{\left( 0\right) }{\omega }_idx^i+\stackrel{\left( 1\right)
}{\omega }_idy^{\left( 1\right) i}+...+\stackrel{\left( k\right) }{\omega }
_idy^{\left( k\right) i},
\]
then
\[
J^{*}\omega =\stackrel{\left( 1\right) }{\omega }_idx^i+...+\stackrel{\left(
k\right) }{\omega }_idy^{\left( k-1\right) i}.
\]
We put $J^{*}f=f$ for any function $f\in \mathcal{F}(T^kM)$ and observe that $J^{*}$ is a
tensor field, of type $(1,1)$ on $T^kM.$ Namely, we have
\begin{equation}
J^{*}=dy^{\left( k-1\right) i}\otimes \displaystyle\frac \partial {\partial
y^{\left( k\right) i}}+dy^{\left( k-2\right) i}\otimes \displaystyle\frac
\partial {\partial y^{\left( k-1\right) i}}+...+dx^i\otimes \displaystyle %
\frac \partial {\partial y^{\left( 1\right) i}}. \tag{1.2.9}
\end{equation}

The rank $\left\| J^{*}\right\| =kn.$ $J^{*}$ can be extended to an
endomorphism of the exterior algebra $\wedge (T^kM),$ by putting
\begin{equation}
(J^{*}\omega )(X_1,...,X_p)=\omega (JX_1,...,JX_p),\forall \omega \in \wedge
^p(T^kM).  \tag{1.2.10}
\end{equation}

The existence of $J^{*}$ allow to introduce the so called \textit{vertical}
differential operators in the exterior algebra $\wedge (T^kM).$

Indeed, let us consider the operators of differentiation $d_k:$
\[
d_k=d=\displaystyle\frac \partial {\partial x^i}dx^i+\displaystyle\frac
\partial {\partial y^{\left( 1\right) i}}dy^{(1)i}+...+\displaystyle\frac
\partial {\partial y^{\left( k\right) i}}dy^{\left( k\right) i}.
\]
We get:
\begin{equation}
\left.
\begin{array}{l}
J^{*}d_k=d_{k-1}=\qquad \displaystyle\frac \partial {\partial y^{\left(
1\right) i}}dx^i+...+\displaystyle\frac \partial {\partial y^{\left(
k\right) i}}dy^{\left( k-1\right) i}, \\
...........................................................................
\\
J^{*}d_1=d_0=\qquad \qquad \qquad \qquad \qquad \displaystyle\frac \partial
{\partial y^{\left( k\right) i}}dx^i.
\end{array}
\right.  \tag{1.2.11}
\end{equation}

It is not difficult to prove that the operators $d_k,d_{k-1},...,d_0$ do
not depend on the changing of coordinates on the manifold $T^kM.$ If $%
L(x,y^{\left( 1\right) },...,y^{\left( k\right) })$ is a 1-form function on $T^kM,$
then $d_kL,d_{k-1}L,...,d_0L$ are differentiable given by (1.2.2').

But $d_0,d_1,...,d_k$ can be extended to the exterior algebra $\wedge (TM)$
giving them restrictions to $\mathcal{F}(T^kM)$ and $\wedge ^1(T^kM).$ So,
we will take $d_0L,...,d_kL$ expressed in (1.2.2') and
\begin{equation}
d_\alpha (dy^{(\beta )i})=0,(\alpha ,\beta =0,1,...k;y^{\left( 0\right) }=x).
\tag{1.2.12}
\end{equation}
Consequently $d_0,...,d_k$ are the antiderivations of degree $1$ in the exterior algebra $\wedge (T^kM).$
For instance, if $\omega \in \wedge ^1(T^kM)$ and it is locally express by
\[
\omega =\stackrel{\left( 0\right) }{\omega }_idx^i+\stackrel{\left( 1\right)
}{\omega }_idy^{\left( 1\right) i}+...+\stackrel{\left( k\right) }{\omega }
_idy^{\left( k\right) i}
\]
then
\[
d_\alpha \omega =\underset{\beta =0}{\stackrel{k}{\sum }}d_\alpha
\stackrel{\left( \beta \right) }{\omega }_i\wedge dy^{\left( \beta
\right) i},\ (\alpha =0,...,k).
\]
It is not so difficult to see that the following properties hold:
\begin{equation}
d_\alpha \circ d_\alpha =0,\ (\alpha =0,...,k).  \tag{1.2.13}
\end{equation}
In the case $k=1$, $\displaystyle\frac 12d_0L=\displaystyle\frac 12 %
\displaystyle\frac{\partial L}{\partial y^i}dx^i=p_idx^i$ is the Liouville $%
1 $- form and $\displaystyle\frac 12(d_1\circ d_2)L=dp_i\wedge dx^i$ is the
symplectic structure on $TM.$

\section{Nonlinear Connections}

The notion of nonlinear connection is also known, [94]. A subbundle $HT^kM$
of the tangent bundle $(TT^kM,d\pi ^k,T^kM)$ which is supplementary to the
vertical subbundle $V_1T^kM:$
\[
TT^kM=HT^kM\oplus V_1T^kM
\]
is called a \textit{nonlinear connection}.

The fibres of $HT^kM,$ determine a horizontal distribution $N:u\in
T^kM\rightarrow N_u=H_uT^kM\subset T_uT^kM,$ $\forall u \in T^kM,$ supplementary to the vertical
distribution $V_1,$ i.e:
\begin{equation}
T_uT^kM=N_u\oplus V_{1,u},\quad \forall u\in T^kM.  \tag{1.3.1}
\end{equation}

If the base manifold $M$ is paracompact then on $T^kM$ there exist nonlinear connections.
The dimension of the horizontal distribution $N$ is $n=\dim M.$

Consider a nonlinear connection $N$ on $T^kM$ and denote by $h$ and $v$ the
horizontal and vertical projectors with respect to the distributions $N$ and
$V_1:$
\[
h+v=I,h^2=h,v^2=v,hv=vh=0.
\]

As usual we denote
\[
X^H=hX,\ X^V=vX,\ \forall X\in \mathcal{X}(T^kM).
\]

An horizontal lift, with respect to $N$ is a $\mathcal{F}(M)-$linear mapping
$l_h:\mathcal{X}(M)\rightarrow \mathcal{X}(T^kM)$ which has the properties:
\[
v\circ l_h=0,d\pi ^k\circ l_h={\rm I_d}
\]

There exists an unique local basis adapted to the horizontal distribution $N.$
It is given by
\begin{equation}
\displaystyle\frac \delta {\delta x^i}=l_h(\displaystyle\frac \partial
{\partial x^i}),(i=1,...,n).  \tag{1.3.3}
\end{equation}

The linearly independent vector fields $\displaystyle\frac \delta {\delta
x^i},(i=1,...,n) $ can be uniquely written in the form
\begin{equation}
\displaystyle\frac \delta {\delta x^i}=\displaystyle\frac \partial
{\partial x^i}-\underset{(1)}{N_i^j}\displaystyle\frac \partial
{\partial y^{\left( 1\right)
j}}-...-\underset{(k)}{N_i^j}\displaystyle\frac \partial {\partial
y^{\left( k\right) j}}.  \tag{1.3.4}
\end{equation}

The system of functions $(\underset{(1)}{N_i^j},...,\underset{(k)}{N_i^j} )$ gives
the\textit{\ coefficients} of the nonlinear connection $N.$

We remark that:

\quad 1) For $X=X^i(x)\displaystyle\frac \partial {\partial x^i}\in \mathcal{%
\ X}(M)$ the horizontal lifts is $l_h$ $X=X^i\displaystyle\frac \delta
{\delta x^i}.$

\quad 2) With respect to (1.1.2), we obtain
\[
\displaystyle\frac \delta {\delta \widetilde{x}^i}=\displaystyle\frac{
\partial x^j}{\partial \widetilde{x}^i}\displaystyle\frac \delta {\delta
x^j}.
\]

\quad 3) With respect to (1.1.2) the coefficients of $N$ are transformed by the rule
\begin{equation}
\left.
\begin{array}{l}
\underset{(1)}{\widetilde{N}_m^i}\displaystyle\frac{\partial
\widetilde{x} ^m}{\partial
x^j}=\underset{(1)}{N_j^m}\displaystyle\frac{\partial
\widetilde{x}^i}{\partial x^m}-\displaystyle\frac{\partial
\widetilde{y}
^{(1)i}}{\partial x^j}, \\
............................................................. \\
\underset{(k)}{\widetilde{N}_m^i}\displaystyle\frac{\partial
\widetilde{x} ^m}{\partial
x^j}=\underset{(k)}{N_j^m}\displaystyle\frac{\partial
\widetilde{x}^i}{\partial
x^m}+...+\underset{(1)}{N_j^m}\displaystyle\frac{
\partial \widetilde{y}^{(k-1)i}}{\partial x^m}-\displaystyle\frac{\partial
\widetilde{y}^{(k)i}}{\partial x^j}.
\end{array}
\right.  \tag{1.3.5}
\end{equation}

The $k$-tangent structure $J$, defined in (1.2.4) applies the horizontal
distribution $N$ into a vertical distribution $N_1\subset V_1$ of dimension $%
n$, supplementary to the distribution $V_2$. Then it applies the distribution
$N_1$ in a distribution $N_2\subset V_2,$ supplementary to the distribution $%
V_3$ and so on. Of course, we have $\dim N_0=\dim N_1=\cdots =\dim N_{k-1}=n$.

Setting $N_0=N$, we can write
\begin{equation}
N_1=J(N_0),N_2=J^2(N_0),....,N_{k-1}=J^{k-1}(N_0)  \tag{1.3.6}
\end{equation}
and we obtain the following direct decomposition:
\begin{equation}
T_uT^kM=N_{0,u}\oplus N_{1,u}\oplus ...\oplus N_{k-1,u}\oplus V_{k,u}\
,\quad \forall u\in T^kM.  \tag{1.3.7}
\end{equation}

An adapted basis to the distributions $N_0$, $N_1$, ..., $N_{k-1}$, $V_k$ at
a point $u\in T^kM$ is given by:
\begin{equation}
\displaystyle\frac \delta {\delta x^i},\displaystyle\frac \delta {\delta
y^{(1)i}},\cdots ,\displaystyle\frac \delta {\delta y^{(k)i}},  \tag{1.3.8}
\end{equation}
where
\[
\displaystyle\frac \delta {\delta y^{(1)i}}=J\left( \displaystyle\frac
\delta {\delta x^i}\right) ,\ \displaystyle\frac \delta {\delta
y^{(2)i}}=J^2\left( \displaystyle\frac \delta {\delta x^i}\right) ,\ ..., %
\displaystyle\frac \delta {\delta y^{(k)i}}=J^k\left( \displaystyle\frac
\delta {\delta x^i}\right) .
\]

Therefore, using (1.2.4) and (1.3.4) we get:
\begin{equation}
\begin{array}{l}
\displaystyle\frac \delta {\delta x^i}=\displaystyle\frac \partial
{\partial x^i}-\underset{(1)}{N_i^j}\displaystyle\frac \partial
{\partial y^{(1)j}}-\cdots
-\underset{(k)}{N_i^j}\displaystyle\frac \partial
{\partial y^{(k)j}}, \\
\\
\displaystyle\frac \delta {\delta y^{(1)i}}=\displaystyle\frac
\partial {\partial
y^{(1)i}}-\underset{(1)}{N_i^j}\displaystyle\frac \partial
{\partial y^{(2)j}}-\cdots
-\underset{(k-1)}{N_i^j}\displaystyle\frac
\partial {\partial y^{(k)j}}, \\
.................................................................. \\
\displaystyle\frac \delta {\delta y^{(k-1)i}}=\displaystyle\frac
\partial {\partial
y^{(k-1)i}}-\underset{(1)}{N_i^j}\displaystyle\frac \partial
{\partial y^{(k)j}}, \\
\\
\displaystyle\frac \delta {\delta y^{(k)i}}=\displaystyle\frac \partial
{\partial y^{(k)i}}.
\end{array}
\tag{1.3.9}
\end{equation}

With respect to (1.1.2) we have:
\begin{equation}
\displaystyle\frac \delta {\delta y^{(\alpha )i}}=\displaystyle\frac{
\partial \widetilde{x}^j}{\partial x^i}\displaystyle\frac \delta {\delta
\widetilde{y}^{(\alpha )j}},\quad (\alpha =0,...,k;y^{(0)}=x).  \tag{1.3.10}
\end{equation}

Taking into account the direct sum (1.3.1) and (1.3.7), it follows that the
vertical distribution $V_1$ at a point $u$ gives rise to the direct
decomposition:
\begin{equation}
V_{1,u}=N_{1,u}\oplus ...\oplus N_{k-1,u}\oplus V_{k,u}\ ,\quad \forall u\in
T^kM.  \tag{1.3.11}
\end{equation}

Let $h,v_1,...,v_k$ be the projectors determined by (1.3.7):
\[
\begin{array}{l}
h+v_1+\cdots +v_k=I,h^2=h,v_\alpha v_\alpha =v_\alpha ,hv_\alpha =v_\alpha
h=0,(\alpha =1,...,k), \\
\\
v_\alpha v_\beta =v_\beta v_\alpha =0,\quad (\alpha \neq \beta ;\alpha
,\beta =1,...,k).
\end{array}
\]

If we denote
\begin{equation}
X^H=hX,\ X^{V_\alpha }=v_\alpha X,\ \forall X\in \mathcal{X}(T^kM)
\tag{1.3.12}
\end{equation}
we have, uniquely,
\begin{equation}
X=X^H+X^{V_1}+...+X^{V_k}.  \tag{1.3.13}
\end{equation}

In the adapted basis (1.3.8) we can write:
\[
X^H=X^{\left( 0\right) i}\displaystyle\frac \delta {\delta
x^i},X^{V\;_\alpha }=X^{\left( \alpha \right) i}\displaystyle\frac \delta
{\delta y^{\left( \alpha \right) i}},(\alpha =1,...,k).
\]

The following properties are important in applications.

1) The distribution $N=N_0$ is integrable if, and only if
\[
\lbrack X^H,Y^H]^{V_\alpha }=0,(\alpha =1,...,k).
\]

2) The distribution $N_\alpha $ is integrable, if and only if:
\[
\lbrack X^{V_\alpha },Y^{V_\alpha }]^H=0,\ [X^{V_\alpha },Y^{V_\alpha
}]^{V_\beta }=0,\ (\alpha \not =\beta \ ,\alpha ,\beta =1,...,k).
\]

The notion of $h$- or $v_\alpha $- lift of a vector fields $X=\mathcal{X}
(M), $ $X=X^i(x)\displaystyle\frac \partial {\partial x^i}$ is obvious. We
have:
\begin{equation}
l_h(X)=X^i(x)\displaystyle\frac \delta {\delta x^i},\ l_{v_\alpha
}(X)=X^i(x) \displaystyle\frac \delta {\delta y^{\left( \alpha \right)
i}},(\alpha =1,...,k).  \tag{1.3.14}
\end{equation}

\section{The Dual Coefficients of a Nonlinear Connection}

Consider a nonlinear connection $N$, having the coefficients $( \underset{\left( 1\right)
}{N_j^i},...,\underset{\left( k\right) }{N_j^i}
).$ The adapted basis $(\displaystyle\frac \delta {\delta x^i},%
\displaystyle
\frac \delta {\delta y^{\left( 1\right) i}},...,\displaystyle\frac \delta
{\delta y^{\left( k\right) i}})$ to the direct decomposition (1.3.7) is
expressed in the formulae (1.3.9). Its dual basis, denoted by
\begin{equation}
\delta x^i,\delta y^{\left( 1\right) i},...,\delta y^{\left( k\right) i},
\tag{1.4.1}
\end{equation}
can be uniquely written in the form:

\begin{equation}
\left.
\begin{array}{l}
\delta x^i=dx^i, \\
\\
\delta y^{\left( 1\right) i}=dy^{\left( 1\right)
i}+\underset{\left(
1\right) }{M_j^i}dx^j, \\
............................................. \\
\delta y^{\left( k\right) i}=dy^{\left( k\right)
i}+\underset{\left( 1\right) }{M_j^i}dy^{\left( k-1\right)
j}+...+\underset{\left( k-1\right) }{M_j^i}dy^{\left( 1\right)
j}+\underset{\left( k\right) }{M_j^i}dx^j,
\end{array}
\right.  \tag{1.4.2}
\end{equation}
where
\begin{equation}
\left.
\begin{array}{c}
\underset{\left( 1\right) }{M_j^i}=\underset{\left( 1\right)
}{N_j^i},\ \underset{\left( 2\right) }{M_j^i}=\underset{\left(
2\right) }{N_j^i}+ \underset{\left( 1\right)
}{N_j^m}\underset{\left( 1\right) }{M_m^i}
,...., \\
\\
\underset{\left( k\right) }{M_j^i}=\underset{\left( k\right)
}{N_j^i}+ \underset{\left( k-1\right) }{N_j^m}\underset{\left(
1\right) }{M_m^i} +....+\underset{\left( 1\right)
}{N_j^m}\underset{\left( k-1\right) }{ M_m^i}.
\end{array}
\right.  \tag{1.4.3}
\end{equation}

The system of functions ($\underset{\left( 1\right) }{M_j^i},...,
\underset{\left( k\right) }{M_j^i})$ is called the system of dual
coefficients of the nonlinear connection $N.$ They are determined
entirely by means of the coefficients $(\underset{\left(
1\right) }{N_j^i} ,...,\underset{\left( k\right) }{N_j^i},).$

Conversely, if the dual coefficients are given, then the coefficients \newline $(\underset{\left(
1\right) }{N_j^i} ,...,\underset{\left( k\right) }{N_j^i},)$
are expressed by:
\begin{equation}
\left.
\begin{array}{c}
\underset{\left( 1\right) }{N_j^i}=\underset{\left( 1\right)
}{M_j^i},\ \underset{\left( 2\right) }{N_j^i}=\underset{\left(
2\right) }{M_j^i}- \underset{\left( 1\right)
}{M_m^i}\underset{\left( 1\right) }{N_j^m}
,...., \\
\\
\underset{\left( k\right) }{N_j^i}=\underset{\left( k\right)
}{M_j^i}- \underset{\left( 1\right) }{M_m^i}\underset{\left(
k-1\right) }{N_j^m} -....-\underset{\left( k-1\right)
}{M_m^i}\underset{\left( 1\right) }{ N_j^m}.
\end{array}
\right.  \tag{1.4.3'}
\end{equation}

It follows, without difficulties, that with respect to (1.1.2) we have:
\begin{equation}
\delta \widetilde{x}^i=\displaystyle\frac{\partial \widetilde{x}^i}{\partial
x^j}\delta x^j,~\delta \widetilde{y}^{(\alpha )i}=\displaystyle\frac{
\partial \widetilde{x}^i}{\partial x^j}\delta y^{(\alpha )j},~(\alpha
=1,...,k).  \tag{1.4.3''}
\end{equation}

Also, with respect to (1.1.2) the dual coefficients are transformed by the
rule:
\begin{equation}
\left\{
\begin{array}{l}
\underset{(1)}{M_j^m}\displaystyle\frac{\partial
\widetilde{x}^i}{\partial
x^m}=\underset{(1)}{\widetilde{M}_m^i}\displaystyle\frac{\partial
\widetilde{x}^m}{\partial x^j}+\displaystyle\frac{\partial
\widetilde{y}
^{(1)i}}{\partial x^j}, \\
............................................. \\
\underset{(k)}{M_j^m}\displaystyle\frac{\partial
\widetilde{x}^i}{\partial
x^m}=\underset{(k)}{\widetilde{M}_m^i}\displaystyle\frac{\partial
\widetilde{x}^m}{\partial x^j}+\underset{(k-1)}{\widetilde{M}_m^i} %
\displaystyle\frac{\partial \widetilde{y}^{(1)m}}{\partial
x^j}+\cdots +
\underset{(1)}{\widetilde{M}_m^i}\displaystyle\frac{\partial
\widetilde{y} ^{(k-1)m}}{\partial x^j}+\displaystyle\frac{\partial
\widetilde{y}^{(k)i}}{
\partial x^j}.
\end{array}
\right.  \tag{1.4.4}
\end{equation}

The relations between the natural basis on the manifold $T^kM$ and the adapted
basis are immediately:
\begin{equation}
\begin{array}{l}
\displaystyle\frac \partial {\partial x^i}=\displaystyle\frac
\delta {\delta x^i}+\underset{(1)}{M_i^j}\displaystyle\frac \delta
{\delta y^{(1)j}}+\cdots +\underset{(k)}{M_i^j}\displaystyle\frac
\delta {\delta
y^{(k)j}}, \\
\displaystyle\frac \partial {\partial y^{(1)i}}=\displaystyle\frac
\delta {\delta y^{(1)i}}+\underset{(1)}{M_i^j}\displaystyle\frac
\delta {\delta y^{(2)j}}+\cdots
+\underset{(k-1)}{M_i^j}\displaystyle\frac \delta {\delta
y^{(k)j}}, \\
................................................................ \\
\displaystyle\frac \partial {\partial y^{\left( k\right) i}}=\displaystyle %
\frac \delta {\delta y^{(k)i}}.
\end{array}
\tag{1.4.5}
\end{equation}

Similarly, we have:
\begin{equation}
\begin{array}{l}
dx^i=\delta x^i, \\
\\
dy^{(1)i}=\delta y^{(1)i}-\underset{(1)}{N_j^i}\delta x^j, \\
...............................................................................
\\
dy^{(k)i}=\delta y^{(k)i}-\underset{(1)}{N_j^i}\delta
y^{(k-1)j}-\cdots - \underset{(k-1)}{N_j^i}\delta
y^{(1)j}-\underset{(k)}{N_j^i}\delta x^j.
\end{array}
\tag{1.4.5'}
\end{equation}

As an application we can prove:

\begin{teo}
1) The Liouville vector fields $\stackrel{1}{\Gamma }$, ..., $\stackrel{k}{%
\Gamma }$ can be expressed in the adapted basis (1.3.8) in the form
\end{teo}
\begin{equation}
\begin{array}{l}
\stackrel{1}{\Gamma }=z^{(1)i}\displaystyle\frac \delta {\delta y^{(k)i}},
\vspace{3mm}\\
\stackrel{2}{\Gamma }=z^{(1)i}\displaystyle\frac \delta {\delta
y^{(k-1)i}}+2z^{(2)i}\displaystyle\frac \delta {\delta y^{(k)i}}, \\
................................................................. \\
\stackrel{k}{\Gamma }=z^{(1)i}\displaystyle\frac \delta {\delta
y^{(1)i}}+2z^{(2)i}\displaystyle\frac \delta {\delta y^{(2)i}}+\cdots
+kz^{(k)i}\displaystyle\frac \delta {\delta y^{(k)i}},
\end{array}
\tag{1.4.6}
\end{equation}
{\it where}
\begin{equation}
\begin{array}{l}
z^{(1)i}=y^{(1)i},~2z^{(2)i}=2y^{(2)i}+\underset{(1)}{M_m^i}y^{(1)m},~...,
\vspace{3mm}\\
kz^{(k)i}=ky^{(k)i}+(k-1)\underset{(1)}{M_m^i}y^{(k-1)m}+\cdots +%
\underset{(k-1)}{M_m^i}y^{(1)m}
\end{array}
\tag{1.4.7}
\end{equation}

{\it 2) With respect to (1.1.2) we have:}
\begin{equation}
\widetilde{z}^{(\alpha )i}=\displaystyle\frac{\partial \widetilde{x}^i}{%
\partial x^j}z^{(\alpha )j},~(\alpha =1,...,k).  \tag{1.4.7'}
\end{equation}

We note that the formulas (1.4.7') express the geometrical meaning
of each entry $z^{(1)i}$, ...,$z^{(k)i}$. So, we call them \textit{the
Liouville distinguished vector fields} (shortly, $d$\textit{-vector fields}
). These vectors are important in the geometry of higher order Lagrangians.

A field of $1$-form $\omega \in \mathcal{X}^{*}(T^kM)$ can be uniquely
written as
\begin{equation}
\omega =\omega ^H+\omega ^{V_1}+\cdots +\omega ^{V_k}  \tag{1.4.8}
\end{equation}
where
\begin{equation}
\omega ^H=\omega \circ h,~\omega ^{V_\alpha }=\omega \circ v_\alpha
,~(\alpha =1,...,k).  \tag{1.4.8'}
\end{equation}

In the adapted cobasis (1.4.1) we get:
\begin{equation}
\omega ^H=\omega _i^{(o)}\delta x^i,~\omega ^{V_\alpha }=\omega _i^{(\alpha
)}\delta y^{(\alpha )i},~(\alpha =1,...,k).  \tag{1.4.8''}
\end{equation}

For any function $f\in \mathcal{F}(T^kM)$, the $1$-form $df$ has the
components:
\[
df=(df)^H+(df)^{V_1}+\cdots +(df)^{V_k}.
\]

Using (1.4.8'') we obtain:
\begin{equation}
(df)^H=\displaystyle\frac{\delta f}{\delta x^i}\delta x^i,~(df)^{V_\alpha }= %
\displaystyle\frac{\delta f}{\delta y^{(\alpha )i}}\delta y^{(\alpha
)i},~(\alpha =1,...,k).  \tag{1.4.8'''}
\end{equation}

Let $\gamma :I\rightarrow T^kM$ be a parametrized curve, locally expressed by
\begin{equation}
x^i=x^i(t),~y^{(\alpha )i}=y^{(\alpha )i}(t),~t\in I,~(\alpha =1,...,k).
\tag{1.4.9}
\end{equation}

The tangent vector field can be expressed as:
$$
\displaystyle\frac{d\gamma }{dt}=\left( \displaystyle\frac{d\gamma }{dt}
\right) ^H+\cdots +\left( \displaystyle\frac{d\gamma }{dt}\right) ^{V_k}=
$$
$$
=
\displaystyle\frac{dx^i}{dt}\displaystyle\frac \delta {\delta x^i}+ %
\displaystyle\frac{\delta y^{(1)i}}{dt}\displaystyle\frac \delta {\delta
y^{(1)i}}+\cdots +\displaystyle\frac{\delta y^{(k)i}}{dt}\displaystyle\frac
\delta {\delta y^{(k)i}},
$$
where, by means of (1.4.2) one can write:
\begin{equation}
\begin{array}{l}
\displaystyle\frac{\delta y^{(1)i}}{dt}=\displaystyle\frac{dy^{(1)i}}{dt}+
\underset{(1)}{M_j^i}\displaystyle\frac{dx^j}{dt},..., \\
............................................. \\
\displaystyle\frac{\delta
y^{(k)i}}{dt}=\displaystyle\frac{dy^{(k)i}}{dt}+
\underset{(1)}{M_j^i}\displaystyle\frac{dy^{(k-1)j}}{dt}+\cdots +
\underset{(k)}{M_j^i}\displaystyle\frac{dx^j}{dt}.
\end{array}
\tag{1.4.10}
\end{equation}

A parametrized curve $\gamma $ is called \textit{horizontal} if $\left( %
\displaystyle\frac{d\gamma }{dt}\right) ^{V_\alpha }=0$, ($\alpha =1,..,k$).
It is characterized by the system of differential equations:
\begin{equation}
\displaystyle\frac{\delta y^{(1)i}}{dt}=\cdots =\displaystyle\frac{\delta
y^{(k)i}}{dt}=0.  \tag{1.4.11}
\end{equation}

A parametrized curve $c:I\rightarrow M$ on the base manifold $M$, analytically
given by $x^i=x^i(t)$, $t\in I$, has its extension $\widetilde{c} :I\rightarrow T^kM$, given by:
$$
x^i=x^i(t), y^{(1)i}=\displaystyle\frac{dx^i}{dt}, ..., y^{(k)i}=\displaystyle\frac{1}{k!}\displaystyle\frac{d^kx^i}{dt^k}.
$$
A horizontal curve $\widetilde c$ is called an \textit{autoparallel curve} of the
nonlinear connection $N$.

\begin{teo}
The autoparallel curves of the nonlinear connection $N$ are characterized by
the system of differential equations
\begin{equation}
\begin{array}{l}
y^{(1)i}=\displaystyle\frac{dx^i}{dt},\cdots ,y^{(k)i}=\displaystyle\frac
1{k!}\displaystyle\frac{d^kx^i}{dt^k}, \vspace{3mm}\\
\displaystyle\frac{\delta y^{(1)i}}{dt}=0,\cdots ,\displaystyle\frac{\delta
y^{(k)i}}{dt}=0.
\end{array}
\tag{1.4.12}
\end{equation}
\end{teo}

\section{The Determination of a Nonlinear Connection}

A nonlinear connection $N$ on the manifold of accelerations of order $k$, $%
T^kM$ can be determined by a $k$-semispray $S$ or a Riemannian structure $%
g(x)$ defined on the base manifold, or by a Finslerian or Lagrangian
structure over the manifold $M$.

A first result is as follows:

\begin{teo}
(R. Miron and Gh. Atanasiu,[94]) If a $k$-semispray $S$, with the coefficients $%
G^i(x,y^{(1)},...,y^{(k)})$ is given on $T^kM$, then the set of functions:
\begin{equation}
\begin{array}{l}
\underset{(1)}{M_j^i}=\displaystyle\frac{\partial G^i}{\partial y^{(k)j}},~%
\underset{(2)}{M_j^i}=\displaystyle\frac 12\left( S\underset{(1)}{M_j^i}+%
\underset{(1)}{M_m^i}\underset{(1)}{M_j^m}\right) ,..., \\
\\
\underset{(k)}{M_j^i}=\displaystyle\frac 1k\left( S\underset{(k-1)}{M_j^i%
}+\underset{(1)}{M_m^i}\underset{(k-1)}{M_j^m}\right)
\end{array}
\tag{1.5.1}
\end{equation}
is the set of dual coefficients of a nonlinear connection $N$ determined
only by the $k$-semispray $S$.
\end{teo}

The proof can be find in the book [94].

Other result obtained by I. Bucataru ([26,27]) is given in:

\begin{teo}
If $S$ is a $k$-semispray with the coefficients $G^i$, then the following
set of functions
\begin{equation}
\underset{(1)}{M_j^{*i}}=\displaystyle\frac{\partial G^i}{\partial y^{(k)j}%
},\underset{(2)}{M_j^{*i}}=\displaystyle\frac{\partial
G^i}{\partial
y^{(k-1)j}},...,\underset{(k)}{M_j^{*i}}=\displaystyle\frac{\partial G^i}{%
\partial y^{(1)j}}  \tag{1.5.2}
\end{equation}
define the dual coefficients of a nonlinear connection $N^{*}$ determined
only by the $k$-semispray $S$.
\end{teo}

The problem is to determine a nonlinear connection $N$ on $T^kM$ from a
Riemannian structure $g_{ij}(x)$, given on the base manifold $M$.

Let us consider $\gamma _{jk}^i(x)$ the Christoffel symbols of the tensor $%
g_{ij}(x)$. Then, we obtain, ([94]), without difficulties:

\begin{teo}
The following set of functions
\begin{equation}
\begin{array}{l}
\underset{(1)}{M_j^i}(x,y^{(1)})=\gamma _{jm}^i(x)y^{(1)m}, \\
\\
\underset{(2)}{M_j^i}(x,y^{(1)},y^{(2)})=\displaystyle\frac
12\left(
\Gamma \underset{(1)}{M_j^i}+\underset{(1)}{M_m^i}\underset{(1)}{M_j^m}%
\right) , \\
............................................................... \\
\underset{(k)}{M_j^i}(x,y^{(1)},...,y^{(k)})=\displaystyle\frac
1k\left(
\Gamma \underset{(k-1)}{M_j^i}+\underset{(1)}{M_m^i}\underset{(k-1)}{%
M_j^m}\right) ,
\end{array}
\tag{1.5.3}
\end{equation}
where $\Gamma $ is the operator (1.2.3), has the properties:

1$^{\circ }$ It defines the dual coefficients of a nonlinear connection $N$
on $\widetilde{T^kM}$, determined only by the Riemannian structure $g_{ij}(x)
$.

2$^{\circ }$ $\underset{(1)}{M_j^i}(x,y^{(1)})$ depend linearly on $y^{(1)i}$, ..., $\underset{(k)}{M_j^i}(x,y^{(1)},...,y^{(k)})$
depend linearly on $y^{(k)i}$.
\end{teo}

In the same manner we can determine a nonlinear connection on $\widetilde{
T^kM}$ by means of a Finsler space. Namely, we can prove:

\begin{teo}
Let $N_j^i(x,y^{(1)})$ be the Cartan nonlinear connection of a Finsler space
$F^n=(M,F(x,y^{(1)}))$. Then, the following set of functions
\begin{equation}
\begin{array}{l}
\underset{(1)}{M_j^i}(x,y^{(1)})=N_j^i(x,y^{(1)}), \\
\\
\underset{(2)}{M_j^i}(x,y^{(1)},y^{(2)})=\displaystyle\frac
12\left(
\Gamma \underset{(1)}{M_j^i}+\underset{(1)}{M_m^i}\underset{(1)}{M_j^m}%
\right) , \\
............................................................... \\
\underset{(k)}{M_j^i}(x,y^{(1)},...,y^{(k)})=\displaystyle\frac
1k\left(
\Gamma \underset{(k-1)}{M_j^i}+\underset{(1)}{M_m^i}\underset{(k-1)}{%
M_j^m}\right) ,
\end{array}
\tag{1.5.4}
\end{equation}
has the properties:

1$^{\circ }$ It defines the dual coefficients of a nonlinear connection $N$
on $\widetilde{T^kM}$, determined only by the fundamental function $%
F(x,y^{(1)})$ of the Finsler space $F^n$.

2$^{\circ }$ $\underset{(2)}{M_j^i}$ depend linearly on $y^{(2)i}$, ..., $\underset{(k)}{M_j^i}$ depend linearly on $y^{(k)i}.$
$\Gamma $ being the operator (1.2.3).
\end{teo}

\begin{rem}
Theorems 1.5.3 and 1.5.4 prove the existence of the nonlinear connections on
$\widetilde{T^kM}$ in the case when the base manifold is paracompact.
\end{rem}

Let us consider the adapted basis (1.3.8) and adapted cobasis (1.4.1)
corresponding to the nonlinear connection with the dual coefficients (1.5.3).

Let $X^i(x)$ be a vector field on the manifold $M$. We obtain
\begin{equation}
\begin{array}{c}
l_hX=X^i(x)\displaystyle\frac \delta {\delta x^i},\ l_{v_\alpha }X=X^i(x) %
\displaystyle\frac \delta {\delta y^{(\alpha )i}}, \vspace{3mm}\\
(\alpha =1,...,n),\ \forall X=X^i(x)\displaystyle\frac \partial {\partial
x^i}.
\end{array}
\tag{1.5.5}
\end{equation}
Of course $l_hX$, $l_{v_\alpha }X$ are $h$- and respectively $v_\alpha $-lifts of the vector
field $X=X^i(x)\displaystyle\frac \partial {\partial x^i}$.

\begin{teo}
If $g_{ij}(x)$ is a Riemannian structure on the base manifold $M$ and $N$ is
a nonlinear connection with the dual coefficients (1.5.3) determined by $%
g_{ij}(x)$, then
\begin{equation}
\Bbb{G}=g_{ij}(x)dx^i\otimes dx^j+g_{ij}(x)\delta y^{(1)i}\otimes \delta
y^{(1)j}+\cdots +g_{ij}(x)\delta y^{(k)i}\otimes \delta y^{(k)j}  \tag{1.5.6}
\end{equation}
is a Riemannian structure on $\widetilde{T^kM}$ induced only by $g_{ij}(x)$.
\end{teo}

The proof is not difficult. The Riemannian structure (1.5.6) is \textit{the
prolongation} to $T^kM$ of the Riemannian structure $g_{ij}(x)$, [94].

Using the same nonlinear connection $N$ with the dual coefficients (1.5.3) we
can consider the $\mathcal{F}(T^kM)$-linear mapping

\begin{equation}
\begin{array}{c}
\Bbb{F}\left( \displaystyle\frac \delta {\delta x^i}\right) =-\displaystyle %
\frac \delta {\delta y^{(k)i}},\Bbb{F}\left( \displaystyle\frac \delta
{\delta y^{(\alpha )i}}\right) =0, (\alpha =1,..,k), \Bbb{F}\left( \displaystyle\frac \delta
{\delta y^{(k)i}}\right) =\displaystyle\frac \delta
{\delta y^{(k)i}}\\
\\
(i=1,...,n).
\end{array}
\tag{1.5.7}
\end{equation}

\newpage
One proves:

1$^{\circ }$ $\Bbb{F}$ is globally defined on the manifold $\widetilde{T^kM}$.

2$^{\circ }$ $\Bbb{F}$ is a tensor field of type $(1,1)$:

\begin{equation}
\Bbb{F=-}\displaystyle\frac \delta {\delta y^{(k)i}}\otimes dx^i+ %
\displaystyle\frac \delta {\delta x^i}\otimes \delta y^{(k)i}.  \tag{1.5.8}
\end{equation}

3$^{\circ }$ $Ker$ $\Bbb{F}=N_1\oplus \cdots \oplus N_{k-1},$ $Im$ $\Bbb{F}
=N_0\oplus V_k$.

4$^{\circ }$ $rank$ $\Bbb{F}=2n$.

5$^{\circ }$ $\Bbb{F}^3+\Bbb{F}=0.$

6$^{\circ }$ $\Bbb{F}$ is determined only by $g_{ij}(x)$.

Concluding we have:

\begin{teo}
A Riemannian structure $g_{ij}(x)$ given on the base manifold determine an
almost $(k-1)n$-contact structure $\Bbb{F}$ on $\widetilde{T^kM}$ depending
only on $g_{ij}(x)$.
\end{teo}

Let $\left( \underset{(1)}{\xi _a},...,\underset{(k-1)}{\xi
_a}\right) $
, $(a=1,...,n)$ be a local basis adapted to the direct decomposition $%
N_1\oplus \cdots \oplus N_{k-1}$ and $\left( \stackrel{(1)}{\eta ^a},...,
\stackrel{(k-1)}{\eta ^a}\right) $, $(a=1,...,n)$ its dual.

In the book [94] is proved the result:

\begin{teo}
The set $\left( \Bbb{G},\Bbb{F},\underset{(1)}{\xi _a},...,\underset{%
(k-1)}{\xi _a},\stackrel{(1)}{\eta ^a},...,\stackrel{(k-1)}{\eta ^a}\right) $
is a Riemannian almost $(k-1)n$-contact structure on the manifold $%
\widetilde{T^kM}$ determined only by the Riemannian structure $g_{ij}(x)$
defined on the base manifold $M$.
\end{teo}

A similar theory we can study for a Finsler space $F^n=(M,F(x,y^{(1)}))$
using the nonlinear connection (1.5.4) defined only by the fundamental
function $F(x,y^{(1)})$.

Also, we can investigate the problem of determination of a nonlinear
connection on $T^kM$ by means of a Riemannian structure $\Bbb{G}$ given on
the manifold $T^kM$.

One shows that: A Riemann structure $\Bbb{G}$ on $T^kM$ determine a
Riemannian almost $(k-1)n$-contact structure depending only on $\Bbb{G}$ (K. Matsumoto and R. Miron, [111]).

\section{$d$-Tensor Fields. $N$-Linear Connections}

Let $N$ be a nonlinear connection on the manifold $T^kM$. Then, the direct
decomposition (1.3.7) holds. With respect to (1.3.7) a vector field $X$ on $T^kM$
and a $1$-form $\omega $ on $T^kM$ can be uniquely written in the form
(1.3.13) and (1.4.8), respectively.

A distinguished tensor field (shortly $d$-tensor) $T$ on $T^kM$ of type $%
(r,s)$ is a tensor field $T$ with the property
\[
T\left( \stackrel{(1)}{\omega },...,\stackrel{(r)}{\omega
},\underset{(1)}{ X},...,\underset{(s)}{X}\right) =T\left(
\stackrel{(1)}{\omega ^H},..., \stackrel{(r)}{\omega
^{V_k}},\underset{(1)}{X^H},...,\underset{(s)}{ X^{V_k}}\right)
\]
for any $1$-forms $\stackrel{(1)}{\omega
},...,\stackrel{(r)}{\omega }$ from $\mathcal{X}^{*}(T^kM)$ and
any vector fields $\underset{(1)}{X},..., \underset{(s)}{X}$ from
$\mathcal{X}(T^kM)$.

In the adapted cobasis (1.3.8) and its dual basis (1.4.1) a $d$-tensor field $T$
has the components:
\[
T_{j_1...j_s}^{i_{1...}i_r}(x,y^{(1)},...,y^{(k)})=T\left( \delta
x^{i_1},...,\delta y^{(k)i_r},\displaystyle\frac \delta {\delta
x^{j_1}},..., \displaystyle\frac \delta {\delta y^{(k)j_s}}\right) .
\]

Using the formulas (1.3.10) and (1.4.3'')we deduce:

A transformation of coordinates (1.2) on $T^kM$ implies a transformation of
the components of the $d$-tensor field $T$ by the classical rule
\begin{equation}
\widetilde{T}_{j_1...j_s}^{i_{1...}i_r}=\displaystyle\frac{\partial
\widetilde{x}^{i_1}}{\partial x^{h_1}}\cdots \displaystyle\frac{\partial
\widetilde{x}^{i_r}}{\partial x^{hr}}\displaystyle\frac{\partial x^{k_1}}{
\partial \widetilde{x}^{j_1}}\cdots \displaystyle\frac{\partial x^{k_s}}{
\partial \widetilde{x}^{j_s}}T_{k_1...k_s}^{h_{1...}h_r}.  \tag{1.6.1}
\end{equation}
But, this fact is possible only for the components of a $d$-tensor in the
adapted basis. The components of $T$ in the natural basis $\left(
\displaystyle\frac \partial {\partial x^i},\displaystyle\frac \partial
{\partial y^{(1)i}},\cdots ,\displaystyle\frac \partial {\partial
y^{(k)i}}\right) $ and natural cobasis $\left(
dx^i,dy^{(1)i},...,dy^{(k)i}\right) $ have very complicated rule of
transformations with respect to the changing of local coordinates (1.1.2) on
the manifold $T^kM$.

If $X\in \mathcal{X}(T^kM),$ its projections $X^H$, $X^{V_1}$, ..., $X^{V_k}$ are $d$
-vector fields and them components $X^{(0)i}$, $X^{(\alpha )i}$, $(\alpha
=1,...,k)$ are called also $d$-vector fields.

As an example, we have that $z^{(1)i}$, ..., $z^{(k)i}$ from (1.4.6), are $d$-vector fields. They are called the Liouville $d$-vector fields.

For a $1$-form $\omega $, $\omega ^H$, $\omega ^{V_\alpha }$ from (1.4.8) are $%
d$-$1$-forms and them components $\omega ^{(0)i}$, ..., $\omega ^{(\alpha
)i} $ from (1.4.8") are $d$-covector fields.

Of course, the set of $d$-tensor fields with respect to ordinary
operation of addition and tensor product determines a tensor
algebra on the ring of functions $\mathcal{F}(T^kM)$.

An $N$-linear connection on the total space of $k$-accelerations bundle $%
T^kM $ is a linear connection $D$ on $T^kM$ with the properties:

(1) $D$ preserves by parallelism the horizontal distribution $N$;

(2) The $k$-tangent structure $J$ is absolutely parallel with respect to $D$.

As a consequence of this definition, the following theorems hold, [94]:

\begin{teo}
A linear connection $D$ on the manifold $T^kM$ is a $N$-linear connection if
and only if the following properties hold:
\begin{equation}
\begin{array}{l}
\left( D_XY^H\right) ^{V_\alpha }=0,\ \left( D_XY^{V_\alpha }\right) ^H=0,
\\
\left( D_XY^{V_\alpha }\right) ^{V_\beta }=0,\ (\alpha \neq \beta ,\alpha
,\beta =1,...,k),
\end{array}
\tag{1.6.2}
\end{equation}
\begin{equation}
D_X(JY^H)=JD_XY^H,D_X(JY^{V_\alpha })=JD_XY^{V_\alpha },(\alpha =1,...,k),
\tag{1.6.2'}
\end{equation}
for any $X,Y\in \mathcal{X}(T^kM)$.
\end{teo}

If $X\in \mathcal{X}(T^kM)$ is written in the form (1.3.13), for any $Y\in
\mathcal{X}(T^kM)$ and $D$ an $N$-linear connection, we have:

\begin{equation}
D_XY=D_{X^H}Y+D_{X^{V_1}}Y+\cdots +D_{X^{V_k}}Y.  \tag{1.6.3}
\end{equation}

Let us denote:

\[
D_X^H=D_{X^H},D_X^{V_1}=D_{X^{V_1}},....,D_X^{V_k}=D_{X^{V_k}}.
\]

Consequently, we can write:

\begin{equation}
D_XY=D_X^HY+D_X^{V_1}Y+\cdots +D_X^{V_k}Y.  \tag{1.6.3'}
\end{equation}

The operators $D_X^H$, ..., $D_X^{V_\alpha }$ are not covariant derivations
in the $d$-tensor algebra, since $D_X^Hf$ $=X^Hf\neq Xf$.

But they have similar properties with the covariant derivative.

So, if $T$ is a tensor field of type $(r,s)$ we have:

\begin{equation}
\begin{array}{l}
\left( D_X^HT\right) \left( \stackrel{(1)}{\omega
},...,\stackrel{(r)}{ \omega
},\underset{(1)}{X},...,\underset{(s)}{X}\right) =X^HT\left(
\stackrel{(1)}{\omega },...,\stackrel{(r)}{\omega
},\underset{(1)}{X},...,
\underset{(s)}{X}\right) - \\
\ \text{-}T\left( D_X^H\stackrel{(1)}{\omega
},...,\stackrel{(r)}{\omega },
\underset{(1)}{X},...,\underset{(s)}{X}\right) \text{-}...\text{-}
T\left( \stackrel{(1)}{\omega },...,D_X^H\stackrel{(r)}{\omega
},\underset{
(1)}{X},...,\underset{(s)}{X}\right) - \\
\ -\left( \stackrel{(1)}{\omega },...,\stackrel{(r)}{\omega },D_X^H
\underset{(1)}{X},...,\underset{(s)}{X}\right) -...-T\left( \stackrel{%
(1) }{\omega },...,\stackrel{(r)}{\omega
},\underset{(1)}{X},...,D_X^H
\underset{(s)}{X}\right), \qquad \\
\\
(D_X^{V_\alpha }T)\left( \stackrel{(1)}{\omega },...,\stackrel{(r)}{\omega }%
, \underset{(1)}{X},...,\underset{(s)}{X}\right) =X^{V_\alpha
}T\left( \stackrel{(1)}{\omega },...,\stackrel{(r)}{\omega
},\underset{(1)}{X},...,
\underset{(s)}{X}\right) - \\
\ \text{-}T\left( D_X^{V_\alpha }\stackrel{(1)}{\omega
},...,\stackrel{(r)}{ \omega
},\underset{(1)}{X},...,\underset{(s)}{X}\right) \text{-}...\text{
-}T\left( \stackrel{(1)}{\omega },...,D_X^{V_\alpha }\stackrel{(r)}{\omega }%
, \underset{(1)}{X},...,\underset{(s)}{X}\right) - \\
\ -\left( \stackrel{(1)}{\omega },...,\stackrel{(r)}{\omega },D_X^{V_\alpha
} \underset{(1)}{X},...,\underset{(s)}{X}\right) -...-T\left( \stackrel{%
(1) }{\omega },...,\stackrel{(r)}{\omega },\underset{(1)}{X}%
,...,D_X^{V_\alpha }\underset{(s)}{X}\right). \qquad
\end{array}
\tag{1.6.4}
\end{equation}

For instance, the formula (1.6.4) for a 1-form $\omega \in \mathcal{X}^{*}(T^kM)$ leads to the following expressions of $D_X^{H}\omega$ and $D_X^{V_\alpha }\omega $:
\begin{equation}
\begin{array}{l}
(D_X^H\omega )(Y)=X^H\omega (Y)-\omega (D_X^HY), \vspace{3mm}\\
(D_X^{V_\alpha }\omega )(Y)=X^{V_\alpha }\omega (Y)-\omega (D_X^{V_\alpha
}Y),~(\alpha =1,...,k).
\end{array}
\tag{1.6.5}
\end{equation}

In the adapted basis (1.3.13) an $N$-linear connection $D$ can be uniquely
represented in the form, [94]:
\begin{equation}
\left\{
\begin{array}{l}
D_{\displaystyle\frac \delta {\delta x^j}}\displaystyle\frac \delta {\delta
x^i}=L_{ij}^m\displaystyle\frac \delta {\delta x^m},\ D_{\displaystyle\frac
\delta {\delta x^j}}\displaystyle\frac \delta {\delta y^{(\alpha
)i}}=L_{ij}^m\displaystyle\frac \delta {\delta y^{(\alpha )m}},(\alpha =
\overline{1,k}), \\
\\
D_{\displaystyle\frac \delta {\delta y^{(\beta
)j}}}\displaystyle\frac \delta {\delta x^i}=\underset{(\beta
)}{C_{ij}^m}\displaystyle\frac \delta
{\delta x^m},\ D_{\displaystyle\frac \delta {\delta y^{(\beta )j}}} %
\displaystyle\frac \delta {\delta y^{(\alpha )i}}=\underset{(\beta
)}{ C_{ij}^m}\displaystyle\frac \delta {\delta y^{(\alpha
)m}},(\alpha ,\beta = \overline{1,k}).
\end{array}
\right.  \tag{1.6.6}
\end{equation}

The system of functions $(L_{ij}^m,\underset{(1)}{C_{ij}^m},...,
\underset{(k)}{C_{ij}^m})$ represents \textit{the coefficients} of
$D$.

With respect to (1.1.2), $L_{ij}^m$ are transformed by the same rule
as the coefficients of a linear connection defined on the base
manifold $M$. Others coefficients $\underset{(\alpha )}{C_{ij}^m},$ $\alpha =1,..,k$ are transformed like $d$ -tensors of type $(1,2)$.

If $T$ is a $d$-tensor field of type $(r,s)$, represented in adapted basis
in the form
\begin{equation}
T=T_{j_1...j_s}^{i_1...i_r}\displaystyle\frac \delta {\delta x^{i_1}}\otimes
\cdots \otimes \displaystyle\frac \delta {\delta y^{(k)i_r}}\otimes
dx^{j_1}\otimes \cdots \otimes \delta y^{(k)j_s}  \tag{1.6.7}
\end{equation}
and $X=X^H=X^i(x,y^{(1)},...,y^{(k)})\displaystyle\frac \delta {\delta x^i}$, then the $h$-covariant derivative $D_X^HT$, by means of (1.6.6), is as follows
\begin{equation}
D_X^HT=X^mT_{j_1...j_s|m}^{i_1...i_r}\displaystyle\frac \delta {\delta
x^{i_1}}\otimes \cdots \otimes \displaystyle\frac \delta {\delta
y^{(k)i_r}}\otimes dx^{j_1}\otimes \cdots \otimes \delta y^{(k)j_s},
\tag{1.6.8}
\end{equation}
where
\begin{equation}
T_{j_1...j_s|m}^{i_1...i_r}=\displaystyle\frac{\delta
T_{j_1...j_s}^{i_1...i_r}}{\delta x^m}
+L_{hm}^{i_1}T_{j_1...j_s}^{hi_2...i_r}+\cdots
-L_{j_sm}^hT_{j_1...h}^{i_1...i_r}.  \tag{1.6.9}
\end{equation}

The operator "$_{\mid}$" will be called \textit{the }$h$\textit{-covariant derivative}, with the same denomination as $D_X^H$ .

Consider the $v_\alpha $-covariant derivatives $D_X^{V_\alpha }$, for $X= \stackrel{(\alpha )}{X^i}\displaystyle\frac \delta {\delta y^{(\alpha )i}}$, \newline $(\alpha =1,...,k)$. Then, using (1.6.6) and (1.6.7) we have:
\begin{equation}
D_X^{V_\alpha }T=\stackrel{(\alpha )}{X^m}T_{j_1...j_s}^{i_1...i_r}\stackrel{
(\alpha )}{\mid }_m\displaystyle\frac \delta {\delta x^{i_1}}\otimes \cdots
\otimes \displaystyle\frac \delta {\delta y^{(k)i_r}}\otimes dx^{j_1}\otimes
\cdots \otimes \delta y^{(k)j_s},  \tag{1.6.8'}
\end{equation}
where
\begin{equation}
T_{j_1...j_s}^{i_1...i_r}\stackrel{(\alpha )}{\mid
}_m=\displaystyle\frac{ \delta T_{j_1...j_s}^{i_1...i_r}}{\delta
y^{(\alpha )m}}+\underset{(\alpha
)}{C_{hm}^{i_1}}T_{j_1...j_s}^{h...i_r}+\cdots -\underset{(\alpha
)}{ C_{j_sm}^h}T_{j_1...h}^{i_1...i_r}, \ (\alpha=1,..,k). \tag{1.6.9'}
\end{equation}

The operators "$\stackrel{(\alpha )}{\mid }$", in number of $k$ are called
$v_\alpha $\textit{-covariant derivatives}, with the same denominations as $%
D_X^{V_\alpha }$. Each of them has similar properties to those of $h$
-covariant derivative. For instance
\newpage
$
\begin{array}{l}
\left( fT_{....}^{....}\right) \stackrel{(\alpha )}{\mid }_m=\displaystyle
\frac{\delta f}{\delta y^{(\alpha )m}}T_{...}^{...}+fT_{...}^{...}\stackrel{
(\alpha )}{\mid }_m, \\

\\
\left( \stackrel{1}{T_{...}^{...}}\otimes \stackrel{2}{T_{...}^{...}}\right)
\stackrel{(\alpha )}{\mid }_m=\stackrel{1}{T_{...}^{...}}\stackrel{(\alpha )
}{\mid }_m\otimes \stackrel{2}{T_{...}^{...}}\stackrel{1}{+T_{...}^{...}}
\otimes \stackrel{2}{T_{...}^{...}}\stackrel{(\alpha )}{\mid }_m.
\end{array}
$

The operators $_{\mid}$ and $\stackrel{(\alpha )}{\mid }$ commute with the
operation of contraction, etc.

These operators applied to the Liouville $d$-vector fields $z^{(1)i}$, ..., $%
z^{(k)i}$ determine 'the deflection tensors' of $D$:

\begin{equation}
\stackrel{(\alpha )}{D_j^i}=z_{\quad |j}^{(\alpha )i};\quad \stackrel{(\beta
\alpha )}{d_j^i}=z^{(\beta )i}\stackrel{(\alpha )}{\mid }_j.  \tag{1.6.10}
\end{equation}

\section{Torsion and Curvature}

\textit{The tensor of torsion} $\Bbb{T}$ of an $N$-linear connection $D$:
\begin{equation}
\Bbb{T}(X,Y)=D_XY-D_YX-[X,Y],\quad \forall X,Y\in \mathcal{X}(T^kM),
\tag{1.7.1}
\end{equation}
for $X=X^H+X^{V_1}+\cdots +X^{V_k}$, $Y=Y^H+Y^{V_1}+\cdots +Y^{V_k}$,
uniquely determines the following vector fields:

\begin{equation}
\Bbb{T}(X^H,Y^H),\Bbb{T}(X^H,Y^{V_\alpha }),\Bbb{T}(X^{V_\alpha },Y^{V_\beta
}),(\alpha ,\beta =1,...,k).  \tag{1.7.2}
\end{equation}

Therefore:

The tensor of torsion $\Bbb{T}$ of the $N$-linear connection $D$ is well
determined by the following components, which are $d$-tensor fields of type $%
(1,2)$:
\begin{equation}
\begin{array}{l}
\Bbb{T}(X^H,Y^H)=h\Bbb{T}(X^H,Y^H)+\sum\limits_{\alpha
=1}^kv_\alpha \Bbb{T}
(X^H,Y^H), \\
\\
\Bbb{T}(X^H,Y^{V_\beta })=h\Bbb{T}(X^H,Y^{V_\beta
})+\sum\limits_{\alpha
=1}^kv_\alpha \Bbb{T}(X^H,Y^{V_\beta }), \\
\\
\Bbb{T}(X^{V_\alpha },Y^{V_\beta })=\sum\limits_{\gamma =1}^kv_\gamma \Bbb{%
T }(X^{V_\alpha },Y^{V_\beta }).
\end{array}
\tag{1.7.3}
\end{equation}
It is not difficult to prove that $h\Bbb{T}(X^{V_\alpha },Y^{V_\beta })=0$, $%
\forall \alpha ,\beta =1$, ..., $k$.
The expressions of $d$-tensors of torsion are the following:
\newpage
\begin{equation}
\begin{array}{l}
\left\{
\begin{array}{l}
h\Bbb{T}(X^H,Y^H)=D_X^HY^H-D_Y^HX^H-[X^H,Y^H]^H, \\
\\
v_\alpha \Bbb{T}(X^H,Y^H)=[X^H,Y^H]^{V_\alpha },
\end{array}
\right. \vspace{3mm}\\
\left\{
\begin{array}{l}
h\Bbb{T}(X^H,Y^{V_\beta })=-\nabla _{Y^{V_\beta }}X^H-[X^H,Y^{V_\beta }]^H,
\\
\\
v_\alpha \Bbb{T}(X^H,Y^{V_\beta })=\left( \nabla _{X^{}}^HY^{V_\beta
}\right) ^{V_\alpha }-[X^H,Y^{V_\beta }]^{V_\alpha }, \\
\\
v_\gamma \Bbb{T}(X^{V\alpha },Y^{V_\beta })=v_\gamma \left( \nabla
_{X^{}}^{V_\alpha }Y^{V_\beta }-\nabla _{Y^{V_\beta }}X^{V_\gamma
}-[X^{V_\alpha },Y^{V_\beta }]\right) .
\end{array}
\right.
\end{array}
\tag{1.7.4}
\end{equation}

In the adapted basis we set:
\begin{equation}
\begin{array}{l}
h\Bbb{T}\left( \displaystyle\frac \delta {\delta x^h},\displaystyle\frac
\delta {\delta x^j}\right) =
\underset{}{\underset{(00)}{T_{jh}^i} \displaystyle\frac \delta {\delta x^i}},
v_\alpha \Bbb{T}\left(\displaystyle \frac \delta {\delta x^h},\displaystyle\frac \delta {\delta
x^j}\right) =
\underset{}{\underset{(0\alpha )}{T_{jh}^i}\displaystyle\frac \delta {\delta y^{(\alpha )i}}}, \\
\\
h\Bbb{T}\left( \displaystyle\frac \delta {\delta x^h},\displaystyle\frac \delta {\delta y^{(\beta )j}}\right)
=\underset{}{\underset{(\beta 0)}{T_{jh}^i}\displaystyle\frac \delta {\delta x^i}},
v_\alpha \Bbb{T}\left( \displaystyle\frac \delta {\delta x^h},\displaystyle\frac \delta {\delta y^{(\beta )j}}\right) =
\underset{}{\underset{(\beta \alpha )}{T_{jh}^i}\displaystyle\frac \delta {\delta y^{(\alpha )i}}}, \\
\\
h\Bbb{T}\left( \displaystyle\frac \delta {\delta y^{(\alpha )h}}, \displaystyle\frac \delta {\delta y^{(\beta )j}}\right)
=\underset{}{0},
v_\gamma \Bbb{T}\left( \displaystyle\frac \delta {\delta y^{(\alpha )h}}, \displaystyle\frac \delta {\delta y^{(\beta )j}}\right) =\underset{}{\underset{(\beta \alpha )}{\stackrel{\gamma }{T_{jh}^i}} \displaystyle\frac \delta {\delta y^{(\gamma )i}}}.
\end{array}
\tag{1.7.5}
\end{equation}

Using (1.7.4) and (1.6.6) as well as the following expressions of the Lie brackets,
we can calculate the components of $d$-tensors of torsion.

The Lie brackets of the vector fields of the adapted basis are given by:
\begin{equation}
\begin{array}{l}
\left[ \displaystyle\frac \delta {\delta x^j},\displaystyle\frac
\delta {\delta x^h}\right] =
\underset{(01)}{R_{jh}^i}
\displaystyle\frac \delta {\delta y^{(1)i}}+\cdots +
\underset{(0k)}{R_{jh}^i}
\displaystyle\frac \delta {\delta y^{(k)i}}, \\
\\
\left[ \displaystyle\frac \delta {\delta x^j},\displaystyle\frac \delta
{\delta y^{(\alpha )h}}\right] =
\underset{(\alpha 1)}{B_{jh}^i}
\displaystyle\frac \delta {\delta y^{(1)i}}+\cdots +
\underset{(\alpha k)}{B_{jh}^i}
\displaystyle\frac \delta {\delta y^{(k)i}}, \\
\\
\left[ \displaystyle\frac \delta {\delta y^{(\alpha)j}},\displaystyle\frac \delta {\delta y^{(\beta )h}}\right]
=\underset{(\alpha \beta )}{
\stackrel{(1)}{C_{jh}^i} }\displaystyle\frac \delta {\delta y^{(1)i}}+\cdots +
\underset{(\alpha \beta )}{\stackrel{(k)}{C_{jh}^i}}\displaystyle \frac \delta {\delta y^{(k)i}},
\end{array}
\tag{1.7.6}
\end{equation}
where the coefficients $R$, $B$, $C$ can be calculated by means of the coefficients of the nonlinear connection $N$, [94].

We remark the following $d$-tensor of torsion
\begin{equation}
\underset{}{\underset{(00)}{T_{jk}^i}=L_{jk}^i-L_{kj}^i},\quad
\underset{(\alpha \alpha)}{\stackrel{\alpha }{T_{jk}^i}}=
\underset{(\alpha)}{C_{jk}^i}- \underset{(\alpha )}{C_{kj}^i}\quad (\alpha =1,...,k).
\tag{1.7.7}
\end{equation}

For simplicity we denote
\begin{equation}
\underset{}{\underset{(00)}{T_{jk}^i}=\underset{(0)}{T_{jk}^i}},\quad
\underset{(\alpha \alpha )}{\stackrel{\alpha }{T_{jk}^i}}=
\underset{(\alpha )}{S_{jk}^i}.
\tag{1.7.7'}
\end{equation}

The notion of curvature can be investigated by the same method.

\textit{The curvature tensor} $\Bbb{R}$ of the $N$-linear connection $D$ is
given by
\begin{equation}
\Bbb{R}(X,Y)Z=\left( D_XD_Y-D_YD_X\right) Z-D_{[X,Y]}Z,\quad \forall
X,Y,Z\in \mathcal{X}(T^kM).  \tag{1.7.8}
\end{equation}

Using the formula (1.6.2) we obtain
\begin{equation}
J\left( \Bbb{R}(X,Y)Z\right) =\Bbb{R}(X,Y)JZ,\ldots ,J^k\left( \Bbb{R}
(X,Y)Z\right) =\Bbb{R}(X,Y)J^kZ.  \tag{1.7.9}
\end{equation}

Setting $Z^{V_\alpha }=J^\alpha Z^H$ we have
\[
\Bbb{R}(X,Y)Z^{V_\alpha }=J^\alpha \left( \Bbb{R}(X,Y)Z^H\right) ,\quad
(\alpha =1,...,k).
\]

The essential components of the curvature tensor field $\Bbb{R}$ are $\Bbb{R}(X,Y)Z^H$. It has the properties:
$$
v_\beta \left( \Bbb{R}(X,Y)Z^H\right) =0,\quad h\left( \Bbb{R}
(X,Y)Z^{V_\beta }\right) =0,
$$
$$
v_\beta \left( \Bbb{R}(X,Y)Z^{V_\alpha }\right) =0,\quad (\alpha \neq \beta
,\alpha ,\beta =1,...,k).
$$

Thus, the curvature tensor $\Bbb{R}$ of the $N$-linear connection $D$ gives
rise to the $d$-vector fields:
\[
\Bbb{R}(X^H,Y^H)Z^H=[D_X^H,D_Y^H]Z^H-D_{[X^H,Y^H]}^HZ^H-\sum\limits_{\gamma
=1}^kD_{[X^H,Y^H]}^{V_\gamma }Z^H,
\]
\[
\Bbb{R}(X^{V_\alpha },Y^H)Z^H=[D_X^{V_\alpha
},D_Y^H]Z^H-D_{[X^{V_\alpha },Y^H]}^HZ^H-\sum\limits_{\gamma
=1}^kD_{[X^{V_\alpha },Y^H]}^{V_\gamma }Z^H,
\]
\begin{equation}
\begin{array}{c}
\Bbb{R}(X^{V_\beta },Y^{V\alpha })Z^H=[D_X^{V_\beta
},D_Y^{V_\alpha }]Z^H-\sum\limits_{\gamma =1}^kD_{[X^{V_\beta
},Y^{V_\alpha }]}^{V_\gamma
}Z^H, \\
(\alpha ,\beta =1,...,k;\ \beta \leq \alpha ).
\end{array}
\tag{1.7.10}
\end{equation}

The $d$-tensor fields (1.7.10) are obtained applying the operators $J$, $J^2$,
..., $J^k$ and taking into account $J^\gamma Z^H=Z^{V_\gamma }$, $(\gamma
=1,...,k)$.

In the applications it is suitable to consider the equalities (1.7.10) as Ricci
identities, [94].

The local expressions of $d$-tensors of curvature in adapted basis are:
\begin{equation}
\begin{array}{l}
\Bbb{R}\left( \displaystyle\frac \delta {\delta x^m},\displaystyle\frac
\delta {\delta x^j}\right) \displaystyle\frac \delta {\delta x^h}=R_{h\ jm}^i\displaystyle\frac \delta {\delta x^i}, \\
\Bbb{R}\left( \displaystyle\frac \delta {\delta y^{(\alpha )m}},
\displaystyle \frac \delta {\delta x^j}\right) \displaystyle\frac \delta {\delta x^h}=
\underset{(\alpha )}{P}{}_{h\ jm}^i \displaystyle\frac \delta {\delta x^i}, \\
\Bbb{R}\left( \displaystyle\frac \delta {\delta y^{(\beta )m}},\displaystyle \frac \delta {\delta y^{(\alpha )j}}\right)
\displaystyle\frac \delta {\delta x^h}=
\underset{(\beta \alpha )}{S}{}_{h\ jm}^i\displaystyle \frac \delta {\delta x^i}.
\end{array}
\tag{1.7.11}
\end{equation}

Note that the other $d$-tensors of curvature (1.7.10) have the same
components. The $d$-tensors (1.7.11) are called also $d$\textit{-tensors of
curvature}. They have the expressions:
\begin{equation}
R_{h\ jm}^i=\displaystyle\frac{\delta L_{hj}^i}{\delta
x^m}-\displaystyle \frac{\delta L_{hm}^i}{\delta
x^j}+L_{hj}^pL_{pm}^i-L_{hm}^pL_{pj}^i+\sum \limits_{\gamma
=1}^k\underset{(\gamma )}{C_{hp}^i}\underset{(0\gamma )}{
R_{jm}^p},  \tag{1.7.12}
\end{equation}
\[
\underset{(\alpha )}{P}{}_{h\ jm}^i=\displaystyle\frac{\delta
L_{hj}^i}{ \delta y^{(\alpha )m}}-\displaystyle\frac{\delta
\underset{(\alpha )}{C_{hm}^i}}{\delta
x^j}+L_{hj}^p\underset{(\alpha )}{C_{pm}^i}- \underset{(\alpha
)}{C_{hm}^p}L_{pj}^i+\sum\limits_{\gamma =1}^k \underset{(\gamma
)}{C_{hp}^i}\underset{(\alpha \gamma )}{B_{jm}^p}\ ,
\]
\[
\underset{(\beta \alpha )}{S}{}_{h\ jm}^i=\displaystyle\frac{\delta
\underset{(\alpha )}{C_{hj}^i}}{\delta y^{(\beta )m}}-%
\displaystyle
\frac{\delta \underset{(\beta )}{C_{hm}^i}}{\delta y^{(\alpha )j}}%
+ \underset{(\alpha )}{C_{hj}^p}\underset{(\beta )}{C_{pm}^i}-%
\underset{ (\beta )}{C_{hm}^p}\underset{(\alpha )}{C_{pj}^i}+\sum%
\limits_{\gamma =1}^k\underset{(\gamma
)}{C_{hp}^i}\underset{(\alpha \beta )}{\stackrel{ \gamma
}{C_{jm}^p}}\
\]

The connection $1$-forms of the $N$-linear connection $D$ given in ch. 7, \S 7 of the book [94] are:
\begin{equation}
\omega _{\ j}^i=L_{jh}^idx^h+\underset{(1)}{C_{jh}^i}\delta
y^{(1)h}+\cdots +\underset{(k)}{C_{jh}^i}\delta y^{(k)h}\ .
\tag{1.7.13}
\end{equation}

Therefore, the structure equations of $D$ are given by the following theorem, [94]:

\begin{teo}
The structure equations of an $N$-linear connection $D$ on the manifold $T^kM
$ are given by:
\begin{equation}
\begin{array}{l}
d(dx^i)-dx^m\wedge \omega _{\ m}^i=-\stackrel{(0)}{\Omega ^i}\ , \\
\\
d(\delta y^{(\alpha )i})-\delta y^{(\alpha )m}\wedge \omega _{\ m}^i=-%
\stackrel{(\alpha )}{\Omega ^i}\ , \\
\\
d\omega _{\ j}^i-\omega _{\ j}^m\wedge \omega _{\ m}^i=-\Omega _{\ j}^i\ ,
\end{array}
\tag{1.7.14}
\end{equation}
where the $2$-forms of torsion $\stackrel{(0)}{\Omega ^i}$, $\stackrel{%
(\alpha )}{\Omega ^i}$ are
\newpage
\begin{equation}
\begin{array}{l}
\stackrel{(0)}{\Omega ^i}\ =\displaystyle\frac 12\underset{(0)}{T_{jh}^i}dx^j\wedge dx^h+
\underset{(1)}{C_{jh}^i}dx^j\wedge \delta y^{(1)h}+\cdots  \vspace{3mm}\\
\qquad +\underset{(k)}{C_{jh}^i}dx^j\wedge \delta x^{(k)h}\ , \\
\\
\stackrel{(\alpha )}{\Omega ^i}\ =\displaystyle\frac 12\underset{(0\alpha )}{R_{jh}^i}dx^j\wedge dx^h+\sum\limits_{\gamma =1}^k\underset{(\gamma \alpha )}{B_{jm}^i}dx^j\wedge \delta y^{(\gamma )m}+ \vspace{3mm}\\
\qquad +\sum\limits_{\gamma =1}^k\underset{(\alpha \gamma )}{\stackrel{\beta }{C_{jm}^i}}\delta y^{(\beta )j}\wedge
\delta y^{(\gamma )m}-\vspace{3mm}
\\
\qquad -\left[ L_{jm}^idx^j+\sum\limits_{\gamma =1}^k\underset{(\gamma )}{C_{jm}^i}\delta y^{(\gamma )j}\wedge \delta y^{(\alpha )m}\right] .
\end{array}
\tag{1.7.15}
\end{equation}
and where the $2$-forms of curvature $\Omega _{\ j}^i$ are:
\begin{equation}
\begin{array}{l}
\Omega _{\ j}^i\ =\displaystyle\frac 12R_{j\ pq}^idx^p\wedge dx^q+ \\
\\
\ \qquad +\sum\limits_{\alpha =1}^k\underset{(\alpha )}{P}{}_{j\ pq}^i%
dx^p\wedge \delta y^{(\alpha )q}+\sum\limits_{\alpha ,\beta =1}^k%
\underset{(\alpha \beta )}{S}{}_{j\ pq}^i\delta y^{(\alpha )p}\wedge
\delta y^{(\beta )q}\ .
\end{array}
\tag{1.7.16}
\end{equation}
\end{teo}

The Bianchi identities of $D$ can be derived from (1.7.14) applying the
operator of exterior differentiation and calculating $d\stackrel{(0)}{\Omega
^i}$, $d\stackrel{(\alpha )}{\Omega ^i}$, $d\Omega _{\ j}^i$ from (1.7.15) and
(1.7.16) modulo the system (1.7.14).

\chapter{Lagrange Spaces of Higher Order}

\markboth{\it{THE GEOMETRY OF HIGHER-ORDER HAMILTON SPACES\ \ \ \ \ }}{\it{Lagrange Spaces of Higher Order}}

We define the notion of higher order Lagrangian and the notion of integral
of action. We investigate the variational problem for autonomous
Lagrangians deriving the Euler-Lagrange equations. The notion of Lagrange
space of order $k$ is introduced by means of regular nondegenerate
Lagrangian defined on the total space of the $k$-accelerations bundle $T^kM.$
 In this case the Craig-Synge equations determine a $k$-semispray, which
depend only on the considered Lagrangian. Therefore the geometry of the

Lagrange space of order $k$ is based on the mentioned $k$-semispray, [94].

\section{Lagrangians of Order $k$}

Let us consider the $k$-accelerations bundle $(T^kM,\pi ^k,M)$. In
Analytical Mechanics $M$ is the space of configurations of a physical
system. A point $x=(x^i)\in U\subset M$ is called a \textit{configuration},
a mapping $c:t\in I\rightarrow (x^i(t))\in U$ is a \textit{law of moving
(evolution)}, $t$ is called \textit{time}, a couple $(t,x)$ is an \textit{\
event} and $\left( \displaystyle\frac{dx^i}{dt},\displaystyle\frac 1{2!} %
\displaystyle\frac{d^2x^i}{dt^2},\cdots ,\displaystyle\frac 1{k!} %
\displaystyle\frac{d^kx^i}{dt^k}\right) $ are the velocity and generalized
accelerations of order $1$, $2$,..., $k-1$ (respectively). The factors $%
\displaystyle\frac 1{h!}$ ($h=1,...,k$) are introduced for convenience.

Throughout this book we omit the word {\it generalized} and say shortly,
accelerations of order $h$ for $\displaystyle\frac 1{h!}\displaystyle\frac{
d^hx^i}{dt^h}$, $h=1,...,k$. A law a moving $c:t\in I\rightarrow (x^i(t))\in
U$ $\subset M$ will be called a \textit{curve} parametrized by time $t$. As
usual the curve $\widetilde{c}:t\in I\rightarrow \widetilde{c}(t)\in (\pi
^k)^{-1}(U)\subset T^kM$,
\begin{equation}
\widetilde{c}:t\in I\rightarrow \left(x^i(t),  \displaystyle\frac{dx^i}{dt}, %
\displaystyle\frac 1{2!}\displaystyle\frac{d^2x^i}{dt^2},\cdots , %
\displaystyle\frac 1{k!}\displaystyle\frac{d^kx^i}{dt^k}\right)  \tag{2.1.1}
\end{equation}
is the extension of the curve $c$ to the total space $T^kM$ of the $k$
-acceleration bundle.

\begin{defi}
A Lagrangian of order $k$, ($k\in \mathbf{N}^{*}$) is a mapping $%
L:T^kM\rightarrow \mathbf{R}$.
\end{defi}

This means that $L$ is a real function $L(x,y^{(1)},...,y^{(k)})$ on $T^kM $.
In the other words, with respect to a change of local coordinates on $T^kM$
(1.1.2), we have
\begin{equation}
\widetilde{L}(\widetilde{x},\widetilde{y}^{(1)},...,\widetilde{y}
^{(k)})=L(x,y^{(1)},...,y^{(k)}).  \tag{2.1.2}
\end{equation}

The previous definition is given for \textit{autonomous Lagrangians}. A
similar definition we have for \textit{nonautonomous Lagrangians}. They are
mappings $L:(t,x,y^{(1)},...,y^{(k)})\in \mathbf{R\times }
T^kM\rightarrow L(t,x,y^{(1)},...,y^{(k)})\in \mathbf{R}$, which are real
functions with respect to a change of coordinates on $\mathbf{R\times }T^kM$
, $(t,x,y^{(1)},...,y^{(k)})\rightarrow (\widetilde{t},\widetilde{x},
\widetilde{y}^{(1)},...,\widetilde{y}^{(k)})$, where $t= \widetilde{t}.$

For us it is preferable to study the autonomous Lagrangians, because the notion
of Lagrange space of order $k$ is a geometrical one. But, one sees that the
nonautonomous Lagrangians can be geometrized by means of the notion of

rheonomic Lagrange space of order $k$. Such kind of geometry can by
constructed by the same methods as in the autonomous case.

In the following we investigate the Lagrangians $L(x,y^{(1)},...,y^{(k)})$
which have the property (2.1.2).

A Lagrangian of order $k$, $L:T^kM\rightarrow \mathbf{R}$ is called \textit{%
\ differentiable} if it is of $C^\infty $-class on $\widetilde{T^kM}$ and
continuous on the null section of the projection $\pi ^k:T^kM\rightarrow M$.

The Hessian of a differentiable Lagrangian $L$, with respect to the
variables $y^{(k)i}$ on $\widetilde{T^kM}$ is the matrix $||g_{ij}||$, where
\begin{equation}
g_{ij}=\displaystyle\frac 12\displaystyle\frac{\partial ^2L}{\partial
y^{(k)i}\partial y^{(k)j}}.  \tag{2.1.3}
\end{equation}

One can prove, [94] that $g_{ij}$ is a $d$-tensor field, on the manifold $%
\widetilde{T^kM}.$ This is  covariant of order $2$ and symmetric.

This property determines the geometrical covariance of the Hessian of $L$.

If
\begin{equation}
rank\ ||g_{ij}||=n,\ \ on\ \ \widetilde{T^kM}  \tag{2.1.4}
\end{equation}
we say that $L(x,y^{(1)},...,y^{(k)})$ is a \textit{regular} (or \textit{\
nondegenerate}) Lagrangian.

The existence of the regular Lagrangians of order $k$ is assured by the
following example.

Let $\gamma _{ij}(x)$ a Riemannian tensor field on the base manifold $M$ and
$z^{(k)i}$ the Liouville vector field on $\widetilde{T^kM}$ determined by
the prolongation of order $k$ of the Riemannian spaces $\mathcal{R}
^n=(M,\gamma _{ij}(x)),$  an arbitrary covector field $b_i$ on $T^{k-1}M$
and a function $b$ on $T^{k-1}M$. Then the Lagrangian of order $k$
\begin{equation}
\begin{array}{lll}
L(x,y^{(1)},...,y^{(k)}) & = & \gamma
_{ij}(x)z^{(k)i}z^{(k)j}+b_i(x,y^{(1)},...,y^{(k-1)})z^{(k)i}+ \\
& + & b(x,y^{(1)},...,y^{(k-1)})
\end{array}
\tag{2.1.5}
\end{equation}
is a regular Lagrangian on $T^kM$. We have $g_{ij}(x,y^{(1)},...,y^{(k)})=
\gamma _{ij}(x)$.

Let us consider the scalar fields
\begin{equation}
I^1(L)=\mathcal{L}_{\stackrel{1}{\Gamma }}L,...,I^k(L)=\mathcal{L}_{%
\stackrel{k}{\Gamma }}L,  \tag{2.1.6}
\end{equation}
where $\stackrel{1}{\Gamma }$, ..., $\stackrel{k}{\Gamma }$ are the
Liouville vector fields on $T^kM$ and $\mathcal{L}$ is the Lie operator of
derivation. $I^1(L)$, ..., $I^k(L)$ will be called the \textit{main invariants }
of the Lagrangian $L$.

For a smooth parametrized curve $c:[0,1]\rightarrow M$
represented in a domain of a local chart by $x^i=x^i(t)$, $t\in [0,1]$. The
parameter $t$ is called {\it time} and its extension to $\widetilde{T^kM}$ is $%
\widetilde{c}$, given by (2.1.1).

The integral of action for the differentiable Lagrangian $%
L(x,y^{(1)},...,y^{(k)})$ along curve $c$ is defined by
\begin{equation}
I(c)=\int\limits_0^1L\left(x(t),
\displaystyle\frac{dx}{dt},\displaystyle\frac
1{2!}\displaystyle\frac{d^2x}{dt^2},\cdots ,\displaystyle\frac 1{k!} %
\displaystyle\frac{d^kx}{dt^k}\right) dt.  \tag{2.1.7}
\end{equation}

Ones proves, [94] the following important results:

\begin{teo}
The necessary conditions for the integral of action be independent on the
parametrization of the curve $c$ are
\begin{equation}
I^1(L)=\cdots =I^{k-1}(L)=0,I^k(L)=L.  \tag{2.1.8}
\end{equation}
\end{teo}

The conditions (2.1.8) are called the \textit{Zermello conditions}, [69, 94].

\begin{teo}

If the differentiable Lagrangian $L$ of order $k$, $k>1$, satisfies the
Zermello conditions (2.1.8), then it is degenerate (singular), i.e.
\begin{equation}
rank\ ||g_{ij}(x,y^{(1)},...,y^{(k)})||<n,\ on\ \widetilde{T^kM}.  \tag{2.1.9}
\end{equation}
\end{teo}

\section{Variational Problem}

The variational problem concerning the functional $I(c)$ from (2.1.7) was
studied in the book [94], ch. 8. So we shall present here only the
corresponding results.

\begin{teo}
In order for the curve $c:t\in [0,1]\rightarrow (x^i(t))\in U\subset M$ to be
an extremal curve for the integral of action $I(c)$ it is necessary that the
following Euler-Lagrange equations hold:
\begin{equation}
\begin{array}{l}
\stackrel{\circ }{E}_i(L):=\displaystyle\frac{\partial L}{\partial x^i}-%
\displaystyle\frac d{dt}\displaystyle\frac{\partial L}{\partial y^{(1)i}}%
+\cdots +(-1)^k\displaystyle\frac 1{k!}\displaystyle\frac{d^k}{dt^k}%
\displaystyle\frac{\partial L}{\partial y^{(k)i}} = 0, \\
\\
y^{(1)i}=\displaystyle\frac{dx^i}{dt},\cdots ,y^{(k)i}=\displaystyle\frac
1{k!}\displaystyle\frac{d^kx^i}{dt^k}.
\end{array}
\tag{2.2.1}
\end{equation}
\end{teo}

An important remark regard the geometrical meaning of the system of
functions $\stackrel{\circ }{E}_i(L)$.
The following property holds:

Every $\stackrel{\circ }{E}_i(L)$ from (2.2.1) determines a
$d$-covector field along curve $c$.

This means that, with respect to a change of local coordinates on the manifold $%
T^kM$, $\stackrel{\circ }{E}_i(L)$ obeys the transformation
\begin{equation}
\stackrel{\circ }{\widetilde{E}}_i(\widetilde{L})\displaystyle\frac{\partial
\widetilde{x}^i}{\partial x^j}=\stackrel{\circ }{E}_j(L).  \tag{2.2.2}
\end{equation}

Let $c:t\in [0,1]\rightarrow (x^i(t))\in U\subset M$ be a smooth curve (or a
law of moving) parametrized by the time $t$ and $V^i(x(t))$ a differentiable
vector field on $c$.

The mapping $S_V:(x(t))\in U\rightarrow S_V^{}(x(t))\in T^kM$, defined by
\begin{equation}
\left\{
\begin{array}{l}
x^i=x^i(t),\ t\in [0,1], \\
\\
y^{(1)i}=V^i(x(t)),\ 2y^{(2)i}=\displaystyle\frac 1{1!}\displaystyle\frac{
dV^i}{dt},...,ky^{(k)i}=\displaystyle\frac 1{(k-1)!}\displaystyle\frac{
d^{k-1}V^i}{dt^{k-1}},
\end{array}
\right.  \tag{2.2.3}
\end{equation}
is a section of the projection $\pi ^k:T^kM\rightarrow M$.

It is not difficult to see that the operator
\begin{equation}
\displaystyle\frac{d_V}{dt}=V^i\displaystyle\frac \partial {\partial x^i}+ %
\displaystyle\frac{dV^i}{dt}\displaystyle\frac \partial {\partial
y^{(1)i}}+\cdots +\displaystyle\frac 1{k!}\displaystyle\frac{d^kV^i}{dt^k} %
\displaystyle\frac \partial {\partial y^{(k)i}}  \tag{2.2.4}
\end{equation}
is invariant with respect to the coordinate transformations (1.1.2).

If $V^i=\displaystyle\frac{dx^i}{dt}$, then $\displaystyle\frac{d_V}{dt}= %
\displaystyle\frac d{dt}$ and for any Lagrangian $L(x,y^{(1)},...,y^{(k)})$,
$\displaystyle\frac{d_VL}{dt}$ is a scalar field.

The action of the $k$-tangent structure $J$ (cf. \S 2.4, ch. 1) on the
operator $\displaystyle\frac{d_V}{dt}$ leads to k new operators:
\begin{equation}
I_V^k=J\left( \displaystyle\frac{d_V}{dt}\right) ,\ I_V^{k-1}=J\left(
I_V^k\right) ,...,I_V^1=J\left( I_V^2\right) ,\ 0=J\left( I_V^1\right) ,
\tag{2.2.5}
\end{equation}
where
\begin{equation}
I_V^k=V^i\displaystyle\frac \partial {\partial y^{(1)i}}+\displaystyle\frac{
dV^i}{dt}\displaystyle\frac \partial {\partial y^{(2)i}}+\cdots + %
\displaystyle\frac 1{(k-1)!}\displaystyle\frac{d^{k-1}V^i}{dt^{k-1}} %
\displaystyle\frac \partial {\partial y^{(k)i}}.  \tag{2.2.5'}
\end{equation}

But the previous operators are vector fields.
Consequently, $I_V^1(L)$, ..., $I_V^k(L)$ are scalar fields.
For $V^i=\displaystyle\frac{dx^i}{dt}$ they coincide with the main invariants $I^1(L)$, ..., $I^k(L)$.

The relations between the operators $\displaystyle\frac{d_V}{dt}$, $I_V^1$,..., $I_V^k$ are expressed cf.[94] by:

\begin{teo}
The following identities hold:
\begin{equation}
\begin{array}{lll}
\displaystyle\frac{d_VL}{dt} & = & V^i\stackrel{\circ }{E_i}(L)+%
\displaystyle \frac{d}{dt}I_V^k(L)-\displaystyle\frac 1{2!}\displaystyle%
\frac{d^2}{dt^2}I_V^{k-1}(L)+\cdots + \\
& + & (-1)^{k-1}\displaystyle\frac 1{k!}\displaystyle\frac{d^k}{dt^k}%
I_V^1(L).
\end{array}
\tag{2.2.6}
\end{equation}
\end{teo}

Besides of the covector field $\stackrel{\circ }{E_i}(L)$, Craig and Synge
introduced other covectors, important in the variational problem of the
integral of action. They are denoted by $\stackrel{1}{E_i},.., \stackrel{k}{E_i}$ and are provided by:

{\bf Lemma 2.2.1}

{\it For any differentiable Lagrangian $L(x,y^{(1)},...,y^{(k)})$ and any
differentiable function $\Phi (t)$ we have:
\begin{equation}
\stackrel{\circ }{E_i}(\Phi L)=\Phi \stackrel{\circ }{E_i}(L)+\displaystyle
\frac{d\Phi }{dt}\stackrel{1}{E_i}(L)+\cdots +\displaystyle\frac{d^k\Phi }{
dt^k}\stackrel{k}{E_i}(L),  \tag{2.2.7}
\end{equation}
where $\stackrel{1}{E_i}(L)$, ..., $\stackrel{k}{E_i}(L)$ are $d$-covector
fields. They are expressed by the actions on $L$ of the following operators:}
\begin{equation}
\begin{array}{l}
\stackrel{\circ }{E_i}=\displaystyle\frac \partial {\partial x^i}- %
\displaystyle\frac d{dt}\displaystyle\frac \partial {\partial
y^{(1)i}}+\cdots +(-1)^k\displaystyle\frac 1{k!}\displaystyle\frac{d^k}{dt^k}
\displaystyle\frac \partial {\partial y^{(k)i}}, \\
\\
\stackrel{1}{E_i}=\sum\limits_{\alpha =1}^k(-1)^\alpha
\displaystyle\frac 1{\alpha !}\left( _{\alpha -1}^\alpha \right)
\displaystyle\frac{d^{k-1}}{
dt^{k-1}}\displaystyle\frac \partial {\partial y^{(\alpha )i}}, \\
\\
\stackrel{2}{E_i}=\sum\limits_{\alpha =2}^k(-1)^\alpha
\displaystyle\frac 1{\alpha !}\left( _{\alpha -2}^\alpha \right)
\displaystyle\frac{d^{k-2}}{
dt^{k-2}}\displaystyle\frac \partial {\partial y^{(\alpha )i}}, \\
...................................................... \\
\stackrel{k}{E_i}=(-1)^k\displaystyle\frac 1{k!}\displaystyle\frac \partial
{\partial y^{(k)i}}.
\end{array}
\tag{2.2.8}
\end{equation}
In the light of the above results we get:

\begin{teo}
For any differentiable Lagrangians $L(x,y^{(1)},...,y^{(k)})$ and any function $%
F(x,y^{(1)},...,y^{(k-1)})$ the following properties hold:
\begin{equation}
\stackrel{0}{E}_i(L+\displaystyle\frac{dF}{dt})=\stackrel{0}{E}_i(L),
\tag{2.2.9}
\end{equation}

\begin{equation}
\stackrel{0}{E}_i(\displaystyle\frac{dF}{dt})=0,\ \stackrel{1}{E}_i(%
\displaystyle\frac{dF}{dt})=-\stackrel{0}{E}_i(F),\ ...,\stackrel{k}{E}_i(%
\displaystyle\frac{dF}{dt})=-\stackrel{k-1}{E}_i(F).  \tag{2.2.10}
\end{equation}
\end{teo}

A consequence of the property (2.2.9) is as follows.

\begin{teo}
The integral of action
\begin{equation}
I(c)=\int\limits_0^1Ldt,\ I^{\prime }(c)=\int\limits_0^1\left( L+%
\displaystyle\frac{dF}{dt}\right) dt  \tag{2.2.11}
\end{equation}
determine the same Euler-Lagrange equations, ($F$ depending on \newline
$(x,y^{(1)},...,y^{(k-1)})$).
\end{teo}

\section{Higher Order Energies}

In the book [94] we introduce the notion of higher order energies of a Lagrangian $L(x,y^{(1)},...,y^{(k)})$.

\begin{defi}
We call energies of order\textit{\ }$k$\textit{, }$k-1$\textit{, ..., }$1$ of
the differentiable Lagrangian $L(x,y^{(1)},...,y^{(k)})$ the following
invariants:
\begin{equation}
\begin{array}{l}
\mathcal{E}^k(L)=I^k(L)-\displaystyle\frac 1{2!}\displaystyle\frac{%
dI^{k-1}(L)}{dt^{}}+\cdots +(-1)^{k-1}\displaystyle\frac 1{k!}\displaystyle
\frac{d^{k-1}I^1(L)}{dt^{k-1}}-L, \\
\\
\mathcal{E}^{k-1}(L)=-\displaystyle\frac 1{2!}I^{k-1}(L)+\displaystyle\frac
1{3!}\displaystyle\frac{dI^{k-2}(L)}{dt}\cdots +(-1)^{k-1}\displaystyle\frac
1{k!}\displaystyle\frac{d^{k-2}I^1(L)}{dt^{k-2}}, \\
............................................................................................
\\
\mathcal{E}^1(L)=(-1)^{k-1}\displaystyle\frac 1{k!}I^1(L).
\end{array}
\tag{2.3.1}
\end{equation}
\end{defi}

As we shall see, the energies $\mathcal{E}^k(L)$, ...,
$\mathcal{E}^1(L)$ are involved in a N\"other theory of symmetries
of the higher order Lagrangians.

The next theorem is well known:

\begin{teo}
For any differentiable Lagrangian $L(x,y^{(1)},...,y^{(k)})$ the following
identity holds
\begin{equation}
\displaystyle\frac{d\mathcal{E}^k(L)}{dt}=-\stackrel{0}{E_i}(L)\displaystyle
\frac{dx^i}{dt}.  \tag{2.3.2}
\end{equation}
\end{teo}

An immediate consequence of the previous theorem is the following \textit{\
law of conservation}:

\begin{teo}
For any differentiable Lagrangian $L(x,y^{(1)},...,y^{(k)})$ the energy of
order $k$, $\mathcal{E}^k(L)$ is conserved along every extremal curve of the
Euler-Lagrange equations $\stackrel{0}{E_i}(L)=0$.
\end{teo}

Theorem 2.2.4 says that the integrals of actions $I(c)$ and $%
I^{\prime }(c)$ from (2.2.11) determine the same Euler-Lagrange equations, if $%
F$ is a Lagrangian with the property $\displaystyle\frac{\partial F}{
\partial y^{(k)i}}=0$.

\begin{defi}
A symmetry of a differentiable Lagrangian $L(x,y^{(1)},...,y^{(k)})$ is a $%
C^\infty $-diffeomorphism $\varphi :M\times \mathbf{R}\rightarrow M\times
\mathbf{R}$ which preserves the variational principle of the integral of
action $I(c)$ from (2.1.7).
\end{defi}

One can consider the notion of local symmetry of the Lagrangian $L$,
taking $\varphi $ as local diffeomorphism. If $U\times (a,b)$ is a domain of
a local chart on the manifold $M\times \mathbf{R}$, then we can express an
infinitesimal diffeomorphism $\varphi $ in the form
\begin{equation}
\begin{array}{l}
x^{\prime i}=x^i+\varepsilon V^i(x,t),\quad (i=1,...,n) \\
\\
t^{\prime }=t+\varepsilon \tau (x,t),
\end{array}
\tag{2.3.3}
\end{equation}
where $(V^i,\tau )$ is a vector field on $U\times (a,b)$ defined
along curve\\$c:t\in [0,1]\rightarrow (x^i(t),t)\in U\times (a,b)$ and $%
\varepsilon $ is a real number, sufficiently small in absolute value, so that $%
Im$ $\varphi \subset U\times (a,b)$ .

The infinitesimal transformation (1.3.3) is a symmetry for the differentiable
Lagrangian $L$ if and only if for any differentiable function\\$%
F(x,y^{(1)},...,y^{(k-1)})$ the following equation holds:

$
\begin{array}{lll}
L\left( x^{\prime },\displaystyle\frac{dx^{\prime }}{dt^{\prime }},\cdots , %
\displaystyle\frac 1{k!}\displaystyle\frac{d^kx^{^{\prime }}}{dt^{\prime k}}
\right) dt^{\prime } & = & \left\{ L\left( x,\displaystyle\frac{dx}{dt}
,\cdots ,\displaystyle\frac 1{k!}\displaystyle\frac{d^kx}{dt^k}\right)
\right. + \\
& + & \left. F\left( x,\displaystyle\frac{dx}{dt},\cdots ,\displaystyle\frac
1{(k-1)!}\displaystyle\frac{d^{k-1}x}{dt^{k-1}}\right) \right\} dt
\end{array}
$

One proves the following N\"other theorem, [94]:

\begin{teo}
For any infinitesimal symmetry (1.3.3) of a Lagrangian

$L(x,y^{(1)},...,y^{(k)})$ and for any function $\Phi
(x,y^{(1)},...,y^{(k-1)})$, the function:
\[
\begin{array}{lll}
\mathcal{F}^k(L,\Phi ) & = & I_V^k(L)-\displaystyle\frac 1{2!}\displaystyle %
\frac d{dt}I_V^{k-1}(L)+\cdots +(-1)^{k-1}\displaystyle\frac 1{k!}%
\displaystyle\frac{d^{k-1}}{dt^{k-1}}I_V^1(L)- \vspace{3mm}\\
& - & \tau \mathcal{E}^k(L)+\displaystyle\frac{d\tau }{dt}\mathcal{E}%
^{k-1}(L)+\cdots +(-1)^k\displaystyle\frac{d^{k-1}\tau }{dt^{k-1}}\mathcal{E}%
^1(L)-\Phi .
\end{array}
\]
is conserved along extremal curves of the Euler-Lagrange equations
\newline $\stackrel{\circ }{E}_i(L)=0$.
\end{teo}

In particular, if the Zermello conditions (2.1.8) are satisfied, then
the energies $\mathcal{E}^1(L)$, ..., $\mathcal{E}^k(L)$ vanish and the
previous theorem has a simpler form.

\section{Jacobi-Ostrogradski Momenta}

We introduce the Jacobi-Ostrogradski momenta and the Hamilton - Jacobi - Ostrogradski equations.
The main results on these are given without proofs, [94].

Consider the energy of order $k$, $\mathcal{E}^k(L)$ of the Lagrangian $L$
from (2.3.1). Noticing that $\mathcal{E}^k(L)$ is a polynomial
function of degree one in $\displaystyle\frac{dx^i}{dt}$, ..., $%
\displaystyle
\frac{d^kx^i}{dt^k}$ we can write:
\begin{equation}
\mathcal{E}^k(L)=p_{(1)i}\displaystyle\frac{dx^i}{dt}+p_{(2)i}\displaystyle
\frac{d^2x^i}{dt^2}+\cdots +p_{(k)i}\displaystyle\frac{d^kx^i}{dt^k}-L,
\tag{2.4.1}
\end{equation}
where
\begin{equation}
\begin{array}{l}
p_{(1)i}=\displaystyle\frac{\partial L}{\partial y^{(1)i}}-\displaystyle %
\frac 1{2!}\displaystyle\frac d{dt}\displaystyle\frac{\partial L}{\partial
y^{(2)i}}+\cdots +(-1)^{k-1}\displaystyle\frac 1{k!}\displaystyle\frac{
d^{k-1}}{dt^{k-1}}\displaystyle\frac{\partial L}{\partial y^{(k)i}}, \\
\\
p_{(2)i}=\displaystyle\frac 1{2!}\displaystyle\frac{\partial L}{\partial
y^{(2)i}}-\displaystyle\frac 1{3!}\displaystyle\frac d{dt}\displaystyle\frac{
\partial L}{\partial y^{(3)i}}+\cdots +(-1)^{k-2}\displaystyle\frac 1{k!} %
\displaystyle\frac{d^{k-2}}{dt^{k-2}}\displaystyle\frac{\partial L}{\partial
y^{(k)i}}, \\
...............................................................................................
\\
p_{(k)i}=\displaystyle\frac 1{k!}\displaystyle\frac{\partial L}{\partial
y^{(k)i}}.
\end{array}
\tag{2.4.2}
\end{equation}
$p_{(1)i}$, ..., $p_{(k)i}$ are called \textit{the Jacobi-Ostrogradski
momenta. }

M. de L\'eon and others, [73, 74] established the following important result.

\begin{teo}
Along extremal curves of Euler-Lagrange equations, \newline
$\stackrel{\circ }{E}_i(L)=0$, the following
Hamilton-Jacobi-Ostrogradski equations hold:
\begin{equation}
\begin{array}{l}
\displaystyle\frac{\partial \mathcal{E}^k(L)}{\partial p_{(\alpha )i}}=%
\displaystyle\frac{d^\alpha x^i}{dt^\alpha },\quad (\alpha =1,...,k), \\
\\
\displaystyle\frac{\partial \mathcal{E}^k(L)}{\partial x^i}=-\displaystyle
\frac{dp_{(1)i}}{dt}, \\
\\
\displaystyle\frac{\partial \mathcal{E}^k(L)}{\partial y^{(\alpha )i}}%
=-\alpha !\displaystyle\frac{dp_{(\alpha +1)i}}{dt},\quad (\alpha
=1,...,k-1).
\end{array}
\tag{2.4.3}
\end{equation}
\end{teo}

The Jacobi-Ostrogradski momenta $p_{(1)i}$, ..., $p_{(k)i}$ allow the
introduction of the $1$-forms:
\begin{equation}
\begin{array}{l}
p_{(1)}=p_{(1)i}dx^i+p_{(2)i}dy^{(1)i}+\cdots +p_{(k)i}dy^{(k-1)i}, \\
\\
p_{(2)}=p_{(2)i}dx^i+p_{(3)i}dy^{(1)i}+\cdots +p_{(k)i}dy^{(k-2)i}, \\
...........................................................................
\\
p_{(k)}=p_{(k)i}dx^i.
\end{array}
\tag{2.4.4}
\end{equation}

The following properties hold:
\[
J^{*}p_{(1)}=p_{(2)},\ ....,\ J^{*}p_{(k-1)}=p_{(k)}.
\]

\section{Higher Order Lagrange Spaces}

The notion of Lagrange space of order $k$ is a natural extension of that of
classical Lagrange space $L^n=(M,L(x,y))$. It was introduced by the
author of this monograph, in the papers quoted in his book [94].
It was introduced by the author in the paper quoted in his book {\it The
Geometry of Higher Order Lagrange Spaces. Applications to Mechanics and
Physics}, (Kluwer, FTPH no. 82) entirely devoted to this subject.

Here we give, without prooofs, the main results from the
geometry of higher order Lagrange spaces, [94].

We call a \textit{Lagrange space of order }$k$ a pair $L^{(k)n}=(M,L)$
formed by a real $n$-dimensional manifold $M$ and a differentiable
Lagrangian of order $k$, $L:(x,y^{(1)},...,y^{(k)})\in T^kM\rightarrow
L(x,y^{(1)},...,y^{(k)})\in \mathbf{R}$ for which the Hessian with the
entries
$$
g_{ij}(x,y^{(1)},...,y^{(k)})=\displaystyle\frac 12\displaystyle\frac{
\partial ^2L}{\partial y^{(k)i}\partial y^{(k)j}},\ on\ \widetilde{T^kM},
$$
has the property
$$
rank\ ||g_{ij}||=n
$$
and the quadratic form $\Psi =g_{ij}\xi ^i\xi ^j$ on $\widetilde{T^kM}$ has
constant signature.

$L$ is called the \textit{fundamental function} and $g_{ij}$ the \textit{\
fundamental (or metric) tensor field} of the space $L^{(k)n}$.

The geometry of the pair $(T^kM,L(x,y^{(1)},...,y^{(k)}))$ is called the
geometry of the space $L^{(k)}=(M,L)$.

Notice that the geometry of Lagrange space of order $k$ is not coincident
with the geometrization of the Lagrangians of order $k$, $%
L(x,y^{(1)},...,y^{(k)})$, which could be degenerate, i.e. $rank\
||g_{ij}||<n$.

We shall study this geometry using the methods suggested by the Higher Order
Lagrangian Mechanics [35] and by the geometry of higher order Finsler spaces [95].
Consequently, we determine a canonical semispray $S$ and derive the main
geometrical object field of the space $L^{(k)n}$ by means of $S$.

Consider the integral of action of the Lagrangian $L$, the fundamental function
of the Lagrange space of order $k$, $L^{(k)n}=(M,L)$ and determine the
covector fields $\stackrel{0}{E_i}(L)$, ..., $\stackrel{k-1}{E_i}(L)$, $%
\stackrel{k}{E_i}(L)$.

Then $\stackrel{0}{E_i}(L)=0$ are just the Euler-Lagrange equations. Thus, $%
\mathcal{E}^k(L)$, ..., $\mathcal{E}^1(L)$ from (2.3.1) give us
the {\it energies} of the space $L^{(k)n}$. Theorems 2.3.1 and
2.3.2 can be applied and Theorem 2.3.3 gives the infinitesimal
symmetries of the considered space.

The following result is known, [94]:

\begin{teo}
The equations $g^{ij}\stackrel{k-1}{E}_i(L)=0$ determine the $k$-semispray
\begin{equation}
S=y^{(1)i}\displaystyle\frac \partial {\partial x^i}+2y^{(2)i}\displaystyle %
\frac \partial {\partial y^{(1)i}}+\cdots +ky^{(k)}\displaystyle\frac
\partial {\partial y^{(k-1)i}}-(k+1)G^i\displaystyle\frac \partial {\partial
y^{(k)i}}  \tag{2.5.1}
\end{equation}
with the coefficients
\begin{equation}
(k+1)G^i=\displaystyle\frac 12g^{ij}\left\{ \Gamma \left( \displaystyle\frac{%
\partial L}{\partial y^{(k)j}}\right) -\displaystyle\frac{\partial L}{%
\partial y^{(k-1)j}}\right\} ,  \tag{2.5.2}
\end{equation}
$\Gamma $ being the operator (1.2.3).
\end{teo}

The $k$-semispray $S$ depends on the fundamental function $L$ of the space $%
L^{(k)n}$. It will be called \textit{canonical}. If $L$ is globally defined
on $T^kM$, then $S$ has the same property on $\widetilde{T^kM}$.

Taking into account Theorem 1.5.1 we find, [94]:

\begin{teo}
The set of functions
\begin{equation}
\begin{array}{l}
\underset{(1)}{M_j^i}=\displaystyle\frac{\partial G^i}{\partial y^{(k)j}}%
,\ \underset{(2)}{M_j^i}=\displaystyle\frac 12\left( S\underset{(1)}{%
M_j^i}+\underset{(1)}{M_m^i}\underset{(1)}{M_j^m}\right) ,\ ..., \\
\\
\underset{(k)}{M_j^i}=\displaystyle\frac 1k\left( S\underset{(k-1)}{M_j^i%
}+\underset{(1)}{M_m^i}\underset{(k-1)}{M_j^m}\right)
\end{array}
\tag{2.5.3}
\end{equation}
are the dual coefficients of a nonlinear connection $N$ determined only by
the canonical semispray $S$.
\end{teo}
$N$ is called canonical, too.

Theorem 1.5.2 give us new dual coefficients $\underset{(1)}{M_{\ j}^{*i}} ,...,\underset{(k)}{M_{\ j}^{*i}}$, (I. Bucataru, [26]).

The coefficients $\underset{(1)}{N_j^i},..., \underset{(k)}{N_j^i}$ of $N$ are obtained
from (1.4.3').

The existence of the spaces $L^{\left( k\right) n}=(M,L),$ when $M$ is a
paracompact manifold is assured by the following examples.

{\it \underline{Example 5.1}}
Let $g_{ij}(x)$ be a Riemannian tensor on $M$ and the nonlinear connection $%
N $ with the dual coefficients $\underset{(1)}{M_j^i}(x,y^{\left(
1\right) }),...\underset{(k)}{M_j^i}(x,y^{\left( 1\right)
},...y^{\left( k\right) }) $ given by (1.5.3). Consider the Liouville
d-vector field, $z^{\left( k\right) i\text{ }}$(Th. 21.4.1):
\begin{equation}
kz^{\left( k\right) i}=ky^{\left( k\right)
i}+(k-1)\underset{(1)}{M_m^i} y^{\left( k-1\right)
m}+...+\underset{(k-1)}{M_m^i}y^{\left( 1\right) m} \tag{2.5.4}
\end{equation}
and remark that $z^{\left( k\right) i}$ is a $d$-vector field, linear in the vertical
variables $y^{\left( k\right) i}.$ Consequence, the function
\begin{equation}
L(x,y^{\left( 1\right) },...y^{\left( k\right) })=g_{ij}(x)z^{\left( k\right)
i}z^{\left( k\right) j}  \tag{2.5.5}
\end{equation}
is a fundamental function of a Lagrange space $L^{\left( k\right) n}.$ Its
fundamental tensor field is exactly $g_{ij}(x).$

{\it \underline{Example 5.2}}
Let $\stackrel{\circ }{L}(x,y^{\left( 1\right) })$
be the Lagrangian of electrodynamics
\begin{equation}
\stackrel{\circ }{L}(x,y^{\left( 1\right) })=mc\gamma _{ij}(x)y^{\left(
1\right) i}y^{\left( 1\right) j}+\displaystyle\frac{2e}mb_i(x)y^{\left(
1\right) i}  \tag{2.5.6}
\end{equation}
$m, c, e$ being the well known physical constants, $\gamma _{ij}(x)$ the
gravitational potentials and $b_i(x)$ the electromagnetic potentials.

Let us consider the nonlinear connection $N$ determined like in previous
example. Then
\begin{equation}
L(x,y^{\left( 1\right) },...,y^{\left( k\right) })=mc\gamma
_{ij}(x)z^{\left( k\right) i}z^{\left( k\right) j}+\displaystyle\frac{2e}
mb_i(x)z^{\left( k\right) i}  \tag{2.5.7}
\end{equation}
is a fundamental function of a Lagrange space of order $k,$ $L^{\left(
k\right) n}=(M,L).$

Evidently $L$ from (2.5.7) is a particular case of the Lagrangian from (2.1.5).
It will be called the \textit{prolongation }of order $k$ of the
electrodynamic Lagrangian in (2.5.6).

In the end of this chapter we will study the geometry of the Lagrange
space $L^{\left( k\right) n},$ with fundamental function (2.5.7).

The adapted basis $\left\{ \displaystyle\frac \delta {\delta x^i}, %
\displaystyle\frac \delta {\delta y^{\left( 1\right)
i}},...,\displaystyle \frac \delta {\delta y^{\left( k\right)
i}}\right\} $ determined by the canonical nonlinear connection $N$
is given by (2.3.9) and its dual basis\\ $\{\delta x^i,\delta
y^{\left( 1\right) i},...,\delta y^{\left( k\right) i}\}$ is
expressed in the formulae (2.4.2).

The horizontal curves of the space $L^{\left( k\right) n}$ with respect to
canonical nonlinear connection $N$ are characterized by the system of
differential equations
\[
\displaystyle\frac{\delta y^{\left( 1\right) i}}{dt}=\cdots =\displaystyle
\frac{\delta y^{\left( k\right) i}}{dt}=0.
\]

In the light of these properties we determine the autoparallel curves of $N$
by adding the conditions to the previous differential equations
\[
y^{\left( 1\right) i}=\displaystyle\frac{dx^i}{dt},...,y^{\left( k\right)
i}= \displaystyle\frac 1{k!}\displaystyle\frac{d^kx^i}{dt^k}.
\]

Using the results from chapter 1 one can proves:

The canonical nonlinear connection $N$ is integrable if, and only if the
following equations hold
\[
\underset{(01)}{R_{jk}^i}=\cdots =\underset{(0k)}{R_{jk}^i}=0.
\]

\section{Canonical Metrical $N$-Connections}

Consider the canonical nonlinear connection $N$ of the Lagrange space of
order $k$, $L^{\left( k\right) n}=(M,L).$ A linear connection $D$ on $T^kM$
is called an $N$-linear connection if:

\vspace*{3mm}

1) $D$ preserves by parallelism the horizontal distribution $N$,

\vspace*{3mm}

2) $DJ=0$.

\vspace*{3mm}

The coefficients of $D$ with respect to the adapted basis denoted by $%
D\Gamma (N)=(L_{jh}^i,\underset{(1)}{C_{jh}^i},...,\underset{(k)}{
C_{jh}^i}),$ can be uniquely determined if $D$ is compatible with
the Riemannian (or pseudoriemannian) metric $\Bbb{G}$ on $T^kM$,
determined by the fundamental tensor field $g_{ij}$ of $L^{(k)n}$.

Namely, $\Bbb{G}$ is expressed in the adapted basis by
\begin{equation}
\Bbb{G}=g_{ij}dx^i\otimes dx^j+g_{ij}\delta y^{\left( 1\right) i}\otimes
\delta y^{\left( 1\right) j}+\cdots +g_{ij}\delta y^{\left( k\right)
i}\otimes \delta y^{\left( k\right) j}.  \tag{2.6.1}
\end{equation}

$D$ is compatible with $\Bbb{G}$ if
\[
D_X\Bbb{G}=0,\forall X\in \mathcal{X}(\widetilde{T^kM}).
\]

\begin{teo}

The following properties hold:

1) There exists a unique $N$-linear connection $D$ on $\widetilde{T^kM}$
verifying the axioms:
\begin{equation}
g_{ij|h}=0,\ g_{ij}\stackrel{(1)}{|}_h=\cdots =g_{ij}\stackrel{(k)}{|}_h=0,
\tag{2.6.2}
\end{equation}
\begin{equation}
L_{jh}^i=L_{hj}^i,\ \underset{(\alpha )}{C_{jh}^i}=\underset{(\alpha )}{%
C_{hj}^i},\ (\alpha =1,...k).  \tag{2.6.2'}
\end{equation}

2) The coefficients $C\Gamma (N)=(L_{jh}^i,\underset{(1)}{C_{jh}^i},...,%
\underset{(k)}{C_{jh}^i})$ of $D$ are given by the generalized
Christoffel symbols:
\begin{equation}
\left\{
\begin{array}{l}
L_{ij}^m=\displaystyle\frac 12g^{ms}(\displaystyle\frac{\delta g_{is}}{%
\delta x^j}+\displaystyle\frac{\delta g_{sj}}{\delta x^i}-\displaystyle\frac{%
\delta g_{ij}}{\delta x^s}), \\
\\
\underset{\left( \alpha \right) }{C_{ij}^m}=\displaystyle\frac 12g^{ms}(%
\displaystyle\frac{\delta g_{is}}{\delta y^{\left( \alpha \right) j}}+%
\displaystyle\frac{\delta g_{sj}}{\delta y^{\left( \alpha \right) i}}-%

\displaystyle\frac{\delta g_{ij}}{\delta y^{\left( \alpha \right) s}}),\
(\alpha =1,...,k).
\end{array}
\right.   \tag{2.6.3}
\end{equation}

3) $D$ depends only on the fundamental function $L$ of the space $L^{\left(
k\right) n}$.
\end{teo}

The connection $D$ from the previous theorem is called the \textit{canonical
metrical} $N$-connection of the Lagrange space of order $k,$ $L^{\left(
k\right) n}.$

The d-tensors of curvature of $D$ satisfy the following identities:

\[
\begin{array}{l}

g_{sj}R_{i\ hm}^s+g_{is}R_{j\ hm}^s=0, \\
\\
g_{sj}\underset{}{\underset{(\alpha )}{P}\text{}_{i\ hm}^s}+g_{is}
\underset{(\alpha )}{P}\text{}_{j\ hm}^s=0, \\
\\
g_{sj}\underset{}{\underset{(\alpha \beta )}{S}\text{}_{i\ hm}^s}+g_{is} \underset{(\alpha \beta )}{S}\text{}_{j\ hm}^s=0.
\end{array}
\]

The connection 1-forms of $D$ are

\begin{equation}
\omega _j^i=L_{jh}^idx^h+\stackrel{k}{\underset{\alpha =1}{\sum }}
\underset{(\alpha )}{C_{jh}^i}\delta y^{\left( \alpha \right) h}.
\tag{2.6.4}
\end{equation}

Finally, the structure equations of the canonical metrical $N$-connection $D$
are given by Theorem 1.7.1 in the conditions (2.6.3).

The Bianchi identities of $D$ are obtained from the structure equations
applying the operator of exterior differentiation and taking into account
the same equations of structure.

Now let us consider the tensor field
\begin{equation}
\Bbb{F}\mathbf{=-}\displaystyle\frac \delta {\delta y^{\left( k\right)
i}}\otimes dx^i+\displaystyle\frac \delta {\delta x^i}\otimes \delta
y^{\left( k\right) i}.  \tag{2.6.5}
\end{equation}
We can see that $\Bbb{F}$ is globally defined on $\widetilde{T^kM}$ and $%
\Bbb{F}^3+\Bbb{F}=0.$ The pair $(\Bbb{G}, \Bbb{F})$ determines a Riemannian almost $%
(k-1)n$-contact structure on $\widetilde{T^kM}$ depending only on then
fundamental function $L$ of the space $L^{\left( k\right) n}.$

{\bf Examples.}

1$^{\circ }$. $\mathcal{R}^{\left( k\right) n}=$Prol$^k\mathcal{R}^n.$

Let $N$ be the nonlinear connection with the dual coefficients (2.5.3).
It is determined only by $g_{ij}(x).$ If $\{\delta x^i,\delta y^{\left(
1\right) i},...,\delta y^{\left( k\right) i}\}$ is the adapted cobasis to $N$
and $V_1,$ then the formula (2.6.1) gives us a Riemannian structure $\Bbb{G}$
on $\widetilde{T^kM}$ which depend on $g_{ij}(x),$ only. The space $\mathcal{%
\ R}^{\left( k\right) n}=(\widetilde{T^kM},\Bbb{G})$ is called the prolongation to $%
\widetilde{T^kM}$ of the Riemannian structure $(M,g_{ij}(x)).$

Consider the $k$-Liouville vector field $z^{\left( k\right) i}$ from (2.5.4)
constructed by means of (2.5.3).

In this case the pair $L^{\left( k\right) n}=(M,L(x,y^{\left( 1\right)
},...,y^{\left( k\right) })$ is a Lagrange space of order $k$, where $L=g_{ij}(x)z^{(k)i}z^{(k)j}$.

The fundamental tensor field of $L^{\left( k\right) n}$ is exactly $%
g_{ij}(x) $, because $z^{\left( k\right) i}$ is of the form $z^{\left(
k\right) i}=y^{\left( k\right) i}+\lambda ^i(x,y^{\left( 1\right)
},...,y^{(k-1)})$. The canonical metrical $N$-connection has the coefficients
\[
L_{jh}^i=\gamma _{jh}^i(x),\ \underset{(\alpha )}{C_{jh}^i}=0,\
(\alpha =1,...,k).
\]

2$^{\circ }.$ $\mathcal{F}^{\left( k\right) n}=$Prol$^kF^n$. The
prolongation of order $k$ of a Finsler space $F^n=(M,F(x,y^{\left( 1\right)
}))$ leads to a second example of Lagrange space of order $k$.

Consider $N$ the nonlinear connection with the dual coefficients (1.5.4). It is determined only by $F^n$. The $N$-lift of the fundamental tensor $%
g_{ij}(x,y^{\left( 1\right) })$ of the space $F^n$ is given by (2.6.1). It is a Riemannian metric on $\widetilde{T^kM}$ determined by $F^n$ only. Thus,
\[
L=g_{ij}(x,y^{\left( 1\right) })z^{\left( k\right) i}z^{\left( k\right) j},
\]
$z^{\left( k\right) i}$ being the $k$-Liouville vector fields corresponding
to $N$, is the fundamental function of a space $L^{\left( k\right) n}.$

The fundamental tensor of space is $g_{ij}(x,y^{\left( 1\right) i})$ and the
canonical metrical $N$-connection has the coefficients
\[
L_{jk}^i=F_{jk}^i(x,y^{\left( 1\right) }),\
\underset{(1)}{C_{jk}^i} =C_{jk}^i,\ \underset{(\alpha
)}{C_{jk}^i}=0,\ (\alpha =2,...,k-1).
\]

\section{Generalized Lagrange Spaces of Order $k$}

The notion of generalized Lagrange spaces of order $k$ is a natural
extension of that of space $L^{(k)n}$. It was used, [94], in the geometrical
theory of the higher order relativistic optics.

\begin{defi}
A generalized Lagrange spaces of order $k$ is a pair $%
GL^{(k)n}=(M,g_{ij}(x,y^{(1)},...,y^{(k)}))$ formed by a real $n$
-dimensional differentiable manifold $M$ and a differentiable symmetric $d$-tensor
field $g_{ij}$ defined on $\widetilde{T^kM}$, having two
properties:

a. $g_{ij}$ has a constant signature on $\widetilde{T^kM}$;

b. $rank$ $||g_{ij}||=n$ on $\widetilde{T^kM}$.
\end{defi}

Of course, $GL^{(k)n}$ can be defined locally, in the case when $g_{ij}$ is
given on an open set $(\pi ^k)^{-1}(U)$, $U\subset M$.

We say that $g_{ij}$ is the \textit{fundamental tensor} of the space $%
GL^{(k)n}$.

Notice that any Lagrange space of order $k$, $L^{(k)n}=(M,L)$ is a
generalized Lagrange space $GL^{(k)n}=(M,g_{ij})$ with
\begin{equation}
g_{ij}=\displaystyle\frac 12\displaystyle\frac{\partial ^2L}{\partial
y^{(k)i}\partial y^{(k)j}}.  \tag{2.7.1}
\end{equation}

But not and conversely. Indeed, is possible that the system of differential
partial equations (2.7.1) does not admit any solution $L(x,y^{(1)},...,y^{(k)})$
when the $d$-tensor field $g_{ij}(x,y^{(1)},...,y^{(k)})$ is apriori given.

Let us consider the tensor field:
\begin{equation}
\underset{(k)}{C}\text{ }_{ijh}=\displaystyle\frac
12\displaystyle\frac{
\partial g_{ij}}{\partial y^{(k)h}}.  \tag{2.7.2}
\end{equation}
It is not hard to see that $\underset{(k)}{C}$ $_{ijh}$ is a covariant of order three
$d$-tensor field .

Easily follows:

{\bf Proposition 2.7.1}

{\it A necessary condition so that the system of differential partial
equation (2.7.1) admits a solution $L(x,y^{(1)},...,y^{(k)})$ is that
the $d$-tensor field $\underset{(k)}{C}$ $_{ijh}$ be
completely symmetric.}

If the tensor $\underset{(k)}{C}$ $_{ijh}$ is not completely
symmetric we say that the space $GL^{(k)n}$ is not reducible to a
Lagrange space of order $k$.

\vspace*{3mm}

\textbf{Example 2.7.1} Let $\mathcal{R}^n=(M,\gamma _{ij}(x))$ be a
Riemannian space and \\
$\sigma (x,y^{(1)},...,y^{(k)})$ a smooth function on $T^kM$ with the
property $\displaystyle\frac{\partial \sigma }{\partial y^{(k)i}}\neq 0$.

Consider the $d$-tensor field on $T^kM$:
\begin{equation}
g_{ij}(x,y^{(1)},...,y^{(k)})=e^{2\sigma (x,y^{(1)},...,y^{(k)})}(\gamma
_{ij}\circ \pi ^k)(x,y^{(1)},...,y^{(k)}).  \tag{2.7.3}
\end{equation}

We can prove that the pair $GL^{(k)n}=(M,g_{ij})$ with $%
g_{ij}$ from (2.7.3) is a generalized Lagrange space of order $k,$ which is not
reducible to a Lagrange space.

This example proves the existence of spaces $GL^{(k)n}$ on the paracompact
manifolds. For $k=1$, this space was introduced by R. Miron and R. Tavakol [100], [114].

\vspace*{3mm}

\textbf{Example 2.7.2} Consider again $\mathcal{R}^n=(M,\gamma _{ij}(x))$ and the Liouville $d$-vector field $z^{(k)i}$ of the
space $Prol^k\mathcal{R}^n$. It is expressed by
\begin{equation}
kz^{(k)i}=ky^{(k)i}+(k-1)\underset{(1)}{M_m^i}y^{(k-1)m}+\cdots +
\underset{(k-1)}{M_m^i}y^{(1)m},  \tag{2.7.4}
\end{equation}
where the dual coefficients $\underset{(1)}{M_j^i}$, ...., $\underset{(k-1)}{M_j^i}$ are
given by the formulae (2.5.3). Evidently $z^{(k)i}$ linearly depends on $%
y^{(k)i}$. So that its covariant $z_i^{(k)}=$ $\gamma _{ij}z^{(k)j}$
linearly depends on $y^{(k)i}$.

Assuming that there exists a function $n(x,y^{(1)},...,y^{(k)})\geq 1$ on $T^kM$,
we can construct the following $d$-tensor field
\begin{equation}
g_{ij}(x,y^{(1)},...,y^{(k)})=\gamma _{ij}(x)+\left( 1-\displaystyle\frac
1{n^2(x,y^{(1)},...,y^{(k)})}\right) z_i^{(k)}z_j^{(k)}.  \tag{2.7.5}
\end{equation}

One proves that the pair $GL^{(k)n}=(M,g_{ij})$, with $g_{ij}$ from
(2.7.5), is a generalized Lagrange space of order $k$, which is not reducible
to a Lagrange space of order $k$.

In the case $k=1$, this space was studied by the author [91, 100] and applied,
together with T. Kawaguchi [104] in the relativistic geometrical optics.

In the previous two examples we have a natural canonical nonlinear
connection, with the dual coefficients (2.5.3), determined by the space
Prol$^k\mathcal{R}^n$.

Returning to the generalized Lagrange spaces of order $k$, $GL^{(k)n}$, we
remark the difficulty to find a nonlinear connection derived only from the
fundamental tensor $g_{ij}$ of the space.

Therefore we assume that a nonlinear connection $N$ on $\widetilde{T^kM}$ is
apriori given. Thus we, shall study the pair $(N,GL^{(k)n})$ using the same
methods like in the geometry of the space $L^{(k)n}$.

Indeed, considering the adapted basis $\left( \displaystyle\frac \delta
{\delta x^i},\displaystyle\frac \delta {\delta y^{(1)i}},\cdots , %
\displaystyle\frac \delta {\delta y^{(k)i}}\right) $ and its dual basis $%
\left( \delta x^i,\delta y^{(1)i},...,\delta y^{(k)i}\right) $, determined
by the nonlinear connection, we define the $N$-lift of the fundamental
tensor $g_{ij}$:

\begin{equation}
\Bbb{G}=g_{ij}dx^i\otimes dx^j+g_{ij}\delta y^{(1)i}\otimes \delta
y^{(1)j}+\cdots +g_{ij}\delta y^{(k)i}\otimes \delta y^{(k)j}.  \tag{2.7.6}
\end{equation}

A $N$-linear connection compatible to $G$ is furnished by the following

theorem:

\begin{teo}
1$^0$ There exists an unique $N$-linear connection $D$ for which
\begin{equation}
\begin{array}{l}
g_{ij|h}=0,\ g_{ij}\stackrel{(\alpha )}{|_h}=0,\ (\alpha =1,...,k), \\
\\
L_{jh}^i=L_{hj}^i,\ \underset{(\alpha )}{C}\text{ }_{jh}^i=\underset{%
(\alpha )}{C}\text{ }_{hj}^i,\ (\alpha =1,...,k).
\end{array}
\tag{2.7.7}
\end{equation}

2$^0$ The coefficients $C\Gamma (N)=\left(
L_{jh}^i,\underset{(1)}{C}\text{
}_{jh}^i,...,\underset{(k)}{C}\text{ }_{jh}^i\right) $ of $D$ are
given by the generalized Christoffel symbols:
\begin{equation}
\begin{array}{l}
L_{ij}^m=\displaystyle\frac 12g^{ms}\left( \displaystyle\frac{\delta g_{is}}{%
\delta x^j}+\displaystyle\frac{\delta g_{sj}}{\delta x^i}-\displaystyle\frac{%
\delta g_{ij}}{\delta x^s}\right) , \\
\\
\underset{(\alpha )}{C}\text{ }_{ij}^m=\displaystyle\frac 12g^{ms}\left( %
\displaystyle\frac{\delta g_{is}}{\delta y^{(\alpha )j}}+\displaystyle\frac{%
\delta g_{sj}}{\delta y^{(\alpha )i}}-\displaystyle\frac{\delta g_{ij}}{%
\delta y^{(\alpha )s}}\right) ,(\alpha =1,...,k).\
\end{array}
\tag{2.7.7'}
\end{equation}
\end{teo}

The previous $N$-linear connection of the space $GL^{(k)n}$ is called
metrical canonical $N$-linear connection.

The structure equations can be written, exactly as in the chapter 1 (theorem
1.7.1).

The tensor $\Bbb{F}$ on $\widetilde{T^kM}$:
\begin{equation}
\Bbb{F}=-\displaystyle\frac \delta {\delta y^{(k)i}}\otimes dx^i+ %
\displaystyle\frac \delta {\delta x^i}\otimes \delta y^{(k)i},  \tag{2.7.8}
\end{equation}
together with the metric tensor $\Bbb{G}$ from (2.7.6) determine a Riemannian $%
(k-1)n$-almost contact structure on $\widetilde{T^kM}$. It is 'the
geometrical model' of the generalized Lagrange space of order $k$, $%
GL^{(k)n}=(M,g_{ij})$.

\chapter{Finsler Spaces of Order $k$}

\markboth{\it{THE GEOMETRY OF HIGHER-ORDER HAMILTON SPACES\ \ \ \ \ }}{\it{Finsler Spaces of Order} $k$}

The geometry of Finsler spaces of order $k$, introduced by the
author and presented in his book 'The geometry of Higher- Order
Finsler Spaces' Hadronic Press, 1998 is a natural extension to
$T^kM$ of the classical theory of Finsler Spaces. The impact of
this geometry in Differential Geometry, Variational Calculus,
Analytical Mechanics or Theoretical Physics is decisive. Finsler
spaces play a role in applications to Biology, Engineering,
Physics or Optimal Control. Also the introduction of the notion of
Finsler space of order $k$ is demanded by the solution of problem
of prolongation to $T^kM$ of the Riemannian or Finslerian
structures defined on the base manifold $M.$

In the present chapter we will develop the geometrical theory of the Finsler
spaces of order $k$, based on the geometry of Lagrange spaces of the same order.
Such that the Finsler spaces $F^{(k)n}$ form a subclass of
the class of spaces $L^{\left( k\right) n}$. Consequently we obtain an
extension of the sequence
$\{R^n\}\subset \{F^n\}\subset \{L^n\}\subset \{GL^n\}$ to the following sequence of the spaces of order $k:$
$\{R^{\left( k\right) n}\}\subset \{F^{\left( k\right) n}\}\subset
\{L^{\left( k\right) n}\}\subset \{GL^{\left( k\right) n}\}.$

In the next chapter we shall investigate the dual of this
sequence, made by the Hamilton spaces of order $k.$

\section{Spaces $F^{\left( k \right) n}$}

In order to introduce the notion of Finsler space of order $k$ there are necessary
some preliminaries. A functions $f:T^kM\rightarrow R,$ of $C^\infty $ class
on $\widetilde{T^kM}$ and continuous on the null sections of the mapping $%
\pi ^k:T^kM\rightarrow M$ is called homogeneous of degree $r\in Z$ on the
fibres on $T^kM$ (briefly $r$-homogeneous) if for any positive constant $a,$ we have:
\begin{equation}
f(x,ay^{\left( 1\right) },a^2y^{\left( 2\right) },...,a^ky^{\left( k\right)
})=a^rf(x,y^{\left( 1\right) },...,y^{(k)}).  \tag{3.1.1}
\end{equation}

An Euler theorem holds:

{\it A function $f\in \mathcal{F}(T^kM)$, differentiable on $\widetilde{T^kM}$
and continuous on the null section of $\pi ^k$ is $r-$homogeneous if and
only if:}
\begin{equation}
\mathcal{L}_{\stackrel{k}{\Gamma }}f=rf , \tag{3.1.2}
\end{equation}
$\mathcal{L}_{\stackrel{k}{\Gamma }}$ being the Lie operator of derivation
with respect to the Lioville vector field $\stackrel{k}{\Gamma },$ i.e.:
\begin{equation}
\mathcal{L}_{\stackrel{k}{\Gamma }}f=\stackrel{k}{\Gamma }f=y^{\left(
1\right) i}\displaystyle\frac{\partial f}{\partial y^{\left( 1\right) i}}
+...+ky^{\left( k\right) i}\displaystyle\frac{\partial f}{\partial y^{\left(
k\right) i}}.  \tag{3.1.2'}
\end{equation}

Notice the following property:

{\it If the function $f\in \mathcal{F}(T^kM)$ is differentiable on $T^kM$ (
inclusive in the points $(x,0,...,0))$ and $k-$homogeneous then $f$ is a polynom of order $k$ in the variable $y^{\left( 1\right) i},y^{\left(2\right) i},...,y^{\left( k\right) i}.$}

The notion of homogeneity can be extended to the vector fields and $1-$form
fields on $T^kM.$ One proves:

A vector field $X\in \mathcal{X}(T^kM)$ is $r$-homogeneous if and only if
\begin{equation}
\mathcal{L}_{\stackrel{k}{\Gamma }}X=(r-1)X . \tag{3.1.3}
\end{equation}

Of course, $\mathcal{L}_{\stackrel{k}{\Gamma }}X=[\stackrel{k}{\Gamma },X].$

For instance, the vector fields $\displaystyle\frac \partial {\partial x^i}, %
\displaystyle\frac \partial {\partial y^{\left( 1\right) i}},..., %
\displaystyle\frac \partial {\partial y^{\left( k\right) i}}$ are $%
1,0,...,1-k$ homogeneous, respectively.

An homogeneous $k$- semispray $S$ is called a $k-$ spray. The following
sentences hold:

1$^{\circ }.$ A $k-$ spray $S$ given by
\begin{equation}
S=y^{\left( 1\right) i}\displaystyle\frac{\partial f}{\partial x^i}
+...+ky^{\left( k\right) i}\displaystyle\frac{\partial f}{\partial y^{\left(
k-1\right) i}}-(k+1)G^i\displaystyle\frac \partial {\partial y^{(k)i}}
\tag{3.1.4}
\end{equation}
is $2$-homogeneous if and only if its coefficients $G^i$ are $k+1$-
homogeneous.

2$^{\circ }.$ The dual coefficients
$\underset{}{\underset{(1)}{M_j^i},}
\underset{}{\underset{(2)}{M_j^i},}...,\underset{}{\underset{(k)}{
M_j^i}}$ of the nonlinear connection $N$, (according to Theorem 2.5.2, ch.2)
determined by a $2$ -homogeneous spray S are $1,2,...,k$-homogeneous, respectively.

3$^{\circ }$. The same property have the coefficients $\underset{%
\left( 1\right) }{N_j^i},...,\underset{\left( k\right) }{N_j^i}$
of the nonlinear connection $N$.

4$^{\circ }.$ The local adapted basis to $N$, $\displaystyle\frac \delta
{\delta x^i},\displaystyle\frac \delta {\delta y^{\left( 1\right) i}},..., %
\displaystyle\frac \delta {\delta y^{\left( k\right) i}}$ has the following
degree of homogeneity 1,0, ...,1- k, respectively.

5$^{\circ }.$ If $X\in \mathcal{X}$ $(\widetilde{T^kM})$ is r-homogeneous
and $f\in \mathcal{F}(\widetilde{T^kM})$ is $s$-homogeneous, then $Xf^{}$is $%
s+r-1$-homogeneous.

Evidently the $1$-forms $dx^i,dy^{\left( 1\right)
i},...,dy^{\left( k\right) i}$ are homogeneous of degree \newline $0,1,...,k,$
respectively. The same degree of homogeneity have \newline  $\delta x^i,\delta
y^{\left( 1\right) i},...,\delta y^{\left( k\right) i}$ (the dual basis
adapted to $N$).

A $q$-form $\omega \in \Lambda ^q(\widetilde{T^kM})$ is $s$-homogeneous, if
and only if
\begin{equation}
\mathcal{L}_{\stackrel{k}{\Gamma }}\omega =s\omega.  \tag{3.1.5}
\end{equation}

Evidently, if $\omega $ is s - homogeneous and $\omega ^{\prime }$ is s'-
homogeneous then $\omega \wedge \omega ^{\prime }$ is s+s'- homogeneous.

For any $s$-homogeneous $\omega \in \Lambda ^q(\widetilde{T^kM})$ and any $r$-homogeneous
$% \underset{1}{X},...,\underset{q}{X}\in \underset{}{\mathcal{X}}(
\widetilde{T^kM})$ the function $\omega
(\underset{1}{X},.., \underset{q}{X})$ is $r+(s-1)q$ homogeneous.

Now we can formulate:

\begin{defi}
A Finsler space of order $k$, $k$ $\geq 1$, is a pair $F^{\left( k\right)
n}=(M,F)$ determined by a real differentiable manifold $M$ of dimension $n$
and a function $F:T^kM\rightarrow R$ having the following properties:

1$^{\circ }$. $F$ is differentiable on $\widetilde{T^kM}$ and continuous on
the null section $0:M\rightarrow T^kM$.

2$^{\circ }$. $F$ is positive.

3$^{\circ }$. $F$ is $k$-homogeneous on the fibres of the bundle $T^kM$,

4$^{\circ }$. The Hessian of $F^2$ with the entries
\begin{equation}
g_{ij}=\displaystyle\frac 12\displaystyle\frac{\partial ^2F^2}{\partial
y^{\left( k\right) i}\partial y^{\left( k\right) j}}  \tag{3.1.6}
\end{equation}
is positively defined on $\widetilde{T^kM}.$
\end{defi}

The $k$-homogeneous means \textit{positively }$k$\textit{- homogeneous}
since \\
$F(x,ay^{\left( 1\right) },...,a^ky^{\left( k\right) })=a^kF(x,y^{\left(
1\right) },...,y^{\left( k\right) })$ holds for any $a>0.$

It is not difficult to see that this definition has a geometrical meaning, $%
g_{ij}$ from (3.1.6) being a d-tensor field on $\widetilde{T^kM}.$ Such that the
function $F$ is called the fundamental or metric function and the $d$-tensor
field $g_{ij}$ is called fundamental or metric tensor of the Finsler space
of order $k$, $F^{\left( k\right) n}.$ Of course, the axiom 4$^{\circ }$
implies:
\begin{equation}
rank\left\| g_{ij}\right\| =n  \tag{3.1.7}
\end{equation}

Clearly, in the case $k=1,F^{(1)n}=(M,F)$ is a Finsler
space [115].

We can see, without difficulties that the following theorem holds:

\begin{teo}
The pair $L^{\left( k\right) n}=(M,F^2(x,y^{\left( 1\right) },...,y^{\left(
k\right) }))$ is a Lagrange space of order $k.$ Conversely, if $L^{\left(
k\right) n}=(M,L(x,y^{\left( 1\right) },...,y^{\left( k\right) })$ is a
Lagrange space of order $k$, having the fundamental function $L$ positively,
$2k$- homogeneous and the fundamental tensor $g_{ij}$ positively defined,
then the pair $F^{\left( k\right) n}=(M,\sqrt{L})$ is a Finsler space of
order $k.\,$
\end{teo}

Consequently, the class of spaces $F^{\left( k\right) n}$ is a subclass of
spaces $L^{\left( k\right) n}.$

Taking into account the $k$-homogeneity of the fundamental function $F$ and $%
2k$-homogeneity of $F^2$ we get:

1$^{\circ }$. $p_i$ given by
\begin{equation}
p_i=\displaystyle\frac 12\displaystyle\frac{\partial F^2}{\partial y^{\left(
k\right) i}}  \tag{3.1.7'}
\end{equation}
is a $k$-homogeneous d-covector fields.

2$^{\circ }$. $p=p_idx^i$ is a $k$-homogeneous 1-form.

3$^{\circ }$. The fundamental tensor $\ g_{ij}$ is $0$-homogeneous. Its contravariant \qquad $g^{ij}$ is $0$-homogeneous ($g_{ij}g^{jh}=\delta _i^h$), too.

These homogeneities imply:
\begin{equation}
\mathcal{L}_{\stackrel{k}{\Gamma }}F^2=2kF^2,\mathcal{L}_{\stackrel{k}{
\Gamma }}\displaystyle\frac{\partial F^2}{\partial y^{\left( k\right) }}=k %
\displaystyle\frac{\partial F^2}{\partial y^k},(or\text{ }\mathcal{L}_{%
\stackrel{k}{\Gamma }}p_i=kp_i),  \tag{3.1.8}
\end{equation}
\begin{equation}
\mathcal{L}_{\stackrel{k}{\Gamma
}}g_{ij}=0,\mathcal{L}_{\stackrel{k}{\Gamma
}}\underset{(k)}{C_{ijh}}=-k \underset{(k)}{C_{ijh}},  \tag{3.1.8'}
\end{equation}
where
\begin{equation}
\underset{(k)}{C_{ijh}}=\displaystyle\frac
12\displaystyle\frac{\partial g_{ij}}{\partial y^{\left( k\right)
h}}=\displaystyle\frac 14\displaystyle \frac{\partial
^3F^2}{\partial y^{\left( k\right) i}\partial y^{\left( k\right)
j}\partial y^{\left( k\right) h}}.  \tag{3.1.8''}
\end{equation}

Of course, the equation  $\mathcal{L}_{\stackrel{k}{\Gamma }}g_{ij}=0$ can be written as follows:
\begin{equation}
\mathcal{L}_{\stackrel{k}{\Gamma }}g_{rs}=(y^{\left( 1\right) i}%
\displaystyle \frac \partial {\partial y^{\left( 1\right) i}}+...+y^{\left(
k\right) i} \displaystyle\frac \partial {\partial y^{\left( k\right)
i}})g_{rs}=0 . \tag{3.1.8'''}
\end{equation}

The integral of action of $F$ along a parametrized curve $c:$
\begin{equation}
I(c)=\int_0^1F(x,\displaystyle\frac{dx}{dt},...,\displaystyle\frac 1{k!} %
\displaystyle\frac{d^kx}{dt^k})dt  \tag{3.1.9}
\end{equation}
can be considered for determining \textit{the length }of $\widetilde{c}
:[0,1]\rightarrow \widetilde{T^kM}$ in the given parametrization.

Theorem 2.1.1 shows that the necessary conditions for the integral of action $%
I(c),$ (3.1.9), to be independent on the parametrization of a curve $c$ are given
by the Zermelo conditions
\[
\mathcal{L}_{\stackrel{1}{\Gamma }}F=...=\mathcal{L}_{\stackrel{k-1}{\Gamma }
}F=0,\ \mathcal{L}_{\stackrel{k}{\Gamma }}F=F.
\]

But $\mathcal{L}_{\stackrel{k}{\Gamma }}F=F$ and $\mathcal{L}_{\stackrel{k}{
\Gamma }}F=kF,$ for $k>1$ are contradictory.

Consequently, in a Finsler space of order $k,$ $k>1,$ the integral
of action (3.1.9) essentially depend on the parametrization of a curve $c.$

\vspace*{3mm}

\textbf{Example 3.5.1}\ Let $F^n=(M,F(x,y^{(1)}))$ \textit{\ be a Finsler
space having } \newline $g_{ij}(x,y^1)$ as the fundamental tensor and $%
\underset{(1)}{M_j^i}(x,y^1)$ as coefficients of the
Cartan nonlinear connection. Then, theorem 1.5.4
gives us the dual coefficients on $T^kM:$
\[
\underset{(1)}{M_j^i}(x,y^{\left( 1\right)
}),\underset{(2)}{M_j^i} (x,y^{\left( 1\right) },y^{\left(
2\right) }),...,\underset{(k)}{M_j^i} (x,y^{\left( 1\right)
},...,y^{\left( k\right) })
\]
of a nonlinear connection, which depends only on the fundamental
function $F(x,y^{\left( 1\right) })$ of the Finsler space
$F^n.$

It is not hard to see that these coefficients are homogeneous of degree $1,2,...,k$
respectively. This property implies that the $d$-Liouville vector field $%
z^{\left( k\right) i}:$
\begin{equation}
kz^{\left( k\right) i}=ky^{\left( k\right)
i}+(k-1)\underset{(1)}{M_j^i} y^{\left( k-1\right)
j}+...+\underset{(k-1)}{M_j^i}y^{\left( 1\right) j} \tag{3.1.10}
\end{equation}
is linear in the variables $y^{\left( k\right) i}$ and it is $k$-homogeneous.

Consider the function
\begin{equation}
F(x,y^{\left( 1\right) },...,y^{\left( k\right) })=\{g_{ij}(x,y^{\left(
1\right) })z^{\left( k\right) i}z^{\left( k\right) j}\}^{1/2},  \tag{3.1.11}
\end{equation}
$g_{ij}(x,y^{\left( 1\right) })$ being a d-tensor positively defined. It
follows that $F$ from (3.1.11) is a positive differentiable function on $%
\widetilde{T^kM}$ and continuous on the null section of $\pi ^k$. It is $k$
-homogeneous and has $g_{ij}(x,y^{(1)})$ as the fundamental tensor.

Consequently, the pair $F^{(k)n}=(M,F(x,y^{\left( 1\right) },...,y^{\left(
k\right) })$ for $F$ from (3.1.11) is a Finsler space of order $k$.

Concluding, we have:
\begin{teo}
If the base manifold is paracompact then there exist a Finsler space of
order $k$, $F^{\left( k\right) n}$.
\end{teo}

The spaces $F^{(k)n}$ constructed in example (3.1.1) is called the
Prolongation of order $k$ of the Finsler space $F^n.$ It is denoted by Prol$%
^kF^n.$

In order to determine the geodesics of the space $F^{(k)n}$ we take the integral
of action of the regular Lagrangian $F^2.$ The variational problem leads to
the Euler - Lagrange equations.
\begin{equation}
\stackrel{\circ }{E_i}(F^2)=\displaystyle\frac{\partial F^2}{\partial x^i}- %
\displaystyle\frac d{dt}\displaystyle\frac{\partial F^2}{\partial y^{\left(
1\right) i}}+...+(-1)^k\displaystyle\frac 1{k!}\displaystyle\frac{d^k}{dt^k} %
\displaystyle\frac{\partial F^2}{\partial y^{\left( k\right) i}}=0.
\tag{3.1.12}
\end{equation}

The integral curves of the previous equations are called the geodesics of the
space $F^{(k)n}.$ Applying the theory from the section 2, ch. 2 we can
determine the infinitesimal symmetries of the spaces $F^{\left( k\right) n}.$

The energies of order $k,k-1,...,1$ of the Finsler space of order $k$, $%
F^{\left( k\right) n}=(M,F(x,y^{\left( 1\right) },...,y^{\left( k\right) }))$
are given by the formulae (2.3.1), for $L=F^2.$

In particular, Theorem 2.3.2 can be applied in order to obtain

\begin{teo}
The energy of order $k$, $\mathcal{E}^k(F^2)$ of the Finsler space $F^{(k)n}$
is conserved along every geodesic of this space.
\end{teo}

Let us consider
\begin{equation}
\stackrel{\vee}{p}_{\left( k\right) i}=\displaystyle\frac{\partial F^2}{
\partial y^{\left( k\right) i}},
\stackrel{\vee}{p}_{\left( k-1\right) i}= %
\displaystyle\frac{\partial F^2}{\partial y^{\left( k-1\right) i}},...,
\stackrel{\vee}{p}_{\left( 1\right) i}=\displaystyle\frac{\partial F^2}{
\partial y^{\left( 1\right) i}},
\stackrel{\vee}{p}_{\left( 0\right) i}= %
\displaystyle\frac{\partial F^2}{\partial x^i} . \tag{3.1.13}
\end{equation}

Then we have

\begin{teo}
1$^{\circ }$. The Cartan differential 1-forms are the followings
\begin{equation}
\begin{array}{l}
d_0F^2=\stackrel{\vee}{p}_{\left( k\right) i}dx^i, \\
\\
d_1F^2=\stackrel{\vee}{p}_{\left( k-1\right) i}dx^i+\stackrel{\vee}{p}%
_{\left( k\right) i}dy^{\left( k\right) i} ,\\
..........................................., \\
d_kF^2=\stackrel{\vee}{p}_{\left( 0\right) i}dx^i+\stackrel{\vee}{p}%
_{\left( 1\right) i}dy^{\left( 1\right) i}+...+\stackrel{\vee}{p}_{\left(
k\right) i}dy^{\left( k\right) i}.
\end{array}
\tag{3.1.14}
\end{equation}

2$^{\circ }$. And the Poincare 2-forms are given by

$dd_0F^2=d\stackrel{\vee}{p}_{\left( k\right) i}\wedge dx^i,$

.........................................,

$dd_{(k-1)}F^2=d\stackrel{\vee}{p}_{\left( 1\right) i}\wedge dx^i+...+d%
\stackrel{\vee}{p}_{\left( k\right) i}dy^{\left( k\right) i}.$\vspace{3mm}\\
Here we have $dd_kF^2=d^2F^2=0.$\vspace{3mm}

3$^{\circ }$. The 1-forms $d_0F^2,d_1F^2,...,d_kF^2$ are $k,k+1,...,2k$
homogeneous, respectively.
\end{teo}

\section{Cartan Nonlinear Connection in $F^{\left( k\right) n}$}

The considerations made in the previous chapter allow us to introduce in a
Finsler space of order $k$, $F^{\left( k\right) n}=(M,F)$ the main
geometrical object fields as: canonical $k$-spray, Cartan nonlinear
connection, canonical $N$-linear connection etc. Canonical mean here that all
these object fields depend only on the fundamental function $F$.

Taking into account the operators $\stackrel{0}{E_i},\stackrel{1}{E_i},...,
\stackrel{k}{E_i}$ given by (2.2.8) we construct the system of $d$
-covector fields
\begin{equation}
\stackrel{0}{E_i}(F^2),\stackrel{1}{E_i}(F^2),...,\stackrel{k}{E_i}(F^2)
\tag{3.2.1}
\end{equation}

All equations
\[
\stackrel{0}{E_i}(F^2)=0,\stackrel{1}{E_i}(F^2)=0,...,\stackrel{k-1}{E_i}
(F^2)=0
\]
have geometrical meanings. The equation
\[
\stackrel{k-1}{E_i}(F^2)=0
\]
is important for us. It will be called the \textit{Craig-Synge equation}.
Using (2.2.8), this equation is expressed as follows
\begin{equation}
\displaystyle\frac{\partial F^2}{\partial y^{(k-1)i}}-\displaystyle\frac
d{dt}\displaystyle\frac{\partial F^2}{\partial y^{\left( k\right) i}}=0.
\tag{3.2.2}
\end{equation}

But $\displaystyle\frac d{dt}\displaystyle\frac{\partial F^2}{\partial
y^{\left( k\right) i}}=\Gamma \displaystyle\frac{\partial F^2}{\partial
y^{\left( k\right) i}}+\displaystyle\frac 2{k!}g_{ij}\displaystyle\frac{
d^{k+1}x^i}{dt^{k+1}},$ where $\Gamma $ is the operator (1.2.3).\vspace{3mm}

Consequently, the Craig-Synge equations (3.2.2) is equivalent to the following
equations
\[
\stackrel{k-1}{g^{ij}E_j(F^2)}=0, y^{(1)i}=\frac{dx^i}{dt}, ..., y^{(k)i}=\frac{1}{k!}\frac{d^kx^i}{dt^k}
\]
or
\begin{equation}
\displaystyle\frac{d^{k+1}x^i}{dt^{k+1}}+(k+1)G^i(x,\displaystyle\frac{dx}{%
dt },...,\displaystyle\frac 1{k!}\displaystyle\frac{d^kx}{dt^k})=0,  \tag{3.2.3}
\end{equation}
where
\begin{equation}
(k+1)G^i(x,y^{\left( 1\right) },...,y^{(k)})=\displaystyle\frac
12g^{ij}\{\Gamma (\displaystyle\frac{\partial F^2}{\partial y^{\left(
k\right) j}})-\displaystyle\frac{\partial F^2}{\partial y^{\left( k-1\right)
j}}\}.  \tag{3.2.4}
\end{equation}

Applying the Theorem 2.5.1 one obtains:

\begin{teo}
The Craig-Synge equations (3.2.3) determines a canonical $k$-spray $S:$
\begin{equation}
S=y^{\left( 1\right) i}\displaystyle\frac \partial {\partial x^i}+2y^{\left(
2\right) i}\displaystyle\frac \partial {\partial y^{\left( 1\right)
i}}...+ky^{\left( k\right) i}\displaystyle\frac \partial {\partial y^{\left(
k-1\right) i}}-(k+1)G^i(x,y^{\left( 1\right) },...,y^{\left( k\right) })%
\displaystyle\frac \partial {\partial y^{\left( k\right) i}}  \tag{3.2.5}
\end{equation}
with the coefficients $G^i$ from (3.2.4).
\end{teo}

Note that $S$ depend only on the fundamental function $F$ of the space $%
F^{\left( k\right) n}.$

It is a $k$-spray, since it is a 2-homogeneous vector field.
The paths of $S$ are given by the differential equations (3.2.3).

By means of the Theorems 1.5.1 and 3.2.1 the dual coefficients of the nonlinear
connection $N$ determined by the canonical $k$-spray $S$ are given by:

\begin{teo}
In a Finsler space of order $k$, $F^{\left( k\right) n}=(M,F)$ there exist
nonlinear connections, depending only on the fundamental function $F$. One
of these, denoted by $N,$ has the dual coefficients:
\begin{equation}
\begin{array}{l}
\underset{(1)}{M_j^i}=\displaystyle\frac
1{2(k+1)}\displaystyle\frac
\partial {\partial y^{\left( k\right) j}}\{g^{im}[\Gamma \displaystyle\frac{%
\partial F^2}{\partial y^{\left( k\right) m}}-\displaystyle\frac{\partial F^2%
}{\partial y^{\left( k-1\right) m}}]\}, \\
\\
\underset{(2)}{M_j^i}=\displaystyle\frac 12\{S\underset{(1)}{M_i^j}+%

\underset{(1)}{M_m^i}\underset{(1)}{M_j^m}]\}, \\
........................................... \\
\underset{(k)}{M_j^i}=\displaystyle\frac 1k\{S\underset{(k-1)}{M_j^i}+%
\underset{(1)}{M_m^i}\underset{(k-1)}{M_j^m}]\},
\end{array}
\tag{3.2.6}
\end{equation}
where $S$ is the canonical $k$-spray of the space $F^{\left( k\right) n}.$
\end{teo}
$N$ is called the Cartan nonlinear connection of the space $F^{\left(
k\right) n}.$

In the case $k=1$, $N$ reduces to the classical Cartan nonlinear
connection of the Finsler space $F^n=(M,F(x,y^{\left( 1\right) })).$

Some properties of $N.$

1$^{\circ }.$ The Cartan nonlinear connection $N$ is globally defined on $%
\widetilde{T^kM}$ (if \newline $F(x,y^{\left( 1\right) },...,y^{\left( k\right) })$
has this property)

2$^{\circ }$. The dual coefficients (3.2.6) of $N$ are homogeneous
of degree $1,2,...,k$ i.e:
\begin{equation}
\mathcal{L}_{\stackrel{k}{\Gamma }}\underset{(\alpha
)}{M_j^i}=\alpha \underset{(\alpha )}{M_j^i},\alpha =1,...,k.
\tag{3.2.7}
\end{equation}

The coefficients
$\underset{(1)}{N_j^i},...,\underset{(k)}{N_j^i}$ of $N,$ (ch.1)
are expressed by
\begin{equation}
\begin{array}{l}
\underset{(1)}{N_j^i}=\underset{(1)}{M_j^i}, \\
\\
\underset{(2)}{N_j^i}=\underset{(2)}{M_j^i}-\underset{(1)}{N_m^i}
\underset{(1)}{M_j^m}, \\
........................................... \\
\underset{(k)}{N_j^i}=\underset{(k)}{M_j^i-}\underset{(k-1)}{N_m^i}
\underset{(1)}{M_j^m}-...-\underset{(1)}{N_m^i}\underset{(k-1)}{M_j^m}.
\end{array}
\tag{3.2.8}
\end{equation}

These coefficients are homogeneous functions of degree $1,...,k$
respectively, i.e
\begin{equation}
\mathcal{L}_{\stackrel{k}{\Gamma }}\underset{(\alpha
)}{N_j^i}=\alpha \underset{(\alpha)}{N_j^i},\ \ \ (\alpha =1,...,k).  \tag{3.2.7'}
\end{equation}

The Cartan nonlinear connection $N$ gives rise to a distribution \newline
$N_u\subset
T_u(\widetilde{T^kM})$ supplementary to the vertical distribution $%
V_{u}\subset T_u(\widetilde{T^kM}),$ $\forall u\in \widetilde{T^kM}$ with
the property:
\[
T_u\widetilde{(T^kM})=N_u\oplus V_u\text{ },\text{ }\forall u\in \widetilde{
T^kM}.
\]

If we consider the distributions
\[
N_0=N,N_1=J(N_0),...,N_{k-1}=J(N_{k-2}),V_k=J(N_{k-1}),
\]
then according to the general theory we obtain the direct decomposition of the
vector spaces:
\begin{equation}
T_u(\widetilde{T^kM})=\ N_{0,u}\oplus N_{1,u}\oplus ...\oplus
N_{k-1,u}\oplus V_{k,u},\forall u\in \widetilde{T^kM} . \tag{3.2.9}
\end{equation}

The local adapted basis to the direct decomposition (3.2.9) is
\begin{equation}
(\displaystyle\frac \delta {\delta x^i},\displaystyle\frac \delta {\delta
y^{\left( 1\right) i}},...,\displaystyle\frac \delta {\delta y^{\left(
k\right) i}}),  \tag{3.2.10}
\end{equation}
where
\begin{equation}
\begin{array}{l}
\displaystyle\frac \delta {\delta x^i}=\displaystyle\frac \partial
{\partial x^i}-\underset{(1)}{N_i^j}\displaystyle\frac \partial
{\partial y^{\left( 1\right)
j}}-...-\underset{(k)}{N_i^j}\displaystyle\frac \partial
{\partial y^{\left( k\right) j}}, \\
\\
\displaystyle\frac \delta {\delta y^{\left( 1\right) i}}=\displaystyle\frac
\partial {\partial y^{\left( 1\right) i}}-\underset{(1)}{N_i^j} %
\displaystyle\frac \partial {\partial y^{\left( 2\right)
j}}-...-\underset{
(k-1)}{N_i^j}\displaystyle\frac \partial {\partial y^{\left( k\right) j}}, \\
........................................................... \\
\displaystyle\frac \delta {\delta y^{\left( k\right) i}}=\displaystyle\frac
\partial {\partial y^{\left( k\right) i}},
\end{array}
\tag{3.2.11}
\end{equation}
$\underset{(1)}{N_i^j},...,\underset{(k)}{N_i^j}$ being the coefficients (3.2.8) of the Cartan nonlinear connection.

Of course we have
\[
\displaystyle\frac \delta {\delta y^{\left( 1\right)
i}}=J(\displaystyle \frac \delta {\delta
x^i}),...,\displaystyle\frac \delta {\delta y^{\left( k\right)
i}}=J(\displaystyle\frac \delta {\delta y^{\left( k-1\right)
i}}),0=J(\displaystyle\frac \delta {\delta y^{\left( k\right)
i}}).
\]

Taking into account section 1.4 ch.1, the dual (adapted) cobasis, of the
basis (3.2.10) is:
\begin{equation}
(\delta x^i,\delta y^{\left( 1\right) i},...,\delta y^{\left( k\right) i})
\tag{3.2.12}
\end{equation}
where
\begin{equation}
\begin{array}{l}
\delta x^i=dx^i, \vspace{3mm}\\
\delta y^{\left( 1\right) i}=dy^{\left( 1\right)
i}+\underset{(1)}{M_j^i}dx^j, \\
....................................... \\
\delta y^{\left( k\right) i}=dy^{\left( k\right)
i}+\underset{(1)}{M_j^i} dy^{\left( k-1\right)
i}+...+\underset{(k)}{M_j^i}dx^j,
\end{array}
\tag{3.2.12a}
\end{equation}
$\underset{\left( 1\right) }{M_j^i},...,\underset{\left( k\right) }{%
M_j^i },$ being the dual coefficients (3.2.6) of the Cartan nonlinear connection $N$.

It is not difficult to see that the following identities hold:
\begin{equation}
\left\{
\begin{array}{l}
J^{*}(\delta y^{\left( k\right) i})=\delta y^{\left( k-1\right) i},\
J^{*}(\delta y^{\left( k-1\right) i})=\delta y^{\left( k-2\right) i},..., \\
\\
J^{*}(\delta y^{\left( 1\right) i})=\delta x^i,\ J^{*}(dx^i)=0.
\end{array}
\right.  \tag{3.2.13}
\end{equation}
$J^{*}$ being the adjoint of the $k$-structure $J$.

Now we can determine the differential operators $d_k,d_{k-1},...,d_0$
defined in (1.2.11), using the expression of the operator of
differentiation $d_k=d$ in the adapted basis:
\begin{equation}
d_k=\displaystyle\frac \delta {\delta x^i}\delta x^i+\displaystyle\frac
\delta {\delta y^{\left( 1\right) i}}\delta y^{\left( 1\right) i}+...+ %
\displaystyle\frac \delta {\delta y^{\left( k\right) i}}\delta y^{\left(
k\right) i} . \tag{3.2.14}
\end{equation}
Applying (3.2.13), we get $d_k$ in (3.2.14), and
\[
d_{k-1}=J^{*}d_k,d_{k-2}=J^{*}d_{k-1},...,d_0=J^{*}d_1.
\]

One obtains:
\begin{equation}
\begin{array}{l}
d_{k-1}=\displaystyle\frac \delta {\delta y^{\left( 1\right) i}}\delta x^i+ %
\displaystyle\frac \delta {\delta y^{\left( 2\right) i}}\delta y^{\left(
1\right) i}+...+\displaystyle\frac \delta {\delta y^{\left( k\right)
i}}\delta y^{\left( k-1\right) i}, \\
\\
d_{k-2}=\displaystyle\frac \delta {\delta y^{\left( 2\right) i}}\delta x^i+ %
\displaystyle\frac \delta {\delta y^{\left( 3\right) i}}\delta y^{\left(
1\right) i}+...+\displaystyle\frac \delta {\delta y^{\left( k\right)
i}}\delta y^{\left( k-2\right) i}, \\
....................................................... \\
d_1=\displaystyle\frac \delta {\delta y^{\left( k-1\right) i}}\delta x^i+ %
\displaystyle\frac \delta {\delta y^{\left( k\right) i}}\delta y^{\left(
1\right) i}, \\
\\
d_0=\displaystyle\frac \delta {\delta y^{(k)i}}\delta x^i.
\end{array}
\tag{3.2.15}
\end{equation}
Consequently, we get:

\begin{teo}
With respect to the direct decomposition (3.2.9), in adapted basis (3.2.11),
(3.2.12), the Cartan 1-forms $d_0F^2,d_1F^2,...,d_kF^2$ of a Finsler
space of order $k$, $F^{\left( k\right) n}=(M,F)$ can be expressed as
follows:
\begin{equation}
\begin{array}{l}
d_0F^2=(d_0F^2)^H, \\
\\
d_1F^2=(d_1F^2)^H+(d_1F^2)^{V_1}, \\
....................................................... \\
d_kF^2=(d_kF^2)^H+(d_kF^2)^{V_1}+...+(d_kF^2)^{V_k}.
\end{array}
\tag{3.2.16}
\end{equation}
Equivalently,
\begin{equation}
\begin{array}{l}
d_0F^2=\displaystyle\frac{\delta F^2}{\delta y^{\left( k\right) i}}\delta x^i,
\\
\\
d_1F^2=\displaystyle\frac{\delta F^2}{\delta y^{\left( k-1\right) i}}\delta
x^i+\displaystyle\frac{\delta F^2}{\delta y^{\left( k\right) i}}\delta
y^{\left( 1\right) i}, \\
....................................................... \\
d_kF^2=\displaystyle\frac{\delta F^2}{\delta x^i}\delta x^i+\displaystyle
\frac{\delta F^2}{\delta y^{\left( 1\right) i}}\delta y^{\left( 1\right)
i}+...+\displaystyle\frac{\delta F^2}{\delta y^{\left( k\right) i}}\delta
y^{\left( k\right) i}.
\end{array}
\tag{3.2.17}
\end{equation}
\end{teo}

In the previous expressions every term is an 1-form field on $\widetilde{%
T^kM }.$ So we have the following main 1-form fields

\begin{equation}
\begin{array}{l}
\theta _0=(d_0F^2)^H=\displaystyle\frac{\delta F^2}{\delta y^{\left(
k\right) i}}\delta x^i=\stackrel{\vee}{p}_{\left( k\right) i}\delta x^i, \\
\\
\theta _1=(d_1F^2)^{V_1}=\displaystyle\frac{\delta F^2}{\delta y^{\left(
k\right) i}}\delta y^{\left( k\right) i}=\stackrel{\vee}{p}_{\left(
k\right) i}\delta y^{\left( 1\right) i}, \\
....................................................... \\
\theta _{k-1}=(d_{k-1}F^2)^{V_{k-1}}=\displaystyle\frac{\delta F^2}{\delta
y^{\left( k\right) i}}\delta y^{\left( k-1\right) i}=\stackrel{\vee}{p}
_{\left( k\right) i}\delta y^{\left( k-1\right) i}, \\
\\
\theta _k=(d_kF^2)^{V_k}=\displaystyle\frac{\delta F^2}{\delta y^{\left(
k\right) i}}\delta y^{\left( k\right) i}=\stackrel{\vee}{p}_{\left(
k\right) i}\delta y^{\left( k\right) i}.
\end{array}
\tag{3.2.18}
\end{equation}

\begin{teo}
1$^{\circ }$. The 1-form fields $\theta _0,...,\theta _k$ depend only on the
fundamental function $F$ of the Finsler space $F^{\left( k\right) n}.$

2$^{\circ }$. The exterior differentials of $\theta _0,...,\theta _k,$
\begin{equation}
\begin{array}{l}
d\theta _0=d\stackrel{\vee}{p}_{\left( k\right) i}\wedge \delta x^i, \\
\\
d\theta _1=d\stackrel{\vee}{p}_{\left( k\right) i}\wedge \delta y^{\left(
1\right) i}+\stackrel{\vee}{p}_{\left( k\right) i}\wedge d\delta y^{\left(
1\right) i}, \\
....................................................... \\
d\theta _k=d\stackrel{\vee}{p}_{\left( k\right) i}\wedge \delta y^{\left(
k\right) i}+\stackrel{\vee}{p}_{\left( k\right) i}\wedge d\delta y^{\left(
k\right) i}
\end{array}
\tag{3.2.19}
\end{equation}
have the same property of homogeneity.
\end{teo}

The second terms of $d\theta _0,...,d\theta _k$, the exterior differentials
of 1-forms $\delta y^{\left( 1\right) i},...,\delta y^{\left( k\right) i}$ are
calculated by means of formulas:
\begin{equation}
\begin{array}{l}
d\delta y^{\left( \alpha \right) i}=\displaystyle\frac
12\underset{\left( 0\alpha \right) }{R_{jm}^i}dx^m\wedge
dx^j+\stackrel{k}{\underset{\gamma =1 }{\sum }}\underset{(\gamma
\alpha )}{B_{jm}^i}dy^{(\gamma )m}\wedge
dx^j+ \\
\qquad \qquad +\stackrel{k}{\underset{\beta ,\gamma =1}{\sum
}}\underset{ (\alpha \gamma )}{\stackrel{(\beta
)}{C_{jm}^i}}dy^{(\gamma )m}\wedge dy^{(\alpha )j},
\end{array}
\tag{3.2.20}
\end{equation}

where $\underset{(\alpha \alpha )}{\stackrel{(\alpha
)}{C_{jm}^i}}=0$, and the coefficients from the right hand side can be
calculated using the following Lie brackets:
\begin{equation}
\begin{array}{c}
\lbrack \displaystyle\frac \delta {\delta x^j},\displaystyle\frac
\delta {\delta x^k}]=\underset{\left( 01\right)
}{R_{jh}^i}\displaystyle\frac \delta {\delta y^{\left( 1\right)
i}}+...+\underset{\left( 0k\right) }{
R_{jh}^i}\displaystyle\frac \delta {\delta y^{\left( k\right) i}}, \\
\\
\lbrack \displaystyle\frac \delta {\delta x^j},\displaystyle\frac
\delta {\delta y^{\left( \alpha \right) h}}]=\underset{\left(
\alpha 1\right) }{ B_{jh}^i}\displaystyle\frac \delta {\delta
y^{\left( 1\right) i}}+...+ \underset{\left( \alpha k\right)
}{B_{jh}^i}\displaystyle\frac \delta
{\delta y^{\left( k\right) i}}, \\
\\
\lbrack \displaystyle\frac \delta {\delta y^{\left( \alpha \right) j}}, %
\displaystyle\frac \delta {\delta y^{\left( \beta \right)
h}}]=\underset{ (\alpha \beta
)}{\stackrel{(1)}{C_{jh}^i}}\displaystyle\frac \delta {\delta
y^{\left( 1\right) i}}+...+\underset{(\alpha \beta
)}{\stackrel{(k)}{
C_{jh}^i}}\displaystyle\frac \delta {\delta y^{\left( k\right) i}}, \\
\\
(\alpha ,\beta =1,...,k).
\end{array}
\tag{3.2.21}
\end{equation}

The whole previous theory can be applied to the following Lagrangians
associated to the Finsler space $F^{\left( k\right) n}$
\begin{equation}
F_1^2=g_{ij}z^{\left( 1\right) i}z^{\left( 1\right)
j},...,F_k^2=g_{ij}z^{\left( k\right) i}z^{\left( k\right) j}  \tag{3.2.22}
\end{equation}

especially in the cases of the particular Finsler spaces of order $k$, \textit{%
Prol }$^k\mathcal{R}^n$ or \textit{Prol}$^kF^n$ (see ch.2). The Lagrangians $%
F_1^2,...,F_k^2$ are positive functions and are $2,4,...,2k$ -homogeneous, respectively.

The autoparallel curves of the Cartan nonlinear connection $N$ are characterized by the system of differential equations
\[
\begin{array}{c}
\displaystyle\frac{\delta y^{\left( 1\right) i}}{dt}=...=\displaystyle\frac{

\delta y^{\left( k\right) i}}{dt}=0, \vspace{3mm} \\
y^{\left( 1\right) i}=\displaystyle\frac{dx^i}{dt},...,y^{\left( k\right)
i}= \frac{1}{k!}\displaystyle\frac{dx^{\left( k\right) i}}{dt^k}.
\end{array}
\]

\section{The Cartan Metrical $N$-Linear Connection}

Let $N$ be the Cartan nonlinear connection of the Finsler space of order $k$, $F^{\left( k\right) n}=(M,F)$ having the adapted basis (3.2.10) and its dual (3.2.12).

The lift of the fundamental tensor field $g_{ij}$ is given by (2.6.1),
\begin{equation}
\Bbb{G}=g_{ij}dx^i\otimes dx^j+g_{ij}\delta y^{\left( 1\right) i}\otimes
\delta y^{\left( 1\right) j}+...+g_{ij}\delta y^{\left( k\right) i}\otimes
\delta y^{\left( k\right) j}.  \tag{3.3.1}
\end{equation}

\begin{teo}
$\Bbb{G}$ from (3.3.1) is a Riemannian structure on $\widetilde{T^kM}$ which
depend only on the fundamental function $F$ of the space $F^{\left( k\right)
n}$. The terms of $\Bbb{G}$ are $0,2,...,2k$ homogeneous, respectively.
\end{teo}

Notice that $\Bbb{G}$ is not homogeneous. We can construct an homogeneous
one using the Lagrangians (3.2.22).

Namely
\begin{equation}
\stackrel{\vee}{\Bbb{G}}=g_{ij}dx^i\otimes dx^j+\displaystyle\frac
1{F_1^2}g_{ij}\delta y^{\left( 1\right) i}\otimes \delta y^{\left( 1\right)
j}+...+\displaystyle\frac 1{F_k^2}g_{ij}\delta y^{\left( k\right) i}\otimes
\delta y^{\left( k\right) j}.  \tag{3.3.1a}
\end{equation}
$\stackrel{\vee}{\Bbb{G}}$ is a Riemannian structure on the manifold $\widetilde{T^kM}$ determined only by the fundamental function $\Bbb{F}$ and it is $0$-homogeneous.

In the following we consider the Riemannian structure $\Bbb{G}$ from (3.3.1).

An $N$-linear connection $D$ is compatible with $\Bbb{G}$ if
\[
D_X\Bbb{G}=0,\ \forall X\in \mathcal{X}(\widetilde{T^kM}).
\]

Applying the Theorem 2.6.1, we have:

\begin{teo}
For a Finsler space of order $k$, $F^{\left( k\right) n}=(M,F)$, the
following properties hold:

1$^{\circ }$. There exists an unique $N$-linear connection $D$ on $\widetilde{%
T^kM}$ verifying the axioms:

A$_1$      $N$ is the Cartan nonlinear connection,

A$_2$      $g_{ij|h}=0$,

A$_3$      $g\stackrel{(\alpha )}{_{ij}|_h}=0$,

A$_4$      $F_{jk}^i=F_{kj}^i$,

A$_5$      $\underset{(\alpha )}{C_{jk}^i}=\underset{(\alpha )}{C_{kj}^i}$ $(\alpha = 1, ..., k)$.

2$^{\circ }$. The coefficients $C\Gamma (N)=(F_{jk}^i,\underset{(1)}{C_{jk}^i},...,\underset{(k)}{C_{jk}^i})$ of $D$ are given by the
generalized Christoffel symbols (2.6.3),
($F_{jk}^i=L_{jk}^i$).

3$^{\circ }$. $D$ depends only on the fundamental function $F$ of the space $F^{\left( k\right) n}$.
\end{teo}

The metrical $N$ linear -connection $D$ from the previous theorem will be called the
\textit{Cartan metrical }$N$\textit{-linear connection} of the space $%
F^{\left( k\right) n}$ and denoted by $C\Gamma (N).$

Of course, the torsion d-tensor fields and the curvature d-tensor fields of $%
C\Gamma (N)$ can be written without difficulties. Such that we have
\begin{equation}
\underset{(0)}{T_{jk}^i}=0,\quad \underset{(\alpha )}{S_{jk}^i}
=0,(\alpha =1,...,k).  \tag{3.3.2}
\end{equation}

Also we can calculate the deflection tensor of $C\Gamma (N):$
\[
\stackrel{(\alpha )}{D_j^i}=z_{\quad |j}^{(\alpha )i},\quad \stackrel{(\beta
\alpha )}{d_j^i}=z^{(\beta )i}\stackrel{(\alpha )}{|_j}.
\]

The coefficients $C\Gamma
(N)=(F_{jh}^i,\underset{(1)}{C_{jk}^i},...,
\underset{(k)}{C_{jk}^i})$ are $0,-1,...,-k$ -homogeneous.

The d-tensors of curvature of $C\Gamma (N)$ satisfy the identities:
\begin{equation}
\begin{array}{l}
g_{sj}R_{i\ hm}^s+g_{is}R_{j\ hm}^s=0, \\
\\
g_{sj}\underset{}{\underset{(\alpha )}{P}\text{}_{i\ hm}^s}+g_{is}
\underset{(\alpha )}{P}\text{}_{j\ hm}^s=0, \\
\\
g_{sj}\underset{}{\underset{(\alpha \beta )}{S}\text{}_{i\ hm}^s}+g_{is} \underset{(\alpha \beta )}{S}\text{}_{j\ hm}^s=0
\end{array}
\tag{3.3.3}
\end{equation}

Notice that the equations $\underset{(0\alpha )}{R}\text{}_{jh}^i=0,(\alpha =1,...,k)$ characterize the integrability of
the Cartan nonlinear connection.

The connection 1-forms $\omega _j^i$ of the Cartan metrical $N$-linear connection $C\Gamma (N)$ are given by:
\begin{equation}
\omega_{\ j}^i=F_{jh}^idx^k+\underset{(1)}{C_{jh}^i}\delta y^{\left(
1\right) h}+...+\underset{(k)}{C_{jh}^i}\delta y^{\left( k\right)h}.  \tag{3.3.4}
\end{equation}

\begin{teo}
The structure equations of the Cartan metrical $N$-linear connection $C\Gamma (N)$ of
the Finsler space $F^{\left( k\right) n\text{ }}$ are given by:
\begin{equation}
\begin{array}{l}
d(dx^i)-dx^m\wedge \omega_{\ m}^i=-\stackrel{(0)}{\Omega ^i}, \\
\\
d(\delta y^{\left( \alpha \right) i})-\delta y^{\left( \alpha \right)
m}\wedge \omega_{\ m}^i=-\stackrel{(\alpha )}{\Omega ^i},\ (\alpha =1,...,k), \\
\\
d\omega_{\ j}^i-\omega_{\ j}^m\wedge \omega_{\ m}^i=-\Omega_{\ j}^i,
\end{array}
\tag{3.3.5}
\end{equation}
where $\stackrel{(0)}{\Omega ^i},\stackrel{(\alpha )}{\Omega ^i},$ are the
2-forms of torsion:
\begin{equation}
\begin{array}{l}
\stackrel{(0)}{\Omega ^i}=\underset{(1)}{C_{jh}^i}dx^j\wedge
\delta y^{\left( 1\right)
h}+...+\underset{(k)}{C_{jh}^i}dx^j\wedge \delta
y^{\left( k\right) h} \\
\\
\stackrel{(\alpha )}{\Omega ^i}=\displaystyle\frac 12\underset{(0\alpha )}{%
R_{jh}^i}dx^j\wedge dx^h+\stackrel{k}{\underset{\gamma =1}{\sum }}%
\underset{(\gamma \alpha )}{B_{jm}^i}dx^j\wedge \delta y^{\left(
\gamma
\right) m}+ \\
\\
+\stackrel{k}{\underset{\gamma =1}{\sum }}\underset{(\alpha \gamma )}{%
\stackrel{(\beta )}{C_{jh}^i}}\delta y^{\left( \beta \right) j}\wedge \delta
y^{\left( \gamma \right) h}-(F_{jh}^idx^j+\stackrel{k}{\underset{\gamma =1%
}{\sum }}\underset{(\gamma )}{C_{jh}^i}\delta y^{\left( \gamma
\right) j})\delta y^{\left( \alpha \right) h}
\end{array}
\tag{3.3.6}
\end{equation}
and $\Omega_{\ j}^i$ are the $2$-forms of curvature:

\begin{equation}
\Omega_{\ j}^i=\displaystyle\frac 12R_{j\ pq}^idx^p\wedge dx^q+\stackrel{k}{%
\underset{\gamma =1}{\sum }}\underset{(\gamma
)}{P}{}_{j\ pq}^i dx^p\wedge
\delta y^{\left( \gamma \right) q}+\stackrel{k}{\underset{\beta ,\gamma =1%
}{\sum }}\underset{(\beta \gamma )}{S}{}_{j\ pq}^i \delta y^{\left(
\beta \right) p}\wedge \delta y^{\left( \gamma \right) q}.
\tag{3.3.7}
\end{equation}
\end{teo}

Now, we can obtain the Bianchi identities of the Cartan metrical $N$-linear
connection $C\Gamma (N)$ if we apply the exterior differential to the system

of equations (3.3.5) and calculate $d\stackrel{(0)}{\Omega ^i}, d\stackrel{
(\alpha )}{\Omega ^i}$ and $d\Omega_{\ j}^i$ from (3.3.6) and (3.3.7), modulo the system (3.3.5).

Finally, consider the tensor field $\Bbb{F}$ determined by the Cartan
nonlinear connection $N$.
\begin{equation}
\Bbb{F}\mathbf{=-}\displaystyle\frac \delta {\delta y^{\left( k\right)
i}}\otimes dx^i+\displaystyle\frac \delta {\delta x^i}\otimes \delta y^{\left(
k\right) i} . \tag{3.3.8}
\end{equation}
It is not difficult to prove that $C\Gamma (N)$ has the property $D_X\Bbb{F}=0.$
Indeed,

$D_{\displaystyle\frac \delta {\delta x^h}}\Bbb{F}=-(D_{\displaystyle\frac \delta
{\delta x^h}}\displaystyle\frac \delta {\delta y^{\left( k\right)
i}})\otimes dx^i-\displaystyle\frac \delta {\delta y^{\left( k\right)
i}}\otimes D_{\displaystyle\frac \delta {\delta x^h}}dx^i+$ \vspace{3mm}\\

$(D_{\displaystyle\frac \delta {\delta x^h}}\displaystyle\frac \delta
{\delta x^i})\otimes \delta y^{\left( k\right) i}+\displaystyle\frac \delta
{\delta x^i}\otimes D_{\displaystyle\frac \delta {\delta x^h}}\delta
y^{\left( k\right) i}=$ \vspace{3mm}\\

$-F_{ih}^m\displaystyle\frac \delta {\delta y^{\left( k\right) m}}\otimes
dx^i+F_{mh}^i\displaystyle\frac \delta {\delta y^{\left( k\right) i}}\wedge
dx^m+F_{ih}^m\displaystyle\frac \delta {\delta x^m}\otimes \delta y^{\left(
k\right) i}-F_{mh}^i\displaystyle\frac \delta {\delta x^i}\otimes \delta
y^{\left( k\right) m}=0.$ \vspace{3mm}\\

We proceed analogously in the case $D_{\displaystyle\frac \delta {\delta
y^{\left( \alpha \right) i}}}\Bbb{F}=0, (\alpha =1,...,k)$. \vspace{3mm}\\

Since $C\Gamma (N)$ has the property $D_X\Bbb{G}=0,$ it follows:

\begin{teo}
The structures $\Bbb{G}$ and $\Bbb{F}$ determine a Riemannian $(k-1)n$-
almost contact structure on $T^kM$, which depends only on the fundamental
Finsler function of Finsler space $F^{\left( k\right) n}=(M,F).$
\end{teo}

In the case $k=1,$ the triple $(\widetilde{TM},\Bbb{F}, \Bbb{G})$ is an almost
K\"ahlerian space.

\chapter{The Geometry of the Dual of $k$-Tangent Bundle}

\markboth{\it{THE GEOMETRY OF HIGHER-ORDER HAMILTON SPACES\ \ \ \ \ }}{\it{The Geometry of the Dual of
$k$-Tangent Bundle}}

In the book [115] one studies the geometry of the dual of the $2$-tangent bundle,
which has been used to investigate the Hamilton spaces of order $2$. Now,
we consider this problem for the general case, $k\geq 1$.

The dual bundle $T^{*k}M$ of the $k$-tangent bundle must have the
same properties as the cotangent bundle $T^{*}M$ with respect to tangent
bundle $TM$. Such that the manifold $T^{*k}M$ should have the same dimension
$(k+1)n$ as the manifold $T^kM$. $T^{*k}M$ should carry a natural
presymplectic structure and at least one Poisson structure. The
manifolds $T^kM$ and $T^{*k}M$ should be locally diffeomorphic.

The dual bundle $T^{*k}M$ plays a main role in construction of the notion of
Hamilton space of order $k$.

In the present chapter we introduce the bundle $T^{*k}M$ and point out the
main geometrical natural object fields, that live on the differentiable
manifold $T^kM$.

\section{The Dual Bundle $(T^{*k}M,\pi ^{*k},M)$}

\begin{defi}
We call \textit{the dual bundle} of $k$-tangent bundle $(T^kM,\pi ^k,M)$ the
differentiable bundle $(T^{*k}M,\pi ^{*k},M)$ whose total
space is the fibered product:
\begin{equation}
T^{*k}M=T^{k-1}M\times _MT^{*}M  \tag{4.1.1}
\end{equation}
and for which the canonical projection $\pi ^{*k}$ is
\begin{equation}
\pi ^{*k}=\pi ^{k-1}\times _M\pi ^{*}.  \tag{4.1.1a}
\end{equation}
\end{defi}

The previous fibered product has a differentiable structure given by that of
the $(k-1)$-tangent bundle $(T^{k-1}M,\pi ^{k-1},M)$ and the cotangent
bundle $(T^{*}M,\pi ^{*},M)$.
For $k=1$, we have $T^{*1}M=T^{*}M$ and $\pi ^{*1}=\pi ^{*}$.
Sometimes we denote $(T^{*k}M,\pi ^{*k},M)$ by $T^{*k}M$.

A point $u\in T^{*k}M$ will be denoted by $u=(x,y^{(1)},...,y^{(k-1)},p)$,
$\pi ^{*k}(u)=x$.

The projections on the factors of the $T^{*k}M$ from (4.1.1) are
\[
\begin{array}{l}
\pi _{k-1}^{*k}:T^{*k}M\rightarrow T^{k-1}M,\quad \pi
_{k-1}^{*k}(x,y^{(1)},...,y^{(k-1)},p)=(x,y^{(1)},...,y^{(k-1)}), \vspace{3mm}\\
\overline{\pi }^{*}:T^{*k}M\rightarrow T^{*}M,\quad \overline{\pi }
^{*}(x,y^{(1)},...,y^{(k-1)},p)=(x,p)
\end{array}
\]
and the canonical projection $\pi ^{*k}:T^{*k}M\rightarrow M$.

Therefore, the following diagram is commutative:
\[
\begin{array}{ccccc}
&  & T^{*k}M &  &  \\
& \stackrel{\pi _{k-1}^{*k}}{\swarrow } & \mid \ \ \text{\quad } & \stackrel{%
\overline{\pi }^{*}}{\searrow } &  \\
T^{k-1}M &  & \mid \pi ^{*k} &  & T^{*}M \\
& \underset{\pi ^{k-1}}{\searrow } & \downarrow \ \quad \  &
\underset{
\pi ^{*}}{\swarrow } &  \\
&  & M &  &
\end{array}
\]

Let $(x^i,y^{(1)i},...,y^{(k-1)i},p_i)$, $(i=1,2,...,n=\dim M)$ be the
coordinates of a point $u=(x,y^{(1)},...,y^{(k-1)},p)$ $\in T^{*k}M$ in a
local chart $\left( \left( \pi ^{*k}\right) ^{-1}(U),\Phi \right) $ on $%
T^{*k}M$.

The change of coordinates on the manifold $T^{*k}M$ is:

\begin{equation}
\begin{array}{l}
\widetilde{x}^i=\widetilde{x}^i(x^1,...,x^n),\ \det \left( \displaystyle
\frac{\partial \widetilde{x}^i}{\partial x^j}\right) \neq 0, \\
\widetilde{y}^{(1)i}=\displaystyle\frac{\partial \widetilde{x}^i}{\partial
x^j}y^{(1)j}, \\
..................................................................................
\\
(k-1)\widetilde{y}^{(k-1)i}=\displaystyle\frac{\partial \widetilde{y}
^{(k-2)i}}{\partial x^j}y^{(1)j}+\cdots +(k-1)\displaystyle\frac{\partial
\widetilde{y}^{(k-2)i}}{\partial y^{(k-2)j}}y^{(k-1)j}, \\
\widetilde{p}_i=\displaystyle\frac{\partial x^j}{\partial \widetilde{x}^i}
p_j.
\end{array}
\tag{4.1.2}
\end{equation}
where the following equalities hold:
\begin{equation}
\displaystyle\frac{\partial \widetilde{y}^{(\alpha )i}}{\partial x^j}= %
\displaystyle\frac{\partial \widetilde{y}^{(\alpha +1)i}}{\partial y^{(1)j}}
=\cdots =\displaystyle\frac{\partial \widetilde{y}^{(k-1)i}}{\partial
y^{(k-1-\alpha )j}};\ (\alpha =0,...,k-2;y^{(0)}=x).  \tag{4.1.3}
\end{equation}
Sometimes (cf. Ch. 1) $y^{(1)i}$, ..., $y^{(k-1)i}$ will be called
\textit{accelerations of order }$1$\textit{, }$2$\textit{, ..., }$k-1$,
respectively and $p_i$ will be called \textit{momenta}. $T^{*k}M$ is a real
manifold of dimension $(k+1)n$. So it has the same dimension as that of the
manifold $T^kM$.

The natural basis of the vector space $T_u(T^{*k}M)$ at the point $u\in
T^{*k}M$, $\left\{ \displaystyle\frac \partial {\partial x^i}|_u, %
\displaystyle\frac \partial {\partial y^{(1)i}}|_u,...,\displaystyle\frac
\partial {\partial y^{(k-1)i}}|_u,\displaystyle\frac \partial {\partial
p_i}|_u\right\} $ is transformed under (4.1.2) as follows:
\begin{equation}
\begin{array}{l}
\displaystyle\frac \partial {\partial x^i}=\displaystyle\frac{\partial
\widetilde{x}^j}{\partial x^i}\displaystyle\frac \partial {\partial
\widetilde{x}^j}+\displaystyle\frac{\partial \widetilde{y}^{(1)j}}{\partial
x^i}\displaystyle\frac \partial {\partial \widetilde{y}^{(1)j}}+\cdots + %
\displaystyle\frac{\partial \widetilde{y}^{(k-1)j}}{\partial x^i} %
\displaystyle\frac \partial {\partial \widetilde{y}^{(k-1)j}}+\displaystyle
\frac{\partial \widetilde{p}_j}{\partial x^i}\displaystyle\frac \partial
{\partial \widetilde{p}_j}, \\
\\
\displaystyle\frac \partial {\partial y^{(1)i}}=\displaystyle\frac{\partial
\widetilde{y}^{(1)j}}{\partial y^{(1)i}}\displaystyle\frac \partial
{\partial \widetilde{y}^{(1)j}}+\displaystyle\frac{\partial \widetilde{y}
^{(2)j}}{\partial y^{(1)i}}\displaystyle\frac \partial {\partial \widetilde{%
y }^{(2)j}}+\cdots +\displaystyle\frac{\partial \widetilde{y}^{(k-1)j}}{
\partial y^{(1)i}}\displaystyle\frac \partial {\partial \widetilde{y}
^{(k-1)j}}, \\
...............................................................................................
\\
\displaystyle\frac \partial {\partial y^{(k-2)i}}=\displaystyle\frac{
\partial \widetilde{y}^{(k-2)j}}{\partial y^{(k-2)i}}\displaystyle\frac
\partial {\partial \widetilde{y}^{(k-2)j}}+\displaystyle\frac{\partial
\widetilde{y}^{(k-1)j}}{\partial y^{(k-2)i}}\displaystyle\frac \partial
{\partial \widetilde{y}^{(k-1)j}}, \\
\\
\displaystyle\frac \partial {\partial y^{(k-1)i}}=\displaystyle\frac{
\partial \widetilde{y}^{(k-1)j}}{\partial y^{(k-1)i}}\displaystyle\frac
\partial {\partial \widetilde{y}^{(k-1)j}}, \\
\\
\displaystyle\frac \partial {\partial p_i}=\displaystyle\frac{\partial x^i}{
\partial \widetilde{x}^j}\displaystyle\frac \partial {\partial \widetilde{p}
_j}
\end{array}
\tag{4.1.4}
\end{equation}
calculated at the point $u\in T^{*k}M$.

The Jacobian matrix of the transformation (4.1.2) at the point $u\in T^{*k}M$
is:
\begin{equation}
J_k=\left(
\begin{array}{cccccc}
\displaystyle\frac{\partial \widetilde{x}^j}{\partial x^i} & 0 & 0 & \cdots
& 0 & 0 \vspace{3mm}\\
\displaystyle\frac{\partial \widetilde{y}^{(1)j}}{\partial x^i} & %
\displaystyle\frac{\partial \widetilde{y}^{(1)j}}{\partial y^{(1)i}} & 0 &
\cdots & 0 & 0 \\

\cdots & \cdots & \cdots & \cdots & \cdots & \cdots \\
\displaystyle\frac{\partial \widetilde{y}^{(k-1)j}}{\partial x^i} & %
\displaystyle\frac{\partial \widetilde{y}^{(k-1)j}}{\partial y^{(1)i}} & %
\displaystyle\frac{\partial \widetilde{y}^{(k-1)j}}{\partial y^{(2)i}} &
\cdots & \displaystyle\frac{\partial \widetilde{y}^{\left( k-1\right) j}}{
\partial y^{\left( k-1\right) i}} & 0 \vspace{3mm}\\
\displaystyle\frac{\partial \widetilde{p}_j}{\partial x^i} & 0 & 0 & \cdots
& 0 & \displaystyle\frac{\partial x^j}{\partial \widetilde{x}^i}
\end{array}
\right) .  \tag{4.1.5}
\end{equation}

It follows

\begin{equation}
\det J_k(u)=\left[ \det \left( \displaystyle\frac{\partial \widetilde{x}^j}{
\partial x^i}(u)\right) \right] ^{k-1}.  \tag{4.1.5a}
\end{equation}

\begin{teo}
1$^0$ If $k$ is an odd number, then $T^{*k}M$ is an orientable manifold.

2$^0$ If $k$ is an even number, the manifold $T^{*k}M$ is orientable if and
only if the base manifold $M$ is orientable.
\end{teo}

The form (4.1.5) of the Jacobian matrix implies the following transformation,
with respect to (4.1.2), of the natural cobasis $\left\{
dx^i,dy^{(1)i},...,dy^{(k-1)i},dp_i\right\} $, at a point $u\in T^{*k}M$:
\begin{equation}
\begin{array}{l}
d\widetilde{x}^i=\displaystyle\frac{\partial \widetilde{x}^i}{\partial x^j}
dx^j, \\
\\
d\widetilde{y}^{(1)i}=\displaystyle\frac{\partial \widetilde{y}^{(1)i}}{
\partial x^j}dx^j+\displaystyle\frac{\partial \widetilde{y}^{(1)i}}{\partial
y^{(1)j}}dy^{(1)j}, \\
.......................................................................................\\
d\widetilde{y}^{(k-1)i}=\displaystyle\frac{\partial \widetilde{y}^{(k-1)i}}{
\partial x^j}dx^j+\displaystyle\frac{\partial \widetilde{y}^{(k-1)i}}{
\partial y^{(1)j}}dy^{(1)j}+\cdots +\displaystyle\frac{\partial \widetilde{y}
^{(k-1)i}}{\partial y^{(k-1)j}}dy^{(k-1)j}, \\
\\
d\widetilde{p}_j=\displaystyle\frac{\partial \widetilde{p}_j}{\partial x^i}
dx^i+\displaystyle\frac{\partial x^i}{\partial \widetilde{x}^j}dp_i.
\end{array}
\tag{4.1.6}
\end{equation}

Also, we can prove without difficulties the following theorem:

\begin{teo}
If the differentiable manifold $M$ is paracompact, then the differentiable
manifold $T^{*k}M$ is paracompact.
\end{teo}

Consider the category $Man$ of differentiable manifolds.

There exists a covariant functor $T^{*k}:Man\rightarrow Man$ in which the
differentiable mappings $f:M\rightarrow N$, analytical expressed by $%
x^{i^{\prime }}=x^{i^{\prime }}(x^1,...,x^n)$, ($i^{\prime }$, $j^{\prime
}=1^{\prime }$, $2^{\prime }$, ..., $n^{\prime }=\dim N$) give the mappings $%
T^{*k}f:T^{*k}M\rightarrow T^{*k}N$, in the form

\begin{equation}
\begin{array}{l}
x^{i^{\prime }}=x^{i^{\prime }}(x^1,...,x^n),\  \\
\\
y^{(1)i^{\prime }}=\displaystyle\frac{\partial x^{i^{\prime }}}{\partial x^j}
y^{(1)j}, \\
.....................................................................................
\\
(k-1)y^{(k-1)i^{\prime }}=\displaystyle\frac{\partial y^{(k-2)i^{\prime }}}{
\partial x^j}y^{(1)j}+\cdots +(k-1)\displaystyle\frac{\partial
y^{(k-2)i^{\prime }}}{\partial y^{(k-2)j}}y^{(k-1)j}, \\
\\
\displaystyle\frac{\partial x^{i^{\prime }}}{\partial x^j}p_{i^{\prime
}}=p_j.
\end{array}
\tag{4.1.7}
\end{equation}

This fact will be used in the theory of subspaces in Hamilton
spaces of order $k$.

\section{Vertical Distributions. Liouville Vector Fields}

The null section $\mathbf{0}:$ $M\rightarrow T^{*k}M$ of the projection $\pi
^{*k}$ is defined by $\mathbf{0}: x \in M\rightarrow (x,0,...,0)\in T^{*k}M$.
As usual we denote $\widetilde{T^{*k}M}=T^{*k}M\setminus \{\mathbf{0}\}$.

The tangent bundle of the manifold $T^{*k}M$, $(TT^{*k}M,d\pi ^{*k},TM)$,
allows to define the vertical subbundle $VT^{*k}M=\ker d\pi ^{*k}$. We get the
vertical distribution $V$, formed by the fibres of $VT^{*k}M$. $V$ is locally
generated by the set of vector fields $\left\{ \displaystyle\frac \partial
{\partial y^{(1)i}},\cdots ,\displaystyle\frac \partial {\partial
y^{(k-1)i}},\displaystyle\frac \partial {\partial p_i}\right\} $ at every
point $u\in T^{*k}M$.

So, $V$ is an integrable distribution and its dimension is $kn$.

It is convenient to adopt the notation:
\begin{equation}
\stackrel{\cdot }{\partial ^i}=\displaystyle\frac \partial {\partial p_i}.
\tag{4.2.1}
\end{equation}

Taking into account the relations (4.1.4) we can consider the following
subdistributions of $V$:

$V_{k-1}$, locally generated by $\left\{ \displaystyle\frac \partial
{\partial y^{(k-1)i}}\right\} $. Its dimension is $n$ and it is integrable.

$V_{k-2}$, locally generated by $\left\{ \displaystyle\frac \partial
{\partial y^{(k-2)i}},\displaystyle\frac \partial {\partial
y^{(k-1)i}}\right\} $. This has dimension $2n$ and it is integrable, too and so
on.

$V_1$, locally generating by $\left\{ \displaystyle\frac \partial {\partial
y^{(1)i}},\cdots ,\displaystyle\frac \partial {\partial y^{(k-1)i}}\right\}. $
Its dimension is $(k-1)n$ and it is also integrable. Of course, we have
the sequence of inclusions:

\[
V_{k-1}\subset V_{k-2}\subset \cdots \subset V_1\subset V.
\]

The transformation of vector field $\stackrel{\cdot }{\partial ^i}$ from
(4.1.4), $\stackrel{\cdot }{\partial ^i}=\displaystyle\frac{\partial x^i}{
\partial \widetilde{x}^j}\stackrel{\cdot }{\widetilde{\partial }^j}$, shows
that we have one more vertical distribution $W_k$, locally generated by the
vector fields $\left\{ \stackrel{\cdot }{\partial ^i}\right\} $ at the
points $u\in \left( \pi ^{*k}\right) ^{-1}(U)$. Its dimension is $n$ and it
is an integrable distribution, too.

We conclude by

\begin{prop}
The following direct sum of vector spaces holds:
\begin{equation}
V_u=V_{1,u}\oplus W_{1,u}\ ,\quad \forall u\in T^{*k}M.  \tag{4.2.2}
\end{equation}
\end{prop}

Using again (4.1.4) we obtain without difficulties

\begin{teo}
1$^0$ The following operators in the algebra of functions \newline $\mathcal{F}%
(T^{*k}M)$:
\begin{equation}
\begin{array}{l}
\stackrel{1}{\Gamma }=y^{(1)i}\displaystyle\frac \partial {\partial
y^{(k-1)i}}, \\
\\
\stackrel{2}{\Gamma }=y^{(1)i}\displaystyle\frac \partial {\partial
y^{(k-2)i}}+2y^{(2)i}\displaystyle\frac \partial {\partial y^{(k-1)i}}, \\
.............................................................. \\
\stackrel{k-1}{\Gamma }=y^{(1)i}\displaystyle\frac \partial {\partial
y^{(1)i}}+\cdots +(k-1)y^{(k-1)i}\displaystyle\frac \partial {\partial
y^{(k-1)i}}
\end{array}
\tag{4.2.3}
\end{equation}
and
\begin{equation}
C^{*}=p_i\stackrel{\cdot }{\partial ^i}  \tag{4.2.4}
\end{equation}
are vector fields on $T^{*k}M$. They are independent vector fields on the
manifold $\widetilde{T^{*k}M}$.

2$^0$ The function
\begin{equation}
\varphi =p_iy^{(1)i}  \tag{4.2.5}
\end{equation}
is a scalar function on the manifold $T^{*k}M$.
\end{teo}

Evidently, $\stackrel{1}{\Gamma }$ belongs to the distribution $V_{k-1}$, $%
\stackrel{2}{\Gamma }$ belongs to $V_{k-2}$, ..., $\stackrel{k-1}{\Gamma }$
belongs to the distribution $V_1$ and the vector field $C^{*}$ belongs to
the distribution $W_k$.

$\stackrel{1}{\Gamma }$, ..., $\stackrel{k-1}{\Gamma }$ are called \textit{\
the Liouville vector fields} and $C^{*}$ is \textit{the Hamilton vector
field } on $T^{*k}M$.

The Liouville vector fields $\stackrel{1}{\Gamma }$, ..., $\stackrel{k-1}{
\Gamma }$ are linearly independent on $\widetilde{T^{*k}M}$ and exactly as
in ch. 1 we can prove:

\begin{teo}

For any differentiable function $H:\widetilde{T^{*k}M}$ $\rightarrow \mathbf{%
\ R}$, $d_0H$, $d_1H$, ..., $d_{k-2}H$ defined by
\begin{equation}
\begin{array}{l}
d_0H=\displaystyle\frac{\partial H}{\partial y^{(k-1)i}}dx^i, \\
\\
d_1H=\displaystyle\frac{\partial H}{\partial y^{(k-2)i}}dx^i+\displaystyle
\frac{\partial H}{\partial y^{(k-1)i}}dy^{(1)i}, \\
....................................................................................
\\
d_{k-2}H=\displaystyle\frac{\partial H}{\partial y^{(1)i}}dx^i+\displaystyle
\frac{\partial H}{\partial y^{(2)i}}dy^{(1)i}+\cdots +\displaystyle\frac{%
\partial H}{\partial y^{(k-1)i}}dy^{(k-2)i}
\end{array}
\tag{4.2.6}
\end{equation}
are fields of $1$-form on $\widetilde{T^{*k}M}$.
\end{teo}

\begin{prop}
1$^0$ $d_{k-1}H$ given by
\begin{equation}
d_{k-1}H=\displaystyle\frac{\partial H}{\partial x^i}dx^i+\displaystyle\frac{%
\partial H}{\partial y^{(1)i}}dy^{(1)i}+\cdots +\displaystyle\frac{\partial H%
}{\partial y^{(k-1)i}}dy^{(k-1)i}  \tag{4.2.7}
\end{equation}
is not a field of $1$-form.

2$^0$ Under a transformation of coordinate on $T^{*k}M$, $d_{k-1}H$ transforms as follows
\begin{equation}
d_{k-1}H=d_{k-1}\widetilde{H}+\stackrel{\cdot }{\widetilde{\partial }^j}%
\widetilde{H}\displaystyle\frac{\partial \widetilde{p}_j}{\partial x^i}dx^i.
\tag{4.2.8}
\end{equation}

3$^0$ If $\stackrel{\cdot }{\partial ^i}H=0$, then $d_{k-1}H$ is a field of $1$-form.
\end{prop}

\begin{prop}
The relation between the differential $dH$ and $d_{k-1}H$ are given by
\begin{equation}
dH=d_{k-1}H+\stackrel{\cdot }{\partial ^i}Hdp_i.  \tag{4.2.9}
\end{equation}
\end{prop}

Indeed,
$$
dH=\displaystyle\frac{\partial H}{\partial x^i}dx^i+\displaystyle\frac{
\partial H}{\partial y^{(1)i}}dy^{(1)i}+\cdots +\displaystyle\frac{\partial
H }{\partial y^{(k-1)i}}dy^{(k-1)i}+\stackrel{\cdot }{\partial ^i}Hdp_i=
$$
$$
=d_{k-1}H+\stackrel{\cdot }{\partial ^i}Hdp_i.
$$
\textit{Remark} The differential $dH$ being invariant under the coordinate
transformations (4.1.2), i.e. $dH=d\widetilde{H}$, from (4.2.9) it
follows the rule of transformation (4.2.8) of $d_{k-1}H$.

If $H=\varphi =p_iy^{(1)i}$, then $d_0H=\cdots =d_{k-3}H=0$ and
\begin{equation}
\omega =d_{k-2}\varphi =p_idx^i.  \tag{4.2.10}
\end{equation}
$\omega $ is called \textit{the Liouville }$1$\textit{-form} on the manifold
$\widetilde{T^{*k}M}$.

The exterior differential $d\omega $ of the Liouville $1$-form $\omega $ is
expressed by
\begin{equation}
\theta =d\omega =dp_i\wedge dx^i.  \tag{4.2.11}
\end{equation}

Using (4.1.6) we can prove the invariance of $2$-form $\theta $ with respect
to (4.1.2).

Now, based on the previous results we obtain:

\begin{teo}
1$^0$ The differential forms $\omega $ and $\theta $ are globally defined on
the manifold $\widetilde{T^{*k}M}$.

2$^0$ $\theta $ is a closed $2$-form, i.e. $d\theta =0$.

3$^0$ $\theta $ is a $2$-form of rank $2n$. It is a presymplectic structure
on $\widetilde{T^{*k}M}$.
\end{teo}

\textit{Proof}\textbf{:} 1$^0$ $\omega $ and $\theta $ are invariant with
respect to a change of coordinates (4.1.2).

2$^0$ $\theta =d\omega $ implies $d\theta =0$.

3$^0$ $\theta =dp_i\wedge dx^i$ is a $2$-form of rank $2n<(k+1)n=\dim
T^{*k}M $, for $k>1$. Consequently $\theta $ is a presymplectic structure on
$T^{*k}M $.

\textit{Remarks} 1$^0$ For $k=1$, $\omega $ and $\theta $ are the

Poincar\'e-Cartan forms on the cotangent bundle $T^{*}M$.

2$^0$ The previous theorem shows the existence of a natural presymplectic
structure on $\widetilde{T^{*k}M}$.

\section{The Structures $J$ and $J^{*}$}

There exists a tangent structure $J$ on $T^{*k}M$ defined as usual by the
endomorphism $J:\mathcal{X}(T^{*k}M)\rightarrow \mathcal{X}(T^{*k}M)$:

\begin{equation}
\begin{array}{l}
J\left( \displaystyle\frac \partial {\partial x^i}\right) =\displaystyle %
\frac \partial {\partial y^{(1)i}},J\left( \displaystyle\frac \partial
{\partial y^{(1)i}}\right) =\displaystyle\frac \partial {\partial
y^{(2)i}},..., \\
\\
J\left( \displaystyle\frac \partial {\partial y^{(k-2)i}}\right) = %
\displaystyle\frac \partial {\partial y^{(k-1)i}},J\left( \displaystyle\frac
\partial {\partial y^{(k-1)i}}\right) =0,J(\stackrel{\cdot }{\partial ^i})=0
\end{array}
\tag{4.3.1}
\end{equation}
at every point $u\in \widetilde{T^{*k}M}$.

By means of (4.1.4) one proves without difficulties

\begin{teo}
1$^0$ $J$ is globally defined on $T^{*k}M$.

2$^0$ $J$ is a tensor field of type $(1,1)$, locally expressed by
\begin{equation}
J=\displaystyle\frac \partial {\partial y^{(1)i}}\otimes dx^i+\displaystyle %
\frac \partial {\partial y^{(2)i}}\otimes dy^{(1)i}+\cdots +\displaystyle %
\frac \partial {\partial y^{(k-1)i}}\otimes dy^{(k-2)i}.  \tag{4.3.2}
\end{equation}

3$^0$ The structure $J$ is integrable.

4$^0$ $J\circ J\circ \cdots \circ J=J^k=0$.

5$^0$ $\ker J=V_{k-1}\oplus W_k$, $Im$ $J=V_1$.

6$^0$ $rank\ ||J||=(k-1)n$.
\end{teo}

According to the above theorem we may call $J$ \textit{the }$k-1$\textit{\
-tangent structure}.

The endomorphism $J$ applied to the Liouville vector fields gives us:
\begin{equation}
J(\stackrel{1}{\Gamma })=0,\ J(\stackrel{2}{\Gamma })=\stackrel{1}{\Gamma }
,\ ...,\ J(\stackrel{k-1}{\Gamma })=\stackrel{k-2}{\Gamma },  \tag{4.3.3}
\end{equation}
and
\begin{equation}
J(C^{*})=0.  \tag{4.3.3a}
\end{equation}

Let us consider a vector field $X\in \mathcal{X}(T^{*k}M)$, locally
expressed by:

\begin{equation}
X=\stackrel{(0)i}{X}\displaystyle\frac \partial {\partial x^i}+\stackrel{%
(1)i }{X}\displaystyle\frac \partial {\partial y^{(1)i}}+\cdots +\stackrel{%
(k-1)i }{X}\displaystyle\frac \partial {\partial y^{(k-1)i}}+X_i\stackrel{%
\cdot }{ \partial ^i}.  \tag{4.3.4}
\end{equation}

\begin{prop}
1$^0$ For any vector field $X\in \mathcal{X}(T^{*k}M)$, $\stackrel{1}{X}$,
..., $\stackrel{k-1}{X}$ \newline given by
\begin{equation}
\stackrel{1}{X}=JX,\ \stackrel{2}{X}=J^2X,\ ...,\ \stackrel{k-1}{X}=J^{k-1}X.
\tag{4.3.4a}
\end{equation}
are vector fields.

2$^0$ If $X$ is given by (4.3.4), then we have:
\[
\begin{array}{l}
\stackrel{1}{X}=\stackrel{(0)i}{X}\displaystyle\frac \partial {\partial
y^{(1)i}}+\cdots +\stackrel{(k-2)i}{X}\displaystyle\frac \partial {\partial
y^{(k-1)i}}, \\
\\
\stackrel{2}{X}=\stackrel{(0)i}{X}\displaystyle\frac \partial {\partial
y^{(2)i}}+\cdots +\stackrel{(k-3)i}{X}\displaystyle\frac \partial {\partial
y^{(k-1)i}}, \\
.................................................. \\
\stackrel{k-1}{X}=\stackrel{(0)i}{X}\displaystyle\frac \partial {\partial
y^{(k-1)i}}.
\end{array}
\]

3$^0$ The vector field $\stackrel{1}{X}$ belongs to the vertical
distribution $V_1$, $\stackrel{2}{X}$ belongs to the distribution $V_2$,
..., $\stackrel{k-1}{X}$ belongs to the distribution $V_{k-1}$.
\end{prop}

Now consider the adjoint of $k-1$-tangent structure $J$. It is the
endomorphism $J^{*}:\mathcal{X}^{*}(T^{*k}M)\rightarrow \mathcal{X}
^{*}(T^{*k}M)$ defined by

\begin{equation}
\begin{array}{l}
J^{*}(dy^{(k-1)i})=dy^{(k-2)i},\ ...,\ J^{*}(dy^{(1)i})=dx^i, \\
\\
J^{*}(dx^i)=0,\ J^{*}(dp_i)=0.
\end{array}
\tag{4.3.5}
\end{equation}

Using (4.3.5) and (4.1.6) we obtain:

\begin{teo}
1$^0$ $J^{*}$ is globally defined on $T^{*k}M$.

2$^0$ $J^{*}$ is a tensor field of type $(1,1)$ on $T^{*k}M$, i.e.
\begin{equation}
J^{*}=dx^i\otimes \displaystyle\frac \partial {\partial
y^{(1)i}}+dy^{(1)i}\otimes \displaystyle\frac \partial {\partial
y^{(2)i}}+\cdots +dy^{(k-2)i}\otimes \displaystyle\frac \partial {\partial
y^{(k-1)i}}.  \tag{4.3.6}
\end{equation}

3$^0$ $rank$ $||J^{*}||=(k-1)n$.

4$^0$ $J^{*}$ is an integrable structure.
\end{teo}

$J^{*}$ is called\textit{\ the }$k-1$\textit{-adjoint tangent structure}.

$J^{*}$ can be extended to an endomorphism of the exterior algebra $\Lambda
(T^{*k}M)$ as follows:

\begin{equation}
\begin{array}{l}
J^{*}f=f,\ \forall f\in \mathcal{F}(T^{*k}M), \\
\\

(J^{*}\omega )(X_1,...,X_q)=\omega (JX_1,...,JX_q),\ \forall \omega \in
\Lambda ^q(T^{*k}M).
\end{array}
\tag{4.3.7}
\end{equation}

Let be $\omega \in \Lambda ^1(T^{*k}M)$ and consider

\begin{equation}
\stackrel{1}{\omega }=J^{*}\omega ,\ ...,\ \stackrel{k-1}{\omega }
=J^{*(k-1)}\omega .  \tag{4.3.7a}
\end{equation}

Then $\stackrel{1}{\omega }$, ..., $\stackrel{k-1}{\omega }$ are $1$-form
fields.

In particular, we get
\begin{equation}
J^{*}dH=d_{k-2}H,\ ...,\ J^{*(k-1)}dH=d_0H.  \tag{4.3.8}
\end{equation}

The $k-1$ adjoint structure $J^{*}$ allows to introduce the vertical
differential operator in the exterior algebra $\Lambda (T^{*k}M)$, [74].

Taking into account the operator of differentiation:

$d=\displaystyle\frac \partial {\partial x^i}dx^i+\displaystyle\frac
\partial {\partial y^{(1)i}}dy^{(1)i}+\cdots +\displaystyle\frac \partial
{\partial y^{(k-1)i}}dy^{(k-1)i}+\stackrel{\cdot }{\partial ^i}dp_i$ \\
we define the following operators:
\begin{equation}
\begin{array}{l}
d_{k-2}=J^{*}d=\displaystyle\frac \partial {\partial y^{(1)i}}dx^i+\cdots + %
\displaystyle\frac \partial {\partial y^{(k-1)i}}dy^{(k-2)i}, \\
\\
d_{k-3}=J^{*2}d=\displaystyle\frac \partial {\partial y^{(2)i}}dx^i+\cdots + %
\displaystyle\frac \partial {\partial y^{(k-1)i}}dy^{(k-3)i}, \\
......................................................................... \\
d_1=J^{*(k-2)}d=\displaystyle\frac \partial {\partial y^{(k-2)i}}dx^i+ %
\displaystyle\frac \partial {\partial y^{(k-1)i}}dy^{(1)i}, \\
\\
d_0=J^{*(k-1)}d=\displaystyle\frac \partial {\partial y^{(k-1)i}}dx^i.
\end{array}
\tag{4.3.9}
\end{equation}

Clearly, these operators $d_0$, ..., $d_{k-2}$ and $d$ do not depend on the
transformation of coordinates on the manifold $T^{*k}M$.

For any $H\in \mathcal{F}(T^{*k}M)$, $d_0H$, ..., $d_{k-2}H$ are given by
Theorem 4.2.2. $d_0$, ..., $d_{k-2}$ are \textit{the vertical
operators of differentiation}.

They can be extended to the exterior algebra $\Lambda (T^{*k}M)$ if we give
their restrictions to $\Lambda ^0(T^{*k}M)$ and $\Lambda ^1(T^{*k}M)$.

As we already have seen $d_0H$, ..., $d_{k-2}H$ are expressed in (4.2.6) and $dH$
is the differential of $H$. The restrictions to $\Lambda ^1(T^{*k}M)$ are
defined by

\begin{equation}
\left\{
\begin{array}{l}
d_\alpha (dx^i)=0,\ d_\alpha (dp_i)=0,\ (\alpha =0,...,k-2), \\
\\
d_\alpha (dy^{(\beta )i})=0,\ \ (\alpha =1,...,k-2;\ \beta =1,...,k-1), \\

\\
d(dx^i)=0,\ d(dy^{(\beta )i})=0,\ d(dp_i)=0.
\end{array}
\right.  \tag{4.3.10}
\end{equation}

In this case $d_0$, ..., $d_{k-2}$ and $d$ are the antiderivations of degree
$1$.

\begin{prop}
The vertical differential operators $d_0$, ..., $d_{k-2}$ and $d$ have the
property
\begin{equation}
d_\alpha \circ d_\alpha =0,\ (\alpha =0,...,k-2),\ d\circ d=0.  \tag{4.3.11}
\end{equation}
\end{prop}

For instance, applying the exterior differential $d$ to the $1$-forms $d_0H$
, $d_1H$, ..., $d_{k-2}H$, written in the form

\begin{equation}
\begin{array}{l}
d_0H=\stackrel{(0)}{p_i}dx^i, \\
\\
d_1H=\stackrel{(1)}{p_i}dx^i+\stackrel{(0)}{p_i}dy^{(1)i}, \\
...........................................................................
\\
d_{k-2}H=\stackrel{(k-2)}{p_i}dx^i+\stackrel{(k-3)}{p_i}dy^{(1)i}+\cdots +
\stackrel{(0)}{p_i}dy^{(k-2)i},
\end{array}
\tag{4.3.12}
\end{equation}
with
\begin{equation}
\stackrel{(0)}{p_i}=\displaystyle\frac{\partial H}{\partial y^{(k-1)i}},\
\stackrel{(1)}{p_i}=\displaystyle\frac{\partial H}{\partial y^{(k-2)i}},\
...,\ \stackrel{(k-2)}{p_i}=\displaystyle\frac{\partial H}{\partial y^{(1)i}},  \tag{4.3.12a}
\end{equation}
we obtain
\begin{equation}
\begin{array}{l}
dd_0H=d\stackrel{(0)}{p_i}\wedge dx^i, \\
\\
dd_1H=d\stackrel{(1)}{p_i}\wedge dx^i+d\stackrel{(0)}{p_i}\wedge dy^{(1)i},
\\
...............................................................................................
\\
dd_{k-2}H=d\stackrel{(k-2)}{p_i}\wedge dx^i+d\stackrel{(k-3)}{p_i}\wedge
dy^{(1)i}+\cdots +d\stackrel{(0)}{p_i}\wedge dy^{(k-2)i}.
\end{array}
\tag{4.3.13}
\end{equation}

The previous formulae (4.3.12) are similar to the Jacobi-Ostrogradski
formulae (2.4.2) from the Lagrange spaces of order $k$, $L^{(k)n}$. The
equalities (4.3.13) are important in the geometrical theory of the
Hamiltonians $H(x,y^{(1)},...,y^{(k-1)},p)$.

\section{Canonical Poisson Structures on $T^{*k}M$}

In this section we state that on the manifold $T^{*k}M$ there exists at least a Poisson structure.
Let us consider the brackets:
\begin{equation}
\left\{f, g \right\}_0=\displaystyle\frac{\partial f}{\partial x^i}
\displaystyle\frac{\partial g}{\partial p_i}-
\displaystyle\frac{\partial f}{\partial p_i}
\displaystyle\frac{\partial g}{\partial x^i}, f, g \in {\cal F}(T^{*k}M),
\tag{4.4.1}
\end{equation}
\begin{equation}
\left\{f, g \right\}_\alpha =\displaystyle\frac{\partial f}{\partial y^{(\alpha )i}}
\displaystyle\frac{\partial g}{\partial p_i}-
\displaystyle\frac{\partial f}{\partial p_i}
\displaystyle\frac{\partial g}{\partial y^{(\alpha )i}}, (\alpha =1,...,k-1).
\tag{4.4.1a}
\end{equation}

\begin{teo} The bracket $\left\{\cdot,\cdot \right\}_\alpha $, $(\alpha =k-1)$ is a canonical Poisson structure on the manifold $T^{*k}M$.
\end{teo}

\textit{Proof:} Indeed, remarking that with respect to a change of local coordinates on $T^{*k}M$,
$\displaystyle\frac{\partial f}{\partial y^{(\alpha )i}}\displaystyle\frac{\partial g}{\partial p_i}$ has a geometrical meaning, it follows

\vspace*{3mm}

1$^0$ For any $f, g \in {\cal F}(T^{*k}M)$, $\left\{ f,g\right\}_{k-1}$ is a differentiable function on $T^{*k}M$.

\vspace*{3mm}

2$^0$ $\{f,g\}_{k-1} =-\{g,f\}_{k-1}.$

\vspace*{3mm}

3$^0$ $\{f,g \}_{k-1}$ is $\mathbf{R}$-linear in every argument.

\vspace*{3mm}

4$^0$ The Jacobi identities are verified:
\[
\left\{ \left\{ f,g\right\}_{k-1} ,h\right\}_{k-1} +\left\{ \left\{
g,h\right\} _{k-1} ,f\right\} _{k-1} +\left\{ \left\{ h,f\right\} _{k-1}, g\right\} _{k-1} =0.
\]

Also we remark here

\vspace*{3mm}

5$^0$ $\left\{ \cdot ,gh\right\} _{k-1} =\left\{ \cdot ,g\right\} _{k-1} h+\left\{ \cdot ,h\right\} _{k-1} g$.\\
We can see, without dificulties that every bracket $\left\{f, g \right\}_0$ and $\left\{f, g \right\}_{\alpha}$, $(\alpha = 1, ..., k-2)$ verifies 2$^0$, 3$^0$, 4$^0$ and 5$^0$ but they do not satisfy the property 1$^0$.

The restrictions of these brackets to the special submanifolds immersed in $T^{*k}M$ can induce Poisson structures. As we shall see in ch. 8, \S 3, the restriction of $\left\{f, g \right\}_0$, $f,g \in {\cal F}(\Sigma_0)$ to the submanifold $\Sigma_0$:
$$
\Sigma_0=\{(x,y^{(1)}, ..., y^{(k-1)}, p) \in T^{(*k)}M| y^{(1)i}= \cdots = y^{(k-1)i}=0\}
$$
has this property.

Using the notion of nonlinear connection $N$, studied in chapter 6 of the present monograph, we can take an adapted basis $\displaystyle\frac{\delta}{\delta x^i}$, $\displaystyle\frac{\delta}{\delta y^{(\alpha)i}}$, $(\alpha=1, ..., k-2)$, $\displaystyle\frac{\partial}{\partial y^{(k-1)i}}$ and remarking that $\displaystyle\frac{\delta f}{\delta x^i}\displaystyle\frac{\partial g}{\partial p_i}$ and $\displaystyle\frac{\delta f}{\delta y^{(\alpha)i}}\displaystyle\frac{\partial g}{\partial p_i}$, $(\alpha=1, ..., k-2)$, have geometrical meaning we can construct the following brackets:
\begin{equation}
\begin{array}{l}
\left\{f, g \right\}^N_0=\displaystyle\frac{\delta f}{\delta x^i}\displaystyle\frac{\partial g}{\partial p_i}-
\displaystyle\frac{\delta g}{\delta x^i}\displaystyle\frac{\partial f}{\partial p^i}, f, g \in {\cal F}(T^{*k}M), \vspace{3mm}\\
\left\{ f,g\right\}^N_\alpha =\displaystyle\frac{\delta f}{\delta y^{(\alpha )i}}\displaystyle\frac{\partial g}{\partial p_i}-
\displaystyle \frac{\partial f}{\partial p_i}\displaystyle\frac{\delta g}{\delta y^{(\alpha )i}},(\alpha =1,...,k-2).
\end{array}
\tag{4.4.2}
\end{equation}
We can check, without difficulties

\begin{prop}
Every bracket from (4.4.2) has the properties 1$^0$, 2$^0$, 3$^0$, 5$^0$ from Theorem 4.4.1.
\end{prop}

We shall see later that for some particular nonlinear connections, the brackets (4.4.2) have also the property 4$^0$.

\section{Homogeneity}

The notion of homogeneity for the functions $H(x,y^{(1)},...,y^{(k-1)},p)$
defined on the manifold $T^{*k}M$ can be considered with respect to the
vertical variables $y^{(1)i}$, ..., $y^{(k-1)i}$, as well as with respect to
the momenta $p_i$, respectively.

Indeed, any homothety $h_a: T^{*k}M\rightarrow T^{*k}M $
$$
h_a(x,y^{(1)},y^{(2)},...,y^{(k-1)},p)=(x,ay^{(1)},a^2y^{(2)},...,a^{k-1}y^{(k-1)},p)
$$
is preserved by the transformations of local coordinates (4.1.2) on $T^{*k}M$.

Let $H_y$ be a group of transformation on $T^{*k}M$:
\[
H_y=\left\{ h_a |a\in \mathbf{R}^{+}\right\} .
\]

The orbit of a point $u_0=(x_0,y_0^{(1)},...,y_0^{(k-1)},p^0)$ by $H_y$ is
given by
\[
x^i=x_0^i,\ y^{(1)i}=ay_0^{(1)i},\ ...,\ y^{(k-1)i}=a^{k-1}y_0^{(k-1)i},\
p_i=p_i^0,\ a\in \mathbf{R}^{+}.
\]

The tangent vector at the point $u_0=h_1(u_0)$ is the Liouville vector field
$\stackrel{k-1}{\Gamma }$ at a point $u_0$:
\[
\stackrel{k-1}{\Gamma }(u_0)=y_0^{(1)i}\displaystyle\frac \partial {\partial
y^{(1)i}}|_{u_0}+2y_0^{(2)i}\displaystyle\frac \partial {\partial
y^{(2)i}}|_{u_0}+\cdots +(k-1)y_0^{(k-1)i}\displaystyle\frac \partial
{\partial y^{(k-1)i}}|_{u_0}.
\]

\begin{defi}
A function $H:T^{*k}M\rightarrow \mathbf{R}$ differentiable on $\widetilde{%
T^{*k}M}$ and continuous on the null section of the projection $\pi ^{*k}$
is called \textit{homogeneous of degree }$r\in \mathbf{Z}$ with respect to $%
(y^{(1)},...,y^{(k-1)})$ if
\begin{equation}
H\circ h_a=a^rH,\quad \forall a\in \mathbf{R}^{+}.  \tag{4.5.1}
\end{equation}
\end{defi}

Exactly as in the section 3.1, ch. 3, it follows:

\begin{teo}
A function $H\in \mathcal{F}(T^{*k}M)$, differentiable on $\widetilde{T^{*k}M%
}$ and continuous on the null section is $r$-homogeneous with respect to $%
(y^{(1)},...,y^{(k-1)})$ if and only if
\begin{equation}
\mathcal{L}_{\stackrel{k-1}{\Gamma }}H=rH,  \tag{4.5.2}
\end{equation}
where $\mathcal{L}_{\stackrel{k-1}{\Gamma }}$ is the Lie derivative with
respect to the Liouville vector field $\stackrel{k-1}{\Gamma }$.
\end{teo}

Evidently, (4.5.2) can be written in the form
\begin{equation}
y^{(1)i}\displaystyle\frac{\partial H}{\partial y^{(1)i}}+2y^{(2)i} %
\displaystyle\frac{\partial H}{\partial y^{(2)i}}+\cdots +(k-1)y^{(k-1)i} %
\displaystyle\frac{\partial H}{\partial y^{(k-1)i}}=rH.  \tag{4.5.2a}
\end{equation}

As is usually (see Ch.3, \S 3.1) the notion of homogeneity can be extended to
the vector fields on $T^{*k}M$.

A vector field $X$ on $\widetilde{T^{*k}M}$ is $r$-homogeneous with respect
to $(y^{(1)},...,y^{(k-1)})$ if
\begin{equation}
X\circ h_a=a^{r-1}h_a^{*}\circ X,\ \forall a\in \mathbf{R}^{+}.  \tag{4.5.3}
\end{equation}

One proves [115]:

\begin{teo}
A vector field $X\in \mathcal{X}(\widetilde{T^{*k}M})$ is $r$-homogeneous
with respect to $(y^{(1)},...,y^{(k-1)})$ if and only if
\begin{equation}
\mathcal{L}_{\stackrel{k-1}{\Gamma }}X=(r-1)X.  \tag{4.5.4}
\end{equation}
\end{teo}

Of course, $\mathcal{L}_{\stackrel{k-1}{\Gamma }}X=[\stackrel{k-1}{\Gamma }
,X]$.

Consequently, $\displaystyle\frac \partial {\partial x^i}$, $\displaystyle
\frac \partial {\partial y^{(1)i}}$, ...,$\displaystyle\frac \partial
{\partial y^{(k-1)i}}$, $\displaystyle\frac \partial {\partial p_i}$ are $0$
, $1$, ..., $2-k$, $0$ homogeneous, respectively, with respect to $%
(y^{(1)},...,y^{(k-1)})$.

The following properties hold:

1$^0$ If $H$ is $s$-homogeneous and $X\in \mathcal{X}(\widetilde{T^{*k}M})$
is $r$-homogeneous, then $HX$ is $s+r$-homogeneous, with respect to $%
(y^{(1)},...,y^{(k-1)})$.

2$^0$ If $H$ is $s$-homogeneous and $X\in \mathcal{X}(\widetilde{T^{*k}M})$
is $r$-homogeneous, then $XH$ is $s+r-1$-homogeneous, with respect to $%
(y^{(1)},...,y^{(k-1)})$.

A $q$-form $\omega \in \Lambda ^q(\widetilde{T^{*k}M})$ is $s$-homogeneous,
with respect to $(y^{(1)},...,y^{(k-1)})$ if

\begin{equation}
\omega \circ h_a^{*}=a^s\omega ,\ \forall a\in \mathbf{R}^{+}.  \tag{4.5.5}
\end{equation}

As we know, [115] the following theorem holds:

\begin{teo}
A $q$-form $\omega $ is $s$-homogeneous with respect to $%
(y^{(1)},...,y^{(k-1)})$ if and only if
\begin{equation}
\mathcal{L}_{\stackrel{k-1}{\Gamma }}\omega =s\omega .  \tag{4.5.6}
\end{equation}
\end{teo}

The $1$-forms $dx^i$, $dy^{(1)i}$, ..., $dy^{(k-1)i}$, $dp_i$ are
homogeneous of degree $0$, $1$, ..., $k-1$, $0$ respectively, with respect
to $(y^{(1)},...,y^{(k-1)})$.

As we remarked above it is important to study the notion of homogeneity with
respect to momenta.

Let $H_p$ be the group of homothety:
$$
H_p=\left\{ h_a^{\prime }: \widetilde{T^{*k}M}
\rightarrow  \widetilde{T^{*k}M}|a\in \mathbf{%
\ R}^{+}\right\}
$$
where $ h_a^{\prime }(x,y^{(1)},...,y^{(k-1)},p)=(x,y^{(1)},...,y^{(k-1)},ap).$
A function $H\in T^{*k}M$, differentiable on $\widetilde{T^{*k}M}$ and
continuous on the null section of $\pi ^{*k}$ is homogeneous of degree $r$,
with respect to the momenta $p_i$, if

\begin{equation}
H\circ h_a^{\prime }=a^rH,\quad \forall a\in \mathbf{R}^{+}.  \tag{4.5.7}
\end{equation}

In other words:

\begin{equation}
H(x,y^{(1)},...,y^{(k-1)},ap)=a^rH(x,y^{(1)},...,y^{(k-1)},p),\quad \forall
a\in \mathbf{R}^{+}.  \tag{4.5.7a}
\end{equation}

We have, [115]:

\begin{teo}
A function $H\in \mathcal{F}(T^{*k}M)$, differentiable on $\widetilde{T^{*k}M%
}$ and continuous on the null section is $r$-homogeneous with respect to $p_i
$ if and only if
\begin{equation}
\mathcal{L}_{C^{*}}H=rH.  \tag{4.5.8}
\end{equation}
\end{teo}

But $\mathcal{L}_{C^{*}}H=C^{*}H=p_i\displaystyle\frac{\partial H}{\partial
p_i}$.

\textit{Remark} It is not difficult to see that if $H$ is differentiable on $%
T^{*k}M$, then the $r$-homogeneous function $H$ is a polinom of degree $r$
in the variables $p_i$.

A vector field $X\in \mathcal{X}(\widetilde{T^{*k}M})$ is homogeneous of
degree $r$ if
\[
X\circ h_a^{\prime }=a^{r-1}h_a^{\prime *}\circ X,\ \forall a\in \mathbf{R}
^{+}.
\]

We have

\begin{teo}
A vector field $X\in \mathcal{X}(\widetilde{T^{*k}M})$ is $r$-homogeneous
with respect to $p_i$ if and only if:
\begin{equation}
\mathcal{L}_{C^{*}}X=(r-1)X.  \tag{4.5.9}
\end{equation}
\end{teo}

This result implies that the vector fields $\displaystyle\frac \partial
{\partial x^i}$, $\displaystyle\frac \partial {\partial y^{(1)i}}$, ...,$%
\displaystyle\frac \partial {\partial y^{(k-1)i}}$, $\displaystyle\frac
\partial {\partial p_i}$ are $1$, $1$, ..., $1$, $0$ homogeneous,
respectively, with respect to $p_i$.

If $H$ is $s$-homogeneous and $X$ is $r$-homogeneous, then $fX$ is $s+r$
-homogeneous and $Xf$ is $s+r-1$-homogeneous, with respect to $p_i$.

For instance, if $H(x,y^{(1)},...,y^{(k-1)},p)$ is a function $r$
-homogeneous with respect to $p_i$, then

1$^0$ $\displaystyle\frac{\partial H}{\partial p_i}$ is $r-1$-homogeneous; \vspace{3mm}\\

2$^0$ $\displaystyle\frac{\partial ^2H}{\partial p_i\partial p_j}$ is $r-2$
-homogeneous.

A $q$-form $\omega \in \Lambda ^q(\widetilde{T^{*k}M})$ is $r$-homogeneous
with respect to momenta $p_i$ if

\[
\omega \circ h_a^{\prime *}=a^r\omega ,\ \forall a\in \mathbf{R}^{+}.
\]

It follows

\begin{teo}
A $q$-form $\omega $ is $r$-homogeneous with respect to $p_i$ if
\begin{equation}
\mathcal{L}_{C^{*}}\omega =r\omega .  \tag{4.5.10}
\end{equation}
\end{teo}

Consequently, $dx^i$, $dy^{(1)i}$, ..., $dy^{(k-1)i}$, $dp_i$ are $0$, $0$,
..., $0$, $1$-homogeneous with respect to $p_i$, respectively.

Finally, we determine the degree of homogeneity of the function $%
\{f,g\}_{k=1}.$

\begin{prop}
If $f$ and $g$ are $r$ and $s$-homogeneous functions with respect to $p_i$,
then the function $ \{f,g\}_{k-1} $ is homogeneous of degree $s+r-1$.
\end{prop}

For applications to the geometry of Cartan spaces of order $k$ is
important to have a special notion of homogeneity on the fibres of the
bundle $T^{*k}M$.

More precisely, the homothety
\[
\overline{h}_a:(x,y^{(1)},...,y^{(k-1)},p)\in T^{*k}M\rightarrow
(x,ay^{(1)},...,a^{k-1}y^{(k-1)},a^kp)\in T^{*k}M
\]
is preserved by the transformation of local coordinates (4.1.2) on $T^{*k}M$.

Let $H_{y,p}$ be the group of transformation on $T^{*k}M$:
\[
H_{y,p}=\left\{ \overline{h}_a:(x,y^{(1)},...,y^{(k-1)},p)\rightarrow
(x,ay^{(1)},...,a^{k-1}y^{(k-1)},a^kp)|a\in \mathbf{R}^{+}\right\} .
\]

The orbit of a point $u_0=(x_0,y_0^{(1)},...,y_0^{(k-1)},p^0)$ by $H_{y,p}$
is given by
\[
x^i=x_0^i,\ y^{(1)i}=ay_0^{(1)i},\ ...,\ y^{(k-1)}=a^{k-1}y_0^{(k-1)i},\
p_i=a^kp_i^0,\ \forall a\in \mathbf{R}^{+}.
\]

The tangent vector at the point $u_0=\overline{h}_1(u_0)$ is the vector
field $\stackrel{k-1}{\Gamma }+kC^{*}$ at the point $u_0$:
\[
\stackrel{k-1}{\Gamma }(u_0)+kC^{*}(u_0)=y_0^{(1)i}\displaystyle\frac
\partial {\partial y^{(1)i}}|_{u_0}+\cdots +(k-1)y_0^{(k-1)i}\displaystyle
\frac \partial {\partial y^{(k-1)i}}|_{u_0}+kp_i^0\stackrel{\cdot }{\partial
^i}|_{u_0}.
\]

\begin{defi}
A function $H:T^{*k}M\rightarrow \mathbf{R}$ differentiable on $\widetilde{%
T^{*k}M}$ and continuous on the null section of the projection $\pi ^{*k}$
is called \textit{homogeneous of degree }$r\in \mathbf{Z}$\textit{\ on the
fibres of the bundle }$T^{*k}M$ if
\begin{equation}
H\circ \overline{h}_a=a^rH,\ \forall a\in \mathbf{R}^{+}.  \tag{4.5.11}
\end{equation}
\end{defi}

Applying the usual methods it follows:

\begin{teo}
A function $H$ on $T^{*k}M$, differentiable on $\widetilde{T^{*k}M}$ and
continuous on the null section is $r$-homogeneous on the fibres of $T^{*k}M$
if and only if
\begin{equation}
\mathcal{L}_{\stackrel{k-1}{\Gamma }+kC^{*}}H=rH.  \tag{4.5.12}
\end{equation}
\end{teo}

If we expand (4.5.12), we can write it in the form
\begin{equation}
y^{(1)i}\displaystyle\frac{\partial H}{\partial y^{(1)i}}+\cdots
+(k-1)y^{(k-1)i}\displaystyle\frac{\partial H}{\partial
y^{(k-1)i}}+kp_i \stackrel{\cdot }{\partial ^i}H=rH.
\tag{4.5.12a}
\end{equation}

A vector field $X$ on $\widetilde{T^{*k}M}$ is $r$-homogeneous on the fibres
of $T^{*k}M$ if

\begin{equation}
X\circ \overline{h}_a=a^{r-1}\overline{h}_a^{*}\circ X,\ \forall a\in
\mathbf{R}^{+}.  \tag{4.5.13}
\end{equation}

It follows

\begin{teo}
A vector field $X\in \mathcal{X}(\widetilde{T^{*k}M})$ is $r$-homogeneous on
the fibres of $T^{*k}M$ if and only if
\begin{equation}
\mathcal{L}_{\stackrel{k-1}{\Gamma }+kC^{*}}X=(r-1)X.  \tag{4.5.14}
\end{equation}
\end{teo}

Of course, (4.5.14) can be given in the form

\begin{equation}
\left[ \stackrel{k-1}{\Gamma },X\right] +k\left[ C^{*},X\right] =(r-1)X.
\tag{4.5.14a}
\end{equation}

\begin{cor}
1$^0$ The vector fields $\displaystyle\frac \partial {\partial x^i}$, $%
\displaystyle\frac \partial {\partial y^{(1)i}}$, ...,$\displaystyle\frac
\partial {\partial y^{(k-1)i}}$ and $\stackrel{\cdot }{\partial ^i}=%
\displaystyle\frac \partial {\partial p_i}$ are $1$, $0$, ..., $2-k$, $1-k$
homogeneous on the fibres of $T^{*k}M$, respectively.

2$^0$ If $H$ $\in \mathcal{F}(\widetilde{T^{*k}M})$ is $s$-homogeneous and $%
X\in \mathcal{X}(\widetilde{T^{*k}M})$ is $r$-homogeneous on the fibres of $%
T^{*k}M$ then

\quad a. $HX$ is $r+s$-homogeneous;

\quad b. $XH$ is $r+s-1$-homogeneous.
\end{cor}

A $q$-form $\omega \in \Lambda ^q(\widetilde{T^{*k}M})$ is $s$-homogeneous
on the fibres of $T^{*k}M$ if
\[
\omega \circ \overline{h}_a^{*}=a^s\omega ,\ \forall a\in \mathbf{R}^{+}.
\]

The following theorem holds:

\begin{teo}
A $q$-form $\omega \in \Lambda ^q(\widetilde{T^{*k}M})$ is $s$-homogeneous
on the fibres of $T^{*k}M$ if and only if
\begin{equation}
\mathcal{L}_{\stackrel{k-1}{\Gamma }+kC^{*}}\omega =s\omega .  \tag{4.5.15}
\end{equation}
\end{teo}

\begin{cor}
The $1$-forms $dx^i$, $dy^{(1)i}$, ..., $dy^{(k-1)i}$, $dp_i$ are $0$, $1$,
..., $k-1$, $k$-homogeneous on the fibres of $T^{*k}M$.
\end{cor}

We will apply these results in the study of the homogeneity on the fibres
of $T^{*k}M$ of $1$-forms $d_0H$, ..., $d_{k-2}H$ and, of course
of the functions $\{f,g\}_{k=1}.$

\chapter{The Variational Problem for the Hamiltonians of Order $k$}

\markboth{\it{THE GEOMETRY OF HIGHER-ORDER HAMILTON SPACES\ \ \ \ \ }}{\it{The Variational Problem for the Hamiltonians of Order} $k$}

The theory of higher order Hamiltonian systems and its applications in
Analytical Mechanics are consistent only if we study the variational problem
for the differentiable Hamiltonians of order $k$, $H(x,y^{(1)},...,y^{(k-1)},p)$, [96, 98].

In this case the integral of action of $H$ must be defined along
curve $c$ on the cotangent manifold $T^{*}M$ by
\[
I(c)=\int_0^1[p_i\displaystyle\frac{dx^i}{dt}-\displaystyle\frac 12H(x, %
\displaystyle\frac{dx}{dt},...,\displaystyle\frac 1{\left( k-1\right) !} %
\displaystyle\frac{d^{k-1}x}{dt^{k-1}},p)]dt.
\]

A local variation of $c$ is a curve
$\widetilde{c}(\varepsilon_1,\varepsilon _2)$ which depend on a
vector field $V^i$ and a covector field $\eta _i$. The integral of
action $I($ $\widetilde{c}(\varepsilon_1,\varepsilon_2))$ depends
on two parameters $\varepsilon _1,\varepsilon _2$. In order for
the functional $I(c)$ to be an extremal value of the functionals
$I($ $\widetilde{c}(\varepsilon_1,\varepsilon_2))$ it is necessary
that
\[
\displaystyle\frac{\partial
I(\widetilde{c}(\varepsilon_1,\varepsilon_2))}{\partial\varepsilon_\alpha}
\Bigg{|}_{\varepsilon_1=\varepsilon_2=0}=0,\ \  (\alpha =1,2)
\]

These conditions allow to determine the Hamilton-Jacoby equations (5.1.17).

Introducing the higher order energies of $H,$ $\mathcal{E}^{k-1}(H),..., \mathcal{E}^1(H)$, a law of conservation of the energy $\mathcal{E}^{k-1}(H)$
is proved and a N\"other type theorem is formulated. This theory is valid in the case when the order $k$ is greater then 1.

\section{The Hamilton-Jacobi Equations}

A function $H:T^{*k}M\rightarrow \mathbf{R}$ differentiable on $\widetilde{
T^{*k}M}$ and continue on the nul section is called a differentiable
Hamiltonian of order $k$. It depends on the variables $%
(x^i,y^{(1)i},...,y^{(k-1)i},p_i)$. So, it will be denoted by \newline
$H(x,y^{(1)},...,y^{(k-1)},p).$

Let us consider a curve
\[
c:t\in [0,1]\rightarrow (x^i(t),p_i(t))\in \widetilde{T^{*}M},
\]
having the image in a local chart of the manifold $\widetilde{T^{*}M}.$
The curve $c$ can be analytically given by the equations:
\begin{equation}
x^i=x^i(t),\ p_i=p_i(t),\ t\in [0,1] . \tag{5.1.1}
\end{equation}

The extension $\stackrel{\vee}{c}$ to the dual bundle $T^{*k}M$ is well
determined. The extension $\stackrel{\vee}{c}$ is given by the equations
\begin{equation}
\begin{array}{l}
x^i=x^i(t), \\
\\
y^{(\alpha )i}(t)=\displaystyle\frac 1{\alpha !}\displaystyle\frac{d^\alpha
x^i}{dt^\alpha }(t),(\alpha =1,...,k-1) ,\\
\\
p_i=p_i(t),\ t\in [0,1].
\end{array}
\tag{5.1.2}
\end{equation}

Also, we consider a vector field $V^i(t)$ and a covector field $\eta _i(t)$
along curve $c$, having the properties:
\begin{equation}
\begin{array}{c}
V^i(0)=V^i(1)=0,\eta _i(0)=\eta _i(1)=0, \\
\\
\displaystyle\frac{d^\alpha V^i}{dt^\alpha }(0)=\displaystyle\frac{d^\alpha
V^i}{dt^\alpha }(1)=0,(\alpha =1,...,k-2).
\end{array}
\tag{5.1.3}
\end{equation}

The variation $\overline{c}(\varepsilon _1,\varepsilon _2)$ of a curve $c$
determined by the pair $(V^i(t),\eta _i(t))$ is defined by
\begin{equation}
\begin{array}{l}
\overline{x}^i=x^i(t)+\varepsilon _1V^i(t), \\
\\
\overline{p}_i=p_i(t)+\varepsilon _2\eta _i(t),t\in [0,1],
\end{array}
\tag{5.1.4}
\end{equation}
where $\varepsilon _1$ and $\varepsilon _2$ are constants, small in the
absolute values, such that the image of the curve $\widetilde{c}(\varepsilon _1,\varepsilon _2)$
belongs to the same domain of chart on $\widetilde{T^{*}M}$ as the image of curve
$c$. The extension of $\overline{c}(\varepsilon _1,\varepsilon _2)$ is the
curve $\stackrel{\stackrel{}{\vee}}{\overline{c}}(\varepsilon
_1,\varepsilon _2)$ given by the equations:
\begin{equation}
\begin{array}{l}
\overline{x}^i=x^i(t)+\varepsilon _1V^i(t), \\
\\
\overline{y}^{(\alpha )i}=\displaystyle\frac 1{\alpha !}(\displaystyle\frac{
d^\alpha x^i}{dt^\alpha }+\varepsilon _1\displaystyle\frac{d^\alpha V^i}{
dt^\alpha }),\ (\alpha =1,...,k-1), \\
\\
\overline{p}_i=p_i(t)+\varepsilon _2\eta _i(t),\ t\in [0,1].
\end{array}
\tag{5.1.5}
\end{equation}

The integral of action for the Hamiltonian
$H(x,y^{(1)},...,y^{(k-1)},p)$ along curve $c$ is defined, like an
extension of the classical form, by
\begin{equation}
I(c)=\int_0^1[p_i(t)\displaystyle\frac{dx^i}{dt}(t)-\displaystyle\frac
12H(x(t),\displaystyle\frac{dx}{dt}(t),...,\displaystyle\frac 1{\left(
k-1\right) !}\displaystyle\frac{d^{k-1}x}{dt^{k-1}}(t),p(t))]dt.  \tag{5.1.6}
\end{equation}

Evidently, $p_i\displaystyle\frac{dx^i}{dt}-\displaystyle\frac 12
H(x, \displaystyle\frac{dx}{dt},...,\displaystyle\frac 1{\left( k-1\right) !}\displaystyle\frac{d^{k-1}x}{dt^{k-1}},p)$
is a differentiable Hamiltonian on the curve $c$.

The integral of action $I(\overline{c}(\varepsilon _1,\varepsilon _2))$ is:
\begin{equation}
\begin{array}{c}
I(\overline{c}(\varepsilon _1,\varepsilon _2))=\displaystyle\int_0^1[(p+\varepsilon
_2\eta )(\displaystyle\frac{dx}{dt}+\varepsilon _1\displaystyle\frac{dV}{dt}
)-\displaystyle\frac 12H(x+\varepsilon _1V,\displaystyle\frac{dx}{dt}
+\varepsilon _{11}\displaystyle\frac{dV}{dt}),... \\
\\
...\displaystyle\frac 1{\left( k-1\right) !}(\displaystyle\frac{d^{k-1}x}{
dt^{k-1}}+\varepsilon _1\displaystyle\frac{d^{k-1}V}{dt^{k-1}}
),p+\varepsilon _2\eta ]dt .
\end{array}
\tag{5.1.7}
\end{equation}

The necessary conditions in order that $I(c)$ is an extremal value of \newline $I(
\overline{c}(\varepsilon _1,\varepsilon _2))$ are:
\begin{equation}
\displaystyle\frac{\partial I(\overline{c}(\varepsilon _1,\varepsilon _2))}{
\partial \varepsilon _1}\left| _{\varepsilon _1=\varepsilon _2=0}\right. =0, %
\displaystyle\frac{\partial I(\overline{c}(\varepsilon _1,\varepsilon _2))}{
\partial \varepsilon _2}\left| _{\varepsilon _1=\varepsilon _2=0}\right. =0.
\tag{5.1.8}
\end{equation}

In our conditions of differentiability, using the equality (5.1.7),
we get the equations:
\begin{equation}
\int_0^1[p_i(t)\displaystyle\frac{dV^i}{dt}(t)-\displaystyle\frac 12( %
\displaystyle\frac{\partial H}{\partial x^i}V^i+\displaystyle\frac{\partial
H }{\partial y^{\left( 1\right) i}}\displaystyle\frac{dV^i}{dt}+...+ %
\displaystyle\frac 1{\left( k-1\right) !}\displaystyle\frac{\partial H}{
\partial y^{(k-1)i}}\displaystyle\frac{d^{k-1}V^i}{dt^{k-1}})]dt=0  \tag{5.1.9}
\end{equation}
and
\begin{equation}
\int_0^1[\displaystyle\frac{dx^i}{dt}-\displaystyle\frac 12\displaystyle
\frac{\partial H}{\partial p_i}]\eta _idt=0 . \tag{5.1.10}
\end{equation}
So, we obtain:

\begin{teo}
The necessary conditions for $I(c)$ to be an extremal value of the functional $%
I(\overline{c}(\varepsilon _1,\varepsilon _2)$ are given by the equations
(5.1.9) and (5.1.10).
\end{teo}

The previous equations imply the Hamilton-Jacobi equations of the
Hamiltonian $H$. To prove this we need to introduce some new notions.

Consider the following main invariants, [96]:
\begin{equation}
I^1(H)=\mathcal{L}_{\stackrel{1}{\Gamma }}H,\ I^2(H)=\mathcal{L}_{%
\stackrel{2}{\Gamma }}H,...,I^{k-1}(H)=\mathcal{L}_{\stackrel{k-1}{
\Gamma }}(H),  \tag{5.1.11}
\end{equation}
in which $\mathcal{L}$ is the Lie operator of derivation and
\begin{equation}
\begin{array}{l}
\stackrel{1}{\Gamma }=y^{\left( 1\right) i}\displaystyle\frac \partial
{\partial y^{\left( k-1\right) i}}, \\
\\
\stackrel{2}{\Gamma }=y^{\left( 1\right) i}\displaystyle\frac \partial
{\partial y^{\left( k-2\right) i}}+2y^{\left( 2\right) i}\displaystyle\frac
\partial {\partial y^{\left( k-1\right) i}}, \\
\\
\stackrel{k-1}{\Gamma }=y^{\left( 1\right) i}\displaystyle\frac \partial
{\partial y^{\left( 1\right) i}}+2y^{\left( 2\right) i}\displaystyle\frac
\partial {\partial y^{\left( 2\right) i}}+...+(k-1)y^{(k-1)i}\displaystyle %
\frac \partial {\partial y^{(k-1)i}}
\end{array}
\tag{5.1.12}
\end{equation}
are the Liouville vector fields on the manifold $T^{*k}M\,.$ Also, we define
the invariants:
\begin{equation}
\begin{array}{l}
I_V^1(H)=V^i\displaystyle\frac{\partial H}{\partial y^{(k-1)i}}, \\
\\
I_V^2(H)=V^i\displaystyle\frac{\partial H}{\partial y^{(k-2)i}}+%
\displaystyle \frac{dV^i}{dt}\displaystyle\frac{\partial H}{\partial
y^{(k-1)i}}, \\
..................................................... \\
I_V^{k-1}(H)=V^i\displaystyle\frac{\partial H}{\partial y^{(1)i}}+ %
\displaystyle\frac{dV^i}{dt}\displaystyle\frac{\partial H}{\partial y^{(2)i}}
+...+\displaystyle\frac 1{\left( k-2\right) !}\displaystyle\frac{d^{k-2}V^i}{
dt^{k-2}}\displaystyle\frac{\partial H}{\partial y^{(k-1)i}}).
\end{array}
\tag{5.1.13}
\end{equation}

For $V^i=\displaystyle\frac{dx^i}{dt},$ the invariants (5.1.13)
are the same with the invariants $I^1(H),...,I^{(k-1)}(H)$ along
curve $c.$ An important notation is as follows:

\begin{equation}
\stackrel{\circ }{E}_i(H)=\displaystyle\frac{dp_i}{dt}+\displaystyle\frac
12[ \displaystyle\frac{\partial H}{\partial x^i}-\displaystyle\frac d{dt} %
\displaystyle\frac{\partial H}{\partial y^{\left( 1\right) i}}+...+\left(
-1\right) ^{k-1}\displaystyle\frac 1{\left( k-1\right) !}\displaystyle\frac{
d^{k-1}}{dt^{k-1}}\displaystyle\frac{\partial H}{\partial y^{\left(
k-1\right) i}}]  \tag{5.1.14}
\end{equation}

Later we prove that $\stackrel{\circ }{E}_i(H)$ is a $d$-covector
field along curve $c.$

By a straighforward calculus we can prove:

\begin{lem}
The following identity holds:

\begin{equation}
\begin{array}{c}
p_i\displaystyle\frac{dV^i}{dt}-\displaystyle\frac 12[\displaystyle\frac{%
\partial H}{\partial x^i}V^i+\displaystyle\frac{\partial H}{\partial
y^{\left( 1\right) i}}\displaystyle\frac{dV^i}{dt}+...+\displaystyle\frac
1{\left( k-1\right) !}\displaystyle\frac{\partial H}{\partial y^{\left(
k-1\right) i}}\displaystyle\frac{d^{k-i}V^i}{dt^{k-1}}]= \\
\\
=-\stackrel{\circ }{E}_i(H)V^i+\displaystyle\frac d{dt}(p_iV^i)-%
\displaystyle \frac 12\displaystyle\frac d{dt}[I_V^{k-1}(H)-\displaystyle%
\frac 1{2!}\displaystyle\frac d{dt}I_V^{k-2}(H)+... \\
\\
+\left( -1\right) ^{k-2}\displaystyle\frac 1{\left( k-1\right) !}%
\displaystyle\frac{d^{k-2}}{dt^{k-2}}I_V^1(H)].
\end{array}
\tag{5.1.15}
\end{equation}
\end{lem}

Using the previous Lemma, we can prove:

\begin{teo}
The equations (5.1.9) and (5.1.10) are equivalent to the equations

\begin{equation}
\int_0^1\stackrel{\circ }{E}_i(H)V^idt=0;\ \int_0^1[\displaystyle\frac{dx^i}{%
dt}-\displaystyle\frac 12\displaystyle\frac{\partial H}{\partial p_i}]\eta
_idt=0.  \tag{5.1.16}
\end{equation}
\end{teo}

\textbf{Proof.} By means of the identity (5.1.15) the equations (5.1.9) can be
written:\\
$$
\int_0^1\{-\stackrel{\circ }{E}_i(H)V^i+\displaystyle\frac d{dt}[p_iV^i- %
\displaystyle\frac 12(I_V^{k-1}(H)-\displaystyle\frac 1{2!}\displaystyle
\frac d{dt}I_V^{k-2}(H)+...+
$$
$ +\left( -1\right) ^{k-2}\displaystyle\frac
1{\left( k-1\right) !}\displaystyle\frac{d^{k-2}}{dt^{k-2}}I_V^1)]\}dt=0.$

Integrating and taking into account the conditions (5.1.3) and the expression
of the invariants (5.1.13) we obtain the announced result. q.e.d.

Now, the equations (5.1.16) in which $V^i$ and $\eta _i$ are arbitrary lead to
the following Hamilton-Jacoby equations:

\begin{teo}
In order for the integral of action $I(c),$ (5.1.6) to be an extremal value for
the functionals $I(\overline{c}(\varepsilon _1,\varepsilon _2))$, (5.1.7) it is
necessary that the curve $c$ to satisfy the following Hamilton-Jacobi
equations
\end{teo}
\begin{equation}
\begin{array}{l}
\displaystyle\frac{dx^i}{dt}=\displaystyle\frac 12\displaystyle\frac{%
\partial H}{\partial p_i}, \\
\\
\displaystyle\frac{dp_i}{dt}=-\displaystyle\frac 12[\displaystyle\frac{%
\partial H}{\partial x^i}-\displaystyle\frac d{dt}\displaystyle\frac{%
\partial H}{\partial y^{\left( 1\right) i}}+...+\left( -1\right) ^{\left(
k-1\right) }\displaystyle\frac 1{\left( k-1\right) !}\displaystyle\frac{%
d^{k-1}}{dt^{k-1}}\displaystyle\frac{\partial H}{\partial y^{\left(
k-1\right) i}}],
\end{array}
\tag{5.1.17}
\end{equation}
where
\begin{equation}
y^{\left( 1\right) i}=\displaystyle\frac{dx^i}{dt},...,y^{\left( k-1\right)
i}=\displaystyle\frac 1{\left( k-1\right) !}\displaystyle\frac{d^{k-1}x^i}{%
dt^{k-1}}.  \tag{5.1.17a}
\end{equation}
Evidently, the second equation (5.1.17) is equivalent to $\stackrel{\circ}{E}_i(H)=0.$

Another important property is expressed in the following theorem:

\begin{teo}
$\stackrel{\circ }{E}_i(H)$ is a covector field.
\end{teo}

This result can be proved by a direct calculation. Another way is as
follows.
With respect to a change of local coordinates on $\widetilde{T^{*k}M}$ we
have
\[
\int_0^1\stackrel{\circ }{\widetilde{E}}_i(\widetilde{H})\widetilde{V}
^idt-\int_0^1\stackrel{\circ }{E}_i(H)V^idt=\int_0^1[\stackrel{\circ }{%
\widetilde{E}}_i(\widetilde{H})\displaystyle\frac{\partial \widetilde{x}^i}{
\partial x^j}-\stackrel{\circ }{E}_j(H)]V^jdt=0
\]
Since $V^i$ is an arbitrary vector field we obtain
\[
\stackrel{\circ }{\widetilde{E}}_i(\widetilde{H})\displaystyle\frac{\partial
\widetilde{x}^i}{\partial x^j}=\stackrel{\circ }{E}_j(H).
\]

The previous property shows that the equation $\stackrel{\circ }{E}_i(H)=0$
has a geometrical meaning.

\section{Zermelo Conditions}

The integral of action (5.1.6) is defined for the parametrized curves $c:t\in
[0,1]\rightarrow (x^i(t),p_i(t))\in \widetilde{T^{*}M}$. The problem is when
it does not depend on the parametrization of curve $c$.

Consider a differentiable diffeomorphism $\widetilde{t}=\widetilde{t}
(t),t\in [0,1]$ which defines a new parametrization of the curve $c$. If $a=
\widetilde{t}(0),b=\widetilde{t}(1),$ then $c$ will be represented by

\[
c:\widetilde{t}\in [a,b]\rightarrow
(\widetilde{x}^i(t),\widetilde{p}_i(\widetilde{t}))\in
\widetilde{T^{*}M}.
\]

In order that the integral of action $I(c)$ does not depend on the
parametrization of curve $c$ is necessary that:
\begin{equation}
\begin{array}{c}
\{\widetilde{p}_i\displaystyle\frac{d\widetilde{x}^i}{d\widetilde{t}}-
\frac{1}{2}\widetilde{H}(\widetilde{x}, \displaystyle\frac{d\widetilde{x}}{d\widetilde{t}
},...,\displaystyle\frac 1{\left( k-1\right) !}\displaystyle\frac{d^{k-1}
\widetilde{x}}{d\widetilde{t}^{k-1}},\widetilde{p})\}\displaystyle\frac{d

\widetilde{t}}{dt}= \\
\\
=p_i\displaystyle\frac{dx^i}{d\widetilde{t}}\displaystyle\frac{d\widetilde{t}
}{dt}-\frac{1}{2} H[x, \displaystyle\frac{dx}{d\widetilde{t}}\displaystyle\frac{d
\widetilde{t}}{dt}, \displaystyle\frac 1{2!}\displaystyle\frac d{dt}( %
\displaystyle\frac{dx}{d\widetilde{t}}\displaystyle\frac{d\widetilde{t}}{dt}
),.... \\
\\
...,\displaystyle\frac 1{\left( k-1\right) !}\displaystyle\frac{d^{k-2}}{
dt^{k-2}}(\displaystyle\frac{dx}{d\widetilde{t}}\displaystyle\frac{d
\widetilde{t}}{dt}),p]. \\
\end{array}
\tag{5.2.1}
\end{equation}

The previous equality holds for any diffeomorpfism $\widetilde{t}=\widetilde{
t}(t).$

If we take the derivatives of (5.2.1) with respect to $\displaystyle\frac{d\widetilde{
t}}{dt}$ and take $\widetilde{t}=t$ we obtain
\begin{equation}
H=y^{\left( 1\right) i}\displaystyle\frac{\partial H}{\partial y^{\left(
1\right) i}}+...+(k-1)y^{\left( k-1\right) i}\displaystyle\frac{\partial H}{
\partial y^{\left( k-1\right) i}}.  \tag{5.2.2}
\end{equation}

If we take again the derivative of (5.2.1) with respect to $\displaystyle\frac{d^2\widetilde{t}}{
dt^2}$ and consider $\widetilde{t}=t$ we have
\begin{equation}
0=y^{\left( 1\right) i}\displaystyle\frac{\partial H}{\partial y^{\left(
2\right) i}}+...+(k-2)y^{\left( k-2\right) i}\displaystyle\frac{\partial H}{
\partial y^{\left( k-2\right) i}} . \tag{5.2.2a}
\end{equation}
And so on.

Therefore we have:

\begin{teo}
The necessary conditions that the integral of action $I(c)$, (5.1.6) does not
depend on the parametrization of the curve $c$ are the following ones:

\begin{equation}
I^{k-1}(H)=H, I^{k-2}(H)=0, ..., I^1(H)=0 . \tag{5.2.3}
\end{equation}
\end{teo}

Indeed, the equation (5.2.2) can be written as $H=\mathcal{L}_{\stackrel{1}{
\Gamma }}H=I^{k-1}(H\dot )$. The equation (5.2.2') is expressed by $0=\mathcal{%
\ L}_{\stackrel{2}{\Gamma }}H=I^{(k-2)}(H),$ etc. Q.E.D.

The equations (5.2.3) will be called the Zermelo conditions. They were introduced
for the higher order Lagrangians by Kazuo Kondo [70] and are fundamental for
the definition of the Kawaguchi spaces [63].

These conditions are very restrictive for the Hamiltonians $H$. Indeed, let
us consider the Hessian of $H,$ with respect to momenta $p_i$. It has the
elements:
\begin{equation}
g^{ij}=\displaystyle\frac 12\displaystyle\frac{\partial ^2H}{\partial
p_i\partial p_j} . \tag{5.2.4}
\end{equation}

In next chapter we shall prove that $g^{ij}$ is a tensor field and its is
called the fundamental tensor field of the Hamiltonian $H$.

Now is not difficult to prove the following result:

\begin{teo}
If the Hamiltonian $H$ satisfies the Zermelo conditions (5.2.3), then its
fundamental tensor $g^{ij}$ has the properties:
\begin{equation}
\mathcal{L}_{\stackrel{1}{\Gamma }}g^{ij}=...=\mathcal{L}_{\stackrel{k-2}{%
\Gamma }}g^{ij}=0,\mathcal{L}_{\stackrel{k-1}{\Gamma }}g^{ij}=g^{ij}.
\tag{5.2.5}
\end{equation}
\end{teo}

\textbf{Proof}. Remarking that the invariants $I^1(H)=\mathcal{L}_{\stackrel{
1}{\Gamma }}H,...,I^{(k-1)}=\mathcal{L}_{\stackrel{k-1}{\Gamma }}H$ have the
properties:
\[
\stackrel{\cdot }{\partial ^i}\stackrel{\cdot }{\partial ^i}
I^{(k-1)}(H)=I^{\left( k-1\right) }(\stackrel{\cdot }{\partial ^i}\stackrel{
\cdot }{\partial ^i}H),...,\stackrel{\cdot }{\partial ^i}\stackrel{\cdot }{
\partial ^i}I^1(H)=I^1(\stackrel{\cdot }{\partial ^i}\stackrel{\cdot }{
\partial ^i}H)
\]
and using the Zermelo conditions (5.2.3), it follows the equations (5.2.5).

The first conditions (5.2.3), $I^{k-1}(H)=H$ and theorem 4.5.1 allow to prove:

\begin{teo}
A necessary condition that the Zermelo conditions (5.2.3) be verified is that
the Hamiltonian $H$ is $1$-homogeneous with respect to the variables $%
y^{\left( 1\right) i},...,y^{\left( k-1\right) i}.$
\end{teo}

\begin{cor}
If the differentiable Hamiltonian $H$ is not $1$-homogeneous with respect to
$y^{\left( 1\right) i},...,y^{\left( k-1\right) i}$ then the Zermelo
conditions (5.2.3) are not verified.
\end{cor}

\begin{cor}
If the differentiable Hamiltonian $H$ is $1$-homogeneous with respect to $%
y^{\left( 1\right) i},...,y^{\left( k-1\right) i}$ and $2$-homogeneous with
respect to $p_i$, then $H$ is $2k+1$-homogeneous on the fibres of $T^{*k}M$.
\end{cor}

Indeed, $\mathcal{L}_{\stackrel{k-1}{\Gamma }+kC^{*}}H=(1+2k)H.$

\begin{cor}
If the differentiable Hamiltonian $H$ has the properties:

a. It satisfies the Zermelo conditions

b. $H$ is $mk+1$-homogeneous on the fibres of $T^{*k}M$, ($m\in Z$), then $H$
is $m$-homogeneous with respect to momenta $p_i$.
\end{cor}

\begin{rem}
The Hamiltonian $H=p_iy^{\left( 1\right) i}$ satisfies the Zermelo conditions
\end{rem}

\section{ Higher Order Energies. Conservation of Energy $%
\mathcal{E}^{k-1}(H)$}

For a differentiable Hamiltonian $H$ the following invariants are important, [98]:

\begin{equation}
\begin{array}{l}
\mathcal{E}^{k-1}(H)= \\
=I^{k-1}(H)-\displaystyle\frac 1{2!}\displaystyle\frac
d{dt}I^{k-2}(H)+...+\left( -1\right) ^{k-2}\displaystyle\frac 1{\left(
k-1\right) !}\displaystyle\frac{d^{k-2}}{dt^{k-2}}I^1(H)-H, \\
\\
\mathcal{E}^{k-2}(H)= \\
= -\displaystyle\frac 1{2!}\displaystyle\frac
d{dt}I^{k-2}(H)+\displaystyle\frac 1{3!}\displaystyle\frac
d{dt}I^{k-3}(H)+\_...+\left( -1\right) ^{k-3}\displaystyle\frac 1{\left(
k-1\right) !}\displaystyle\frac{d^{k-3}}{dt^{k-3}}I^1(H), \\
.......................................... \\
\mathcal{E}^1(H)=\left( -1\right) ^{k-2}\displaystyle\frac 1{\left(
k-1\right) !}I^1(H).
\end{array}
\tag{5.3.1}
\end{equation}

The invariants $\mathcal{E}^{k-1}(H), \mathcal{E}^{k-2}(H),..., \mathcal{E}^{1}(H)$ are called {\it the energies of order} $k-1,k-2,...,1$ of the Hamiltonian $H$,
respectively.

Their expressions justify the invariant character for each. And these
energies are essential for studying the Nother symmetries of $H$.

In order to prove the law of conservation for the energy $\mathcal{E}
^{k-1}(H)$ we need some preliminary considerations.

\begin{lem}
If the differentiable Hamiltonian $H$ has the property
\begin{equation}
\displaystyle\frac{dx^i}{dt}=\displaystyle\frac 12\displaystyle\frac{%
\partial H}{\partial p_i},  \tag{5.3.2}
\end{equation}
\end{lem}
then we have
\begin{equation}
\begin{array}{c}
\displaystyle\frac{dH}{dt}=2\stackrel{\circ }{E}_i(H)\displaystyle\frac{dx^i
}{dt}+\displaystyle\frac d{dt}[I^{k-1}(H)-\displaystyle\frac 1{2!} %
\displaystyle\frac d{dt}I^{k-2}(H)+... \\
\\
...+\left( -1\right) ^{k-2}\displaystyle\frac 1{(k-1)!}\displaystyle\frac{
d^{k-2}}{dt^{k-2}}I^1(H)].
\end{array}
\tag{5.3.3}
\end{equation}

\textbf{Proof}. The condition (5.3.2) implies:

$\displaystyle\frac{dH}{dt}=\displaystyle\frac{\partial H}{\partial x^i} %
\displaystyle\frac{dx^i}{dt}+\displaystyle\frac{\partial H}{\partial p_i} %
\displaystyle\frac{dp_i}{dt}+(\displaystyle\frac{\partial H}{\partial
y^{(1)i}}\displaystyle\frac{dy^{\left( 1\right) i}}{dt}+...+\displaystyle
\frac{\partial H}{\partial y^{(k-1)i}}\displaystyle\frac{dy^{\left(
k-1\right) i}}{dt})=$\vspace{3mm}\\

$=2\displaystyle\frac{dp_i}{dt}\displaystyle\frac{dx^i}{dt}+\{\displaystyle
\frac{\partial H}{\partial x^i}\displaystyle\frac{dx^i}{dt}+\displaystyle
\frac{\partial H}{\partial y^{(1)i}}\displaystyle\frac{dy^{\left( 1\right)
i} }{dt}+...+\displaystyle\frac{\partial H}{\partial y^{(k-1)i}}%
\displaystyle
\frac{dy^{\left( k-1\right) i}}{dt}\}=$ \vspace{3mm}\\

$$
=-2p_i\displaystyle\frac{d^2x^i}{dt^2}+\{\displaystyle\frac{\partial H}{
\partial x^i}\displaystyle\frac{dx^i}{dt}+\displaystyle\frac{\partial H}{
\partial y^{(1)i}}\displaystyle\frac{dy^{\left( 1\right) i}}{dt}+...+
$$
$$
+\displaystyle\frac{\partial H}{\partial y^{(k-1)i}}\displaystyle\frac{
dy^{\left( k-1\right) i}}{dt}\}+2\displaystyle\frac d{dt}(p_i\displaystyle
\frac{dx^i}{dt}).
$$

Using Lemma 4.1.1 for $V^i=\displaystyle\frac{dx^i}{dt},\displaystyle\frac{
dV^i}{dt}=\displaystyle\frac{d^2x^i}{dt^2},...$ the previous equality leads
to the formula (5.3.3). Q.E.D.

Now, we can prove without difficulties

\begin{teo}
If the differentiable Hamiltonian $H$ satisfies the first \newline Hamilton-Jacobi
equations (5.3.2) then the variation of the energy of order $k-1$, $\mathcal{E}%
^{k-1}(H)$ is given by:

\begin{equation}
\displaystyle\frac{d\mathcal{E}^{k-1}(H)}{dt}=-\stackrel{\circ }{E}_i(H)%
\displaystyle\frac{dx^i}{dt} . \tag{5.3.4}
\end{equation}
\end{teo}

Clearly, this formula leads to an interesting result

\begin{teo}
Along every solution curve $c$ of the Hamilton-Jacobi equations (5.1.17)
the energy of order $k-1,\ \mathcal{E}^{k-1}(H)\,$, is conserved.
\end{teo}

Finally, we notice the following theorems:

\begin{teo}
If the differentiable Hamiltonian $H$ verifies the Zermelo conditions (5.2.3),
then the energies of $H$, $\mathcal{E}^{k-1}(H),...,\mathcal{E}^1(H)$ vanish.
\end{teo}

\begin{teo}
If the differentiable Hamiltonian $H$ verifies the Zermelo conditions (5.2.3),
then along any curve $c$ which satisfies the first Hamilton-Jacobi
equations (5.1.17) we have
\end{teo}
\[
\stackrel{\circ }{E}_i(H)\displaystyle\frac{dx^i}{dt}=0.
\]

\section{The Jacobi-Ostrogradski Momenta}

The theory of Jacoby-Ostrogradski momenta of the Lagrangians of order $k$,
briefly presented in the Chapter 2, can be extended to the Hamiltonians of
order $k$.

Indeed, from (5.3.1) we see that the energy of order $k-1,$ $\mathcal{E}^{k-1}(H)$ is a polynomial function of degree one in the
higher order accelerations $\displaystyle\frac{dx^i}{dt},...,\displaystyle
\frac{d^{k-1}x^i}{dt^{k-1}}$. So, we have

\begin{equation}
\mathcal{E}^{k-1}(H)=p_{\left( 1\right) i}\displaystyle\frac{dx^i}{dt}
+...+p_{\left( k-1\right) i}\displaystyle\frac{d^{k-1}x^i}{dt^{k-1}}-H
\tag{5.4.1}
\end{equation}

A straighfoward calculus shows that $p_{\left( 1\right) i},...,p_{\left(
k-1\right) i}$ are given by

\begin{equation}
\begin{array}{l}
p_{\left( 1\right) i}=\displaystyle\frac{\partial H}{\partial y^{\left(
1\right) i}}-\displaystyle\frac 1{2!}\displaystyle\frac d{dt}\displaystyle
\frac{\partial H}{\partial y^{(2)i}}+...+\left( -1\right) ^{k-2}%
\displaystyle \frac 1{\left( k-1\right) !}\displaystyle\frac{d^{k-2}}{%
dt^{k-2}} \displaystyle\frac{\partial H}{\partial y^{\left( k-1\right) i}},
\\
\\
p_{\left( 2\right) i}=\displaystyle\frac 1{2!}\displaystyle\frac{\partial H}{
\partial y^{(2)i}}-\displaystyle\frac 1{3!}\displaystyle\frac d{dt} %
\displaystyle\frac{\partial H}{\partial y^{(3)i}}+...+\left( -1\right)
^{k-3} \displaystyle\frac 1{\left( k-1\right) !}\displaystyle\frac{d^{k-3}}{%
dt^{k-3} }\displaystyle\frac{\partial H}{\partial y^{\left( k-1\right) i}},
\\
...................................... \\
p_{\left( k-1\right) i}=\displaystyle\frac 1{\left( k-1\right) !} %
\displaystyle\frac{\partial H}{\partial y^{\left( k-1\right) i}}.
\end{array}
\tag{5.4.2}
\end{equation}

The entries $p_{\left( 1\right) i},...p_{\left( k-1\right) i}$ are called the
Jacobi-Ostrogradski momenta of the Hamiltonian $H.$

Using the rule of transformation (4.1.4) of the natural basis of the
module $\mathcal{X}(T^{*k}M)$ and the definition (5.4.2) of the Jacobi-Ostrogradski momenta $%
p_{\left( 1\right) i},...,p_{\left( k-1\right) i}\,$we obtain

\begin{prop}
With respect to the transformations of local coordinates (4.1.2) on the
manifold $T^{*k}M$ the Jacobi-Ostrogradski momenta \newline  $p_{\left( 1\right)
i},...,p_{\left( k-1\right) i}$ are transformed as follows:

\begin{equation}
\begin{array}{l}

p_{\left( 1\right) i}=\displaystyle\frac{\partial \widetilde{y}^{\left(
1\right) m}}{\partial y^{\left( 1\right) i}}\widetilde{p}_{\left( 1\right)
m}+...+\displaystyle\frac{\partial \widetilde{y}^{\left( k-1\right) m}}{%
\partial y^{\left( 1\right) i}}\widetilde{p}_{\left( k-1\right) m}, \\
\\
p_{\left( 2\right) i}=\displaystyle\frac{\partial \widetilde{y}^{\left(
2\right) m}}{\partial y^{\left( 2\right) i}}\widetilde{p}_{\left( 2\right)
m}+...+\displaystyle\frac{\partial \widetilde{y}^{\left( k-1\right) m}}{%
\partial y^{\left( 2\right) i}}\widetilde{p}_{\left( k-1\right) m}, \\
...................................... \\
p_{\left( k-1\right) i}=\displaystyle\frac{\partial \widetilde{y}^{\left(
k-1\right) m}}{\partial y^{\left( k-1\right) i}}\widetilde{p}_{\left(
k-1\right) m}.
\end{array}
\tag{5.4.3}
\end{equation}
\end{prop}

Let us consider the following $1$-form fields on $T^{*k}M:$

\begin{equation}
\begin{array}{l}
p_{\left( 1\right) }=p_{\left( 1\right) i}dx^i+p_{\left( 2\right)
i}dy^{\left( 1\right) i}+...+p_{\left( k-1\right) i}dy^{\left( k-2\right) i},
\\
\\
p_{\left( 2\right) }=p_{\left( 2\right) i}dx^i+p_{\left( 3\right)
i}dy^{\left( 1\right) i}+...+p_{\left( k-1\right) i}dy^{\left( k-3\right) i},
\\
...................................... \\
p_{\left( k-1\right) }=p_{\left( k-1\right) i}dx^i.
\end{array}
\tag{5.4.4}
\end{equation}

\begin{prop}
With respect the transformations of coordinates on the manifold $T^{*k}M$ we
have

\begin{equation}
p_{(\alpha )}=\widetilde{p}_{\left( \alpha \right) },\ \ (\alpha =1,...,k-1).
\tag{5.4.5}
\end{equation}
\end{prop}

Indeed, (5.4.3), (5.4.4) and (4.1.6), have as consequence the equalities
(5.4.5).

The relations between the momenta $p_i$ and the Jacobi-Ostrogradski momenta $%
p_{\left( 1\right) i}$ are given by:

\begin{lem}
The following identity holds:
\begin{equation}
\displaystyle\frac 12\displaystyle\frac{dp_{\left( 1\right) i}}{dt}=-%
\stackrel{\circ }{E}_i(H)+(\displaystyle\frac{dp_i}{dt}+\displaystyle\frac 12%
\displaystyle\frac{\partial H}{\partial x^i}) . \tag{5.4.6}
\end{equation}
\end{lem}

Indeed, (5.1.14) and (5.4.2) imply the last equality. Now, we can formulate

\begin{teo}
Along the solution curves of the Hamilton-Jacobi equations (5.1.17), the
following Hamilton-Jacobi-Ostrogradski equations hold:

\begin{equation}
\begin{array}{l}
\displaystyle\frac{\partial \mathcal{E}^{k-1}(H)}{\partial p_{\left( \alpha
\right) i}}=\displaystyle\frac{d^\alpha x^i}{dt^\alpha },\ \ (\alpha =1,...,k-1),
\\
\\
\displaystyle\frac{\partial \mathcal{E}^{k-1}(H)}{\partial x^i}=-%
\displaystyle\frac{dp_{\left( 1\right) i}}{dt}+2\displaystyle\frac{dp_i}{dt},
\\
\\
\displaystyle\frac{\partial \mathcal{E}^{k-1}(H)}{\partial y^{\left( \alpha
\right) i}}=-\alpha !\displaystyle\frac{dp_{(\alpha +1)i}}{dt},\ \ \alpha
=(1,...,k-1).
\end{array}
\tag{5.4.7}
\end{equation}
\end{teo}

\textbf{Proof}. The energy $\mathcal{E}^{k-1}(H),$ (5.4.1) can be written in
the form

\begin{equation}
\mathcal{E}^{k-1}(H)=p_{\left( 1\right) i}y^{\left( 1\right) i}+2!p_{\left(
2\right) i}y^{\left( 2\right) i}+...+\left( k-1\right) !p_{\left( k-1\right)
i}y^{\left( k-1\right) i}-H,  \tag{5.4.8}
\end{equation}
where

\[
y^{\left( \alpha \right) i}=\displaystyle\frac 1{\alpha !}\displaystyle\frac{
d^\alpha x^i}{dt^\alpha },(\alpha =1,...,k-1).
\]

Therefore the equations (5.4.7)$_1$ and (5.4.7)$_3$ hold. In order to prove
(5.4.7) $_2$ we remark that (5.4.8) implies $\displaystyle\frac{\partial
\mathcal{E} ^{k-1}(H)}{\partial x^i}=-\displaystyle\frac{\partial H}{%
\partial x^i}$ and using Lemma 5.4.1 we obtain the equations (5.4.7)$_2$
.Q.E.D.

>From Lemma 5.4.1 we deduce also

\begin{teo}
Along the solution curves of the Hamilton-Jacobi equations (5.1.17) we have

\begin{equation}
\displaystyle\frac{dx^i}{dt}=\displaystyle\frac 12\displaystyle\frac{%
\partial H}{\partial p_i};\displaystyle\frac{dp_i}{dt}=-\displaystyle\frac 12%
\displaystyle\frac{\partial H}{\partial x^i}+\displaystyle\frac 12%
\displaystyle\frac{\partial p_{\left( 1\right) i}}{dt}.  \tag{5.4.9}
\end{equation}
\end{teo}

\section{N\"other Type Theorems}

In this section we will define the notion of symmetry of a differentiable
Hamiltonian $H$ and will prove two N\"other type theorems, using the model of
symmetries from the Lagrangian theory, (Ch.2).

Consider a differentiable Hamiltonian $H_0$ with the properties:

\begin{equation}
\displaystyle\frac{\partial H_0}{\partial y^{\left( k-1\right) i}}=0, %
\displaystyle\frac{\partial H_0}{\partial p_i}=0  \tag{5.5.1}
\end{equation}

\begin{prop}
If $H_0$ is a differentiable Hamiltonian, having the properties (5.5.1), then
the following identities hold:

\begin{equation}
\begin{array}{l}

\displaystyle\frac \partial {\partial p_i}\displaystyle\frac{dH_0}{dt}=%
\displaystyle\frac d{dt}\displaystyle\frac{\partial H_0}{\partial p_i}=0, \\
\\
\displaystyle\frac \partial {\partial x^i}\displaystyle\frac{dH_0}{dt}=%
\displaystyle\frac d{dt}\displaystyle\frac{\partial H_0}{\partial x^i}, \\
\\
\displaystyle\frac \partial {\partial y^{\left( 1\right) i}}\displaystyle
\frac{dH_0}{dt}=\displaystyle\frac d{dt}\displaystyle\frac{\partial H_0}{%
\partial y^{\left( 1\right) i}}+\displaystyle\frac{\partial H_0}{\partial x^i%
}, \\
\\
\displaystyle\frac \partial {\partial y^{\left( 2\right) i}}\displaystyle
\frac{dH_0}{dt}=\displaystyle\frac d{dt}\displaystyle\frac{\partial H_0}{%
\partial y^{\left( 2\right) i}}+2\displaystyle\frac{\partial H_0}{\partial
y^{\left( 1\right) i}}, \\
..................................................... \\
\displaystyle\frac \partial {\partial y^{\left( k-2\right) i}}\displaystyle
\frac{dH_0}{dt}=\displaystyle\frac d{dt}\displaystyle\frac{\partial H_0}{%
\partial y^{\left( k-2\right) i}}+\left( k-2\right) \displaystyle\frac{%
\partial H_0}{\partial y^{\left( k-3\right) i}}, \\
\\
\displaystyle\frac \partial {\partial y^{\left( k-1\right) i}}\displaystyle
\frac{dH_0}{dt}=\left( k-1\right) \displaystyle\frac{\partial H_0}{\partial
y^{\left( k-2\right) i}}.
\end{array}
\tag{5.5.2}
\end{equation}
\end{prop}

Indeed, along curve $c:t\in [0,1]\rightarrow (x^i(t),p_i(t))\in
\widetilde{ T^{*}M}$ we have

\[
\displaystyle\frac{dH_0}{dt}=\displaystyle\frac{\partial H_0}{\partial x^i}
y^{\left( 1\right) i}+2\displaystyle\frac{\partial H_0}{\partial y^{\left(
1\right) i}}y^{\left( 2\right) i}+...+\left( k-1\right) \displaystyle\frac{
\partial H_0}{\partial y^{\left( k-2\right) i}}y^{\left( k-1\right) i}.
\]

By a straightforward calculus we deduce the previous identities.

\begin{prop}
If $H_0$ is a differentiable Hamiltonian which satisfies the equations
(5.5.1), then the following identities hold:

\begin{equation}
\begin{array}{c}
\displaystyle\frac 12\displaystyle\frac \partial {\partial p_i}(H+%
\displaystyle\frac{dH_0}{dt})=\displaystyle\frac 12\displaystyle\frac{%
\partial H}{\partial p_i}, \\
\\
\stackrel{\circ }{E}_i(H+\displaystyle\frac{dH_0}{dt})=\stackrel{\circ }{E}%
_i(H).
\end{array}
\tag{5.5.3}
\end{equation}
\end{prop}

Indeed, taking into account the proposition 5.5.1 and the expression (5.1.14)
of the covector $\stackrel{\circ }{E}_i(H)$ we have (5.5.3)$_1$ and\vspace{3mm}\\

$\stackrel{\circ }{E}_i(H+\displaystyle\frac{dH_0}{dt})=\stackrel{\circ }{E}
_i(H)+\displaystyle\frac 12[\displaystyle\frac \partial {\partial x^i} %
\displaystyle\frac{\partial H_0}{dt}-\displaystyle\frac d{dt}\displaystyle
\frac \partial {\partial y^{\left( 1\right) i}}\displaystyle\frac{dH_0}{dt}+$

$+\displaystyle\frac 1{2!}\displaystyle\frac{d^2}{dt^2}\displaystyle\frac
\partial {\partial y^{\left( 2\right) i}}\displaystyle\frac{dH_0}{dt}
+...+\left( -1\right) ^{k-1}\displaystyle\frac 1{\left( k-1\right) !} %
\displaystyle\frac{d^{k-1}}{dt^{k-1}}\displaystyle\frac \partial {\partial
y^{\left( k-1\right) i}}\displaystyle\frac{dH_0}{dt}]=\stackrel{\circ }{E}
_i(H).$\vspace{3mm}\\

Now, it is easy to prove

\begin{teo}
$H_0$ being an arbitrary differentiable Hamiltonian with the
properties (5.5.1), then the Hamiltonians $H$ and $H+\displaystyle\frac{dH_0}{dt}$ have
the same Hamilton-Jacobi equations:

\begin{equation}
\displaystyle\frac{dx^i}{dt}=\displaystyle\frac
12\displaystyle\frac{
\partial H}{\partial p_i},\ \ \ \ \stackrel{\circ }{E}_i(H)=0  \tag{5.5.4}
\end{equation}
\end{teo}

Indeed, (5.5.4) are the consequence of the equations (5.5.3) Q.E.D.

\begin{cor}
If $H_0$ is an arbitrary differentiable Hamiltonian with the \vspace{3mm}\\
properties $\displaystyle\frac{\partial H_0}{\partial y^{(k-1)i}}=0,\displaystyle\frac{
\partial H_0}{\partial p_i}=0,$ then the integrals of action: \vspace{3mm}\\
\begin{equation}
I\left( c\right) =\int_0^1[p_i\displaystyle\frac{dx^i}{dt}-H(x,\displaystyle
\frac{dx}{dt},...,\displaystyle\frac 1{\left( k-1\right) !}\displaystyle
\frac{d^{k-1}x}{dt^{k-2}},p)]dt  \tag{5.5.5}
\end{equation}
and
\begin{equation}
\begin{array}{l}
I^{\prime }\left( c\right)
=\int_0^1[p_i\displaystyle\frac{dx^i}{dt}- \displaystyle\frac
12H(x,\displaystyle\frac{dx}{dt},...,\displaystyle\frac 1{\left(
k-1\right) !}\displaystyle\frac{d^{k-1}x}{dt^{k-2}},p)]-\\
\ -\displaystyle\frac
12\displaystyle\frac{dH_0}{dt}(x,\displaystyle\frac{dx
}{dt},...,\displaystyle\frac 1{\left(
k-2\right)!}\displaystyle\frac{ d^{k-2}x}{dt^{k-2}})]dt
\end{array} \tag{5.5.6}
\end{equation}
determine the same Hamilton-Jacobi equations (5.5.4).
\end{cor}

Be means of the previous result we can formulate

\begin{defi}
A symmetry of the differentiable Hamiltonian $H(x, y^{(1)}, ...,$ $y^{(k-1)} ,p)$
is a $C^{\infty \,}$-diffeomorphism $\varphi :\widetilde{T^{*}M}\times
R\rightarrow \widetilde{T^{*}M}\times R$ which preserves the variational
principle of the integral of action, expressed in the Corollary 5.5.1.
\end{defi}

We consider the local symmetries only.
Therefore we study the infinitesimal
symmetries defined on an open set $\pi ^{*-1}(U)\times (a,b)$ in the
infinitesimal form:

\begin{equation}
\begin{array}{l}
x^{\prime i}(t^{\prime })=x^i(t)+\varepsilon _1V^i(t), \\
\\
p_i^{\prime }(t^{\prime })=p_i(t)+\varepsilon _2\eta _i(t), \\
\\
t^{\prime }=t+\varepsilon _1\tau _1(t)+\varepsilon _2\tau _2(t),
\end{array}
\tag{5.5.7}
\end{equation}
where $\varepsilon _1,\varepsilon _2$ are real numbers, sufficiently small
in absolute values so that the points $(x,p,t)$ and $(x^{\prime },p^{\prime
},t^{\prime })$ belong to the same set $\pi ^{*-1}(U)\times (a,b),$ where
the curve
\[
c:t\in [0,1]\rightarrow (x^i(t),p_i(t),t)\in \pi ^{*-1}(U)\times (a,b)
\]
is given. By $V^i(t)$ we mean a vector field $V^i(x(t))$ and $\eta
_i(t) $ denotes a covector field $\eta _i(x(t))$, along $c$. The pair $%
(V^i(t),\eta _i(t))$ satisfies (5.1.3) and the differentiable functions $\tau
_1(t),\tau _2(t)$ satisfy the conditions $\tau _\alpha (0)=\tau _\alpha
(1)=0,(\alpha =1,2).$

In the following considerations the terms of higher order in $\varepsilon _{1,}\varepsilon _2,$ will be neglected$.$

The infinitesimal transformation (5.5.7) is a symmetry of the differentiable
Hamiltonian $H(x,y^{(1)},...,y^{(k-1)},p)$ if for any
differentiable Hamiltonian \newline  $H_0(x,y^{(1)},...,y^{(k-2)})$ the following
equations hold:

\begin{equation}
\begin{array}{l}
\lbrack p_i^{\prime }\displaystyle\frac{dx^{\prime i}}{dt^{\prime }}- %
\displaystyle\frac 12H(x^{\prime },\displaystyle\frac{dx^{\prime }}{
dt^{\prime }},...,\displaystyle\frac 1{\left( k-1\right) !}\displaystyle
\frac{d^{k-1}x^{\prime }}{dt^{^{\prime }k-1}},p^{\prime })]dt^{\prime }= \\
\\
=p_i\displaystyle\frac{dx^i}{dt}-\displaystyle\frac 12H(x,\displaystyle\frac{
dx}{dt},...,\displaystyle\frac 1{\left( k-1\right) !}\displaystyle\frac{
d^{k-1}}{dt^{k-1}},p)- \\
\\
-\displaystyle\frac 12\displaystyle\frac{dH_0}{dt}(x,\displaystyle\frac{dx}{
dt},...,\displaystyle\frac 1{\left( k-2\right) !}\displaystyle\frac{d^{k-2}x
}{dt^{k-2}})]dt.
\end{array}
\tag{5.5.8}
\end{equation}

From (5.5.7) we deduce
\begin{equation}
\begin{array}{l}

\displaystyle\frac{dt^{\prime }}{dt}=1+\varepsilon _1\displaystyle\frac{
d\tau _1}{dt}+\varepsilon _2\displaystyle\frac{d\tau _2}{dt}, \\
\\
\displaystyle\frac{dx^{\prime i}}{dt^{\prime }}=\displaystyle\frac{dx^i}{dt}
+\varepsilon _1(\displaystyle\frac{dV^i}{dt}-\varphi _1^i)+\varepsilon
_2\Psi _1^i ,\\
........................................................, \\
\displaystyle\frac{d^{k-1}x^{\prime i}}{dt^{\prime k-1}}=\displaystyle\frac{
d^{k-1}x^i}{dt^{k-1}}+\varepsilon _1(\displaystyle\frac{d^{k-1}V^i}{dt^{k-1}}
-\varphi _{k-1}^i)+\varepsilon _2\Psi _{k-1}^i,
\end{array}
\tag{5.5.9}
\end{equation}
where
\begin{equation}
\begin{array}{l}
\varphi _1^i=\displaystyle\frac{dx^i}{dt}\displaystyle\frac{d\tau _1}{dt},
\vspace{3mm}\\
\varphi _2^i=\left( _1^2\right) \displaystyle\frac{d^2x^i}{dt^2}%
\displaystyle \frac{d\tau _1}{dt}+\left( _2^2\right) \displaystyle\frac{dx^i%
}{dt} \displaystyle\frac{d^2\tau _1}{dt^2} ,\\
......................................................., \\

\varphi _{k-1}^i=\left( _1^{k-1}\right) \displaystyle\frac{d^{k-1}x^i}{
dt^{k-1}}\displaystyle\frac{d\tau _1}{dt}+\left( _2^{k-1}\right) %
\displaystyle\frac{d^{k-2}x^i}{dt^{k-2}}\displaystyle\frac{d^2\tau _1}{dt^2}
+... \vspace{3mm}\\
\qquad +\left( _{k-1}^{k-1}\right) \displaystyle\frac{dx^i}{dt^{}} %
\displaystyle\frac{d^{k-1}\tau _1}{dt^{k-1}}
\end{array}
\tag{5.5.10}
\end{equation}
and similarly, for $\Psi _1^i,...,\Psi _{k-1}^i$:
\begin{equation}
\begin{array}{l}
\Psi _1^i=\displaystyle\frac{dx^i}{dt}\displaystyle\frac{d\tau _2}{dt}, \\
\\
\Psi _2^i=\left( _1^2\right) \displaystyle\frac{d^2x^i}{dt^2}\displaystyle
\frac{d\tau _2}{dt}+\left( _2^2\right) \displaystyle\frac{dx^i}{dt} %
\displaystyle\frac{d^2\tau _2}{dt^2} ,\\
......................................................., \\
\Psi _{k-1}^i=\left( _1^{k-1}\right) \displaystyle\frac{d^{k-1}x^i}{dt^{k-1}}
\displaystyle\frac{d\tau _2}{dt}+\left( _2^{k-1}\right) \displaystyle\frac{
d^{k-2}x^i}{dt^{k-2}}\displaystyle\frac{d^2\tau _2}{dt^2}+... \\
\qquad +\left( _{k-1}^{k-1}\right) \displaystyle\frac{dx^i}{dt^{}} %
\displaystyle\frac{d^{k-1}\tau _2}{dt^{k-1}}.
\end{array}
\tag{5.5.10a}
\end{equation}

The previous formulas are proved starting from (5.5.7) and taking into account
the following expressions:

$\displaystyle\frac{dx^{\prime i}}{dt^{\prime }}=\left( \displaystyle\frac{
dx^i}{dt}+\varepsilon _1\displaystyle\frac{dV^i}{dt}\right) \displaystyle
\frac{dt}{dt^{\prime }}=\left( \displaystyle\frac{dx^i}{dt}+\varepsilon _1 %
\displaystyle\frac{dV^i}{dt}\right) \left( 1-\varepsilon _1\displaystyle
\frac{d\tau _1}{dt}-\varepsilon _2\displaystyle\frac{d\tau _2}{dt}\right) =$\vspace{3mm}\\

$=\displaystyle\frac{dx^i}{dt}+\varepsilon _1\left( \displaystyle\frac{dV^i}{
dt}-\varphi _1^i\right) -\varepsilon _2\Psi _1^i,$etc.

By means of the formulas (5.5.9) the equality (5.5.8) becomes:\\
$\{(p_i+\varepsilon _2\eta _i)[\displaystyle\frac{dx^i}{dt}+\varepsilon
_1\left( \displaystyle\frac{dV^i}{dt}-\varphi _1^i\right) $-$\varepsilon
_2\Psi _1^i]-\displaystyle\frac 12H[x+\varepsilon _1V,\displaystyle\frac{dx}{
dt}+\varepsilon _1\left( \displaystyle\frac{dV}{dt}-\varphi _1\right) $- \vspace{3mm}\\
$ \varepsilon _2\Psi _1,...,\displaystyle\frac 1{\left( k-1\right) !}\left( \displaystyle\frac{
d^{k-1}x}{dt^{k-1}}+\varepsilon _1\left( \displaystyle\frac{d^{k-1}V}{
dt^{k-1}}-\varphi _{k-1}\right) -\varepsilon _2\Psi _{k-1}\right)
,p+\varepsilon _2\eta ]\}\vspace{3mm}\\
\left( 1+\varepsilon _1\displaystyle\frac{d\tau _1}{dt}+\varepsilon _2\displaystyle\frac{d\tau _2}{dt}\right)=
p_i\displaystyle\frac{dx^i}{dt}-\displaystyle\frac 12H(x,\displaystyle
\frac{dx}{dt},...,\displaystyle\frac 1{\left( k-1\right) !}\displaystyle
\frac{d^{k-1}x}{dt^{k-1}},p)-\vspace{3mm}\\
\displaystyle\frac 12\displaystyle\frac{dH_0}{%
dt }(x,\displaystyle\frac{dx}{dt},...,\displaystyle\frac 1{\left( k-2\right)
!} \displaystyle\frac{d^{k-2}x}{dt^{k-2}}).$\vspace{3mm}\\

Using the Taylor expansion with respect to $\varepsilon _1,\varepsilon _2$
and neglecting the terms of higher order in $\varepsilon _1,\varepsilon _2$ we
obtain
\begin{equation}
\begin{array}{l}
p_i\displaystyle\frac{dV}{dt}-\displaystyle\frac 12[\displaystyle\frac{
\partial H}{\partial x^i}V^i+\displaystyle\frac{\partial H}{\partial
y^{\left( 1\right) i}}\displaystyle\frac{dV^i}{dt},...,\displaystyle\frac
1{\left( k-1\right) !}\displaystyle\frac{\partial H}{\partial y^{\left(
k-1\right) i}}\displaystyle\frac{d^{k-1}V^i}{dt^{k-1}}]+ \\
\\
+\displaystyle\frac 12[\displaystyle\frac{\partial H}{\partial y^{\left(
1\right) i}}\varphi _1^i+\displaystyle\frac 1{2!}\displaystyle\frac{\partial
H}{\partial y^{\left( 2\right) i}}\varphi _2^i+...+\displaystyle\frac
1{\left( k-1\right) !}\displaystyle\frac{\partial H}{\partial y^{\left(
k-1\right) i}}\varphi _{k-1}^i]- \\
\\
-\displaystyle\frac 12H\displaystyle\frac{d\tau _1}{dt}=\displaystyle\frac
12 \displaystyle\frac{d\Phi _1}{dt}
\end{array}
\tag{5.5.11}
\end{equation}
and
\begin{equation}
\begin{array}{c}
\left( \displaystyle\frac{dx^i}{dt}-\displaystyle\frac 12\displaystyle\frac{
\partial H}{\partial p_i}\right) \eta _i+\displaystyle\frac 12[\displaystyle
\frac{\partial H}{\partial y^{\left( 1\right) i}}\Psi _1^i+\displaystyle %
\frac 1{2!}\displaystyle\frac{\partial H}{\partial y^{\left( 2\right) i}}
\Psi _2^i+... \\
\\
...+\displaystyle\frac 1{\left( k-1\right) !}\displaystyle\frac{\partial H}{
\partial y^{\left( k-1\right) i}}\Psi _{k-1}^i]-\displaystyle\frac 12H %
\displaystyle\frac{d\tau _2}{dt}=\displaystyle\frac{d\Phi _2}{dt}
\end{array}
\tag{5.5.11a}
\end{equation}
where
\begin{equation}
\varepsilon _1\Phi _1+\varepsilon _2\Phi _2=H_0  \tag{5.5.12}
\end{equation}

As a consequence we have that the functions $\Phi _1$ and $\Phi _2$ are depending only on the variables
$(x^i,y^{\left( 1\right) i},...,y^{\left( k-2\right) i})$ and, of course, these can be arbitrary chosen.

Taking into account of Lemma 5.1.1 and remarking the identity
\[
\begin{array}{c}
\displaystyle\frac{\partial H}{\partial y^{\left( 1\right) i}}\varphi _1^i+ %
\displaystyle\frac 1{2!}\displaystyle\frac{\partial H}{\partial y^{\left(
2\right) i}}\varphi _2^i+...+\displaystyle\frac 1{\left( k-1\right) !} %
\displaystyle\frac{\partial H}{\partial y^{\left( k-1\right) i}}\varphi
_{k-1}^i= \\
\\
=\displaystyle\frac{d\tau _1}{dt}I^{k-1}(H)+\displaystyle\frac 1{2!} %
\displaystyle\frac{d^2\tau _1}{dt^2}I^{k-2}(H)+...+\displaystyle\frac
1{(k-1)!}\displaystyle\frac{d^{k-1}\tau _1}{dt^{k-1}}I^1(H)
\end{array}
\]
and similar identities for $\Psi _\alpha $ and $\tau _2,$ the equations (5.5.11)
and (5.5.11') can be expressed more simple.
\begin{prop}
The equations (5.5.11) and (5.5.11') are echivalent to
\begin{equation}
\begin{array}{l}
-\stackrel{\circ }{E}_i(H)V^i+\displaystyle\frac d{dt}(p_iV^i)-\displaystyle %
\frac 12\displaystyle\frac d{dt}[I_V^{k-1}(H)-\displaystyle\frac 1{2!}%
\displaystyle\frac d{dt}I_V^{k-2}(H)+... \\
\\
...+\left( -1\right) ^{k-2}\displaystyle\frac 1{\left( k-1\right) !}%
\displaystyle\frac{d^{k-2}}{dt^{k-1}}I_V^1(H)]+\\
\displaystyle\frac 12[%
\displaystyle\frac{d\tau _1}{dt}I^{k-1}(H)+\displaystyle\frac 1{2!}%
\displaystyle\frac{d^2\tau _1}{dt^2}I^{k-2}(H)+... \\
\\
...+\displaystyle\frac 1{(k-1)!}\displaystyle\frac{d^{k-1}\tau _1}{dt^{k-1}}%
I^1(H)]-\displaystyle\frac 12H\displaystyle\frac{d\tau _1}{dt}=\displaystyle %
\frac 12\displaystyle\frac{d\Phi _1}{dt}
\end{array}
\tag{5.5.13}
\end{equation}
and

\end{prop}
\begin{equation}
\begin{array}{l}
(\displaystyle\frac{dx^i}{dt}-\displaystyle\frac 12\displaystyle\frac{
\partial H}{\partial p_i})\eta _i+\displaystyle\frac 12[\displaystyle\frac{
d\tau _2}{dt}I^{k-1}(H)+\displaystyle\frac 1{2!}\displaystyle\frac{d^2\tau
_2 }{dt^2}I^{k-2}(H)+...+ \\
\\
...+\displaystyle\frac 1{(k-1)!}\displaystyle\frac{d^{k-1}\tau _1}{dt^{k-1}}
I^1(H)]-\displaystyle\frac 12H\displaystyle\frac{d\tau _2}{dt}=\displaystyle %
\frac 12\displaystyle\frac{d\Phi _2}{dt}.
\end{array}
\tag{5.5.13a}
\end{equation}

Now, we can prove by a straightforward calculus

\begin{lem}
The following identities hold
\begin{equation}
\begin{array}{l}
\displaystyle\frac{d\tau _\alpha }{dt}I^{k-1}(H)+\displaystyle\frac 1{2!}%
\displaystyle\frac{d^2\tau _\alpha }{dt^2}I^{k-2}(H)+...+\displaystyle\frac
1{(k-1)!}\displaystyle\frac{d^{k-1}\tau _\alpha }{dt^{k-1}}I^1(H)-H%
\displaystyle\frac{d\tau _\alpha }{dt}= \\
\\
=-\tau _\alpha \displaystyle\frac d{dt}\mathcal{E}^{k-1}(H)+\displaystyle %
\frac d{dt}[\tau _\alpha \mathcal{E}^{k-1}(H)-\displaystyle\frac{d\tau
_\alpha }{dt}\mathcal{E}^{k-2}(H)+...+\\
+\left( -1\right) ^{k-2}\displaystyle
\frac{d^{k-2}\tau _\alpha }{dt^{k-2}}\mathcal{E}^1(H)], \ \ \ (\alpha =1,2).
\end{array}
\tag{5.5.14}
\end{equation}
\end{lem}
In the previous equality $\mathcal{E}^{k-1}(H),...\mathcal{E}^1(H)$ are the
energies of order $(k-1,...1)$ respectively, expressed in (5.3.1).
\begin{teo}
For any infinitesimal symmetry (5.5.7) the left hand sides of the following
equalities do not depend on the variables $y^{\left( k-1\right)i}$ and
momenta $p_i.$
\begin{equation}
\begin{array}{l}
-\stackrel{\circ }{E}_i(H)V^i+\displaystyle\frac d{dt}(p_iV^i)-\displaystyle %
\frac 12\displaystyle\frac d{dt}[I_V^{k-1}(H)-\displaystyle\frac 1{2!}%
\displaystyle\frac{dI_V^{k-2}(V)}{dt}+... \\
\\
...+\left( -1\right) ^{k-2}\displaystyle\frac 1{\left( k-1\right) !}%
\displaystyle\frac{d^{k-2}I_V^1(H)}{dt^{k-2}}]-\displaystyle\frac 12\tau _1%
\displaystyle\frac{d\mathcal{E}^{k-1}(H)}{dt}+\displaystyle\frac 12%
\displaystyle\frac d{dt}[\tau _1\mathcal{E}^{k-1}(H)-\vspace{3mm}\\
\displaystyle\frac{d\tau _1}{dt}\mathcal{E}^{k-2}(H)+... +(-1)^{k-2}\displaystyle\frac{d^{k-2}\tau _1}{dt^{k-2}}\mathcal{E}^1(H)]=%
\displaystyle\frac 12\displaystyle\frac{d\Phi _1}{dt}. \\
\\
(\displaystyle\frac{dx^i}{dt}-\displaystyle\frac 12\displaystyle\frac{%
\partial H}{\partial p_i})\eta _i-\displaystyle\frac 12\tau _2\displaystyle
\frac{d\mathcal{E}^{k-1}(H)}{dt}+\displaystyle\frac 12\displaystyle\frac
d{dt}[\tau _2\mathcal{E}^{k-1}(H)-\displaystyle\frac{d\tau _2}{dt}\mathcal{E}%
^{k-2}(H)+... \\
\\
+\left( -1\right) ^{k-2}\displaystyle\frac{d^{k-2}\mathcal{E}^1(H)}{dt^{k-2}}%
]=\displaystyle\frac 12\displaystyle\frac{d\Phi _2}{dt}.
\end{array}
\tag{5.5.15}
\end{equation}
\end{teo}
\textbf{Proof}. Indeed, the equalities (5.5.15) result from proposition 5.5.3
and Lemma 5.5.1. But the functions $\Phi _1$ and $\Phi _2$ satisfies the
conditions \vspace{3mm}\\
$\displaystyle\frac{\partial \Phi _\alpha }{\partial y^{(k-1)i}}=0, %
\displaystyle\frac{\partial \Phi _\alpha }{\partial p_i}=0,(\alpha =1,2)$
Q.E.D. \vspace{3mm}

Now, we can prove a N\"other type theorem:
\begin{teo}
For any infinitesimal symmetry (5.5.7) of a differentiable Hamiltonian $H(x,y^{(1)},...,y^{(k-1)},p)$ and for any differentiable functions $\Phi_1,\Phi _2$ \vspace{1mm}\\
with the properties $\displaystyle\frac{\partial \Phi _\alpha
}{\partial y^{(k-1)i}}=\displaystyle\frac{\partial \Phi _\alpha
}{\partial p_i} =0,(\alpha =1,2)$, the following functions
$\mathcal{F}_1^k(H,\Phi _1), \mathcal{F}_2^k(H,\Phi _2)$ are
conserved along the solutions
curves of the Hamilton-Jacobi equations \vspace{3mm}\\
$\displaystyle\frac{dx^i}{dt}=\displaystyle\frac 12\displaystyle\frac{\partial H}{\partial p_i},\stackrel{\circ }{E}_i(H)=0:$

\begin{equation}
\begin{array}{c}
\mathcal{F}_1^k(H,\Phi _1)= \hfill \\
=p_iV^i-\displaystyle\frac 12[I_V^{k-1}(H)-%
\displaystyle\frac 1{2!}\displaystyle\frac{dI_V^{k-2}(H)}{dt}+...+(-1)^{k-2}%
\displaystyle\frac 1{\left( k-1\right) !}\displaystyle\frac{dI_V^{k-2}(H)}{%
dt^{k-2}}]+ \\
\\
+\displaystyle\frac 12[\tau _1\mathcal{E}^{k-1}(H)-\displaystyle\frac{d\tau
_1}{dt}\mathcal{E}^{k-2}(H)+...+\left( -1\right) ^{k-2}\displaystyle\frac{%
d^{k-2}\tau _1}{dt^{k-2}}\mathcal{E}^1(H)-\Phi _1], \\
\\
\mathcal{F}_2^k(H,\Phi _2)=\tau _2\mathcal{E}^{k-1}(H)-\displaystyle\frac{%
d\tau _2}{dt}\mathcal{E}^{k-2}(H)+...+\left( -1\right) ^{k-2}\displaystyle
\frac{d^{k-2}\tau _2}{dt^{k-2}}\mathcal{E}^1(H)-\Phi _2.
\end{array}
\tag{5.5.16}
\end{equation}
\end{teo}

\textbf{Proof}. Taking into account the Theorem 5.5.2 and the expression of the
function $\mathcal{F}_1^k(H,\Phi _1)$ and $\mathcal{F}_2^k(H,\Phi _2)$ it
follows that the following equations hold:
\begin{equation}
\begin{array}{l}
\stackrel{\circ }{E}_i(H)V^i=\displaystyle\frac{d\mathcal{F}_1^k(H,\Phi _1)}{
dt}-\displaystyle\frac 12\tau _1\displaystyle\frac{d\mathcal{E}^{k-1}(H)}{dt}
, \\
(\displaystyle\frac{dx^i}{dt}-\displaystyle\frac 12\displaystyle\frac{
\partial H}{\partial p_i})\eta _i=\displaystyle\frac 12\displaystyle\frac{d
\mathcal{F}_2^k(H,\Phi _2)}{dt}-\displaystyle\frac 12\tau _2\displaystyle
\frac{d\mathcal{E}^{k-1}(H)}{dt}.
\end{array}
\tag{5.5.17}
\end{equation}
According to the Hamilton-Jacobi equations (5.1.17) and the law of
conservation of energy $\mathcal{E}^{k-1}(H),$ the last equations (5.5.17) implies
\[
\displaystyle\frac{d\mathcal{F}_1(H,\Phi _1)}{dt}=0, \displaystyle\frac{d
\mathcal{F}_2(H,\Phi _{2})}{dt}=0.
\]

The functions $\mathcal{F}_1^k(H,\Phi _1)$ and $\mathcal{F}_2^k(H,\Phi _2)$
depend on the invariants $I_V^1(H),$ $...,I_V^{k-1}(H)$, the energies of order $1,...,k-1,$ $\mathcal{E}^1(H),...,\mathcal{E}^{k-1}(H)$ and the arbitrary functions $\Phi _\alpha (x,y^{(1)},...,y^{(k-2)}),$ ($\alpha =1,2$).

In particular, if the Zermelo conditions (5.2.3) are verified then the
energies $\mathcal{E}^1(H),...,\mathcal{E}^{k-1}(H)=0.$ Assuming $\Phi _2=0,$
the previous N\"other Theorem, becomes:
\begin{teo}
For any infinitesimal symmetry (5.5.7) of a differentiable Hamiltonian $%
H(x,y^{(1)},...,y^{(k-1)},p)$ , which satisfies the Zermelo conditions
(5.2.3), and for any differentiable function $\Phi (x,y^{(1)},...,y^{(k-1)})$
, along the solution curves of the Hamilton-Jacobi equations $\displaystyle
\frac{dx^i}{dt}=\displaystyle\frac 12\displaystyle\frac{\partial H}{p_i},%
\stackrel{\circ }{E}_i(H)=0,$ the following function is constant:

\begin{equation}
\begin{array}{l}
\mathcal{F}^k(H,\Phi )=p_iV^i-\displaystyle\frac 12[I_V^{k-1}(H)-%
\displaystyle\frac 1{2!}\displaystyle\frac{dI_V^{k-2}(H)}{dt}+...+\vspace{3mm}\\
+(-1)^{k-2}%
\displaystyle\frac 1{\left( k-1\right) !}\displaystyle\frac{d^{k-2}I_V^1(H)}{%
dt^{k-2}}+\Phi ] .
\end{array}
\tag{5.5.18}
\end{equation}
\end{teo}

The theory from this chapter will be applied to the geometrical study of the
Hamilton space of order $k$.

\chapter{Dual Semispray. Nonlinear Connections}

\markboth{\it{THE GEOMETRY OF HIGHER-ORDER HAMILTON SPACES\ \ \ \ \ }}{\it{Dual Semispray. Nonlinear Connections}}

The notion of semispray from the geometry of higher order Lagrange spaces has a dual correspondent in the geometrical theory of the Hamilton spaces of order $k$. The same remark is true concerning the concept of nonlinear connection which is canonical related with that of semispray.

\section{Dual Semispray}

\begin{defi}
\textit{A dual }$k$\textit{-semispray} on $T^{*k}M$ is a vector field $S\in
\mathcal{X}(\widetilde{T^{*k}M})$ with the property
\begin{equation}
JS=\stackrel{k-1}{\Gamma },  \tag{6.1.1}
\end{equation}
where $\stackrel{k-1}{\Gamma }$ is the Liouville vector field
\begin{equation}
\stackrel{k-1}{\Gamma }=y^{(1)i}\displaystyle\frac \partial {\partial
y^{(1)i}}+\cdots +(k-1)y^{(k-1)i}\displaystyle\frac \partial {\partial
y^{(k-1)i}}.  \tag{6.1.1a}
\end{equation}
\end{defi}
We have:
\begin{prop}
A dual $k$-semispray $S$ can be locally represented by
\begin{equation}
\begin{array}{lll}
S & = & y^{(1)i}\displaystyle\frac \partial {\partial x^i}+2y^{(2)i}%
\displaystyle\frac \partial {\partial y^{(1)i}}+\cdots +(k-1)y^{(k-1)i}%
\displaystyle\frac \partial {\partial y^{(k-2)i}}+ \vspace{3mm}\\
& + & k\xi ^i(x,y^{(1)},...,y^{(k-1)},p)\displaystyle\frac \partial
{\partial y^{(k-1)i}}+\eta _i(x,y^{(1)},...,y^{(k-1)},p)\displaystyle\frac
\partial {\partial p_i}.
\end{array}
\tag{6.1.2}
\end{equation}
\end{prop}

Indeed, $S$ from (6.1.2) satisfies the equation (6.1.1) for an arbitrary system of functions $\left\{ \xi
^i\right\} $ and $\left\{ \eta _i\right\} $ ($i=1$, ..., $n$), $J$ being the $k-1$-tangent endomorphism (\S 4.3).
The systems of functions $\left\{ \xi ^i(x,y^{(1)},...,y^{(k-1)},p)\right\} $
and $\left\{ \eta _i(x,y^{(1)},...,y^{(k-1)},p)\right\} $ are called \textit{%
\ the coefficients of the dual }$k$\textit{-semispray }$S$.

Since $S$ is a vector field on $T^{*k}M$ it follows that the functions $\xi
^i$ and $\eta _i$, defined on every local chart of the manifold $T^{*k}M$
are important geometrical object fields on the manifold $T^{*k}M$.

So, we have
\begin{teo}
With respect to the transformation (4.1.2) of the local coordinates on $%
T^{*k}M$, the systems of functions $\left\{ \xi ^i\right\} $ and $\left\{
\eta _i\right\} $ transform as follows:
\begin{equation}
k\widetilde{\xi }^i=k\displaystyle\frac{\partial \widetilde{y}^{(k-1)i}}{%
\partial y^{(k-1)j}}\xi ^j+(k-1)\displaystyle\frac{\partial \widetilde{y}%
^{(k-1)i}}{\partial y^{(k-2)j}}y^{(k-1)j}+\cdots +2\displaystyle\frac{%
\partial \widetilde{y}^{(k-1)i}}{\partial y^{(1)j}}y^{(2)j}+\displaystyle
\frac{\partial \widetilde{y}^{(k-1)i}}{\partial x^j}y^{(1)j},  \tag{6.1.3}
\end{equation}
\begin{equation}
\widetilde{\eta }_i=\displaystyle\frac{\partial x^j}{\partial \widetilde{x}^i%
}\eta _j+\displaystyle\frac{\partial \widetilde{p}_i}{\partial x^j}y^{(1)j}.
\tag{6.1.4}
\end{equation}
\end{teo}

\textit{Proof:} $S$ being a vector field, from (6.1.2), (4.1.4) we
have:
$$
S= y^{(1)i}\left\{ \displaystyle\frac{\partial
\widetilde{x}^m}{\partial x^i} \displaystyle\frac{\partial}
{\partial \widetilde{x}^m}+\displaystyle\frac{\partial
\widetilde{y}^{(1)m}}{\partial x^i}\displaystyle\frac{\partial}
{\partial \widetilde{y}^{(1)m}}+\cdots +
\displaystyle\frac{\partial\widetilde{y}^{(k-1)m}}{\partial
x^i}\displaystyle\frac{\partial}{\partial\widetilde{y}^{(k-1)m}}+\displaystyle\frac{\partial\widetilde{p}_m}{\partial
x^i}\displaystyle\frac{\partial}{\partial\widetilde{p}_m}\right\}+
$$
\vspace{3mm}
$$
\begin{array}{l}
+2y^{(2)i}\left\{\displaystyle\frac{\partial
\widetilde{y}^{(1)m}}{\partial
y^{(1)i}}\displaystyle\frac{\partial}{\partial\widetilde{y}^{(1)m}}+\displaystyle\frac{\partial
\widetilde{y}^{(2)m}}{\partial
y^{(1)i}}\displaystyle\frac{\partial}{\partial\widetilde{y}^{(2)m}}+\cdots
+\displaystyle\frac{\partial\widetilde{y}^{(k-1)m}}{\partial
y^{(1)i}}\displaystyle\frac{\partial}{\partial\widetilde{y}^{(k-1)m}}\right\} + \\
.............................................................................................\\
+(k-1)y^{(k-1)i}\left\{\displaystyle\frac{\partial
\widetilde{y}^{(k-2)m}}{\partial
y^{(k-2)i}}\displaystyle\frac{\partial}{\partial\widetilde{y}^{(k-2)m}}+\displaystyle\frac{\partial
\widetilde{y}^{(k-1)m}}{\partial
y^{(k-2)i}}\displaystyle\frac{\partial}{\partial\widetilde{y}^{(k-1)m}}\right\} + \vspace{3mm}\\
 + k\xi ^i\displaystyle\frac{\partial \widetilde{y}^{(k-1)m}}{\partial
y^{(k-1)i}}\displaystyle\frac{\partial}{\partial\widetilde{y}^{(k-1)m}}+
\eta_i\displaystyle\frac{\partial x^i}{\partial
\widetilde{x}^m}\displaystyle\frac{\partial}{\partial\widetilde{p}_m}.
\end{array}
$$
\vspace*{3mm}

Taking into account (6.1.4) one obtains:\\
$
\begin{array}{lll}
S & = & \widetilde{y}^{(1)m}\displaystyle\frac \partial {\partial \widetilde{
x}^m}+2\widetilde{y}^{(2)m}\displaystyle\frac \partial {\partial \widetilde{%
y }^{(1)m}}+\cdots +(k-1)\widetilde{y}^{(k-1)m}\displaystyle\frac \partial
{\partial \widetilde{y}^{(k-2)m}}+ \vspace{3mm}\\
& + & \left( y^{(1)i}\displaystyle\frac{\partial \widetilde{y}^{(k-1)m}}{
\partial x^i}+2y^{(2)i}\displaystyle\frac{\partial \widetilde{y}^{(k-1)m}}{
\partial y^{(1)i}}+\cdots +(k-1)y^{(k-1)i}\displaystyle\frac{\partial
\widetilde{y}^{(k-1)m}}{\partial y^{(k-1)i}}\right. + \vspace{3mm}\\
& + & \left. k\xi ^i\displaystyle\frac{\partial \widetilde{y}^{(k-1)m}}{
\partial y^{(k-1)i}}\right) \displaystyle\frac \partial {\partial \widetilde{
y}^{(k-1)m}}+\left( y^{(1)i}\displaystyle\frac{\partial \widetilde{p}_m}{
\partial x^i}+\eta _i\displaystyle\frac{\partial x^i}{\partial \widetilde{x}
^m}\right) \displaystyle\frac \partial {\partial \widetilde{p}_m}.
\end{array}
$\vspace{3mm}

Now, writting $S$ in the form
$$
\begin{array}{l}
S=\widetilde{y}^{(1)m}\displaystyle\frac \partial {\partial \widetilde{x}
^m}+2\widetilde{y}^{(2)m}\displaystyle\frac \partial {\partial \widetilde{y}
^{(1)m}}+\cdots +  \\
 +(k-1)\widetilde{y}^{(k-1)m}\displaystyle\frac \partial
{\partial \widetilde{y}^{(k-2)m}}+k\widetilde{\xi }^m\displaystyle\frac
\partial {\partial y^{(k-1)m}}+\widetilde{\eta }_m\displaystyle\frac
\partial {\partial \widetilde{p}_m}
\end{array}
$$
and identifying with the previous expression of $S$ we obtain the relations
(6.1.3) and (6.1.4). q.e.d.

\begin{teo}
If on every domain of local chart of the manifold $T^{*k}M$  the
systems of functions $\left\{ \xi ^i\right\} $ and $\left\{ \eta
_i\right\}$, ($i=1$, ..., $n$) are given, such that, with respect
to (4.1.2) the formulae (6.1.3) and (6.1.4) hold, then $S$ from
(6.1.2) is a dual $k$-semispray on $T^{*k}M $.
\end{teo}

The proof follow the usual way, [115].

\begin{prop}
The integral curves of the dual $k$-semispray $S$, from (6.1.2) are the
solution curves of the system of differential equations:
\begin{equation}
\left\{
\begin{array}{l}
\displaystyle\frac{dx^i}{dt}=y^{(1)},\ \displaystyle\frac{dy^{(1)i}}{dt}%
=2y^{(2)},\ ...,\ \displaystyle\frac{dy^{(k-2)i}}{dt}=(k-1)y^{(k-1)i}, \\
\\
\displaystyle\frac{dy^{(k-1)i}}{dt}=k\xi ^i(x,y^{(1)},...,y^{(k-1)},p),\ %
\displaystyle\frac{dp_i}{dt}=\eta _i(x,y^{(1)},...,y^{(k-1)},p).
\end{array}
\right.   \tag{6.1.5}
\end{equation}
\end{prop}

Notice that this system is equivalent to the following:

\begin{equation}
\left\{
\begin{array}{l}
\displaystyle\frac{dx^i}{dt}=y^{(1)},\ \displaystyle\frac 1{2!}\displaystyle
\frac{d^2x^i}{dt^2}=y^{(2)},\ ...,\ \displaystyle\frac 1{(k-1)!}%
\displaystyle \frac{d^{k-1}x^i}{dt^{k-1}}=y^{(k-1)i}, \\
\\
\displaystyle\frac 1{k!}\displaystyle\frac{d^kx^i}{dt^k}=\xi ^i(x, %
\displaystyle\frac{dx}{dt},...,\displaystyle\frac 1{(k-1)!}\displaystyle
\frac{d^{k-1}x}{dt^{k-1}},p), \\
 \displaystyle\frac{dp_i}{dt}=\eta _i(x, %
\displaystyle\frac{dx}{dt},...,\displaystyle\frac 1{(k-1)!}\displaystyle
\frac{d^{k-1}x}{dt^{k-1}},p).
\end{array}
\right.  \tag{6.1.5a}
\end{equation}

The problem of integration of this system of differential equations is solved by
usual methods. The solutions $\{(x^i(t),p_i(t)),\ t\in (a,b) \}$ are curves on the cotangent manifold $T^{\ast }M.$

Concerning the homogeneity of $S$, we see that every term which does not
include $\xi ^i$ or $\eta _i$ has the degree of homogeneity $2$ on the
fibres of $T^{*k}M$. Therefore we have:

\begin{prop}
The dual $k$-semispray $S$, from (6.1.2) is $2$-homogeneous on the fibres of $%
T^{*k}M$ if and only if the coefficients $\xi ^i$ are $k$-homogeneous and $%
\eta _i$ are $k+1$-homogeneous on the fibres of $T^{*k}M$.
\end{prop}

A dual $k$-semispray $S$, $2$-homogeneous on the fibres of $T^{*k}M$ is
called \textit{a }$k$\textit{-spray}.

Remarking that the rule (6.1.3) of transformation of the coefficients $\xi ^i$
of a $k$-semispray $S$ is same rule with that of the coordinates $y^{(k)i}$ on the
total space of bundle of accelerations $T^kM$ we can give a remarkable
geometrical meaning of the coefficients $\xi ^i$:

\begin{prop}
Every dual $k$-semispray $S$ on the manifold $T^{*k}M$ with the coefficients
$(\xi ^i,\eta _i)$ determine a bundle morphism
\[
\xi :(x,y^{(1)},...,y^{(k-1)},p)\in T^{*k}M\rightarrow
(x,y^{(1)},...,y^{(k-1)},y^{(k)})\in T^kM
\]
defined by
\begin{equation}
x^i=x^i,\ y^{(1)i}=y^{(1)i},\ ...,\ y^{(k-1)i}=y^{(k-1)i},\ y^{(k)i}=\xi
^i(x,y^{(1)},...,y^{(k-1)},p).  \tag{6.1.6}
\end{equation}

Moreover, the bundle morphism $\xi $ is a local diffeomorphism if and only if
\[
rank\left\| \stackrel{\cdot }{\partial ^i}\xi ^j\right\| =n.
\]
\end{prop}

Indeed, the mapping $\xi :T^{*k}M\rightarrow T^kM$ does not depend on the
transformation of local coordinates. It is local invertible if and only if
\\
$rank\left\| \stackrel{\cdot }{\partial ^i}\xi ^j\right\| =n$. \hfill Q.E.D.

We shall see in Chapter 8, that the bundle morphism $\xi $ defined in
(6.1.6) is uniquely determined by the Legendre transformation between a
Lagrange spaces of order $k$, $L^{(k)n}=(M,L)$ and a Hamilton spaces of
order $k$, \\
$H^{(k)n}=(M,H)$.

Consequently, if the bundle morphism $\xi $ from (6.1.6) is apriori given we
can consider the dual $k$-semispray
\begin{equation}
S_\xi =y^{(1)i}\displaystyle\frac \partial {\partial x^i}+2y^{(2)i} %
\displaystyle\frac \partial {\partial y^{(1)i}}+\cdots +(k-1)y^{(k-1)i} %
\displaystyle\frac \partial {\partial y^{(k-2)i}}+k\xi ^i\displaystyle\frac
\partial {\partial y^{(k-1)i}}+\eta _i\stackrel{\cdot }{\partial ^i}.
\tag{6.1.7}
\end{equation}
In this case $S_\xi $ is characterized only by the coefficients $\eta
_i(x,y^{(1)},...,y^{(k-1)},p)$.

The values of $1$-forms $d_0H$, ..., $d_{k-2}H$ from (4.3.14), (4.3.5) on
the vector fields $S_\xi $ are as follows:
\begin{equation}
\begin{array}{l}
d_0H(S_\xi )=\stackrel{(0)}{p_i}y^{(1)i}, \vspace{3mm}\\
d_1H(S_\xi )=\stackrel{(1)}{p_i}y^{(1)i}+2\stackrel{(0)}{p_i}y^{(2)i}, \\
..............................................................................................
\\
d_{k-2}H(S_\xi )=\stackrel{(k-2)}{p_i}y^{(1)i}+2\stackrel{(k-3)}{p_i}
y^{(2)i}+\cdots +(k-1)\stackrel{(0)}{p_i}y^{(k-1)i}.
\end{array}
\tag{6.1.8}
\end{equation}

These scalar fields does not depend on the coefficients $\xi ^i$ and $\eta
_i $ of $S_\xi $. Also the $1$-form $dH$ leads to the formula
\begin{equation}
\begin{array}{l}
dH(S_\xi )=\displaystyle\frac{\partial H}{\partial x^i}y^{(1)i}+2 %
\displaystyle\frac{\partial H}{\partial y^{(1)i}}y^{(2)i}+\cdots +(k-1) %
\displaystyle\frac{\partial H}{\partial y^{(k-2)i}}y^{(k-1)i}+ \\
\\
+k\displaystyle\frac{\partial H}{\partial y^{(k-1)i}}\xi ^i+\stackrel{\cdot
}{\partial ^i}H\eta _i=S_\xi (H).
\end{array}
\tag{6.1.8a}
\end{equation}

The existence of a $k$-semispray on the manifold $T^{*k}M$ is assured  by the
following theorem
\begin{teo}
If the base manifold $M$ is paracompact, then on the manifold $T^{*k}M$ there
exists the dual $k$-semispray $S_\xi $, with apriori given bundle morphism $%
\xi $.
\end{teo}

\textbf{Proof:} Assuming that the manifold $M$ is paracompact by the Theorem 4.1.2,
the manifold $T^{*k}M$ is paracompact, too. We
shall see (Ch. 8) that a bundle morphism $\xi $, defined in (6.1.6) exists.
Now, let $\gamma _{ij}(x)$, $x\in M$, be a Riemann metric on $M$ and $\gamma
_{jk}^i(x)$ its Christoffel symbols.

Setting
\begin{equation}
\eta _j=\gamma _{jh}^i(x)p_iy^{(1)h}  \tag{6.1.9}
\end{equation}
we can prove that the rule of transformation of the system of
functions $\{\eta _i\}$, with respect to (4.1.2) is just (6.1.4).
Applying Theorem 6.1.2 we obtain a $k$-semispray $S_\xi $, with the
coefficients $\eta _i$.

Another properties of $S$ are given in \S 4.3.
We have
\[
\stackrel{k-1}{\Gamma }=JS,\ \stackrel{k-2}{\Gamma }=J^2S,\ ...,\ \stackrel{%
1 }{\Gamma }=J^{(k-1)}S.
\]

\section{Nonlinear Connections}

The notion of nonlinear connection on the total space of the dual bundle \newline
$(T^{*k}M,\pi ^{*k},M)$ can be introduced by the classical method, [115].

\begin{defi}
A nonlinear connection on the manifold $T^{*k}M$ is a regular distribution $N
$ on the $T^{*k}M$ supplementary to the vertical distribution $V$, i.e
\end{defi}
\begin{equation}
T_u(T^{*k}M)=N_u\oplus V_u,\ \forall u\in T^{*k}M.  \tag{6.2.1}
\end{equation}

Taking into account the Proposition 4.2.1 it follows that the distribution $%
N $ has the property
\begin{equation}
T_u(T^{*k}M)=N_u\oplus V_{1,u}\oplus W_{k,u} . \tag{6.2.1a}
\end{equation}

Locally $V_1$ is generated by the system of vector fields
$(\displaystyle\frac{\partial}{\partial y^{\left(1\right)
i}},...,\displaystyle\frac{\partial}{\partial y^{\left(k-1\right)
i}})$ and $W_k$ is generated locally by $(\displaystyle\frac
\partial {\partial p_i}).$

As usual, we shall write these systems of tangent vectors
$\left(\displaystyle\frac{\partial}{\partial y^{(1)i}},...,
\displaystyle\frac{\partial}{\partial y^{(k-1)i}}\right)$ as
$(\stackrel{(1)}{\partial_i},...,\stackrel{(k-1)}{\partial_i})$
and $(\displaystyle\frac{\partial}{\partial p_i})$ as
$(\stackrel{\cdot}{\partial_i}),$ respectively.

It follows that the local dimension of the distribution $N$ is $n,$ local
dimension of the distribution $V_1$ is $(k-1)n$ and that of distribution $W_k$ is $n$.

Consider a nonlinear connection $N$ on $T^{*k}M$ and denote by $h$ and $v$ the
projectors determined by direct decomposition (6.2.1). Then we have:
\[
h+v=Id, h^2=h, v^2=v, hv=vh=0.
\]

As usual we denote
\begin{equation}
X^H=hX,\ X^V=vX,\ \forall X\in \mathcal{X}(T^{*k}M).  \tag{6.2.2}
\end{equation}

A horizontal lift, with respect to $N$ is a $\mathcal{F}(M)$-linear mapping\\
$l_h:\mathcal{X}(M)\rightarrow \mathcal{X}(T^{*k}M)$ which has the properties
\[
v\circ l_h=0,\ d\pi ^{*k}\circ l_h=I_d.
\]

There exists an unique local basis adapted to the \textit{horizonta}l
distribution $N$. It is given by

\begin{equation}
\displaystyle\frac \delta {\delta x^i}=l_h(\displaystyle\frac \partial
{\partial x^i}),(i=1,...,n).  \tag{6.2.3}
\end{equation}

The linearly independent vector fields $\displaystyle\frac \delta {\delta
x^i},$ $(i=1,...,n)$ can be uniquely written in the form:

\begin{equation}
\displaystyle\frac \delta {\delta x^i}=\displaystyle\frac \partial
{\partial x^i}-\underset{(1)}{N_i^j}\displaystyle\frac \partial
{\partial y^{\left( 1\right)
j}}-...-\underset{(k-1)}{N_i^j}\displaystyle\frac \partial
{\partial y^{\left( k-1\right) j}}+N_{ij}\displaystyle\frac
\partial {\partial p_j}.  \tag{6.2.4}
\end{equation}

The systems of functions
$\underset{(1)}{N_i^j},...,\underset{(k-1)}{ N_i^j},N_{ij}$ are
called the \textit{coefficients} of the nonlinear connection $N$.

They determine an important object fields on the manifold $T^{*k}M.$
Indeed, a change of local coordinates (4.1.2) and (6.2.3) imply:
\begin{equation}
\displaystyle\frac \delta {\delta x^i}=\displaystyle\frac{\partial
\widetilde{x}^j}{\partial x^i}\displaystyle\frac \delta {\delta \widetilde{x}
^j} , \tag{6.2.5}
\end{equation}
since
\[
\displaystyle\frac \delta {\delta x^i}\left| _u\right. =l_h(\displaystyle
\frac \delta {\delta x^i}\left| _x\right. )=l_h(\displaystyle\frac{\partial
\widetilde{x}^j}{\partial x^i}\displaystyle\frac \partial {\partial
\widetilde{x}^j})_x=(\displaystyle\frac{\partial \widetilde{x}^j}{\partial
x^i}\displaystyle\frac \delta {\delta \widetilde{x}^j})_u,x=\pi ^{*k}(u).
\]

It is not difficult to prove that the formula (6.2.5) has the following
consequence:

\begin{teo}
1$^{\circ }.$ With respect to changes of coordinate on $T^{*k}M$, (4.1.2),
 the coefficients of a nonlinear connection $N$ are transformed by the
rule
\begin{equation}
\begin{array}{l}
\underset{(1)}{\widetilde{N}_m^i}\displaystyle\frac{\partial \widetilde{x}%
^m}{\partial x^j}=\underset{(1)}{N_j^m}\displaystyle\frac{\partial
\widetilde{x}^i}{\partial x^m}-\displaystyle\frac{\partial \widetilde{y}%
^{(1)i}}{\partial x^j}, \\
\\
\underset{(2)}{\widetilde{N}_m^i}\displaystyle\frac{\partial \widetilde{x}%
^m}{\partial x^j}=\underset{(2)}{N_j^m}\displaystyle\frac{\partial
\widetilde{x}^i}{\partial x^m}+\underset{(1)}{N_j^m}\displaystyle\frac{%
\partial \widetilde{y}^{(1)i}}{\partial x^m}-\displaystyle\frac{\partial
\widetilde{y}^{(2)i}}{\partial x^j}, \\
.............................................................. \\
\underset{(k-1)}{\widetilde{N}_m^i}\displaystyle\frac{\partial \widetilde{x%
}^m}{\partial
x^j}=\underset{(k-1)}{N_j^m}\displaystyle\frac{\partial
\widetilde{x}^i}{\partial x^m}+\underset{(k-2)}{N_j^m}\displaystyle\frac{%
\partial \widetilde{y}^{(1)i}}{\partial x^m}+...+\underset{(1)}{N_j^m}%
\displaystyle\frac{\partial \widetilde{y}^{(k-2)i}}{\partial x^m}-%
\displaystyle\frac{\partial \widetilde{y}^{(k-1)i}}{\partial x^j}, \\

\\
\widetilde{N}_{ij}=\displaystyle\frac{\partial x^r}{\partial \widetilde{x}^i}%
\displaystyle\frac{\partial x^s}{\partial \widetilde{x}^j}N_{rs}+p_r%
\displaystyle\frac{\partial ^2x^r}{\partial \widetilde{x}^i\partial
\widetilde{x}^j}.
\end{array}
\tag{6.2.6}
\end{equation}

2$^{\circ }.$ Conversely, if the system of functions $($ $\underset{(1)}{%
N_i^j},...,\underset{(k-1)}{N_i^j},N_{ij})$ are given on every
domain of local chart of the manifold $T^{*k}M$ such that the
equations (6.2.6) are verified, then $($ $\underset{(1)}{N_i^j},...,\underset{(k-1)}{N_i^j}, N_{ij})$ are the coefficients of a nonlinear connection $N$ on $T^{*k}M.$
\end{teo}

The $(k-1)$-tangent structure $J,$ defined in the section 3, ch.4, maps the horizontal distribution $N$ into a vertical distribution $N_1=J(N_0),$ $N_0=N$ which is supplementary to the distribution $V_2$ in $V_1.$
If we continue this process, we obtain:
\begin{equation}
\begin{array}{l}
N_0=N,N_1=J(N_0),N_2=J^2(N_0),...,N_{k-2}=J^{k-2}(N_0), \\
\\
V_k=J^{k-1}(N_0).
\end{array}
\tag{6.2.7}
\end{equation}
Therefore, we have the direct sum of vector spaces:
\begin{equation}
V_{1,u}=N_{1,u}\oplus N_{2,u}\oplus ...\oplus N_{k-2,u}\oplus V_{k-1,u},\
\forall u\in T^{*k}M  \tag{6.2.8}
\end{equation}
and
\begin{equation}
V_u=N_{1,u}\oplus N_{2,u}\oplus ...\oplus N_{k-2,u}\oplus V_{k-1,u}\oplus
W_{k,u},\forall u\in T^{*k}M.  \tag{6.2.8a}
\end{equation}
Taking into account the formula (6.1.1), we get:
\begin{teo}
1$^{\circ }.$ A nonlinear connection $N$ on the manifold $T^{*k}M$ gives
rise to the direct sum of linear spaces for the tangent space $%
T_u(T^{*k}M):$
\begin{equation}
T_u(T^{*k}M)=N_{0,u}\oplus N_{1,u}\oplus ...\oplus N_{k-2,u}\oplus
V_{k-1,u}\oplus W_{k,u},\forall u\in T^{*k}M  \tag{6.2.9}
\end{equation}

2$^{\circ }.$ The adapted basis of every term of the previous direct sum,
respectively is:

\begin{equation}
\begin{array}{l}
\displaystyle\frac \delta {\delta x^i}=\displaystyle\frac \partial
{\partial x^i}-\underset{(1)}{N_i^j}\displaystyle\frac \partial
{\partial y^{(1)j}}-\cdots
-\underset{(k-1)}{N_i^j}\displaystyle\frac \partial
{\partial y^{(k-1)j}}+N_{ij}\displaystyle\frac \partial {\partial p_j}, \\
\\
\displaystyle\frac \delta {\delta y^{(1)i}}=\displaystyle\frac
\partial {\partial
y^{(1)i}}-\underset{(1)}{N_i^j}\displaystyle\frac \partial
{\partial y^{(2)j}}-\cdots
-\underset{(k-2)}{N_i^j}\displaystyle\frac
\partial {\partial y^{(k-1)j}}, \\
................................................................ \\
\displaystyle\frac \delta {\delta y^{\left( k-1\right) i}}=\displaystyle %
\frac \partial {\partial y^{(k-1)i}}, \\
\\
\displaystyle\frac \delta {\delta p_i}=\displaystyle\frac \partial {\partial
p_i}.
\end{array}
\tag{6.2.10}
\end{equation}
\end{teo}

\textbf{Proof:} 1$^{\circ }.$ The direct sum (6.2.9) is obtained from (6.2.1)
and (6.2.8').

2$^{\circ }.$ $J(\displaystyle\frac \delta {\delta x^i})=\displaystyle\frac
\delta {\delta y^{\left( 1\right) i}}$ is obtained by means of the
definition of the $k-1$-tangent structure $J.$ Also, $\displaystyle\frac
\delta {\delta y^{\left( 2\right) i}}=J(\displaystyle\frac \delta {\delta
y^{\left( 1\right) i}}),$ ....... leads to the formulae (6.2.10).

Remarking that (6.2.5) holds, it follows:

\begin{prop}
Under a change of local coordinates on $T^{*k}M$, the vector fields of the adapted basis:
\begin{equation}
\{\displaystyle\frac \delta {\delta x^i},\displaystyle\frac \delta {\delta
y^{\left( 1\right) i}}...,\displaystyle\frac \delta {\delta y^{\left(
k-1\right) i}},\displaystyle\frac \delta {\delta p_i}\}  \tag{6.2.11}
\end{equation}
transform by the rule:
\begin{equation}
\begin{array}{l}
\displaystyle\frac \delta {\delta x^i}=\displaystyle\frac{\partial
\widetilde{x}^j}{\partial x^i}\displaystyle\frac \delta {\delta \widetilde{x}%
^j},\displaystyle\frac \delta {\delta y^{\left( 1\right) i}}=\displaystyle
\frac{\partial \widetilde{x}^j}{\partial x^i}\displaystyle\frac \delta
{\delta \widetilde{y}^{\left( 1\right) j}},..., \\
\\

\displaystyle\frac \delta {\delta y^{\left( k-1\right) i}}=\displaystyle
\frac{\partial \widetilde{x}^j}{\partial x^i}\displaystyle\frac \delta
{\delta \widetilde{y}^{\left( k-1\right) j}},\displaystyle\frac \delta

{\delta p_i}=\displaystyle\frac{\partial x^i}{\partial \widetilde{x}^j}%
\displaystyle\frac \delta {\delta \widetilde{p}_j}.
\end{array}
\tag{6.2.11a}
\end{equation}
\end{prop}

According to (6.2.8') the vertical projector $v$ is expressed by means of
the projectors $v_1,...,v_{k-1}, w_k$ determined by the distributions $%
N_1,N_2,...,N_{k-2},V_k$ and $W_k$ as follows
\[
v=v_1+v_2+...+v_{k-1}+w_k.
\]

If we consider also the horizontal projector $h,$ determined by the
distribution $N=N_0$ we have, for $\alpha =1,...,k-1$
\begin{equation}
h+v_1+...+v_{k-1}+w_k=Id , \tag{6.2.12}
\end{equation}
\[
h^2=h,\ v_\alpha ^2=v_\alpha ,\ w_k^2=w_k,
\]
\[
h\circ v_\alpha =v_\alpha \circ h=0,\ h\circ w_\alpha =w_\alpha \circ h=0,\
v_\beta \circ v_\alpha =v_\alpha \circ v_\beta =0,\alpha \not =\beta.
\]
As usual we put:
\begin{equation}
X^H=hX,X^{V_\alpha }=v_\alpha X,\ X^{W_k}=w_kX,\ \forall X\in \mathcal{X}
(T^{*k}M).  \tag{6.2.13}
\end{equation}

In adapted basis we get
\begin{equation}
X^H=\stackrel{\left( 0\right) }{X^i}\displaystyle\frac \delta {\delta x^i},\
X^{V_\alpha }=\stackrel{\left( \alpha \right) }{X^i}\displaystyle\frac
\delta {\delta y^{\left( \alpha \right) i}},\ X^{W_k}=X_i\displaystyle\frac
\partial {\partial p_i}, (\alpha =1,...,k-1).  \tag{6.2.13a}
\end{equation}

Therefore, the formulae (6.2.11') show that one has

\begin{prop}
With respect to (4.1.2), the coordinates of the vectors $X^H,$ $X^{V_\alpha },X^{W_k}$ are changed by the rule:
\begin{equation}
\stackrel{\left( 0\right) }{X^i}=\displaystyle\frac{\partial \widetilde{x}^i%
}{\partial x^j}\stackrel{\left( 0\right) }{X^j},\ \stackrel{\left( \alpha
\right) }{X^i}=\displaystyle\frac{\partial \widetilde{x}^i}{\partial x^j}%
\stackrel{\left( \alpha \right) }{X^j},\ \widetilde{X}_i=\displaystyle\frac{%
\partial x^j}{\partial \widetilde{x}^i}X_j.  \tag{6.2.13b}
\end{equation}
\end{prop}

The following result is important:

\begin{prop}
1$^{\circ }$. The distribution $N_0$ is integrable if and only if for any $%
X,Y\in \mathcal{X}(T^{*k}M):$

\[
[X^H,Y^H]^{V_\alpha }=0,\ [X^H,Y^H]^{W_k}=0,\ (\alpha =1,...,k-1).
\]

2$^{\circ }$. The distributions $N_1,...,N_{k-2}$ are integrable if and only
if
\[
[X^{V_\alpha },Y^{V_\alpha }]^H=[X^{V_\alpha },Y^{V_\alpha }]^{V_\beta
}=[X^{V_\alpha },Y^{V_\alpha }]^{W_k}=0,
\]
for $\alpha =1,...,k-2,$ $\beta \not =\alpha ,$ $\beta =1,...,k-1,$
respectively.
\end{prop}

Indeed, the previous equations appear if we express that $[X^H,Y^H] $ is a vector field which belongs to the distribution $N_0,$ $[X^{V_\alpha },Y^{V_\alpha }]$ is a vector field which belongs to the distribution $N_\alpha,$ for $\alpha =1,.., k-2.$

Notices that the distributions $V_{k-1}$and $W_k$ are integrable. q.e.d

The notions of $h$-lift, $v_\alpha $-lift or $w_k$-lift of a vector field $%
X=X^i(x)\displaystyle\frac \partial {\partial x^i}\in \mathcal{X}(M)$ are
not difficult to define. If $l_h$ is the horizontal lift to $N_0$, $l_{v_\alpha} $ are the lifts to $N_\alpha $ and $l_{v_{k-1}}$ is the lift to $V_{k-1},$ then we have
\begin{equation}
X^H=l_hX=X^i(x)\displaystyle\frac \delta {\delta x^i},X^{V_\alpha
}=l_{v_\alpha }X=X^i\left( x\right) \displaystyle\frac \delta {\delta
y^{(\alpha)i}}  \tag{6.2.14}
\end{equation}
at the point $u\in T^{*k}M,$ $\pi ^{*k}(u)=x.$

In the case of an 1-form $\omega =\omega _i(x)dx^i$ from $\mathcal{X}^{*}(M)$
the lift $l_{w_k}(\omega )$ to the distribution $W_k$ can be defined, at
every point $u\in T^{*k}M,$with $\pi ^{*k}(u)=x,$by

\begin{equation}
\omega ^{W_k}=l_{w^{\alpha}}\omega =\omega _i\stackrel{\cdot }{\partial ^i}.
\tag{6.2.14a}
\end{equation}

\section{The Dual Coefficients of the Nonlinear Connection $N$}

Throughout this chapter, we consider the coefficients $(\underset{%
\left( 1\right) }{N_j^i},...,\underset{\left( k-1\right)
}{N_j^i},N_{ij})$ of the nonlinear connection $N$. By means of these
coefficients, the vector fields from the adapted
basis (6.2.11), are expressed . The dual basis of (6.2.11) will be denoted by
\begin{equation}
(\delta x^i,\delta y^{\left( 1\right) i},...,\delta y^{\left( k-1\right)
i},\delta p_i).  \tag{6.3.1}
\end{equation}

Remark that the coefficients of the basis (6.3.1), called the dual coefficients of
$N,$ are expressed by the coefficients of the nonlinear connection but are not coincident with them.

First of all, the conditions of duality between the adapted basis (6.2.11) and
its dual basis (6.3.1) impose the following form of the 1-form fields (6.3.1)

\begin{equation}
\begin{array}{l}
\delta x^i=dx^i, \\
\\
\delta y^{\left( 1\right) i}=dy^{\left( 1\right)
i}+\underset{\left(
1\right) }{M_j^i}dx^j, \\
..................................... \\
\delta y^{\left( k-1\right) i}=dy^{\left( k-1\right)
i}+\underset{\left( 1\right) }{M_j^i}dy^{\left( k-2\right)
i}+...+\underset{\left( k-2\right)
}{M_j^i}dy^{\left( 1\right) j}+\underset{\left( k-1\right) }{M_j^i}dx^j, \\
\\
\delta p_i=dp_i-N_{ji}dx^j.
\end{array}
\tag{6.3.2}
\end{equation}

We obtain, without difficulties, [115]

\begin{prop}
The dual coefficients $(\underset{\left( 1\right) }{M_j^i},...,\underset{%
\left( k-1\right) }{M_j^i},N_{ji})$ are uniquely determined by the coefficients $(\underset{\left( 1\right)
}{N_j^i},...,\underset{\left( k-1\right) }{N_j^i},N_{ij})$ through the following formulae:

\begin{equation}
\begin{array}{c}
\underset{\left( 1\right) }{M_j^i}=\underset{\left( 1\right)
}{N_j^i},\
\underset{\left( 2\right) }{M_j^i}=\underset{\left( 2\right) }{N_j^i}+%
\underset{\left( 1\right) }{N_m^i}\underset{\left( 1\right) }{M_j^m}, \\
.................................................. \\
\underset{\left( k-1\right) }{M_j^i}=\underset{\left( k-1\right) }{N_j^i}%
+\underset{\left( k-2\right) }{N_m^i}\underset{\left( 1\right) }{M_j^m}%
+...+\underset{\left( 1\right) }{N_m^i}\underset{\left( k-2\right) }{%
M_j^m}.
\end{array}
\tag{6.3.3}
\end{equation}
\end{prop}

The same formulae determine the coefficients $(\underset{\left( 1\right)
}{N_j^i},...,\underset{\left( k-1\right) }{N_j^i},N_{ij})$ as functions by the dual
coefficients.

These two dual basis (6.2.11) and (6.3.1) and the Proposition 6.1.1 lead to:
\begin{prop}
If we change the local coordinates on $T^{*k}M$, then the 1-form of dual basis
\[
(\delta x^i,\delta y^{\left( 1\right) i},...,\delta y^{\left( k-1\right)
i},\delta p_i)
\]
are transformed as follows:
\begin{equation}
\begin{array}{l}
\delta \widetilde{x}^i=\displaystyle\frac{\partial \widetilde{x}^i}{\partial
x^j}\delta x^j,\delta \widetilde{y}^{\left( 1\right) i}=\displaystyle\frac{%
\partial \widetilde{x}^i}{\partial x^j}\delta y^{\left( 1\right)
j},...,\delta \widetilde{y}^{\left( k-1\right) i}=\displaystyle\frac{%
\partial \widetilde{x}^i}{\partial x^j}\delta y^{\left( k-1\right) j}, \\
\\
\delta \widetilde{p}_i=\displaystyle\frac{\partial x^i}{\partial \widetilde{x%
}^i}\delta p_j.
\end{array}
\tag{6.3.4}
\end{equation}
\end{prop}

Evidently, the properties (6.3.4) imply some special rules of transformation
of dual coefficients. In this respect we remark the formulae, which can be
easily deduced
\begin{equation}
\begin{array}{l}
\displaystyle\frac \partial {\partial x^i}=\displaystyle\frac
\delta {\delta x^i}+\underset{\left( 1\right)
}{M_i^j}\displaystyle\frac \delta {\delta
y^{\left( 1\right) j}}+...+\underset{\left( k-1\right) }{M_i^j} %
\displaystyle\frac \delta {\delta y^{\left( k-1\right) j}}+N_{ij} %
\displaystyle\frac \delta {\delta p_j}, \\
\\
\displaystyle\frac \partial {\partial y^{\left( 1\right) i}}=\displaystyle %
\frac \delta {\delta y^{\left( 1\right) i}}+\underset{\left(
1\right) }{ M_i^j}\displaystyle\frac \delta {\delta y^{\left(
2\right) j}}+...+ \underset{\left( k-2\right)
}{M_i^j}\displaystyle\frac \delta {\delta
y^{\left( k-1\right) j}}, \\
..................................... \\
\displaystyle\frac \partial {\partial y^{\left( k-1\right) i}}=\displaystyle %
\frac \delta {\delta y^{\left( k-1\right) i}}, \\
\\
\displaystyle\frac \partial {\partial p_i}=\displaystyle\frac \delta {\delta
p_i},
\end{array}
\tag{6.3.5}
\end{equation}
and
\begin{equation}
\begin{array}{l}
dx^i=\delta x^i, \\
\\
dy^{\left( 1\right) i}=\delta y^{\left( 1\right)
i}-\underset{\left(
1\right) }{N_j^i}\delta x^j, \\
\\
dy^{\left( 2\right) i}=\delta y^{\left( 2\right)
i}-\underset{\left( 1\right) }{N_j^i}\delta y^{\left( 1\right)
i}-\underset{\left( 2\right) }{
N_j^i}\delta x^j, \\
...................................................... \\
dy^{\left( k-1\right) i}=\delta y^{\left( k-1\right)
i}-\underset{\left( 1\right) }{N_j^i}\delta y^{\left( k-2\right)
i}-...-\underset{\left( k-2\right) }{N_j^i}\delta y^{\left(
1\right) j}-\underset{\left(
k-1\right) }{N_j^i}\delta x^j, \\
\\
dp_i=\delta p_i+N_{ji}\delta x^j.
\end{array}
\tag{6.3.5a}
\end{equation}
Therefore, we obtain
\begin{teo}

1$^{\circ }.$ A change of local coordinates on $T^{*k}M$ implies for the
dual coefficients of the nonlinear connection $N$, the following rule of
transformations
\begin{equation}
\begin{array}{l}
\underset{\left( 1\right) }{M_j^m}\displaystyle\frac{\partial \widetilde{x}%
^i}{\partial x^m}=\underset{\left( 1\right) }{\widetilde{M}_m^i}%
\displaystyle\frac{\partial \widetilde{x}^m}{\partial x^j}+\displaystyle
\frac{\partial \widetilde{y}^{\left( 1\right) i}}{\partial x^j}, \\
\\
\underset{\left( 2\right) }{M_j^m}\displaystyle\frac{\partial \widetilde{x}%
^i}{\partial x^m}=\underset{\left( 2\right) }{\widetilde{M}_m^i}%
\displaystyle\frac{\partial \widetilde{x}^m}{\partial x^j}+\underset{%
\left( 1\right) }{\widetilde{M}_m^i}\displaystyle\frac{\partial \widetilde{y}%
^{\left( 1\right) m}}{\partial x^j}+\displaystyle\frac{\partial \widetilde{y}%
^{\left( 2\right) i}}{\partial x^j}, \\
........................................................................ \\
\underset{\left( k-1\right) }{M_j^m}\displaystyle\frac{\partial \widetilde{%
x}^i}{\partial x^m}=\underset{\left( k-1\right) }{\widetilde{M}_m^i}%
\displaystyle\frac{\partial \widetilde{x}^m}{\partial x^j}+\underset{%
\left( k-2\right) }{\widetilde{M}_m^i}\displaystyle\frac{\partial \widetilde{%
y}^{\left( 1\right) m}}{\partial x^j}+...+\underset{\left( 1\right) }{%
\widetilde{M}_m^i}\displaystyle\frac{\partial \widetilde{y}^{\left(
k-2\right) m}}{\partial x^j}+\displaystyle\frac{\partial \widetilde{y}%
^{\left( k-1\right) i}}{\partial x^j}, \\
\\
\widetilde{N}_{ij}=\displaystyle\frac{\partial x^r}{\partial \widetilde{x}^i}%
\displaystyle\frac{\partial x^r}{\partial \widetilde{x}^i}N_{rs}+p_r%
\displaystyle\frac{\partial ^2x^r}{\partial \widetilde{x}^i\partial
\widetilde{x}^j}.
\end{array}
\tag{6.3.6}
\end{equation}

2$^{\circ }.$ Conversely, if the system of function
$(\underset{\left( 1\right) }{M_j^i},...,\underset{\left(
k-1\right) }{M_j^i},N_{ij})$ are given on every domain a local
chart of the manifold $T^{*k}M$ such that the
equations (6.3.6) are verified, then $(\underset{\left( 1\right) }{M_j^i}%
,...,\underset{\left( k-1\right) }{M_j^i},N_{ji})$ are the dual
coefficients of a nonlinear connection $N$ on $T^{*k}M.$
\end{teo}

The proof can be done as in the usual cases, [115].

Some applications:

1) The vector fields $\stackrel{1}{\Gamma },\stackrel{2}{\Gamma },...,
\stackrel{k-1}{\Gamma }$ can be expressed in the adapted basis,  in the
form:
\begin{equation}
\begin{array}{l}
\stackrel{1}{\Gamma }=z^{\left( 1\right) i}\displaystyle\frac \delta {\delta
y^{\left( k-1\right) i}}, \\
\\
\stackrel{2}{\Gamma }=z^{\left( 1\right) i}\displaystyle\frac \delta {\delta
y^{\left( k-2\right) i}}+2z^{\left( 2\right) i}\displaystyle\frac \delta
{\delta y^{\left( k-1\right) i}}, \\
............................................. \\
\stackrel{k-1}{\Gamma }=z^{\left( 1\right) i}\displaystyle\frac \delta
{\delta y^{\left( 1\right) i}}+2z^{\left( 2\right) i}\displaystyle\frac
\delta {\delta y^{\left( 2\right) i}}+...+\left( k-1\right) z^{\left(
k-1\right) }\displaystyle\frac \delta {\delta y^{\left( k-1\right) i}},
\end{array}
\tag{6.3.7}
\end{equation}
where
\begin{equation}
\begin{array}{l}
z^{\left( 1\right) i}=y^{\left( 1\right) i}, \vspace{3mm}\\
2z^{\left( 2\right) i}=2y^{\left( 2\right)
i}+\underset{(1)}{M_m^i}
y^{\left( 1\right) m} \\
............................................. \\
\left( k-1\right) z^{\left( k-1\right) i}=\left( k-1\right)
y^{\left( k-1\right) i}+\left( k-2\right)
\underset{(1)}{M_m^i}y^{\left( k-2\right)
m}+... \\
...+\underset{(k-2)}{M_m^i}y^{\left( 1\right) m}.
\end{array}
\tag{6.3.7a}
\end{equation}
With respect to (4.1.2), we have
\begin{equation}
\widetilde{z}^{\left( \alpha \right) i}=\displaystyle\frac{\partial
\widetilde{x}^i}{\partial x^j}z^{\left( \alpha \right) j},\left( \alpha
=1,...,k-1\right)  \tag{6.3.7b}
\end{equation}

>From (6.3.7b) it follows that $z^{\left( 1\right) i},...,z^{\left(
k-1\right) i}$ have a geometrical meaning. They will be called the \textit{\
Liouville d-vector fields}.

\vspace{3mm}

2) The operators $d_0,...,d_{k-2}$ and $d$ are represented in the adapted
basis by the formulae
\begin{equation}
\begin{array}{l}
d_0=\displaystyle\frac \delta {\delta y^{\left( k-1\right) i}}\delta x^i, \\
\\
d_1=\displaystyle\frac \delta {\delta y^{\left( k-2\right) i}}\delta x^i+ %
\displaystyle\frac \delta {\delta y^{\left( k-1\right) i}}\delta y^{\left(
1\right) i} ,\\
................................................., \\
d_{k-2}=\displaystyle\frac \delta {\delta y^{\left( 1\right) i}}\delta x^i+ %
\displaystyle\frac \delta {\delta y^{\left( 2\right) i}}\delta y^{\left(
1\right) i}+...+\displaystyle\frac \delta {\delta y^{\left( k-1\right)
i}}\delta y^{\left( k-2\right) i}
\end{array}
\tag{6.3.8}
\end{equation}
and
\begin{equation}
d=\displaystyle\frac \delta {\delta x^i}\delta x^i+\displaystyle\frac \delta
{\delta y^{\left( 1\right) i}}\delta y^{\left( 1\right) i}+...+\displaystyle %
\frac \delta {\delta y^{\left( k-1\right) i}}\delta y^{\left( k-1\right) i}+ %
\displaystyle\frac \delta {\delta p_i}\delta p_i.  \tag{6.3.9}
\end{equation}

If $H\in \mathcal{F}(T^{*k}M),$ then we have the 1-forms
\begin{equation}
\begin{array}{l}
d_0H=\displaystyle\frac{\delta H}{\delta y^{\left( k-1\right) i}}\delta x^i,
\\
\\
d_1H=\displaystyle\frac{\delta H}{\delta y^{\left( k-2\right) i}}\delta x^i+ %
\displaystyle\frac{\delta H}{\delta y^{\left( k-1\right) i}}\delta y^{\left(
1\right) i} ,\\
................................................., \\
d_{k-2}H=\displaystyle\frac{\delta H}{\delta y^{\left( 1\right) i}}\delta
x^i+\displaystyle\frac{\delta H}{\delta y^{\left( 2\right) i}}\delta
y^{\left( 1\right) i}+...+\displaystyle\frac{\delta H}{\delta y^{\left(
k-1\right) i}}\delta y^{\left( k-2\right) i}
\end{array}
\tag{6.3.10}
\end{equation}
and
\begin{equation}
dH=\displaystyle\frac{\delta H}{\delta x^i}\delta x^i+\displaystyle\frac{
\delta H}{\delta y^{\left( 1\right) i}}\delta y^{\left( 1\right) i}+...+ %
\displaystyle\frac{\delta H}{\delta y^{\left( k-1\right) i}}\delta y^{\left(
k-1\right) i}+\displaystyle\frac{\delta H}{\delta p_i}\delta p_i.  \tag{6.3.11}
\end{equation}
In the previous formulas
\begin{equation}
\displaystyle\frac{\delta H}{\delta x^i},...,\displaystyle\frac{\delta H}{
\delta y^{\left( k-1\right) i}}  \tag{6.3.12}
\end{equation}
are d-covector and $\displaystyle\frac{\delta H}{\delta p_i}$ is a d-vector.

\qquad

3) In the adapted basis (6.2.11) a dual $k$-semispray $S_\xi $ can be written
as follows
\begin{equation}
\begin{array}{l}
S_\xi =z^{\left( 1\right) i}\displaystyle\frac \delta {\delta
x^i}+2z^{\left( 2\right) i}\displaystyle\frac \delta {\delta y^{\left(
1\right) i}}+...+\left( k-1\right) z^{\left( k-1\right) i}\displaystyle\frac
\delta {\delta y^{\left( k-2\right) i}}+ \\
\qquad +k\stackrel{\vee}{\xi }^i\displaystyle\frac \delta {\delta y^{\left(
k-1\right) i}}+\stackrel{\vee}{\eta }_i\displaystyle\frac \delta {\delta
p_i},
\end{array}
\tag{6.3.13}
\end{equation}

where
\begin{equation}
\begin{array}{l}
k\stackrel{\vee}{\xi }^i=k\xi ^i+(k-1)y^{\left( k-1\right) s}\underset{%
(1) }{M_s^i}+...+y^{\left( 1\right) s}\underset{(k-1)}{M_s^i} ,\\
\\
\stackrel{\vee}{\eta }_i=\eta _i+y^{\left( 1\right) s}N_{si}.
\end{array}
\tag{6.3.14}
\end{equation}

Evidently, $\stackrel{\vee}{\xi }^i$is a $d$-vector field, and $\stackrel{
\vee}{\eta }_i$ is a $d$-covector field. It is remarkable that the equation $%
\stackrel{\vee}{\eta }_i=0$ has a geometrical meaning. Especially the
sprays $S_\xi $ with this property will be considered.

The interior products of the 1-forms $d_0H,...,d_{k-2}H$ and $dH$ with the
vector field $S_\xi $ are given by

\begin{equation}
\begin{array}{l}
d_0H(S_\xi )=z^{\left( 1\right) i}\displaystyle\frac{\delta H}{\delta
y^{\left( k-1\right) i}} ,\\
\\
d_1H(S_\xi )=z^{\left( 1\right) i}\displaystyle\frac{\delta H}{\delta
y^{\left( k-2\right) i}}+2z^{\left( 2\right) i}\displaystyle\frac{\delta H}{
\delta y^{\left( k-1\right) i}}, \\
............................................................, \\
d_{k-2}H(S_\xi )=z^{\left( 1\right) i}\displaystyle\frac{\delta H}{\delta
y^{\left( 1\right) i}}+2z^{\left( 2\right) i}\displaystyle\frac{\delta H}{
\delta y^{\left( 2\right) i}}+...+(k-1)z^{\left( k-1\right) i}\displaystyle
\frac{\delta H}{\delta y^{\left( k-1\right) i}}
\end{array}
\tag{6.3.15}
\end{equation}
and
\begin{equation}
\begin{array}{l}
dH(S_\xi )=z^{\left( 1\right) i}\displaystyle\frac{\delta H}{\delta x^i}
+2z^{\left( 2\right) i}\displaystyle\frac{\delta H}{\delta y^{\left(
1\right) i}}+...+(k-1)z^{\left( k-1\right) i}\displaystyle\frac{\delta H}{
\delta y^{\left( k-2\right) i}}+ \\
\\
\qquad +k\stackrel{\vee}{\xi }^i\displaystyle\frac \delta {\delta y^{\left(
k-1\right) i}}+\stackrel{\vee}{\eta }_i\displaystyle\frac{\delta H}{\delta
p_i}.
\end{array}
\tag{6.3.16}
\end{equation}

Consequently, we have

\begin{equation}
dH(S_\xi )=S_\xi (H) . \tag{6.3.16a}
\end{equation}

We shall use all these formulae in the theory of higher order Hamilton
spaces.

An 1-form fields $\omega \in \mathcal{X}^{*}(T^{*k}M)$ can be uniquely
written as

\begin{equation}
\omega =\omega ^H+\omega ^{V_1}+...+\omega ^{V_{k-1}}+\omega ^{W_k},
\tag{6.3.17}
\end{equation}
where
\begin{equation}
\omega ^H=\omega \circ h,\ \omega ^{V_\alpha }=\omega \circ v_\alpha ,\
(\alpha =1,...,k-1),\ \omega ^{W_k}=\omega \circ w_k . \tag{6.3.18}
\end{equation}

These components can be easily written in the adapted basis.

For a function $H\in F(T^{*k}M)$ we deduce

\begin{equation}
dH=(dH)^H+(dH)^{V_1}+...+(dH)^{V_{k-1}}+(dH)^{W_k} , \tag{6.3.19}
\end{equation}
with
\begin{equation}
(dH)^H=\displaystyle\frac{\delta H}{\delta x^i}dx^i,(dH)^{V_1}=\displaystyle
\frac{\delta H}{\delta y^{\left( 1\right) i}}\delta y^{\left( 1\right)
i},...,(dH)^{W_k}=\displaystyle\frac{\delta H}{\delta p_i}\delta p_i.
\tag{6.3.19a}
\end{equation}

In the case of 1-forms $d_0H,...,d_{k-2}H$ we have:
\[
\begin{array}{l}
(d_0H)^{V_\alpha }=0,(d_0H)^{W_k}=0,(\alpha =1,...,k-1) ,\\
\\
(d_1H)^{V_\alpha }=0,(d_1H)^{W_k}=0,(\alpha =2,...,k-1) ,\\
................................................................, \\
(d_{k-2}H)^{V_{k-1}}=0,(d_{k-2}H)^{W_k}=0.
\end{array}
\]

4) The horizontal curves with respect to a nonlinear connection $N,$ can be
studied exactly as in Ch. 10, from the book [115].

Let $\gamma :I\rightarrow T^{*k}M$ be a parametrized curve, locally
expressed by
\begin{equation}
x^i=x^i(t),y^{\left( 1\right) i}=y^{\left( 1\right) i}(t),...,y^{\left(
k-1\right) i}=y^{\left( k-1\right) i}(t),p_i=p_i(t),t\in I . \tag{6.3.20}
\end{equation}

For the tangent vector field $\displaystyle\frac{d\gamma }{dt},t\in I$ one can
write:

\begin{equation}
\displaystyle\frac{d\gamma }{dt}=(\displaystyle\frac{d\gamma }{dt})^H+( %
\displaystyle\frac{d\gamma }{dt})^{V_1}+...+(\displaystyle\frac{d\gamma }{dt}
)^{V_{k-1}}+(\displaystyle\frac{d\gamma }{dt})^{W_k}.  \tag{6.3.20a}
\end{equation}

Using the basis (6.2.11), adapted to the direct decomposition (6.2.9), we obtain

\begin{equation}
\displaystyle\frac{d\gamma }{dt}=\displaystyle\frac{\delta x^i}{dt} %
\displaystyle\frac \delta {\delta x^i}+\displaystyle\frac{\delta y^{\left(
1\right) i}}{dt}\displaystyle\frac \delta {\delta y^{\left( 1\right)
i}}+...+ \displaystyle\frac{\delta y^{\left( k-1\right) i}}{dt}\displaystyle%
\frac \delta {\delta y^{\left( k-1\right) i}}+\displaystyle\frac{\delta p_i}{%
dt} \displaystyle\frac \delta {\delta p_i}.  \tag{6.3.20b}
\end{equation}

Taking into account the formulas (6.3.2) the coefficients of the tangent
vector field $\displaystyle\frac{d\gamma }{dt\text{ }}$ from (6.3.20b) have
the following expressions
\begin{equation}
\begin{array}{l}
\displaystyle\frac{\delta x^i}{dt}=\displaystyle\frac{dx^i}{dt}, \\
\\
\displaystyle\frac{\delta y^{\left( 1\right) i}}{dt}=\displaystyle\frac{
dy^{\left( 1\right) i}}{dt}+\underset{\left( 1\right) }{M_m^i}%
\displaystyle \frac{dx^m}{dt}, \\
..................................................... \\
\displaystyle\frac{\delta y^{\left( k-1\right) i}}{dt}=\displaystyle\frac{
dy^{\left( k-1\right) i}}{dt}+\underset{\left( 1\right) }{M_m^i} %
\displaystyle\frac{dy^{(k-2)i}}{dt}+...+\underset{(k-2)}{M_m^i} %
\displaystyle\frac{dy^{\left( 1\right) i}}{dt}+\underset{(k-1)}{M_m^i} %
\displaystyle\frac{dx^m}{dt}, \\
\\
\displaystyle\frac{\delta p_i}{dt}=\displaystyle\frac{dp_i}{dt}-N_{mi} %
\displaystyle\frac{dx^m}{dt}.
\end{array}
\tag{6.3.21}
\end{equation}

An \textit{horizontal curve} $\gamma :I\rightarrow T^{*k}M$, with respect
to the nonlinear connection $N$ is defined by the conditions
\[
(\displaystyle\frac{d\gamma }{dt})^{V_1}=...=(\displaystyle\frac{d\gamma }{%
dt })^{V_{k-1}}=0,(\displaystyle\frac{d\gamma }{dt})^{W_k}=0
\]
It follows:
\begin{teo}
A parametrized curve $\gamma :I\rightarrow T^{*k}M$ is horizontal if and
only if the following system of ordinary differential equations is verified:
\begin{equation}
\displaystyle\frac{\delta y^{\left( 1\right) i}}{dt}=...=\displaystyle\frac{%
\delta y^{\left( k-1\right) i}}{dt}=0,\displaystyle\frac{\delta p_i}{dt}=0.
\tag{6.3.22}
\end{equation}
\end{teo}

The horizontal curves with the property
\begin{equation}
y^{\left( 1\right) i}=\displaystyle\frac{dx^i}{dt},...,y^{\left( k-1\right)
i}=\displaystyle\frac 1{\left( k-1\right) !}\displaystyle\frac{d^{k-1}x^i}{
dt^{k-1}}  \tag{*}
\end{equation}
are called \textit{the autoparalel curves} of the nonlinear connection $N.$
These curves are characterized by (6.3.22) in the conditions (*) .

\section{The Determination of the Nonlinear Connection by a Dual $k$%
-Semispray}

The main problem concerning the notion of nonlinear connection $N$ is wether a
dual $k$-semispray $S_\xi $ determines a nonlinear connection $N$ or not. The answer is yes. We have the following important result.

\begin{teo}
Any dual $k$-semispray $S_\xi $ with the coefficients $(\xi ^i,\eta _i)$
determines:

1$^{\circ }.$ The dual coefficients $\underset{(1)}{M_j^i},...,\underset{%
(k-1)}{M_j^i},$ of a nonlinear connection $N$, by the formulas:

\begin{equation}
\underset{(1)}{M_j^i},=-\displaystyle\frac{\partial \xi
^i}{\partial
y^{\left( k-1\right) j}},\underset{(2)}{M_j^i},=-\displaystyle\frac{%
\partial \xi ^i}{\partial y^{\left( k-2\right) j}},...,\underset{(k-1)}{%
M_j^i}=-\displaystyle\frac{\partial \xi ^i}{\partial y^{\left( 1\right) j}}.
\tag{6.4.1}
\end{equation}

2$^{\circ }$. The coefficients $N_{ij}$ by the formula

\begin{equation}
N_{ij}=\displaystyle\frac{\delta \eta _i}{\delta y^{\left( 1\right) j}},
\tag{6.4.2}
\end{equation}
where the operators
\[
\begin{array}{l}
\displaystyle\frac \delta {\delta y^{\left( 1\right) j}}=\displaystyle\frac
\partial {\partial y^{\left( 1\right) j}}-\underset{(1)}{N_j^i}%
\displaystyle\frac \partial {\partial y^{\left( 2\right) j}}-...-\underset{%
(k-2)}{N_j^i}\displaystyle\frac \partial {\partial y^{\left( k-1\right)
j}},\ \underset{(\alpha )}{N_j^i}, \\
\\
(\alpha =1,...,k-2)
\end{array}
\]
are determined by
$\underset{(1)}{M_j^i},...,\underset{(k-2)}{M_j^i}$ from (6.4.1).
\end{teo}

\textit{Proof}: We prove that with respect to a change of local
coordinates on the manifold $T^{*k}M,$ the system of functions
$\underset{(1)}{M_j^i} ,...,\underset{(k-1)}{M_j^i}$ from (6.4.1)
obey the transformation (6.3.6).

By means of (6.1.3), we have
\[
k\widetilde{\xi }^i=k\displaystyle\frac{\partial \widetilde{x}^i}{\partial
x^j}\xi ^j+(k-1)\displaystyle\frac{\partial \widetilde{y}^{\left( 1\right)
i} }{\partial x^j}y^{\left( k-1\right) j}+...+\displaystyle\frac{\partial
\widetilde{y}^{\left( k-1\right) i}}{\partial x^j}y^{\left( 1\right) j}.
\]
Applying the formula $\displaystyle\frac \partial {\partial y^{\left(
k-1\right) m}}=\displaystyle\frac{\partial \widetilde{x}^s}{\partial x^m} %
\displaystyle\frac \partial {\partial \widetilde{y}^{\left( k-1\right) }},$
we deduce the first formula (6.3.6). By the same method, we
establish inductively the other formula (6.3.6), for the coefficients (6.4.1).

In order to prove (6.4.2), we remark that the vector fields $\displaystyle\frac \delta {\delta
y^{\left( 1\right) i}},$ constructs by means of coefficients $\underset{(1)}{M_j^i},...,\underset{%
(k-1)}{M_j^i},$ from (6.4.1) has the law of
transformation:
\[
\displaystyle\frac \delta {\delta \widetilde{y}^{\left( 1\right) i}}= %
\displaystyle\frac{\partial x^j}{\partial \widetilde{x}^i}\displaystyle\frac
\delta {\delta y^{\left( 1\right) j}}.
\]
The rule of transformation of the coefficients $\eta _i$ of $S_\xi $ is
\[
\widetilde{\eta }_i=\displaystyle\frac{\partial x^r}{\partial \widetilde{x}%
^i }\eta _r+\displaystyle\frac{\partial \widetilde{p}_i}{\partial x^r}%
y^{\left( 1\right) r}.
\]
These formula have as a consequence:
\[
\displaystyle\frac{\delta \widetilde{\eta }_i}{\delta \widetilde{y}^{\left(
1\right) j}}=\displaystyle\frac{\partial x^r}{\partial \widetilde{x}^i} %
\displaystyle\frac{\partial x^m}{\partial \widetilde{x}^j}\displaystyle\frac{
\delta \eta _r}{\delta y^{\left( 1\right) m}}+p_r\displaystyle\frac{\partial
^2x^r}{\partial \widetilde{x}^r\partial \widetilde{x}^j}.
\]

This is the rule of transformation of the coefficients $N_{ij}$ of a
nonlinear connection $N$ on the manifold $T^{*k}M.$

\section{Lie Brackets. Exterior Differential}

In the following it is important to determine the Lie brackets of the vector
fields of the adapted basis (6.2.11) and the exterior differentials of the covector
fields of adapted cobasis (6.3.1).

By a direct calculus we obtain

\begin{prop}
The following expressions of the Lie brackets hold true:
\[
[\displaystyle\frac \delta {\delta x^j},\displaystyle\frac \delta
{\delta x^h}]=\underset{(01)}{R_{jh}^i}\displaystyle\frac \delta
{\delta y^{\left( 1\right)
i}}+...+\underset{(0,k-1)}{R_{jh}^i}\displaystyle\frac \delta
{\delta y^{\left( k-1\right)
i}}+\underset{(0)}{R_{ijh}}\displaystyle\frac \delta {\delta p_i},
\]
\[
[\displaystyle\frac \delta {\delta x^j},\displaystyle\frac \delta {\delta
y^{\left( \alpha \right) h}}]=\underset{(\alpha 1)}{B_{jh}^i}%
\displaystyle
\frac \delta {\delta y^{\left( 1\right) i}}+...+\underset{(\alpha ,k-1)}{%
B_{jh}^i}\displaystyle\frac \delta {\delta y^{\left( k-1\right) i}}+%
\underset{(\alpha )}{B_{ijh}}\displaystyle\frac \delta {\delta
p_i},
\]
\begin{equation}
[\displaystyle\frac \delta {\delta y^{\left( \alpha \right) j}},%
\displaystyle \frac \delta {\delta y^{\left( \beta \right) h}}]=\stackrel{(1)%
}{\underset{(\alpha \beta )}{C_{jh}^i}}\displaystyle\frac \delta
{\delta
y^{\left( 1\right) i}}+...+\stackrel{(k-1)}{\underset{(\alpha \beta )}{%
C_{jh}^i}}\displaystyle\frac \delta {\delta y^{\left( k-1\right) i}}+%
\underset{(\alpha \beta )}{B_{ijh}}\displaystyle\frac \delta
{\delta p_i}, \tag{6.5.1}
\end{equation}
\[
[\displaystyle\frac \delta {\delta y^{\left( \alpha \right) j}},%
\displaystyle \frac \delta {\delta p_h}]=\underset{(\alpha
1)}{C_j^{ih}}\displaystyle
\frac \delta {\delta y^{\left( 1\right) i}}+...+\underset{(\alpha ,k-1)}{%
C_j^{ih}}\displaystyle\frac \delta {\delta y^{\left( k-1\right) i}}+%
\underset{(\alpha )}{C_{ij}^h}\displaystyle\frac \delta {\delta
p_i},
\]
\[
[\displaystyle\frac \delta {\delta x^j},\displaystyle\frac \delta
{\delta p_h}]=\underset{(1)}{B_j^{ih}}\displaystyle\frac \delta
{\delta y^{\left( 1\right)
i}}+...+\underset{(k-1)}{B_j^{ih}}\displaystyle\frac \delta
{\delta y^{\left( k-1\right)
i}}+\underset{(0)}{B_{ij}^h}\displaystyle \frac \delta {\delta
p_i},
\]
$\alpha ,\beta =1,...k;\beta \leq \alpha $
in which:
\begin{equation}
\begin{array}{l}
\stackrel{(1)}{\underset{(\alpha \beta )}{C_{jh}^i}}=0,\ \ \ 1\leq
\alpha \leq
\beta,  \\
\\
\stackrel{(2)}{\underset{(\alpha \beta )}{C_{jh}^i}}=0,\ \ \ 2\leq
\alpha \leq
\beta,  \\
............................... \\
\stackrel{(\gamma )}{\underset{(\alpha \beta
)}{C_{jh}^i}}=0,\ \ \ \gamma \leq \alpha \leq \beta
\end{array}
\tag{6.5.2}
\end{equation}
\end{prop}
{\it and the coefficients $R,B,C$ can be obtained by a straiforward calculus, $%
\underset{(0\alpha )}{R_{jh}^i}$ being expressed by:}
\begin{equation}
\begin{array}{l}
\underset{(01)}{R_{jh}^i}=\underset{(01)}{r_{jh}^i}, \\
\\
\underset{(02)}{R_{jh}^i}=\underset{(02)}{r_{jh}^i}+\underset{(1)}{
M_s^i}\underset{(01)}{r_{jh}^s}, \\
......................... \\
\underset{(0,k-1)}{R_{jh}^i}=\underset{(0,k-1)}{r_{jh}^i}+\underset{%
(1)
}{M_s^i}\underset{(0,k-2)}{r_{jh}^s}+...+\underset{(k-2)}{M_s^i}
\underset{(0,k-2)}{r_{jh}^s}
\end{array}
\tag{6.5.3}
\end{equation}
and
\begin{equation}
R_{ijh}=\displaystyle\frac{\delta N_{ji}}{\delta x^h}-\displaystyle\frac{
\delta N_{ki}}{\delta x^j}, \tag{6.5.3a}
\end{equation}
with
\begin{equation}
\underset{(0\alpha )}{r_{jh}^i}=\displaystyle\frac{\delta
\underset{ (\alpha )}{N_j^i}}{\delta
x^h}-\displaystyle\frac{\delta \underset{(\alpha )}{N_h^i}}{\delta
x^j}, (\alpha =1,...k-1).  \tag{6.5.3b}
\end{equation}

Taking into account the conditions of integrability of the nonlinear
connection $N,$ (cf Prop. 6.2.3) we have:

\begin{teo}
The nonlinear connection $N$ is integrable if and only if the following
equations hold:
\begin{equation}
\underset{(0\alpha )}{R_{jh}^i}=0,(\alpha =1,...,k-1),\underset{(0)}{%
R_{ijh}}=0 . \tag{6.5.4}
\end{equation}
\end{teo}

Indeed, the distribution $N$ is integrable if and only if the vector fields $[\displaystyle
\frac \delta {\delta x^j},\displaystyle\frac \delta {\delta x^h}]$ belong
to $N$. So, the equations (6.5.4) (which have a geometrical meaning) express
the necessary and sufficient conditions for $N$ to be integrable.

Similar characterizations hold for the case when each distribution $%
N_1,...,N_{k-2}$ is integrable. (cf. Prop.6.2.3, 2$^{\circ }$)

In order to determine the exterior differentials of 1-forms \\
$(\delta x^i,\delta y^{\left( 1\right) i},...,\delta y^{\left(
k-1\right) i},\delta p_i),$ we start from the formulae (6.3.2) and
express the mentioned differentials with respect to the base of
$\wedge ^2(T^{*k}M)$ determined by $(\delta x^i,\delta y^{\left(
1\right) i},...,\delta y^{\left( k-1\right)i},\delta p_i).$

>From (6.3.2) one obtains:
\begin{equation}
\begin{array}{l}
d\delta x^i=0, \\
\\
d\delta y^{\left( 1\right) i}=d\underset{(1)}{M_j^i}\wedge dx^j, \\
\\
d\delta y^{\left( 2\right) i}=d\underset{(1)}{M_j^i}\wedge
dy^{\left(
1\right) j}+d\underset{(2)}{M_j^i}\wedge dx^j, \\
........................................................ \\
d\delta y^{\left( k-1\right) i}=d\underset{(1)}{M_j^i}\wedge
dy^{\left( k-2\right) j}+d\underset{(2)}{M_j^i}\wedge dy^{\left(
k-3\right) j}+...+d \underset{(k-1)}{M_j^i}\wedge dx^j
\end{array}
\tag{6.5.5}
\end{equation}
and
\begin{equation}
d\delta p_i=-dN_{ji}\wedge dx^j.  \tag{6.5.5a}
\end{equation}

Substituting $(dx^i,dy^{\left( 1\right) i},...,dy^{\left( k-1\right) i})$
from (6.3.5') we can write (6.5.5) in the following form
\begin{equation}
\begin{array}{l}
d\delta y^{\left( \alpha \right) i}=\underset{(\alpha
0)}{P_j^i}\wedge dx^j+\underset{(\alpha 1)}{P_j^i}\wedge \delta
y^{\left( 1\right) j}+...+ \underset{(\alpha ,\alpha
-1)}{P_j^i}\wedge \delta y^{\left( \alpha
-1\right) j}, \\
\\
(\alpha =1,...,k-1), \\
\\
d\delta p_i=P_{ij}\wedge dx^j,
\end{array}
\tag{6.5.6}
\end{equation}
where $\underset{(\alpha \beta )}{P_j^i}$ are 1-forms, which should be calculated by means of formula
\[
d\omega (X,Y)=X\omega (Y)-Y\omega (X)-\omega ([X,Y]),\forall \omega \in
\wedge ^1(T^{*k}M)
\]
and using the Lie brackets (6.5.1).

\begin{teo}
The exterior differentials of the 1-forms\\
$(\delta x^i,\delta y^{\left( 1\right) i},...,\delta y^{\left( k-1\right)
i},\delta p_i)$ are given by
\begin{equation}
\begin{array}{l}
d\delta x^i=0, \\
\\
d\delta y^{\left( 1\right) i}=\underset{(10)}{P_j^i}\wedge dx^j, \\
\\
d\delta y^{\left( 2\right) i}=\underset{(20)}{P_j^i}\wedge dx^j+%
\underset{(21)}{P_j^i}\wedge \delta y^{\left( 1\right) j}, \\
........................................................ \\
d\delta y^{\left( k-1\right) i}=\underset{(k-1,0)}{P_j^i}\wedge dx^j+%
\underset{(k-1,1)}{P_j^i}\wedge \delta y^{\left( 1\right) j}+...+%
\underset{(k-1,k-2)}{P_j^i}\wedge \delta y^{\left( k-2\right) j}
\end{array}
\tag{6.5.7}
\end{equation}
and
\begin{equation}
d\delta p_i=P_{ij}\wedge dx^j , \tag{6.5.7a}
\end{equation}
where the 1-forms $\underset{(\alpha \beta )}{P_j^i}$ and $P_{ij}$ are
given by
\begin{equation}
\underset{(10)}{P_j^i}=\displaystyle\frac 12\underset{(01)}{R_{jm}^i}%
dx^m+\underset{\gamma =1}{\stackrel{k-1}{\sum }}\underset{(\gamma 1)}{%
B_{jm}^i}\delta y^{\left( \gamma \right)
m}+\underset{(1)}{B_j^{im}}\delta p_m,  \tag{6.5.8$_1$}
\end{equation}
\begin{equation}
\left\{
\begin{array}{c}
\underset{(20)}{P_j^i}=\displaystyle\frac 12\underset{(02)}{R_{jm}^i}%
dx^m+\underset{\gamma =1}{\stackrel{k-1}{\sum }}\underset{(\gamma 2)}{%
B_{jm}^i}\delta y^{\left( \gamma \right)
m}+\underset{(2)}{B_j^{im}}\delta
p_m, \\
\\
\underset{(21)}{P_j^i}=-\underset{(12)}{B_{mj}^i}dx^m-\underset{\gamma
=1}{\stackrel{k-1}{\sum }}\stackrel{\left( 2\right) }{\underset{(\gamma 1)%
}{C_{mj}^i}}\delta y^{\left( \gamma \right) m}+\underset{(12)}{C_j^{im}}%
\delta p_m,
\end{array}
\right.   \tag{6.5.8$_2$}
\end{equation}

\begin{center}
.......................................................................................
\end{center}

\begin{equation}
\begin{array}{l}
\underset{(\alpha 0)}{P_j^i}=\displaystyle\frac 12\underset{(0\alpha )}{%
R_{jm}^i}dx^m+\underset{\gamma =1}{\stackrel{k-1}{\sum }}\underset{%
(\gamma \alpha )}{B_{jm}^i}\delta y^{\left( \gamma \right) m}+\underset{%
(\alpha )}{B_j^{im}}\delta p_m, \\
\\
\underset{(\alpha 1)}{P_j^i}=-\underset{(1\alpha )}{B_{mj}^i}dx^m-%
\underset{\gamma =1}{\stackrel{k-1}{\sum }}\stackrel{\left( \alpha
\right)
}{\underset{(\gamma 1)}{C_{mj}^i}}\delta y^{\left( \gamma \right) m}+%
\underset{(1\alpha )}{C_j^{im}}\delta p_m, \\
................................................................... \\
\underset{(\alpha ,\alpha -1)}{P_j^i}=-\underset{(\alpha -1,\alpha )}{%
B_{jm}^i}dx^m-\underset{\gamma =1}{\stackrel{k-1}{\sum
}}\stackrel{\left( \alpha \right) }{\underset{(\gamma ,\alpha
-1)}{C_{mj}^i}}\delta y^{\left(
\gamma \right) m}+\underset{(\alpha -1,\alpha )}{C_j^{im}}\delta p_m, \\
\\
(\alpha =3,...,k-1)
\end{array}
\tag{6.5.8$_\alpha $}
\end{equation}
and
\begin{equation}
P_{ij}=\displaystyle\frac
12\underset{(0)}{R_{ijm}}dx^m+\underset{\gamma
=1}{\stackrel{k-1}{\sum }}\underset{(\gamma )}{B_{ijm}}\delta
y^{\left( \gamma \right) m}+\underset{(0)}{B_{ij}^m}\delta p_m.
\tag{6.5.9}
\end{equation}
\end{teo}

\begin{rem}
In order to give a proof of the equations (6.5.2) we remark that
\[
d\delta y^{\left( \gamma \right) i}(\displaystyle\frac \delta {\delta
y^{(\alpha )j}},\displaystyle\frac \delta {\delta y^{(\beta )h}})=-\stackrel{%
\left( \gamma \right) }{\underset{(\alpha \beta )}{C_{jh}^i}}
\]
\end{rem}

>From (6.5.6) it follows that the left hand side of the previous formula vanishes
for $\gamma <\alpha \leq \beta .$

\begin{rem}
By this method we can calculate the coefficients for the expression of the Lie brackets (6.5.1).
\end{rem}

Indeed, $d\delta ^{\left( \alpha \right) i}$ from (6.5.5) can be written in
the adapted basis
\begin{equation}
\begin{array}{l}
d\delta y^{\left( \alpha \right) i}=\{d\underset{\left( \alpha \right)
}{M_j^i}- \underset{\left( 1\right) }{N_j^m}d\underset{\left(
\alpha -1\right) }{ M_m^i}-...-\underset{\left( \alpha -1\right)
}{N_j^i}d\underset{\left(
1\right) }{M_j^i}\}\wedge dx^j+ \\
\\
\qquad \{d\underset{\left( \alpha -1\right)
}{M_j^i}-\underset{\left( 1\right) }{N_j^m}d\underset{\left(
\alpha -2\right) }{M_m^i}-...- \underset{\left( \alpha -2\right)
}{N_j^i}d\underset{\left( 1\right) }{
M_j^i}\}\wedge dy^{\left( 1\right) j}+ \\
\\
\qquad +...+\{d\underset{\left( 2\right) }{M_j^i}-\underset{\left(
1\right) }{N_j^m}d\underset{\left( 1\right) }{M_m^i}\}\wedge
dy^{\left( \alpha -2\right) j}+d\underset{\left( 1\right)
}{M_j^i}\wedge dy^{\left( \alpha \right) j}.
\end{array}
\tag{6.5.10}
\end{equation}

Identifying to $d\delta y^{\left( \alpha \right) i}$ for (6.5.6) we have

\begin{equation}
\begin{array}{l}
\underset{(\alpha 0)}{P_j^i}=d\underset{\left( \alpha \right)
}{M_j^i}- \underset{\left( 1\right) }{N_j^m}d\underset{\left(
\alpha -1\right) }{ M_m^i}-...-\underset{\left( \alpha -1\right)
}{N_j^m}d\underset{\left(
1\right) }{M_m^i}, \\
\\
\underset{(\alpha 1)}{P_j^i}=d\underset{\left( \alpha -1\right) }{M_j^i}%
- \underset{\left( 1\right) }{N_j^m}d\underset{\left( \alpha
-2\right) }{ M_m^i}-...-\underset{\left( \alpha -2\right)
}{N_j^m}d\underset{\left(
1\right) }{M_m^i}, \\
.................................................................. \\
\underset{\left( \alpha ,\alpha -1\right)
}{P}=d\underset{(1)}{M_j^i}.
\end{array}
\tag{6.5.11}
\end{equation}

Taking into account that, with respect to the adapted basis the 1-form
$d\underset{ (\alpha )}{M_j^i}$ are given by
\begin{equation}
d\underset{(\alpha )}{M_j^i}=\displaystyle\frac{\delta
\underset{(\alpha )}{M_j^i}}{\delta x^h}\delta
x^h+\underset{\gamma =1}{\stackrel{k-1}{\sum }
}\displaystyle\frac{\delta \underset{(\alpha )}{M_j^i}}{\delta
y^{\left( \gamma \right) h}}\delta y^{\left( \gamma \right)
h}+\displaystyle\frac{ \delta \underset{(\alpha )}{M_j^i}}{\delta
p_h}\delta p_h.  \tag{6.5.12}
\end{equation}

Substituting (6.5.11) and identifying to (6.5.8) we completely determine the
coefficients of the Lie brackets.

For instance, from $P_{\alpha ,\alpha -1}$, (6.5.8) $_\alpha $) and $\underset{%
\left( \alpha ,\alpha -1\right) }{P}$ (6.5.11) we obtain
\[
\underset{\left( \alpha -1,\alpha \right)
}{B_{mj}^i}=-\displaystyle\frac{ \delta \underset{(\alpha
)}{M_j^i}}{\delta x^m};\stackrel{\left( \alpha
\right) }{\underset{\left( \gamma ,\alpha -1\right) }{C_{mj}^i}}=- %
\displaystyle\frac{\delta \underset{(\alpha )}{M_j^i}}{\delta
y^{(\gamma )h}};\underset{\left( \alpha -1,\alpha \right)
}{C_j^{im}}=-\displaystyle \frac{\delta \underset{(\alpha
)}{M_j^i}}{\delta p_m}.
\]

\section{The Almost Product Structure $\Bbb{P}$. The Almost Contact
Structure $\Bbb{F}$}

The $\mathcal{F}(T^{*k}M)$-linear mapping $\Bbb{P}:\mathcal{X}
(T^{*p}M)\rightarrow \mathcal{X}(T^{*p}M)$ defined by

\begin{equation}
\Bbb{P}(X^H)=X^H,\Bbb{P}(X^{V_\alpha })=-X^{V\alpha },\Bbb{P}
(X^{W_k})=-X^{W_k},(\alpha =1,...,k-1)  \tag{6.6.1}
\end{equation}
determines an almost product structure on the manifold $T^{*k}M.$ It is
given by means of a nonlinear connection $N.$

We have

\begin{equation}
\begin{array}{l}
\Bbb{P}\circ \Bbb{P}=I, \\
\Bbb{P}=I-2(v_1+...+v_{k-1}+w_k), \\
{\rm rank} \Bbb{P}=(k+1)n.
\end{array}
\tag{6.6.2}
\end{equation}

\begin{teo}
A nonlinear connection $N$ on $T^{*k}M$ is characterized by the existence of
an almost product structure $\Bbb{P}$ on $T^{*k}M$ whose eingenspaces
corresponding to the eingenvalue $-1$ coincides with the linear space of the
vertical distribution $V$ on $T^{*k}M.$
\end{teo}

The proof is same as in the case $k=2$ (see the book [115]).

\begin{teo}
The almost product structure $\Bbb{P},$ defined by (6.6.1) is integrable if
and only if the horizontal distribution $N$ is integrable.
\end{teo}

The proof is exactly as in Prop. 9.8.1 of the book [115].

Another important structure on $T^{*k}M$ is determined by the $\mathcal{F}(T^{*k}M)-$ linear
mapping

\begin{equation}
\begin{array}{c}
\Bbb{F}(\displaystyle\frac \delta {\delta x^i})=-\displaystyle\frac \partial
{\partial y^{\left( k-1\right) i}},
\Bbb{F}(\displaystyle\frac \delta {\delta y^{\left( \alpha \right) i}})=0,
(\alpha =1,..,k-2),\\
 \Bbb{F}(\displaystyle\frac \partial
{\partial y^{\left( k-1\right) }})=\displaystyle\frac \delta {\delta x^i},\Bbb{F}(\displaystyle\frac \delta {\delta p_i})=0,  \tag{6.6.3}
\end{array}
\end{equation}
where $(\displaystyle\frac \delta {\delta x^i}...,\displaystyle\frac \delta
{\delta p_i})$ is the adapted basis of a nonlinear connection $N$ and of the
vertical distribution $V.$ Now, is not difficult to prove

\begin{teo}
The mapping $\Bbb{F}$ has the following properties:

1$^{\circ }$ $\Bbb{F}$ is globally defined on $\widetilde{T^{*k}M}.$

2$^{\circ }$ $\Bbb{F}$ is a tensor field of type $(1,1):$

\begin{equation}
F=-\displaystyle\frac \partial {\partial y^{\left( k-1\right) i}}\otimes
dx^i+\displaystyle\frac \delta {\delta x^i}\otimes \delta y^{\left(
k-1\right) i},  \tag{6.6.4}
\end{equation}

3$^{\circ }$ $Ker\Bbb{F}=N_1\oplus ...\oplus N_{k-2}\oplus W_k,\ Im\Bbb{F}%
=N_0\oplus V_{k-1}$,

4$^{\circ }$ ${\rm rank} \Bbb{F}=2n$,

5$^{\circ }$ $\Bbb{F}^3+\Bbb{F}=0.$
\end{teo}

Looking at 5$^{\circ }$, we can say that $\Bbb{F}$ is an \textit{almost }$(k-1)n$\textit{-contact structure} determined by the nonlinear
connection $N$.

The Nijenhuis tensor of structure $\Bbb{F}$ is expressed by:

\[
\mathcal{N}_{\Bbb{F}}(X,Y)=\Bbb{F}^2[X,Y]+[\Bbb{F}X,\Bbb{F}Y]-\Bbb{F[F}X,Y]-
\Bbb{F[}X,\Bbb{F}Y]
\]
and the condition of normality of the structure $\Bbb{F}$ is as follows:

\begin{equation}
\begin{array}{lr}
\mathcal{N}_{\Bbb{F}}(X,Y)+\underset{i=1}{\stackrel{n}{\sum
}}[\underset{ \alpha =1}{\stackrel{k-1}{\sum }}d(\delta y^{\left(
\alpha \right)
i}(X,Y)+d(\delta p_i)(X,Y)]=0, & \\
 & \forall X,Y\in \mathcal{X}(T^{*k}M).
\end{array}
\tag{6.6.5}
\end{equation}

Using the formulas (6.4.2) and (6.4.5) we can obtain the explicit form of the
last equation.

\section{The Riemannian Structure $\Bbb{G}$ on $T^{*k}M$}

Let $\Bbb{G}$ be a Riemannian structure on the manifold $T^{*k}M.$ $\Bbb{G}$
determines uniquely a nonlinear connection $N$ on $T^{*k}M.$ $N$ is the
ortoghonal distribution, with respect to $\Bbb{G}$, to the vertical
distribution $V.$

In the case when the base manifold $M$ is paracompact, Theorem 4.1.2 affirms
that the manifolds $T^{*k}M$ is paracompact, too. So on $T^{*k}M$ there
exists a Riemannian structure $\Bbb{G}.$ Consequently, we have:

\begin{teo}

If the base manifold $M$ is paracompact then on the manifold $T^{*k}M$ there
exist nonlinear connections $N.$
\end{teo}

Let $\Bbb{G}$ be a Riemannian structure on $T^{*k}M$ and $N$ the
nonlinear connection, whose distribution is orthogonal to the
vertical distribution $V$. The problem is to determine the local
coefficients $\underset{\left( 1\right) }{N_j^i},...,$
$\underset{\left( k-1\right) }{N_j^i},N_{ij}$ of $N$ by means
of the local coefficients of $\Bbb{G}$:

\begin{equation}
\begin{array}{c}
\stackrel{\left( 00 \right) }{g_{ij}} = \Bbb{G}(\displaystyle\frac
\partial {\partial x^i},\displaystyle\frac \partial {\partial
x^j}),\stackrel{\left( 0\alpha \right)
}{g_{ij}}=\Bbb{G}(\displaystyle\frac \partial {\partial x^i},
\displaystyle\frac \partial {\partial y^{\left( \alpha \right)
j}}),\\
\stackrel{\left( 0k\right) }{g_i^j}=\Bbb{G}(\displaystyle\frac
\partial {\partial x^i},\displaystyle\frac \partial {\partial
p_j}), ,...,\stackrel{\left( k,k\right) }{g^{ij}}
=\Bbb{G}(\displaystyle\frac
\partial {\partial p_i},\displaystyle\frac \partial {\partial p_j}).
\end{array} \tag{6.7.1}
\end{equation}

It follows that $\stackrel{\left( k-1,k-1\right) }{g_{ij}}$ and $\stackrel{%
\left( k,k\right) }{g^{ij}}$ are d-tensor fields symmetric and positively
defined. They are coefficients of the restrictions of $\Bbb{G}$ to the
distributions $V_{k-1}$ and $W_k.$

The coefficients of the nonlinear connection $N$ enter in the adapted basis $%
\{\displaystyle\frac \delta {\delta x^i}\}$ to the distribution $N=N_0,$ (6.2.4),
in the adapted basis \newline $\{\displaystyle\frac \delta {\delta y^{\left( 1\right) i}}\},...,\{ %
\displaystyle\frac \delta {\delta y^{\left( k-1\right) i}}\}$ to the
distribution $N_1,...,V_{k-1}$ and in the adapted basis $\{\displaystyle
\frac \delta {\delta p_i}\}$ to the vertical distributions $W_k$. They are
uniquely determined by the conditions that each of the distributions $%
\{N_1,...,N_{k-2}\}$ is orthogonal to $V_{k-1}$ with respect to $\Bbb{G}$
and $N_0$ is orthogonal to $V.$

Indeed, the conditions of orthogonality
\[
\Bbb{G}(\displaystyle\frac \delta {\delta y^{\left( k-2\right) i}}, %
\displaystyle\frac \delta {\delta y^{\left( k-1\right) j}})=0
\]
give us
\[
\stackrel{\left( k-2,k-1\right) }{g_{ij}}-\underset{\left( 1\right) }{%
N_i^m }\stackrel{\left( k-1,k-1\right) }{g_{mj}}=0.
\]

But $rank(\stackrel{\left( k-1,k-1\right) }{g_{ij}})=n$. So the previous
equation uniquely gives the coefficients $\underset{\left( 1\right) }{%
N_i^m }.$

The following equations
\[
\Bbb{G}(\displaystyle\frac \delta {\delta y^{\left( k-3\right) i}}, %
\displaystyle\frac \delta {\delta y^{\left( k-1\right) j}})=0
\]
lead to the equations
\[
\stackrel{\left( k-3,k-1\right) }{g_{ij}}-\underset{\left( 1\right) }{%
N_i^m }\stackrel{\left( k-2,k-1\right) }{g_{mj}}-\underset{\left(
2\right) }{ N_i^m}\stackrel{\left( k-1,k-1\right) }{g_{mj}}=0.
\]
These equations uniquely determine the coefficients
$\underset{\left( 2\right) }{N_i^m}$, etc.

Now we prove the following result.

\begin{teo}
Any Riemannian structure $G$ on the manifold $T^{*k}M$ determines on this
manifold a Riemannian almost contact structure $(\stackrel{\circ }{\Bbb{G}},%
\Bbb{F}).$

\end{teo}

\textit{Proof}. The structure $\Bbb{G}$ determine a nonlinear
connection $N$ with the local coefficients $\underset{\left(
1\right) }{(N_i^m},..., \underset{\left( k-1\right)
}{N_i^m},N_{ij}).$ With these coefficients we construct the adapted basis

$(\displaystyle\frac \delta {\delta x^i},\displaystyle\frac \delta {\delta
y^{\left( 1\right) i}},...,\displaystyle\frac \delta {\delta y^{\left(
k-1\right) j}},\displaystyle\frac \delta {\delta p_j})$ to $N$ and $V_1.$
The restrictions of $G$ to the distribution $V_{k-1} $ and $W_k$ give us the
symmetric, positively defined d-tensor fields

\[
g_{ij}=\Bbb{G}(\displaystyle\frac \delta {\delta y^{\left( k-1\right) i}}, %
\displaystyle\frac \delta {\delta y^{\left( k-1\right) j}}),\ h^{ij}=\Bbb{G}%
( \displaystyle\frac \partial {\partial p_i},\displaystyle\frac \partial
{\partial p_j})
\]

Therefore $\Bbb{G}$ determines on the manifold $T^{\ast k}M$ the Riemannian structure
\begin{equation}
\begin{array}{lll}
\stackrel{\circ }{\Bbb{G}} & = &g_{ij}dx^i\otimes
dx^j+g_{ij}\delta y^{\left( 1\right) i}\otimes \delta y^{\left(
1\right) j}+...+ \vspace{3mm}\\
& + & g_{ij}\delta y^{\left( k-1\right)i}\otimes \delta y^{\left( k-1\right) i}+ h^{ij}\delta p_i\otimes \delta p_j.
\end{array} \tag{6.7.3}
\end{equation}

Consider the almost $(k-1)n$-contact structure $\Bbb{F}$ determined by $N.$
It is given by (6.6.3). Consequently, $\Bbb{F}$ is determined only by the
Riemannian structure $\Bbb{G}$.

Therefore the pair $(\stackrel{\circ }{\Bbb{G}},\Bbb{F})$ is a Riemannian almost
contact structure determined only by $\Bbb{G.}$

Of course, the equation
\[
\stackrel{\circ }{\Bbb{G}}(\Bbb{F}X,Y)=-\stackrel{\circ }{\Bbb{G}}(X,\Bbb{F}Y),
\]
is verified on the adapted basis to $N$ and $V.$ q.e.d

\begin{rem}
We shall use the Riemannian structure $\stackrel{\circ }{\Bbb{G}},$ in the
case of Hamilton space of order $k$, for $h^{ij}=g^{ij}$
\end{rem}

\section{The Riemannian Almost Contact Structure $(\stackrel{\vee}{\Bbb{G}},%
\stackrel{\vee}{\Bbb{F}})$}

Consider a d-tensor field $g_{ij}$ on $\widetilde{T^{*k}M}$, symmetric and
positively defined. It follows $\det \left\| g_{ij}\right\| >0$ on $%
\widetilde{T^{*k}M}$. Let $N$ be an apriori given nonlinear connection on
the manifold $T^{*k}M$ and the adapted basis with respect to $N$ and $V:$ $( %
\displaystyle\frac \delta {\delta x^j},\displaystyle\frac \delta {\delta
y^{\left( 1\right) j}},...,\displaystyle\frac \delta {\delta y^{\left(
k-1\right) j}},\displaystyle\frac \delta {\delta p_j}).$

As usual, the dual
basis is $(dx^i\delta y^{\left( 1\right) i},...,\delta y^{\left( k-1\right)
i},\delta p_i).$ The $N$-lift of the d-tensor fields $g_{ij}$ at every point $%
u\in T^{*k}M$ is defined by
\begin{equation}
\begin{array}{l}
\stackrel{\vee}{\Bbb{G}}=g_{ij}dx^i\otimes dx^j+g_{ij}\delta y^{\left(
1\right) i}\otimes \delta y^{\left( 1\right) j}+...+g_{ij}\delta y^{\left(
k-1\right) i}\otimes \delta y^{\left( k-1\right) i}+ \\
\\
\qquad +g^{ij}\delta p_i\otimes \delta p_j,
\end{array}
\tag{6.8.1}
\end{equation}
where $g^{ij}$ is the contravariant tensor of $g_{ij}.$ We have:
\begin{prop}
$\stackrel{\vee}{\Bbb{G}}$ is a Riemannian structure on $T^{*k}M.$
\end{prop}

Indeed, $\stackrel{\vee}{\Bbb{G}}$ is a tensor field on the manifold $%
\widetilde{T^{*k}M}$ covariant symmetric and positively defined.

The geometrical object fields $N$ and $g_{ij}$ allows to define the $%
\mathcal{F}(\widetilde{T^{*k}M})$-linear mapping
\[
\stackrel{\vee}{\Bbb{F}}:\mathcal{X}(\widetilde{T^{*k}M})\rightarrow
\mathcal{X}(\widetilde{T^{*k}M})
\]
by
\begin{equation}
\stackrel{\vee}{\Bbb{F}}(\displaystyle\frac \delta {\delta x^i})=-g_{ij} %
\displaystyle\frac \delta {\delta p_j},\stackrel{\vee}{\Bbb{F}}( %
\displaystyle\frac \delta {\delta y^{\left( \alpha \right) i}})=0,\alpha
=1,...,k,\stackrel{\vee}{\Bbb{F}}(\displaystyle\frac \delta {\delta
p_i})=g^{ij}\displaystyle\frac \delta {\delta x^j}  \tag{6.8.2}
\end{equation}
Taking into account the rule of transformations of the d-tensor $g_{ij}$ and
of the vector fields $\displaystyle\frac \delta {\delta x^j},..., %
\displaystyle\frac \delta {\delta p_j} $ it follows that $\stackrel{\vee}{%
\Bbb{F}}$ has a geometrical meaning.

We have also:

\begin{teo}
1$^{\circ }$ $\stackrel{\vee}{\Bbb{F}}$ is a tensor field on $\widetilde{%
T^{*k}M}$ of type (6.1.1)

2$^{\circ }$ In adapted basis $\stackrel{\vee}{\Bbb{F}}$ is expressed by
\begin{equation}
\stackrel{\vee}{\Bbb{F}}=-g_{ij}\displaystyle\frac \delta {\delta
p_j}\otimes dx^i+g^{ij}\displaystyle\frac \delta {\delta x^i}\otimes \delta
p_j,  \tag{6.8.3}
\end{equation}

3$^{\circ }$ $Ker\stackrel{\vee}{\Bbb{F}}=N_1\oplus ...\oplus N_{(k-1)},\ Im%
\stackrel{\vee}{\Bbb{F}}=N_0\oplus W_k$,

4$^{\circ }$ rank$\stackrel{\vee}{\Bbb{F}}=2n$,

5$^{\circ }$ $\stackrel{\vee}{\Bbb{F}}^3+\stackrel{\vee}{\Bbb{F}}=0$.
\end{teo}

By means of (6.8.2) the proof is immediate.

Consequently, $\stackrel{\vee}{\Bbb{F}}$ is an almost $(k-1)n$-contact
structure on $T^{*k}M$ determined by the nonlinear connection $N$ and
d-tensor $g_{ij}$. Of course, it is useful in the case of Hamilton spaces of
order $k.$

The conditions of normality of the structure $\stackrel{\vee}{\Bbb{F}}$ can
be written exactly as in (6.6.5).

Between the Riemannian structure $\stackrel{\vee}{\Bbb{G}}$ and the
structure $\stackrel{\vee}{\Bbb{F}}$ there is a strongly relation.

\begin{teo}
1$^{\circ }$ The pair $(\stackrel{\vee}{\Bbb{G}}$ ,$\stackrel{\vee}{\Bbb{F}%
}$ $)$ is a Riemannian almost contact structure determined only by the pair $%
(N,g_{ij})$

2$^{\circ }$ The associated 2-form is
\[
\theta =\delta p_i\wedge dx^i.
\]
3$^{\circ }$ If the coefficients $N_{ij}$ of $N$ are symmetric, then $\theta
$ reduces to the presymplectic structure
\[
\theta =dp_i\wedge dx^i.
\]
\end{teo}

\textit{Proof}:

1$^{\circ }$ The following formula
$\stackrel{\vee}{\Bbb{G}}(\stackrel{\vee}{\Bbb{F}}X,Y)=-\stackrel{\vee}{\Bbb{G}}(X,\stackrel{\vee
}{\Bbb{F}}Y)$ can be verified on the adapted basis, using (6.8.2).

2$^{\circ }$ $\theta (X,Y)=\stackrel{\vee}{\Bbb{G}}(\stackrel{\vee}{\Bbb{F}}X,Y)$ is
satisfied, too.

3$^{\circ }$ From $N_{ij}=N_{ji}$ and $\delta
p_i=dp_i-N_{ji}dx^j$ we deduce $\delta p_i\wedge dx^i=dp_i\wedge dx^i,$
q.e.d.

\chapter{Linear Connections on the Manifold $T^{*k}M$}

\markboth{\it{THE GEOMETRY OF HIGHER-ORDER HAMILTON SPACES\ \ \ \ \ }}{\it{Linear Connection on the Manifold} $T^{*k}M$}

In the chapter 10 of the book 'The Geometry of Hamilton and
Lagrange Spaces', [115], we studied the notions of $d$-tensor
algebra and $N$ -linear connections on $T^{*2}M$. The
corresponding theory will be extended in the present section to
the case $k>2.$

\section{The Algebra of Distinguished Tensor Fields}

Let $N$ be a nonlinear connection on the manifold $T^{*k}M.$ It determines,
at every point $u\in T^{*k}M,$ the direct decomposition of the linear space $%
T_u(T^{*k}M)$:

\begin{equation}
T_u(T^{*k}M)=N_{0,u}\oplus N_{1,u}\oplus ...\oplus N_{k-2,u}\oplus
V_{k-1, u}\oplus W_{k,u}.  \tag{7.1.1}
\end{equation}

A vector field $X\in \mathcal{X}(T^{*k}M)$ and an one form $\omega \in
\mathcal{X}^{*}(T^{*k}M)$ can be uniquely written in the form

\begin{equation}
\left\{
\begin{array}{c}
X=X^H+X^{V_1}+...+X^{V_{k-1}}+X^{W_k}, \\
\\
\omega =\omega ^H+\omega ^{V_1}+...+\omega ^{V_{k-1}}+\omega ^{W_k}.
\end{array}
\right.  \tag{7.1.2}
\end{equation}

Clearly if $h,v_1,v_2,...,v_{k-1},w_k$ are the projectors determined by the
decomposition (7.1.1), we have

\begin{equation}
\begin{array}{c}
X^H=hX,\ X^{V_\alpha }=v_\alpha X,\ (\alpha =1,...,k-1),\ X^{W_k}=w_kX, \\
\\
\omega ^H=\omega \circ h,\ \omega ^{V_\alpha }=\omega \circ v_\alpha ,\
(\alpha =1,...,k-1),\ \omega ^{W_k}=\omega \circ w_k.
\end{array}
\tag{7.1.3}
\end{equation}

\begin{defi}
A distinguished tensor field (briefly $d$-tensor field) on the manifold $%
T^{*k}M$ of type $(r,s)$ is a $d$-tensor field $T$ of type $(r,s)$ on $%
T^{*k}M$ with the property:
\begin{equation}
\begin{array}{l}
T(\stackrel{1}{\omega },...,\stackrel{r}{\omega },.\underset{1}{X},...,%
\underset{s}{X})=T(\stackrel{1}{\omega }^H,...,\stackrel{1}{\omega }%
^{W_k},.\underset{1}{X}^H,...,\underset{s}{X^{W_k}}) \\
\\
\forall \stackrel{1}{\omega },...,\stackrel{r}{\omega }\in \mathcal{X}%
^{*}(T^{*k}M),\forall \underset{1}{X},...,\underset{s}{X}\in \mathcal{X}%
(T^{*k}M)
\end{array}
\tag{7.1.4}
\end{equation}
\end{defi}

For instance, every components $X^H,X^{V_1},X^{V_{k-1}},X^{W_k}$ of a vector field $X$ is a $d$
-vector field and every component of an 1-form $\omega ^H,\omega ^{V_1},\omega
^{V_{k-1}},\omega ^{W_k}$ is a $d$-1-form field.

In the adapted basis $(\displaystyle\frac \delta {\delta x^i},\displaystyle
\frac \delta {\delta y^{\left( 1\right) i}},...,\displaystyle\frac \delta
{\delta y^{\left( k-1\right) i}},\displaystyle\frac \delta {\delta p_i})$ to
the decomposition (7.1.1) and in its dual basis $(dx^i,\delta
y^{(1)i},...,\delta y^{\left( k-1\right) i},\delta p_i)$, given by (4.1.10)
and (4.2.2), a $d$-tensor field $T$ of type $(r,s)$ can be written in the
form

\begin{equation}
T\left( u\right) =T_{j_1...j_s}^{i_1...i_r}(u)\displaystyle\frac \delta
{\delta x^{i_1}}\otimes ...\otimes \displaystyle\frac \delta {\delta
p_{j_s}}\otimes dx^{j_1}\otimes ...\otimes \delta p_{i_r},\forall u\in
T^{*k}M  \tag{7.1.5}
\end{equation}

It follows that the set $\{1,\displaystyle\frac \delta {\delta x^i}, %
\displaystyle\frac \delta {\delta y^{\left( 1\right) i}},...,\displaystyle
\frac \delta {\delta y^{\left( k-1\right) i}},\displaystyle\frac \delta
{\delta p_i}\}$ generates the algebra of the $d$-tensor fields over the ring
of functions $\mathcal{F}(T^{*k}M).$

With respect to a local transformation of the coordinates on $T^{*k}M,$ the
local coefficients $T_{j_1...j_s}^{i_1...i_r}$ of $T$ are transformed by the
classical rule

\begin{equation}
\widetilde{T}_{j_1...j_s}^{i_1...i_r}=\displaystyle\frac{\partial \widetilde{
x}^{i_1}}{\partial x^{h_1}}...\displaystyle\frac{\partial \widetilde{x}%
^{i_r} }{\partial x^{h_r}}\displaystyle\frac{\partial x^{k_1}}{\partial
\widetilde{x }^{j_1}}...\displaystyle\frac{\partial x^{k_s}}{\partial
\widetilde{x}^{j_s}} T_{k_1...k_s}^{h_1...h_r}.  \tag{7.1.6}
\end{equation}

For instance, if $f\in {\cal F}(T^{*k}M),$ then each set of
functions $\displaystyle\frac{\delta f}{\delta
x^i},\displaystyle\frac{\delta f}{\delta y^{\left( 1\right)
i}},...,$ $\displaystyle\frac{\delta f}{\delta y^{\left(
k-1\right) i}},$ $(i=1,...,n)$ is a $d$-covector field, and
$\displaystyle \frac{\delta f}{\delta p_i}$ is a $d$-vector field.

\section{$N$-Linear Connections}

The $N$-linear connections, for $k=2,$ were studied in chapter 10, of the book
[115]. Their definition for $k\geq 2$ is as follows:

\begin{defi}
A linear connection $D$ on the manifold $T^{*k}M$ is called an $N$-linear
connection if:

(1) $D$ preserves by parallelism the distributions $%
N_0,N_1,...,N_{k-2}, V_{k-1}, W_k.$

(2) The $k-1$ tangent structure $J$ is absolutely parallel with respect to $D.$

(3) The presymplectic structure $\theta $ is absolutely parallel with respect $%
D$.
\end{defi}

Directly from the definition we can establish, without difficulties the
following characterization of an $N$-linear connection

\begin{teo}
A linear connection $D$ is an $N$-linear connection on the manifold $T^{*k}M$
if and only if:

(1) $D$ preserves by parallelism the distributions $%
N_0,N_1,...,N_{k-2},V_{k-1}$ and $W_k.$

(2)
\begin{equation}
\begin{array}{l}
D_X(JY^H)=J(D_XY^H),\ D_X(JY^{V_\alpha })=J(D_XY^{V_\alpha }),\  \\

\qquad \qquad \qquad (\alpha =1,...,k-1) \\
D_X(JY^{W_k})=J(D_XY^{W_k}),\ \forall X,Y\in X(T^{*k}M).
\end{array}
\tag{7.2.1}
\end{equation}

(3)
\begin{equation}
D\theta =0 . \tag{7.2.2}
\end{equation}
\end{teo}

We remark that the equalities $D_X(JY^{V_{k-1} })=J(D_XY^{V_{k-1} })$ and $%
D_X(JY^{W_k})=J(D_XY^{W_k})$ are trivial, since $J(Y^{V_{k-1} })=0,$ and $%
J(Y^{W_k})=0.$

We obtain also

\begin{teo}
For any $N$- linear connection $D$ we have
\begin{equation}
D_Xh=D_Xv_\alpha =D_Xw_k=0,(\alpha =1,...,k-1),  \tag{7.2.3}
\end{equation}
\begin{equation}
D_X\Bbb{P}=0,D_X\Bbb{F}=0.  \tag{7.2.4}
\end{equation}
\end{teo}
Indeed, from $(D_Xh)(Y)=D_X(hY)-hD_XY$ if $Y=Y^H$ we obtain
\[
(D_Xh)(Y^H)=D_XY^H-hD_XY^H=0
\]
and for $Y=Y^{V_\alpha },$
\[
(D_Xh)(Y^{V_\alpha })=D_X(hY^{V_\alpha })-hD_XY^{V_\alpha }=0.
\]
Also $(D_Xh)(Y^{W_k})=0.$ That means $(D_Xh)(Y)=0,\forall Y\in \chi ((T^{*k}M)$.
Similarly, we prove  the other equalities (7.2.3).

Now, taking into account the formula (6.6.2), we deduce $D_X\Bbb{P}=0$. The last equality (7.2.4) can be proved by the formula $(D_X\Bbb{F})(Y)=D_X\Bbb{F}(Y)-\Bbb{F}(D_XY),$ using the local expression of $D_XY$ (see the section 5, from this chapter).

Let us consider the vector field $X$ written in the form (7.1.2). From the linearity of the operator $D_XY$ with respect to $X$ we deduce
\begin{equation}
D_XY=D_{X^H}Y+D_{X^{V_1}}Y+...+D_{X^{V_{k-1}}}Y+D_{X^{W_k}}Y.  \tag{7.2.5}
\end{equation}

Here appears $(k+1)$ new operators of derivation in the $d-$tensor algebra, defined by
\begin{equation}
D_X^H=D_{X^H},D_X^{V_1}=D_{X^{V_1}},...,D_X^{V_{k-1}}=D_{X^{V_{k-1}}},D_X^{W_k}=D_{X^{W_k}}.
\tag{7.2.6}
\end{equation}

These operators are not covariant derivations in the $d$-tensor algebra, since
$$
\begin{array}{l}
D_X^Hf=X^Hf\not =Xf,D_X^{V_\alpha }f=X^{V_\alpha }f\not =Xf,(\alpha
=1,...,k),\\
D_X^{W_k}f=X^{W_k}f\not =Xf.
\end{array}
$$
However, they have similar properties with the covariant derivatives.

>From (7.2.5) and (7.2.6) we deduce
\begin{equation}
D_XY=D_X^HY+D_X^{V_1}Y+...+D_X^{V_{k-1}}V+D_X^{W_k}Y.  \tag{7.2.7}
\end{equation}

By means of Theorem 7.2.2, we have:

\begin{teo}
The operators $D_X^H,D_X^{V_1},...,D_X^{V_{k-1}}$ and $D_X^{W_k}$ have the
following properties:

1) Each operator $D_X^H,D_X^{V_1},...,D_X^{V_{k-1}}$ and $D_X^{W_k}$ maps
a vector field that belongs to one of the distributions $N_0,N_1,...,N_{k-2}, V_{k-1}
$ and $W_k$ into a vector field that belongs to the same distribution,

\qquad

2) $D_X^Hf=X^Hf,\ D_X^{V_\alpha }f=X^{V_\alpha }f,\ (\alpha =1,...,k),\
D_X^{W_k}f=X^{W_k}f,$

\qquad

3) $D_X^H(fY)=X^H(fY)+fD_X^HY,\ $ $D_X^{V_\alpha }(fY)=X^{V_\alpha
}(fY)+fD_X^{V_\alpha }Y,$

\qquad

$(\alpha =1,...,k-1),$ $D_X^{W_k}(fY)=X^{W_k}(fY)+fD_X^{W_k}Y,$

\qquad

4) $D_X^H(Y+Z)=D_X^HY+D_X^HZ,\ D_X^{V_\alpha }(Y+Z)=D_X^{V_\alpha
}Y+D_X^{V_\alpha }Z,$

\qquad

$(\alpha =1,...,k-1),$ $D_X^{W_k}(Y+Z)=D_X^{W_k}Y+D_X^{W_k}Z,$

\qquad

5) $D_{X+Y}^H=D_X^H+D_Y^H,\ D_{X+Y}^{V_\alpha }=D_X^{V_\alpha
}+D_Y^{V_\alpha },\ (\alpha =1,...k-1),$

\qquad

$D_{X+Y}^{W_k}=D_X^{W_k}+D_Y^{W_k},$

\qquad

6) $D_{fX}^H=fD_X^H;D_{fX}^{V_\alpha }=fD_Y^{V_\alpha },(\alpha
=1,...,k-1),D_{fX}^{W_k}=fD_Y^{W_k},$

\qquad

7) $D_X^H(JY)=JD_X^HY,$ $D_X^{V_\alpha }(JY)=JD_X^{V_\alpha }Y,\ (\alpha
=1,...,k-1)$

\qquad

$D_X^{W_k}(JY)=JD_X^{W_k}Y,$

\qquad

8) $D_X^H\theta =0,D_X^{V_\alpha }\theta =0,(\alpha =1,...,k-1),D_X^{W_k}\theta =0.$

\qquad

9) For any open set $U\subset T^{*k}M$ the following properties hold:

\quad

$(D_X^HY)_{\left| U\right. }=D_{X\left| _U\right. }^HY_{\left| U\right.
},(D_X^{V_\alpha }Y)_{\left| U\right. }=D_{X\left| _U\right. }^{V_\alpha
}Y_{\left| U\right. },(\alpha =1,...,k-1),$

\quad

$(D_X^{W_k}Y)\left| _U\right. =D_{X\left| _U\right. }^{W_k}Y_{\left|
U\right. }.$
\end{teo}
The proof of this theorem can be done by the classical methods, [115].

The operators $D_X^H,D_X^{V_\alpha },D_X^{W_k}$ will be called the operators
of $h$- , $v_\alpha $- and $w_k$- covariant derivation.

The action of these operators over the 1-form fields $\omega $ are given by
\begin{equation}
\begin{array}{l}
(D_X^H\omega )(Y)=X^H\omega (Y)-\omega (D_X^HY), \\
\\
(D_X^{V_\alpha }\omega )(Y)=X^{V_\alpha }\omega (Y)-\omega (D_X^{V_\alpha
}Y),(\alpha =1,...k-1), \\
\\
(D_X^{W_k}\omega )(Y)=X^{W_k}\omega (Y)-\omega (D_X^{W_k}Y).
\end{array}
\tag{7.2.8}
\end{equation}

Of course, the action of the operators $D_X^H,D_X^{V_\alpha },D_X^{W_k}$ can
be extended to any tensor fields, particularly to any d-tensor field on $%
T^{*k}M.$

For instance, if the d-tensor $T$ verifies (7.1.4) we have
\begin{equation}
\begin{array}{l}
(D_X^HT)(\stackrel{1}{\omega }^H,...,\stackrel{1}{\omega }^{W_k},.
\underset{1}{X^H},...,\underset{s}{X^{W_k}})=X^HT(\stackrel{1}{\omega
} ^H,...,\stackrel{1}{\omega
}^{W_k},\underset{1}{X}^H,...,\underset{s}{
X^{W_k}}) \\
\\
\quad -T(D_X^H\stackrel{1}{\omega }^H,...,\stackrel{1}{\omega }^{W_k},
\underset{1}{X^H},...,\underset{s}{X^{W_k}})-...- \\
\\
\quad -T(\stackrel{1}{\omega ^H},...,\stackrel{r}{\omega ^{W_k}},.
\underset{1}{X^H},...,D_X^H\underset{s}{X^{W_k}}).
\end{array}
\tag{7.2.9}
\end{equation}

Now, let us consider a parametrized smooth curve $\gamma :t\in I\rightarrow
\gamma (t)\in \widetilde{T^{*k}M}$ having the image in a domain of local chart.

Its tangent vector field $\stackrel{\cdot }{\gamma }=\displaystyle\frac{
d\gamma }{dt}$ is uniquely written in the form

\begin{equation}
\stackrel{\cdot }{\gamma }=\stackrel{\cdot }{\gamma }^H+\stackrel{\cdot }{
\gamma }^{V_1}+...+\stackrel{\cdot }{\gamma }^{V_{k-1}}+\stackrel{\cdot }{
\gamma }^{W_k} . \tag{7.2.10}
\end{equation}

If the curve $\gamma $ is analytically given by (6.3.2), then $\stackrel{
\cdot }{\gamma },$ in the adapted basis is given by (6.3.20''), (6.3.21).
The horizontal curves are defined by Theorem 6.3.2 and the autoparallel curves
of the nonlinear connection $N$ are given by the equations (6.3.22) in the
conditions
\[
y^{(1)i}=\displaystyle\frac{dx^i}{dt},...,y^{(k-1)i}=\displaystyle\frac
1{(k-1)!}\displaystyle\frac{d^{k-1}x^i}{dt^{k-1}}.
\]

A vector field $Y(\gamma(t))$, along curve $\gamma $ has the
covariant derivative
\begin{equation}
D_{\stackrel{\cdot }{\gamma }}Y=D_{\stackrel{\cdot }{\gamma }}^HY+D_{%
\stackrel{\cdot }{\gamma }}^{V_1}Y+...+D_{\stackrel{\cdot }{\gamma }
}^{V_{k-1}}Y+D_{\stackrel{\cdot }{\gamma }}^{W_k}Y.  \tag{7.2.11}
\end{equation}

The vector field $Y(\gamma(t))$ is \textit{parallel} along curve $%
\gamma $ if
\begin{equation}
D_{\stackrel{\cdot }{\gamma }}Y=0.  \tag{7.2.12}
\end{equation}

The curve $\gamma $ is autoparallel with respect to the $N$- linear
connection $D$ if $D_{\stackrel{\cdot }{\gamma }}\stackrel{\cdot }{\gamma }=0$.
This equation is equivalent to
\begin{equation}
D_{\stackrel{\cdot }{\gamma }}^H\stackrel{\cdot }{\gamma }+D_{\stackrel{
\cdot }{\gamma }}^{V_1}\stackrel{\cdot }{\gamma }+...+D_{\stackrel{\cdot }{
\gamma }}^{V_{k-1}}\stackrel{\cdot }{\gamma }+D_{\stackrel{\cdot }{\gamma }
}^{W_k}\stackrel{\cdot }{\gamma }=0  \tag{7.2.13}
\end{equation}

In a next section we shall express this equation in an adapted basis.

\section{The Torsion and Curvature of an $N$-Linear Connection}

The torsion tensor field $\Bbb{T}$ of an $N$-linear connection $D$ is
expressed as usually by

\begin{equation}
\Bbb{T}(X,Y)=D_XY-D_YX-[X,Y],\forall X,Y\in \mathcal{X}(T^{*k}M).  \tag{7.3.1}
\end{equation}
$\Bbb{T}$ can by characterized by the vector fields
\begin{eqnarray*}
\Bbb{T}(X^H,Y^H)
&=&h\Bbb{T}(X^H,Y^H)+v_1\Bbb{T}(X^H,Y^H)+...+v_{k-1}\Bbb{T}
(X^H,Y^H)+\\
& & \\
&+&w_k\Bbb{T}(X^H,Y^H),
\end{eqnarray*}

\begin{eqnarray*}
\Bbb{T}(X^H,Y^{V_\alpha}) & = & h \Bbb{T}(X^H,Y^{V_\alpha}) +
v_1\Bbb{T}
(X^H,Y^{V_\alpha })+...+v_{k-1}\Bbb{T}(X^H,Y^{V_\alpha })+ \\
& & \\
&+&w_k\Bbb{T}(X^H,Y^{V_\alpha}),
\end{eqnarray*}

\begin{equation}
\begin{array}{lll}
\Bbb{T}(X^H,Y^{W_k})
&=&h\Bbb{T}(X^H,Y^{W_k})+v_1\Bbb{T}(X^H,Y^{W_k})+...+\\
& & \\
&+&v_{k-1}\Bbb{T}(X^H,Y^{W_k})+w_k\Bbb{T}(X^H,Y^{W_k}),
\end{array} \tag{7.3.2}
\end{equation}
\begin{eqnarray*}
\Bbb{T}(X^{V_\alpha },Y^{V_\beta }) &=&h\Bbb{T}(X^{V_\alpha },Y^{V_\beta
})+v_1\Bbb{T}(X^{V_\alpha },Y^{V_\beta })+...+v_{k-1}\Bbb{T}(X^{V_\alpha
},Y^{V_\beta })+ \\
& & \\
& + & w_k\Bbb{T}(X^{V_\alpha },Y^{V_\beta }),\quad
\end{eqnarray*}

\[
\alpha \leq \beta ;\alpha ,\beta =1,2,...,k-1,
\]

\begin{eqnarray*}
\Bbb{T}(X^{V_\alpha },Y^{W_k}) &=&h\Bbb{T}(X^{V_\alpha },Y^{W_k})+v_1\Bbb{T}
(X^{V_\alpha },Y^{W_k})+...+v_{k-1}\Bbb{T}(X^{V_\alpha },Y^{W_k})+ \\
& & \\
& + & w_k\Bbb{T}(X^{V_\alpha },Y^{W_k}),
\end{eqnarray*}

\begin{eqnarray*}
\Bbb{T}(X^{W_k},Y^{W_k}) = h\Bbb{T}(X^{W_k},Y^{W_k})+v_1\Bbb{T}
(X^{W_k},Y^{W_k})+...+ \\
+v_{k-1}\Bbb{T}(X^{W_k},Y^{W_k})+ w_k\Bbb{T}(X^{W_k},Y^{W_k}).
\end{eqnarray*}

Since $D$ preserves by parallelism the distributions $%
N_0,N_1,...,N_{k-2},V_{k-1}$ and $W_k$ we deduce

\begin{prop}
The following properties of the torsion $\Bbb{T}$ of the $N$-linear
connection $D$ hold:
\begin{equation}
hT(X^{V_{k-1}},Y^{V_{k-1}})=0,hT(X^{W_k},X^{W_k})=0. \tag{7.3.3}
\end{equation}
\end{prop}

Indeed, the distributions $V_{k-1}$ and $W_k$ being integrable, the equations
(7.3.3) are verified.

Now, using the formula (7.3.1), we can write the expression of $\Bbb{T}(X,Y)$
in the form
\begin{equation}
\Bbb{T}(X,Y)=h\Bbb{T}(X,Y)+v_1\Bbb{T}(X,Y)+...+v_{k-1}\Bbb{T}(X,Y)+w_k\Bbb{T}
(X,Y) , \tag{7.3.4}
\end{equation}
taking into account the components of $X$ and $Y$ from the decomposition
(7.1.2).

The curvature of the $N$-linear connection $D$ is given by
\begin{equation}
\Bbb{R}(X,Y)Z=(D_XD_Y-D_YD_X)Z-D_{[X,Y]}Z,\ \forall X,Y,Z\in \mathcal{X}
(T^{*k}M)  \tag{7.3.5}
\end{equation}

Taking into account the decomposition of vector fields $X,Y,Z$ in the form
(7.1.2) we can write the curvature $\Bbb{R}$ in a similar manner as the torsion $%
\Bbb{T}$.

The definition 7.2.1 allows to prove, without difficulties
\begin{prop}
The following properties hold:
\begin{equation}
J(\Bbb{R}(X,Y)Z)=\Bbb{R}(X,Y)JZ;\ D_X\theta =0.  \tag{7.3.5a}
\end{equation}
\end{prop}
Consequently, we have

\begin{teo}
The curvature $\Bbb{R}$ has the properties:

1$^{\circ }$ The essential components of $R$ are
\begin{equation}
\Bbb{R}(X,Y)Z^H,\Bbb{R}(X,Y)Z^{V_\alpha },(\alpha =1,...k-1),\Bbb{R}%
(X,Y)Z^{W_k}.  \tag{7.3.6}
\end{equation}

2$^{\circ }$ The vector field $\Bbb{R}(X,Y)Z^H$ belongs to the horizontal
distribution $N=N_0$.

3$^{\circ }$ The vector field $\Bbb{R}(X,Y)Z^{V_\alpha },(\alpha =1,...k-1)$
belongs to the distribution $N_\alpha .$

4$^{\circ }$ The vector field $\Bbb{R}(X,Y)Z^{W_k\text{ }}$belongs to the
distribution $W_k$.

5$^{\circ }$ The following equations hold:
\begin{equation}
\begin{array}{l}
v_\alpha \{\Bbb{R}(X,Y)Z^H\}=0,(\alpha =1,...k-1),w_k\{\Bbb{R}\left(
X,Y\right) Z^H\}=0, \\
\\
v_\alpha \{\Bbb{R}(X,Y)Z^{V_\beta }\}=0,(\alpha \not =\beta ;\alpha ,\beta
=1,...k-1), \\
\\
v_\alpha \{\Bbb{R}\left( X,Y\right) Z^{W_k}\}=0, \\
\\
h\{\Bbb{R}(X,Y)Z^{V_\alpha }\}=0,(\alpha =1,...k-1),h\{\Bbb{R}\left(
X,Y\right) Z^{W_k}\}=0.
\end{array}
\tag{7.3.7}
\end{equation}
\end{teo}

Of course we can express the $d$-tensor of curvature by means of the \newline
operators $D_X^H,D_X^{V_\alpha },D_X^{W_k}.$ They will be written in the
adapted basis in a next section.
\begin{prop}
The Ricci identities of the $N$-linear connection $D$ are:
\end{prop}
\[
\lbrack D_X,D_Y]Z^H=\Bbb{R}(X,Y)Z^H+D_{[X,Y]}Z^H,
\]
\begin{equation}
\lbrack D_X,D_Y]Z^{V_\alpha }=\Bbb{R}(X,Y)Z^{V_\alpha }+D_{[X,Y]}Z^{V_\alpha
},(\alpha =1,...,k-1) , \tag{7.3.8}
\end{equation}
\[
\lbrack D_X,D_Y]Z^{W_k}=\Bbb{R}(X,Y)Z^{W_k}+D_{[X,Y]}Z^{W_k}.
\]

Let us consider the Liouville vector fields $\stackrel{1}{\Gamma },...,$ $%
\stackrel{k-1}{\Gamma }$ and $C^{*}$ from Theorem 4.2.1 and let us apply the
previous Proposition.

\begin{teo}
For any $N$-linear connection $D$ the following identities hold:
\end{teo}
\begin{equation}
\begin{array}{l}
\lbrack D_X,D_Y]\stackrel{\alpha }{\Gamma }=\Bbb{R}(X,Y)\stackrel{\alpha }{
\Gamma }+D_{[X,Y]}\stackrel{\alpha }{\Gamma },(\alpha =1,...,k-1), \\
\\
\lbrack D_X,D_Y]C^{*}=\Bbb{R}(X,Y)C^{*}+D_{[X,Y]}C^{*}.
\end{array}
\tag{7.3.9}
\end{equation}

Using the considerations from this chapter we can establish the Bianchi
identities of an $N$-linear connection $D,$ by means of the operators $%
D_X^H, $ $D_X^{V_\alpha },$ $D_X^{W_k}$ taking into account the classical
identities
\begin{equation}
\begin{array}{l}
\underset{(X,Y,Z)}{\sum
}\{(D_X\Bbb{T})(Y,Z)-\Bbb{R}(X,Y)Z+\Bbb{T}(\Bbb{T}
(X,Y),Z)\}=0, \\
\\
\underset{(X,Y,Z)}{\sum
}\{(D_X\Bbb{R})(U,Y,Z)-\Bbb{R}(\Bbb{T}(X,Y),Z)U\}=0,
\end{array}
\tag{7.3.10}
\end{equation}
where $\underset{(X,Y,Z)}{\sum }$ means the cyclic sum.

\section{The Coefficients of a $N$-Linear Connection}

A $N$-linear connection $D$ is characterized by its coefficients in the
adapted basis $(\displaystyle\frac \delta {\delta x^i},\displaystyle\frac
\delta {\delta y^{\left( \alpha \right) i}},\displaystyle\frac \delta
{\delta p_i})$. As we shall see these coefficients obey particular rules of
transformations with respect to the change of local coordinates on the
manifold $T^{*k}M.$

Taking into account Theorem 7.2.1 we can prove:
\begin{teo}
We have:

1$^{\circ }$ An $N$-linear connection $D$ can be uniquely represented in the
adapted basis in the following form:
\begin{equation}
\begin{array}{l}
D_{\displaystyle\frac \delta {\delta x^j}}\displaystyle\frac \delta {\delta
x^i}=H_{ij}^s\displaystyle\frac \delta {\delta x^s};D_{\displaystyle\frac
\delta {\delta x^j}}\displaystyle\frac \delta {\delta y^{\left( \alpha
\right) i}}=H_{ij}^s\displaystyle\frac \delta {\delta y^{\left( \alpha
\right) s}},(\alpha =1,...k-1); \\
D_{\displaystyle\frac \delta {\delta x^j}}\displaystyle\frac \delta {\delta
p_i}=-H_{sj}^i\displaystyle\frac \delta {\delta p_s}; \\
\\
D_{\displaystyle\frac \delta {\delta y^{\left( \alpha \right) j}}}%
\displaystyle\frac \delta {\delta x^i}=\underset{(\alpha )}{C_{ij}^s}%

\displaystyle\frac \delta {\delta x^s};D_{\displaystyle\frac
\delta {\delta y^{\left( \alpha \right) j}}}\displaystyle\frac
\delta {\delta y^{(\beta )i}}=\underset{(\alpha
)}{C_{ij}^s}\displaystyle\frac \delta {\delta
y^{(\beta )s}},(\alpha ,\beta =1,...k-1); \\
D_{\displaystyle\frac \delta {\delta y^{\left( \alpha \right) j}}}%
\displaystyle\frac \delta {\delta p_i}=\underset{(\alpha )}{C_{sj}^i}%
\displaystyle\frac \delta {\delta p_s}; \\
\\
D_{\displaystyle\frac \delta {\delta p_j}}\displaystyle\frac \delta {\delta
x^i}=C_i^{js}\displaystyle\frac \delta {\delta x^s};D_{\displaystyle\frac
\delta {\delta p_j}}\displaystyle\frac \delta {\delta y^{(\alpha
)i}}=C_i^{js}\displaystyle\frac \delta {\delta y^{(\alpha )s}},(\alpha
=1,...k-1); \\
D_{\displaystyle\frac \delta {\delta p_i}}\displaystyle\frac \delta {\delta
p_j}=-C_s^{ij}\displaystyle\frac \delta {\delta p_s}.
\end{array}
\tag{7.4.1}
\end{equation}

2$^{\circ }$ With respect to the transformation (4.1.2), the coefficients $H_{jk}^i$ obey the
rule of transformation
\begin{equation}
\widetilde{H}_{rs}^{i\,}\displaystyle\frac{\partial \widetilde{x}^r}{%
\partial x^j}\displaystyle\frac{\partial \widetilde{x}^s}{\partial x^h}=%
\displaystyle\frac{\partial \widetilde{x}^i}{\partial x^r}H_{jh}^r-%
\displaystyle\frac{\partial ^2\widetilde{x}^i}{\partial x^j\partial x^h}
\tag{7.4.2}
\end{equation}

3$^{\circ }$ The system of functions $\underset{(\alpha
)}{C_{jh}^i,}C_i^{jh},(\alpha =1,...,k-1)$ are $d$-tensor fields
of type (1,2) and (2,1), respectively.
\end{teo}

\textbf{Proof}: According to the definition 2.1, we can write

$D_{\displaystyle\frac \delta {\delta x^j}}\displaystyle\frac
\delta {\delta x^i}=\underset{\left( 0\right)
}{H_{ij}^s}\displaystyle\frac \delta {\delta
x^s};D_{\displaystyle\frac \delta {\delta x^j}}\displaystyle\frac
\delta {\delta y^{\left( 1\right) i}}=\underset{\left( 1\right)
}{H_{ij}^s} \displaystyle\frac \delta {\delta y^{\left( 1\right)
s}},...,$

\vspace{3mm}

$D_{\displaystyle\frac \delta {\delta x^j}}\displaystyle\frac \delta {\delta
y^{\left( k-1\right) i}}=\underset{\left( k-1\right) }{H_{ij}^s} %
\displaystyle\frac \delta {\delta y^{\left( k-1\right) s}};D_{\displaystyle
\frac \delta {\delta x^j}}\displaystyle\frac \delta {\delta p_i}=\underset{%
\left( k\right) }{H_{js}^i}\displaystyle\frac \delta {\delta p_s};$

\vspace{3mm}

Applying the mapping $J$ and looking to (4.3.1) and to Theorem 7.2.3, we obtain

$D_{\displaystyle\frac \delta {\delta x^j}}(J\displaystyle\frac
\delta {\delta x^i})=J(\underset{\left( 0\right)
}{H_{ij}^s}\displaystyle\frac
\delta {\delta x^s}),_{}D_{\displaystyle\frac \delta {\delta x^j}}(J %
\displaystyle\frac \delta {\delta y^{\left( 1\right)
i}})=\underset{\left( 1\right) }{J(H_{ij}^s}\displaystyle\frac
\delta {\delta y^{\left( 1\right) s}}),...,$

$D_{\displaystyle\frac \delta {\delta x^j}}(J\displaystyle\frac
\delta {\delta y^{\left( k-2\right) i}})=J(\underset{\left(
k-1\right) }{H_{ij}^s} \displaystyle\frac \delta {\delta y^{\left(
k-1\right) s}}).$

Consequently,
\[
\underset{\left( 1\right) }{H_{ij}^s}=\underset{\left( 0\right) }{
H_{ij}^s},\underset{\left( 2\right) }{H_{ij}^s}=\underset{\left(
1\right) }{H_{ij}^s},...,\underset{\left( k-1\right) }{H_{ij}^s}=
\underset{\left( k-2\right) }{H_{ij}^s}
\]
and, from $D_{\displaystyle\frac \delta {\delta x^j}}\theta =0,$ we have $%
\underset{\left( 0\right) }{H_{ij}^s}=\underset{\left( k\right)
}{ H_{ij}^s}.$ It follows that the formulae of the first line of
(7.4.1) are valid. Similarly we prove the other equalities of (7.4.1).
By a direct calculus we prove that 2$^{\circ }$ and
3$^{\circ }$ from this theorem hold. q.e.d

The systems of functions
\begin{equation}
D\Gamma (N)=\{H_{jh}^i,\underset{\left( \alpha \right) }{C_{jh}^i}
,C_i^{jh}\},(\alpha =1,...,k-1)  \tag{7.4.3}
\end{equation}
is the system of coefficients of the $N$-linear connection $D$ in the
adapted basis.

The converse of the previous theorem holds,too.
\begin{teo}
If the system of functions (7.4.3) are apriori given over every domain of a local
chart on the manifold $T^{*k}M,$ having the rule of transformation mentioned
in the previous theorem, then there exists an unique $N$-linear connection $D$
whose local coefficients are just the system of given functions.
\end{teo}

\begin{cor}
The following formulae hold:
\begin{equation}
\begin{array}{l}
D_{\displaystyle\frac \delta {\delta x^j}}dx^i=-H_{js}^idx^s;D_{%
\displaystyle \frac \delta {\delta x^j}}\delta y^{\left( \alpha \right)
i}=-H_{js}^i\delta y^{\left( \alpha \right) s},(\alpha =1,...k-1); \\
D_{\displaystyle\frac \delta {\delta x^j}}\delta p_i=H_{ij}^s\delta p_s, \\
\\
D_{\displaystyle\frac \delta {\delta y^{\left( \alpha \right) j}}}dx^i=-%
\underset{(\alpha )}{C_{js}^i}dx^s;D_{\displaystyle\frac \delta
{\delta y^{\left( \alpha \right) j}}}\delta y^{(\beta
)i}=-H_{js}^i\delta y^{(\beta
)s},(\alpha ,\beta =1,...k-1); \\
D_{\displaystyle\frac \delta {\delta y^{\left( \alpha \right) j}}}\delta p_i=%
\underset{(\alpha )}{C_{ij}^s}\delta p_s; \\
\\
D_{\displaystyle\frac \delta {\delta p_j}}dx^i=-C_s^{ij}dx^s;D_{%
\displaystyle \frac \delta {\delta p_j}}\delta y^{\left( \alpha \right)
i}=C_s^{ij}\delta y^{\left( \alpha \right) s},(\alpha =1,...k-1); \\
D_{\displaystyle\frac \delta {\delta p_j}}\delta p_i=C_i^{js}\delta p_s.
\end{array}
\tag{7.4.4}
\end{equation}
\end{cor}

Indeed, the formula (7.4.1) and the conditions of duality between the vector
fields from the adapted basis and its dual 1-form basis lead to the formula
(7.4.4).

\section{The $h$-, $v_\alpha$- and $w_k$-Covariant Derivatives in Local
Adapted Basis}

In the adapted basis a tensor field $T$ can be written in the form (7.1.5):

\begin{equation}
T=T_{j_1...j_s}^{i_1...i_r}\displaystyle\frac \delta {\delta x^{i_1}}\otimes
...\otimes \displaystyle\frac \delta {\delta p_{j_s}}\otimes dx^{j_1}\otimes
...\otimes \delta p_{i_r}.  \tag{7.5.1}
\end{equation}

Applying the operator of covariant derivation $D_X$ for $X=X^H=X^i %
\displaystyle\frac \delta {\delta x^i}$ and taking into account the formulae
(7.4.1), (7.4.4) and the properties of the operator $D_X^HT=X^iD_{\displaystyle\frac
\delta {\delta x^j}}^HT,$ expressed in Theorem 7.2.3, we deduce
\begin{equation}
D_X^HT=X^mT_{j_1...j_s|m}^{i_1...i_r}\displaystyle\frac \delta {\delta
x^{i_1}}\otimes ...\otimes \displaystyle\frac \delta {\delta p_{j_s}}\otimes
dx^{j_1}\otimes ...\otimes \delta p_{i_r},  \tag{7.5.2}
\end{equation}
where
\begin{equation}
\begin{array}{l}
T_{j_1...j_s|m}^{i_1...i_r}=\displaystyle\frac{\delta
T_{j_1...j_s}^{i_1...i_r}}{\delta x^m}
+T_{j_1...j_s}^{hi_2...i_r}H_{hm}^{i_1}+...+T_{j_1...j_s}^{i_1...h}H_{hm}^{i_r}-
\\
\\
-T_{h...j_s}^{i_1...i_r}H_{j_1m}^h-...-T_{j_1...h}^{i_1...i_r}H_{j_sm}^h.
\end{array}
\tag{7.5.2a}
\end{equation}

The operator $"_{|}"$ is called the $h$-covariant derivative with respect to $%
D\Gamma (N).$

Now, we put $X=X^{V_\alpha }=X^i\displaystyle\frac \delta {\delta
y^{\left( \alpha \right) i}},(\alpha =1,...,k).$ From (7.5.1) we deduce
\begin{equation}
\begin{array}{l}
D_X^{V_\alpha }T=X^mT_{j_1...j_s}^{i_1...i_r}\stackrel{(\alpha )}{|_m} %
\displaystyle\frac \delta {\delta x^{i_1}}\otimes ...\otimes \displaystyle %
\frac \delta {\delta p_{j_s}}\otimes dx^{j_1}\otimes ...\otimes \delta
p_{i_r} \\
\alpha =1,...,k-1.
\end{array}
\tag{7.5.3}
\end{equation}
where
\begin{equation}
\begin{array}{l}
T_{j_1...j_s}^{i_1...i_r}\stackrel{(\alpha
)}{|_m}=\displaystyle\frac{\delta
T_{j_1...j_s}^{i_1...i_r}}{\delta y^{\left( \alpha \right) m}}
+T_{j_1...j_s}^{hi_2...i_r}\underset{\left( \alpha \right)
}{C_{hm}^{i_1}} +...+T_{j_1...j_s}^{i_1...h}\underset{\left(
\alpha \right) }{C_{hm}^{i_r}}
- \\
\\
-T_{h...j_s}^{i_1...i_r}\underset{\left( \alpha \right)
}{C_{j_1m}^h} -...-T_{j_1...h}^{i_1...i_r}\underset{\left( \alpha
\right) }{C_{j_sm}^h} ,\left( \alpha =1,...,k-1\right).
\end{array}
\tag{7.5.3a}
\end{equation}

The operator $\stackrel{\left( \alpha \right) }{|}$ is called the $v_\alpha $
-covariant derivative with respect to $D\Gamma (N).$

Finally, taking $X=X^{W_k}=X_i\displaystyle\frac \delta {\delta p_i},$ then $%
D_X^{W_k}T$ has the form:
\begin{equation}
D_X^{W_k}T=X_m T_{j_1...j_s}^{i_1...i_r}|^m\displaystyle\frac \delta {\delta
x^{i_1}}\otimes ...\otimes \displaystyle\frac \delta {\delta p_{j_s}}\otimes
dx^{j_1}\otimes ...\otimes \delta p_{i_r},  \tag{7.5.4}
\end{equation}
where
\begin{equation}
\begin{array}{l}
T_{j_1...j_s}^{i_1...i_r}|^m=\displaystyle\frac{\delta
T_{j_1...j_s}^{i_1...i_r}}{\delta p_m}
+T_{j_1...j_s}^{hi_2...i_r}C_h^{i_1m}+...+T_{j_1...j_s}^{i_1...h}C_h^{i_rm}-
\\
\\
-T_{h...j_s}^{i_1...i_r}C_{j_1}^{hm}-...-T_{j_1...h}^{i_1...i_r}C_{j_s}^{hm}.
\end{array}
\tag{7.5.4a}
\end{equation}

It is not hard to prove
\begin{prop}
The following properties hold:
\[
T_{j_1...j_s|m}^{i_1...i_r},\ T_{j_1...j_s}^{i_1...i_r}\stackrel{(\alpha )}{%
|_m},\ T_{j_1...j_s}^{i_1...i_r}|^m,\left( \alpha =1,...,k-1\right)
\]
are $d$-tensor fields. The first two are of type $(r,s+1)$ and the last one
is of type $(r+1,s).$
\end{prop}

As an application, the $d$-tensor field $g_{ij}$ has the $h$-, $v_\alpha $-
and $w_k$- covariant derivatives with respect to the $N$-linear connection
with the coefficients $D\Gamma (N),$ given by:
\begin{equation}
\left\{
\begin{array}{c}
g_{ij|k}=\displaystyle\frac{\delta g_{ij}}{\delta x^m}
-H_{im}^hg_{hj}-H_{jm}^hg_{ih}, \\
\\
g_{ij}\stackrel{\left( \alpha \right)
}{|}_m=\displaystyle\frac{\delta g_{ij} }{\delta y^{\left( \alpha
\right) m}}-\underset{\left( \alpha
\right) }{ C_{im}^h}g_{hj}-\underset{\left( \alpha \right) }{C_{jm}^h}%
g_{ih}, \ \ (\alpha =1,..,k-1),\\
\\
g_{ij}|^m=\displaystyle\frac{\delta g_{ij}}{\delta p_m}
-C_i^{hm}g_{hj}-C_j^{hm}g_{ih}.
\end{array}
\right.  \tag{7.5.5}
\end{equation}

In a next chapter we shall determine the coefficients $D\Gamma \left(
N\right) $ from the conditions that $g_{ij}$ is covariant constant with
respect $D$.

\begin{prop}
The operators $_{|}$ , $\stackrel{\left( \alpha \right) }{|}$ and $|$ have
the properties:

\vspace{3mm}

1$^{\circ }$ $f_{|m}=\displaystyle\frac{\delta f}{\delta x^m},f\stackrel{%
\left( \alpha \right) }{|}_m=\displaystyle\frac{\delta f}{\delta y^{\left(
\alpha \right) m}},\left( \alpha =1,...,k-1\right) ,f|^m=\displaystyle\frac{%
\delta f}{\delta p_m}$

\quad

2$^{\circ }$ These operators are distributive with respect to the addition
of the $d$-tensor of the same type.

3$^{\circ }$ They commute with the operation of contraction

4$^{\circ }$ They verify the Leibniz rule to the tensor product.
\end{prop}

As an application we study the $(\stackrel{(\alpha )}{z})$ -deflection
tensor of $D\Gamma \left( N\right) .$ They are defined by:
\begin{equation}
\stackrel{\left( \alpha \right) }{D_j^i}=\stackrel{\left( \alpha \right) }{
z_{|j}^i},\ \stackrel{\left( \alpha \beta \right) }{D_j^i}=\stackrel{\left(
\alpha \right) }{z_{}^i}\stackrel{\left( \beta \right) }{|}_j,(\alpha ,\beta
=1,...,k-1),\ \stackrel{\left( \alpha \right) }{D^{ij}}=\stackrel{\left(
\alpha \right) }{z^i}|^j,  \tag{7.5.6}
\end{equation}
where $\stackrel{\left( \alpha \right) }{z^i},(\alpha =1,...,k-1)\,$are the
Liouville $d$-vector fields.

Evidently, they have the following expressions
\begin{equation}
\begin{array}{l}
\stackrel{\left( \alpha \right) }{D_j^i}=\displaystyle\frac{\delta \stackrel{%
\left( \alpha \right) }{z^i}}{\delta x^j}+\stackrel{\left( \alpha
\right) }{ z^m}H_{mj}^i,\ \stackrel{\left( \alpha \beta \right)
}{D_j^i}=\displaystyle \frac{\delta \stackrel{\left( \alpha
\right) }{z^i}}{\delta y^{(\beta )j}}+ \stackrel{\left( \alpha
\right) }{z^m}\underset{\left( \beta \right) }{
C_{mj}^i}, \\
\\
\stackrel{\left( \alpha \right) }{D^{ij}}=\displaystyle\frac{\delta
\stackrel{\left( \alpha \right) }{z^i}}{\delta p_j}+\stackrel{\left( \alpha
\right) }{z^m}C_m^{ij}.
\end{array}
\tag{7.5.6a}
\end{equation}

Similarly, we introduce the $(p)$-deflection tensors by

\begin{equation}
\Delta _{ij}=p_{i|j},\ \stackrel{(\alpha )}{\not \delta _{ij}}=p_i\stackrel{
(\alpha )}{|_j},\ \not \delta _i^j=p_i|^j.  \tag{7.5.7}
\end{equation}

We deduce:
\[
\Delta _{ij}=-p_hH_{ij}^h,\ \stackrel{(\alpha )}{\not \delta
_{ij}}=-p_h \underset{\left( \alpha \right) }{C_{ij}^h},\ \not
\delta _i^j=\delta _i^j-p_hC_i^{hj}.
\]

The deflection tensors will be used in some important identities determined
by the Ricci identities, applied to the Liouville $d$-tensor fields $%
\stackrel{\left( \alpha \right) }{z^i},$ and to the $d$-covector $p_i.$

\textit{Some remarks}:

1$^{\circ }$ In the adapted basis we can prove the equation

\[
D_X\Bbb{F}=0,\text{ }\forall X\in X(T^{*k}M)
\]

2$^{\circ }$ A Berwald connection is an $N$-linear connection $D$ with the

coefficients
\begin{equation}
B\Gamma \left( N\right) =(B_{jk}^i,0,...,0,0)  \tag{7.5.8}
\end{equation}
where $B_{jh}^i$ has the same rule of transformation as $H_{jh}^i$ from
(7.4.2) and is determined by the nonlinear connection $N$.

We have:
\begin{teo}
Any nonlinear connection $N,$ with the coefficients \newline
$(\underset{\left( 1\right) }{N_j^i},...,\underset{\left(
k-1\right) }{N_j^i},N_{ij})$ determines the following Berwald
connections:
\begin{equation}
B\Gamma \left( N\right) =(\stackrel{\cdot }{\partial^i} N_{ij},0,...,0,0)
\tag{7.5.9}
\end{equation}
\begin{equation}
\underset{\left( 1\right) }{B}\Gamma \left( N\right) =(\displaystyle\frac{%
\delta \underset{\left( 1\right) }{N_j^i}}{\delta y^{\left( 1\right) h}}%
,0,...,0,0)  \tag{7.5.9a}
\end{equation}
\end{teo}

\textbf{Proof}: The first one is given by the last formula (7.1.6) applying
the derivation $\stackrel{\cdot }{\widetilde{\partial }^i}=\displaystyle
\frac{\partial \widetilde{x}^i}{\partial x^s}\stackrel{\cdot }{\partial ^s}$
and (7.5.9') is obtained from the first formula (7.1.6), applying the derivation

$\displaystyle\frac \delta {\delta \widetilde{y}^{(1)h}}=\displaystyle\frac{
\partial x^m}{\partial \widetilde{x}^h}\displaystyle\frac \delta {\delta
y^{(1)m}}.$ q.e.d.

If the base manifold $M$ is paracompact, then the manifold $T^{*k}M$ is
paracompact, too. Consequently, on $T^{*k}M$ there exist nonlinear
connections. Therefore, we have:

\begin{teo}
If the base manifold $M$ is paracompact, then on the manifold $T^{*k}M$
there exist $N$-linear connections $D$.
\end{teo}

Indeed, on $T^{*k}M$ there exist nonlinear connections $N.$ Applying the
Theorem 7.5.1 the conclusion follows. \hfill  Q.E.D.

\section{Ricci Identities. Local Expressions of $d$-Tensor of Curvature and
Torsion. Bianchi Identities.}

Let $D$ be an $N$-linear connection with the local coefficients
\begin{equation}
D\Gamma \left( N\right) =(H_{jh}^i,\underset{\left( 1\right)
}{C_{jh}^i} ,...,\underset{\left( k-1\right) }{C_{jh}^i},C_i^{jh})
\tag{7.6.1}
\end{equation}

The Ricci identities for a $d$-vector field $X^i$ can be deduced from the
formulae (7.3.8) written in the adapted basis. But we can obtain them by a
straighforward calculus.

\begin{teo}
For any $N$-linear connection $D$ and any $d$-vector field $X^i$ the
following Ricci formulae hold:
\begin{equation}
\begin{array}{l}
X_{|j|h}^i-X_{|h|j}^i=X^mR_{m\ jh}^i-X_{|m}^iT_{\ jh}^m-\{X^i\stackrel{\left(
1\right) }{|_m}\underset{\left( 01\right) }{R_{\ jh}^m}+...+ \\
\\
+X^i\stackrel{\left( k-1\right) }{|_m}\underset{\left( 0,k-1\right) }{%
R_{\ jh}^m}+X^i|^m\underset{\left( 0\right) }{R_{mjh}}\},
\end{array}
\tag{7.6.2$_1$}
\end{equation}
\begin{equation}
\begin{array}{l}
\stackrel{\left( \alpha \right) }{X_{|j}^i|_h}-\stackrel{\left( \alpha
\right) }{X^i|_{h|j}}=X^m\underset{\left( \alpha \right) }{\ P_{m\ jh}^i}%
-X_{|m}^i\underset{\left( \alpha \right)
}{C_{jh}^m}-\{X^i\stackrel{\left( 1\right) }{|_m}\underset{\left(
01\right) }{R_{\ jh}^m}+X^i\stackrel{\left(
\alpha \right) }{|_m}H_{jh}^m\}- \\
\\
-\{X^i\stackrel{\left( 1\right) }{|_m}\underset{\left( \alpha ,1\right) }{%
B_{jh}^m}+...+X^i\stackrel{\left( k-1\right)
}{|_m}\underset{\left( \alpha ,k-1\right)
}{B_{jh}^m}+X^i|^m\underset{\left( \alpha \right) }{B_{mjh}}\},
\end{array}
\tag{7.6.2$_2$}
\end{equation}
\begin{equation}
\begin{array}{l}
X_{|j}^i|_h-X^i|_{\ |j}^h=X^mP_{m\ \ j}^{\ i\ h}-X_{|m}^iC_j^{mh}-X^i|^mH_{mj}^{\ \ h}- \\
\\
-\left\{ X^i\stackrel{(1)}{|}_m\underset{(1)}{B\text{ }}\text{ }%
_j^{mh}+\cdots +X^i\stackrel{(k-1)}{|}_m\underset{(k-1)}{B\text{}}\text{}%
_j^{mh}+X^i|_m\underset{(0)}{B\text{}}\text{}_{mj}^{\ \ h}\right\},
\end{array}
\tag{7.6.2$_3$}
\end{equation}
\begin{equation}
\begin{array}{l}
X^i\stackrel{(\alpha )}{|}_j\stackrel{(\beta )}{|}_h-X^i\stackrel{(\beta )}{|%
}_h\stackrel{(\alpha )}{|}_j=X^m\underset{(\alpha \beta
)}{S}\text{}_{m\ jh}^{\ i}-X^i\stackrel{(\alpha
)}{|}_m\underset{(\beta )}{C}\text{}_{\ jh}^m+X^i\stackrel{(\beta )}{|}_m\underset{(\alpha )}{C}\text{}_{\ hj}^m- \\
\\
-\left\{ X^i\stackrel{(1)}{|}_m\underset{(\alpha \beta )}
{\stackrel{(1)}{C}\text{}}\text{}_{\ jh}^m+\cdots +X^i\stackrel{(k-1)}{|}_m\underset{(\alpha \beta )}
{\stackrel{(k-1)}{C}\text{}}\text{}_{\ jh}^m+X^i|_m
\underset{(\alpha \beta )}{B\text{}}\text{}_{mjh}\right\},
\end{array}
\tag{7.6.2$_4$}
\end{equation}

\begin{equation}
\begin{array}{l}
X^i\stackrel{(\alpha )}{|}_j|^h-X^i|^h\stackrel{(\alpha
)}{|}_j=X^m\underset{(\alpha)}{S}\textrm{ }_{m\ \ j}^{\ i\
h}-X^i|^m\underset{(\alpha)}{C}\textrm{
}_{mj}^h+X^i\stackrel{(\alpha
)}{|}_mC\textrm{ }_j^{mh}-\\
-\left\{X^i\stackrel{(1)}{|}_m\underset{(\alpha 1)}{C\textrm{
}}\textrm{ }_{\
j}^{mh}+\cdots+X^i\stackrel{(k-1)}{|}_m\underset{(\alpha,k-1)}{C\textrm{
}}\textrm{ }_{j}^{mh}+X^i|^m\underset{(\alpha )}{C\text{}}\textrm{
}_{mj}^{h}\right\},
\end{array}\tag{7.6.$2_5$}
\end{equation}

\begin{equation}
X^i|^j|^h-X^i|^h|^j=X^mS_m^{\ \ \ ijh}+X^i|^mS_m^{jh},
\tag{7.6.$2_6$}
\end{equation}
where all terms in $\underset{(0\alpha )}{R}\text{}_{\ jh}^i$ , $\underset{(0)}{R}\text{}_{mjh}$,
$\underset{(\alpha \beta )}{B}\text{}_{\ jh}^i$, $\underset{(\alpha )}{B}\text{}_{ijh}$,
$\underset{(0)}{B}\text{}_{mj}^{\ \ h}$ are known from the Lie brackets (6.5.1).
\end{teo}
The coefficients $D\Gamma (N)$ are given in (7.6.1). Supplementary we put:
\begin{equation}
T_{\ jh}^i=H_{jh}^i-H_{hj}^i,\quad S_i^{\ jh}=C_i^{jh}-C_i^{hj}. \tag{7.6.3}
\end{equation}
Here $R_{\ mjh}^i$, ..., are called\textit{\ }$d$\textit{-tensor of
curvature of }$D$ and they have the expressions:
\begin{equation}
\left\{
\begin{array}{lll}
R_{m\ jh}^{\ \ i} & = & \displaystyle\frac{\delta H_{mj}^i}{\delta x^h}-%
\displaystyle \frac{\delta H_{mh}^i}{\delta x^j}%
+H_{mj}^sH_{sh}^i-H_{mh}^sH_{sj}^i+ \underset{(1)}{C}\text{ }_{ms}^i%
\underset{(01)}{R}\text{}_{\ jh}^s+\text{ } \\
& & \\
& + & \cdots +\underset{(k-1)}{C}\text{}_{ms}^i\underset{(0,k-1)}{R}\text{}_{\ jh}^s+C_m^{is}\underset{(0)}{R}\text{}_{sjh}, \\
& & \\
\underset{(\alpha )}{P}\text{}_{m \ jh}^{\ \ i} & = &
\displaystyle \frac{\delta H_{mj}^i}{\delta y^{(\alpha
)h}}-\displaystyle\frac{\delta \underset{(\alpha )}{C}\text{}_{mh}^i}{\delta x^j}+H_{mj}^s\underset{ (\alpha )}{C}\text{}_{sh}^i-\underset{(\alpha )}{C}\text{}_{mh}^sH_{sj}^i+\underset{(1)}{C}\text{}_{ms}^i\underset{(\alpha1)}{B}\text{}_{jh}^s+ \\
& & \\
& + & \cdots +\underset{(k-1)}{C}\text{ }_{ms}^i\underset{(\alpha,k-1)}{B}\text{}_{jh}^s+C_m^{is}\underset{(\alpha )}{B}\text{}_{sjh}, \\
& & \\
P_{m\ \ j}^{\ \ i\ \ \ h} & = & \displaystyle\frac{\delta H_{mj}^i}{\delta p_h}-\displaystyle\frac{\delta C_m^{ih}}{\delta
x^j} +H_{mj}^sC_s^{ih}-C_m^{sh}H_{sj}^i+\underset{(1)}{C}\text{}_{ms}^i
\underset{(1)}{B}\text{}_j^{sh}+ \\
& & \\
& + & \cdots +\underset{(k-1)}{C}\text{}_{ms}^i\underset{(k-1)}{B}\text{}_j^{sh}+C_m^{is}\underset{(0)}{B}\text{}_{s\ j}^{\ h}.
\end{array}
\right.  \tag{7.6.4}
\end{equation}
and
\begin{equation}
\left\{
\begin{array}{lll}
\underset{(\alpha \beta )}{S}\text{}_{m\ jh}^{\ \ i} & = & %
\displaystyle\frac{\delta \underset{(\alpha )}{C}\text{}_{mj}^i}{\delta y^{(\beta )h}}-\displaystyle\frac{\delta
\underset{(\beta )}{C}\text{}_{mh}^i}{\delta y^{(\alpha)j}}+\underset{(\alpha )}{C}\text{}_{mj}^s \underset{(\beta
)}{C}\text{}_{sh}^i-\underset{(\beta )}{C}\text{}_{mh}^s\underset{(\alpha )}{C}\text{}_{sj}^i+\underset{(1)}{C}\text{}_{ms}^i\underset{(\alpha \beta )}{\stackrel{(1)}{C}}\text{}_{jh}^s+ \\
& & \\
& + & \cdots +\underset{(k-1)}{C}\text{}_{ms}^i\underset{(\alpha \beta) }{\stackrel{(k-1)}{C}}\text{}_{jh}^s+
C_m^{is}\underset{(\alpha \beta )}{B} \text{}_{sjh},\ (\alpha \leq \beta ,\alpha ,\beta =1,...,k-1) \\
& & \\
\underset{(\alpha )}{S}\text{}_{m\ j}^{\ \ i\ \ h} & = & %
\displaystyle\frac{\delta \underset{(\alpha)}{C}\text{}_{mj}^i}{\delta p_h}-\displaystyle\frac{\delta C_m^{ih}}{\delta
y^{(\alpha )j}}+\underset{ (\alpha )}{C}\text{}_{mj}^sC_s^{ih}-C_{mh}^s\underset{(\alpha )}{C}\text{}_{sj}^i+\underset{(1)}{C}\text{}_{ms}^i\underset{(\alpha
1)}{C}\text{}_j^{sh}+ \\
& & \\
& + & \cdots +\underset{(k-1)}{C}\text{}_{ms}^i\underset{(\alpha,k-1)}{C}\text{}_{sj}^i+C_m^{is}\underset{(\alpha )}{C}\text{}_{sj}^h, \\
& & \\
S_m^{\ \ \ ijh} & = & \displaystyle\frac{\delta C_m^{ij}}{\delta p_h}- %
\displaystyle\frac{\delta C_m^{ih}}{\delta p_j}
+C_m^{sj}C_s^{ih}-C_m^{sh}C_s^{ij}.
\end{array}\right.
\tag{7.6.5}
\end{equation}

As usually, we extend the Ricci identities (7.6.2) for any $d$-tensor field $%
T_{j_1...j_s}^{i_1...i_r}$.

For instance, if $g^{ij}(x,y^{(1)},...,y^{(k-1)},p)$ is a $d$-tensor field,
the Ricci identities of $g^{ij}$, with respect to the $N$-linear connection $D\Gamma (N)$, are:
\begin{equation}
\begin{array}{lll}
g_{\ \ |h|m}^{ij}-g_{\ |m|h}^{ij} &= & g^{sj}R_{s\ hm}^{\ i}-g^{is}R_{s\
hm}^{\ j}-g_{\ |s}^{ij}T_{\ hm}^s- \{g^{ij}\stackrel{(1)}{|}_s\underset{(01)}{R}\text{}_{\ hm}^s+\\
& & \\
& & +\cdots
+g^{ij}\stackrel{(k-1)}{|}_s\underset{(0,k-1)}{R}\text{}_{\ hm}^s+g^{ij}|^s\underset{(0)}{R}\text{}_{shm}\} ,  \\
.................... &  & \\
& & \\
g^{ij}|^h|^m-g^{ij}|^m|^h & = &
g^{sj}S_s^{\ ihm}+g^{is}S_s^{\ jhm}+g^{ij}|^sS_s^{\ hm}.
\end{array}
\tag{7.6.6}
\end{equation}

In particular if $D\Gamma (N)$ satisfies the supplementary conditions:
\begin{equation}
g_{\ \ |h}^{ij}=0,\quad g^{ij}\stackrel{(\alpha )}{|}_h=0,\quad
g^{ij}|^h=0,\quad (\alpha =1,...,k-1),  \tag{7.6.7}
\end{equation}
then the Ricci identities (7.6.6) give us:
\begin{equation}
\begin{array}{l}
g^{sj}R_{s\ \ hm}^{\ \ i}+g^{is}R_{s\ \ hm}^{\ \ j}=0, \\
\\
g^{sj}\underset{(\alpha )}{P}\text{}_{s\ \ hm}^{\ \ i}+g^{is}\underset{(\alpha )}{P}\text{}_{s\ \ hm}^{\ \ j}=0 ,\\
........................................, \\
g^{sj}S_s^{\ \ ihm}+g^{is}S_s^{\ \ jhm}=0.
\end{array}
\tag{7.6.8}
\end{equation}

Some important identities are obtained applying the Ricci identities to the $d$
-covector $p_i$ and to the Liouville vector fields $z^{(\alpha )i}$, ($\alpha =1$, ..., $k-1$).
\begin{teo}
Any $N$-linear connection $D\Gamma (N)$ satisfies the following identities:
\begin{equation}
\begin{array}{lll}
\Delta _{ij|h}-\Delta _{ih|j} & = & -p_sR_{i\ \ jh}^{\ \ s}-\Delta
_{is}T_{\ jh}^s- \\
& - & \left\{ \stackrel{(1)}{\not \delta }_{is}\underset{(01)}{R}\text{}_{jh}^s+\cdots +\stackrel{(k-1)}{\not \delta }_{is} \underset{(0,k-1)}{R}\text{}_{jh}^s+\not \delta _i^s\underset{(0)}{R}\text{}_{sjh}\right\} ,
\\
\Delta _{ij}\stackrel{(\alpha )}{|}_h-\stackrel{(\alpha )}{\not \delta }%
\text{}_{ih|j}\text{} & = & -p_s\underset{(\alpha )}{P}\text{}_{i\ \
jh}^{\ s}-\Delta _{is}\underset{(\alpha )}{C}\text{}_{jh}^s-\stackrel{(\alpha )}{\not \delta }\text{}_{is}H_{jh}^s- \\
& - & \left\{ \stackrel{(1)}{\not \delta }_{is}\underset{(\alpha 1)}{B}%
\text{}_{\ jh}^s+\cdots +\stackrel{(k-1)}{\not \delta }_{is}\underset{%
(\alpha ,k-1)}{B}\text{}_{\ jh}^s+\not \delta _i^s\underset{(\alpha )}{B}%
\text{}_{sjh}\right\} , \\
\Delta _{ij}|^h-\not \delta _{i\ |j}^{\ h} & = & -p_sP_{i\ \ j}^{\ s\ \
h}-\Delta _{is}C_j^{sh}-\not \delta _i^sH_{sj}^h- \\
& - & \left\{ \stackrel{(1)}{\not \delta }_{is}\underset{(1)}{B}\text{}%
_j^{sh}+\cdots +\stackrel{(k-1)}{\not \delta }_{is}\underset{(k-1)}{B}%
\text{}_j^{sh}+\not \delta _i^s\underset{(0)}{B}\text{}_{sj}^{\
\
h}\right\} , \\
\stackrel{(\alpha )}{\not \delta }_{ij}\stackrel{(\beta )}{|}_h-\stackrel{%
(\beta )}{\not \delta }_{ih}\stackrel{(\alpha )}{|}_j & = & -p_s\underset{%
(\alpha \beta )}{S}\text{}_{i\ \ jh}^{\ s}-\stackrel{(\alpha
)}{\not \delta }_{is}\underset{(\beta )}{C}\text{}_{\
jh}^s-\stackrel{(\beta )}{\not
\delta }\text{}_{is}\underset{(\alpha )}{C}\text{}_{\ hj}^s- \\
& - & \left\{ \stackrel{(1)}{\not \delta }_{is}\underset{(\alpha \beta )}{%
\stackrel{(1)}{C}}\text{}_{\ jh}^s+\cdots +\stackrel{(k-1)}{\not \delta }%
_{is}\underset{(\alpha \beta )}{\stackrel{(k-1)}{C}}\text{}_{\
jh}^s+\not
\delta _i^s\underset{(\alpha \beta )}{B}\text{}_{sjh}\right\} , \\
\stackrel{(\alpha )}{\not \delta }_{ij}|^h-\not \delta _i^h\stackrel{(\beta )%
}{|}_j & = & -p_s\underset{(\alpha )}{S}\text{ }_{i\ \ j}^{\ s\ \
h}-\not
\delta _i^s\underset{(\alpha )}{C}\text{}_{s\ j}^{\ h}-\stackrel{(\alpha )%
}{\not \delta }\text{}_{is}C_j^{sh}- \\
& - & \left\{ \stackrel{(1)}{\not \delta }_{is}\underset{(\alpha 1)}{C}%
\text{}_{\ j}^{sh}+\cdots +\stackrel{(k-1)}{\not \delta }_{is}\underset{%
(\alpha ,k-1)}{C}\text{}_{\ j}^{sh}+\not \delta _i^s\underset{(\alpha )}{C%
}\text{}_{sj}^{\ h}\right\} , \\
\not \delta _i^{\ j}|^h-\not \delta _i^{\ h}|^j & = & -p_sS_i^{\ sjh}+\not
\delta _i^sS_s^{jh}.
\end{array}
\tag{7.6.9}
\end{equation}

\end{teo}

The similar identities are obtained applying the Ricci identities to the
Liouville vector fields $z^{(\alpha )i}$.
\begin{teo}
Any $N$-linear connection $D\Gamma (N)$ satisfies the following identities,
obtained from (7.6.2) for $X^i=z^{(\alpha )i}$:
\begin{equation}
\begin{array}{lll}
\stackrel{(\alpha )}{D}\text{}_{\ j|h}^i-\stackrel{(\alpha )}{D}\text{}_{\
h|j}^i & = & z^{(\alpha )s}R_{s\ jh}^{\ i}-\stackrel{(\alpha )}{D}\text{}%
_s^iT_{\ jh}^s- \\
& - & \{ \stackrel{(\alpha 1)}{D}\text{}_{\ s}^i\underset{(0)}{R}%
\text{}_{\ jh}^s+\cdots +\stackrel{(\alpha ,k-1)}{D}\text{}_{\ s}^i%
\underset{(0,k-1)}{R}\text{}_{\ jh}^s+\stackrel{(\alpha )}{D}\text{}^{is}%
\underset{(0)}{R}\text{}_{sjh}\} , \\
..................... & ... &
.......................................................................................
\\
\stackrel{(\alpha )}{D}\text{}^{ij}|^h-\stackrel{(\alpha )}{D}\text{}%
^{ih}|^j & = & z^{(\alpha )s}S_s^{\ ijh}-\stackrel{(\alpha )}{\not \delta }%
\text{}^{is}S_s^{\ \ jh},\ (\alpha =1,...,k-1).
\end{array}
\tag{7.6.10}
\end{equation}
\end{teo}

Evidently, the whole theory from the present section is simplified if $D\Gamma
(N)=B\Gamma (N)$=($B_{jh}^i$, $0$, ..., $0$, $0$) is a Berwald connection,
taking $X_{\ ||j}^i=\displaystyle\frac{\delta X^i}{\delta x^j}+X^hB_{hj}^i$,
$X^i\stackrel{(\alpha )}{||}_j=\displaystyle\frac{\delta X^i}{\delta
y^{(\alpha )j}}$, $X^i||^j=\displaystyle\frac{\delta X^i}{\delta p_j}$, '$%
_{||}$', '$\stackrel{(\alpha )}{||}$' and '$||$' being the $h-$, $%
v_\alpha -$ and $w_k$-covariant derivatives with respect to $B\Gamma (N)$.

The Bianchi identities for the $N$-linear connection $D\Gamma (N)$ can be
obtained using the methods described above. Namely, applying the Ricci
identities to $d$-tensor field $X_{\ |j}^i$ we obtain:

$\left( X_{\ |j}^i\right) {}_{|h|m}-\left( X_{\ |j}^i\right) {}_{|m|h}=X_{\
|j}^rR_{r\ hm}^i-X_{\ |r}^iR_{j\ hm}^i-X_{\ |r}^iT_{hm}^r-$

$\left( X_{\ |j}^i\stackrel{(1)}{|}_r\underset{(01)}{R}\text{ }_{\
hm}^r+\cdots +X_{\
|j}^i\stackrel{(k-1)}{|}_r\underset{(0,k-1)}{R}\text{ } _{\
hm}^r+X_{\ |j}^i|^r\underset{(01)}{R}\text{ }_{rhm}\right) .$

If we  cyclically permute the indices $j$, $h$, $m$ and add the identitiessuch obtained
we determine a first set of Bianchi identities:
\begin{equation}
\begin{array}{l}
\underset{(jhm)}{S}\left\{ R_{j\ hm}^{\
i}-T_{\ jh|m}^i-T_{\ jr}^iT_{\ hm}^r+\sum\limits_{\alpha
=1}^{k-1}\underset{ (\alpha )}{C}\text{}_{jr}^i\underset{(0\alpha
)}{R}\text{}_{\
hm}^r+C_j^{ir}\underset{(0)}{R}\text{}_{rhm}\right\} =0, \\
\\
\underset{(jhm)}{S}\left\{ R_{s\ jh|m}^i+R_{s\ jr}^{\
i}T_{\ hm}^r-\sum\limits_{\alpha =1}^{k-1}\underset{(\alpha
)}{P}\text{}_{s\ jr}^{\ i}\underset{(0\alpha )}{R}\text{}_{\
hm}^r-P_{s\ j}^{\ i\ \
r}\underset{(0)}{R}\text{}_{rhm}\right\} =0, \\
\\
\underset{(jhm)}{S}\left\{ \underset{(0\alpha )}{R}\text{}_{\
jh|m}^i- \underset{(0\alpha )}{R}\text{}_{\ jr}^iT_{\ hm}^r+\underset{(\alpha )}{P} \text{}_{jr}^{is}\underset{(0\alpha )}{R}\text{}_{\ hm}^r+\underset{
(\alpha )}{B}\text{}_j^{ir}\underset{(0)}{R}\text{}_{rhm}\right\} =0, \\
\\
\underset{(jhm)}{S}\left\{ \underset{(0)}{R}\text{}_{sjh|m}-\underset{ (0)}{R}\text{}_{sjr}T_{\ hm}^r+\sum\limits_{\alpha =1}^{k-1}\underset{ (\alpha
)}{B}\text{}_{sjr}\underset{(0\alpha )}{R}\text{}_{\ hm}^r+
\underset{(0)}{P}\text{}_{\ sj}^r\underset{(0)}{R}\text{}_{rhm}\right\} =0.
\end{array}
\tag{7.6.11}
\end{equation}
\begin{equation}
\underset{(\alpha )}{P}_{\ jr}^i=H_{jr}^i+\underset{(\alpha
)}{B}\text{ } _{\ jr}^i,\ (\alpha =1,...,k-1).  \tag{7.6.12}
\end{equation}

In a similar way we have the second set of Bianchi identities:

\begin{equation}
\begin{array}{l}
\underset{(jhm)}{S}\left\{ S_s^{\ ijh}|^m-S_s^{\ ijr}S_r^{\
hm}\right\} =0,
\\
\\
\underset{(ihm)}{S}\left\{ S_s^{\ ihm}-S_s^{\ ih}|^m-S_s^{\
ir}S_r^{\ hm}\right\} =0,
\end{array}
\tag{7.6.13}
\end{equation}
where $\underset{(jhm)}{S}$ is the cyclic sum.

Similarly,  we can get the other Bianchi identities.

\section{Parallelism of the Vector Fields on the Manifold $T^{*k}M$}

Consider a $N$-linear connection $D$ on the manifold $T^{*k}M$
with the coefficients $D\Gamma (N)$=($H_{\ jh}^i$,
$\underset{(\alpha )}{C}$ $_{\
jh}^i$, $C_i^{\ jh}$), ($\alpha =1,...,k-1$) in the adapted basis \newline
$\left(\displaystyle\frac \delta {\delta x^i},\displaystyle\frac \delta {\delta
y^{(\alpha )i}},\displaystyle\frac \delta {\delta p_i}\right) $.

A smooth curve $\gamma :I\rightarrow T^{*k}M$ having the image in a domain
of a local chart is given by
\begin{equation}
x^i=x^i(t),\ y^{(\alpha )i}=y^{(\alpha )i}(t),\ p_i=p_i(t),(\alpha
=1,...,k-1),\ t\in I.  \tag{7.7.1}
\end{equation}

The tangent vector field $\stackrel{\cdot }{\gamma }=\displaystyle\frac{
d\gamma }{dt}$ can be written by means of (7.2.10), in the form:
\begin{equation}
\stackrel{\cdot }{\gamma }=\displaystyle\frac{dx^i}{dt}\displaystyle\frac
\delta {\delta x^i}+\displaystyle\frac{\delta y^{(1)i}}{dt}\displaystyle %
\frac \delta {\delta y^{(1)i}}+\cdots +\displaystyle\frac{\delta y^{(k-1)i}}{
dt}\displaystyle\frac \delta {\delta y^{(k-1)i}}+\displaystyle\frac{\delta
p_i}{dt}\displaystyle\frac \partial {\partial p_i}  \tag{7.7.2}
\end{equation}
where, from (6.3.2),
\begin{equation}
\left\{
\begin{array}{l}
\displaystyle\frac{\delta y^{(1)i}}{dt}=\displaystyle\frac{dy^{(1)i}}{dt}+
\underset{(1)}{M}\text{}_j^i\displaystyle\frac{dx^j}{dt}, \\
......................................... \\
\displaystyle\frac{\delta
y^{(k-1)i}}{dt}=\displaystyle\frac{dy^{(k-1)i}}{dt}
+\underset{(1)}{M}\text{}_j^i\displaystyle\frac{dy^{(k-2)j}}{dt}+\cdots +
\underset{(k-2)}{M}\text{}_j^i\displaystyle\frac{dy^{(1)j}}{dt}+
\underset{(k-1)}{M}\text{}_j^i\displaystyle\frac{dx^j}{dt}, \\
\\
\displaystyle\frac{\delta p_i}{dt}=\displaystyle\frac{dp_i}{dt}-N_{ji} %
\displaystyle\frac{dx^j}{dt}.
\end{array}
\right.  \tag{7.7.3}
\end{equation}

Let us denote
\begin{equation}
D_{\stackrel{\cdot }{\gamma }}X=\displaystyle\frac{DX}{dt},\ DX=%
\displaystyle \frac{DX}{dt}dt,\ \forall X\in \mathcal{X}(T^{*k}M).  \tag{7.7.4}
\end{equation}

$DX$ is called \textit{the covariant differential of the vector field }$X$
and $\displaystyle\frac{DX}{dt}$ is \textit{the covariant differential of }$%
X $\textit{\ along curve }$\gamma $.

If $X$ is written in the form
\begin{equation}
\begin{array}{lll}
X & = & X^H+X^{V_1}+\cdots +X^{V_{k-1}}+X^{W_k}= \\
& & \\
& = & \stackrel{(0)i}{X}\displaystyle\frac \delta {\delta x^i}+\stackrel{%
(1)i }{X}\displaystyle\frac \delta {\delta y^{(1)i}}+\cdots +\stackrel{(k-1)i%
}{X} \displaystyle\frac \delta {\delta y^{(k-1)i}}+X_i\displaystyle\frac
\delta {\delta p_i}
\end{array}
\tag{7.7.4a}
\end{equation}
and $\stackrel{\cdot }{\gamma }$ from (7.7.2) is $\stackrel{\cdot }{\gamma }=
\stackrel{\cdot }{\gamma }^H+\stackrel{\cdot }{\gamma }^{V_1}+\cdots +
\stackrel{\cdot }{\gamma }^{V_{k-1}}+\stackrel{\cdot }{\gamma }^{W_k}$ we
have

\[
D_{\stackrel{\cdot }{\gamma }}=D_{\stackrel{\cdot }{\gamma }}^H+D_{\stackrel{
\cdot }{\gamma }}^{V_1}+\cdots +D_{\stackrel{\cdot }{\gamma }}^{V_{k-1}}+D_{%
\stackrel{\cdot }{\gamma }}^{W_k}.
\]
Then $DX$ has the final form:
\begin{equation}
\begin{array}{lll}
DX & = & \left( d\stackrel{(0)i}{X}+\stackrel{(0)s}{X}\omega _s^i\right) %
\displaystyle\frac \delta {\delta x^i}+\left( d\stackrel{(1)i}{X}+\stackrel{
(1)s}{X}\omega _s^i\right) \displaystyle\frac \delta {\delta
y^{(1)i}}+\cdots + \\
& + & \left( d\stackrel{(k-1)i}{X}+\stackrel{(k-1)s}{X}\omega _s^i\right) %
\displaystyle\frac \delta {\delta y^{(k-1)i}}+\left( dX_i-X_s\omega
_i^s\right) \displaystyle\frac \delta {\delta p_i},
\end{array}
\tag{7.7.5}
\end{equation}
where
\begin{equation}
\omega _j^i=H_{js}^idx^s+\underset{(1)}{C}\text{}_{js}^i\delta
y^{(1)s}+\cdots +\underset{(k-1)}{C}\text{}_{js}^i\delta
y^{(k-1)s}+C_j^{is}\delta p_s.  \tag{7.7.6}
\end{equation}

The differential forms $\omega _j^i$ are called \textit{the }$1$\textit{\
-forms connection of the connection }$D$.

Putting
\begin{equation}
\displaystyle\frac{\omega
_j^i}{dt}=H_{js}^i\displaystyle\frac{dx^s}{dt}+
\underset{(1)}{C}\text{}_{js}^i\displaystyle\frac{\delta
y^{(1)s}}{dt}+\cdots +
\underset{(k-1)}{C}\text{}_{js}^i\displaystyle\frac{\delta
y^{(k-1)s}}{dt} +C_j^{is}\displaystyle\frac{\delta p_s}{dt}.
\tag{7.7.7}
\end{equation}
the covariant differential $\displaystyle\frac{DX}{dt}$ along curve $%
\gamma $ is
\begin{equation}
\begin{array}{lll}
\displaystyle\frac{DX}{dt} & = & \left( \displaystyle\frac{d\stackrel{(0)i}{%
X }}{dt}+\stackrel{(0)s}{X}\displaystyle\frac{\omega _s^i}{dt}\right) %
\displaystyle\frac \delta {\delta x^i}+\left( \displaystyle\frac{d\stackrel{
(1)i}{X}}{dt}+\stackrel{(1)s}{X}\displaystyle\frac{\omega _s^i}{dt}\right) %
\displaystyle\frac \delta {\delta y^{(1)i}}+\cdots + \\
& + & \left( \displaystyle\frac{d\stackrel{(k-1)i}{X}}{dt}+\stackrel{(k-1)s}{
X}\displaystyle\frac{\omega _s^i}{dt}\right) \displaystyle\frac \delta
{\delta y^{(k-1)i}}+\left( \displaystyle\frac{dX_i}{dt}-X_s\displaystyle
\frac{\omega_i^s}{dt}\right) \displaystyle\frac \delta {\delta p_i}.
\end{array}
\tag{7.7.8}
\end{equation}

The theory of the parallelism of vector fields along curve $\gamma
$, presented in the chapter 10, section 7 of the book [115], in
the case $k=2$, can be extended without difficulties for $k>2$.

Consequently, we define the notion of parallelism of a vector field $X(\gamma
(t)) $ along curve $\gamma $, by the differential equation $%
\displaystyle
\frac{DX}{dt}=0$. We obtain:
\begin{teo}
The vector field $X$, given by (7.7.4') is parallel along the parametrized
curve $\gamma $, with respect to the $N$-linear connection $D$, if and only
if its components $\stackrel{(0)i}{X}$, $\stackrel{(\alpha )i}{X}$, $X_i$ ($%
\alpha =1$, ..., $k-1$) are solutions of the differential equations
\begin{equation}
\begin{array}{l}
\displaystyle\frac{d\stackrel{(0)i}{X}}{dt}+\stackrel{(0)s}{X}\displaystyle
\frac{\omega _s^i}{dt}=0, \\
\\
\displaystyle\frac{d\stackrel{(\alpha )i}{X}}{dt}+\stackrel{(\alpha )s}{X}%
\displaystyle\frac{\omega _s^i}{dt}=0,\ (\alpha =1,...,k-1), \\
\\
\displaystyle\frac{dX_i}{dt}-X_s\displaystyle\frac{\omega _i^s}{dt}=0.
\end{array}
\tag{7.7.9}
\end{equation}
\end{teo}

By means of the formula (7.7.8), the proof of the previous theorem is
immediate.

The vector field $X\in X(T^{*k}M)$ is called \textit{absolute parallel} with
respect to the $N$-linear connection $D$ if the equation $DX=0$ holds for
any curve $\gamma $. This equations $DX=0$, $\forall \gamma $ is equivalent
to the integrability of the following system of Pffaf equations
\begin{equation}
d\stackrel{(0)i}{X}+\stackrel{(0)s}{X}\omega _s^i=0,\ d\stackrel{(\alpha )i}{
X}+\stackrel{(\alpha )s}{X}\omega _s^i=0,\ (\alpha =1,...,k-1),\
dX_i-X_s\omega _i^s=0.  \tag{7.7.10}
\end{equation}

But the previous system is equivalent to
\begin{equation}
\begin{array}{l}
\stackrel{(0)i}{X}_{|j}=\stackrel{(\alpha )i}{X_{|j}}=0,\ (\alpha
=1,...,k-1),\ X_{i|j}=0, \\
\\
\stackrel{(0)i}{X}\stackrel{(\beta )}{|}_j=\stackrel{(\alpha )i}{X}\stackrel{
(\beta )}{|}_j=0,\ (\alpha =1,...,k-1),\ X_i\stackrel{(\beta )}{|}_j=0,\
(\beta =1,...,k-1), \\
\\
\stackrel{(0)i}{X}|^j=\stackrel{(\alpha )i}{X}|^j=0,\ (\alpha =1,...,k-1),\
X_i|^j=0,
\end{array}
\tag{7.7.10a}
\end{equation}
which must be integrable.

Using the Ricci identities (7.6.2) the system (7.7.10') is integrable if and
only if the coordinates ($\stackrel{(0)i}{X}$, $\stackrel{(\alpha )i}{X}$, $%
X_i$) of the vector field $X$ satisfy the following equations:
\begin{equation}
\begin{array}{l}
\stackrel{(\alpha )s}{X}R_{s\ \ jh}^{\ \ i}=0,\ \stackrel{(\alpha
)s}{X} \underset{(\beta )}{P}\text{}_{s\ \ jh}^{\ \ i}=0,\
\stackrel{(\alpha )s}{ X}P_{s\ \ j}^{\ \ i\ \ h}=0,\
\stackrel{(\alpha )s}{X}\underset{(\beta
\gamma )}{S}\text{}_{s\ \ jh}^{\ \ i}=0, \\
\\
\stackrel{(\alpha )s}{X}\underset{(\beta )}{S}\text{}_{s\ \ j}^{\
\ i\ \ h}=0,\ \stackrel{(\alpha )s}{X}S_s^{\ \ ijh}=0,\ (\alpha
=0,1,...,k-1;\beta ,\gamma =1,...,k-1)
\end{array}
\tag{7.7.11}
\end{equation}
and
\begin{equation}
\begin{array}{l}
X^sR_{s\ \ jh}^{\ \ i}=0,\ X_s\underset{(\beta )}{P}\text{}_{i\ \
jh}^{\ \ s}=0,\ X_sP_{i\ \ j}^{\ \ s\ \ h}=0,\ X_s\underset{(\beta
\gamma )}{S}\text{}_{i\ \ jh}^{\ \ s}=0, \\
\\
X_sS\text{}_{i\ \ j}^{\ \ s\ \ h}=0,\ X_sS_i^{\ \ sjh}=0,(\beta ,\gamma
=1,...,k-1).
\end{array}
\tag{7.7.11a}
\end{equation}

The manifold $T^{*k}M$ is called \textit{with absolute parallelism of
vectors, with respect to }$D$ if any vector field on $T^{*k}M$ is absolute
parallel.

In this case the systems of equations (7.7.11), (7.7.11') are verified for any
vector field $X$ with the coefficients ($\stackrel{(0)i}{X}$, $\stackrel{
(\alpha )i}{X}$, $X_i$) in the adapted cobasis.

We obtain:
\begin{teo}
The manifold $T^{*k}M$ is with absolute parallelism of vectors, with respect
to the $N$-linear connection $D$ if and only if all curvature $d$-tensors of
$D$ vanish, i.e.

$R_{m\ jh}^{\ \ \ i}=0,\ \underset{(\alpha )}{P}$ $_{m\ jh}^{\ \ \
i}=0,\ P_{m\ \ j}^{\ \ \ i\ \ h}=0,\ \underset{(\alpha \beta
)}{S}$ $_{m\ \ jh}^{\ \ \ i}=0,\ \underset{(\alpha )}{S}$ $_{m\ \
j}^{\ \ \ i\ \ h}=0,S_m^{\ \ ijh}=0$, $\alpha $, $\beta =1$, ...,
$k-1$.
\end{teo}

The previous theory can be applied to investigate the autoparallel
curves with respect to a $N$-linear connection $D$.

The parametrized curve $\gamma :t\in I\rightarrow \gamma (t)\in T^{*k}M$, is
an autoparellel curve with respect to $D$ if $D_{\stackrel{\cdot }{\gamma }}
\stackrel{\cdot }{\gamma }=0$.

By means of (7.7.2), (7.7.8) we obtain
\begin{equation}
\begin{array}{l}
D_{\stackrel{\cdot }{\gamma }}\stackrel{\cdot }{\gamma }=\displaystyle\frac{%
D \stackrel{\cdot }{\gamma }}{dt}  = \\
 = \left( \displaystyle\frac{d^2x^i}{%
dt^2} +\displaystyle\frac{dx^s}{dt}\displaystyle\frac{\omega _s^i}{dt}%
\right) \displaystyle\frac \delta {\delta x^i}+\left( \displaystyle\frac
d{dt} \displaystyle\frac{\delta y^{(1)i}}{dt}+\displaystyle\frac{\delta
y^{(1)s}}{ dt}\displaystyle\frac{\omega _s^i}{dt}\right) \displaystyle\frac
\delta {\delta y^{(1)i}}+\cdots + \vspace{3mm}\\
+ \left( \displaystyle\frac d{dt}\displaystyle\frac{\delta y^{(k-1)i}}{%
dt }+\displaystyle\frac{\delta y^{(k-1)s}}{dt}\displaystyle\frac{\omega _s^i%
}{dt }\right) \displaystyle\frac \delta {\delta y^{(k-1)i}}+\left( %
\displaystyle \frac d{dt}\displaystyle\frac{\delta p_i}{dt}-\displaystyle%
\frac{\delta p_s}{ dt}\displaystyle\frac{\omega _i^s}{dt}\right) %
\displaystyle\frac \delta {\delta p_i}.
\end{array}
\tag{7.7.12}
\end{equation}

\begin{teo}
A smooth parametrized curve $\gamma $, (7.7.1) is an autoparallel curve with
respect to the $N$-linear connection $D$ if and only if the functions $x^i(t)
$, $y^{(\alpha )i}(t)$, $p_i(t)$, ($\alpha =1$, ..., $k-1$), $t\in I$ verify
the following system of differential equations:
\begin{equation}
\left\{
\begin{array}{l}
\displaystyle\frac{d^2x^i}{dt^2}+\displaystyle\frac{dx^s}{dt}\displaystyle
\frac{\omega _s^i}{dt}=0, \\
\\
\displaystyle\frac d{dt}\displaystyle\frac{\delta y^{(\alpha )i}}{dt}+%
\displaystyle\frac{\delta y^{(\alpha )s}}{dt}\displaystyle\frac{\omega _s^i}{%
dt}=0,\ (\alpha =1,...,k-1), \\
\\
\displaystyle\frac d{dt}\displaystyle\frac{\delta p_i}{dt}-\displaystyle
\frac{\delta p_s}{dt}\displaystyle\frac{\omega _i^s}{dt}=0.
\end{array}
\right.   \tag{7.7.13}
\end{equation}
\end{teo}

Of course, the theorem of existence and uniqueness for the autoparallel
curves can be formulated taking into account the system of differential
equations (7.7.13).

We recall that $\gamma $ is an horizontal curve if $\stackrel{\cdot }{\gamma
}=\stackrel{\cdot }{\gamma }^H$. The horizontal curves are characterized by
\begin{equation}
x^i=x^i(t),\ \displaystyle\frac{\delta y^{(\alpha )i}}{dt}=0,\ (\alpha
=1,...,k-1),\ \displaystyle\frac{\delta p_i}{dt}=0.  \tag{7.7.14}
\end{equation}

\begin{defi}
A horizontal path of an $N$-linear connection $D$ is an horizontal
autoparallel curve $\gamma $, with respect to $D$.
\end{defi}

So, a horizontal path $\gamma $ is characterized by $D_{\stackrel{\cdot }{
\gamma }^H}\stackrel{\cdot }{\gamma }^H=0$.
Taking into account (7.7.13) we get:
\begin{teo}
The horizontal paths of an $N$-linear connection $D$ are characterized by
the system of differential equations
\begin{equation}
\displaystyle\frac{d^2x^i}{dt^2}+H_{jh}^i\displaystyle\frac{dx^j}{dt}%
\displaystyle\frac{dx^h}{dt}=0,\ \displaystyle\frac{\delta y^{(\alpha )i}}{dt%
}=0,\ (\alpha =1,...,k-1),\ \displaystyle\frac{\delta p_i}{dt}=0.  \tag{7.7.15}
\end{equation}
\end{teo}

Indeed, (7.7.14) implies $\displaystyle\frac{\omega _j^i}{dt}=H_{jh}^i %
\displaystyle\frac{dx^h}{dt}$. But (7.7.13) gives us the mentioned equations
(7.7.15).

A parametrized curve $\gamma :I\rightarrow T^{*k}M$ is called\textit{\ }$%
v_\alpha $\textit{-vertical at the point }$x_0\in M$ if $\stackrel{\cdot }{
\gamma }=\stackrel{\cdot }{\gamma }^{V_\alpha }$. It is analytically given by
\[
x^i=x_0^i,\ y^{(\alpha )i}=y^{(\alpha )i}(t),\ y^{(\beta )i}(t)=0,\ \beta
\neq \alpha ,\ p_i=0,\ t\in I.
\]

A $v_\alpha $-vertical path $\gamma $, with respect to $D$ is defined by $D_{%
\stackrel{\cdot }{\gamma }^{V_\alpha }}\stackrel{\cdot }{\gamma
}^{V_\alpha}=0.$

In this case, the equations (7.7.13) are as follows
\begin{equation}
\begin{array}{l}
\displaystyle\frac{dx^i}{dt}=0,\ \displaystyle\frac{dy^{(\beta )i}}{dt}=0,\
(\beta \neq \alpha ),\ \displaystyle\frac{dp_i}{dt}=0\text{ and} \\
\\
\displaystyle\frac d{dt}\displaystyle\frac{\delta y^{(\alpha
)i}}{dt}+ \underset{(\alpha )}{C}\text{}_{sj}^i\displaystyle\frac{dy^{(\alpha )s}}{
dt}\displaystyle\frac{dy^{(\alpha )j}}{dt}=0
\end{array}
\tag{7.7.16}
\end{equation}

Similarly, a $w_k$-vertical curve $\gamma $ at the point $x_0\in M$ is
defined by the condition $\stackrel{\cdot }{\gamma }=\stackrel{\cdot }{\gamma }
^{W_k}$. Analytically it is expressed by
\[
x^i=x_0^i,\ y^{(\alpha )i}=y^{(\alpha )i}(t)=0,\ (\alpha =1,...,k-1),\
p_i=p_i(t),\ t\in I.
\]

A $w_k$-path $\gamma $, with respect to $D$ has the property $D_{\stackrel{
\cdot }{\gamma }^{W_k}}\stackrel{\cdot }{\gamma }^{W_k}=0$.

The $w_k$-paths, with respect to the $N$-linear connection $D$ are
characterized by
\begin{equation}
\begin{array}{l}
\displaystyle\frac{dx^i}{dt}=\displaystyle\frac{dy^{(1)i}}{dt}=\cdots = %
\displaystyle\frac{dy^{(k-1)i}}{dt}=0, \\
\\
\displaystyle\frac{d^2p_i}{dt^2}-C_i^{jm}(x_0,0,...,0,p)\displaystyle\frac{
dp_j}{dt}\displaystyle\frac{dp_m}{dt}=0.
\end{array}
\tag{7.7.17}
\end{equation}

In the case when $D$ is a Berwald $N$-linear connection the previous theory
is a simple one.

\section{Structure Equations of a $N$-Linear Connection}

Let us consider a $N$-linear connection $D$ with the coefficients $D\Gamma
(N)$ in the adapted basis $\left( \displaystyle\frac \delta {\delta x^i}, %
\displaystyle\frac \delta {\delta y^{(1)i}},...,\displaystyle\frac \delta
{\delta y^{(k-1)i}},\displaystyle\frac \delta {\delta p_i}\right) $.

It is not difficult to prove:
\begin{lem}
1$^0$ Each of the following object fields
\[
d(dx^i)-dx^m\wedge \omega _m^i;\ d(\delta y^{(\alpha )i})-\delta y^{(\alpha
)m}\wedge \omega _m^i,\ (\alpha =1,...,k-1);\ d(\delta p_i)+\delta p_m\wedge
\omega _i^m
\]
is a $d$-vector field, except the last one which is a $d$-covector field.

2$^0$ The geometrical object field
\[
d\omega _j^i-\omega _j^m\wedge \omega _m^i
\]
is a $d$-tensor field of type $(1,1)$.
\end{lem}
Using this lemma we obtain a fundamental result in the geometry of the
manifold $T^{*k}M$ and implicitly in the geometry of higher order Hamilton
spaces.
\begin{teo}
For any $N$-linear connection $D$ with the coefficients \newline $D\Gamma (N)=\left(
H_{jh}^i,\underset{(1)}{C}_{jh}^i,...,\underset{(k-1)}{C}%
_{jh}^i,C_i^{jh}\right) $ the following structure equations hold:
\begin{equation}
\begin{array}{l}
d(\delta x^i)-dx^m\wedge \omega _m^i=-\stackrel{(0)}{\Omega ^i}, \\
\\
d(\delta y^{(\alpha )i})-dy^{(\alpha )m}\wedge \omega _m^i=-\stackrel{%
(\alpha )}{\Omega ^i},(\alpha =1,...,k-1), \\
\  \\
d(\delta p_i)+\delta p_m\wedge \omega _i^m=-\Omega _i
\end{array}
\tag{7.8.1}
\end{equation}
and
\begin{equation}
d\omega _j^i-\omega _j^m\wedge \omega _m^i=-\Omega _j^i,  \tag{7.8.2}
\end{equation}
where $\stackrel{(0)}{\Omega ^i}$, $\stackrel{(1)}{\Omega ^i}$, $\stackrel{%
(k-1)}{\Omega ^i}$ and $\Omega _i$ are the $2$-forms of torsion
\begin{equation}
\begin{array}{lll}
\stackrel{(0)}{\Omega ^i} & = & dx^j\wedge \left\{
\displaystyle\frac 12T_{jm}^idx^m+\sum\limits_{\alpha
=1}^{k-1}\underset{(\alpha )}{C}\text{}_{jm}^i\delta y^{(\alpha )m}+C_j^{im}\delta p_m\right\} , \\
\stackrel{(\alpha )}{\Omega ^i} & = & dx^j\wedge \underset{(\alpha 0)}{P}%
\text{}_{\ j}^i+\sum\limits_{\gamma =1}^{k-1}\delta y^{(\gamma
)j}\wedge
\underset{(\alpha \gamma )}{P}\text{}_{\ j}^i+ \\
& + & \delta y^{(\alpha )j}\wedge \left\{
H_{jm}^idx^m+\sum\limits_{\gamma =1}^{k-1}\underset{(\gamma
)}{C}\text{}_{jm}^i\delta y^{(\gamma
)m}+C_j^{im}\delta p_m\right\} ,\ (\alpha =1,...,k-1) \\
\Omega _i & = & dx^j\wedge \left\{ \displaystyle\frac 12\underset{(0)}{R}%
\text{}_{ijm}dx^m+\sum\limits_{\gamma =1}^{k-1}\underset{(\gamma )}{B}%
\text{}_{ijm}\delta y^{(\gamma )m}+\underset{(0)}{B}\text{}_{ij}^{\ \
m}\delta p_m\right\} - \\
& - & \delta p_j\wedge \left\{ H_{im}^jdx^m+\sum\limits_{\gamma =1}^{k-1}%
\underset{(\gamma )}{C}\text{}_{im}^j\delta y^{(\gamma
)m}+C_i^{jm}\delta p_m\right\} ,
\end{array}
\tag{7.8.3}
\end{equation}
$\underset{(\alpha 0)}{P}$ $_j^i$, ..., $\underset{(\alpha ,k-1)}{P}$ $%
_j^i$ being given by (5.4.8) and where $\Omega _j^i$ are the $2$-forms
of curvature:
\begin{equation}
\begin{array}{lll}
\Omega _j^i & = & \displaystyle\frac 12R_{j\ hm}^{\ i}dx^h\wedge
dx^m+\sum\limits_{\gamma =1}^{k-1}\underset{(\gamma )}{P}\text{}_{j\
hm}^{\ i}dx^h\wedge \delta y^{(\gamma )m}+ \\
& + & P_{j\ h}^{\ i\ \ m}dx^h\wedge \delta p_m+\stackrel{k-1}{%
\sum\limits_{\alpha \leq \beta }}\sum\limits_{\beta =1}^{k-1}\underset{%
(\alpha \beta )}{S}\text{}_{j\ hm}^{\ i}\delta y^{(\alpha )h}\wedge \delta
y^{(\beta )m}+ \\
& + & \sum\limits_{\gamma =1}^{k-1}\underset{(\gamma )}{S}\text{}_{j\ h}^{\ i\ \ m}\delta y^{(\gamma )h}\wedge \delta
p_m+\displaystyle\frac 12S_i^{\ jhm}\delta p_h\wedge \delta p_m.
\end{array}
\tag{7.8.4}
\end{equation}
\end{teo}

Indeed, by means of the exterior differential $d\delta y^{(\alpha )i}$
from (7.4.7), (7.4.8) and $\omega _j^i$ from (7.7.6) we get the formulas (7.8.3) and (7.8.4).

These formulas have a very simple form in the case of Berwald
connection, where $\underset{(\alpha )}{C}$ $_{jh}^i=0$,
$C_j^{ih}=0$.

The structure equations will be used in a theory of submanifold of the
Hamilton spaces, studied in chapter 9.

\chapter{Hamilton Spaces of Order $k\geq 1$}

\markboth{\it{THE GEOMETRY OF HIGHER-ORDER HAMILTON SPACES\ \ \ \ \ }}{\it{Hamilton Spaces of Order} $k\geq 1$}

The Hamilton spaces of order $1$ and $2$ were investigated in the chapter $5$
and $12$ of the book [115]. In the present chapter we study the natural
extension of this notion to order $k\geq 1$.

A Hamilton space of order $k$ is a pair $%
H^{(k)n}=(M,H(x,y^{(1)},...,y^{(k-1)},p))$ in which $M$ is a real $n$
-dimensional manifold and $H:T^{*k}M\rightarrow \mathbf{R}$ is a regular
Hamiltonian function on the manifold $T^{*k}M=T^{(k-1)}M\times _MT^{*}M$.

The geometry of the spaces $H^{(k)n}$ can be developed step by step following the
same ideas as in the cases $k=1$ or $k=2$ and using the geometrical theory of
the manifold $T^{*k}M$ described in the last three chapters. Of course, $%
T^{*k}M$ being the dual of $T^kM$, the geometry of the Hamilton spaces of order
$k$, $H^{(k)n}=(M,H),$ appears as dual of the geometry of Lagrange spaces of
order $k$, $L^{(k)n}=(M,L)$, via a Legendre mapping.

Therefore, in this chapter we study the notion of Hamilton space $%
H^{(k)n}=(M,H)$, the canonical presymplectic structure and canonical Poisson
structure, Legendre mappings, the nonlinear connection and canonical
metrical connection. We end with the Riemannian almost contact model of this
space.

\section{The Spaces $H^{(k)n}$}

Let us consider the dual bundle $(T^{*k}M,\pi ^{*k},M)$. The local
coordinates of a point $u=(x,y^{(1)},...,y^{(k-1)},p)$, $u\in T^{*k}M$, will

be denoted as usually by $(x^i,y^{(1)i},...,y^{(k-1)i},p_i)$; $(x^i)$ being
the coordinates of the particle $x$, $y^{(1)i}$, ..., $y^{(k-1)i}$ are
seen as the coordinates of the accelerations of order $1$, ..., $k-1$,
respectively and $p_i$ are the momenta. The coordinate transformations on $%
T^{*k}M$ are given by (4.1.2), (4.1.3).

On the manifold $T^{*k}M$ there are the vertical distributions $%
V_{k-1}\subset V_{k-2}\subset \cdots \subset V_1\subset V$ and a vertical
distribution $W_k$ such that $V_u=V_{1,u}\oplus W_{k,u}$, $\forall u\in
T^{*k}M$.

Also, on the manifold $T^{*k}M$ there exist the Liouville vector fields $%
\stackrel{(1)}{\Gamma }$, ..., $\stackrel{(k-1)}{\Gamma }$ and the Hamilton
vector field $C^{*}$, linearly independent, expressed by
\begin{equation}
\left\{
\begin{array}{l}
\stackrel{(1)}{\Gamma }=y^{(1)i}\displaystyle\frac \partial {\partial
y^{(k-1)i}}, \\
\\
\stackrel{(2)}{\Gamma }=y^{(1)i}\displaystyle\frac \partial {\partial
y^{(k-2)i}}+2y^{(2)i}\displaystyle\frac \partial {\partial y^{(k-1)i}}, \\
......................................................... \\
\stackrel{(k-1)}{\Gamma }=y^{(1)i}\displaystyle\frac \partial {\partial
y^{(1)i}}+2y^{(2)i}\displaystyle\frac \partial {\partial y^{(2)i}}+\cdots
+(k-1)y^{(k-1)i}\displaystyle\frac \partial {\partial y^{(k-1)i}}.
\end{array}
\right.  \tag{8.1.1}
\end{equation}
and
\begin{equation}
C^{*}=p_i\displaystyle\frac \partial {\partial p_i}.  \tag{8.1.1a}
\end{equation}

Theorem 4.2.1 stipulates that these vector fields are globally defined on the
total space of the dual bundle.

The function
\begin{equation}
\varphi =p_iy^{(1)i}  \tag{8.1.2}
\end{equation}
is a scalar function on $T^{*k}M$.

A Hamiltonian is a scalar function
$H:(x,y^{(1)},...,y^{(k-1)},p)\in T^{*k}M\rightarrow
H(x,y^{(1)},...,y^{(k-1)},p)\in \mathbf{R}$. 'Scalar' means that
$H$ does not depend on the changing of coordinates on $T^{*k}M$.

As we know, the Hamiltonian $H$ is differentiable if it is differentiable
on the manifold $\widetilde{T^{*k}M}=T^{*k}M\setminus \{\mathbf{0}\}$ (where
$\mathbf{0}$ is the null section of the projection $\pi ^{*k}$) and $H$ is
continuous on the null section. Evidently,
\[
\widetilde{T^{*k}M}=\left\{ (x,y^{(1)},...,y^{(k-1)},p)\in
T^{*k}M|y^{(1)},...,y^{(k-1)},p\text{ are not all null}\right\} .
\]

The null section $\mathbf{0}:M\rightarrow T^{*k}M$, having the property $\pi
^{*k}\circ \mathbf{0}=1_M$ can be identified with the manifold $M$.

\begin{defi}
A regular Hamiltonian $H:T^{*k}M\rightarrow \mathbf{R}$ is a differentiable
Hamiltonian whose Hessian with respect to the momenta $p_i$, with the
entries:
\begin{equation}
g^{ij}(x,y^{(1)},...,y^{(k-1)},p)=\displaystyle\frac 12\displaystyle\frac{%
\partial ^2H}{\partial p_i\partial p_j}  \tag{8.1.3}
\end{equation}
is nondegenerate on the manifold $\widetilde{T^{*k}M}$.
\end{defi}

Of course, $g^{ij}$ from (8.1.3) is a $d$-tensor field, contravariant of order
$2$, symmetric.

The condition of regularity is expressed by
\begin{equation}
rank||g^{ij}||=n,\quad on\ \widetilde{T^{*k}M}.  \tag{8.1.3a}
\end{equation}

If the base manifold $M$ is paracompact, then the manifold $T^{*k}M$ is
paracompact, too and on $T^{*k}M$ there exist regular Hamiltonians.

The $d$-tensor field $g^{ij}$ being nonsingular on $\widetilde{T^{*k}M}$
there exists a $d$-tensor field $g_{ij}$ covariant of order $2$,
symmetric, uniquely determined, at every point $u\in \widetilde{T^{*k}M}$, by
\begin{equation}
g_{ij}g^{jk}=\delta _i^k.  \tag{8.1.4}
\end{equation}
\begin{defi}
An Hamilton space of order $k$ is a pair\\
$H^{(k)n}=(M,H(x,y^{(1)},...,y^{(k-1)},p))$, where $M$ is a real $n$
-dimensional manifold and $H$ is a differentiable regular Hamiltonian having
the property that the $d$-tensor field $g^{ij}$ has a constant signature on $%
\widetilde{T^{*k}M}$.
\end{defi}

As usually, $H$ is called \textit{the fundamental function} and $g^{ij}$
\textit{the fundamental tensor} of the space $H^{(k)n}$.

In the case when the fundamental tensor $g^{ij}$ is positively defined, then
the condition of regularity (8.1.3') is verified.

\begin{teo}
Assuming that the base manifold $M$ is paracompact, then there exists on $\widetilde{%
T^{*k}M}$ a regular Hamiltonian $H$ such that the pair $(M,H)$ is a Hamilton
space of order $k$.
\end{teo}

\textit{Proof:} Let $F^{(k-1)n}=(M,F(x,y^{(1)},...,y^{(k-1)}))$ be a Finsler
space of order $k-1$ on the manifold $T^{k-1}M$, where $T^{k-1}M=\pi
_{k-1}^{*k}(T^{*k}M)$, having $\gamma _{ij}(x,y^{(1)},...,y^{(k-1)})$ as
fundamental tensor. The manifold $M\,$ being paracompact, the space $%
F^{(k-1)n}$ exists.

Then, the function
\[
H(x,y^{(1)},...,y^{(k-1)},p)=\alpha \gamma
^{ij}(x,y^{(1)},...,y^{(k-1)})p_ip_j,\ (\alpha \in \mathbf{R},\alpha >0),
\]
is well defined in every point $(x,y^{(1)},...,y^{(k-1)},p)\in \widetilde{
T^{*k}M}$ and it is a fundamental function for a Hamilton space of order $k$. Its fundamental tensor field is $\alpha \gamma ^{ij}$. Q.E.D.

One of the important $d$-tensor field derived from the fundamental function $%
H$ of the space $H^{(k)n}$ is:
\begin{equation}
C^{ijh}=-\displaystyle\frac 12\stackrel{\cdot }{\partial }^hg^{ij}=- %
\displaystyle\frac 14\stackrel{\cdot }{\partial }^i\stackrel{\cdot }{
\partial }^j\stackrel{\cdot }{\partial }^hH,\quad \left( \stackrel{\cdot }{
\partial }^i=\displaystyle\frac \partial {\partial p_i}\right) .  \tag{8.1.5}
\end{equation}

\begin{prop}
We have:

1$^0$ $C^{ijh}$ is a totally symmetric $d$-tensor field;

2$^0$ $C^{ijh}$ vanishes if and only if the fundamental tensor $g^{ij}$ does
not depend on the momenta $p_i$.
\end{prop}

Other geometrical object fields which are entirely determined by the
fundamental function $H$ are the coefficients $C_i^{\ jh}$ of the $w_k$
-covariant derivation, given by
\begin{equation}
C_i^{\ jh}=-\displaystyle\frac 12g_{is}\left( \stackrel{\cdot }{\partial }
^j\ g^{sh}+\stackrel{\cdot }{\partial }^h\ g^{js}-\stackrel{\cdot }{\partial }
^s\ g^{jh}\right) .  \tag{8.1.6}
\end{equation}

\begin{prop}
1$^0$ $C_i^{\ jh}$ are the components of a $d$-tensor fields of type $(2,1)$.

2$^0$ They depend on the fundamental function $H$ only.

3$^0$ They are symmetric in the indices $j,h$.

4$^0$ The formula
\[
C_i^{\ jh}=g_{is}C^{sjh}
\]
holds.

5$^0$ The $w_k$-covariant derivative of the fundamental tensor $g^{ij}$,
vanishes:
\begin{equation}
g^{ij}|^h=0.  \tag{8.1.7}
\end{equation}
\end{prop}
The proof is not difficult.

\section{The $k$-Tangent Structure $J$ and the Adjoint $k$-Tangent Structure $J^{*}$}

For the Hamilton space of order $k$, $H^{(k)n}=(M,H),$ the structures $J$ and
$J^{*}$ defined on the manifold $T^{*k}M$ in the section 3, Ch.4 have some
special properties.

The $k$-tangent structure is the mapping:

\begin{equation}
\begin{array}{l}
J\left( \displaystyle\frac \partial {\partial x^i}\right) =\displaystyle %
\frac \partial {\partial y^{(1)i}},\cdots ,J\left( \displaystyle\frac
\partial {\partial y^{(k-2)i}}\right) =\displaystyle\frac \partial {\partial
y^{(k-1)i}},\\
J\left( \displaystyle\frac \partial {\partial y^{(k-1)i}}\right)
=0,J\left( \stackrel{\cdot }{\partial }^i\right) =0.  \tag{8.2.1}
\end{array}
\end{equation}

Locally, it is expressed by (4.3.2):

\begin{equation}
J=\displaystyle\frac \partial {\partial y^{(1)i}}\otimes dx^i+\displaystyle %
\frac \partial {\partial y^{(2)i}}\otimes dy^{(1)i}+\cdots +\displaystyle %
\frac \partial {\partial y^{(k-1)i}}\otimes dy^{(k-2)i}.  \tag{8.2.2}
\end{equation}

The main properties of $J$ are explicitly given in Theorem 4.3.1.

Let $X$ be a vector field on $T^{*k}M$, locally expressed by

\begin{equation}
X=\stackrel{(0)i}{X}\displaystyle\frac \partial {\partial x^i}+\stackrel{%
(1)i }{X}\displaystyle\frac \partial {\partial y^{(1)i}}+\cdots +\stackrel{%
(k-1)i }{X}\displaystyle\frac \partial {\partial y^{(k-1)i}}+X_i\stackrel{%
\cdot }{ \partial }^i.  \tag{8.2.3}
\end{equation}

Here, $\stackrel{\cdot }{\partial }^i=\displaystyle\frac \partial
{\partial p_i}$.

Consider the following vector fields $\stackrel{0}{X}$, $\stackrel{1}{X}$,
..., $\stackrel{k-1}{X}$

\begin{equation}
\stackrel{1}{X}=JX,\ \stackrel{2}{X}=J^2X,\ ...,\ \stackrel{k-1}{X}=J^{k-1}X.
\tag{8.2.3a}
\end{equation}

Taking into account (8.2.3), these vector fields have the form:

\begin{equation}
\begin{array}{l}
\stackrel{1}{X}=\stackrel{(0)i}{X}\displaystyle\frac \partial {\partial
y^{(1)i}}+\cdots +\stackrel{(k-2)i}{X}\displaystyle\frac \partial {\partial
y^{(k-1)i}}, \\
\\
\stackrel{2}{X}=\stackrel{(0)i}{X}\displaystyle\frac \partial {\partial
y^{(2)i}}+\cdots +\stackrel{(k-3)i}{X}\displaystyle\frac \partial {\partial
y^{(k-1)i}}, \\
.................................................. \\
\stackrel{k-1}{X}=\stackrel{(0)i}{X}\displaystyle\frac \partial {\partial
y^{(k-1)i}}.
\end{array}
\tag{8.2.4}
\end{equation}

Now, the adjoint $J^{*}$ of $J$ is defined by

\begin{equation}
J^{*}(dx^i)=0,J^{*}(dy^{(1)})=dx^i,...,J^{*}(dy^{(k-1)i})=dy^{(k-2)i},J^{*}(dp_i)=0.
\tag{8.2.5}
\end{equation}

$J^{*}$ is the following $d$-tensor field of type $(1,1)$:

\begin{equation}
J^{*}=dx^i\otimes \displaystyle\frac \partial {\partial
y^{(1)i}}+dy^{(1)i}\otimes \displaystyle\frac \partial {\partial
y^{(2)i}}+\cdots +dy^{(k-2)i}\otimes \displaystyle\frac \partial {\partial
y^{(k-1)i}}.  \tag{8.2.6}
\end{equation}

$J^{*}$ is an integrable structure and $rank$ $J^{*}=(k-1)n$.

If $\omega $ is an $1$-form field on the manifold $T^{*k}M$ and

\begin{equation}
\omega =\stackrel{(0)}{\omega }_idx^i+\cdots +\stackrel{(k-1)}{\omega }
_idy^{(k-1)i}+\omega ^idp_i,  \tag{8.2.7}
\end{equation}
then by means of $J^{*}$ we obtain a number of $k-1$ $1$-forms on $T^{*k}M$:
\begin{equation}
\stackrel{1}{\omega }=J^{*}\omega ,\ ...,\ \stackrel{k-1}{\omega }
=J^{*(k-1)}\omega .  \tag{8.2.8}
\end{equation}

The vertical differential operators $d_0$, ..., $d_{k-2}$ are introduced in
\S 4.3 by
\begin{equation}
d_0=J^{*(k-1)}d,...,d_{k-2}=J^{*}d,  \tag{8.2.9}
\end{equation}
where $d$ is the operator of differentiation on the manifold $T^{*k}M$. We
know from  \S 4.3 that these operators are the antiderivations of degree $1$ in
the exterior algebra $\Lambda (T^{*k}M)$.

The following formula hold:
\begin{equation}
d\circ d=0,\ d_\alpha \circ d_\alpha =0,\ (\alpha =0,1,...,k-2).  \tag{8.2.10}
\end{equation}

We get:

\begin{prop}
For any Hamilton space of order $k$, $H^{(k)n}=(M,H)$ the $1$-forms (8.2.8) of
the form $dH$ are given by
\begin{equation}
\begin{array}{l}

d_0H=\stackrel{(0)}{p_i}dx^i, \\
\\
d_1H=\stackrel{(1)}{p_i}dx^i+\stackrel{(0)}{p_i}dy^{(1)i}, \\
........................................................................ \\
d_{k-2}H=\stackrel{(k-2)}{p_i}dx^i+\stackrel{(k-3)}{p_i}dy^{(1)i}+\cdots +%
\stackrel{(0)}{p_i}dy^{(k-2)i}.
\end{array}
\tag{8.2.11}
\end{equation}
where
\begin{equation}
\stackrel{(0)}{p_i}=\displaystyle\frac{\partial H}{\partial y^{(k-1)i}},\
\stackrel{(1)}{p_i}=\displaystyle\frac{\partial H}{\partial y^{(k-2)i}},\
...,\stackrel{(k-2)}{p_i}=\displaystyle\frac{\partial H}{\partial y^{(1)i}}.
\tag{8.2.11a}
\end{equation}
\end{prop}

\begin{prop}
The following $2$-forms depend only on the fundamental function $H$ of the
Hamilton space $H^{(k)n}$:
\begin{equation}
\begin{array}{l}
dd_0H=d\stackrel{(0)}{p_i}\wedge dx^i, \\
\\
dd_1H=d\stackrel{(1)}{p_i}\wedge dx^i+d\stackrel{(0)}{p_i}\wedge dy^{(1)i},
\\
.................................................................... \\
dd_{k-2}H=d\stackrel{(k-2)}{p_i}\wedge dx^i+\cdots +d\stackrel{(0)}{p_i}%
\wedge dy^{(k-2)i}.
\end{array}
\tag{8.2.12}
\end{equation}
\end{prop}

We have, also:

\begin{equation}
d_\alpha \circ d_\alpha H=0,\ \forall \alpha =0,1,...,k-2 . \tag{8.2.13}
\end{equation}
Clearly, $dd_\alpha H=0,\ \alpha =0,...,k-2$ are closed $2$-forms.

As we know, the operator
\begin{equation}
d_{k-1}=\displaystyle\frac \partial {\partial x^i}dx^i+\displaystyle\frac
\partial {\partial y^{(1)i}}dy^{(1)i}+\cdots +\displaystyle\frac \partial
{\partial y^{(k-1)i}}dy^{(k-1)i}  \tag{8.2.14}
\end{equation}
is not a vertical differentiation.

\begin{prop}
We have:

1$^0$ $d_{k-1}H$ is not an $1$-form;

2$^0$ Under a change of local coordinate on $T^{*k}M$ we have
\[
d_{k-1}H=d_{k-1}\widetilde{H}+\widetilde{\stackrel{\cdot }{\partial }^j}%
\widetilde{H}\displaystyle\frac{\partial \widetilde{p}_j}{\partial x^i}dx^i;
\]

3$^0$ $J^{*}\circ d_{k-1}=d_{k-2}$.
\end{prop}

\section{The Canonical Poisson Structure of the Hamilton Space $H^{(k)n}$}

Consider a Hamilton space $H^{(k)n}=(M,H)$. As we know from the section 4 ch.4 on the manifold $T^{*k}M$ there exists a canonical
Poisson structure $\{,\}_{k-1}$. Besides this natural Poisson structure on the submanifold $\Sigma_0$ of $T^{*k}M$ there exist a remarkable Poisson structure of the space $H^{(k)n}=(M,H)$.

Let $(T^{*k}M,\overline{\pi }^{*},T^{*}M)$ be the bundle introduced in \S 1, ch.4. The projection $\overline{\pi }^{*}$ is given by $\overline{
\pi }^{*}(x,y^{(1)},...,y^{(k-1)},p)=(x,p)$. The canonical section $\sigma
_0:(x,p)\in T^{*}M\rightarrow (x,0,...,0,p)\in T^{*k}M$ has the image $%
\Sigma _0=Im$ $\sigma _0$, a submanifold of the manifold $T^{*k}M$. The
canonical presymplectic structure $\theta =dp_i\wedge dx^i$ has its
restriction $\theta _0$ to $\Sigma _0$, given by $\theta _0=dp_i\wedge dx^i$
in every point $(x,p)\in \Sigma _0$. The equations of the submanifold $%
\Sigma _0$ being $y^{(\alpha )i}=0$, $(\alpha =1,...,k-1)$, then $(x^i,p_i)$
are the coordinate of the points $(x,p)\in \Sigma _0$.

\begin{teo}
The pair $(\Sigma _0,\theta _0)$ is a symplectic manifold.
\end{teo}

\textit{Proof:} Indeed, $\theta _0=dp_i\wedge dx^i$ is a closed $2$-form and
$rank$ $||\theta _0||=2n=\dim \Sigma _0$. Q.E.D.

In a point $u=(x,p)\in \Sigma _0$ the tangent space $T_u\Sigma _0$ has the
natural basis $\left( \displaystyle\frac \partial {\partial x^i}, %
\displaystyle\frac \partial {\partial p_i}\right) _u$, $(i=1,...,n)$ and
natural cobasis $\left( dx^i,dp_i\right) _u$.

Let us consider $\mathcal{F}(\Sigma _0)$-module $\mathcal{X}(\Sigma _0)$ of
vector fields and $\mathcal{F}(\Sigma _0)$-module $\mathcal{X}^{*}(\Sigma
_0) $ of covector fields on the submanifold $\Sigma _0$.

The following $\mathcal{F}(\Sigma _0)$-linear mapping

$S_{\theta _0}:\mathcal{X}(\Sigma _0)\rightarrow \mathcal{X}^{*}(\Sigma _0)$
defined by
\begin{equation}
S_{\theta _0}(X)=i_X\theta _0,\forall X\in \mathcal{X}(\Sigma _0)  \tag{8.3.4}
\end{equation}
gives us
\begin{equation}
S_{\theta _0}\left( \displaystyle\frac \partial {\partial x^i}\right)
=-dp_i,\ S_{\theta _0}\left( \displaystyle\frac \partial {\partial
p_i}\right) =dx^i.  \tag{8.3.4a}
\end{equation}

These equalities have as a consequence:

\begin{prop}
The mapping $S_{\theta _0}$ is an isomorphism.
\end{prop}

The Hamilton space $H^{(k)n}=(M,H)$, allows to consider the restriction $H_0$
of the fundamental function $H$ to the submanifold $\Sigma _0$, $%
H_0(x,p)=H(x,0,...,0,p)$. Therefore the pair $\left( M,H_0(x,p)\right) $ is
a classical Hamilton space (cf. [115]) with fundamental tensor field
\[
g^{ij}(x,0,...,0,p)=\displaystyle\frac 12\stackrel{\cdot }{\partial }^i
\stackrel{\cdot }{\partial }^jH_0.
\]

By means of the last proposition it follows:

\begin{prop}
1$^0$ There exists a unique vector field $X_{H_0}\in \mathcal{X}(\Sigma _0)$
such that
\[
S_{\theta _0}(X_{H_0})=i_{X_{H_0}}\theta _0=-dH_0.
\]

2$^0$ $X_{H_0}$ is given by
\begin{equation}
X_{H_0}=\displaystyle\frac{\partial H_0}{\partial p_i}\displaystyle\frac
\partial {\partial x^i}-\displaystyle\frac{\partial H_0}{\partial x^i}%
\displaystyle\frac \partial {\partial p_i}.  \tag{8.3.5}
\end{equation}
\end{prop}

\begin{teo}
The integral curves of the vector field $X_{H_0\text{ }}$ are given by the $%
\Sigma _0$-canonical equations
\begin{equation}
\displaystyle\frac{dx^i}{dt}=\displaystyle\frac{\partial H_0}{\partial p_i}%
,\ \displaystyle\frac{dp_i}{dt}=-\displaystyle\frac{\partial H_0}{\partial
x^i},\ y^{(\alpha )i}=0,\ (\alpha =1,...,k-1).  \tag{8.3.6}
\end{equation}
\end{teo}

For two functions $f$, $g\in \mathcal{F}(\Sigma _0)$, let $X_f$ and $X_g$ be
the corresponding Hamilton vector fields given by
\[
i_{X_f}\theta _0=-df,\ i_{X_g}\theta _0=-dg.
\]

\begin{teo}
The following formula holds
\begin{equation}
\left\{ f,g\right\} _0=\theta _0(X_f,X_g),\ \forall f,g\in \mathcal{F}%
(\Sigma _0)  \tag{8.3.7}
\end{equation}
\end{teo}

\textit{Proof}: We have \\
$\theta (X_f,X_g)=(i_{X_f}\theta _0)(X_g)=S_{\theta
_0}(X_f)(X_g)=-df(X_g)=-X_gf=$ \vspace{3mm}\\
$=\left( \displaystyle\frac{\partial f}{\partial x^i}\displaystyle\frac{
\partial g}{\partial p_i}-\displaystyle\frac{\partial f}{\partial p_i} %
\displaystyle\frac{\partial g}{\partial x^i}\right) =\left\{ f,g\right\} _0$
. Q.E.D.

\textbf{Remark} 1$^0$ The previous theory can be extended to the other
Poisson structures $\left\{ ,\right\} _\alpha $ , $(\alpha =1,...,k-1)$ (cf. [4]).

2$^0$ The triple $\left( T^{*k}M,H(x,y^{(1)},...,y^{(k-1)},p),\theta \right)
$ is an Hamiltonian system in which $\theta $ is a presymplectic structure.
Therefore we can apply Gotay's method (cf. M. de Leon and Gotay, [115])
taking into account and the considerations from the previous section, \S 8.2.

The equations (8.3.6) are particular. For a Hamilton space $H^{(k)n}=(M, H),$ the integral of action, (see Ch.5):
$$
I(c)=\int^1_0[p_i\frac{dx^i}{dt} -\frac{1}{2} H(x, \frac{dx}{dt}, ..., \frac{1}{(k-1)!}\frac{d^{k-1}}{dt^{k-1}}, p)]dt
$$
leads, via the variational problem, to the fundamental equations of the space $H^{(k)n}$, i.e. the Hamilton-Jacobi equations
\begin{equation}
\begin{array}{l}
\displaystyle\frac{dx^i}{dt}= \displaystyle\frac{1}{2}\displaystyle\frac{\partial H}{\partial p_i}, \vspace{3mm}\\
\displaystyle\frac{dp_i}{dt}=-\displaystyle\frac{1}{2}[\displaystyle\frac{\partial H}{\partial x^i} - \displaystyle\frac{d}{dt}\displaystyle\frac{\partial H}{\partial y^{(1)i}} + \cdots + (-1)^{k-1}
\displaystyle\frac{d^{k-1}}{dt^{k-1}}\displaystyle\frac{\partial H}{\partial y^{(k-1)i}}].
\end{array}
\tag{8.3.8}
\end{equation}
The energies of order k-1, ${\cal E}^{k-1}(H)$ of the considered space has the expression
$$
{\cal E}^{k-1}(H) = I^{k-1}(H) - \frac{1}{2!}\frac{d}{dt}I^{k-2}(H) + \cdots + (-1)^{k-2}\frac{1}{(k-1)!}\frac{d^{k-2}}{dt^{k-2}}I^1(H) - H
$$
and we have:
\begin{teo}
For a Hamilton space $H^{(k)n}=(M,H)$ the energy of order $k-1$, ${\cal E}^{k-1}(H)$ is constant along every solution curve of the Hamilton-Jacobi equations.
\end{teo}

\section{Legendre Mapping Determined by a Lagrange Space $L^{(k)n}=(M,L)$}

Let $L^{(k)n}=(M,L(x,y^{(1)},...,y^{(k)}))$ be a Lagrange space of order $k.$ It
determines a local diffeomorphism $\varphi :\widetilde{T^kM}\rightarrow
\widetilde{T^{*k}M}$ which preserves the fibres. The mapping $\varphi $
transforms the canonical $k$-semispray $S$, (1.2.5), of $L^{(k)n}$ in the
dual $k$-semispray $S_\xi $, where $\xi =\varphi ^{-1}$ and determines a
nonlinear connection $N^{*}$ on $T^{*k}M$. But $\varphi $ does not transform
the regular Lagrangian $L$ in a regular Hamiltonian, in the case $k>1$. We
need some supplementary geometrical object fields for getting a regular
Hamiltonian $H$ from the fundamental function $L$ of the space $L^{(k)n}$.
In this section we investigate the above mentioned problems.

A point $(x,y^{(1)},...,y^{(k)})$ of the manifold $T^kM$ will be denoted by \\
$(y^{(0)},y^{(1)},...,y^{(k)})$ and its coordinates by $(y^{(0)i},...,y^{(k)i})$.

The fundamental tensor of the space $L^{(k)n}=(M,L(y^{(0)},...,y^{(k)}))$
will be given by
\begin{equation}
a_{ij}=\displaystyle\frac 12\displaystyle\frac{\partial ^2L}{\partial
y^{(k)i}\partial y^{(k)j}}.  \tag{8.4.1}
\end{equation}

\begin{prop}
The mapping $\varphi :u=(y^{(0)},y^{(1)},...,y^{(k)})\in \widetilde{T^kM}%
\rightarrow u^{*}=(x,y^{(1)},...,y^{(k-1)},p)\in \widetilde{T^{*k}M}$ given
by
\begin{equation}
\left\{
\begin{array}{l}
x^i=y^{(0)i},\ y^{(1)i}=y^{(1)i},\ ...,\ y^{(k-1)i}=y^{(k-1)i}, \vspace{3mm}\\
p_i=\displaystyle\frac 12\displaystyle\frac{\partial L}{\partial y^{(k)i}}.
\end{array}
\right.   \tag{8.4.2}
\end{equation}
is a local diffeomorphism, which preserves the fibres.
\end{prop}

\textit{Proof:} The mapping $\varphi $ is differentiable on the manifold $%
\widetilde{T^kM}$ and its Jacobian has the determinant equal to $\det
||a_{ij}||\neq 0$. Of course, $\pi ^k(y^{(0)},...,y^{(k)})=\pi ^{*k}\circ
\varphi (y^{(0)},...,y^{(k)})=y^{(0)}$.

The diffeomorphism $\varphi $ is called \textit{the Legendre mapping} (or
\textit{Legendre transformation}).

We denote
\begin{equation}
p_i=\displaystyle\frac 12\displaystyle\frac{\partial L}{\partial y^{(k)i}}
=\varphi _i(y^{(0)},y^{(1)},...,y^{(k)}).  \tag{8.4.2a}
\end{equation}
Clearly, $\varphi _i$ is a $d$-covector field on $\widetilde{T^kM}$.

The local inverse diffeomorphism $\xi =\varphi ^{-1}:\widetilde{T^{*k}M}
\rightarrow \widetilde{T^kM}$ is expressed by
\begin{equation}
\left\{
\begin{array}{l}
y^{(0)i}=x^i,\ y^{(1)i}=y^{(1)i},\ ...,\ y^{(k-1)i}=y^{(k-1)i}, \\
y^{(k)i}=\xi ^i(x,y^{(1)},...,y^{(k-1)},p).
\end{array}
\right.  \tag{8.4.3}
\end{equation}

With respect to a change of local coordinates on the manifold $T^{*k}M$, $%
\xi ^i$ is transformed exactly as the variables $y^{(k)i}$ from the formulas (1.1.2).

The mappings $\varphi $ and $\xi $ satisfy the conditions
\[
\xi \circ \varphi =1_{\stackrel{\wedge }{U}},\varphi \circ \xi =1_{\stackrel{
\vee}{U}},\stackrel{\wedge }{U}=(\pi ^{*k})^{-1}(U),\stackrel{\vee}{U}
=(\pi ^k)^{-1}(U),U\subset M.
\]

Therefore we have the following identities
\begin{equation}
a_{ij}(y^{(0)},...,y^{(k)})=\displaystyle\frac{\partial \varphi _i}{\partial
y^{(k)j}};\ a^{ij}(x,y^{(1)},...,y^{(k-1)},\xi (x,y^{(1)},...,y^{(k-1)},p))= %
\displaystyle\frac{\partial \xi ^i}{\partial p_j}  \tag{8.4.4}
\end{equation}
and
\begin{equation}
\left\{
\begin{array}{l}
\displaystyle\frac{\partial \varphi _i}{\partial x^j}=-a_{is}\displaystyle
\frac{\partial \xi ^s}{\partial x^j};\ \displaystyle\frac{\partial \varphi
_i }{\partial y^{(\alpha )j}}=-a_{is}\displaystyle\frac{\partial \xi ^s}{
\partial y^{(\alpha )j}},\ (\alpha =1,...,k-1);\ \displaystyle\frac{\partial
\varphi _i}{\partial y^{(k)j}}=a_{ij}; \\
\\
\displaystyle\frac{\partial \xi ^i}{\partial x^j}=-a^{is}\displaystyle\frac{
\partial \varphi _s}{\partial x^j};\ \displaystyle\frac{\partial \xi ^i}{
\partial y^{(\alpha )j}}=-a^{is}\displaystyle\frac{\partial \varphi _s}{
\partial y^{(\alpha )j}},\ (\alpha =1,...,k-1);\ \displaystyle\frac{\partial
\xi ^i}{\partial p_j}=a^{ij}.\
\end{array}
\right.  \tag{8.4.5}
\end{equation}

The differential $d\varphi $ of the mapping $\varphi $,
$\varphi_{*}:T_u(T^kM)\rightarrow T_{\varphi (u)}(T^{*k}M)$ in the natural basis is
expressed by
\begin{equation}
\begin{array}{l}
\varphi _{*}\left( \displaystyle\frac \partial {\partial y^{(0)i}}|_u\right)
=\displaystyle\frac \partial {\partial x^i}|_{u^{*}}+\displaystyle\frac{
\partial \varphi _m}{\partial x^i}\displaystyle\frac \partial {\partial
p_m}|_{u^{*}}\ ,\quad u^{*}=\varphi (u) \\
\\
\varphi _{*}\left( \displaystyle\frac \partial {\partial y^{(\alpha
)i}}|_u\right) =\displaystyle\frac \partial {\partial y^{(\alpha
)i}}|_{u^{*}}+\displaystyle\frac{\partial \varphi _m}{\partial y^{(\alpha
)i} }\displaystyle\frac \partial {\partial p_m}|_{u^{*}}\ ,\quad (\alpha
=1,...,k-1), \\
\\
\varphi _{*}\left( \displaystyle\frac \partial {\partial y^{(k)i}}|_u\right)
=a_{im}\displaystyle\frac \partial {\partial p_m}|_{u^{*}}\ .
\end{array}
\tag{8.4.6}
\end{equation}

\begin{teo}
The Legendre mapping $\varphi :T^kM\rightarrow T^{*k}M$ transforms the
canonical $k$-semispray $S$ of the Lagrange space $L^{(k)n}$
\begin{equation}
S=y^{(1)i}\displaystyle\frac \partial {\partial y^{(0)i}}+2y^{(2)i}%
\displaystyle\frac \partial {\partial y^{(1)i}}+\cdots +ky^{(k)i}%
\displaystyle\frac \partial {\partial
y^{(k-1)i}}-(k+1)G^i(y^{(0)},...,y^{(k)})\displaystyle\frac \partial
{\partial y^{(k)i}}  \tag{8.4.7}
\end{equation}
in the dual $k$-semispray $S_\xi ^{*}$ on $T^{*k}M$:
\begin{equation}
\begin{array}{l}
S_\xi ^{*}=y^{(1)i}\displaystyle\frac \partial {\partial x^i}+2y^{(2)i}%
\displaystyle\frac \partial {\partial y^{(1)i}}+\cdots +(k-1)y^{(k-1)i}%
\displaystyle\frac \partial {\partial y^{(k-2)i}}+ \\
\\
+k\xi ^i(x,y^{(1)},...,y^{(k-1)},p)\displaystyle\frac \partial {\partial
y^{(k-1)}}+\eta _i(x,y^{(1)},...,y^{(k-1)},p)\displaystyle\frac \partial
{\partial p_i},
\end{array}
\tag{8.4.8}
\end{equation}
with the coefficients $\xi ^i$ from (8.4.3) and
\begin{equation}
\begin{array}{l}
\eta _i=-a_{is}\left( \displaystyle\frac{\partial \xi ^s}{\partial x^r}%
y^{(1)r}+\cdots +(k-1)\displaystyle\frac{\partial \xi ^s}{\partial y^{(k-1)r}%
}\xi ^r\right. + \\
\\
+\left. (k+1)G^s(x,y^{(1)},...,y^{(k-1)},\xi
(x,y^{(1)},...,y^{(k-1)},p))\right)
\end{array}
\tag{8.4.9}
\end{equation}
\end{teo}

\textit{Proof}: If $X\in \mathcal{X}(\widetilde{T^kM})$ have the local
expression

$X=X^{(0)i}\displaystyle\frac \partial {\partial y^{(0)i}}+X^{(1)i} %
\displaystyle\frac \partial {\partial y^{(1)i}}+\cdots +X^{(k)i}%
\displaystyle
\frac \partial {\partial y^{(k)i}}$ in every point $u\in \widetilde{T^kM}$,
then

$X^{*}=\left( \varphi _{*}X\right) (u^{*})=X^{(0)i}(u^{*})\varphi _{*} %
\displaystyle\frac \partial {\partial x^i}+\cdots +X^{(k)i}(u^{*})\varphi
_{*}\displaystyle\frac \partial {\partial y^{(k)i}}$ $=$

$=X^{(0)i}(u^{*})\displaystyle\frac \partial {\partial x^i}|_{u^{*}}+\cdots
+X^{(k-1)i}(u^{*})\displaystyle\frac \partial {\partial y^{(k-1)i}}|_{u*}+$

$+\left( X^{(0)m}(u^{*})\displaystyle\frac{\partial \varphi _i}{\partial x^m}
+\cdots +X^{(k-1)m}(u^{*})\displaystyle\frac{\partial \varphi _i}{\partial
y^{(k-1)m}}\right) \displaystyle\frac \partial {\partial p_i}|_{u^{*}}+$

$+X^{(k)m}(u^{*})a_{mi}(u^{*})\displaystyle\frac \partial {\partial
p_i}|_{u^{*}}$ in the points $u^{*}=\varphi (u)$.

Consequently, $S_\xi ^{*}=\varphi _{*}S$ holds.

Now, applying Theorem 5.3.1 we obtain without difficulties:

\begin{teo}
The Legendre mapping $\varphi :T^kM\rightarrow T^{*k}M$ determined by a
Lagrange space of order $k$, $L^{(k)n}$, transforms the canonical $k$
-semispray $S$ in the dual $k$-semispray $S_\xi ^{*}$ with the coefficients $\xi ^i$, $\eta_i$.
$S_\xi ^{*}$ determines locally a nonlinear connection $N^{*}$ on $%
T^{*k}M$ with dual coefficients:
\begin{equation}
\underset{(1)}{M}^{*}\text{ }_j^i=-\displaystyle\frac{\partial \xi ^i}{%
\partial y^{(1)j}},...,\underset{(k-1)}{M^{*}}\text{ }_j^i=-\displaystyle
\frac{\partial \xi ^i}{\partial y^{(k-1)j}}  \tag{8.4.10}
\end{equation}
and
\begin{equation}
N_{ij}^{*}=\displaystyle\frac{\delta \eta _i}{\delta y^{(1)j}},  \tag{8.4.10a}
\end{equation}
where the operators $\displaystyle\frac \delta {\delta y^{(1)i}}$ are
obtained by means of the coefficients, (8.4.10).
\end{teo}

We may ask wether by means of the Legendre mapping $\varphi $ we
can transfer the fundamental function $L$ of the Lagrange space $L^{(k)n}$
to a fundamental function $H$ of a Hamilton space of order $k$, $H^{(k)}$,
like in the classical case of $k=1$. Remarking that the energy of order $k$,
$\mathcal{E}_c^k(L)$ of the space $L^{(k)n}$ is not a regular Lagrangian on $T^kM$, the Legendre
transformation $\varphi $ cannot apply $\mathcal{E}_c^k(L)$ in a regular
Hamiltonian on $T^{*k}M$. Therefore it is necessary to look for another way to
solve this problem.

Let us consider a fixed nonlinear connection $\stackrel{\circ }{N}$ on the
submanifold $T^{k-1}M$ in $T^kM$. The dual coefficients of $\stackrel{\circ
}{N}$ are denoted by $\stackrel{\circ }{\underset{(1)}{M}}$ $%
_j^i(x,y^{(1)},...,y^{(k-1)})$, ..., $\stackrel{\circ
}{\underset{(k-1)}{M} }$ $_j^i(x,y^{(1)},...,y^{(k-1)})$. It is
not difficult (cf. (1.4.7)) to prove that
\begin{equation}
kz^{(k)i}=ky^{(k)i}+(k-1)\stackrel{\circ
}{\underset{(1)}{M}}\text{ } _m^iy^{(k-1)m}+\cdots
+\stackrel{\circ }{\underset{(k-1)}{M}}\text{ } _m^iy^{(1)m}
\tag{8.4.11}
\end{equation}
is a $d$-vector field at every point $u=(x,y^{(1)},...,y^{(k)})$.

The Legendre mapping
$$
\varphi :u\in T^kM\rightarrow u^{*}\in T^{*k}M,
u^{*}=(x,y^{(1)},...,y^{(k-1)},p) \textrm{ transform } z^{(k)i}
$$
in a $d$-vector
field $\stackrel{\vee}{z}^{(k)i}$at $u^{*}$, given by

\begin{equation}
k\stackrel{\vee}{z}^{(k)i}=k\xi ^i+(k-1)\stackrel{\circ }{\underset{(1)}{%
M }}\text{ }_m^iy^{(k-1)m}+\cdots +\stackrel{\circ
}{\underset{(k-1)}{M}} \text{ }_m^iy^{(1)m}  \tag{8.4.12}
\end{equation}

Let us consider the following Hamiltonian
\begin{equation}
H(x,y^{(1)},...,y^{(k-1)},p)=2p_i\stackrel{\vee}{z}
^{(k)i}-L(x,y^{(1)},...,y^{(k-1)},\xi (x,y^{(1)},...,y^{(k-1)},p)).
\tag{8.4.13}
\end{equation}

Clearly, $H$ is a differentiable Hamiltonian function on $T^{*k}M$.
\begin{teo}
The Hamiltonian function $H$, (8.4.13), is the fundamental function of a
Hamilton space of order $k$, $H^{(k)n}$ and its fundamental tensor field is
\begin{equation}
g^{ij}(x,y^{(1)},...,y^{(k-1)},p)=a^{ij}(x,y^{(1)},...,y^{(k-1)},\xi
(x,y^{(1)},...,y^{(k-1)},p)),  \tag{8.4.14}
\end{equation}
$a_{ij}$ being the fundamental tensor field of the Lagrange space of order $k
$, $L^{(k)n}=(M,L)$.
\end{teo}

\textit{Proof: }From (8.4.13) we have \\
$\displaystyle\frac 12\stackrel{\cdot }{\partial }^jH=\stackrel{\vee}{z}
^{(k)j}+p_m\stackrel{\cdot }{\partial }^j\stackrel{\vee}{z}^{(k)m}- %
\displaystyle\frac 12\displaystyle\frac{\partial L}{\partial y^{(k)m}}
\stackrel{\cdot }{\partial }^j\xi ^{(k)m}=$\vspace{3mm}\\
$=\stackrel{\vee}{z}^{(k)j}+p_m\stackrel{\cdot }{\partial }^j\stackrel{\vee
}{z}^{(k)m}-p_m\stackrel{\cdot }{\partial }^j\stackrel{\vee}{z}^{(k)m}=
\stackrel{\vee}{z}^{(k)j}$. \vspace{3mm}\\
Consequently,\\
$$
g^{ij}(x,y^{(1)},...,y^{(k-1)},p)=\displaystyle\frac 12\stackrel{\cdot }{
\partial }^i\stackrel{\cdot }{\partial }^jH=\displaystyle\frac 12\stackrel{
\cdot }{\partial }^i\stackrel{\vee}{z}^{(k)j}=$$
$$=\displaystyle\frac 12\stackrel{\cdot }{\partial }^i\xi
^{(k)j}=a^{ij}(x,y^{(1)},...,y^{(k-1)},\xi ).$$

It follows that $H$ is a regular Hamiltonian and its fundamental tensor $%
g^{ij}$ has a constant signature on $T^{*k}M$. Q.E.D.

The Hamilton space of order $k$, $H^{(k)n}=(M,H)$ is called the \textit{dual} of
the Lagrange space of order $k$, $L^{(k)n}=(M,L)$.

\textbf{Remarks.} All geometrical object fields of the space $H^{(k)n}$,
which are derived from the fundamental tensor field $g^{ij}$ do not depend on the
apriori fixed nonlinear connection $\stackrel{\circ }{N}$. Consequently, the
geometry of the Hamilton space $H^{(k)n}=(M,H)$ can be constructed by means of
the dual $k$-semispray $S_\xi ^{*}$ and the fundamental tensor field $g^{ij}$
from (8.4.14).

\section{Legendre Mapping Determined by a Hamilton Space of Order $k$}

Let us consider a converse for the previous problem: being given a Hamilton space of order $k$, $%
H^{(k)n}=(M,H(x,y^{(1)},...,y^{(k-1)},p))$ let us determine its \textit{dual}
like a Lagrange space of order $k$, $L^{(k)n}=(M,L(x,y^{(1)},...,y^{(k)}))$.

In this case, we start from the space $H^{(k)n}$ and try to determine the
local diffeomorphism in the form of $\varphi ^{-1}$ from (8.4.3). But $%
y^{(k)i} $ is not a vector field on $T^kM$ and we can not define it only by
the $d$-vector field $\stackrel{\cdot }{\partial }^iH$. As in the previous
section we assume that a nonlinear connection $\stackrel{\circ }{N}$ on the
submanifold $T^{k-1}M$ in $T^kM$ is apriori fixed. The dual coefficients of $%
\stackrel{\circ }{N}$ being $\stackrel{\circ }{\underset{(1)}{M}}$
$_j^i$, ..., $\stackrel{\circ }{\underset{(k-1)}{M}}$ $_j^i$ they are
functions of the points $(x,y^{(1)},...,y^{(k-1)})$ from $T^{k-1}M$.

Then the $d$-vector field $z^{(k)i}$, given by (8.4.11) on $T^kM$ can be
considered.

The mapping

$\stackrel{*}{\xi }:u^{*}=(x,y^{(1)},...,y^{(k-1)},p)\in T^{*k}M\rightarrow
u=(y^{(0)},y^{(1)},...,y^{(k-1)},y^{(k)})\in T^kM$,
defined as follows:
\begin{equation}
\begin{array}{l}
y^{(0)i}=x^i,y^{(1)i}=y^{(1)i},...,y^{(k-1)i}=y^{(k-1)i}, \\
\\
y^{(k)i}=\stackrel{*}{\xi }^i(x,y^{(1)},...,y^{(k-1)},p),
\end{array}
\tag{8.5.1}
\end{equation}
where $\stackrel{*}{\xi }^i$is expressed from the formula:
\begin{equation}
\begin{array}{l}
k\stackrel{*}{\xi }^i(x,y^{(1)},...,y^{(k-1)},p)+(k-1)\stackrel{\circ }{%
\underset{(1)}{M}}_m^i(x,y^{(1)},...,y^{(k-1)})y^{(k-1)m}+ \\
\\
+\cdots +\stackrel{\circ }{\underset{(k-1)}{M}}
_m^i(x,y^{(1)},...,y^{(k-1)})y^{(1)m}=\displaystyle\frac
k2\stackrel{\cdot }{
\partial }^iH.
\end{array}
\tag{8.5.2}
\end{equation}

\begin{teo}
The mapping $\stackrel{*}{\xi }$ is a local diffeomorphism, which preserves
the fibres of $T^{*k}M$ and $T^kM$.
\end{teo}

\textit{Proof}: The determinant of the Jacobian of $\stackrel{*}{\xi }$ is \newline
$%
\det ||g^{ij}(x,y^{(1)},...,y^{(k-1)},p)||\neq 0$, since $\stackrel{\cdot }{
\partial }^j\stackrel{*}{\xi }^i=g^{ij}$ and $\pi ^{*k}=\pi ^k\circ
\stackrel{*}{\xi }$. q.e.d.

The formulae (8.5.1) and (8.5.2) imply:
\begin{equation}
\stackrel{\cdot }{\partial }^i\stackrel{*}{\xi }^j=g^{ij}  \tag{8.5.3}
\end{equation}
\begin{equation}
z^{(k)i}(x,y^{(1)},...,y^{(k-1)},\stackrel{*}{\xi })=\displaystyle\frac 12
\stackrel{\cdot }{\partial }^iH(x,y^{(1)},...,y^{(k-1)},p),  \tag{8.5.4}
\end{equation}
$g^{ij}$ being the fundamental tensor field of space $H^{(k)n}$.

The inverse mapping $\stackrel{*}{\varphi }
:u=(y^{(0)},y^{(1)},...,y^{(k-1)},y^{(k)})\in T^kM\rightarrow
u^{*}=(x,y^{(1)},...,y^{(k-1)},p)\in T^{*k}M$ of the Legendre transformation
$\stackrel{*}{\xi }$ can be written as follows:

\begin{equation}
\begin{array}{l}
x^i=y^{(0)i},y^{(1)i}=y^{(1)i},...,y^{(k-1)i}=y^{(k-1)i}, \\
\\
p_i=\stackrel{*}{\varphi }^i(y^{(0)},y^{(1)},...,y^{(k-1)},y^{(k)}).
\end{array}
\tag{8.5.5}
\end{equation}
>From here we deduce
\begin{equation}
\displaystyle\frac{\partial \stackrel{*}{\varphi _i}}{\partial y^{(k)j}}
=g_{ij}(x,y^{(1)},...,y^{(k-1)},\stackrel{*}{\varphi }(u)),  \tag{8.5.6}
\end{equation}
where $g_{ij}$ is the covariant tensor of the fundamental tensor $g^{ij}$ of
$H^{(k)n}$.

The Legendre transformation $\stackrel{*}{\xi }$ maps the regular Hamiltonian $H$ into the differentiable Lagrangian $L=\stackrel{*}{
\xi }(H)$:
\begin{equation}
\begin{array}{l}
L(x,y^{(1)},...,y^{(k-1)},y^{(k)})=2p_iz^{(k)i}-H(x,y^{(1)},...,y^{(k-1)},p),
\\
\\
p_i=\stackrel{*}{\varphi }_i(y^{(0)},y^{(1)},...,y^{(k-1)},y^{(k)}).
\end{array}
\tag{8.5.7}
\end{equation}

\begin{teo}
The Lagrangian $L$ from (8.5.7) is a regular one. Its fundamental tensor field
is $g_{ij}(x,y^{(1)},...,y^{(k-1)},\stackrel{*}{\varphi }%
^i(x,y^{(1)},...,y^{(k-1)},y^{(k)}))$.
\end{teo}

Indeed, $\displaystyle\frac{\partial z^{(k)i}}{\partial y^{(k)j}}=\delta
_j^i $ implies $\displaystyle\frac 12\displaystyle\frac{\partial L}{\partial
y^{(k)i}}=\displaystyle\frac{\partial p_m}{\partial y^{(k)i}}z^{(k)m}+p_i- %
\displaystyle\frac 12\stackrel{\cdot }{\partial }^mH\displaystyle\frac{
\partial p_m}{\partial y^{(k)i}}$.

Taking into account (8.5.4) it follows
\begin{equation}
p_i=\displaystyle\frac 12\displaystyle\frac{\partial L}{\partial y^{(k)i}},\
\ p_i=\stackrel{*}{\varphi }_i(x,y^{(1)},...,y^{(k-1)},y^{(k)}).  \tag{8.5.8}
\end{equation}

By means of (8.5.6) it follows
\begin{equation}
a_{ij}(u)=\displaystyle\frac 12\displaystyle\frac{\partial ^2L}{\partial
y^{(k)i}\partial y^{(k)j}}(u)=g_{ij}(x,y^{(1)},...,y^{(k-1)},\stackrel{*}{
\varphi }(u)),  \tag{8.5.9}
\end{equation}
where $u=(x,y^{(1)},...,y^{(k)})\in T^kM$. Q.E.D.$\quad $

The Lagrange space of order $k$, $L^{(k)n}=(M,L)$ is called the \textit{dual}
of the Hamilton space of order $k$.

\textbf{Remarks} 1$^0$ Of course, the space $L^{(k)n}=(M,L)$ is locally
determined.

2$^0$ The relations (8.5.8) give us the geometrical meaning of the mapping $%
\stackrel{*}{\varphi }$, the inverse of the bundle morphism $\stackrel{*}{
\xi }$.

It follows

\begin{teo}
We have the following local properties:

Consider the Legendre transformation determined by the Lagrange space of order $k$, $%
L^{(k)n}=(M,L)$. Then, the dual Hamilton space of order $k$, $H^{(k)n}=(M,H)
$, $L=\stackrel{*}{\xi }(H)$ is defined by the local diffeomorphism $%
\stackrel{*}{\varphi }$, the inverse of the local diffeomorphism $\stackrel{*}{%
\xi }$ and $H=\stackrel{*}{\varphi }(L)$ is locally given by
\begin{equation}
\begin{array}{l}
H(x,y^{(1)},...,y^{(k-1)},p)=2p_iz^{(k)}(x,y^{(1)},...,y^{(k-1)},y^{(k)})-\\
\hfill -L(x,y^{(1)},...,y^{(k-1)},y^{(k)}),
\\
\\
\ y^{(k)i}=\stackrel{*}{\xi }^i(x,y^{(1)},...,y^{(k-1)},p).
\end{array}
\tag{8.5.10}
\end{equation}
\end{teo}

Indeed, as the Legendre transformation $\varphi :T^kM\rightarrow T^{*k}M$ is
defined by \\
$\varphi :(x,y^{(1)},...,y^{(k-1)},y^{(k)})\in T^kM\rightarrow
(x,y^{(1)},...,y^{(k-1)},p)\in T^{*k}M$, with $p_i=\displaystyle\frac 12 %
\displaystyle\frac{\partial L}{\partial y^{(k)i}}=\stackrel{*}{\varphi }
(x,y^{(1)},...,y^{(k)})$. It follows, locally, $\varphi =\stackrel{*}{
\varphi }$, and $\varphi ^{-1}=\stackrel{*}{\xi }$. But $H$ from (8.5.10) is $%
H=\stackrel{*}{\varphi }(L)$. So, locally, $H=\stackrel{*}{\varphi }(
\stackrel{*}{\xi }H)$ and $L=\stackrel{*}{\xi }(\stackrel{*}{\varphi }L)$
.Q.E.D.

\section{The Canonical Nonlinear Connection of the Space $H^{(k)n}$}

There exists a nonlinear connection $N^{*}$ on the manifold $T^{*k}M$, which
is determined only by a Hamilton space of order $k$, $H^{(k)n}=(M,H)$.

Such a property holds in the case of Lagrange spaces of order $k$, $%
L^{(k)n}=(M,L)$, (see ch. 2). Namely, $L^{(k)n}$ determines a canonical $k$
-semispray $S$ and $S$ allows to construct a nonlinear connection $N$, which
depends only on the fundamental function $L$ (Theorem 2.5.2.). Since $L^{(k)n}$
defines the Legendre mapping $\varphi $ from (8.4.2), (8.4.2'), then $\varphi
^{-1}=\xi $ determines a bundle morphism and $S$ is transformed in a $k$
-dual semispray $S_\xi ^{*}=\varphi _{*}(S)$. Theorem 8.4.2. shows that $%
S_\xi ^{*}$ determine a nonlinear connection $N^{\ast}$ on $T^{*k}M$, which depends
on the fundamental function $L$ of the space $L^{(k)n}$.

Considering the Hamilton space $H^{(k)n}=(M,H)$, the dual of $L^{(k)n}$, its
fundamental function $H=\varphi (L)$ is built locally by means of $L$ and
by an apriori given nonlinear connection $\stackrel{\circ }{N}$ on the
submanifold $T^{k-1}M$. The connection $N^{*}$ can be considered the $%
\stackrel{\circ }{N}$-canonical nonlinear connection of the space $H^{(k)n}$.

Assuming now that a Hamilton space of order $k$, $H^{(k)n}=(M,H)$ is
considered, we construct a bundle morphism $\stackrel{*}{\xi }$, (8.5.1),
(8.5.2), by means of $H$ and $\stackrel{\circ }{N}$, and we determine locally the
Lagrange space $L^{(k)n}=(M,L)$ , where $L=\stackrel{*}{\xi }(H)$ is from
(8.5.7). It is the \textit{dual} space of $H^{(k)n}$. The Legendre transformation $%
\stackrel{*}{\varphi }=\stackrel{*}{\xi }^{-1}$, with $\stackrel{*}{\varphi }
_i=\displaystyle\frac 12\displaystyle\frac{\partial L}{\partial y^{(k)i}}$,
transforms the $k$-semispray $S$ of $L$ from (8.4.7) in the dual $k$-semispray $%
S_\xi ^{*}$ from (8.4.8):

\begin{equation}
S_\xi ^{*}=y^{(1)i}\displaystyle\frac \partial {\partial x^i}+\cdots
+(k-1)y^{(k-1)i}\displaystyle\frac \partial {\partial y^{(k-2)i}}+k\stackrel{
*}{\xi }^i\displaystyle\frac \partial {\partial y^{(k-1)i}}+\eta _i %
\displaystyle\frac \partial {\partial p_i},  \tag{8.6.1}
\end{equation}
where $\stackrel{*}{\xi }^i$ is given by (8.5.2) and its coefficients are

\begin{equation}
\begin{array}{l}
\eta _i=-g_{is}\big[ \displaystyle\frac{\partial
\stackrel{*}{\xi}^s}{\partial
x^r}y^{(1)r}+\cdots+(k-1)\displaystyle\frac{\partial \stackrel{*}{\xi }^s}{\partial y^{(k-1)r}}\stackrel{*}{\xi }^r+\\
(k+1)G^s(x,...,y^{(k-1)},\stackrel{*}{\xi}(u^{*})) \big]
\end{array}
\tag{8.6.2}
\end{equation}
So, we have obtained:

\begin{teo}
The dual coefficients of the $\stackrel{\circ }{N}$-canonical nonlinear
connection $N^{*}$ of the Hamilton space of order $k$, $H^{(k)n},$ are
expressed as follows:
\begin{equation}
\underset{(1)}{\stackrel{*}{M}}\text{}_j^i=-\displaystyle\frac{\partial
\stackrel{*}{\xi }^i}{\partial y^{(1)j}},...,\underset{(k-1)}{\stackrel{*}{%
M}}\text{ }_j^i=-\displaystyle\frac{\partial \stackrel{*}{\xi }^i}{\partial
y^{(k-1)j}}  \tag{8.6.3}
\end{equation}
and
\begin{equation}
\stackrel{\ast }{N}_{ij}=\displaystyle\frac{\delta \eta _i}{\delta y^{(1)j}},
\tag{8.6.3a}
\end{equation}
where the operator $\displaystyle\frac \delta {\delta y^{(1)j}}$ is
constructed by means of the coefficients \newline  $\underset{(1)}{N}\text{}_j^i, ..., \underset{(k)}{N}\text{}_j^i$ given by the dual coefficients (8.6.3).
\end{teo}

The previous theory was presented in the book [115], ch. 12, in the case $k=2$.

It is useful to prove the existence of a nonlinear connection canonically
determined by the fundamental function of the Hamilton space $H^{(k)n}$.

\section{Canonical Metrical $N$-Linear Connection of the Space $H^{(k)n}$}

For a Hamilton space of order $k$, $H^{(k)n}=(M,H)$, we consider the
canonical nonlinear connection $N$ determined in the previous section.

We consider also the adapted basis $\left( \displaystyle\frac \delta
{\delta x^i},\displaystyle\frac \delta {\delta y^{(1)i}},\cdots , %
\displaystyle\frac \delta {\delta y^{(k-1)i}},\stackrel{\cdot }{\partial }
^i\right) $ and its dual basis $\left( \delta x^i,\delta y^{(1)i},...,\delta
y^{(k-1)i},\delta p_i\right) $ determined by $N$ and by the distribution $%
W_k$ (cf. ch. 6).

We can define a $N$-linear connection $D$ by its coefficients
$D\Gamma (N)=(H_{jh}^i,\underset{(\alpha )}{C}$
$_{jh}^i,C_i^{jh})$, ($\alpha =1,2,...,k-1$). Therefore the $h$-,
$v_\alpha $-,$w_k$-operators of
covariant derivations with respect to $D$, denoted by '$_{|}$', '$%
\stackrel{(\alpha )}{|}$', '$|$', can be applied to the $d$-tensor
fields on the manifold $T^{*k}M$.

Therefore if $g^{ij}=\displaystyle\frac 12\stackrel{\cdot }{\partial }^i
\stackrel{\cdot }{\partial }^jH$ is the fundamental tensor of the Hamilton
space $H^{(k)n}$, we have
\begin{equation}
\begin{array}{l}
g_{\quad |h}^{ij}=\displaystyle\frac{\delta g^{ij}}{\delta x^h}
+g^{mj}H_{mh}^i+g^{im}H_{mh}^j, \\
\\
g^{ij}\stackrel{(\alpha )}{|}_h=\displaystyle\frac{\delta
g^{ij}}{\delta y^{(\alpha )h}}+g^{mj}\underset{(\alpha )}{C}\text{}_{mh}^i+g^{im}
\underset{(\alpha )}{C}\text{}_{mh}^j,\ (\alpha =1,...,k-1), \\
\\
g^{ij}|^h=\displaystyle\frac{\partial g^{ij}}{\partial p_h}
+g^{mj}C_m^{ih}+g^{im}C_m^{jh}.
\end{array}
\tag{8.7.1}
\end{equation}

A $N$-linear connection $D\Gamma (N)$ is called \textit{compatible} with the
fundamental tensor $g^{ij}$ of the Hamilton space of order $k$, $%
H^{(k)n}=(M,H)$, (or it is \textit{metrical}) if $g^{ij}$ is covariant constant
(or absolute parallel) with respect to $D\Gamma (N)$, i.e.

\begin{equation}
g_{\quad |h}^{ij}=0,\ g^{ij}\stackrel{(\alpha )}{|}_h=0,\ g^{ij}|^h=0.
\tag{8.7.2}
\end{equation}

The previous equations have a geometrical meaning. Indeed, if the tensor field $g^{ij}$ is positively defined and $%
\omega _i$ is a $d$-covector field, then
\[
||\omega ||^2=g^{ij}\omega _i\omega _j
\]

is a scalar field. Then, locally, $\displaystyle\frac d{dt}||\omega ||^2=0$,
along any curve $\gamma $ and for any parallel covector $\omega _i$ along
$\gamma $, if and only if the equations (8.7.2) hold.

As usual, we can prove:

\begin{teo}
1) In a Hamilton space of order $k$, $H^{(k)n}=(M,H)$, there exists a
unique $N$-linear connection $D$, with the coefficients $D\Gamma
(N)=$ \newline $=(H_{jh}^i,\underset{(1)}{C}_{jh}^i,...,\underset{(k-1)}{C}$ $%
_{jh}^i,C_i^{jh})$ verifying the axioms:

1$^0$ $N$ is the canonical nonlinear connection of $H^{(k)n}$.

2$^0$ The fundamental tensor $g^{ij}$ is $h$-covariant constant:
\begin{equation}
g_{\quad |h}^{ij}=0.  \tag{8.7.3}
\end{equation}

3$^0$ $g^{ij}$ is $v_\alpha $-covariant constant:
\begin{equation}
g^{ij}\stackrel{(\alpha )}{|}_h=0,\ (\alpha =1,...,k-1).  \tag{8.7.4}
\end{equation}

4$^0$ $g^{ij}$ is $w_k$-covariant constant:
\begin{equation}
\ g^{ij}|^h=0.  \tag{8.7.4a}
\end{equation}

5$^0$ $D\Gamma (N)$ is $h$-torsion free:
\begin{equation}
T_{jh}^i=H_{jh}^i-H_{hj}^i=0.  \tag{8.7.5}
\end{equation}

6$^0$ $D\Gamma (N)$ is $v_\alpha $-torsion free:
\begin{equation}
\underset{(\alpha )}{S}\text{ }_{jh}^i=\underset{(\alpha )}{C}\text{ }%
_{jh}^i-\underset{(\alpha )}{C}\text{ }_{hj}^i=0,\ (\alpha
=1,...,k-1). \tag{8.7.6}
\end{equation}

7$^0$ $D\Gamma (N)$ is $w_k$-torsion free:
\begin{equation}
S_i^{jh}=C_i^{jh}-C_i^{hj}=0.  \tag{8.7.7}
\end{equation}

2) The connection $D\Gamma (N)$ has the coefficients given by the
generalized Christoffel symbols:
\begin{equation}
\begin{array}{l}
H_{jh}^i=\displaystyle\frac 12g^{is}\left( \displaystyle\frac{\delta g_{sh}}{%
\delta x^j}+\displaystyle\frac{\delta g_{js}}{\delta x^h}-\displaystyle\frac{%
\delta g_{jh}}{\delta x^s}\right) , \\
\\
\underset{(\alpha )}{C}\text{ }_{jh}^i=\displaystyle\frac 12g^{is}\left( %
\displaystyle\frac{\delta g_{sh}}{\delta y^{(\alpha )j}}+\displaystyle\frac{%
\delta g_{js}}{\delta y^{(\alpha )h}}-\displaystyle\frac{\delta g_{jh}}{%
\delta y^{(\alpha )s}}\right) ,\ (\alpha =1,...,k-1), \\
\\
C_i^{jh}=-\displaystyle\frac 12g_{is}\left( \stackrel{\cdot }{\partial }%
^jg^{sh}+\stackrel{\cdot }{\partial }^hg^{js}-\stackrel{\cdot }{\partial }%
^sg^{jh}\right) .
\end{array}
\tag{8.7.8}
\end{equation}

3) This connection depends only by the fundamental function $H$ of the space
$H^{(k)n}$ and by the canonical nonlinear connection $N$.
\end{teo}

The connection $D\Gamma (N)$ with the coefficients (8.7.8) will be denoted by $%
C\Gamma (N)$ and called \textit{canonical }for the space $H^{(k)n}$.

The canonical $N$-metrical connection $C\Gamma (N)$ has zero torsions $%
T_{jh}^i$, $\underset{(\alpha )}{S}$ $_{jh}^i$, $S_i^{jh}$. Of
course, we can determine all $N$-linear metrical connections of the space
$H^{(k)n}$. Therefore we consider the Obata's operator, [115]:
\begin{equation}
\Omega _{hm}^{ij}=\displaystyle\frac 12\left( \delta _h^i\delta
_m^j-g_{hm}g^{ij}\right) ,\ \Omega _{hm}^{*ij}=\displaystyle\frac 12\left(
\delta _h^i\delta _m^j+g_{hm}g^{ij}\right)  \tag{8.7.9}
\end{equation}
and prove without difficulties:

\begin{teo}
The set of all $N$-linear metrical connections \\
$D\overline{\Gamma }(N)=(\overline{H}_{jh}^i,\underset{(\alpha )}{%
\overline{C}}$ $_{jh}^i,\overline{C}_i^{jh})$ of the space $H^{(k)n}$ are
expressed by
\begin{equation}
\begin{array}{l}
\overline{H}_{jh}^i=H_{jh}^i+\Omega _{rj}^{is}X_{sh}^r, \\
\\
\underset{(\alpha )}{\overline{C}}\text{ }_{jh}^i=\underset{(\alpha )}{C}%
\text{ }_{jh}^i+\Omega _{rj}^{is}\underset{(\alpha )}{Y}\text{
}_{sh}^r,\
(\alpha =1,...,k-1), \\
\\
\overline{C}_i^{jh}=C_i^{jh}+\Omega _{rj}^{js}Z_s^{rh}.
\end{array}
\tag{8.7.10}
\end{equation}
where $C\Gamma (N)=(H_{jh}^i,\underset{(\alpha )}{C}$
$_{jh}^i,C_i^{jh})$
is the canonical $N$-metrical connection and $X_{sh}^r$, $Y$ $_{sh}^r$ and $%
Z_s^{rh}$ are arbitrary $d$-tensors.
\end{teo}

It is important to remark that a triple $\left(
X_{jh}^i,\underset{(\alpha )}{Y}\text{ }_{jh}^i,Z_i^{jh}\right) $
determines by means of (8.7.10) a transformation of the $N$-metrical
connections $D\Gamma (N)\rightarrow D\overline{\Gamma } (N) $.

The set of transformations $\left\{ D\Gamma (N)\rightarrow D\overline{\Gamma
}(N)\right\} $ and the composition of these mappings is an Abelian group.

The last theorems lead to the following result:
\begin{teo}
There exists an unique $N$-linear connection $D\overline{\Gamma }(N)=$ \newline
$=(%
\overline{H}_{jh}^i,\underset{(\alpha )}{\overline{C}}_{jh}^i,\overline{%
C}_i^{jh})$, $(\alpha =1,...,k-1)$, metrical with respect to the fundamental
tensor $g^{ij}$ of the space $H^{(k)n}$ having as torsion the $d$
-tensor fields $T_{\ jh}^i$ $(=-T_{\ hj}^i)$, $\underset{(\alpha )}{S}$ $%
_{jh}^i$ $(=-\underset{(\alpha )}{S}$ $_{hj}^i)$, $S_i^{jh}$
$(=-S_i^{hj})$, apriori given. The coefficients of $D\overline{\Gamma }(N)$ have
the following expressions:
\begin{equation}
\begin{array}{l}
\overline{H}_{jh}^i=\displaystyle\frac 12g^{is}\left( \displaystyle\frac{%
\delta g_{sh}}{\delta x^j}+\displaystyle\frac{\delta g_{js}}{\delta x^h}-%
\displaystyle\frac{\delta g_{jh}}{\delta x^s}\right) +\displaystyle\frac
12g^{is}\left( g_{sm}T_{\ jh}^m-g_{jm}T_{\ sh}^m+g_{hm}T_{\ sj}^m\right) , \\
\\
\underset{(\alpha )}{\overline{C}}\text{}_{jh}^i=\displaystyle\frac
12g^{is}\left( \displaystyle\frac{\delta g_{sh}}{\delta y^{(\alpha )j}}+%
\displaystyle\frac{\delta g_{js}}{\delta y^{(\alpha )h}}-\displaystyle\frac{%
\delta g_{jh}}{\delta y^{(\alpha )s}}\right) +\displaystyle\frac
12g^{is}( g_{sm}\underset{(\alpha )}{S}\text{}_{\ jh}^m-g_{jm}%
\underset{(\alpha )}{S}\text{}_{\ sh}^m+\\
\\
+g_{hm}\underset{(\alpha )}{S}%
\text{}_{\ sj}^m) , \\
\\
\overline{C}_i^{jh}=\displaystyle\frac 12g_{is}\left( \stackrel{\cdot}{\partial^j}g^{sh}+\stackrel{\cdot}{\partial^h}g^{js}-\stackrel{\cdot}{\partial^s}g^{jh}\right) -\displaystyle\frac 12g_{is}\left(g^{sm}S_m^{\ jh}-g^{jm}S_m^{\ sh}+g^{hm}S_m^{\ js}\right) .
\end{array}
\tag{8.7.11}
\end{equation}
\end{teo}

Now, consider the canonical $N$-metrical connection $C\Gamma (N)$
only. Then, we can write the Ricci identities, from Theorem 6.6.1,
where $T_{\ jh}^i=0$, $\underset{(\alpha )}{S}\text{}_{jh}^i=0$,
($\alpha =1,...,k-1$), $S_i^{jh}=0$. Since the metricity conditions (6.6.7) are verified it follows:

\begin{prop}
The curvature $d$-tensor fields of $C\Gamma (N)$ satisfy the identities
(6.6.8).
\end{prop}

\begin{prop}
The tensors of deflection of the canonical $N$-metrical connection $C\Gamma
(N)$ satisfy the identities (6.6.9), with $T_{\ jh}^i=0$, $\underset{(\alpha )}{S}\text{}_{jh}^i=0$, ($\alpha =1,...,k-1$),
$S_i^{\ jh}=0$.
\end{prop}

Let $\omega _{\ j}^i$ be the $1$-forms connection of $C\Gamma (N)$

\begin{equation}
\omega _{\ j}^i=H_{js}^idx^s+\underset{(1)}{C}\text{
}_{js}^i\delta y^{(1)s}+\cdots +\underset{(k-1)}{C}\text{
}_{js}^i\delta y^{(k-1)s}+C_i^{js}\delta p_s.  \tag{8.7.12}
\end{equation}

Theorem 6.7.1 is valid for $C\Gamma (N)$.

In particular, Theorem 6.7.2 for $C\Gamma (N)$ holds, too. Namely:

\begin{teo}
The manifold $T^{*k}M$ is with absolute parallelism of vectors with respect
to the canonical $N$-metrical connection $C\Gamma (N)$ if and only if the
curvature $d$-tensors of $C\Gamma (N)$ vanishes.
\end{teo}

Also, we have:

\begin{teo}
A smooth parametrized curve $\gamma $, given by $x^i=x^i(t)$, $y^{(\alpha
)i}=y^{(\alpha )i}(t)$, $(\alpha =1,...,k-1)$, $p_i=p_i(t)$, $t\in I$ is an
autoparallel curve with respect to $C\Gamma (N)$ if and only if the
equations (6.7.13), are verified, $\omega _{\ j}^i$ being the $1$-forms
connection (8.7.12).

The horizontal curves of the space $H^{(k)n}$ are characterized by the
equations (6.7.14).

The horizontal paths of $H^{(k)n}$ are expressed by the equations (6.7.15) for the connection $C\Gamma (N)$.

The $v_\alpha $-vertical paths of $H^{(k)n}$ are characterized by (6.7.16) and $w_k$-vertical paths are given by (6.7.17).
\end{teo}

The canonical $N$-metrical connection $C\Gamma (N)$ of the Hamilton space of
order $k$, $H^{(k)n}=(M,H)$ is $h$-, $v_\alpha $- and $w_k$-torsion
free if the $d$-tensors of torsion vanish:

\begin{equation}
\begin{array}{l}
T_{jh}^i=H_{jh}^i-H_{hj}^i=0,\ \underset{(\alpha)}{S}\text{}_{jh}^i=
\underset{(\alpha )}{C}\text{}_{jh}^i-\underset{(\alpha)}{C}\text{}_{hj}^i=0,
\ (\alpha =1,...,k-1),\vspace{3mm}\\
 S_i^{jh}=C_i^{jh}-C_i^{hj}.
\end{array}
\tag{8.7.13}
\end{equation}

Therefore, we have:

\begin{teo}
The canonical $N$-metrical connection $C\Gamma (N)$, with the coefficients
(8.7.8), of the Hamilton space of order $k$, $H^{(k)n}$, has the structure
equations given by the equations (7.8.1), (7.8.2), (7.8.3) and (7.8.4) from Theorem 7.8.1, in the conditions (8.7.13) and with $1$-forms $\omega _{\ j}^i$ from (8.7.12).
\end{teo}

The Bianchi identities of $C\Gamma (N)$ can be obtained from (7.6.11), (7.6.12),
(7.6.13) by means of the conditions (8.7.13).

\section{The Hamilton Space $H^{(k)n}$ of Electrodynamics}

The classical Lagrangian of electrodynamics has been considered in section 5 of
chapter 2:
\[
\stackrel{\cdot }{L}(x,y^{(1)})=mc\gamma _{ij}(x)y^{(1)i}y^{(1)j}+ %
\displaystyle\frac{2e}mb_i(x)y^{(1)i}
\]
in which $m$, $c$, $e$ are known physical constants, $\gamma _{ij}(x)$ are
gravitational potentials and $b_i(x)$ are electromagnetic potentials. This can be
extended to the manifold of accelerations of order $k$, $T^kM$ as follows:
\begin{equation}
L(x,y^{(1)},...,y^{(k)})=mc\gamma _{ij}(x)z^{(k)i}z^{(k)j}+\displaystyle
\frac{2e}mb_i(x)z^{(k)i},  \tag{8.8.1}
\end{equation}
$z^{(k)i}$ being the Liouville $d$-vector field
\begin{equation}
kz^{(k)i}=ky^{(k)i}+(k-1)\underset{(1)}{M}\text{}_j^iy^{(k-1)j}+\cdots + \underset{(k-1)}{M}\text{}_j^iy^{(1)j}.
\tag{8.8.2}
\end{equation}

The dual coefficients $\underset{(1)}{M}\text{}_j^i$, ...,
$\underset{(k-1)}{M}\text{}_j^i$ are given by the prolongation to
$T^kM$ of the Riemannian structure $\mathcal{R}^n=(M,\gamma_{ij}(x))$.

Consequently, $\underset{(1)}{M}$ $_j^i$ depends on the variables
$(x,y^{(1)})$; $\underset{(2)}{M}$ $_j^i$ depends on the variables
$(x,y^{(1)},y^{(2)})$; ...; $\underset{(k-1)}{M}$ $_j^i$ depends on
the variables $(x,y^{(1)},...,y^{(k-1)})$. This is an important property, since
$\left( \underset{(1)}{M}\text{}_j^i,...,\underset{(k-1)}{M}\text{} _j^i\right) $
can be considered as {\it the dual coefficients
of a nonlinear connection} $\stackrel{\circ }{N}$ on the manifold $ T^{*k}M,$
determined only by the tensor field $\gamma _{ij}(x)$.

It follows that:

1) The Lagrangian $L$, (8.8.1), is well determined only by $\gamma _{ij}(x)$ and $b_i(x)$;

2) $L$ is differentiable on $T^kM$;

3) $L$ is regular:
\begin{equation}
g_{ij}=\displaystyle\frac 12\displaystyle\frac{\partial ^2L}{\partial
y^{(k)i}\partial y^{(k)j}}=mc\gamma _{ij}(x);  \tag{8.8.3}
\end{equation}

4) The pair $L^{(k)n}=(M,L)$ is a Lagrange space of order $k$.

It is called \textit{the Lagrange space of order }$k$\textit{\ of
electrodynamics}. The geometry of the space $L^{(k)n}$ gives us a good
geometrical model for the Analytical Mechanics based on the Lagrangian of
electrodynamics. So, the integral of action $I(c)=\int\limits_0^1L(x, %
\displaystyle\frac{dx}{dt},\cdots ,\displaystyle\frac
1{k!}\displaystyle \frac{d^kx}{dt^k})dt$ leads to the
Euler-Lagrange equations $\stackrel{\circ }{ E_i}(L)=0$ (2.2.1).
The energy of higher order can be determined. The energy of order
$k$, $\mathcal{E}_c^k(L),$ satisfies the law of conservation. A
N\"other theorem, holds.

In the same time the canonical $k$-semispray $S$ is given by
\begin{equation}
S=y^{(1)}\displaystyle\frac \partial {\partial x^i}+2y^{(2)i}\displaystyle %
\frac \partial {\partial y^{(1)i}}+\cdots +ky^{(k)i}\displaystyle\frac
\partial {\partial y^{(k-1)i}}-(k+1)G^i\displaystyle\frac \partial {\partial
y^{(k)i}}  \tag{8.8.4}
\end{equation}
whose coefficients $G^i$, calculated by means of (8.8.1) have the expressions
\begin{equation}
(k+1)G^i=\displaystyle\frac 1{2mc}\gamma ^{ij}\left\{ \Gamma \left( %
\displaystyle\frac \partial {\partial y^{(k)j}}\right) -\displaystyle\frac
\partial {\partial y^{(k-1)j}}\right\}  \tag{8.8.4a}
\end{equation}
with
\begin{equation}
\Gamma =y^{(1)}\displaystyle\frac \partial {\partial x^i}+2y^{(2)i} %
\displaystyle\frac \partial {\partial y^{(1)i}}+\cdots +ky^{(k)i} %
\displaystyle\frac \partial {\partial y^{(k-1)i}}  \tag{8.8.4b}
\end{equation}
and
\begin{equation}
\displaystyle\frac{\partial L}{\partial y^{(k)j}}=2mc\gamma
_{ih}(x)z^{(k)h}+ \displaystyle\frac{2e}mb_j(x).  \tag{8.8.4c}
\end{equation}

The Legendre transformation $\xi :T^kM\rightarrow T^{*k}M$ determined by the
Lagrange space of electrodynamics $L^{(k)n}$, is defined by (8.4.2), (8.4.2').
This is:
\begin{equation}
p_i=\displaystyle\frac 12\displaystyle\frac{\partial L}{\partial y^{(k)j}}
=mc\gamma _{ih}(x)z^{(k)h}+\displaystyle\frac emb_i(x).  \tag{8.8.5}
\end{equation}

The mapping $\varphi $ is a local diffeomorphism. Its inverse $\varphi
^{-1}=\xi $ is given by (8.4.3), where

\begin{equation}
\begin{array}{lll}
\xi ^i(x,y^{(1)},...,y^{(k-1)},p) & = & \displaystyle\frac 1{mc}\gamma
^{ih}(x)\left\{ p_h-\displaystyle\frac emb_h(x)\right\} - \\
& & \\
& - & \displaystyle\frac 1k\left\{ (k-1)\underset{(1)}{M}\text{}
_h^iy^{(k-1)h}+\cdots +\underset{(k-1)}{M}\text{}_h^iy^{(1)h}\right\} .
\end{array}
\tag{8.8.5a}
\end{equation}

The canonical $k$-semispray $S$, (8.8.4) is transformed by Legendre mapping $%
\varphi $ in the dual $k$-semispray $S_\xi ^{*}$, from (8.4.8), with the
coefficients $\eta_i$ from (8.4.9).

Theorem 8.4.2 allows to determine, locally, the canonical nonlinear connection
$N^{*}$ on the manifold $T^{*k}M$.

As we remarked above the connection $N$ of the prolongation of the
Riemannian structure $(M,\gamma _{ij}(x))$, with the dual coefficients $%
\underset{(1)}{M}\text{}_j^i(x,y^{(1)})$, ..., \\
$\underset{(k-1)}{M}\text{}_j^i(x,y^{(1)},...,y^{(k-1)})$ defines a
nonlinear connection on the submanifold $T^{k-1}M$ in $T^kM$.

Thus, the Liouville vector field $z^{(k)i}$ from (8.8.2) is transformed by the
Legendre mapping $\varphi $ in the vector field $\stackrel{\vee}{z}^{(k)i}$:

\begin{equation}
k\stackrel{\vee}{z}^{(k)i}=k\xi ^i+(k-1)\underset{(1)}{M}\text{}
_m^iy^{(k-1)m}+\cdots +\underset{(k-1)}{M}\text{}_m^iy^{(1)m}
\tag{8.8.6}
\end{equation}

Let us consider the following Hamiltonian $H$, the dual of $L$:

\begin{equation}
H(x,y^{(1)},...,y^{(k-1)},p)=2p_i\stackrel{\vee}{z}
^{(k)i}-L(x,y^{(1)},...,y^{(k-1)},\xi (x,y^{(1)},...,y^{(k-1)},p)).
\tag{8.8.7}
\end{equation}

We can state:

\begin{teo}
We have:

1$^0$ The dual of the Lagrange space of electrodynamics $L^{(k)n}=(M,L)$, (8.8.1),
is a Hamilton space $H^{(k)n}=(M,H)$, (8.8.7).

2$^0$ The fundamental function $H$ of the space $H^{(k)n}$ is given by
\begin{equation}
H=\displaystyle\frac 1{mc}\left\{ \gamma ^{ij}(x)p_ip_j-\displaystyle\frac{2e%
}m\gamma ^{ij}p_ib_j+\displaystyle\frac e{m^2}\gamma ^{ij}b_ib_j\right\} .
\tag{8.8.8}
\end{equation}

3$^0$ $H$ depends only on the point $x$ and the momenta $p$.

4$^0$ $H$ is the fundamental function of a Hamilton space of order $1$, $%
H^n=(M,H)$.
\end{teo}

\textit{Proof}: 1$^0$, 2$^0$. Applying theorem 7.4.3 and the formulas (8.8.5),
(8.8.6) we obtain:
\[
\stackrel{\vee}{z}^{(k)i}=\displaystyle\frac 1{mc}\gamma ^{ij}(x)\left[
p_i- \displaystyle\frac emb_i(x)\right]
\]
and
$$
L(x,y^{(1)},...,y^{(k-1)},\xi )=mc\gamma _{ij}\displaystyle\frac 1{mc}\gamma
^{ij}\left[ p_i-\displaystyle\frac emb_i\right] \left[ p_j-\displaystyle
\frac emb_j\right]+$$
$$+\displaystyle\frac{2e}{m^2c}\gamma ^{ij}b_j\left[ p_i- %
\displaystyle\frac emb_i\right] .
$$
Substituting in (8.8.7) we have the equality (8.8.8).

3$^0$ Looking at the function $H$ one remarks that $H$ does not depend on
the variables $y^{(1)i}$, ..., $y^{(k-1)i}$.

4$^0$ $H$ is a Hamilton function on the cotangent bundle $T^{*}M$ and its
Hessian, with respect to the momenta $p_i$, has the matrix $\left\| %
\displaystyle\frac 1{mc}\gamma ^{ij}(x)\right\| $. So, the pair $H^n=(M,H)$
is a Hamilton space. Q.E.D.

\textbf{Remarks}. 1$^0$ The Hamilton space $H^n=(M,H(x,p))$ is the dual of
the Lagrange space of electrodynamics $L^n=(M,\stackrel{\circ }{L}
(x,y^{(1)}))$.

2$^0$ $H^{(k)n}$ does not depend on the nonlinear connection $%
N $ of the prolongation of the Riemann space $\mathcal{R}^n=(M,\gamma _{ij}(x))$.

For the Hamilton space $H^n=(M,H)$ with the fundamental function $H(x,p)$
from (8.8.8) we have a canonical nonlinear connection $N$, with the coefficients
\begin{equation}
N_{ij}=\gamma _{ij}^h(x)p_h+\displaystyle\frac ec\left(
b_{i|j}+b_{j|i}\right) ,  \tag{8.8.9}
\end{equation}
where $b_{i|j}=\displaystyle\frac{\partial b_i}{\partial x^j}-b_s\gamma
_{ij}^s$, $\gamma _{jh}^i(x)$ being the Christoffel symbols of the metric $%
\displaystyle\frac 1{mc}\gamma _{ij}(x)$.

The canonical $N$-metrical connection $C\Gamma (N)=(H_{jh}^i,C_i^{jh})$ has
the following coefficients
\begin{equation}
H_{jh}^i=\gamma _{jh}^i,\quad C_i^{jh}=0.  \tag{8.8.10}
\end{equation}

The geometrical object fields $H$, $g_{ij}=\displaystyle\frac 1{mc}\gamma
_{ij}$, $N_{ij}$ and $H_{jh}^i$ allow to develop the geometry of
the space $H^n=(M,H)$.

\section{The Riemannian Almost Contact Structure Determined by the Hamilton
Space $H^{(k)n}$}

Consider a Hamilton space of order $k$, $H^{(k)n}=(M,H)$ having $g^{ij}$ as
fundamental tensor field. Let $N$ be a nonlinear connection canonical
associated to $H^{(k)n},$ according to Theorem 8.6.1. As usually, we
consider the adapted basis $\left( \displaystyle\frac \delta {\delta x^i}, %
\displaystyle\frac \delta {\delta y^{(1)i}},\cdots ,\displaystyle\frac
\delta {\delta y^{(k-1)i}},\stackrel{\cdot }{\partial }^i\right) $ and its
dual basis $\left( \delta x^i,\delta y^{(1)i},...,\delta y^{(k-1)i},\delta
p_i\right) $. Then, at every point $u^{*}=(x,y^{(1)},...,y^{(k-1)},p)\in
T^{*k}M$ we can define the tensor
\begin{equation}
\stackrel{\vee}{\Bbb{G}}=g_{ij}dx^i\otimes dx^j+g_{ij}\delta
y^{(1)i}\otimes \delta y^{(1)j}+\cdots +g_{ij}\delta y^{(k-1)i}\otimes
\delta y^{(k-1)j}+g^{ij}\delta p_i\otimes \delta p_j  \tag{8.9.1}
\end{equation}

$\stackrel{\vee}{\Bbb{G}}$ is the $N$-lift of the fundamental tensor $%
g^{ij} $ of the Hamilton space of order $k$, $H^{(k)n}=(M,H)$, (cf. \S 6.7, ch. 6).

\begin{teo}
1$^0$ $\stackrel{\vee}{\Bbb{G}}$ is a (pseudo)-Riemannian structure on the
manifold $\widetilde{T^{*k}M}$, determined only by the space $H^{(k)n}$ and the nonlinear connection $N$.

2$^0$ The distributions $N_0$, $N_1$, ..., $N_{k-2}$, $V_{k-1}$, $W_k$ are mutual
orthogonal with respect to $\stackrel{\vee}{\Bbb{G}}$.
\end{teo}

Indeed: 1$^0$. Every term of $\stackrel{\vee}{\Bbb{G}}$ is defined on $%
\widetilde{T^{*k}M}$ because $g_{ij}$ is a $d$-tensor field and $\delta x^i$
, ..., $\delta p_i$ have geometrical meaning, i.e. they are transformed as
in (6.3.4), $\det ||\stackrel{\vee}{\Bbb{G}}||\neq 0$. 2$^0$. This
property is obvious, taking into account the form of $\stackrel{\vee}{\Bbb{%
G }}$. Q.E.D.

The tensor $\stackrel{\vee}{\Bbb{G}}$ is of the form:

\begin{equation}
\begin{array}{l}
\stackrel{\vee}{\Bbb{G}}\Bbb{=G}^H+\Bbb{G}^{V_1}+\cdots +\Bbb{G}^{V_{k-1}}+
\Bbb{G}^{W_k}, \vspace{3mm}\\
\Bbb{G}^H=g_{ij}dx^i\otimes dx^j,\Bbb{G}^{V_\alpha }=g_{ij}\delta y^{(\alpha
)i}\otimes \delta y^{(\alpha )j},(\alpha =1,...,k-1),\vspace{3mm}\\
\Bbb{G}^{W_k}=g^{ij}\delta p_i\otimes \delta p_j.
\end{array}
\tag{8.9.2}
\end{equation}

Here $\Bbb{G}^H$ is the restriction of $\stackrel{\vee}{\Bbb{G}}$ to the
distribution $N_0=N^{*}$, $\Bbb{G}^{V_\alpha }$ is the restriction of $%
\stackrel{\vee}{\Bbb{G}}$ to the distribution $N_\alpha $ , $\alpha
=1,...,k-1$ ($N_{k-1}=V_{k-1}$) and $\Bbb{G}^{W_k}$ is the restriction of $%
\stackrel{\vee}{\Bbb{G}}$ to the distribution $W_k$. Moreover, $\Bbb{G}^H$,
$\Bbb{G}^{V_\alpha }$, $\Bbb{G}^{W_k}$ are $d$-tensor fields.

\begin{teo}
The tensor fields $\Bbb{G}^H$, $\Bbb{G}^{V_\alpha }$, ($\alpha =1,...,k-1$),
$\Bbb{G}^{W_k}$ are covariant constant with respect to any metrical $N$
-connection $D\Gamma (N)$.
\end{teo}

Indeed, using the considerations of this section, it follows $D_X\Bbb{G}
^H=0$, $D_X\Bbb{G}^{V_\alpha }=0$, $D_X\Bbb{G}^{W_k}=0$ and, consequently $%
D_X\stackrel{\vee}{\Bbb{G}}=0$. Q.E.D.

The geometry of the (pseudo)-Riemannian space $\left( \widetilde{T^{*k}M},
\stackrel{\vee}{\Bbb{G}}\right) $ can be studied by means of a metrical
$N$-linear connection, presented in the section 7 of this chapter.

As we know (\S 6, Ch.6) on $\widetilde{T^{*k}M}$, endowed with the
nonlinear connection $N$, there exists a natural almost $(k-1)n$-contact
structure $\Bbb{F}$ , defined by (6.6.3), having the properties
from in Theorem 6.5.2. The condition of normality of $\Bbb{F}$ is
(6.6.5). Exactly as in Theorem 6.7.2 we can deduce:

\begin{teo}
The pair $\left( \stackrel{\vee}{\Bbb{G}},\Bbb{F}\right) $ is a Riemannian $%
(k-1)n$-almost contact structure determined by the space $H^{(k)n}$ and
the nonlinear connection $N$.
\end{teo}

The considerations from section 6.8, ch. 8 show that the nonlinear
connection $N$ and the fundamental tensor $g^{ij}$ of the Hamilton space of
order $k$, $H^{(k)n},$ determine on the manifold $\widetilde{T^{*k}M}$, two
important geometric object fields. One is the $N$-lift $\stackrel{\vee}{%
\Bbb{G}}$ from (8.9.1) and another one is $\stackrel{\vee}{\Bbb{F}}$ the almost
$(k-1)n$-contact structure, given by (6.8.5).

Namely,
\[
\stackrel{\vee}{\Bbb{F}}\left( \displaystyle\frac \delta {\delta
x^i}\right) =-g_{ij}\stackrel{\cdot }{\partial }^j,\ \stackrel{\vee}{\Bbb{F}
}\left( \displaystyle\frac \delta {\delta y^{(\alpha )i}}\right) =0,\
(\alpha =1,...,k-1),\ \stackrel{\vee}{\Bbb{F}}\left( \stackrel{\cdot }{
\partial }^i\right) =g^{ij}\displaystyle\frac \delta {\delta x^j}
\]

By means of Theorem 7.9.1 it follows that $\stackrel{\vee}{\Bbb{F}}$ is a
tensor field of type $(1,1)$ on the manifold $T^{*k}M$ given in (6.8.3), $rank$ $\stackrel{\vee}{\Bbb{F}}=2n$ and $\stackrel{\vee}{\Bbb{F}}^3+
\stackrel{\vee}{\Bbb{F}}=0$.

The following version of the Theorem 6.8.2 holds:

\begin{teo}
For any Hamilton space of order $k$, $H^{(k)n}=(M,H),$ endowed with a
canonical nonlinear connection $N$ the following properties hold:

1$^0$ The pair $\left( \stackrel{\vee}{\Bbb{G}},\stackrel{\vee}{\Bbb{F}}%
\right) $ is a Riemannian almost $(k-1)n$-contact structure.

2$^0$ The associated $2$-form is
\[
\theta =\delta p_i\wedge dx^i.
\]

3$^0$ If the coefficients $N_{ij}$ of $N$ are symmetric then
\[
\theta =dp_i\wedge dx^i
\]
and $\theta $ is the canonical presymplectic structure on $\widetilde{T^{*k}M}$.

4$^0$ The condition of normality of structure $\stackrel{\vee}{\Bbb{F}}$ is
\[
\mathcal{N}_{\stackrel{\vee
}{\Bbb{F}}}(X,Y)+\sum\limits_{i=1}^n\left[ \sum\limits_{\alpha
=1}^{k-1}d(\delta y^{(\alpha )i})(X,Y)+d\delta p_i(X,Y)\right]
=0,\ \forall X,Y\in \mathcal{X}(\widetilde{T^{*k}M}),
\]
where $\mathcal{N}_{\stackrel{\vee}{\Bbb{F}}}$ is the Nijenhuis tensor of $%
\stackrel{\vee}{\Bbb{F}}$.
\end{teo}

Now, taking the local expression of the tensor field $\stackrel{\vee}{\Bbb{F}}$:
\[
\stackrel{\vee}{\Bbb{F}}=-g_{ij}\stackrel{\cdot }{\partial }^j\otimes
dx^i+g^{ij}\displaystyle\frac \delta {\delta x^i}\otimes \delta p_j
\]
and applying the theory of $h$-, $v_\alpha $- and $w_k$-covariant derivative
from \S 7.5, ch. 7, it follows:

\begin{teo}
With respect to a metrical $N$-linear connection $D\Gamma (N)$, the
Riemannian almost $(k-1)n$-contact structure $\left( \stackrel{\vee}{\Bbb{G}%
},\stackrel{\vee}{\Bbb{F}}\right) $ is covariant constant, i.e.
\[
D_X\stackrel{\vee}{\Bbb{G}}=0,\ D_X\stackrel{\vee}{\Bbb{F}}=0.
\]
\end{teo}

It follows that the geometry of the Riemannian almost $(k-1)n$-contact space
$\left( \widetilde{T^{*k}M},\stackrel{\vee}{\Bbb{G}},\stackrel{\vee}{\Bbb{%
F }}\right) $ can be studied by means of the metrical $N$-linear connection $%
D\Gamma (N)$. This space is called \textit{the Riemannian almost }$(k-1)n$
\textit{-contact model }of the Hamilton space of order $k$, $H^{(k)n}$.

\chapter{Subspaces in Hamilton Spaces of Order $k$}

\markboth{\it{THE GEOMETRY OF HIGHER-ORDER HAMILTON SPACES\ \ \ \ \ }}{\it{Subspaces in Hamilton Spaces of Order} $k$ }

In this chapter we shall study the geometry of subspaces in a Hamilton space
$H^{(k)n}=(M,H)$. A submanifold $\stackrel{\vee}{M}$ of the manifold $M$
determines a dual manifold $T^{*k}\!\!\stackrel{\vee}{M}$. But the immersion $%
i: \!\!\stackrel{\vee}{M}\rightarrow M$ does not automatically imply an
immersion of $T^{*k}\!\!\stackrel{\vee}{M}$, into the dual manifold $T^{*k}M$,
like in the case of manifold $T^{k}\!\!\stackrel{\vee}{M}$ into the total
space of $k$ tangent bundle $T^{k}M$, [94], since $T^{*k}\!\!\stackrel{\vee}{M}
=T^{k-1}\!\!\stackrel{\vee}{M}\times _{\stackrel{\vee}{M}}T^{*}\!\!\stackrel{\vee}{M}.$

Therefore, by means of the immersion $i$ and by an immersion of the cotangent
manifold $T^{*}M\!\!\stackrel{\vee}{M}$ into $T^{*}M$ we can define $T^{*k}\!\!\stackrel{\vee}{M}$
as an immersed submanifold in $T^{*k}M.$ Thus, the main geometrical objects
fields on $T^{*k}M$ induce the coresponding geometrical object fields on
submanifold $T^{*k}\!\!\stackrel{\vee}{M}.$ The Hamilton space $H^{\left( k\right) n}=(M,H)$
induces a Hamilton subspaces $\stackrel{\vee}{H}^{\left( k\right) m}=(
\stackrel{\vee}{M},\stackrel{\vee}{H}).$ So we study the intrinsic
geometrical object fields on $\stackrel{\vee}{H}^{\left( k\right)m}$ and
the induced geometrical object fields, as well as the relations between
them. These problems are approached by means of the method of moving frame, used
in the case of subspaces $\stackrel{\vee}{L}^{(k)m}=(\stackrel{\vee}{M},
\stackrel{\vee}{L})$ in the Lagrange spaces of order $k.$

\section{ Submanifolds $T^{*k}\!\!\stackrel{\vee}{M}$ in the Manifold $T^{*k}M$}

Let $M$ be a $C^\infty -$ real, $n$-dimensional manifold and $\stackrel{\vee
}{M}$ be a $C^\infty $-real, $m$ -dimensional manifold, $1<m<n,$ immersed
in $M$ through the immersion $i:\stackrel{\vee}{M}\rightarrow M.$ Locally $i$ can be given in the form
\begin{equation}
x^i=x^i(u^1,....,u^m),\ rank||\displaystyle\frac{\partial x^i}{\partial
u^\alpha }\parallel =m,\forall (u^\alpha )\in \stackrel{\vee}{U}\subset
\!\!\stackrel{\vee}{M}  \tag{9.1.1}
\end{equation}

The indices $i,j,h,r,s,p,q$ run over the set $\{1,2,...,n\};$ the indices $\alpha
,\beta ,\gamma ,...$ run over the set $\{1,2...m\}$ and $\overline{\alpha },
\overline{\beta },\overline{\gamma },...$run over the set $\{1,..n-m\}$.

If $i:\stackrel{\vee}{M}\rightarrow M$ is an embeding, then we identify $\stackrel{\vee}{M}$ with $i(\stackrel{\vee}{M})\subset M$ and say that $\stackrel{\vee}{M}$ is a submanifold of $M$. In this case (9.1.1) are called the parametric equations of the submanifold $\stackrel{\vee}{M}$ in the manifold $M$.

Let us consider the manifold $T^{*k}\!\!\!\!\stackrel{\vee}{M}$ determined by $\stackrel{\vee}{M}$, i.e. \\
$T^{*k}\stackrel{\vee}{M\text{}}=T^{k-1}\!\!\stackrel{\vee}{M} \times _{\stackrel{\vee}{M}} T^{*}\!\!\stackrel{\vee}{M}$ . A point $\stackrel{\vee}{u} \in T^{*k}\stackrel{\vee}{M\text{}}$ will be denoted by $\stackrel{\vee}{u }=(u,v^{(1)},...v^{(k-1)},\stackrel{\vee}{p})$ and its coordinates by $(u^\alpha ,v^{(1)\alpha },...,v^{(k-1)\alpha },\stackrel{\vee}{p}_{\alpha
}).$ A change of local coordinates on the manifold $T^k\!\!\stackrel{\vee}{M}$ is given by
\begin{equation}
\left\{
\begin{array}{l}
\overline{u}^\alpha =\overline{u}^\alpha (u^1,...,u^m),rank\left\| %
\displaystyle\frac{\partial \overline{u}^\alpha }{\partial u^\beta }\right\|
=m ,\\
\\
\overline{v}^{\left( 1\right) \alpha }=\displaystyle\frac{\partial \overline{
u}^\alpha }{\partial u^\beta }v^{\left( 1\right) \beta }, \\
....................................................., \\
(k-1)\overline{v}^{(k-1)\alpha }=\displaystyle\frac{\partial \overline{v}
^{(k-2)\alpha }}{\partial u^\beta }v^{(1)\beta }+...+(k-1)\displaystyle\frac{
\partial \overline{v}^{(k-2)\alpha }}{\partial v^{(k-2)\beta }}v^{(k-1)\beta
}, \\
\\
\stackrel{\vee}{\overline{p}}_{\alpha}=\displaystyle\frac{\partial u^\beta }{
\partial \overline{u}^\alpha }\stackrel{\vee}{p}_\beta,
\end{array}
\right.  \tag{9.1.2}
\end{equation}
with
\[
\displaystyle\frac{\partial \overline{u}^\alpha }{\partial u^\beta }= %
\displaystyle\frac{\partial \overline{v}^{(1)\alpha }}{\partial v^{\left(
1\right) \beta }}=...=\displaystyle\frac{\partial \overline{v}^{(k-1)\alpha
} }{\partial v^{\left( k-1\right) \beta }}.
\]

Remarking that $T^{*k}:Man\rightarrow Man$ is a covariant functor from the
category of differentiable manifolds Man to itself, it follows that the
immersion $i:\stackrel{\vee}{M}\rightarrow M$ uniquelly determines the
mapping $T^{*k}i:T^{*k}\!\!\stackrel{\vee}{M}\rightarrow T^{*k}M$ (see Ch.4)
analytically given by the equations
\begin{equation}
\left\{
\begin{array}{l}
x^i=x^i(u^1,...,u^m),rank\left\| \displaystyle\frac{\partial x^i}{\partial
u^\alpha }\right\| =m,  \\
y^{\left( 1\right) i}=\displaystyle\frac{\partial x^i}{\partial u^\alpha }
v^{\left( 1\right) \alpha },  \\
........................................  \\
\left( k-1\right) y^{\left( k-1\right) i}=\displaystyle\frac{\partial
y^{\left( k-2\right) i}}{\partial u^\alpha }v^{\left( 1\right) \alpha }+2 %
\displaystyle\frac{\partial y^{\left( k-2\right) i}}{\partial v^{\left(
1\right) \alpha }}v^{\left( 2\right) \alpha }+...+  \\
\hfill  +\left( k-1\right) %
\displaystyle\frac{\partial y^{\left( k-2\right) i}}{\partial v^{\left(
k-2\right) \alpha }}v^{\left( k-1\right) \alpha },
\end{array}
\right.  \tag{9.1.3}
\end{equation}
and
\begin{equation}
\displaystyle\frac{\partial x^i}{\partial u^\alpha }p_i=\stackrel{\vee}{p}
_\alpha , \tag{9.1.3a}
\end{equation}
where
\begin{equation}
\begin{array}{l}
\displaystyle\frac{\partial x^i}{\partial u^\alpha }=\displaystyle\frac{
\partial y^{\left( 1\right) i}}{\partial v^{\left( 1\right) \alpha }}=...= %
\displaystyle\frac{\partial y^{\left( k-1\right) i}}{\partial v^{\left(
k-1\right) \alpha }}, \\
\\
\displaystyle\frac{\partial y^{(a)i}}{\partial u^\alpha }=\displaystyle\frac{
\partial y^{\left( a+1\right) i}}{\partial v^{\left( 1\right) \alpha }}=...= %
\displaystyle\frac{\partial y^{\left( k-1\right) i}}{\partial v^{\left(
k-1-a\right) \alpha }},(a=1,...,k-2).
\end{array}
\tag{9.1.3b}
\end{equation}

We shall denote
\begin{equation}
B_\alpha ^i(u)=\displaystyle\frac{\partial x^i}{\partial u^\alpha }
,B_{\alpha \beta }^i(u)=\displaystyle\frac{\partial ^2x^i}{\partial u^\alpha
\partial u^\beta },...  \tag{9.1.4}
\end{equation}

In order to obtain an immersion $i^{*}$ of $T^{*k}\!\!\stackrel{\vee}{M}=T^{k-1}
\!\!\stackrel{\vee}{M}\times _{\stackrel{\vee}{M}}T^{*}\!\!\stackrel{\vee}{M}$ in $%
T^{*k}M=T^{k-1}M\times _MT^{*}M,$ it is necessary that $i^{*}$ to be of the
form $i^{*}=T^{k-1}i\times f,$ with $f:T^{*}\!\!\stackrel{\vee}{M}\rightarrow
T^{*}M$. Analytically $T^{k-1}i$ is defined by the equations (9.1.3).
Consequently, taking into account (9.1.3'), it follows that the mapping $f$
can be given in the form
\begin{equation}
p_i=A_i^\alpha (u)\stackrel{\vee}{p}_\alpha ,rank\Vert A_i^\alpha (u)\Vert
=m,\forall u=(u^1,...,u^m)\in \stackrel{\vee}{U}\subset \!\!\stackrel{\vee}{M}
\tag{9.1.5}
\end{equation}

We obtain
\begin{teo}
The equations (9.1.3), (9.1.5) define a local immersion

$i^{*}: T^{*k}\!\!\stackrel{\vee}{M}\rightarrow T^{*k}M,$

$i^{*}: \stackrel{\vee}{u}=(u,v^{\left( 1\right) },...,v^{\left( k-1\right)
},\stackrel{\vee}{p}) \rightarrow u^{*}=(x,y^{\left( 1\right) },...,y^{(k-1)}, p).$
\end{teo}

\textbf{Proof}. The Jacobian matrix of the mapping $i^{*}$ given by the
equations (9.1.3), (9.1.5) at a point $\stackrel{\vee}{u}\in \stackrel{\vee}{U}$ is as follows:
\begin{equation}
J(i^{*})_{\stackrel{\vee}{u}}=\left\|
\begin{array}{cccccc}
\displaystyle\frac{\partial x^i}{\partial u^\alpha } & 0 & 0 & \cdots & 0 & 0
\vspace{3mm}\\
\displaystyle\frac{\partial y^{(1)j}}{\partial u^\alpha } & \displaystyle
\frac{\partial x^i}{\partial u^\alpha } & 0 & \cdots & 0 & 0 \\
\cdots & \cdots & \cdots & \cdots & \cdots & \cdots \\
\displaystyle\frac{\partial y^{(k-1)i}}{\partial u^\alpha } & \displaystyle
\frac{\partial y^{(k-1)i}}{\partial v^{(1)\alpha }} & \displaystyle\frac{
\partial y^{(k-1)i}}{\partial v^{(2)\alpha }} & \cdots & \displaystyle\frac{
\partial x^i}{\partial u^\alpha } & 0 \vspace{3mm}\\
\displaystyle\frac{\partial A_i^\beta }{\partial u^\alpha }\stackrel{\vee}{%
p }_\beta & 0 & 0 & \cdots & 0 & A_i^\alpha
\end{array}
\right\| _{\stackrel{\vee}{u}}  \tag{9.1.6}
\end{equation}

Since the matrices $\left\| \displaystyle\frac{\partial x^i}{\partial
u^\alpha }\right\| $ and $\left\| A_i^\alpha (u)\right\| $ have the rank $m$
it follows that $i^{*}$ defines a local immersion.    Q.E.D.

The functions $A_i^\alpha (u)$ from (9.1.5) are not arbitrary. They must satisfy some conditions.

\begin{prop}
The following properties hold:

1$^{\circ}$
\begin{equation}
B_\alpha ^i(u)A_i^\beta (u)=\delta _\alpha ^\beta,
\tag{9.1.7}
\end{equation}

2$^{\circ}$ $A_i^\alpha (u)$ is a d-covector field with respect to index $i$,

$A_i^\alpha (u)$ is a d-vector field with respect to index $\alpha$.
\end{prop}

Indeed, (9.1.7) follows from (9.1.3') and (9.1.5'), 2$^{\circ }$ $p_i$ being a $d$
-covector it follows that $A_i^\alpha (u)$ has the same quality with respect
to index $i.$ The same remark is valid for the index $\alpha $ of $%
A_i^\alpha (u).$

Of course $B_\alpha ^i(u)$ is a $d$-vector with respect to index $i$ and it is a
covector with respect to index $\alpha$.

The previous properties allow to call $B_\alpha ^i(u)$ and $A_i^\alpha (u)$
the \textit{mixed} $d$-tensor fields. Along $\stackrel{\vee}{M}$ a mixed $d$
-tensor will be given by the components $T_{\beta _1...\beta
_qj_1...j_s}^{\alpha _1...\alpha _pi_1...i_r}$ having the property that it a
$d$-tensor of type $(p,q)$ with respect \vspace{3mm}\\
to indices $\alpha _1...\alpha _p,$ and $\beta _1...\beta _q$ (to a change of coordinates on the submanifold
$T^{*k}\!\!\stackrel{\vee}{M}$ and it is a $d$-tensor of type $(r,s)$ with respect to indices $
i_1...i_r$ and $j_1...j_s,$ to a change of coordinates on the enveloping
manifold $T^{*k}M$.

\textbf{Remark.}The notion of mixed $d$-tensor will be extended in the section 3.

As usual we set
\[
\stackrel{\cdot }{\partial ^i}=\displaystyle\frac \partial {\partial p_i},
\stackrel{\cdot }{\partial ^\alpha }=\displaystyle\frac \partial {\partial
\stackrel{\vee}{p}_\alpha }.
\]

The differential $di^{*}$ of the immersion $i^{*}$ acts on the natural
basis \vspace{3mm}\\
$(\displaystyle\frac \partial {\partial u^\alpha },\displaystyle\frac
\partial {\partial v^{\left( 1\right) \alpha }},...,\displaystyle\frac
\partial {\partial v^{\left( k-1\right) \alpha }},\stackrel{\cdot }{\partial
^\alpha })$ by the rule:
\begin{equation}
di^{*}\left\| \displaystyle\frac \partial {\partial u^\alpha },\displaystyle %
\frac \partial {\partial v^{\left( 1\right) \alpha }},...,\displaystyle\frac
\partial {\partial v^{\left( k-1\right) \alpha }},\stackrel{\cdot }{\partial
^\alpha }\right\| =\left\| \displaystyle\frac \partial {\partial x^i},..., %
\displaystyle\frac \partial {\partial y^{\left( k-1\right) \alpha }},
\stackrel{\cdot }{\partial ^i}\right\| J(i^{*})  \tag{9.1.8}
\end{equation}
and for the natural cobasis:
\begin{equation}
di^{*}\left\| dx^i,...,dy^{(k-1)i},dp_i\right\| ^T=J(i^{*})\left\| dx^\alpha
,...,dv^{\left( k-1\right) \alpha },d\stackrel{\vee}{p}_\alpha \right\|
\tag{9.1.8a}
\end{equation}

Therefore, it is not difficult to see that:

1$^{\circ }$. We have

\begin{equation}
\stackrel{\cdot }{\partial ^\alpha }=A_i^\alpha \stackrel{\cdot }{\partial
^i };dp_i=dA_i^\alpha \stackrel{\vee}{p}_\alpha +A_i^\alpha d\stackrel{\vee
}{p }_\alpha  \tag{9.1.9}
\end{equation}

2$^{\circ }$ Along submanifold $T^{*k}\!\!\stackrel{\vee}{M}$ the vertical
distributions $v_1,...,v_{k-1},w_k$ are subdistributions of the vertical distributions $%
V_1,...,V_{k-1},W_k $ from the manifold $T^kM$.

3$^{\circ }$ We have for the Liouville vector fields $\stackrel{1}{\gamma }
,...,\stackrel{k-1}{\gamma }$ from $T^{*k}M$ the relations

\begin{equation}
di^{*}(\stackrel{1}{\gamma })=\stackrel{1}{\Gamma },...,di^{*}(\stackrel{k-1
}{\gamma })=\stackrel{k-1}{\Gamma }  \tag{9.1.10}
\end{equation}

These equations will be applied in the theory of Hamilton subspaces.

\section{Hamilton Subspaces $\stackrel{\vee}{H}^{(k)m}$ in $H^{(k)n}$.
Darboux Frames}

Let $H^{(k)n}=(M,H)$ be a Hamilton spaces of order $k,$ its fundamental
function $H$ being defined on the manifold $T^{*k}M$ and an immersion $%
i^{*}:T^{*k}\!\!\stackrel{\vee}{M}\rightarrow T^{*k}M$ given locally by the
equations (9.1.3), (9.1.5). For a point $u^{*}=(x,y^{(1)},...,y^{(k-1)},p)\in
T^{*k}M$ and a point $\stackrel{\vee}{u}=(u,v^{(1)},...,v^{(k-1)},\stackrel{
\vee}{p})\in T^{*k}\!\!\stackrel{\vee}{M}$ with the property $u^{*}=i^{*}(\stackrel{\vee}{u})$ we
obtain the restriction of the fundamental function $H$ expressed by
\begin{equation}
H\circ i^{*}(\stackrel{\vee}{u})=\stackrel{\vee}{H}(\stackrel{\vee}{u}
),\forall \stackrel{\vee}{u}\in T^{*k}\!\!\stackrel{\vee}{M}.
\tag{9.2.1}
\end{equation}

It follows that $\stackrel{\vee}{H}(\stackrel{\vee}{u})$ is a
differentiable Hamiltonian on the submanifold $T^{*k}\!\!\stackrel{\vee}{M}.$

Consider the $d$-tensor field
\begin{equation}
\stackrel{\vee}{g}^{\alpha \beta }=\displaystyle\frac 12\stackrel{\cdot }{
\partial ^\alpha }\stackrel{\cdot }{\partial ^\beta }\stackrel{\vee}{H}
on\quad \widetilde{T^{*k}\!\!\stackrel{\vee}{M}}.  \tag{9.2.2}
\end{equation}
>From (9.2.1) we deduce that at the points $\stackrel{\vee}{u}\in \widetilde{T^{*k}\!\!\stackrel{\vee}{M}}$ the fundamental tensor field $\stackrel{\vee}{g}^{\alpha \beta }$ is given by
\begin{equation}
\stackrel{\vee}{g}^{\alpha \beta }(\stackrel{\vee}{u})=A_i^\alpha
(u)A_j^\beta (u)g^{ij}(i^{*}\stackrel{\vee}{u}) . \tag{9.2.3}
\end{equation}

It follows from previous equality that

\begin{equation}
rang\left\| \stackrel{\vee}{g}^{\alpha \beta }\right\| =m.  \tag{9.2.3a}
\end{equation}

Therefore we have:
\begin{teo}
The pair $\stackrel{\vee}{H}^{(k)m}=(\!\!\stackrel{\vee}{M},\stackrel{\vee}{H}%
)$ is a Hamilton space of order $k.$
\end{teo}

Indeed, this property is a consequences of the equations (9.2.1), (9.2.2) and
(9.2.3) and the fact that $A_i^\alpha (u^1,...,u^m)$ is a mixed $d$-tensor,
with rank $\left\| A_i^\alpha \right\| =m.$

The space $\stackrel{\vee}{H}^{(k)m}$ will be called the Hamilton subspace of the
Hamilton space $H^{(k)n}=(M,H).$

In the point $\stackrel{\vee}{u}\in \stackrel{\vee}{U}\subset T^{*k}
\!\!\stackrel{\vee}{M}$, the $d$-covector fields $A_i^\alpha$, $(\alpha =1,...,m)$ are defined by $A_i^\alpha (\stackrel{
\vee}{u})=A_i^\alpha (u)$. So $A_i^1(\stackrel{\vee}{u})$,..., $A_i^m
(\stackrel{\vee}{u})$ are linearly independent. They are called tangent
covectors to $T^{*k}\!\!\stackrel{\vee}{M}.$ A $d$-covector $\omega
_i(u^{*})\,, $ $u^{*}\in T^{*k}M$ is called normal to $T^{*k}\!\!\stackrel{\vee}{M}
$ in $T^{*k}M $ if $u^{*}=i^{*}(\stackrel{\vee}{u})$ and
\[
g^{ij}(i^{*}(\stackrel{\vee}{u}))A_i^\alpha (\stackrel{\vee}{u})\omega
_j(i^{*}(\stackrel{\vee}{u}))=0,\forall \stackrel{\vee}{u}\in \stackrel{
\vee}{U}\subset T^{*k}\!\!\stackrel{\vee}{M}.
\]

Then we can determine a Darboux coframe $\mathcal{R}^{*}$, at every point $%
\stackrel{\vee}{u}\in \stackrel{\vee}{U},$ of the form
\begin{equation}
\mathcal{R}^{*}=\{\stackrel{\vee}{u},A_i^1(\stackrel{\vee}{u}),...,A_i^m(
\stackrel{\vee}{u});A_i^{\overline{1}}(\stackrel{\vee}{u}),...,A_i^{%
\overline{n-m}}(\stackrel{\vee}{u})\}  \tag{9.2.4}
\end{equation}
formed by $m$-tangent $d$-covectors $A_i^\alpha (\stackrel{\vee}{u})$ and by
$n-m$ normal unit $d$-covectors $A_i^{\overline{\alpha }}(\stackrel{\vee}{u}
),$ ($\overline{\alpha }=\overline{1},...,\overline{n-m})$, which verify the
conditions:
\begin{equation}
\begin{array}{l}
g^{ij}(i^{*}\stackrel{\vee}{u})A_i^\alpha (\stackrel{\vee}{u})A_j^{%
\overline{\alpha }}(\stackrel{\vee}{u})=0, \\
\\
g^{ij}(i^{*}\stackrel{\vee}{u})A_i^{\overline{\alpha }}(\stackrel{\vee}{u}
)A_j^{\overline{\beta }}(\stackrel{\vee}{u})=\delta ^{\overline{\alpha }
\overline{\beta }},\text{ }\forall \stackrel{\vee}{u}\in \stackrel{\vee}{U}
\end{array}.
\tag{9.2.4a}
\end{equation}

Of course $\mathcal{R}^{*}$ exists and it has a geometrical meaning with
respect to the change of local coordinates on the manifold $T^{*k}
\!\!\stackrel{\vee}{M}$ and $T^{*k}M$ and with respect to the orthogonal
(transformation of the normal $d$-covectors $A_i^{\overline{\alpha }})$
given by

\begin{equation}
A_i^{\overline{\alpha ^{\prime }}}(\stackrel{\vee}{u})=C_{\overline{\beta }
}^{\overline{\alpha }^{\prime }}(\stackrel{\vee}{u})A_i^{\overline{\beta }%
}( \stackrel{\vee}{u}),\quad \left\| C_{\overline{\beta }}^{\overline{%
\alpha } ^{\prime }}(\stackrel{\vee}{u})\right\| \in O(n-m) . \tag{9.2.5}
\end{equation}

$\mathcal{R}^{*}\,$ is called a \textit{moving coframe} on the
submanifold $T^{*k}\!\!\stackrel{\vee}{M}$. The dual frame $\mathcal{R}$ of the coframe $\mathcal{%
\ R}^{*}$ is given by
\begin{equation}
\mathcal{R}=\{\stackrel{\vee}{u},A_\alpha ^i(\stackrel{\vee}{u}),A_{%
\overline{\alpha }}^i(\stackrel{\vee}{u})\} , \tag{9.2.4b}
\end{equation}
where
\begin{equation}
\left\{
\begin{array}{l}
A_i^\alpha A_\beta ^i=\delta _\beta ^\alpha ,A_i^\alpha A_{\overline{\beta }
}^i=0,A_i^{\overline{\alpha }}A_\beta ^i=0,A_i^{\overline{\alpha }}A_{%
\overline{\beta }}^i=\delta _{\overline{\beta }}^{\overline{\alpha }}, \\
\\
A_i^\alpha A_\alpha ^j+A_i^{\overline{\alpha }}A_{\overline{\alpha }
}^j=\delta _i^j,
\end{array}
\right.  \tag{9.2.6}
\end{equation}
at every point $\stackrel{\vee}{u}\in T^{*k}\!\!\stackrel{\vee}{M}.$

So, $\mathcal{R}$ will be called the {\it moving frame or the Darboux frame} along the Hamilton
subspaces $\stackrel{\vee}{H}^{(k)m}$ in the Hamilton space $H^{(k)n}.$

The conditions (9.2.4'), (9.2,6) imply
\begin{equation}
\begin{array}{c}
\stackrel{\vee}{g}\stackrel{}{^{\alpha \beta }}A_\beta ^i=g^{ij}A_j^\alpha
\\
\delta ^{\overline{\alpha }\overline{\beta }}A_{\overline{\beta }
}^i=g^{ij}A_j^{\overline{\alpha }}
\end{array}
\tag{9.2.7}
\end{equation}

Indeed,

$\stackrel{\vee}{g}\stackrel{}{^{\alpha \beta }}A_\beta ^i=g^{rs}A_r^\alpha
A_s^\beta A_\beta ^i=g^{rs}A_r^\alpha (\delta _s^i-A_s^{\overline{\beta }}A_{%
\overline{\beta }}^i)=g^{ij}A_j^\alpha .$

In the following investigation we use the Darboux frame $\mathcal{R}$ in
order to represent in $\mathcal{R}$ the $d$-tensor from the Hamilton spaces $%
H^{(k)n}.$

We get:

\begin{prop}
The fundamental tensor $g^{ij}$ of the Hamilton space $H^{(k)n}$ and
its covariant $g_{ij}$ are represented in the Darboux frame $\mathcal{R}$ by
\begin{equation}
\left\{
\begin{array}{c}
g^{ij}=\stackrel{\vee}{g}\stackrel{}{^{\alpha \beta }}A_\alpha ^iA_\beta
^j+\delta ^{\overline{\alpha }\overline{\beta }}A_{\overline{\alpha }}^iA_{%
\overline{\beta }}^j \\
\\
g_{ij}=\stackrel{\vee}{g}_{\alpha \beta }A_i^\alpha A_j^\beta +\delta _{%
\overline{\alpha }\overline{\beta }}A_i^{\overline{\alpha }}A_j^{\overline{%
\beta }}
\end{array}
\right.   \tag{9.2.8}
\end{equation}
\end{prop}

Indeed, the formulas (9.2.8) are consequence of (9.2.7)

The moving frame $\mathcal{R}$ will be used in the next section in order to
derive the Gauss-Weingarden formulae and Gauss-Codazzi equations of the
Hamilton subspaces $\stackrel{\vee}{H}^{\left( k\right) m}$ in $H^{\left(k\right) n}.$

\section{Induced Nonlinear Connection}

Now, let us consider the canonical nonlinear connection $N$ of the
Hamilton space of order $k,$ $H^{(k)n}=(M,H)$. $N$ has the dual
coefficients $(
\underset{(1)}{M_j^i},...,\underset{(k-1)}{M_j^i},N_{ij}).$ We
shall prove that the restriction of $N$ to the Hamilton subspaces
$\stackrel{\vee}{H}^{(k)m}=(\!\!\stackrel{\vee}{M},\!\!\stackrel{\vee}{H})$ determines an
induced nonlinear connection $\stackrel{\vee}{N}$ on the manifold
$T^{*k}\!\!\stackrel{\vee}{M}.$

The nonlinear connection $\stackrel{\vee}{N}$ is well determined by its
dual coefficients or by means of its adapted cobasis $(du^\alpha ,\delta
v^{(1)\alpha },....,\delta v^{(k-1)\alpha },\delta \stackrel{\vee}{p}
_\alpha ).$

\begin{defi}
A nonlinear connection $\stackrel{}{\stackrel{\vee}{N}}$ \thinspace is
called \textbf{induced} by the canonical nonlinear connection $N$ if the following
conditions hold:
\begin{equation}
\begin{array}{l}
\delta v^{(1)\alpha }=A_i^\alpha \delta y^{(1)i},...,\delta v^{(k-1)\alpha
}=A_i^\alpha \delta y^{(k-1)i} \\
\\
\delta \stackrel{\vee}{p}_\alpha =B_\alpha ^i\delta p_i,
\end{array}
\tag{9.3.1}
\end{equation}
\end{defi}

The previous conditions have a geometrical meaning. In fact, to a change
of coordinates on $T^{*k}\!\!\stackrel{\vee}{M}$ the previous formulae are
preserved. Tha same happens if we change the coordinates on the manifold $%
T^{*k}M.$

\begin{teo}
The dual coefficients of the induced nonlinear connection $\stackrel{\vee}{N}
$ are given by the following formulae
\begin{equation}
\stackrel{\vee}{\underset{\left( 1\right) }{M_\beta ^\alpha
}}=A_i^\alpha
\stackrel{\vee}{\underset{\left( 1\right) }{M_\beta ^i}},...,\stackrel{%
\vee}{\underset{\left( k-1\right) }{M_\beta ^\alpha }}=A_i^\alpha
\stackrel{\vee}{\underset{\left( k-1\right) }{M_\beta ^i}},\quad
\stackrel{\vee}{N_{\alpha \beta }}=\stackrel{\vee}{N_{\alpha
i}}B_\beta ^i, \tag{9.3.2}
\end{equation}
\end{teo}
where
\begin{equation}
\begin{array}{l}
\underset{(1)}{\stackrel{\vee}{M_\beta ^i}}=\underset{\left(
1\right) }{ M_j^i}\displaystyle\frac{\partial y^{\left( 1\right)
j}}{\partial v^{\left( 1\right) \beta
}}+\displaystyle\frac{\partial y^{\left( 1\right) i}}{
\partial u^\beta }, \\
\\
\underset{(2)}{\stackrel{\vee}{M_\beta ^i}}=\underset{\left(
2\right) }{ M_j^i}\displaystyle\frac{\partial y^{\left( 2\right)
j}}{\partial v^{\left( 2\right) \beta
}}+\underset{(1)}{M_j^i}\displaystyle\frac{\partial y^{\left(
2\right) j}}{\partial v^{\left( 1\right) \beta }}+\displaystyle
\frac{\partial y^{\left( 2\right) i}}{\partial u^\beta }, \\
...................................................... \\
\underset{(k-1)}{\stackrel{\vee}{M_\beta ^i}}=\underset{\left(
k-1\right) }{M_j^i}\displaystyle\frac{\partial y^{\left(
k-1\right) j}}{
\partial v^{\left( k-1\right) \beta }}+\underset{(k-2)}{M_j^i}%
\displaystyle \frac{\partial y^{\left( k-1\right) j}}{\partial
v^{\left( k-2\right) \beta }
}+...+\underset{(1)}{M_j^i}\displaystyle\frac{\partial
y^{\left( k-1\right) j}}{\partial v^{\left( 1\right) \beta }}+\displaystyle%
\frac{ \partial y^{\left( k-1\right) i}}{\partial u^\alpha }
\end{array}
\tag{9.3.3}
\end{equation}
and
\begin{equation}
\stackrel{\vee}{N_{\alpha i}}=N_{ji}B_\alpha ^j-\displaystyle
\frac{\partial p_i}{\partial u^\alpha }.  \tag{9.3.4}
\end{equation}

\textbf{Proof}. The first equations (9.3.1) lead to
$dv^{\left( 1\right) \alpha }+\underset{(1)}{\stackrel{\vee
}{M_\beta
^\alpha }}du^\beta =A_i^\alpha (dy^{\left( 1\right) i}+\underset{(1)}{%
M_j^i }dx^j)=A_i^\alpha (\displaystyle\frac{\partial y^{\left( 1\right) i}}{\partial
u^\beta }+\displaystyle\frac{\partial y^{\left( 1\right) i}}{\partial
v^{\left( 1\right) \beta }}dv^{\left( 1\right) \beta }+\underset{(1)}{%
M_j^i }B_\beta ^\alpha du^\beta ). $

Because of $\displaystyle\frac{\partial y^{\left( 1\right) i}}{\partial
v^{\left( 1\right) \alpha }}=B_\alpha ^i$ and $A_i^\alpha B_\beta ^i=\delta_\beta ^\alpha $ the last equality leads to (9.3.2):
$\stackrel{\vee}{\underset{(1)}{M_\beta ^\alpha }}=A_i^\alpha
\stackrel{ \vee}{\underset{(1)}{M_\beta ^i}},$ where
$\stackrel{\vee}{\underset{ (1)}{M_\beta ^i}}\,$ is given in (9.3.3), etc.

Remark that
\begin{equation}
\delta y^{\left( a\right) i}=dy^{\left( a\right)
i}+\underset{(1)}{M_j^i}
dy^{(a-1)i}+...+\underset{(a)}{M_j^i}dx^j,(a=1,...,k-1) \tag{9.3.5}
\end{equation}
and
\begin{equation}
dy^{\left( a\right) i}=\displaystyle\frac{\partial y^{\left( a\right) i}}{
\partial u^\alpha }du^\alpha +...+\displaystyle\frac{\partial y^{\left(
a\right) i}}{\partial v^{\left( a-1\right) \alpha }}dv^{(a-1)\alpha
}+B_\alpha ^idv^{\left( a\right) \alpha }.  \tag{9.3.5a}
\end{equation}

We can prove:
\begin{prop}
The following formulae hold:
\begin{equation}
\begin{array}{l}
dx^i=B_\beta ^idu^\beta,  \\
\\
\delta y^{\left( 1\right) i}=B_\beta ^i\delta v^{\left( 1\right) \beta }+A_{%
\overline{\alpha }}^i\underset{(1)}{H_\beta ^{\overline{\alpha
}}}du^\beta,
\\
................................................ \\
\delta y^{\left( k-1\right) i}=B_\beta ^i\delta v^{\left( k-1\right) \beta
}+A_{\overline{\alpha }}^i\underset{(1)}{(H_\beta ^{\overline{\alpha }}}%
dv^{\left( k-2\right) \beta }+...+\underset{(k-1)}{H_\beta ^{\overline{%
\alpha }}}du^\beta ), \\
\\
\delta p_i=A_i^\beta \delta \stackrel{\vee}{p}_\beta -A_i^{\overline{\alpha
}}H_{\overline{\alpha }\beta }du^\beta,
\end{array}
\tag{9.3.6}
\end{equation}
where
\begin{equation}
A_{\overline{\alpha }}^i\underset{(1)}{H_\beta ^{\overline{\alpha }}}=%
\stackrel{\vee}{\underset{(1)}{M_\beta ^i}}-B_\alpha ^i\stackrel{\vee}{%
\underset{(1)}{M_\beta ^\alpha }};...;A_{\overline{\alpha }}^i\underset{%
(k-1)}{H_\beta ^{\overline{\alpha }}}=\stackrel{\vee}{\underset{(k-1)}{%
M_\beta ^i}}-B_\alpha ^i\stackrel{\vee}{\underset{(k-1)}{M_\beta ^\alpha }%
}, \tag{9.3.7}
\end{equation}
\begin{equation}
A_i^{\overline{\alpha }}H_{\overline{\alpha }\beta }=\stackrel{\vee}{N}%
_{\beta i}.  \tag{9.3.7a}
\end{equation}
\end{prop}

\textbf{Proof}. The formulae (9.3.5), (9.3.5') and (9.3.3) have as consequence the formulae (9.3.6).
Also, we have $A_i^{\overline{\alpha }}H_{\overline{\alpha }
\beta }=-A_\beta ^i\delta p_i=\stackrel{\vee}{N}_{\beta i}.$ q.e.d.

The previous expressions (9.3.6) is not convenient for us, because
in the right hand sides we have the natural cobasis $du^\alpha
,dv^{\left( 1\right) i},...,dv^{\left( k-2\right) i}.$ This means
that the corresponding coefficients $A_{\overline{\alpha
}}^i\underset{(a)}{H_\beta ^{\overline{ \alpha }}},(a=1,...,k-1)$
do not have a geometrical meaning.

So, we prove:
\begin{teo}
In every point $\stackrel{\vee}{u}\in T^{*k}\!\!\stackrel{\vee}{M}$ the
following formulas hold:
\begin{equation}
\begin{array}{l}
dx^i=B_\beta ^idu^\beta,  \\
\\
\delta y^{\left( 1\right) i}=B_\beta ^i\delta y^{\left( 1\right) \beta }+A_{%
\overline{\alpha }}^i\underset{(1)}{K_\beta ^{\overline{\alpha
}}}du^\beta,
\\
\\
\delta y^{\left( 2\right) i}=B_\beta ^i\delta v^{\left( 2\right) \beta }+A_{%
\overline{\alpha }}^i\underset{(1)}{(K_\beta ^{\overline{\alpha
}}}\delta
v^{\left( 1\right) \beta }+\underset{(2)}{K_\beta ^{\overline{\alpha }}}%
du^\beta ), \\
............................................................................
\\
\delta y^{\left( k-1\right) i}=B_\beta ^i\delta v^{\left( k-1\right) \beta
}+A_{\overline{\alpha }}^i\underset{(1)}{(K_\beta ^{\overline{\alpha }}}%
\delta v^{\left( k-2\right) \beta }+\underset{(2)}{K_\beta ^{\overline{%
\alpha }}}\delta v^{\left( k-3\right) \beta
}+...+\underset{(k-1)}{K_\beta ^{\overline{\alpha }}}du^\beta )
\end{array}
\tag{9.3.8}
\end{equation}
and
$$
\delta p_i=A_i^\beta \delta \stackrel{\vee}{p}_\beta -A_i^{\overline{%
\alpha }}H_{\overline{\alpha }\beta }du^\beta ,
$$
where $A_i^{\overline{\alpha }}H_{\overline{\alpha }\beta }$ is given by
(9.3.7') and
\begin{equation}
\begin{array}{l}
\underset{(1)}{K_\beta ^{\overline{\alpha }}}=\underset{(1)}{H_\beta ^{%
\overline{\alpha }}}, \\
\\
\underset{(2)}{K_\beta ^{\overline{\alpha }}}=\underset{(2)}{H_\beta ^{%
\overline{\alpha }}}-\underset{(1)}{H_\gamma ^{\overline{\alpha }}}%
\stackrel{\vee}{\underset{(1)}{N_\beta ^\alpha }}, \\
............................ \\
\underset{(k-1)}{K_\beta ^{\overline{\alpha
}}}=\underset{(k-1)}{H_\beta
^{\overline{\alpha }}}-\underset{(1)}{H_\gamma ^{\overline{\alpha }}}%
\stackrel{\vee}{\underset{(k-2)}{N_\beta ^\gamma }}-...-\underset{(k-2)%
}{H_\gamma ^{\overline{\alpha }}}\stackrel{\vee
}{\underset{(1)}{N_\beta ^\gamma }}.
\end{array}
\tag{9.3.9}
\end{equation}
\end{teo}

Of course, $\stackrel{\vee}{\underset{(1)}{N_\beta ^\alpha
}},..., \stackrel{\vee}{\underset{(k-1)}{N_\beta ^\alpha }}$ are
the coefficients of the induced connection $\stackrel{\vee
}{N.}$

\textbf{Proof}. Taking into account the formulas (6.3.5) we have
\[
dv^{\left( a\right) \alpha }=\delta v^{\left( a\right) \alpha
}-\stackrel{ \vee}{\underset{(1)}{N_\beta ^\alpha }}\delta
v^{\left( a-1\right) \beta }-...-\stackrel{\vee
}{\underset{(a)}{N_\beta ^\alpha }}\delta u^\beta ,(a=1,...,k-1).
\]
So, $\delta y^{\left( a\right) i}$ from (9.3.6) is as follows

$\delta y^{\left( a\right) i}=B_\beta ^i\delta v^{\left( a\right) \beta }+A_{%
\overline{\alpha }}^i[\underset{(1)}{H_\beta ^\alpha (}\delta
v^{\left( a-1\right) \gamma }-\stackrel{\vee
}{\underset{(1)}{N_\gamma ^\beta }} \delta v^{\left( a-2\right)
\gamma }-...-\stackrel{\vee}{\underset{(a-1)}{ N_\gamma ^\beta
}}\delta u^\gamma )$

$\qquad +\underset{(2)}{H_\beta ^{\overline{\alpha }}(}\delta
v^{\left( a-2\right) \beta }-\stackrel{\vee
}{\underset{(1)}{N_\gamma ^\beta }} \delta v^{\left( a-3\right)
\gamma }-...-\stackrel{\vee}{\underset{(a-2)}{ N_\gamma ^\beta
}}\delta u^\gamma )+...+$

\qquad $+\underset{(a-1)}{H_\beta ^{\overline{\alpha }}(}\delta
v^{\left( 1\right) \beta }-\stackrel{\vee
}{\underset{(1)}{N_\gamma ^\beta }}\delta u^\gamma
)+\stackrel{\vee}{\underset{(a)}{N_\beta ^{\overline{\alpha }}}}
\delta u^\beta ).$

Identifying with $\delta y^{\left( a\right) i}$ from (9.3.6): \vspace{3mm}

$\delta y^{\left( a\right) i}=B_\beta ^i\delta v^{\left( a\right) \beta }+A_{%
\overline{\alpha }}^i[\underset{(1)}{K_\beta ^{\overline{\alpha
}}}\delta v^{\left( a-1\right) \beta }+\underset{(2)}{K_\beta
^{\overline{\alpha }}} \delta v^{\left( a-2\right) \beta
}+...+\underset{(a)}{K_\beta ^{\overline{ \alpha }}}du^{\beta}]$\\
we obtain the formulae (9.3.9). q.e.d

The previous theorem has an important consequence.

\begin{cor}
With respect to the transformations of coordinates on the submanifold $%
T^{*k}\!\!\stackrel{\vee}{M}$ and to the transformation (9.2.5),
$\underset{(1)}{K_\beta
^{\overline{\alpha}}},...,\underset{(k-1)}{K_\beta
^{\overline{\alpha }}}$ and $H_{\overline{\alpha }\beta }\,$ are
the mixed $d$-tensor fields.
\end{cor}

Generally a set of functions $T_{j...\beta ...\overline{\beta }}^{i...\alpha
...\overline{\alpha }}(\stackrel{\vee}{u})$ which are the components of a $d$-tensor in the
indices $i,j,...$and $d$-tensor in the indices $\alpha ,\beta ,...$ and tensor
with respect to the transformation (9.2.5) in the indices $\overline{\alpha },
\overline{\beta },...$ is called a \textit{mixed} $d$-tensor field on the
submanifolds $T^{*k}\!\!\stackrel{\vee}{M}.$

For instance $B_\alpha ^i,A_\alpha ^i,A_{\overline{\alpha
}}^i,A_i^\alpha ,A_i^{\overline{\alpha }},\stackrel{\vee
}{g}_{\alpha \beta },g_{ij},\delta _{\overline{\alpha
}\overline{\beta }},\underset{(1)}{K_\beta ^{\overline{
\alpha }}},...,\underset{(k-1)}{K_\beta ^{\overline{\alpha }}}$ and $H_{%
\overline{\alpha }\beta }$ are \newline mixed $d$-tensor fields. The previous
definition of mixed $d$-tensor fields can be extended to any geometrical $d$
-object fields.

\section{The Relative Covariant Derivative}

In this section we shall construct an operator $\nabla $ of relative
covariant derivation in the algebra of mixed $d$-tensor fields. It is
clear that $\nabla $ will be well determined if we know its action on the
mixed $d$-vector field:

\begin{equation}
X^i(\stackrel{\vee}{u}),X^\alpha (\stackrel{\vee}{u}),X^{\overline{\alpha }
}(\stackrel{\vee}{u}),\forall \stackrel{\vee}{u}=(u,v^{(1)},...,v^{\left(
k-1\right) },\stackrel{\vee}{p})\in T^{*k}\!\!\stackrel{\vee}{M}.  \tag{9.4.1}
\end{equation}

\begin{defi}
We call a \textit{coupling} of the canonical metrical $N$-connection $D$ of
the Hamilton spaces of order $k,$ $H^{\left( k\right) n}=(M,H)$ \thinspace
to the induced nonlinear connection $\stackrel{\vee}{N}$, of the submanifold $\stackrel{%
\vee}{M}$, in the manifold $M$ the operator $\stackrel{\vee}{D}$ with the
property
\begin{equation}
\stackrel{\vee}{D}X^i=DX^i\text{ modulo \{(9.3.8), (9.3.8')\}},  \tag{9.4.2}
\end{equation}
where
\begin{equation}
DX^i=dx^i+X^i\omega _j^i  \tag{9.4.2a}
\end{equation}
and where $\omega^i_j$ are 1-forms connection:
\begin{equation}
\omega _j^i=H_{jh}^idx^h+\underset{(1)}{C_{jh}^i}\delta y^{\left(
1\right) h}+...+\underset{(k-1)}{C_{jh}^i}\delta y^{\left(
k-1\right) h}+C_j^{ih}\delta p_h.  \tag{9.4.2b}
\end{equation}
\end{defi}

Thus, in every point $\stackrel{\vee}{u}$ $\in T^{*k}\stackrel{\vee
}{M}$, $\stackrel{\vee}{D}X^i$ has the form
\begin{equation}
\stackrel{\vee}{D}X^i=dx^i+X^j\stackrel{\vee}{\omega }_j^i , \tag{9.4.3}
\end{equation}
where
\begin{equation}
\stackrel{\vee}{\omega }_j^i=\stackrel{\vee}{H_{j\beta
}^i}du^\beta + \stackrel{\vee}{\underset{(1)}{C_{j\beta
}^i}}\delta v^{\left( 1\right) \beta
}+...+\underset{(k-1)}{\stackrel{\vee}{C_{j\beta }^i}}\delta
v^{\left( k-1\right) \beta }+\stackrel{\vee}{C_j^{i\beta }}\delta \stackrel{\vee
}{p}_\beta. \tag{9.4.4}
\end{equation}

\begin{teo}
The connection one forms $\stackrel{\vee}{\omega }_j^i$ of $\stackrel{%
\vee}{D}$ are given by (9.4.4), where:
\begin{equation}
\begin{array}{l}
\stackrel{\vee}{H_{j\beta}^i}=H_{j\beta}^iB_\beta^h+(\underset{(1)}{C_{jh}^i}\underset{\left(1\right)}
{K_\beta^{\overline{\alpha}}}+...+{\underset{(k-1)}{C_{jh}^i}}\underset{\left(k-1\right)}
{K_\beta^{\overline{\alpha}}})A_{\overline{\alpha}}^h-C_i^{jh}A_h^{\overline{\alpha}}H_{\overline{\alpha}\beta}, \\
\\
\stackrel{\vee}{\underset{(1)}{C_{j\beta }^i}}=\underset{(1)}{C_{jh}^i}%
B_\beta ^h+(\underset{(2)}{C_{jh}^i}\underset{\left( 1\right)
}{K_\beta ^{\overline{\alpha
}}}+...+\underset{(k-1)}{C_{jh}^i}\underset{\left(
k-2\right) }{K_\beta ^{\overline{\alpha }}})A_{\overline{\alpha }}^h, \\
...................................................... \\
\stackrel{\vee}{\underset{(k-1)}{C_{j\beta }^i}}=\underset{(k-1)}{%
C_{jh}^i}B_\beta ^h+\underset{(k-2)}{C_{jh}^i}\underset{\left(
1\right)
}{K_\beta ^{\overline{\alpha }}}A_{\overline{\alpha }}^h, \\
\\
\stackrel{\vee}{C}_i^{j \beta }=C_i^{jh}A_h^\beta.
\end{array}
\tag{9.4.5}
\end{equation}
\end{teo}

\textbf{Proof}. Using the formulae (9.3.8), (9.3.8'), $\omega _j^i$ modulo
(9.3.8), (9.3.8') from (9.4.2'') leads to
\begin{equation}
\begin{array}{ll}
\stackrel{\vee}{\omega }_j^i=H_{jh}^iB_\beta ^hdu^\beta + &
\underset{(1) }{C_{jh}^i}(B_\beta ^h\delta v^{\left( 1\right)
\beta }+A_{\overline{\alpha }
}^h\underset{(1)}{K_\beta ^{\overline{\alpha }}}du^\beta) +...+ \\
& \\
& \underset{(k-1)}{C_{jh}^i}[B_\beta ^h\delta v^{\left( k-1\right)
\beta }+A_{\overline{\alpha }}^h\underset{(1)}{(K_\beta
^{\overline{\alpha }}}
\delta v^{\left( k-2\right) \beta }+...+\underset{(k-1)}{K_\beta ^{%
\overline{\alpha }}}du^\beta )]+ \\
& \\
& C_i^{jh}[A_\beta ^h\delta \stackrel{\vee}{p}_\beta -A_h^{\overline{\alpha
}}H_{\overline{\alpha }\beta }du^\beta ].
\end{array}
\tag{9.4.4a}
\end{equation}

Identifying with $\stackrel{\vee}{\omega }_j^i$ from (9.4.4) we obtain the
formulae (9.4.5). Evidently the coupling connection $\stackrel{\vee}{D}$ of the
canonical metrical $N$-connection $D$, depend only on the fundamental
function $H$ of the space $H^{(k)n}$ and by the immersion $i^{*}:T^{*k}
\!\!\stackrel{\vee}{M}$ in $T^{*k}M.$

Of course, we can write $\stackrel{\vee}{D}X^i$ in the form
\begin{equation}
\stackrel{\vee}{D}X^i=X_{|\alpha }^idu^\alpha +X^i\stackrel{\left( 1\right)
}{|_\alpha }\delta v^{\left( 1\right) \alpha }+...+X^i\stackrel{\left(k-1\right) }
{|_\alpha }\delta v^{\left( k-1\right) \alpha }+X^i|^\alpha
\delta \stackrel{\vee}{p}_\alpha , \tag{9.4.6}
\end{equation}
where
\begin{equation}
\begin{array}{l}
X_{|\alpha }^i=\displaystyle\frac{\delta X^i}{\delta u^\alpha }+X^j\stackrel{
\vee}{H}_{j\alpha }^i, \\
\\
X^i\stackrel{\left( a\right) }{|}_\alpha
=\displaystyle\frac{\delta X^i}{ \delta v^{(a)\alpha
}}+X^j\underset{\left( a\right) }{\stackrel{\vee}{C}
_{j\alpha }^i},(a=1,...,k-1), \\
\\
X^i|^\alpha =\stackrel{\cdot}{\partial ^\alpha} X^i+X^j\stackrel{\vee}{C}\text{}_j^{i\alpha }.
\end{array}
\tag{9.4.7}
\end{equation}

We can extend, without difficulties the action of the linear connection $%
\stackrel{\vee}{D}$ to the $d$-tensors $T_{j_i...j_s}^{i_1...i_r}(\stackrel{
\vee}{u})$, $\forall \stackrel{\vee}{u}\in T^{*k}\!\!\stackrel{\vee}{M}.$

\begin{defi}
We call the induced tangent connection on $T^{*k}\!\!\stackrel{\vee}{M}$ by the
canonical metrical $N$-connection $D$ the operator \thinspace $D^{T\text{ }}$
given by
\begin{equation}
D^TX^\alpha =A_i^\alpha \stackrel{\vee}{D}X^i,\quad for\quad X^i=A_\alpha
^iX^\alpha   \tag{9.4.8}
\end{equation}
\end{defi}

We can see that:
\begin{equation}
D^TX^\alpha =dX^\alpha +X^\beta \omega _\beta ^\alpha , \tag{9.4.9}
\end{equation}
where $\omega _\beta ^\alpha $ are the connection one forms of $D^T:$
\begin{equation}
\omega _\beta ^\alpha =H_{\beta \gamma }^\alpha du^\gamma
+\underset{(1)}{ C_{\beta \gamma }^\alpha }\delta v^{\left(
1\right) \gamma }+...+\underset{ (k-1)}{C_{\beta \gamma }^\alpha
}\delta v^{\left( k-1\right) \gamma }+C_\beta ^{\alpha \gamma
}\delta \stackrel{\vee}{p}_\gamma  \tag{9.4.10}
\end{equation}

\begin{teo}
1$^{\circ }$ The coefficients ($H_{\beta \gamma }^\alpha ,\underset{(1)}{%
C_{\beta \gamma }^\alpha },...,\underset{(k-1)}{C_{\beta \gamma }^\alpha }%
,C_\beta ^{\alpha \gamma }$) of the one forms connection $\omega _\beta
^\alpha $ have the following expressions
\begin{equation}
\begin{array}{c}
H_{\beta \gamma }^\alpha =A_i^\alpha (\displaystyle\frac{\delta A_\beta ^i}{%
\delta u^\gamma }+A_\beta ^j\stackrel{\vee}{H_{j\gamma }^i}), \\
\\
\underset{(a)}{C_{\beta \gamma }^\alpha }=A_i^\alpha (\displaystyle\frac{%
\delta A_\beta ^i}{\delta v^{(a)\gamma }}+A_\beta ^h\underset{(a)}{%
\stackrel{\vee}{C_{h\gamma }^\alpha }}), (a=1,..., k-1), \\
\\
C_\beta ^{\alpha \gamma }=A_i^\alpha (\stackrel{\cdot }{\partial ^\gamma }%
A_\beta ^i+A_\beta ^h\stackrel{\vee}{C_h^{i\gamma }}).
\end{array}
\tag{9.4.11}
\end{equation}

2$^{\circ }$ With respect of the change of local coordinates on the
manifold $T^{*k}\!\!\stackrel{\vee}{M},$ the coefficients
\[
H_{\beta \gamma }^\alpha ,\ \underset{(1)}{C_{\beta \gamma }^\alpha },...,
\underset{(k-1)}{\ C_{\beta \gamma }^\alpha },..., C_\beta ^{\alpha
\gamma }
\]
have the rule of transformation as the coefficients of $\stackrel{\vee}{N}$
-linear connection on $T^{*k}\!\!\stackrel{\vee}{M}.$
\end{teo}

\textbf{Proof}. 1$^{\circ }$ The equalities (9.4.8) and (9.4.9) imply
$$
dX^\alpha +X^\beta \omega _\beta ^\alpha =A_i^\alpha (dX^i+X^j\stackrel{
\vee}{\omega }_j^i)=A_i^\alpha [(dA_\beta ^i+A_\beta ^h\stackrel{\vee}{
\omega }_h^i)X^\beta +A_\beta ^idx^\beta ]
$$

Consequently,
\begin{equation}
\omega _\beta ^\alpha =A_i^\alpha (dA_\beta ^i+\stackrel{\vee}{\omega }
_j^iA_\beta ^j).  \tag{9.4.12}
\end{equation}

Thus, the equalities (9.4.4), (9.4.10) and (9.4.12) lead to the form (9.4.11) of the coefficients of induced tangent connection $D^T.$

\qquad 2$^{\circ }$ Using the rule of transformations of the coefficients of
$\stackrel{\vee}{\omega}\text{}_j^i$ it follows the sentence 2$^{\circ }$ of
theorem . q.e.d.

As in the case of $\stackrel{\vee}{D}$ we can write:
\begin{equation}
D^TX^\alpha =X_{|\beta }^\alpha du^\beta +X^\alpha \stackrel{\left( 1\right)
}{|}_\beta \delta v^{\left( 1\right) \beta }+...+X^\alpha \stackrel{\left(
k-1\right) }{|}_\beta \delta v^{\left( k-1\right) \beta }+X^\alpha |^\beta
\delta \stackrel{\vee}{p}_\beta,  \tag{9.4.13}
\end{equation}
where
\begin{equation}
\begin{array}{l}
X_{|\gamma }^\alpha =\displaystyle\frac{\delta X^\alpha }{\delta u^\gamma }
+X^\beta H_{\beta \gamma }^\alpha, \\
\\
X^\alpha \stackrel{\left( a\right) }{|}_\gamma
=\displaystyle\frac{\delta X^\alpha }{\delta v^{(a)\gamma
}}+X^\beta \underset{\left( a\right) }{
C_{\beta \gamma }^\alpha },\ \ (a=1,...k-1), \\
\\
X^\alpha |^\gamma =\stackrel{\cdot }{\partial ^\gamma }X^\alpha +X^\beta
C_\beta ^{\alpha \gamma }.
\end{array}
\tag{9.4.14}
\end{equation}

Evidently, the induced tangent connection $D^T$ depends on the canonical $%
\stackrel{\vee}{N}$- metrical connection $D$ and depends on the immersion $%
i^{*}:T^{*k}\!\!\stackrel{\vee}{M}\rightarrow T^{*k}M.$

This operator $D^T$, as well the operators ,,$_{|\gamma }"$ , $_{,,}\stackrel{%
\left( a\right) }{|}_\gamma "$ and $_{,,}|^\gamma "$ can be extended to the
tensor fields $T_{\beta _1...\beta _q}^{\alpha _1...\alpha _a}(\stackrel{
\vee}{u})$.

\begin{defi}
We call the \textit{induced normal connection} on $T^{*k}\!\!\stackrel{\vee}{M}$
by the canonical metrical $N$-connection D, the operator $D^{\perp \text{ }}$
given by
\begin{equation}
D^{\perp \text{}}X^{\overline{\alpha }}=A_i^{\overline{\alpha }}\stackrel{%
\vee}{D}X^i,\ for\ X^i=A_{\overline{\alpha }}^iX^{\overline{\alpha }}.
\tag{9.4.15}
\end{equation}
\end{defi}

As before we set
\begin{equation}
D^{\perp \text{ }}X^{\overline{\alpha }}=dX^{\overline{\alpha }}+X^{%
\overline{\beta }}\omega _{\overline{\beta }}^{\overline{\alpha }},
\tag{9.4.16}
\end{equation}
with the connection $1$-forms of $D^{\perp \text{}}.$
\begin{equation}
\omega _{\overline{\beta }}^{\overline{\alpha
}}=H_{\overline{\beta }\gamma }^{\overline{\alpha }}du^\gamma
+\underset{(1)}{C_{\overline{\beta }\gamma }^{\overline{\alpha
}}}\delta v^{\left( 1\right) \gamma }+...+\underset{
(k-1)}{C_{\overline{\beta }\gamma }^{\overline{\alpha }}}\delta
v^{\left( k-1\right) \gamma }+C_{\overline{\beta
}}^{\overline{\alpha }\gamma }\delta \stackrel{\vee}{p}_\gamma.
\tag{9.4.17}
\end{equation}

Applying the same method as in the case of $\stackrel{\vee}{D}$ and $D^T$
we can calculate the coefficients of the $1$-forms connection $\omega _{%
\overline{\beta }}^{\overline{\alpha }}.$

\begin{teo}
The coefficients $(H_{\overline{\beta }\gamma }^{\overline{\alpha }},%
\underset{(1)}{C_{\overline{\beta }\gamma }^{\overline{\alpha }}},...,%
\underset{(k-1)}{C_{\overline{\beta }\gamma }^{\overline{\alpha }}},C_{%
\overline{\beta }}^{\overline{\alpha }\gamma })$ of the induced normal
connection $D^{\perp }$ are given by
\begin{equation}
\begin{array}{l}
H_{\overline{\beta }\gamma }^{\overline{\alpha }}=A_i^{\overline{\alpha }}(%
\displaystyle\frac{\delta A_{\overline{\beta }}^i}{\delta u^\gamma }+A_{%
\overline{\beta }}^j\stackrel{\vee}{H}\text{}_{j\gamma}^i), \\
\\
\underset{(a)}{C_{\overline{\beta }\gamma }^{\overline{\alpha }}}=A_i^{%
\overline{\alpha }}(\displaystyle\frac{\delta A_{\overline{\beta }}^i}{%
\delta v^{(\alpha )\gamma }}+A_{\overline{\beta }}^j\underset{(a)}{%
\stackrel{\vee}{C}\text{}_{j\gamma}^i}), (a=1,...,k-1), \\
\\
C_{\overline{\beta }}^{\overline{\alpha }\gamma }=A_i^{\overline{\alpha }}(%
\stackrel{\cdot }{\partial ^\gamma }A_{\overline{\beta }}^i+A_\gamma ^i%
\stackrel{\vee}{C}_{\overline{\beta }}^{\overline{\alpha }\gamma }).
\end{array}
\tag{9.4.18}
\end{equation}
\end{teo}

\textbf{Proof}. Indeed, (9.4.15) and (9.4.16) have as consequence:
\begin{equation}
\omega _{\overline{\beta }}^{\overline{\alpha }}=A_i^{\overline{\alpha }
}(dA_{\overline{\beta }}^i+A_{\overline{\beta }}^j\stackrel{\vee}{\omega }\text{}_j^i).
\tag{9.4.19}
\end{equation}

Therefore, $\stackrel{\vee}{\omega}\text{}_j^i$ from (9.4.4) and $\omega _{%
\overline{\beta }}^{\overline{\alpha }}$ expressed in (9.4.17) substituted in
the equality (9.4.19) leads to (9.4.18). q.e.d.

As usual, we may set for the induced normal connection
\begin{equation}
D^{\perp }X^{\overline{\alpha }}=X_{|\beta }^{\overline{\alpha }}du^\beta
+X^{\overline{\alpha }}\stackrel{\left( 1\right) }{|_\beta }\delta v^{\left(
1\right) \beta }+...+X^{\overline{\alpha }}\stackrel{\left( k-1\right) }{
|_\beta }\delta v^{\left( k-1\right) \beta }+X^{\overline{\alpha }}|^\beta
\delta \stackrel{\vee}{p}_\beta  \tag{9.4.20}
\end{equation}
with
\begin{equation}
\begin{array}{l}
X_{|\beta }^{\overline{\alpha }}=\displaystyle\frac{\delta X^{\overline{
\alpha }}}{\delta u^\beta }+X^{\overline{\gamma }}H_{\overline{\gamma }\beta
}^{\overline{\alpha }}, \\
\\
X^{\overline{\alpha }}\stackrel{\left( a\right) }{|}_\beta
=\displaystyle \frac{\delta X^{\overline{\alpha }}}{\delta
v^{(a)\beta }}+X^\gamma \underset{\left( a\right)
}{C_{\overline{\gamma }\beta }^{\overline{\alpha
}}},(a=1,...k-1), \\
\\
X^{\overline{\alpha }}|^\beta =\stackrel{\cdot }{\partial ^\beta }X^{%
\overline{\alpha }}+X^\gamma C_{\overline{\gamma }}^{\overline{\alpha }\beta
}.
\end{array}
\tag{9.4.21}
\end{equation}

These are the induced normal covariant derivations.

Now, we can define the relative (or mixed) derivation $\nabla $ introduced at
the begining of this section.
\begin{defi}
A \textit{relative} (mixed) \textit{covariant} derivation in the algebra of
mixed $d$-tensor field is an operator $\nabla $ which has the following
properties.
\begin{equation}
\begin{array}{l}
\nabla f=df,\forall f\in \mathcal{F}(T^{*k}\!\!\stackrel{\vee}{M}) \\
\\
\nabla X^i=\stackrel{\vee}{D}X^i,\nabla X^\alpha =D^TX^\alpha ,\nabla X^{%
\overline{\alpha }}=D^{\perp }X^{\overline{\alpha }}
\end{array}
\tag{9.4.22}
\end{equation}
for any mixed vector fields $X^i(\stackrel{\vee}{u}),X^\alpha (\stackrel{%
\vee}{u})$ and $X^{\overline{\alpha }}(\stackrel{\vee}{u}).$
\end{defi}

The connection $1$-forms $\stackrel{\vee}{\omega }_j^i,\omega _\beta
^\alpha $ and $\omega _{\overline{\beta }}^{\overline{\alpha }}$ will be
called the \textit{connections} $1$-forms of the relative covariant
derivation $\nabla $.
Of course, the operator $\nabla $ can be extended to any mixed $d$-tensor $%
T_{j...\beta ...\overline{\beta }}^{i...\alpha ...\overline{\alpha }}(
\stackrel{\vee}{u})$, $\stackrel{\vee}{u}\in T^{*k}\!\!\stackrel{\vee}{M}.$

We can write the Ricci identities for the relative covariant derivative $%
\nabla ,$ and its torsion and curvature tensor, taking into account every
components $\stackrel{\vee}{D},D^T,$and $D^{\perp }$ of $\nabla .$

Then, it is easy to prove:

\begin{teo}
The structure equations of the mixed covariant derivation $\nabla $ are as
follows:
\begin{equation}
\begin{array}{l}
d(du^\alpha )-du^\beta \wedge \omega _\beta ^\alpha =-\stackrel{\left(
0\right) }{\Omega ^\alpha } ,\\
\\
d(dv^{(a)\alpha })-\delta v^{(a)\beta }\wedge \omega _\beta ^\alpha =-%
\stackrel{\left( a\right) }{\Omega ^\alpha },(\alpha =1,...,k-1), \\
\\
d(\delta \stackrel{\vee}{p}_\alpha )-\delta \stackrel{\vee}{p}_\beta
\wedge \omega _\alpha ^\beta =-\Omega _\alpha
\end{array}
\tag{9.4.23}
\end{equation}
\end{teo}
and
\begin{equation}
\begin{array}{l}
d(\stackrel{\vee}{\omega _j^i})-\stackrel{\vee}{\omega }_j^h\wedge
\stackrel{\vee}{\omega }_h^j=-\stackrel{\vee}{\Omega _j^i}, \\
\\
d(\omega _\beta ^\alpha )-\omega _\beta ^\gamma \wedge \omega _\gamma
^\alpha =-\Omega _\beta ^\alpha, \\
\\
d(\omega _{\overline{\beta }}^{\overline{\alpha }})-\omega _{\overline{\beta
}}^{\overline{\gamma }}\wedge \omega _{\overline{\gamma }}^{\overline{\alpha
}}=-\Omega _{\overline{\beta }}^{\overline{\alpha }},
\end{array}
\tag{9.4.24}
\end{equation}
in which $\Omega ^\alpha ,\stackrel{\left( a\right) }{\Omega ^\alpha }$ and $%
\Omega _\alpha $ are the $2$-forms of torsion:
\begin{equation}
\begin{array}{l}
\stackrel{\left( 0\right) }{\Omega ^\alpha }=du^\beta \wedge \{\displaystyle %
\frac 12T_{\beta \gamma }^\alpha du^\gamma
+\underset{(1)}{C_{\beta \gamma }^\alpha }\delta v^{(1)\gamma
}+...+\underset{(k-1)}{C_{\beta \gamma }^\alpha }\delta
v^{(k-1)\gamma }+C_\beta ^{\alpha \gamma }\delta \stackrel{
\vee}{p}_\gamma\} ,\\
\\
\stackrel{\left( a\right) }{\Omega ^\alpha }=du^\beta \wedge \underset{%
(a0) }{P_\beta ^\alpha }+\underset{b=1}{\stackrel{k-1}{\sum
}}\delta v^{\left( b\right) \beta }\wedge \underset{(ab)}{P_\beta
^\alpha }+\delta
v^{(a)\beta }\wedge \{H_{\beta \gamma }^\alpha du^\gamma + \\
\\
\qquad +\underset{b=1}{\stackrel{k-1}{\sum
}}\underset{(b)}{C_{\beta \gamma }^\alpha }\delta v^{(b)\gamma
}+C_\beta ^{\alpha \gamma }\delta \stackrel{\vee}{p}_\gamma
\},(a=1,...,k-1),
\end{array}
\tag{9.4.25}
\end{equation}
with $\underset{(ab)}{P_\beta ^\alpha }$ from (5.4.8)
and where the $2$-forms of curvature are given by

\[
\begin{array}{c}
\stackrel{\vee}{\Omega }_j^i=\displaystyle\frac 12\stackrel{\vee}{R}\text{}_{j\ \alpha \beta }^i du^\alpha \wedge du^\beta
+\stackrel{k-1}{\underset{a=1 }{\sum \text{ }}}\underset{(a)\quad
\quad }{\stackrel{\vee}{P}\text{}_{j\ \alpha \beta }^i}du^\alpha \wedge
\delta v^{\left( a\right) \beta }+\stackrel{\vee}{P}\text{}_{j\ \alpha
}^{\ i\ \ \beta }du^\alpha \wedge \delta \stackrel{\vee}{p}_\beta +
\\
\\
\quad +\stackrel{k-1}{\underset{a,b=1}{\sum \text{ }}}\underset{
(ab)\quad \quad }{\stackrel{\vee}{S}\text{}_{j\ \alpha \beta
}^i}dv^{\left( a\right) \alpha }\wedge \delta v^{\left( b\right)
\beta }+\stackrel{k-1}{\underset{ a=1}{\sum \text{
}}}\stackrel{\vee}{S}_{j\ \alpha }^{\ i\ \ \beta }dv^{\left( a\right)
\alpha }\wedge \delta \stackrel{\vee}{p}_\beta
+\displaystyle\frac 12\stackrel{\vee}{S}\text{}_j^{\ i\alpha\beta }\delta
\stackrel{\vee}{p}_\alpha \wedge \stackrel{\vee}{p}_\beta ....,
\end{array}
\]

\begin{equation}
\begin{array}{c}
\Omega _\beta ^\alpha =\displaystyle\frac 12R_{\beta\  \gamma
\varphi }^\alpha du^\gamma \wedge du^\varphi
+\stackrel{k-1}{\underset{a=1}{\sum \text{ }}} \underset{\left(
a\right) \quad }{P_{\beta\ \gamma \varphi }^\alpha } du^\gamma
\wedge \delta v^{\left( a\right) \varphi }+P_{\beta \gamma
}^{\alpha \varphi }du^\gamma \wedge \delta \stackrel{\vee}{p}_\varphi + \\
\\
\quad +\stackrel{k-1}{\underset{a,b=1}{\sum \text{ }}}S_{\beta\
\gamma \varphi }^\alpha \delta v^{\left( a\right) \gamma }\wedge
\delta v^{\left( b\right) \varphi
}+\stackrel{k-1}{\underset{a=1}{\sum \text{ }}}S_{\beta\ \gamma
}^{\ \alpha \ \varphi }dv^{\left( a\right) \gamma }\wedge \delta
\stackrel{\vee}{p}_\varphi +\displaystyle\frac 12S_\beta ^{\ \alpha
\gamma \varphi }\delta \stackrel{\vee}{p}_\gamma \wedge \delta \stackrel{\vee}{p}_\varphi,
\end{array}
\tag{9.4.26}
\end{equation}
\[
\begin{array}{c}
\Omega _{\overline{\beta }}^{\overline{\alpha }}=\displaystyle\frac 12R_{%
\overline{\beta }\ \gamma \varphi }^{\overline{\alpha }}du^\gamma \wedge
du^\varphi +\stackrel{k-1}{\underset{a=1}{\sum \text{ }}}\underset{%
\left( a\right) \quad }{P_{\overline{\beta }\ \gamma \varphi }^{\overline{
\alpha }}}du^\gamma \wedge \delta v^{\left( a\right) \varphi }+P_{\overline{
\beta }\ \gamma }^{\ \overline{\alpha }\ \varphi }du^\gamma \wedge \delta
\stackrel{\vee}{p}_\varphi + \\
\\
\quad +\stackrel{k-1}{\underset{a,b=1}{\sum \text{
}}}S_{\overline{\beta }\  \gamma \varphi }^{\overline{\alpha
}}\delta v^{\left( a\right) \gamma }\wedge \delta v^{\left(
b\right) \varphi }+\stackrel{k-1}{\underset{a=1}{ \sum \text{
}}}S_{\overline{\beta }\ \gamma\varphi}^{\ \overline{\alpha }
}dv^{\left( a\right) \gamma }\wedge \delta \stackrel{\vee}{p}_\varphi + %
\displaystyle\frac 12S_{\overline{\beta }}^{\ \overline{\alpha }\gamma \varphi
}\delta \stackrel{\vee}{p}_\gamma \wedge \delta \stackrel{\vee}{p}_\varphi.
\end{array}
\]

The previous equations can be particularized in the case when the induced
tangent connection $D^T$ has vanishing tensors of torsion $T_{\ \beta \gamma
}^\alpha ,\underset{(a)}{S}\text{}_{\ \beta \gamma }^\alpha$ and $S_\alpha ^{\ \beta \gamma }.$

In the following we shall adopt the notations
\begin{equation}
\stackrel{\vee}{\Omega }_{ij}=\stackrel{\vee}{\Omega }_i^hg_{hj},\Omega
_{\alpha \beta }=\Omega _\alpha ^\gamma \stackrel{\vee}{g}_{\gamma \beta
},\Omega _{\overline{\alpha }\overline{\beta }}=\Omega _{\overline{\alpha }
}^{\overline{\gamma }}\delta _{\overline{\gamma }\overline{\beta }}.
\tag{9.4.27}
\end{equation}

These covariant $2$-forms of curvature of the mixed covariant connection $%
\nabla $ will be useful to write the fundamental equations of the
immersion $i^{*}:T^{*k}\!\!\stackrel{\vee}{M}\rightarrow T^kM.$

As a direct consequence of the Theorems from the present section we get:
\begin{teo}
The mixed covariant derivation $\nabla $ is a metrical one, i.e.:
\begin{equation}
\nabla g^{ij}=0,\nabla \delta ^{\overline{\alpha }\overline{\beta }%
}=0,\nabla \stackrel{\vee}{g}^{\alpha \beta }=0.  \tag{9.4.28}
\end{equation}
\end{teo}

\textbf{Proof}. The first formula is immediate, because of $\nabla g^{ij}=
\stackrel{\vee}{D}g^{ij}=0$ while $\nabla \delta ^{\overline{\alpha }
\overline{\beta }}=0,\nabla \stackrel{\vee}{g}^{\alpha \beta }=0$ can be
proved by means of Gauss-Weingarten formulas, given in next section.

\section{The Gauss-Weingarten Formula}

We need to study the moving equations of the Darboux frame:
\begin{equation}
\mathcal{R}=\{\stackrel{\vee}{u},A_\alpha ^i(\stackrel{\vee}{u}),A_{%
\overline{\alpha }}^i(\stackrel{\vee}{u})\},\forall \stackrel{\vee}{u}\in
T^{*k}\!\!\stackrel{\vee}{M} . \tag{9.5.1}
\end{equation}
So, we obtain
\begin{teo}
The following Gauss-Weingarten formulae hold good:
\begin{equation}
\bigtriangledown A_\alpha ^i=A_{\overline{\alpha }}^i\Pi ^{\overline{\alpha }}_\alpha ,\quad \bigtriangledown A_{\overline{\alpha }}^i=-A_\alpha ^i\Pi _{%
\overline{\alpha }}^\alpha  , \tag{9.5.2}
\end{equation}
where
\begin{equation}
\begin{array}{l}
\prod_\alpha ^{\overline{\beta }}=\underset{\left( 0\right) \quad }{%
H_{\alpha \gamma }^{\overline{\beta }}}du^\gamma +\underset{\left(
1\right) \quad }{H_{\alpha \gamma }^{\overline{\beta }}}\delta
v^{\left( 1\right) \gamma }+...+\underset{\left( k-1\right) \quad
}{H_{\alpha \gamma
}^{\overline{\beta }}}\delta v^{\left( k-1\right) \gamma }+H_\alpha ^{%
\overline{\beta }\gamma }\delta \stackrel{\vee}{p}_\gamma , \\
\\
\prod_{\overline{\beta }}^\alpha =\stackrel{\vee\quad }{g^{\alpha \gamma }}%
\delta _{\overline{\beta }\overline{\gamma }}\prod_\gamma ^{\overline{\gamma
}}
\end{array}
\tag{9.5.3}
\end{equation}
and where
\begin{equation}
\begin{array}{l}
\underset{\left( 0\right) \quad }{H_{\alpha \gamma }^{\overline{\beta }}}%
=A_h^{\overline{\beta }}(\displaystyle\frac{\delta A_\alpha ^h}{\delta
u^\gamma }+A_\alpha ^j\stackrel{\vee\quad }{H_{j\gamma }^h}), \\
\\
\underset{\left( a\right) \quad }{H_{\alpha \gamma }^{\overline{\beta }}}%
=A_h^{\overline{\beta }}(\displaystyle\frac{\delta A_\alpha
^h}{\delta v^{(a)\gamma }}+A_\alpha ^j\stackrel{\vee\quad
}{\underset{\left(
a\right) \quad }{C_{j\gamma }^h}}),(a=1,...,k-1), \\
\\
H_\alpha ^{\overline{\beta }\gamma }=A_h^{\overline{\beta }}(\stackrel{\cdot
}{\partial^\gamma} A_\alpha ^h+A_\alpha ^j\stackrel{\vee\quad }{%
C_j^{h\gamma }}),
\end{array}
\tag{9.5.4}
\end{equation}
with $\stackrel{\vee\quad }{H_{j\gamma }^h},\stackrel{\vee\quad }{%
\underset{\left( a\right) \quad }{C_{j\gamma }^h}}$ and
$\stackrel{\vee\quad }{C_j^{h\gamma }}$ are from (9.4.5).
\end{teo}

\textbf{Proof}. Taking into account that $A_\alpha ^i$ is a mixed $d$-tensor we
have for its relative covariant derivation:

$
\begin{array}{ll}
\bigtriangledown A_\alpha ^i= & dA_\alpha ^i+A_\alpha ^j\stackrel{\vee\quad
}{\omega _j^i}-A_\beta ^i\omega _\alpha ^\beta =dA_\alpha ^i+A_\alpha ^j
\stackrel{\vee\quad }{\omega _j^i}- \\
& -A_\gamma ^iA_h^\gamma (dA_\alpha ^h+\stackrel{\vee\quad }{\omega _j^h}
A_\alpha ^j)=dA_\alpha ^i+A_\alpha ^j\stackrel{\vee\quad }{\omega _j^i}- \\
& -(\delta _h^i-A_{\overline{\gamma }}^iA_h^{\overline{\gamma }})(dA_\alpha
^h+\stackrel{\vee\quad }{\omega _j^h}A_\alpha ^j)=A_{\overline{\alpha }
}^iA_h^{\overline{\alpha }}(dA_\alpha ^h+\stackrel{\vee\quad }{\omega _j^h}
A_\alpha ^j).
\end{array}
$

So, we get
\begin{equation}
\Pi _\alpha ^{\overline{\alpha }}=A_h^{\overline{\alpha }}(dA_\alpha ^h+
\stackrel{\vee\quad }{\omega _j^hA_\alpha ^j)}  \tag{9.5.5}
\end{equation}

Remarking that
$$
dA_\alpha ^h=\displaystyle\frac{\delta A_\alpha ^h}{\delta u^\beta }
du^\beta +\displaystyle\frac{\delta A_\alpha ^h}{\delta v^{(1)\beta }}\delta
v^{(1)\beta }+...+\displaystyle\frac{\delta A_\alpha ^h}{\delta
v^{(k-1)\beta }}\delta v^{(k-1)\beta }+\stackrel{\cdot }{\partial ^\beta }
A_\alpha ^h\delta \stackrel{\vee}{p}_\beta
$$
and $\stackrel{\vee}{\omega _j^i}$ is from (9.4.4) we obtain
$$
\begin{array}{ll} \Pi _\alpha ^{\overline{\alpha }}= &
A_h^{\overline{\alpha }}[(\displaystyle \frac{\delta A_\alpha
^h}{\delta u^\beta }+A_\alpha ^j\stackrel{\vee\quad }{ H_{j\beta
}^h})du^\beta +(\displaystyle\frac{\delta A_\alpha ^h}{\delta
v^{(1)\beta }}+A_\alpha ^j \stackrel{\vee\quad }{\underset{(1)}{
C_{j\beta }^h}})\delta v^{(1)\beta }+... \\
& ...+(\displaystyle\frac{\delta A_\alpha ^h}{\delta v^{(k-1)\beta
}} +A_\alpha ^j \stackrel{\vee\quad
}{\underset{(k-1)}{C_{j\beta }^h}} )\delta v^{(k-1)\beta
}+\stackrel{\cdot }{(\partial ^\beta }A_\alpha ^h+A_\alpha
^j\stackrel{\vee\quad }{\underset{(1)}{C_j^{h\beta }}})\delta
\stackrel{\vee}{p}_\beta .
\end{array}
$$

The last expression of the $1$-forms $\prod_\alpha ^{\overline{\alpha }}$ can
be written as in (9.5.3) with the coefficients (9.5.4).

In order to prove the second formula (9.5.2) we remark, using the same method,
the formula
\begin{equation}
\bigtriangledown A_{\overline{\alpha }}^i=A_\alpha ^i[A_j^\alpha (dA_{%
\overline{\alpha }}^j+A_{\overline{\alpha }}^r\stackrel{\vee}{\omega _r^j}
)]=-A_\alpha ^i\Pi _{\overline{\alpha }}^\alpha . \tag{9.5.6}
\end{equation}

So, we have
\begin{equation}
\Pi _{\overline{\alpha }}^\alpha =-A_j^\alpha (dA_{\overline{\alpha }}^j+A_{%
\overline{\alpha }}^r\stackrel{\vee}{\omega _r^j}).  \tag{9.5.6a}
\end{equation}

The second relation between the $1$-forms $\Pi _\alpha ^{\overline{\alpha }}$
and $\Pi _{\overline{\alpha }}^\alpha $ is proved without difficulties by means
of the equations $\bigtriangledown g_{ij}=0,\bigtriangledown \stackrel{\vee
}{g}_{\alpha \beta }=0$ and $\bigtriangledown \delta _{\overline{\alpha }
\overline{\beta }}$ $=0,$ which will be proved in the next Lemma. q.e.d

\begin{lem}
The following properties of the relative covariant derivation $%
\bigtriangledown $ hold:
\begin{equation}
\bigtriangledown g_{ij}=0, \ \bigtriangledown g_{\alpha \beta
}=0, \ \bigtriangledown \delta _{\overline{\alpha }\overline{\beta }}=0.
\tag{9.5.7}
\end{equation}
\end{lem}

Evidently, the equation $\bigtriangledown g_{ij}=0$ is immediate. After the
formula (9.2.8') we have
\begin{equation}
\stackrel{\vee}{g}_{\alpha \beta }=g_{ij}A_\alpha ^iA_\beta ^j, \ \delta _{%
\overline{\alpha }\overline{\beta }}=g_{ij}A_{\overline{\alpha }}^iA_{%
\overline{\beta }}^j  \tag{9.5.8}
\end{equation}

So, using (9.5.2) and $\bigtriangledown g_{ij}=0,$ we obtain $\bigtriangledown
\stackrel{\vee}{g}_{\alpha \beta }=0,\bigtriangledown \delta _{\overline{
\alpha }\overline{\beta }}=0$. q.e.d.

\begin{rem}
1$^{\circ }$ Theorem 8.4.5 is a consequences of the previous Lemma.

2$^{\circ }$ Because of $g_{ij}A_\alpha ^iA_{\overline{\beta }}^j=0,$ applying
the operator $\bigtriangledown $ we get by means of (9.5.2) that
$g_{ij}(\Pi _\alpha ^{\overline{\alpha }}A_{\overline{\alpha }}^iA_{%
\overline{\beta }}^j-A_\alpha ^iA_\beta ^j\Pi _{\overline{\alpha }}^\beta
)=0.$ But, this is $\delta _{\overline{\alpha }\overline{\beta }}\Pi _\alpha
^{\overline{\alpha }}=\stackrel{\vee}{g}_{\alpha \beta }\Pi _{\overline{
\alpha }}^\beta ,$ which imply the second formula (9.5.3).
\end{rem}

The formulae (9.5.4) allows to prove that the coefficients
\begin{equation}
\underset{\left( 0\right) \quad }{H_{\alpha \gamma
}^{\overline{\beta }}}, \underset{\left( 1\right) \quad
}{H_{\alpha \gamma }^{\overline{\beta }}} ,...,\underset{\left(
k-1\right) \quad }{H_{\alpha \gamma }^{\overline{ \beta
}}},H_\alpha ^{\overline{\beta }\gamma }  \tag{9.5.9}
\end{equation}
are mixed $d$-tensor. They will be called the \textit{second fundamental
tensors} of the Hamilton subspace $\stackrel{\vee\quad }{H^{\left( k\right)
m}}=(\stackrel{\vee}{M}, \stackrel{\vee}{H})$ of the Hamilton space $%
H^{\left( k\right) n}=(M,H).$

As an application of the previous considerations we get:
\begin{prop}
We have
\begin{equation}
\stackrel{\vee}{\Omega }_{ij}=-\stackrel{\vee}{\Omega }_{ji},\Omega
_{\alpha \beta }=-\Omega _{\beta \alpha },\Omega _{\overline{\alpha }%
\overline{\beta }}=-\Omega _{\overline{\beta }\overline{\alpha }}  \tag{9.5.10}
\end{equation}
\end{prop}

The proof follows the same way as for the covariant curvature $2$-forms $%
\Omega _{ij}$ of a metrical connection.

Let us consider a parametrized smooth curve $\stackrel{\vee}{c}$ on the
Hamilton subspace $\stackrel{\vee\quad }{H^{\left( k\right) m}}=(\stackrel{
\vee}{M}, \stackrel{\vee}{H}).$ Locally $\stackrel{\vee}{c}$ can be given
by
\[
u^\alpha =u^\alpha (t),v^{\left( 1\right) \alpha }=v^{\left( 1\right) \alpha
}(t),...,v^{(k-1)\alpha }=v^{\left( k-1\right) \alpha }\left( t\right) ,
\stackrel{\vee}{p}_\alpha =\stackrel{\vee}{p_\alpha }(t)
\]

A tangent vector $X^i$ along $\stackrel{\vee}{c}$ is given by $X^i=A_\alpha
^iX^\alpha $.

But, along curve $\stackrel{\vee}{c}$ $,$we have
\begin{equation}
\displaystyle\frac{\bigtriangledown X^i}{dt}=A_{\overline{\alpha }}^i\Pi
_\alpha ^{\overline{\alpha }}X^\alpha +A_\alpha ^i\displaystyle\frac{
\bigtriangledown X^\alpha }{dt}  \tag{9.5.11}
\end{equation}
$\displaystyle\frac{\bigtriangledown X^\alpha }{dt}=0$ along $\stackrel{\vee
}{c},$ implies $\displaystyle\frac{\bigtriangledown X^i}{dt}=A_{\overline{
\alpha }}^i\Pi _\alpha ^{\overline{\alpha }}X^\alpha .$ This means that $%
\displaystyle\frac{\bigtriangledown X^i}{dt}$ is normal to the subspace $%
\stackrel{\vee\quad }{H^{\left( k\right) m}}.$

We say that the subspace $\stackrel{\vee\quad }{H^{\left( k\right) n}}$ is
\textit{totally geodesic }in the space $H^{\left( k\right) n}$ if along any
curve $\stackrel{\vee}{c}$ , $\displaystyle\frac{\bigtriangledown X^\alpha
}{dt}=0$ implies $\displaystyle\frac{\bigtriangledown X^i}{dt}=0.$ The
geometrical meaning of the condition is evident.

\begin{teo}
The Hamilton subspaces of order $k,$ $\stackrel{\vee\quad }{H^{\left(
k\right) n}}=(\stackrel{\vee}{H},\!\!\stackrel{\vee}{M})$ is totally geodesic
in the Hamilton space of order $k,$ $H^{\left( k\right) n}=(H,M)$ if and
only if, the second fundamental tensors $\underset{\left( 0\right) \quad }{%
H_{\alpha \gamma }^{\overline{\beta }}},\underset{\left( 1\right) \quad }{%
H_{\alpha \gamma }^{\overline{\beta }}},...,\underset{\left(
k-1\right)
\quad }{H_{\alpha \gamma }^{\overline{\beta }}},H_\alpha ^{\overline{\beta }%
\gamma }$ vanish.
\end{teo}

\textbf{Proof}. If the second fundamental tensors vanish then $\Pi _\beta ^{%
\overline{\alpha }}$ identically vanish and by (9.5.11) it follows that $%
\displaystyle\frac{\bigtriangledown X^\alpha }{dt}=0$ imply $\displaystyle
\frac{\bigtriangledown X^i}{dt}=0,$ for any $X^i=A_\alpha ^iX^\alpha .$

Conversely, the condition $\displaystyle\frac{\bigtriangledown X^\alpha }{dt}
=0\Rightarrow \displaystyle\frac{\bigtriangledown X^i}{dt}=0$ for any $%
X^i=A_\alpha ^i$ $X^\alpha $ and the formula (9.5.11) leads to $\Pi _\alpha ^{%
\overline{\alpha }} $ $=0.$ Taking into account (9.5.3) and the fact that $%
\stackrel{\vee}{c}$ is arbitrary, it follows that all second fundamental tensors
vanish.

\section{The Gauss-Codazzi Equations}

The Gauss-Codazzi equations of a Hamilton subspaces of order $k$, $\stackrel{
\vee\quad }{H^{\left( k\right) m}}=(\stackrel{\vee}{M}, \stackrel{\vee}{H}
)\,$ in the Hamilton spaces of order $k$, $H^{\left( k\right) n}=(H,M)$
endowed with a canonical metrical $N$-connection $D$, are obtained from the integrability
conditions of the system of equalities (9.5.2). We can deduce these equations
using the structure equations (9.4.23) and (9.4.24).

\begin{teo}
The Gauss-Codazzi equations of the Hamilton subspaces $\stackrel{\vee\quad
}{H^{\left( k\right) m}}$ in the Hamilton space$H^{\left( k\right) n}$ are as
follows:
\begin{equation}
\begin{array}{l}
A_\alpha ^iA_\beta ^j\stackrel{\vee}{\Omega }_{ij}-\Omega _{\alpha \beta
}=\Pi _{\beta \overline{\gamma }}\wedge \Pi _\alpha ^{\overline{\gamma }}, \\
\\
A_{\overline{\alpha }}^iA_{\overline{\beta }}^j\stackrel{\vee}{\Omega }%
_{ij}-\Omega _{\overline{\alpha }\overline{\beta }}=\Pi _{\gamma \overline{%
\beta }}\wedge \Pi _{\overline{\alpha }}^\gamma, \\
\\
-A_\alpha ^iA_{\overline{\beta }}^\alpha \stackrel{\vee}{\Omega }%
_{ij}=\delta _{\overline{\alpha }\overline{\gamma }}(d\Pi _\alpha ^{%
\overline{\gamma }}+\Pi _\beta ^{\overline{\gamma }}\wedge \omega _\alpha
^\beta -\Pi _\alpha ^{\overline{\alpha }}\wedge \omega _{\overline{\alpha }%
}^{\overline{\gamma }}),
\end{array}
\tag{9.6.1}
\end{equation}
where $\Pi _{\alpha \overline{\beta }}=\stackrel{\vee}{g}_{\alpha \gamma
}\Pi _{\overline{\beta }}^\gamma .$
\end{teo}

\textbf{Proof}. Consider the first equations (9.5.2), written in the form
\begin{equation}
dA_\alpha ^i+A_\alpha ^j\stackrel{\vee}{\omega _j^i}-A_\beta ^i\omega
_\alpha ^\beta =A_{\overline{\alpha }}^i\Pi _\alpha ^{\overline{\alpha }}.
\tag{*}
\end{equation}
By exterior differentiation of the both side of this equation we get:
\begin{equation}
dA_\alpha ^j\wedge \stackrel{\vee}{\omega _j^i}+A_\alpha ^jd\stackrel{\vee
}{\omega _j^i}-dA_\beta ^i\wedge \omega _\alpha ^\beta -A_\beta ^id\omega
_\alpha ^\beta =dA_{\overline{\alpha }}^i\wedge \Pi _\alpha ^{\overline{
\alpha }}+A_{\overline{\alpha }}^id\Pi _\alpha ^{\overline{\alpha }}.
\tag{**}
\end{equation}
Looking to the equality (*) and at:
\begin{equation}
dA_{\overline{\alpha }}^i+A_{\overline{\alpha }}^j\stackrel{\vee}{\omega }
_j^i-A_{\overline{\beta }}^i\omega _{\overline{\alpha }}^{\overline{\beta }
}=-A_\beta ^i\Pi _{\overline{\alpha }}^\beta ,  \tag{***}
\end{equation}
the equality (**) becomes
\begin{equation}
-A_{\overline{\alpha }}^j\stackrel{\vee}{\Omega }_j^i+A_\beta ^i(\Omega
_\alpha ^\beta +\Pi _{\overline{\alpha }}^\beta \wedge \Pi _\alpha ^{%
\overline{\alpha }})=A_{\overline{\alpha }}^i(d\Pi _\alpha ^{\overline{
\alpha }}+\Pi _\beta ^{\overline{\alpha }}\wedge \omega _\alpha ^\beta -\Pi
_\alpha ^{\overline{\beta }}\wedge \omega _{\overline{\beta }}^{\overline{
\alpha }}).  \tag{9.6.2}
\end{equation}

Now, multiplying (9.6.2) by $g_{ih}A_\beta ^h=\stackrel{\vee}{g}_{\beta
\gamma }A_i^\gamma $ we obtain the first equations (9.6.1). The same
operation, taking the factor $g_{ij}A_{\overline{\beta }}^j=\delta _{%
\overline{\beta }\overline{\alpha }}A_i^{\overline{\alpha }}$ leads to third
equation (9.6.1). Of course, to deduce the second equations (9.6.1) we will take
the exterior differential of both sides of equations (***) and apply the
same method. q.e.d.

In order to obtain the system of all fundamental equations of the subspaces $%
\stackrel{\vee\quad }{H^{\left( k\right) m}}$in $H^{\left( k\right) n}$, we
must find the relations between the torsion $2$-forms $\stackrel{\left(
0\right) }{\Omega ^i},\stackrel{\left( a\right) }{\Omega ^i}
,(a=1,...,k-1),\Omega _i$ of the canonical metrical $N$-connection $D$ of $%
H^{\left( k\right) n}$ and the torsion $2$-forms $\stackrel{\left( 0\right)
}{\Omega ^\alpha },\stackrel{\left( a\right) }{\Omega ^\alpha }
,(a=1,...,k-1),\Omega _\alpha $ of the relative connection $\nabla .$

Therefore, we obtain:
\begin{teo}
The fundamental equation of the Hamilton subspaces of order $k,$ $\stackrel{%
\vee\quad }{H^{\left( k\right) m}}=(\stackrel{\vee}{M}, \stackrel{\vee}{H})
$ in a Hamilton space of order $k,$ $H^{\left( k\right) n}=(M,H)$ endowed
with the canonical metrical $N$-connection $D$ are: the Gauss-Codazzi
equations (9.6.1), as well as the following equations:

$d(dx^i)-dx^j\wedge \omega _j^i=\stackrel{\left( 0\right) }{\Omega ^i},$
modulo (9.3.8),

$d(\delta y^{\left( a\right) i})-\delta y^{\left( a\right) j}\wedge \omega
_j^i=-\stackrel{\left( a\right) }{\Omega ^i}$ modulo (9.3.8), $(a=1,...,k-1)$, \vspace{3mm}

$d(\delta p_i)+\delta p_j\wedge \omega _i^j=-\Omega _i,$ modulo (9.3.8').
\end{teo}

We end here the theory of Hamilton subspaces of order $k$ in a Hamilton
space of the same order.

We underline the importance of this theory for applications in the Higher Order
Hamiltonian Mechanics.

The particular case $m=n-1$ of the hypersurfaces $\!\!\stackrel{\vee}{M}$ in $M$
can be obtained from the previous study without difficulties. If $k=1$ we
have the theory of Hamilton subspaces $\stackrel{\vee}{H^n}=(\stackrel{\vee}{M}, \stackrel{\vee}{
H}(x,\stackrel{\vee}{p})$ in a Hamilton space $H^n=(M,H(x,p)).$

\chapter{The Cartan Spaces of Order $k$ as Dual of Finsler Spaces of Order $k$}

\markboth{\it{THE GEOMETRY OF HIGHER-ORDER HAMILTON SPACES\ \ \ \ \ }}{\it{The Cartan Spaces of Order} $k$}

The geometry of Finsler spaces of order $k$, $%
F^{(k)n}=(M,F(x,y^{(1)},...,y^{(k)}))$, studied in chapter 3, is a
particular case of the geometry of Lagrange spaces of order $k$, $%
L^{(k)n}=(M,L(x,y^{(1)},...,y^{(k)}))$. The fundamental function $F$ being a
$k$-homogeneous regular Lagrangian on the fibres of the bundle $T^kM$.

As we know, the 'dual' of the space $L^{(k)n}$, via Legendre
transformation $Leg$, is a Hamilton space of order $k$, $%
H^{(k)n}=(M,H(x,y^{(1)},...,y^{(k-1)},p))$. In this chapter we will prove
that the restriction of the mapping $Leg$ to the Finsler spaces of order $k$
determines a new class of Hamilton spaces of order $k$, that are called the Cartan
spaces of order $k$ and are denoted by $\mathcal{C}^{(k)n}$ $= (M, K(x,y^{(1)}, ..., y^{(k-1)}, p))$.

Also, the spaces $\mathcal{C}^{(k)n}$ are the Hamilton spaces $%
H^{(k)n}=(M,H) $ in which the fundamental function $K(x,y^{(1)},...,y^{(k-1)},p)$ is regular and $k$-homogeneous on the fibres
of the dual bundle $T^{*k}M$.

For the spaces $\mathcal{C}^{(k)n}$ it is important to determine the
fundamental geometrical object fields which are important in their differential geometry.

\section{$\mathcal{C}^{(k)n}$-Spaces}

\begin{defi}
A Cartan space of order $k\geq 1$ is a pair $\mathcal{C}^{(k)n} =$ \newline
$= (M, K(x,y^{(1)},...,y^{(k-1)},p))$ for which the following axioms
hold:

1$^0$ $K$ is a real function on the manifold $T^{*k}M$, differentiable on $%
\widetilde{T^{*k}M}$ and continuous on the null section of the projection $\pi^{*k}:T^{*k}M\rightarrow M$.

2$^0$ $K>0$ on $T^{*k}M$.

3$^0$ $K$ is positively $k$-homogeneous on the fibres of bundle $T^{*k}M$,
i.e.
\begin{equation}
K(x,ay^{(1)},...,a^{k-1}y^{(k-1)},a^kp)=a^kK(x,y^{(1)},...,y^{(k-1)},p),\
\forall a\in \mathbf{R}^{+}.  \tag{10.1.1}
\end{equation}

4$^0$ The Hessian of $K^2$, with respect to the momenta $p_i$, having the
elements
\begin{equation}
g^{ij}(x,y^{(1)},...,y^{(k-1)},p)=\displaystyle\frac 12\stackrel{\cdot }{%
\partial }^i\stackrel{\cdot }{\partial }^jK^2  \tag{10.1.2}
\end{equation}
is positively defined.
\end{defi}

It follows that the $k$-homogeneity of the function $K$ is not $k$
-homogeneity of $K$ with respect to $p_i$. It is considered as in section 1,
ch. 3. In the book [115], ch. 13 is considered only the case of homogeneity of
$K$ with respect to $p_i$.

We deduce from (10.1.2) that $g^{ij}$ is contravariant of order two, symmetric
and nondegenerate, i.e.

\begin{equation}
rank\left\| g^{ij}\right\| =n,\quad on\quad \widetilde{T^{*k}M}.  \tag{10.1.3}
\end{equation}

The function $K$ is called the \textit{\ fundamental} (or \textit{metric})
function of $\mathcal{C}^{(k)n}$ and $g^{ij}$ is \textit{the fundamental
tensor} of this space.

The first problem is if exist the Cartan spaces $\mathcal{C}^{(k)n}$.

\begin{teo}
If the base manifold $M$ is paracompact, then on $T^{*k}M$ there exist
functions $K(x,y^{(1)},...,y^{(k-1)},p)$ such that the pair $(M,K)$ is a
Cartan space of order $k$.
\end{teo}

\textit{Proof}: The manifold $M$ being paracompact there exists a Finsler
space of order $k-1$, $F^{(k-1)n}=(M,F(x,y^{(1)},...,y^{(k-1)}))$. Its
fundamental tensor $a_{ij}$:
\[
a_{ij}(x,y^{(1)},...,y^{(k-1)})=\displaystyle\frac 12\displaystyle\frac{
\partial ^2F^2}{\partial y^{(k-1)i}\partial y^{(k-1)j}}
\]
is positively defined on $T^{k-1}M$ and $0$-homogeneous on the fibres of
bundle $(T^{k-1}M,\pi ^{k-1},M)$.

Now, we can construct on the dual bundle $T^{*k}M=T^{k-1}M\times _MT^{*}M$
the following Hamiltonian

\begin{equation}
K(x,y^{(1)},...,y^{(k-1)},p)=\left\{
a^{ij}(x,y^{(1)},...,y^{(k-1)})p_ip_j\right\} ^{1/2},  \tag{10.1.4}
\end{equation}
where $a^{ij}$ is the contravariant tensor of $a_{ij}$.

It follows, without difficulties, that $K$ is a scalar function (i.e. it
does not depend of the transformations of coordinates on $T^{*k}M$) which
satisfies the axioms 1$^0$-4$^0$ from Definition 1.1.

The fact that $K$ is positively $k$-homogeneous on the fibres of the bundle $T^{*k}M$ follows directly. As $a^{ij}$ are $0$-homogeneous we have that
$$
K(x, ay^{(1)}, ..., a^{k-1}y^{(k-1)}, a^kp) = a^k K(x, y^{(1)}, ..., y^{(k-1)}, p).
$$
The fundamental tensor $g^{ij}(x,y^{(1)},...,y^{(k-1)},p)$ coincides with \newline $a^{ij}(x,y^{(1)},...,y^{(k-1)})$. The conclusion is that the pair $(M,K)$, with $K$ from (10.1.4) is a Cartan space of order $k$.

\begin{rem}
{\rm The Cartan space $\mathcal{C}^{(k)n}$ with the fundamental function $K$ from
(10.1.4) is a particular one, but it is an important and useful example of Cartan space of
order $k$.}
\end{rem}

\begin{rem}
{\rm If $\mathcal{C}^{(k)n}=(M,K)$ is a Cartan space of order $k$, then $%
H^{(k)n}=(M,K^2)$ is a Hamilton space of order $k$ having the same
fundamental tensor $g^{ij}$ as the space $\mathcal{C}^{(k)n}$. $H^{(k)n}$
will be called the Hamilton space associated to the Cartan space of order $k$
, $\mathcal{C}^{(k)n}$. All geometrical properties of $H^{(k)n}=(M,K^2)$ are
geometrical properties of $\mathcal{C}^{(k)n}=(M,K)$.}
\end{rem}

\section{Geometrical Properties of the Cartan Spaces of Order $k$}

First of all we shall study those properties of the spaces $\mathcal{C}
^{(k)n}=(M,K)$ which result from the $k$-homogeneity of the fundamental
function $K$ expressed in the identity (10.1.1).

Consider the vector field on $T^{*k}M$:
\begin{equation}
\stackrel{k-1}{\Gamma }+kC^{*}=y^{(1)i}\displaystyle\frac \partial {\partial
y^{(1)i}}+\cdots +(k-1)y^{(k-1)i}\displaystyle\frac \partial {\partial
y^{(k-1)i}}+kp_i\stackrel{\cdot }{\partial }^i  \tag{10.2.1}
\end{equation}
and the Lie derivation with respect to this vector field $\mathcal{L}_{%
\stackrel{k-1}{\Gamma }+kC^{*}}$ .

The functions $K$, $K^2$, $\stackrel{\cdot }{\partial }^jK^{2k}$, $%
g^{ij}$ ($i,j=1,...,n$) and $C^{ijh}$:
\begin{equation}
C^{ijh}=\displaystyle\frac 12\stackrel{\cdot }{\partial }^hg^{ij}=- %
\displaystyle\frac 14\stackrel{\cdot }{\partial }^i\stackrel{\cdot }{
\partial }^j\stackrel{\cdot }{\partial }^hK^2  \tag{10.2.2}
\end{equation}
are $k$, $2k$, $k$, $0$ and $-k$ homogeneous, respectively on the fibres of
the bundle $T^{*k}M$. Taking into account the Theorems 4.5.7 and 4.5.8 we have:

\begin{prop}
The following identities hold:
\begin{equation}
\mathcal{L}_{\stackrel{k-1}{\Gamma }+kC^{*}}K=kK,  \tag{10.2.3}
\end{equation}
\begin{equation}
\mathcal{L}_{\stackrel{k-1}{\Gamma }+kC^{*}}K^2=2kK,  \tag{10.2.4}
\end{equation}
\begin{equation}
\mathcal{L}_{\stackrel{k-1}{\Gamma }+kC^{*}}\stackrel{\cdot }{\partial }%
^hK^2=k\stackrel{\cdot }{\partial }^hK^2,  \tag{10.2.5}
\end{equation}
\begin{equation}
\mathcal{L}_{\stackrel{k-1}{\Gamma }+kC^{*}}g^{ij}=0,  \tag{10.2.6}
\end{equation}
\begin{equation}
\mathcal{L}_{\stackrel{k-1}{\Gamma }+kC^{*}}C^{ijh}=-kC^{ijh}.  \tag{10.2.7}
\end{equation}
\end{prop}
We remark, especially, $0$-homogeneity of the functions $
g^{ij}(x,y^{(1)},...,y^{(k-1)},$ \newline $p)$. By virtue of (10.2.1) and (10.2.6) we have the
identity:
\begin{equation}
y^{(1)i}\displaystyle\frac{\partial g^{jh}}{\partial y^{(1)i}}+\cdots
+(k-1)y^{(k-1)i}\displaystyle\frac{\partial g^{jh}}{\partial y^{(k-1)i}}%
+kp_i \stackrel{\cdot }{\partial }^ig^{jh}=0.  \tag{10.2.6a}
\end{equation}

Of course, the covariant tensor field $g_{ij}$ of the fundamental tensor $%
g^{ij}$ of the Cartan space $\mathcal{C}^{(k)n}$ is $0$-homogeneous on the
fibres of $T^{*k}M$. Indeed, the equality (10.2.6) implies $\mathcal{L}_{%
\stackrel{k-1}{\Gamma }+kC^{*}}g_{ij}=0$. Therefore we have the identity:
\begin{equation}
y^{(1)i}\displaystyle\frac{\partial g_{jh}}{\partial y^{(1)i}}+\cdots
+(k-1)y^{(k-1)i}\displaystyle\frac{\partial g_{jh}}{\partial y^{(k-1)i}}%
+kp_i \stackrel{\cdot }{\partial }^ig_{jh}=0.  \tag{10.2.6b}
\end{equation}

In any Cartan space of order $k$, $\mathcal{C}^{(k)n}=(M,K)$,
there exists two important tensors $\underset{(k-1)}{C}$
$_{jh}^i$, $C_i^{jh}$ which are the $v_{k-1}$- and $w_k$-coefficients
of a canonical metrical connection.

We have:
\begin{teo}
For any Cartan space of order $k$, $\mathcal{C}^{(k)n}=(M,K)$, $\underset{(k-1)}{C}\text{}_{jh}^i$, $C_i^{jh}$ given by
\begin{equation}
\underset{(k-1)}{C}\text{ }_{jh}^i=\displaystyle\frac 12g^{is}\left( %
\displaystyle\frac{\partial g_{sh}}{\partial y^{(k-1)j}}+\displaystyle\frac{%
\partial g_{js}}{\partial y^{(k-1)h}}-\displaystyle\frac{\partial g_{jh}}{%
\partial y^{(k-1)s}}\right) ,  \tag{10.2.8}
\end{equation}
\begin{equation}
C_i^{jh}=-\displaystyle\frac 12g_{is}\left( \displaystyle\frac{\partial
g^{sh}}{\partial p_j}+\displaystyle\frac{\partial g^{js}}{\partial p_h}-%
\displaystyle\frac{\partial g^{jh}}{\partial p_s}\right)   \tag{10.2.9}
\end{equation}
have the following properties:

1$^0$ They are $d$-tensor fields of type (1,2) and (2,1), respectively.

2$^0$ $\underset{(k-1)}{C}$ $_{jh}^i$ are $1-k$ homogeneous and
$C_i^{jh}$ are $-k$ homogeneous on the fibres of $T^{*k}M$.

3$^0$ They are the $v_{k-1}$- and $w_k$-coefficients of a canonical metrical
connection $D$, i.e.:
\begin{equation}
g^{ij}\stackrel{(k-1)}{|}_h=0,\ g^{ij}|^h=0.  \tag{10.2.10}
\end{equation}

4$^0$ We have
\begin{equation}
C_i^{jh}=g_{is}C^{sjh}.  \tag{10.2.11}
\end{equation}

5$^0$ $\underset{(k-1)}{S}$ $_{jh}^i=\underset{(k-1)}{C}$ $_{jh}^i-%
\underset{(k-1)}{C}$ $_{hj}^i=0$, $S_i^{jh}=C_i^{jh}-C_i^{hj}=0.$
\end{teo}

The proofs of these affirmations are not difficult.

\section{Canonical Presymplectic Structures, Variational Problem of the Space $\mathcal{C}^{(k)n}$}

Consider a Cartan space of order $k$, $\mathcal{C}
^{(k)n}=(M,K(x,y^{(1)},...,y^{(k-1)},p))$, and the canonical presymplectic
structure
\begin{equation}
\theta =dp_i\wedge dx^i,  \tag{10.3.1}
\end{equation}
where
\begin{equation}
\theta =d\omega ,\omega =p_idx^i.  \tag{10.3.1a}
\end{equation}
Of course the structure $\theta$ can be directly studied using the property of integrability, given by $d\theta=0$. The relations between the structure $\theta $ and the Poisson structure $\{f, g\}_0:$
\begin{equation}
\left\{ f,g\right\} _0=\displaystyle\frac{\partial f}{\partial x^i} %
\displaystyle\frac{\partial g}{\partial p_i}-\displaystyle\frac{\partial f}{
\partial p_i}\displaystyle\frac{\partial g}{\partial x^i},  \forall f, g, \in {\cal F}(\Sigma_0), \tag{10.3.2}
\end{equation}
and
\begin{equation}
\Sigma _0=\left\{ (x,y^{(1)},...,y^{(k-1)},p)|y^{(1)i}=\cdots
=y^{(k-1)i}=0\right\},  \tag{10.3.3}
\end{equation}
lead to some important geometrical results.

Recall that the submanifold $\Sigma _0$ has been introduced in the section 3 of chapter
8 as a section of the canonical projection of the differentiable bundle $%
(T^{*k}M,\overline{\pi }^{*},T^{*}M)$.

Let us consider the restriction $K_0$ of the fundamental function $K$ of
space $\mathcal{C}^{(k)n}$ to $\Sigma_0$:

\begin{equation}
K_0=K_{|\Sigma_0}.  \tag{10.3.4}
\end{equation}

Thus we have
\[
K_0(x,p)=K(x,0,...,0,p).
\]

It follows that $K_0$ is $1$-homogeneous with respect to $p_i.$

Therefore the pair $\mathcal{C}_0^{(1)n}=(M,K_0)$ is a classical Hamilton
space, [115], with the fundamental tensor field

\begin{equation}
g_0^{ij}(x,p)=g^{ij}(x,0,...,0,p)=\displaystyle\frac 12\stackrel{\cdot }{
\partial }^i\stackrel{\cdot }{\partial }^jK_0^2.  \tag{10.3.5}
\end{equation}

Theorem 7.3.1 is valid in the particular case of Cartan space $\mathcal{C}
_0^{(1)n}$:

\begin{teo}
The pair $(\Sigma _0,\theta _0)$, with $\theta _0=\theta _{|\Sigma _0}$, is
a symplectic manifold.
\end{teo}

Since $\theta _0=dp_i\wedge dx^i$ in every point $(x,p)\in \Sigma _0$, it is
a closed $2$-form of rank $2n=\dim \Sigma _0$.

The tangent space $T_u\Sigma _0$ at a point $u\in \Sigma _0$ has the natural
basis $\left( \displaystyle\frac \partial {\partial x^i},\displaystyle\frac
\partial {\partial p_i}\right) _u$, ($i=1,...,n$) and the natural cobasis $%
\left( dx^i,dp_i\right) $, ($i=1,...,n$).

Consider the $\mathcal{F}(\Sigma _0)$-module $\mathcal{X}(\Sigma _0)$ of vector
fields and the $\mathcal{F}(\Sigma _0)$-module $\mathcal{X}^{*}(\Sigma _0)$ of $%
1 $-form fields on the submanifold $\Sigma _0$.

Taking into account the theory from ch. 8, section 3 we have:

1$^0$ The $\mathcal{F}(\Sigma _0)$-linear mapping $S_{\theta _0}:\mathcal{X}
(\Sigma _0)\rightarrow \mathcal{X}^{*}(\Sigma _0)$ given by
\[
S_{\theta _0}(X)=i_X\theta _0,\ \forall X\in \mathcal{X}(\Sigma _0)
\]
is an isomorphism.

2$^0$ There exists an unique vector field $X_{K_0^2}\in \mathcal{X}(\Sigma
_0)$ such that
\[
S_{\theta _0}(X_{K_0^2})=i_{X_{K_0^2}}\theta _0=-dK_0^2.
\]

3$^0$ The Hamiltonian vector field $X_{K_0^2}$ is given by

\[
X_{K_0^2}=\displaystyle\frac{\partial K_0^2}{\partial p_i}\displaystyle\frac
\partial {\partial x^i}-\displaystyle\frac{\partial K_0^2}{\partial x^i} %
\displaystyle\frac \partial {\partial p_i}.
\]
Consequently, we can formulate:
\begin{teo}
The integral curves of the Hamiltonian vector field $X_{K_0^2}$ are given by
the $\Sigma _0$-canonical equations:
\begin{equation}
\displaystyle\frac{dx^i}{dt}=\displaystyle\frac{\partial K_0^2}{\partial p_i}%
,\ \displaystyle\frac{dp_i}{dt}=-\displaystyle\frac{\partial K_0^2}{\partial
x^i},\ y^{(\alpha )i}=0,\ (\alpha =1,...,k-1).  \tag{10.3.6}
\end{equation}
\end{teo}

Now, let us consider two functions $f,g\in \mathcal{F}(\Sigma _0)$ and the
vector fields $X_f$, $X_g$ given by
\[
i_{X_f}\theta _0=-df,\ i_{X_g}\theta _0=-dg.
\]
\begin{prop}
The following relations between the structures $\theta _0$ and $\left\{
,\right\} _0$ on the submanifold $\Sigma _0$, hold:
\begin{equation}
\left\{ f,g\right\} _0=\theta _0(X_f,X_g),\ \forall f,g\in \mathcal{F}%
(\Sigma _0).  \tag{10.3.7}
\end{equation}
\end{prop}

\begin{rem}
{\rm The triple $(T^{*k}M,K^2(x,y^{(1)},...,y^{(k-1)},p),\theta )$ is a
particular Hamiltonian system. It can be studied directly applying the
method of Gotay, [115].}
\end{rem}
But, the equations (10.3.6) are extremely particular. For Cartan spaces ${\cal C}^{(k)n}=(M, K(x, y^{(1)}, ..., y^{(k-1)}, p))$ the integral of action (see Ch.5)
\begin{equation}
I(c)=\int^1_0 [p_i\frac{dx^i}{dt} - K(x, \frac{dx}{dt}, ..., \frac{1}{(k-1)!}\frac{d^{k-1}}{dt^{k-1}}, p)]dt
\tag{10.3.8}
\end{equation}
leads to the fundamental equations of the space:
\begin{equation}
\begin{array}{l}
\displaystyle\frac{dx^i}{dt} = \displaystyle\frac{1}{2} \displaystyle\frac{\partial K}{\partial p_i}, \vspace{3mm}\\
\displaystyle\frac{dp_i}{dt} = -\displaystyle\frac{1}{2}[\displaystyle\frac{\partial K}{\partial x^i} - \displaystyle\frac{d}{dt} \displaystyle\frac{\partial K}{\partial y^{(1)i}} + \cdots + (-1)^{k-1} \displaystyle\frac{d^{k-1}}{dt^{k-1}}\displaystyle\frac{\partial K}{\partial y^{(k-1)i}}]
\end{array}
\tag{10.3.9}
\end{equation}

The integrand of the integral of action is $k$-homogeneous on the fibres of $T^{*k}M$. The Hamilton-Jacobi equations (10.3.9) are homogeneous on the fibres of $T^{*k}M$, too.

The energy of order $k-1$, ${\cal E}^{k-1}$, of the Cartan space ${\cal C}^{(k)n}= (M, K)$, by means of formula (5.3.1), is given by:
\begin{equation}
{\cal E}^{k-1}(K) = I^{k-1}(K) - \frac{1}{2!}\frac{d}{dt} I^{k-2}(K) + \cdots + (-1)^{k-2} \frac{1}{(k-1)!}\frac{d^{k-2}}{dt^{k-2}}I^1(K) - K.
\tag{10.3.10}
\end{equation}
We have:

\begin{teo}
For a Cartan space ${\cal C}^{(k)n}= (M, K)$, the energy of order $k-1$, ${\cal E}^{k-1}(K)$ is constant along every solution curve of the Hamilton-Jacobi equations.
\end{teo}

\section{The Cartan Spaces $\mathcal{C}^{(k)n}$ as Dual of Finsler Spaces $%
F^{(k)n}$}

Let $F^{(k)n}=(M,F(x,y^{(1)},...,y^{(k)}))$ be a Finsler space of order $k$
having \newline $F(x,y^{(1)},...,y^{(k)})$ as fundamental function and

\begin{equation}
a_{ij}(x,y^{(1)},...,y^{(k)})=\displaystyle\frac 12\displaystyle\frac{
\partial ^2F^2}{\partial y^{(k)i}\partial y^{(k)j}}  \tag{10.4.1}
\end{equation}
as fundamental tensor.

Remembering the Definition 3.1.1, $F$ is a function from $T^kM$ to $%
\mathbf{R}$, differentiable on the manifold $\widetilde{T^kM}=T^kM\setminus
\{\mathbf{0}\}$ and continuous on the null section. $F$ is a positive
function, $k$-homogeneous:
\[
F(x,ay^{(1)},...,a^ky^{(k)})=a^kF(x,y^{(1)},...,y^{(k)}),\ \forall a\in
\mathbf{R}^{+}
\]
and the tensor field $a_{ij}$ is positively defined on $\widetilde{T^kM}$.

The Legendre mapping $\varphi :\widetilde{T^kM}\rightarrow \widetilde{%
T^{*k}M }$ defined in section 4, ch.8, by
\begin{equation}
\left\{
\begin{array}{l}
x^i=y^{(0)i},y^{(1)i}=y^{(1)i},...,y^{(k-1)i}=y^{(k-1)i}, \\
p_i=\displaystyle\frac 12\displaystyle\frac{\partial F^2}{\partial y^{(k)i}}
\end{array}
\right.  \tag{10.4.2}
\end{equation}
is a local diffeomorphism.

We denote
\begin{equation}
p_i=\varphi _i(y^{(0)},y^{(1)},...,y^{(k)}).  \tag{10.4.2a}
\end{equation}
Evidently $y^{(0)i}=x^i$ is the point of the base manifold $M$, $x=\pi
^k(x,y^{(1)},...,y^{(k)})$.

The local inverse $\xi =\varphi ^{-1}:\widetilde{T^{*k}M}\rightarrow
\widetilde{T^kM}$ is expressed by (8.4.3):
\begin{equation}
\left\{
\begin{array}{l}
y^{(0)i}=x^i,y^{(1)i}=y^{(1)i},...,y^{(k-1)i}=y^{(k-1)i}, \\
y^{(k)i}=\xi ^i(x,y^{(1)},...,y^{(k-1)},p)
\end{array}
\right. \tag{10.4.3}
\end{equation}

Remarking that the Legendre transformation $\varphi $, from (10.4.2) is \newline
$k$-homogeneous, it follows that its local inverse $\xi ^{-1}$ is $k$
-homogeneous on the fibres of bundle $T^{*k}M$. So, we have:
\begin{equation}
\xi ^i(x,ay^{(1)},...,a^{k-1}y^{(k-1)},a^kp)=a^k\xi
^i(x,y^{(1)},...,y^{(k-1)},p),\ \forall a\in \mathbf{R}^{+}.  \tag{10.4.4}
\end{equation}

Applying the Theorem 8.4.1 we have:

\begin{teo}
The Legendre mapping $\varphi $ transforms the canonical $k$-spray of the
Finsler space $F^{(k)n}$:

\begin{equation}
S=y^{(1)i}\displaystyle\frac \partial {\partial x^i}+\cdots +ky^{(k)i}%
\displaystyle\frac \partial {\partial y^{(k-1)i}}-(k+1)G^i\displaystyle\frac
\partial {\partial y^{(k)i}}  \tag{10.4.5}
\end{equation}
with the coefficients
\[
(k+1)G^i(x,y^{(1)},...,y^{(k)})=\displaystyle\frac 12a^{ij}\left\{ \Gamma
\left( \displaystyle\frac{\partial F^2}{\partial y^{(k)j}}\right) -%
\displaystyle\frac{\partial F^2}{\partial y^{(k-1)j}}\right\}
\]
in the $\xi $-dual $k$-spray
\begin{equation}
\begin{array}{l}
S_\xi ^{*}=y^{(1)i}\displaystyle\frac \partial {\partial x^i}+\cdots
+(k-1)y^{(k-1)i}\displaystyle\frac \partial {\partial y^{(k-2)i}}+ \\
\quad +k\xi ^i(x,y^{(1)},...,y^{(k-1)},p)\displaystyle\frac \partial
{\partial y^{(k-1)i}}+\eta _i(x,y^{(1)},...,y^{(k-1)},p)\displaystyle\frac
\partial {\partial p_i}
\end{array}
\tag{10.4.6}
\end{equation}
with $\xi ^i$ from (10.4.3) and
\begin{equation}
\begin{array}{l}
\eta _i= -a_{is}\big[ \displaystyle\frac{\partial \xi ^s}{\partial x^r}%
y^{(1)r}+\cdots +\\
+(k-1)\displaystyle\frac{\partial \xi ^s}{\partial y^{(k-1)r}%
}\xi ^r+(k+1)G^s(x,y^{(1)},...,y^{(k-1)},\xi )\big]
\tag{10.4.7}
\end{array}
\end{equation}
\end{teo}

\textbf{Proof}: By means of theorem 8.4.1, the $k$-spray $S$, which is $2$
-homogeneous on the fibres of $T^kM$, is transformed in the dual $k$
-semispray $S_\xi ^{*}$ from (10.4.6), (10.4.7). We must prove that $S_\xi ^{*}$
is $2$-homogeneous on the fibres of $T^{*k}M$. But every term of $S_\xi ^{*}$
is $2$-homogeneous on the fibres of $T^{*k}M$. Thus, $S_\xi ^{*}$ is a dual $%
k$-spray.

\begin{teo}
The dual $k$-spray $S_\xi ^{*}$ from (10.4.6) determines a local nonlinear
connection $N^{*}$ on the manifold $\widetilde{T^{*k}M}$, which depends on
the Finsler space of order $k$, $F^{(k)n}$, and $N^{*}$ has the following
dual coefficients
\begin{equation}
\begin{array}{l}
\underset{(1)}{M}^{*}\text{ }_j^i=-\displaystyle\frac{\partial \xi ^i}{%
\partial y^{(1)j}},\ ...,\ \underset{(k-1)}{M}^{*}\text{ }_j^i=-%
\displaystyle\frac{\partial \xi ^i}{\partial y^{(k-1)j}}, \\
N_{ij}^{*}=\displaystyle\frac{\delta \eta _i}{\delta y^{(1)j}},
\end{array}
\tag{10.4.8}
\end{equation}
where the operators $\displaystyle\frac \delta {\delta y^{(1)j}}$ are
constructed by means of the coefficients $\underset{(1)}{N}^{*}$ $%
_j^i$, ..., $\underset{(k-1)}{N}^{*}$ $_j^i$ from the dual coefficients
(10.4.8).
\end{teo}

Indeed, this theorem is just the Theorem 8.4.2 applied to Finsler spaces of
order $k$.

We remark the following property of homogeneity of the coefficients $M^{*}$,
$N^{*}$:

\begin{prop}
The coefficients $\underset{(1)}{M}^{*}$ $_j^i$, ..., $\underset{(k-1)}{M%
}^{*}$ $_j^i$, $N_{ij}^{*}$ are homogeneous on the fibres of $T^{*k}M$ of
degree $k-1$, ..., $1$, $k$, respectively.
\end{prop}

Indeed, $\xi ^i$ being $k$-homogeneous and $\eta _i$ from (10.4.7) being $k+1$
homogeneous, by means of (10.4.8) the property follows.

Now, let us consider $\stackrel{\circ }{N}$ a fixed nonlinear connection on
the submanifold $T^{k-1}M$ of $T^{*k}M=T^{k-1}M\times _MT^{*}M$. We assume
that the coefficients $\stackrel{\circ }{\underset{(1)}{M}}$ $_j^i$, ..., $%
\stackrel{\circ }{\underset{(k-1)}{M}}$ $_j^i$ are homogeneous of degree $%
1 $, ..., $k-1$ on the fibres of $T^{k-1}M$.

Thus, the Legendre mapping $\varphi :u=(x,y^{(1)},...,y^{(k-1)},y^{(k)})\in
T^kM\rightarrow u^{*}=(x,y^{(1)},...,y^{(k-1)},p)\in T^{*k}M$ transforms the
Liouville $d$-vector field $z^{(k)i}$ at the point $u$:

\begin{equation}
kz^{(k)i}=ky^{(k)i}+(k-1)\stackrel{\circ }{\underset{(1)}{M}}
_m^iy^{(k-1)m}+\cdots +\stackrel{\circ
}{\underset{(k-1)}{M}}_m^iy^{(1)m} \tag{10.4.9}
\end{equation}
in the $d$-vector field $\stackrel{\vee}{z}^{(k)i}$ at the point $u^{*}$.
The vector field $\stackrel{\vee}{z}^{(k)i}$ is:
\begin{equation}
k\stackrel{\vee}{z}^{(k)i}=k\xi ^i+(k-1)\stackrel{\circ }{\underset{(1)}{%
M }}_m^iy^{(k-1)m}+\cdots +\stackrel{\circ }{\underset{(k-1)}{M}}
_m^iy^{(1)m}.  \tag{10.4.9a}
\end{equation}
Of course, $\stackrel{\vee}{z}^{(k)i}$is $k$-homogeneous on the fibres of $%
T^{*k}M$.
Consider the function
\begin{equation}
K^2(x,y^{(1)},...,y^{(k-1)},p)=2p_i\stackrel{\vee}{z}
^{(k)i}-F^2(x,y^{(1)},...,y^{(k-1)},\xi (x,y^{(1)},...,y^{(k-1)},p)).
\tag{10.4.10}
\end{equation}
\begin{teo}
We have:

1$^0$ The pair $(M,K^2(x,y^{(1)},...,y^{(k-1)},p))$ is a Hamilton space. Its
fundamental function $K^2$ is $2k$-homogeneous on the fibres of $T^{*k}M$.

2$^0$ The fundamental tensor of the space $(M,K^2)$ is positively defined
and it does not depend on the apriori given nonlinear connection $\stackrel{%
\circ }{N}$. It is given by:
\[
g^{ij}(x,y^{(1)},...,y^{(k-1)},p)=a^{ij}(x,y^{(1)},...,y^{(k-1)},\xi
^i(x,y^{(1)},...,y^{(k-1)},p)).
\]
\end{teo}

The proof follows from Theorem 8.4.3 in which $H=K^2$.

We shall say that $(M,K^2)$ is a Cartan spaces of order $k$, dual (via
Legendre transformation) to the Lagrange space $L^{(k)n}=(M,F^2)$ associated
to the Finsler space $F^{(k)n}=(M,F)$.

The inverse problem: being given a Cartan space of order $k$, \newline
$\mathcal{C}^{(k)n}=(M,K(x,y^{(1)},...,y^{(k-1)},p))$ let us determine its dual as a
Finsler space of order $k$, $F^{(k)n}=(M,F(x,y^{(1)},...,y^{(k-1)},y^{(k)}))$.
We follow the theory that has been done in the section 5, ch. 8.

Let $\stackrel{\circ }{N}$ be an apriori given nonlinear connection on $%
T^{k-1}M$, with the dual coefficients $\left( \stackrel{\circ
}{\underset{ (1)}{M}}\text{}_j^i,...,\stackrel{\circ
}{\underset{(k-1)}{M}}\text{}_j^i\right) $ and the mapping
\[
\xi ^{*}:u^{*}=(x,y^{(1)},...,y^{(k-1)},p)\in T^{*k}M\rightarrow
u=(x,y^{(1)},...,y^{(k-1)},y^{(k)})\in T^kM
\]
defined by
\begin{equation}
\left\{
\begin{array}{l}
y^{(0)i}=x^i,y^{(1)i}=y^{(1)i},...,y^{(k-1)i}=y^{(k-1)i}, \vspace{3mm}\\
y^{(k)i}=\xi ^{*i}(x,y^{(1)},...,y^{(k-1)},p),
\end{array}
\right.  \tag{10.4.11}
\end{equation}
where $\xi ^{*i}$ is expressed from the formula
\begin{equation}
k\xi ^{*i}+(k-1)\stackrel{\circ }{\underset{(1)}{M}}\text{ }
_j^iy^{(k-1)j}+\cdots +\stackrel{\circ
}{\underset{(k-1)}{M}}\text{ } _j^iy^{(1)j}=\displaystyle\frac
k2\stackrel{\cdot }{\partial }^iK^2. \tag{10.4.12}
\end{equation}

\begin{teo}
The mapping $\xi ^{*}$ is a local diffeomorphism which preserves the fibres
of $T^{*k}M$ and $T^kM$.
\end{teo}

Indeed, by means of theorem 7.5.1, the Jacobian of $\xi ^{*}$ is $\det
||g^{ij}||$, $g^{ij}$ being the fundamental tensor of the Cartan space $%
\mathcal{C}^{(k)n}$. q.e.d.

The formulas (10.4.11) and (10.4.12) imply:
\begin{equation}
\begin{array}{l}
\stackrel{\cdot }{\partial }^i\xi ^{*j}=g^{ij}, \\
z^{(k)i}(x,y^{(1)},...,y^{(k-1)},\xi ^{*}(u^*))=\displaystyle\frac 12\stackrel{
\cdot }{\partial }^iK^2(x,y^{(1)},...,y^{(k-1)},p),
\end{array}
\tag{10.4.13}
\end{equation}
where $z^{(k)i}$ is the Liouville vector field (10.4.9).

Let $F^2=\xi ^{*}(K^2)$ be the Lagrangian

\begin{equation}
F^2(x,y^{(1)},...,y^{(k)})=2p_iz^{(k)i}(x,y^{(1)},...,y^{(k-1)},\xi
^{*})-K^2(x,y^{(1)},...,y^{(k-1)},\varphi ^{*}).  \tag{10.4.14}
\end{equation}
where $\varphi ^{*}$ is the local inverse of $\xi ^{*}$.

Theorem 8.5.2 allows to state:

\begin{teo}
The pair $(M,F)$, with $F$ from (10.4.14), has the properties:

1$^0$ It is a Finsler space, having the fundamental function $F^2$, $2k$
-homogeneous on the fibres of $T^kM$.

2$^0$ Its fundamental tensor field is given by
\[
a_{ij}(u)=g_{ij}(x,y^{(1)},...,y^{(k-1)},\varphi ^{*}(u)).
\]
\end{teo}

So, the space $(M,F)$ is called \textit{the Finsler space of order }$k$
\textit{\ dual to the Cartan space of order }$k$\textit{, }$\mathcal{C}
^{(k)n}=(M,K)$.

>From (10.4.14) we deduce

\begin{equation}
p_i=\displaystyle\frac 12\displaystyle\frac{\partial F^2}{\partial y^{(k)i}}
,\ p_i=\varphi _i^{*}(x,y^{(1)},...,y^{(k)}).  \tag{10.4.15}
\end{equation}

Therefore, by means of the Theorem 8.5.3, we have:

\begin{teo}
The Legendre transformation determined by the Lagrange space $%
L^{(k)n}=(M,F^2)$, with $F^2$ from (10.4.14) is defined by the local
diffeomorphism $\varphi ^{*}$, the inverse of the local diffeomorphism $\xi ^{*}$
and $K^2=\varphi ^{*}(F^2)$ is given by
\[
K^2(x,y^{(1)},...,y^{(k-1)},p)=2p_iz^{(k)i}(x,y^{(1)},...,y^{(k-1)},\xi
^{*})-F^2(x,y^{(1)},...,y^{(k-1)},\xi ^{*}).
\]
\end{teo}

\begin{rem}
It follows $\varphi =\varphi ^{*}$ and $\varphi ^{-1}=\xi ^{*}$. Therefore,
locally, we get $K^2=\varphi ^{*}(\xi ^{*}(K^2))$ and $F^2=\xi ^{*}(\varphi
^{*}(F^2))$.
\end{rem}

\section{Canonical Nonlinear Connection. $N$-Linear Connections}

As we know, from section 6, Ch. 8 and from Theorem 10.4.6 we can determine a
nonlinear connection $N^{*}$ of the Cartan space $\mathcal{C}^{(k)n}=(M,K)$
by means of an apriori given nonlinear connection $\stackrel{\circ }{N}$ on
the manifold $T^{k-1}M$. We construct a bundle morphism and determine the
space $L^{(k)n}=(M,F^2)$, with $F^2=\xi ^{*}(K^2)$.

The Legendre transformation $\varphi ^{*}=\xi ^{*-1}$ (with $\varphi _i^{*}= %
\displaystyle\frac 12\displaystyle\frac{\partial F^2}{\partial y^{(k)i}}$)
transforms the $k$-spray $S$ of $L^{(k)n}$ in the dual $k$-spray $S_{\xi
^{*}}^{*} $:

\begin{equation}
S_{\xi ^{*}}^{*}=y^{(1)i}\displaystyle\frac \partial {\partial x^i}+\cdots
+(k-1)y^{(k-1)i}\displaystyle\frac \partial {\partial y^{(k-2)i}}+k\xi ^{*i} %
\displaystyle\frac \partial {\partial y^{(k-1)i}}+\eta _i\displaystyle\frac
\partial {\partial p_i}  \tag{10.5.1}
\end{equation}

\begin{teo}
The following systems of functions
\begin{equation}
\underset{(1)}{M}_j^{*i}=-\displaystyle\frac{\partial \xi
^{*i}}{\partial y^{(1)j}},\ ...,\
\underset{(k-1)}{M}_j^{*i}=-\displaystyle\frac{\partial \xi
^{*i}}{\partial y^{(k-1)j}},  \tag{10.5.2}
\end{equation}
and
\begin{equation}
N_{ij}^{*}=\displaystyle\frac{\delta \eta _i}{\delta y^{(1)j}},  \tag{10.5.2a}
\end{equation}
where the operators $\displaystyle\frac \delta {\delta y^{(1)i}}$ are
constructed by means of (10.5.2), give the dual coefficients of a nonlinear
connection $N^{*}$ which depends on the fundamental function $K$ of Cartan
space $\mathcal{C}^{(k)n}=(M,K)$ and on the apriori given nonlinear
connection $\stackrel{\circ }{N}$ of the manifold $T^{k-1}M$.
\end{teo}

Indeed, this is a particular case of the Theorem 8.6.1.

The nonlinear connection $N^{*}$ is called \textit{canonical} for the Cartan
space $\mathcal{C}^{(k)n}=(M,K)$.

In the following we denote $N^{*}$ by $N$ and consider the adapted basis and
adapted cobasis determined by $N$ and by the vertical distribution $W_k$:

\begin{equation}
\left\{ \displaystyle\frac \delta {\delta x^i},\displaystyle\frac \delta
{\delta y^{(1)i}},...,\displaystyle\frac \delta {\delta y^{(k-1)i}}, %
\displaystyle\frac \partial {\partial p_i}=\stackrel{\cdot }{\partial }
^i\right\}  \tag{10.5.3}
\end{equation}
and
\begin{equation}
\left\{ dx^i,\delta y^{(1)i},...,\delta y^{(k-1)i},\delta p_i\right\} .
\tag{10.5.3a}
\end{equation}

Let $D$ be a $N$-linear connection with the coefficients $D\Gamma
(N)=$ \newline $=(H_{jh}^i,\underset{(\alpha )}{C}\text{}_{jh}^i,C_i^{jh})$,
($\alpha =1,...,k-1$).

The fundamental tensor field $g^{ij}$ of $\mathcal{C}^{(k)n}$ is absolute
parallel with respect to $D$ if

\begin{equation}
g_{\ \ |h}^{ij}=0,\ g^{ij}\stackrel{(\alpha )}{|}_h=0,\ (\alpha
=1,...,k-1),\ g^{ij}|^h=0.  \tag{10.5.4}
\end{equation}

Assume that the equations (10.5.4) hold. Then $D$ is called the \textit{metrical }$%
N $\textit{-linear connection}.

In the case when $h$-, $v_\alpha $-, and $w_k$-torsions of $D$ vanish, $D$
is called \textit{canonical metrical }$N$\textit{-linear connection} of the
Cartan space $\mathcal{C}^{(k)n}$.

Theorem 8.7.1 leads to the following important result:

\begin{teo}
1$^0$ There exists an unique canonical metrical $N$-linear connection $D$ of
the space $\mathcal{C}^{(k)n}=(M,K)$. Its coefficients are given by the
generalized Christoffel symbols:
\begin{equation}
\begin{array}{l}
H_{jh}^i=\displaystyle\frac 12g^{is}\left( \displaystyle\frac{\delta g_{sh}}{%
\delta x^j}+\displaystyle\frac{\delta g_{js}}{\delta x^h}-\displaystyle\frac{%
\delta g_{jh}}{\delta x^s}\right) , \vspace{3mm}\\
\underset{(\alpha )}{C}\text{}_{jh}^i=\displaystyle\frac 12g^{is}\left( %
\displaystyle\frac{\delta g_{sh}}{\delta y^{(\alpha )j}}+\displaystyle\frac{%
\delta g_{js}}{\delta y^{(\alpha )h}}-\displaystyle\frac{\delta g_{jh}}{%
\delta y^{(\alpha )s}}\right) ,\ (\alpha =1,...,k-1), \vspace{3mm}\\
C_i^{jh}=-\displaystyle\frac 12g_{is}\left( \displaystyle\frac{\partial
g^{sh}}{\partial p_j}+\displaystyle\frac{\partial g^{js}}{\partial p_h}-%
\displaystyle\frac{\partial g^{jh}}{\partial p_s}\right) .
\end{array}
\tag{10.5.5}
\end{equation}

2$^0$ This connection depends only on the canonical nonlinear connection $N$ and
on the fundamental function $K$ of the Cartan space of order $k$, $\mathcal{C%
}^{(k)n}$.

3$^0$ The coefficients $(H_{jh}^i,\underset{(1)}{C}\text{}_{jh}^i,..., \underset{(k-1)}{C}\text{}_{jh}^i,C_i^{jh})$ are homogeneous on the
fibres of bundle $T^{*k}M$ of degree: $0$, $-1$, ..., $1-k$, $-k$,
respectively.
\end{teo}

\begin{cor}
The following identities hold
\begin{equation}
\begin{array}{l}
\mathcal{L}_{\stackrel{k-1}{\Gamma }+kC^{*}}H_{jh}^i=0, \vspace{3mm}\\
\mathcal{L}_{\stackrel{k-1}{\Gamma }+kC^{*}}\underset{(\alpha )}{C}\text{}%
_{jh}^i=-\alpha \underset{(\alpha )}{C}\text{ }_{jh}^i,\ (\alpha
=1,...,k-1), \vspace{3mm}\\
\mathcal{L}_{\stackrel{k-1}{\Gamma }+kC^{*}}C_i^{jh}=-kC_i^{jh}.
\end{array}
\tag{10.5.6}
\end{equation}
\end{cor}

\begin{cor}
We have
\begin{equation}
C_i^{jh}=g_{is}C^{sjh}.  \tag{10.5.7}
\end{equation}
\end{cor}

Consider the $d$-Liouville vector fields of the space $\mathcal{C}^{(k)n}$: $%
z^{(1)i}$, ..., $z^{(k-1)i}$, (cf. (6.2.7), (6.2.7')). They are
homogeneous on the fibres of $T^{*k}M$ of degree $1$, $2$, ..., $k-1$,
respectively.

The deflection tensor fields of the canonical metrical $N$-linear connection $D$
are given by (see ch. 7, section 5):
\begin{equation}
\begin{array}{l}
\stackrel{(\alpha )}{D}\text{ }_j^i=\stackrel{(\alpha )}{z}\text{ }_{|j}^i= %
\displaystyle\frac{\delta \stackrel{(\alpha )}{z}\text{ }^i}{\delta x^j}+
\stackrel{(\alpha )}{z}\text{ }^mH_{mj}^i\text{,} \\
\stackrel{(\alpha \beta )}{D}\text{ }_j^i=\stackrel{(\alpha )}{z}\text{ }^i
\stackrel{(\beta )}{|}_j=\displaystyle\frac{\delta \stackrel{(\alpha )}{z}
\text{ }^i}{\delta y^{(\beta )j}}+\stackrel{(\alpha )}{z}\text{ }^m
\underset{(\beta )}{C}\text{ }_{mj}^i\text{,} \\
\stackrel{(\alpha )}{D}\text{ }^{ij}=\stackrel{(\alpha )}{z}\text{ }^i|^j= %
\displaystyle\frac{\delta \stackrel{(\alpha )}{z}\text{ }^i}{\delta p_j}+
\stackrel{(\alpha )}{z}\text{ }^mC_m^{ij}\text{,\ (}\alpha ,\beta =1,...,k-1
\text{)}
\end{array}
\tag{10.5.8}
\end{equation}
and
\begin{equation}
\begin{array}{l}
\Delta _{ij}=p_{i|j}=-p_mH_{ij}^m, \\
\stackrel{(\alpha )}{\not \delta }_{ij}=p_i\stackrel{(\alpha )}{|}_j=-p_m
\underset{(\alpha )}{C}\text{ }_{ij}^m,\ (\alpha =1,...,k-1), \\
\not \delta _i^j=p_i|^j=\delta _i^j-p_mC_i^{mi}.
\end{array}
\tag{10.5.8a}
\end{equation}

The degrees of homogeneity on the fibres of $T^{*k}M$ for these tensor fields are
easily determined.

The deflection tensors (10.5.8) and (10.5.8a) lead to important identities for the
canonical metrical $N$-linear connection $D$, derived from the Ricci
identities applied to the Liouville $d$-vector fields $z^{(1)i}$, ..., $%
z^{(k-1)i}$.

Indeed, Theorem 7.6.2 and 7.6.3 give us:
\begin{teo}
The canonical metrical $N$-linear connection of Cartan space $\mathcal{C}%
^{(k)n}$ satisfies the following identities:
\[
\stackrel{(\alpha )}{D}_{\ j|h}^i-\stackrel{(\alpha )}{D}_{\
h|j}^i=z^{(\alpha )s}R_{s\ jh}^{\ i}-\sum\limits_{\beta
=1}^{k-1}\left\{ \stackrel{(\alpha \beta )}{D}\text{}_{\
s}^i\underset{(0\beta )}{R}\text{}_{\ jh}^s\right\}
+\stackrel{(\alpha )}{D}^{is}\underset{(0)}{R}\text{}_{sjh},
\]
$$\begin{array}{l}
\stackrel{(\alpha )}{D}_{\ j}^i\stackrel{(\beta )}{|}_h-\stackrel{(\alpha
\beta )}{D}_{\ h|j}^i=z^{(\alpha )s}\underset{(\beta )}{P}\text{}_{\ jh}^i-%
\stackrel{(\alpha )}{D}_s^i\underset{(\beta
)}{C}\text{}_{jh}^s-\stackrel{(\alpha
\beta )}{D}_s^iH_{jh}^s- \\
-\sum\limits_{\gamma =1}^{k-1}\left\{ \stackrel{%
(\alpha \gamma )}{D}\text{}_{\ s}^i\underset{(\beta \gamma )}{B}\text{}%
_{\ jh}^s\right\} -\stackrel{(\alpha )}{D}\text{}^{is}\underset{(\beta
)}{B}\text{}_{sjh},\\
\stackrel{(\alpha )}{D}\text{}_{\ j}^i|^h-\stackrel{(\alpha )}{D}\text{}_{\ \
|j}^{ih}=z^{(\alpha )s}P_{s\ j}^{\ i\ h}-\stackrel{(\alpha
)}{D}\text{}_s^iC_{\ \ j}^{sh}-\stackrel{(\alpha )}{D}\text{}^{is}H_{s\ j}^{\
h}-\\
-\sum\limits_{\gamma
=1}^{k-1}\left\{ \stackrel{(\alpha \gamma )}{D}\text{}_{\ s}^i\underset{%
(\gamma )}{B}\text{}_{\ \ j}^{s\ \ h}\right\} -\stackrel{(\alpha )}{D}\text{}^{is}%
\underset{(0)}{B}\text{}_{sj}^{\ h},
\end{array}
$$

\begin{equation}
\begin{array}{c}
\stackrel{(\alpha \beta )}{D}\text{}_j^i\stackrel{(\gamma )}{|}_h-\stackrel{%
(\alpha \gamma )}{D}\text{}_h^i\stackrel{(\beta )}{|}_j=z^{(\alpha )s}%
\underset{(\beta \gamma )}{S}\text{}_{s\ jh}^{\
i}-\stackrel{(\alpha
\beta )}{D}\text{}_{is}\underset{(\gamma )}{C}\text{}_{jh}^s+\stackrel{%
(\alpha \gamma )}{D}\text{}^i_s\underset{(\beta )}{C}\text{}_{hj}^s- \\
 \\
\qquad \qquad \qquad \qquad \qquad -\sum\limits_{\sigma
=1}^{k-1}\left\{
\stackrel{(\alpha \sigma )}{D}\text{}_{\ s}^i\underset{(\beta \gamma )}{%
\stackrel{(\sigma )}{C}}\text{}_{\ jh}^s\right\} -\stackrel{(\alpha )}{D}%
\text{}^{is}\underset{(\beta \gamma )}{B}\text{}_{\ sjh},
\end{array}
\tag{10.5.9}
\end{equation}
\[
\begin{array}{c}
\stackrel{(\alpha \beta )}{D}\text{}_{\ j}^i|^h-\stackrel{(\alpha )}{D}%
\text{}^{ih}\stackrel{(\beta )}{|}_j=z^{(\alpha )s}\underset{(\beta )}{S}%
\text{}_{s\ j}^{\ i\ h}-\stackrel{(\alpha )}{D}\text{}^{is}\underset{%
(\beta )}{C}\text{}_{s\ j}^{\ \ h}-\stackrel{(\alpha \beta )}{D}\text{}%
_s^iC_j^{sh}- \\
 \\
\qquad \qquad \qquad \qquad \qquad -\sum\limits_{\sigma
=1}^{k-1}\left\{
\stackrel{(\alpha \sigma )}{D}\text{}_{\ s}^i\underset{(\beta \sigma )}{B}%
\text{}_{\ j}^{s\ h}\right\} -\stackrel{(\alpha )}{D}\text{}^{is}%
\underset{(\beta )}{C}\text{}_{sj}^{\ \ h},
\end{array}
\]
\[
\stackrel{(\alpha )}{D}^{ij}|^h-\stackrel{(\alpha )}{D}^{ih}|^j=z^{(\alpha
)s}S_s^{ijh}
\]
\textit{and}
\[
\Delta _{ij|h}-\Delta _{ih|j}=-p_sR_{i\ jh}^{\ s}-\sum\limits_{\gamma
=1}^{k-1}\left( \stackrel{(\gamma )}{\not \delta }_{is}\underset{(0\gamma )%
}{R}\text{}_{jh}^s\right) -\not \delta
_i^s\underset{(0)}{R}\text{}_{sjh},
\]
\[
\Delta _{ij}\stackrel{(\beta )}{|}_h-\not \delta _{ih|j}=-p_s\underset{%
(\beta )}{P}\text{}_{i\ jh}^{\ s}-\Delta _{is}\underset{(\beta )}{C}\text{}_{jh}^s-%
\stackrel{(\beta )}{\not \delta }_{is}H_{jh}^s-\sum\limits_{\gamma
=1}^{k-1}\left( \stackrel{(\gamma )}{\not \delta
}_{is}\underset{(\alpha
\gamma )}{B}\text{}_{\ jh}^s\right) -\not \delta _i^s\underset{(\beta )}{B%
}\text{}_{sjh},
\]
\[
\Delta _{ij}|^h-\not \delta _{i\ |j}^{\ h}=-p_sP_{i\ \ j}^{\ s\ \
h}-\Delta _{is}C_j^{sh}-\not \delta _i^sH_{s\ j}^{\
h}-\sum\limits_{\gamma
=1}^{k-1}\left( \stackrel{(\gamma )}{\not \delta }_{is}\underset{(\gamma )%
}{B}\text{}_j^{sh}\right) -\not \delta
_i^s\underset{(0)}{B}\text{}_{sj}^{\ h},
\]
\begin{equation}
\begin{array}{c}
\stackrel{(\beta )}{\not \delta }_{ij}\stackrel{(\gamma )}{|}_h-\stackrel{%
(\gamma )}{\not \delta }_{ih}\stackrel{(\beta
)}{|}_j=-p_s\underset{(\beta\gamma)}{S}\text{}_{i\ \ jh}^{\ s}-\stackrel{(\beta )}{\not \delta }_{is} \underset{(\gamma )}{C}\text{}_{\ jh}^s- \\
-\stackrel{(\gamma )}{\not \delta }_{is}\underset{(\beta )}{C}\text{}_{hj}^s-\\
-\sum\limits_{\sigma =1}^{k-1}\left( \stackrel{(\sigma )}{\not \delta }%
_{is}\underset{(\beta \gamma )}{\stackrel{(\sigma )}{C}}\text{}_{\ jh}^s\right)
-\not \delta _i^s\underset{(\beta \gamma)}{B}\text{}_{sjh},
\end{array}
\tag{10.5.9a}
\end{equation}
\[
\stackrel{(\beta )}{\not \delta }_{ij}|^h-\not \delta _i^h\stackrel{(\beta )%
}{|}_j=-p_s\underset{(\beta )}{S}\text{}_{i\ \ j}^{\ s\ \ h}-\not \delta _i^s%
\underset{(\beta )}{C}\text{}_{s\ j}^{\ h}-\stackrel{(\beta )}{\not \delta }\text{}_{is}C_j^{sh}-\sum\limits_{\sigma =1}^{k-1}\left(
\stackrel{(\sigma )}{\not \delta }_{is}\underset{(\beta \sigma
)}{C}\text{}_{\ j}^{sh}\right) -\not \delta _i^s\underset{(\beta)}{C}\text{}_{sj}^{\ h},
\]
\[
\not \delta _i^{\ j}|^h-\not \delta _i^{\ h}|^j=-p_sS_i^{\ sjh}.
\]
\end{teo}

The Ricci identities and Bianchi identities of the canonical metrical $N$
-linear connection of the Cartan space $\mathcal{C}^{(k)n}$ can be written
using the corresponding indentities of the $N$-linear connections of $%
T^{*k}M $, given in the chapter 6.

Applying the Ricci identities to the fundamental tensor $g^{ij}$ of the
space $\mathcal{C}^{(k)n}$, with respect to the canonical metrical $N$
-linear connection, we obtain the identities (7.6.8).

So, we have:
\begin{teo}
In a Cartan space $\mathcal{C}^{(k)n}$, with respect to the canonical
metrical $N$-linear connection, the following identities hold:
\[
\begin{array}{l}
g^{sj}R_{s\ \ hm}^{\ \ i}+g^{is}R_{s\ \ hm}^{\ \ j}=0, \\
\\
g^{sj}\underset{(\alpha )}{P}\text{}_{s\ \ hm}^{\ \ i}+g^{is}\underset{%
(\alpha )}{P}\text{}_{s\ \ hm}^{\ \ j}=0,\ (\alpha =1,...,k-1), \\
\\
g^{sj}S_s^{\ ihm}+g^{is}S_s^{\ jhm}=0.
\end{array}
\]
\end{teo}

\section{Parallelism of Vector Fields in Cartan Space $\mathcal{C}^{(k)n}$}

Consider a Cartan space of order $k$, $\mathcal{C}^{(k)n}=(M,K)$,
endowed with the canonical metrical $N$-linear connection $C\Gamma
(N)=(H_{jm}^i, \underset{(\alpha )}{C}$ $_{jm}^i,C_i^{jm})$,
($\alpha =1,...,k-1$), the
coefficients being given by (10.5.5). The local vector fields $\left( %
\displaystyle\frac \delta {\delta x^i},\displaystyle\frac \delta {\delta
y^{(\alpha )i}},\displaystyle\frac \delta {\delta p_i}\right) $, ($\alpha
=1,..,k-1$) determine an adapted basis and \newline $\left( dx^i,\delta y^{(\alpha
)i},\delta p_i\right) $, ($\alpha =1,...,k-1$) is its dual basis.

Along a smooth parametrized curve $\gamma :I\rightarrow T^{*k}M$, having
the image in a domain of a local chart:
\begin{equation}
x^i=x^i(t),\ y^{(\alpha )i}=y^{(\alpha )i}(t),\ p_i=p_i(t),\ t\in I,\
(\alpha =1,..,k-1)  \tag{10.6.1}
\end{equation}
the tangent vector field $\stackrel{\cdot }{\gamma }$ is given by (7.7.2), i.e.
\begin{equation}
\stackrel{\cdot }{\gamma
}=\displaystyle\frac{dx^i}{dt}\displaystyle\frac \delta {\delta
x^i}+\sum\limits_{\alpha =1}^{k-1}\displaystyle\frac{\delta
y^{(\alpha )i}}{dt}\displaystyle\frac \delta {\delta y^{(\alpha )i}}+ %
\displaystyle\frac{\delta p_i}{dt}\stackrel{\cdot }{\partial }^i,  \tag{10.6.2}
\end{equation}
where
\begin{equation}
\begin{array}{l}
\displaystyle\frac{\delta y^{(\alpha
)i}}{dt}=\displaystyle\frac{dy^{(\alpha
)i}}{dt}+\underset{(1)}{M}\text{
}_j^i\displaystyle\frac{dy^{(\alpha -1)j} }{dt}+\cdots
+\underset{(\alpha -1)}{M}\text{ }_j^i\displaystyle\frac{
dy^{(1)j}}{dt}+\underset{(\alpha )}{M}\text{
}_j^i\displaystyle\frac{dx^j}{
dt}, \vspace{3mm}\\
\displaystyle\frac{\delta p_i}{dt}=\displaystyle\frac{dp_i}{dt}-N_{ji} %
\displaystyle\frac{dx^j}{dt},\ \ \ \ \ (\alpha =1,...,k-1).
\end{array}
\tag{10.6.3}
\end{equation}

Consider the $1$-forms of metrical canonical $N$-linear connection $C\Gamma
(N)$:

\begin{equation}
\omega _{\ j}^i=H_{js}^idx^s+\sum\limits_{\alpha
=1}^{k-1}\underset{ (\alpha )}{C}_{js}^i\delta y^{(\alpha
)s}+C_j^{is}\delta p_s.  \tag{10.6.4}
\end{equation}

Then the vector field $X\in \mathcal{X}(T^{*k}M)$:

\begin{equation}
X=\stackrel{(0)}{X}^i\displaystyle\frac \delta {\delta
x^i}+\sum\limits_{\alpha =1}^{k-1}\stackrel{(\alpha )}{X}^i\displaystyle %
\frac \delta {\delta y^{(\alpha )i}}+X_i\stackrel{\cdot }{\partial }^i
\tag{10.6.5}
\end{equation}
has the covariant differential along curve $\gamma $:
\begin{equation}
\begin{array}{l}
\displaystyle\frac{DX}{dt}=\left( \displaystyle\frac{d\stackrel{(0)}{X}^i}{%
dt }+\stackrel{(0)}{X}^s\displaystyle\frac{\omega _s^i}{dt}\right) %
\displaystyle \frac \delta {\delta x^i}+\sum\limits_{\alpha =1}^{k-1}\left( %
\displaystyle \frac{d\stackrel{(\alpha )}{X}^i}{dt}+\stackrel{(\alpha )}{X}^s%
\displaystyle \frac{\omega _s^i}{dt}\right) \displaystyle\frac \delta
{\delta y^{(\alpha )i}}+\\
\qquad +\left( \displaystyle\frac{dX_i}{dt}-X_s\displaystyle%
\frac{\omega _i^s}{ dt}\right) \stackrel{\cdot }{\partial }^i.  \tag{10.6.6}
\end{array}
\end{equation}

Theorem 7.7.1 takes the form:
\begin{teo}
The vector $X$ from (10.6.5) is parallel along curve $\gamma $ if
and only if $(\stackrel{(0)}{X}^i,\stackrel{(\alpha )}{X}^i,X_i)$
are the solutions of the system of differential equations:
\begin{equation}
\left\{
\begin{array}{l}
\displaystyle\frac{d\stackrel{(0)}{X}^i}{dt}+\stackrel{(0)}{X}^s%
\displaystyle \frac{\omega _s^i}{dt}=0, \vspace{3mm}\\
\displaystyle\frac{d\stackrel{(\alpha )}{X}^i}{dt}+\stackrel{(\alpha )}{X}^s%
\displaystyle\frac{\omega _s^i}{dt}=0,\ (\alpha =1,...,k-1), \vspace{3mm}\\
\displaystyle\frac{dX_i}{dt}-X_s\displaystyle\frac{\omega _i^s}{dt}=0.
\end{array}
\right.   \tag{10.6.7}
\end{equation}
\end{teo}
We obtain, also
\begin{teo}
The Cartan spaces $\mathcal{C}^{(k)n}$ is with absolute parallelism of
vectors, with respect to $C\Gamma (N)$, if and only if all curvature $d$
-tensors of $C\Gamma (N)$ vanish.
\end{teo}

In the case $X=\stackrel{\cdot }{\gamma }$, the equation $\displaystyle\frac{
D\stackrel{\cdot }{\gamma }}{dt}=0$ says that $\gamma $ is an autoparallel curve
of $\mathcal{C}^{(k)n}$ with respect to $C\Gamma (N)$.

>From Theorem 7.7.3 it follows:

\begin{teo}
A smooth parametrized curve $\gamma $, (10.6.1), is an autoparallel curve of
the Cartan space $\mathcal{C}^{(k)n}$, endowed with metrical canonical $N$
-linear connection $C\Gamma (N)$, if and only if the following system of
differential equations is verified:
\begin{equation}
\left\{
\begin{array}{l}
\displaystyle\frac{d^2x^i}{dt^2}+\displaystyle\frac{dx^s}{dt}\displaystyle
\frac{\omega _s^i}{dt}=0, \vspace{3mm}\\
\displaystyle\frac d{dt}\displaystyle\frac{\delta y^{(\alpha )i}}{dt}+%
\displaystyle\frac{\delta y^{(\alpha )s}}{dt}\displaystyle\frac{\omega _s^i}{%
dt}=0, \vspace{3mm}\\
\qquad \qquad \ (\alpha =1,...,k-1), \vspace{3mm}\\
\displaystyle\frac d{dt}\displaystyle\frac{\delta p_s}{dt}-\displaystyle
\frac{\delta p_s}{dt}\displaystyle\frac{\omega _i^s}{dt}=0.
\end{array}
\right.   \tag{10.6.8}
\end{equation}
\end{teo}

As we know from section 7, ch. 7, a curve $\gamma $ is horizontal if $%
\stackrel{\cdot }{\gamma }=\stackrel{\cdot }{\gamma }^H$, i.e.

\begin{equation}
x^i=x^i(t),\ \displaystyle\frac{\delta y^{(\alpha )i}}{dt}=0,\ (\alpha
=1,...,k-1),\ \displaystyle\frac{\delta p_i}{dt}=0,\ t\in I.  \tag{10.6.9}
\end{equation}

Therefore, taking into account the definitions of horizontal paths, $v_\alpha $-paths and $w_k$-paths (\S 7, ch. 6), we obtain:
\begin{teo}
The Cartan space $\mathcal{C}^{(k)n}$ endowed with the metrical canonical $N$
-linear connection $C\Gamma (N)$ has the following properties:

a. The horizontal paths are characterized by the system of differential
equations:
\begin{equation}
\displaystyle\frac{d^2x^i}{dt^2}+H_{jh}^i\displaystyle\frac{dx^j}{dt}%
\displaystyle\frac{dx^h}{dt}=0,\ \displaystyle\frac{\delta y^{(\alpha )i}}{dt%
}=0,\ \displaystyle\frac{\delta p_i}{dt}=0.  \tag{10.6.10}
\end{equation}

b. The $v_\alpha $-paths at the point $x=x_0$ are characterized by the
system of differential equations:
\begin{equation}
\begin{array}{l}
\displaystyle\frac{dx^i}{dt}=0,\ \displaystyle\frac{dy^{(\beta )i}}{dt}=0,\
(\beta \neq \alpha ),\ \displaystyle\frac{dp_i}{dt}=0, \vspace{3mm}\\
\displaystyle\frac d{dt}\displaystyle\frac{dy^{(\alpha )i}}{dt}+\underset{%
(\alpha )}{C}\text{}_{sj}^{\ i}\displaystyle\frac{dy^{(\alpha )s}}{dt}%
\displaystyle\frac{dy^{(\alpha )j}}{dt}=0.
\end{array}
\tag{10.6.11}
\end{equation}

c. The $w_k$-paths at the point $x=x_0$ are characterized by the system of
differential equations:
\begin{equation}
\begin{array}{l}
\displaystyle\frac{dx^i}{dt}=0,\ \displaystyle\frac{dy^{(\alpha )i}}{dt}=0,\
(\alpha =1,...,k-1), \vspace{3mm}\\
\displaystyle\frac{d^2p_i}{dt^2}-C_i^{js}(x,0,...,0,p)\displaystyle\frac{dp_j%
}{dt}\displaystyle\frac{dp_s}{dt}=0.

\end{array}
\tag{10.6.12}
\end{equation}
\end{teo}

\section{Structure Equations of Metrical Canonical $N$-Connection}

In a Cartan space of order $k$, $\mathcal{C}^{(k)n}=(M,K)$, endowed with the
metrical canonical $N$-linear connection $C\Gamma (N)$, lemma 7.8.1 (ch. 7)
holds:

The following object fields
\[
\begin{array}{l}
d(dx^i)-dx^m\wedge \omega _m^i,\ d(\delta y^{(\alpha )i})-dy^{(\alpha
)m}\wedge \omega _m^i,\ (\alpha =1,...,k-1), \vspace{3mm}\\
d(\delta p_i)+\delta p_m\wedge \omega _m^i
\end{array}
\]
are $d$-vector fields and
\[
d\omega _j^i-\omega _j^m\wedge \omega _m^i
\]
is a $d$-tensor field of type $(1,1)$.

Consequently, we have:

\begin{teo}
A Cartan space of order $k$, $\mathcal{C}^{(k)n}=(M,K)$, has the following
structure equations of the metrical canonical $N$-linear connection $C\Gamma
(N)$:
\begin{equation}
\begin{array}{l}
d(dx^i)-dx^m\wedge \omega _m^i=-\stackrel{(0)}{\Omega }_i, \\
\\
d(\delta y^{(\alpha )i})-dy^{(\alpha )m}\wedge \omega _m^i=-\stackrel{%
(\alpha )}{\Omega }_i,\ (\alpha =1,...,k-1), \\
\\
d(\delta p_i)+\delta p_m\wedge \omega _m^i=-\Omega _i.
\end{array}
\tag{10.7.1}
\end{equation}
and
\begin{equation}
d\omega _j^i-\omega _j^m\wedge \omega _m^i=-\Omega _j^i,  \tag{10.7.2}
\end{equation}
where the $2$-forms of torsion are:
\begin{equation}
\begin{array}{l}
\stackrel{(0)}{\Omega }_i=dx^j\wedge \left( \sum\limits_{\alpha =1}^{k-1}%
\underset{(\alpha )}{C}\text{}_{jm}^i\delta y^{(\alpha
)m}+C_j^{im}\delta
p_m\right) , \\
\\
\stackrel{(\alpha )}{\Omega }_i=dx^j\wedge \underset{(\alpha 0)}{P}\text{}%
_j^i+\sum\limits_{\gamma =1}^{k-1}\delta y^{(\gamma )j}\wedge \underset{%
(\alpha \gamma )}{P}\text{}_j^i+ \\
\qquad +\delta y^{(\alpha )j}\wedge \left\{
H_{jm}^idx^m+\sum\limits_{\gamma =1}^{k-1}\underset{(\gamma )}{C}\text{}%
_{jm}^i\delta y^{(\gamma )m}+C_j^{im}\delta p_m\right\} ,\ (\alpha
=1,...,k-1), \\
\\
\Omega _i=dx^j\wedge \left\{ \displaystyle\frac 12\underset{(0)}{R}\text{}%
_{ijm}dx^m+\sum\limits_{\gamma =1}^{k-1}\underset{(\gamma )}{B}\text{}%
_{ijm}\delta y^{(\gamma )m}+\underset{(0)}{B}\text{}_{ij}^m\delta
p_m\right\} - \\
\qquad -\delta p_j\wedge \left\{ H_{im}^jdx^m+\sum\limits_{\gamma =1}^{k-1}%
\underset{(\gamma )}{C}\text{}_{im}^j\delta y^{(\gamma
)m}+C_i^{jm}\delta p_m\right\}
\end{array}
\tag{10.7.3}
\end{equation}
and where the $2$-forms of curvature are
\begin{equation}
\begin{array}{l}
\Omega _j^i=\displaystyle\frac 12R_{jhm}^idx^h\wedge dx^m+ \vspace{3mm}\\
\qquad +\sum\limits_{\gamma =1}^{k-1}\underset{(\gamma )}{P}\text{}_{j\ \ hm}^{\ \ i}dx^h\wedge \delta y^{(\gamma )m}+P_{j\ h}^{\
i\ m}dx^h\wedge \delta p_m+ \vspace{3mm}\\
\qquad +\sum\limits_{\alpha \leq \beta ,\alpha
=1}^{k-1}\sum\limits_{\beta =1}^{k-1}\underset{(\alpha \beta
)}{S}\text{}_{j\ hm}^{\ i}\delta
y^{(\alpha )h}\wedge \delta y^{(\beta )m}+ \vspace{3mm}\\
\qquad +\sum\limits_{\gamma =1}^{k-1}\underset{(\gamma )}{S}\text{}_{j\ hm}^{\ i}\delta y^{(\gamma )h}\wedge \delta
p_m+\displaystyle\frac 12S_{j}^{\ ihm}\delta p_h\wedge \delta
p_m,
\end{array}
\tag{10.7.4}
\end{equation}
where, according to \S 5, ch. 6, $\underset{(\alpha ,0)}{P}\text{}_j^i$, ..., $\underset{(\alpha ,\alpha -1)}{P}\text{}_j^i$ are as follows:
\begin{equation}
\begin{array}{l}
\underset{(\alpha ,0)}{P}\text{}_j^i=d\underset{(\alpha )}{M}\text{}%
_j^i-\underset{(1)}{N}\text{}_j^md\underset{(\alpha -1)}{M}\text{}%
_m^i-\cdots -\underset{(\alpha -1)}{N}\text{}_j^md\underset{(1)}{M}%
\text{}_m^i, \\
\\
\underset{(\alpha ,1)}{P}\text{}_j^i=d\underset{(\alpha -1)}{M}\text{}%
_j^i-\underset{(1)}{N}\text{}_j^md\underset{(\alpha -2)}{M}\text{}%
_m^i-\cdots -\underset{(\alpha -2)}{N}\text{}_j^md\underset{(1)}{M}%
\text{}_m^i, \\
.......................................................................................
\\
\underset{(\alpha ,\alpha -1)}{P}\text{}_j^i=d\underset{(1)}{M}\text{}_j^i.
\end{array}
\tag{10.7.5}
\end{equation}
\end{teo}

Of course this theorem is important for the theory of metrical canonical
connection. Also, it allows to determine the Bianchi identities of the
spaces $\mathcal{C}^{(k)n}$.

We can use the previous results in a theory of Cartan subspaces of order $k$
in the Cartan space $\mathcal{C}^{(k)n},$ cf. Ch.9.

\section{Riemannian Almost Contact Structure of the Space $\mathcal{C}%
^{(k)n} $}

Consider a Cartan space of order $k$, $\mathcal{C}
^{(k)n}=(M,K(x,y^{(1)},...,y^{(k-1)},p))$ and its canonical nonlinear
connection $N$.

The adapted basis $(\displaystyle\frac \delta {\delta x^i},
\displaystyle\frac \delta {\delta y^{(1)i}}, ..., \displaystyle\frac
\delta {\delta y^{(k-1)i}},\displaystyle\frac \delta {\delta p_i}) $,
$\displaystyle\frac \delta {\delta p_i}=\displaystyle\frac \partial
{\partial p_i}=\stackrel{\cdot }{\partial }^i $ and its dual basis
$ \left(\delta x^i,\delta y^{(1)i},...,\delta y^{(k-1)i},\delta p_i\right) $,
where $ \delta x^i=dx^i $ are well determined by $N$. As we know from section 1, the
fundamental tensor of $\mathcal{C}^{(k)n}$ is:

\begin{equation}
g^{ij}=\displaystyle\frac 12\stackrel{\cdot }{\partial }^i\stackrel{\cdot }{
\partial }^jK^2.  \tag{10.8.1}
\end{equation}

Taking into account the associated Hamilton space $H^{(k)n}=(M,K^2)$ of the
space $\mathcal{C}^{(k)n}$, we can study the Riemannian almost contact
structure of the Cartan space of order $k$ by means of the corresponding
structure of the space $H^{(k)n}=(M,K^2)$. Such that, using the covariant
tensor $g_{ij}$ of the fundamental tensor $g^{ij}$ we define the tensor

\begin{equation}
\stackrel{\vee}{\Bbb{G}}=g_{ij}dx^i\otimes
dx^j+\sum\limits_{\alpha =1}^{k-1}g_{ij}\delta y^{(\alpha
)i}\otimes \delta y^{(\alpha )j}+g^{ij}\delta p_i\otimes \delta
p_j,  \tag{10.8.2}
\end{equation}
$\stackrel{\vee}{\Bbb{G}}$ is the $N$-lift of the fundamental tensor $%
g^{ij} $ of the space $\mathcal{C}^{(k)n}$.
Since $g^{ij}$ is positively defined on $\widetilde{T^{*k}M}$ and $N$ is
given on $T^{*k}M$, it follows:

\begin{teo}
1$^0$ $\stackrel{\vee}{\Bbb{G}}$ is a Riemannian structure on the manifold $%
\widetilde{T^{*k}M}$ determined only by the fundamental tensor $g^{ij}$ of
the Cartan space $\mathcal{C}^{(k)n}$ and by the canonical nonlinear
connection $N$.

2$^0$ The distributions $N_0$, $N_1$, ..., $N_{k-2}$, $V_{k-1}$, $W_k$ are mutual
orthogonal with respect to $\stackrel{\vee}{\Bbb{G}}$ .
\end{teo}

\begin{prop}
The tensor $\stackrel{\vee}{\Bbb{G}}$ is not homogeneous on the fibres of
the bundle $T^{*k}M$.
\end{prop}

Indeed, the first term in $\stackrel{\vee}{\Bbb{G}}$ is $0$-homogeneous,
the second term is $2$-homogeneous, ..., the last term is $2k$-homogeneous.
So, the whole $\stackrel{\vee}{\Bbb{G}}$ is not homogeneous.

Let us consider the following invariants:

\begin{equation}
K_1^2=g_{ij}z^{(1)i}z^{(1)j},\ K_2^2=g_{ij}z^{(2)i}z^{(2)j},\ ...,\
K_{k-1}^2=g_{ij}z^{(k-1)i}z^{(k-1)j},\ K_0^2=g^{ij}p_ip_j,  \tag{10.8.3}
\end{equation}
where $z^{(1)i}$, ..., $z^{(k)i}$ are the Liouville vector fields determined
by the canonical nonlinear connection of space $\mathcal{C}^{(k)n}$.

Of course, all invariants $K_0^2$, $K_1^2$, ..., $K_{k-1}^2$ are positive.

Thus we can construct a new Riemannian structure on $\widetilde{T^{*k}M}$:

\begin{equation}
\stackrel{\circ }{\Bbb{G}}=g_{ij}dx^i\otimes
dx^j+\sum\limits_{\alpha =1}^{k-1}\displaystyle\frac 1{K_\alpha
^2}g_{ij}\delta y^{(\alpha )i}\otimes \delta y^{(\alpha
)j}+\displaystyle\frac 1{K_0^2}g^{ij}\delta p_i\otimes \delta p_j.
\tag{10.8.4}
\end{equation}

\begin{teo}
1$^0$ $\stackrel{\circ }{\Bbb{G}}$ is a Riemannian structure on the manifold
$\widetilde{T^{*k}M}$ determined only by $g^{ij}$ and $N$.

2$^0$ $\stackrel{\circ }{\Bbb{G}}$ is $0$-homogeneous on the fibres of the
bundle $T^{*k}M$.

3$^0$ The distributions $N_0$, $N_1$, ..., $N_{k-2}$, $V_{k-1}$, $W_k$ are
mutual orthogonal with respect to $\stackrel{\circ }{\Bbb{G}}$.
\end{teo}

The proof follows without difficulties.

The Riemannian structure $\stackrel{\vee}{\Bbb{G}}$ is of the form

\begin{equation}
\stackrel{\vee}{\Bbb{G}}=\stackrel{H}{\Bbb{G}}+\stackrel{V_1}{\Bbb{G}}
+\cdots +\stackrel{V_{k-1}}{\Bbb{G}}+\stackrel{W_k}{\Bbb{G}}  \tag{10.8.5}
\end{equation}
with
\begin{equation}
\begin{array}{l}
\stackrel{H}{\Bbb{G}}=g_{ij}dx^i\otimes dx^j,\ \stackrel{V_1}{\Bbb{G}}
=g_{ij}\delta y^{(1)i}\otimes \delta y^{(1)j},\ ..., \\
\\
\stackrel{V_{k-1}}{\Bbb{G}}=g_{ij}\delta y^{(k-1)i}\otimes \delta
y^{(k-1)j},\ \stackrel{W_k}{\Bbb{G}}=g^{ij}\delta p_i\otimes \delta p_j.
\end{array}
\tag{10.8.6}
\end{equation}

The tensors $\stackrel{H}{\Bbb{G}}$ , $\stackrel{V_\alpha }{\Bbb{G}}$, ($%
\alpha =1,...,k-1$) and $\stackrel{W_k}{\Bbb{G}}$ are $d$-tensor fields on
the manifold $\widetilde{T^{*k}M}$.

Exactly as in the chapter 8 we can prove:
\begin{teo}
The $d$-tensor fields $\stackrel{H}{\Bbb{G}}$, $\stackrel{V_\alpha }{\Bbb{G}}
$, ($\alpha =1,...,k-1$) and $\stackrel{W_k}{\Bbb{G}}$ are covariant
constant with respect to the canonical metrical $N$-connection of the Cartan
space of order $k$, $\mathcal{C}^{(k)n}$.
\end{teo}

Now, let us consider the natural almost $(k-1)n$-contact structure $\Bbb{F}$
determined by the canonical nonlinear connection $N$. It is defined by
(6.6.3):

\begin{equation}
\begin{array}{l}
\Bbb{F}\left( \displaystyle\frac \delta {\delta x^i}\right) =-\displaystyle %
\frac \partial {\partial y^{(k)i}},\ \Bbb{F}\left( \displaystyle\frac
\partial {\partial y^{(k-1)i}}\right) =\displaystyle\frac \delta {\delta
x^i},\Bbb{F}\left( \displaystyle\frac \delta {\delta y^{(\alpha )i}}\right)
=0,\ (\alpha =\overline{1,k-1}),\vspace{3mm}\\
\ \Bbb{F}\left( \displaystyle\frac \delta
{\delta p_i}\right) =0.
\end{array}
\tag{10.8.7}
\end{equation}
$\Bbb{F}$ is a tensor field of type $(1,1)$ and in adapted basis it is
expressed by (3.5.4).

The condition of normality of $\Bbb{F}$ is given by the equation (6.6.5).

But the pair of structures $(\stackrel{\vee}{\Bbb{G}},\Bbb{F})$ is a
Riemannian almost $(k-1)n$ contact structure on $\widetilde{T^{*k}M}$
determined only by $N$ and by the fundamental function $K$ of the Cartan
space $\mathcal{C}^{(k)n}$. So, we have:

\begin{teo}
On the manifold $T^{*k}M$ there exists a natural Riemannian almost $(k-1)n$
-contact structure $(\stackrel{\vee}{\Bbb{G}},\Bbb{F})$,
determined only by the canonical nonlinear connection $N$ and the fundamental
function $K$ of the Cartan space $\mathcal{C}^{(k)n}=(M,K)$.
\end{teo}

The canonical nonlinear connection $N$ and the fundamental tensor $g^{ij}$
of the Cartan space $\mathcal{C}^{(k)n}$ determine the almost $(k-1)n$
-contact structure $\stackrel{\vee}{\Bbb{F}}$:

\begin{equation}
\stackrel{\vee}{\Bbb{F}}\left( \displaystyle\frac \delta {\delta
x^i}\right) =-g_{ij}\displaystyle\frac \partial {\partial p_j},\ \stackrel{
\vee}{\Bbb{F}}\left( \displaystyle\frac \delta {\delta y^{(\alpha
)i}}\right) =0,\ (\alpha =\overline{1,k-1}),\ \stackrel{\vee}{\Bbb{F}}
\left( \displaystyle\frac \delta {\delta p_i}\right) =g^{ij}\displaystyle %
\frac \delta {\delta x^j}.  \tag{10.8.8}
\end{equation}

Theorem 6.7.1, reads:

\begin{teo}
The structure $\stackrel{\vee}{\Bbb{F}}$ of a Cartan space $\mathcal{C}%
^{(k)n}$ has the following properties:
\begin{equation}
\begin{array}{l}
\stackrel{\vee}{\Bbb{F}}=-g_{ij}\displaystyle\frac \delta {\delta
p_j}\otimes dx^i+g^{ij}\displaystyle\frac \delta {\delta x^i}\otimes \delta
p_j, \\
\\
Ker\stackrel{\vee}{\Bbb{F}}=N_1\oplus \cdots \oplus N_{k-1},\ Im\stackrel{%
\vee}{\Bbb{F}}=N_0\oplus W_k, \\
\\
rank\stackrel{\vee}{\Bbb{F}}=2n,
\end{array}
\tag{10.8.9}
\end{equation}
\begin{equation}
\stackrel{\vee}{\Bbb{F}}^3+\stackrel{\vee}{\Bbb{F}}=0.  \tag{10.8.10}
\end{equation}
\end{teo}

Consequently, $\stackrel{\vee}{\Bbb{F}}$ is an almost $(k-1)n$-contact
structure on the manifold $\widetilde{T^{*k}M}$.

We have, also:

\begin{teo}
For a Cartan space $\mathcal{C}^{(k)n}$ the following properties hold:

1$^0$ The pair $(\stackrel{\vee}{\Bbb{G}},\stackrel{\vee}{\Bbb{F}})$ is a
Riemannian almost $(k-1)n$-contact structure determined only by the canonical
nonlinear connection $N$ and the fundamental tensor $g^{ij}$.

2$^0$ The associated $2$-form is
\[
\theta =\delta p_i\wedge dx^i
\]
and if the coefficients $N_{ij}$ of $N$ are symmetric then $\theta $ is the
canonical presymplectic structure
\[
\theta =dp_i\wedge dx^i.
\]
\end{teo}

Concluding, the space $(\widetilde{T^{*k}M},\stackrel{\vee}{\Bbb{G}},
\stackrel{\vee}{\Bbb{F}})$ is the \textit{geometrical model} of the Cartan
space of order $k$, $\mathcal{C}^{(k)n}$.

\chapter{Generalized Hamilton and Cartan Spaces of Order $k$.
Applications to Hamiltonian Relativistic Optics}

\markboth{\it{THE GEOMETRY OF HIGHER-ORDER HAMILTON SPACES\ \ \ \ \ }}{\it{Generalized Hamilton and Cartan Spaces of Order} $k$}

On the total space of the dual bundle $(T^{*k}M,\pi ^{*},M)$ there exist
some geometrical structures defined by a general $d$-tensor field $%
g^{ij}(x,y^{(1)},...,y^{(k-1)},p),$ symmetric and nondegenerate, which are
suggested by the Relativistic Optics. Generally the $d$-tensor $g^{ij}$ is
not the metric tensor of a Hamilton space $H^{\left( k\right) n}$ or of a
Cartan space $\mathcal{C}^{\left( k\right) n}.$ In this case the pair $%
GH^{\left( k\right) n}=(M,g^{ij}(x,y^{\left( 1\right) },...,y^{\left(
k-1\right) },p))$ defines a 'Generalized Hamilton space of order $k$'. If $%
g^{ij}$ is $0$-homogeneous on the fibres of the bundle $T^{*k}M$
we say that the pair $G\mathcal{C}^{(k)n}$ is a 'Generalized
Cartan space of order $k$'.

We study, in this chapter, the geometry of these spaces and apply it to the theory of Hamiltonian Relativistic Optics.

\section{The Space $GH^{(k)n}$}

\begin{defi}
A generalized Hamilton space of order $k$ is a pair $%
GH^{(k)n}=(M,g^{ij}(x,y^{\left( 1\right) },...,y^{\left( k-1\right) },p),$
where

1$^{\circ }$ $g^{ij}$ is a $d$-tensor field of type $(2,0),$ symmetric and
nondegenerate on the manifold $\widetilde{T^{*k}M}:$
\begin{equation}
rank\left\| g^{ij}\right\| =n  \tag{11.1.1}
\end{equation}

2$^{\circ }$ The quadratic form $g^{ij}X_iX_j$ has a constant signature on $%
\widetilde{T^{*k}M}.$
\end{defi}

The tensor $g^{ij}$ is called fundamental for the space $GH^{(k)n}$.

In the case when the base manifold $M$ is paracompact then $\widetilde{
T^{*k}M}$ is paracompact. Thus on $\widetilde{T^{*k}M}$ there exists the $d$
-tensor $g^{ij},$ with the property that $GH^{(k)n}=(M,g^{ij})$ is a
generalized Hamilton space of order $k$.

If $g^{ij}$ is positively defined on $\widetilde{T^{*k}M},$ then the
conditions (11.1.1) is verified.

\begin{defi}
The space $GH^{(k)n}=(M,g^{ij})$ is called reducible to a Hamilton space
of order $k,$ if there exists an Hamiltonian $H(x,y^{\left( 1\right)
},...,y^{\left( k-1\right) },p)$ such that the following equality holds:
\begin{equation}
g^{ij}=\displaystyle\frac 12\stackrel{\cdot }{\partial ^i}\stackrel{\cdot }{%
\partial ^j}H  \tag{11.1.2}
\end{equation}
\end{defi}

Let us consider the $d$-tensor field:

\begin{equation}
C^{ijh}=-\displaystyle\frac 12\stackrel{\cdot }{\partial ^h}g^{ij}  \tag{11.1.3}
\end{equation}

We have:

\begin{prop}
A necessary condition for a generalized Hamilton space $%
GH^{(k)n}=(M,g^{ij})\,$ be reducible to a Hamilton space of order $k$ is that the $d-$tensor field $C^{ijk}$ be totally symmetric.
\end{prop}

Indeed, if (11.1.2) holds, then $C^{ijh}=-\displaystyle\frac
14\stackrel{\cdot }{\partial ^i}\stackrel{\cdot }{\partial
^i}\stackrel{\cdot }{\partial ^i}H$ is totally symmetric.

\textbf{Example 1.} 1$^{\circ }$ Let
$g^{ij}(x,y^{(1)},...,y^{(k-1)},p)=\gamma
^{ij}(x,y^{(1)},...,y^{(k-1)})$ be a $d$-tensor on the
$\widetilde{T^{*k}M}$ determined by the fundamental tensor $\gamma
_{ij}$ of a Finsler space of order $k-1$.

The space $GH^{(k)n}=(M,g^{ij})$ is reducible to the Hamilton space $%
H^{(k)n}=(M,g^{ij}p_ip_j).$

2$^{\circ }$ Consider the fundamental tensor $\gamma _{ij}(x,y^{\left(
1\right) },...,y^{(k-1)})$ of a Finsler space of order $k-1,$ $%
F^{(k-1)n}=(M,y^{\left( 1\right) },...,y^{(k-1)})).$ The $d$-tensor field

\begin{equation}
g^{ij}(x,y^{(1)},...,y^{(k-1)},p)=e^{-2\sigma
(x,y^{(1)},...,y^{(k-1)},p)}\gamma ^{ij}(x,y^{(1)},...,y^{(k-1)})  \tag{11.1.4}
\end{equation}
with the property $\sigma \in \mathcal{F}(\widetilde{T^{*k}M})$ and $\gamma
_{ij}\,$is a fundamental tensor of $F^{(k-1)n}$, determines a generalized
Hamilton space of order $k,$ $GH^{\left( k\right) n}=(M,g^{ij}).$

This space is reducible to a Hamilton space $H^{(k)n}$ only if $\stackrel{\cdot }{%
\partial ^i}\sigma =0.$

Let $g_{ij}$ be the covariant tensor of the fundamental tensor $g^{ij}$ of
the space $GH^{(k)n}=(M,g^{ij}).$ Then we have:

\begin{equation}
g_{ih}g^{hj}=\delta _i^j  \tag{11.1.5}
\end{equation}

Also, we consider the $d$-tensor field

\begin{equation}
C_i^{jh}=-\displaystyle\frac 12g_{is}(\stackrel{\cdot }{\partial ^j}g^{sh}+
\stackrel{\cdot }{\partial ^h}g^{js}-\stackrel{\cdot }{\partial ^s}g^{jh})
\tag{11.1.6}
\end{equation}

Evidently:
\begin{equation}
S_i^{jh}=C_i^{jh}-C_i^{hj}=0.  \tag{11.1.7}
\end{equation}

If the space $GH^{(k)n}$ is reducible to a Hamilton space then:

\begin{equation}
C_i^{jh}=-g_{is}C^{shj}.  \tag{11.1.6a}
\end{equation}

The $d$-tensor field $C_i^{jh}$ are the coefficients of the $w_k$-covariant
derivatives. Indeed, we have:

\begin{equation}
g^{ij}|^h=\stackrel{\cdot }{\partial ^h}
g^{ij}+C_s^{ih}g^{sj}+C_s^{jh}g^{is}=0.  \tag{11.1.8}
\end{equation}

\section{Metrical $N$-Linear Connections}

If we consider a generalized Hamilton space of order $k$, $%
GH^{(k)n}=(M,g^{ij}), $ in general we can not determine a nonlinear connection
only by means of the fundamental tensor $g^{ij}.$

But there are some particular cases when this is possible. For
instance, in examples 1$^{\circ }$ and 2$^{\circ }$ we can consider the
canonical nonlinear connection $\stackrel{\circ }{N}$ of the Finsler space $%
F^{\left( k-1\right) n}=(M,F(x,y^{(1)},...,y^{(k-1)})$ with the
coefficients $(\underset{(1)}{N_j^i},...\underset{(k-1)}{N_j^i}).$
Then the system of functions
$(\underset{(1)}{N_j^i},...\underset{(k-1)}{N_j^i},N_{ij})$ with
$N_{ij}=\displaystyle\frac \delta {\delta
y^{(1)j}}(\underset{(1)}{
N_j^h}p_h),$ (if $\underset{(1)}{N_j^i}(x,y^{(1)}),$ does not depend on $%
y^{(2)i},...,y^{(k-1)i})$ determines a nonlinear connection only by means of
the fundamental tensor $g^{ij}$ of the space.

Now, let us consider an apriori fixed nonlinear connection $N$,
with the coefficients
$(\underset{(1)}{N_j^i},...\underset{(k-1)}{N_j^i},N_{ij})$ on the
manifold $\widetilde{T^{*k}M}$. We will study the geometry of the
space $GH^{(k)n}$ endowed with the nonlinear connection $N$.

The adapted basis to the direct decomposition (6.2.9), \\
$(\displaystyle\frac \delta {\delta x^i},\displaystyle\frac \delta {\delta
y^{(1)i}},...,\displaystyle\frac \partial {\partial y^{(k-1)i}},%
\displaystyle
\frac \partial {\partial p_i})$ is expressed in (6.2.10) and the
adapted cobasis $(dx^i,\delta y^{\left( 1\right) i},...,\delta y^{\left(
k-1\right) i},\delta p_i)$ is written in (6.3.2).

Now, as usual, we can prove:
\begin{teo}
1) A generalized Hamilton space of order $k$ endowed with a nonlinear
connection $N$, has an unique $N$-linear connection $C\Gamma (N)=$
\newline $=(H_{jh}^i,%
\underset{(1)}{C_{jh}^i},...,\underset{(k-1)}{C_{jh}^i},C_i^{jh})$
satisfying the following axioms:\vspace{3mm}

1$^{\circ }$ The nonlinear connection $N$ is apriori given.

2$^{\circ }$ $g^{ij}{}_{|h}=0,\ g^{ij}\stackrel{(\alpha )}{|}_h=0,\ (\alpha
=1,...,k-1),\ g^{ij}|^h=0.$

3$^{\circ }$ $T_{\ jh}^i=0,\ \underset{\left( \alpha \right)
}{S}\text{}_{\ jh}^i=0,\ (\alpha =1,...,k-1), \ S_i^{\ jh}=0.$

2) The coefficients of $C\Gamma (N)$ are the following generalized
Christoffel symbols:
\begin{equation}
\begin{array}{l}
H_{jh}^i=\displaystyle\frac 12g^{is}(\displaystyle\frac{\delta g_{sh}}{%
\delta x^j}+\displaystyle\frac{\delta g_{js}}{\delta x^h}-\displaystyle\frac{%
\delta g_{jh}}{\delta x^s}), \\
\\
\underset{(\alpha )}{C_{jh}^i}=\displaystyle\frac
12g^{is}(\displaystyle
\frac{\delta g_{sh}}{\delta y^{\left( \alpha \right) j}}+\displaystyle\frac{%
\delta g_{js}}{\delta y^{\left( \alpha \right) h}}-\displaystyle\frac{\delta
g_{jh}}{\delta y^{\left( \alpha \right) s}}),\ (\alpha =1,...,k-1), \\
\\
C_i^{jh}=-\displaystyle\frac 12g_{is}(\stackrel{\cdot }{\partial ^j}g^{sh}+%
\stackrel{\cdot }{\partial ^h}g^{js}-\stackrel{\cdot }{\partial ^s}g^{jh}).
\end{array}
\tag{11.2.1}
\end{equation}
\end{teo}
The previous connection $C\Gamma (N)$ is called canonical metrical $N$
-connection.

More general, one proves

\begin{teo}
1) A generalized Hamilton space $GH^{(k)n},$ endowed with a nonlinear
connection $N,$ has an unique $N$- linear connection \newline
$\overline{D}\Gamma (N)=(%
\overline{H}_{jh}^i,\underset{(1)}{\overline{C}_{jh}^i},...,\underset{%
(k-1)}{\overline{C}_{jh}^i},\overline{C}_i^{jh})$ satisfying the axioms:

1$^{\circ }$ $N$ is apriori given on $\widetilde{T^{*k}M}.$

2$^{\circ }$ $h$-$,v_\alpha $- and $w_k$- covariant derivation of $g^{ij}$
vanish:
\begin{equation}
g_{\ \ |h}^{ij}=0,\ g^{ij}\stackrel{(\alpha )}{|}_h=0,\ \ g^{ij}|^h=0,  \tag{11.2.2}
\end{equation}

3$^{\circ }$ The skewsymmetric tensors of torsion
\begin{equation}
\overline{T}_{\ jh}^i=\overline{H}_{jh}^i-\overline{H}_{hj}^i,\underset{%
\left( \alpha \right) }{\overline{S}_{\ jh}^i}=\underset{\left(
\alpha
\right) }{\overline{C}_{jh}^i}-\underset{\left( \alpha \right) }{\overline{%
C}_{hj}^i},(\alpha =1,...,k-1),\ \overline{S}_i^{\ jh}=\overline{C}_i^{jh}-%
\overline{C}_i^{hj}  \tag{11.2.2a}
\end{equation}
are apriori given.

2) This connection has the following coefficients:
\begin{equation}
\begin{array}{l}
\overline{H}_{jh}^i=H_{jh}^i+\displaystyle\frac 12g^{is}(g_{sr}\overline{T}%
_{jh}^r-g_{jr}\overline{T}_{sh}^r+g_{hr}\overline{T}_{js}^r), \\
\\
\underset{\left( \alpha \right)
}{\overline{C}_{jh}^i}=\underset{\left(
\alpha \right) }{C_{jh}^i}+\displaystyle\frac 12g^{is}(g_{sr}\underset{%
\left( \alpha \right)
}{\overline{S}_{jh}^r}-g_{jr}\underset{\left( \alpha
\right) }{\overline{S}_{sh}^r}+g_{hr}\underset{\left( \alpha \right) }{%
\overline{S}_{js}^r}),\ (\alpha =1,...,k-1), \\
\\
\overline{C}_i^{jh}=C_i^{jh}-\displaystyle\frac 12g_{is}(g^{sh}\overline{S}%
_r^{jh}-g^{jr}\overline{S}_r^{jh}+g^{hr}\overline{S}_r^{js}),
\end{array}
\tag{11.2.2b}
\end{equation}
where $(H_{jh}^i,\underset{(1)}{C}\text{}_{jh}^i,...,%
\underset{(k-1)}{C}\text{}_{jh}^i, C_i^{jh})$ are the coefficients of the canonical metrical $N$-connection $C\Gamma
(N).$
\end{teo}

If we are interested on the all $N$-linear connections which verify the
equations (11.2.2), we can prove:

\begin{teo}
In a generalized Hamilton space $GH^{(k)n}=(M,g^{ij})$ the set of all $N$
-linear connections $\overline{D}\Gamma (N)$ which satisfy the equations
(11.2.2) is given by

\begin{equation}
\begin{array}{l}
\overline{H}_{jh}^i=H_{jh}^i+\Omega _{rj}^{is}X_{sh}^r, \\
\\
\underset{\left( \alpha \right)
}{\overline{C}_{jh}^i}=\underset{\left( \alpha \right)
}{C_{jh}^i}+\Omega _{rj}^{is}\underset{\left( \alpha
\right) }{X_{sh}^r},\ (\alpha =1,...,k-1) \\
\\
\overline{C}_i^{jh}=C_i^{jh}+\Omega _{ri}^{js}X_s^{rh}
\end{array}
\tag{11.2.3}
\end{equation}
where $\Omega _{rj}^{is}=\displaystyle\frac 12(\delta _r^i\delta
_j^s-g_{rj}g^{is})$ are Obata's operators and $C\Gamma (N)=(H_{jh}^i,%
\underset{(1)}{C_{jh}^i},...,$ \newline
$\underset{(k-1)}{C_{jh}^i},C_i^{jh})$
is the canonical metrical $N$-connection and
$X_{sh}^r,\underset{\left( \alpha \right) }{\text{
}X_{sh}^r},(\alpha =1,...,k-1),$ $X_s^{rh}$ are arbitrary
$d$-tensor fields.
\end{teo}

\begin{cor}
The mappings $D\Gamma (N)\rightarrow \overline{D}\Gamma (N)$ determined by
(11.2.3), together with the composition of these mappings is an abelian group.
\end{cor}

Now, we can repeat all considerations from the section 7 of the chapter 8.
So, we have:

\begin{prop}
The curvature $d$-tensor fields of the canonical metrical $N$-connection $%
C\Gamma (N),$ (11.2.1) satisfy the identities (7.6.8).
\end{prop}

\begin{prop}
The tensors of deflection of $C\Gamma (N)$, (11.2.1) satisfy the
identities (7.6.9), with $T_{\ jk}^i=0,$ $S_i^{\ jh}=0,$
$\underset{\left( \alpha \right) }{S_{\ jh}^i}=0,\ (\alpha
=1,...,k-1).$
\end{prop}

Let $\omega _j^i$ be the $1$-forms connection of $C\Gamma (N),$

\begin{equation}
\omega _j^i=H_{js}^idx^s+\underset{\alpha =1}{\stackrel{k-1}{\sum
}} \underset{\left( \alpha \right) }{C_{js}^i}\delta y^{\left(
\alpha \right) s}+C_i^{js}\delta p_s  \tag{11.2.4}
\end{equation}
and a curve $\gamma :t\in I\rightarrow \gamma (t)\in T^{*k}M$ expressed by
(7.7.1).

Thus, we have (cf. Th. 7.7.1):

\begin{teo}
In a space $GH^{\left( k\right) n}$ endowed with the canonical metrical $N$
-connection $C\Gamma (N)$ a vector field $X=X^{\left( 0\right) i}%
\displaystyle\frac \delta {\delta x^i}+...+X^{\left( k-1\right)i }%
\displaystyle\frac \delta {\delta y^{\left( k-1\right)
i}}+X_i\displaystyle \frac \partial {\partial p_i}$ is parallel
along curve $\gamma $ if and only if $(X^{\left( 0\right)
i},...,X^{\left( k-1\right) },X_i)$ are the solutions of the
system of differential equations
\begin{equation}
\displaystyle\frac{dX^{\left( 0\right) i}}{dt}+X^{\left( 0\right) s}%
\displaystyle\frac{\omega _s^i}{dt}=0,...,\displaystyle\frac{dX^{\left(
k-1\right) i}}{dt}+X^{\left( k-1\right) s}\displaystyle\frac{\omega _s^i}{dt}%
=0,\ \displaystyle\frac{dX_i}{dt}-X_s\displaystyle\frac{\omega _i^s}{dt}=0.
\tag{11.2.5}
\end{equation}
\end{teo}

\begin{teo}
The space $GH^{(k)n}$ endowed with the $N$-linear connection $C\Gamma (N)$
is with absolute parallelism of vectors if and only if all curvature $d$
-tensor of $C\Gamma (N)$ vanish.

\end{teo}

A smooth parametrized curve $\gamma :t\in I\rightarrow
(x(t),y^{(1)}(t),...,y^{(k-1)}(t),p(t))\in T^{*k}M$ is an autoparallel curve of
$C\Gamma (N)$ if $D_{\stackrel{\cdot }{\gamma }}\stackrel{\cdot }{\gamma }
=0. $

Thus, applying Theorem 7.7.3 we have:

\begin{teo}
The curve $\gamma :I\rightarrow \widetilde{T^{*k}M}$ is autoparallel for the
space $GH^{(k)n},$ with respect to the canonical $N$-linear connection $%
C\Gamma (N)\,$ if and only if the following system of differential equations
is verified:
\end{teo}

\begin{equation}
\begin{array}{l}
\displaystyle\frac{d^2x^i}{dt^2}+\displaystyle\frac{dx^s}{dt}\displaystyle
\frac{\omega _s^i}{dt}=0, \\
\\
\displaystyle\frac d{dt}\left( \displaystyle\frac{\delta y^{\left( \alpha

\right) i}}{dt}\right) +\displaystyle\frac{\delta y^{\left( \alpha \right)
s} }{dt}\displaystyle\frac{\omega _s^i}{dt}=0,(\alpha =1,...,k-1), \\
\\
\displaystyle\frac d{dt}\left( \displaystyle\frac{\delta p_i}{dt}\right) - %
\displaystyle\frac{\delta p_s}{dt}\displaystyle\frac{\omega _i^s}{dt}=0.
\end{array}
\tag{11.2.6}
\end{equation}

Recall that $\gamma $ is an horizontal curve if and only if:
\[
x^i=x^i(t),\displaystyle\frac{\delta p_i}{dt}=0,\displaystyle\frac{\delta
y^{\left( \alpha \right) i}}{dt}=0,(\alpha =1,...,k-1).
\]

Therefore, we have:

\begin{teo}
For a generalized Hamilton space of order $k,GH^{(k)n}$ endowed with the
canonical metrical $N$-connection $C\Gamma (N)$ the following properties
hold:

1$^{\circ }$ The horizontal paths are characterized by the system of
differential equations:
\[
\displaystyle\frac{d^2x^i}{dt^2}+H_{jh}^i\displaystyle\frac{dx^j}{dt}%
\displaystyle\frac{dx^h}{dt}=0,\ \displaystyle\frac{\delta y^{\left( \alpha
\right) i}}{dt}=0,\ \displaystyle\frac{\delta p_i}{dt}=0,(\alpha =1,...,k-1).
\]

2$^{\circ }$ The $v_\alpha $- paths in a point $x_0\in M$ are characterized
by
\[
x^i=x_0^i,\displaystyle\frac{dy^{(\beta )i}}{dt}=0,(\beta \neq \alpha ),%
\displaystyle\frac{dp_i}{dt}=0,
\]
\[
\displaystyle\frac d{dt}\left( \displaystyle\frac{\delta y^{\left( \alpha
\right) i}}{dt}\right) +\underset{\left( \alpha \right) }{C_{sj}^i}%
\displaystyle\frac{dy^{\left( \alpha \right) s}}{dt}\displaystyle\frac{%
dy^{\left( \alpha \right) j}}{dt}=0
\]

3$^{\circ }$ The $w_k$-path are characterized by
\[
\displaystyle\frac{dx^i}{dt}=0,\displaystyle\frac{dy^{\left( 1\right) i}}{dt}%
=...=\displaystyle\frac{dy^{\left( k-1\right) i}}{dt}=0,
\]
\[
\displaystyle\frac{d^2p_i}{dt^2}-C_i^{jm}(x_0,0,...0,p)\displaystyle\frac{%
dp_j}{dt}\displaystyle\frac{dp_m}{dt}=0.
\]
\end{teo}

Finally, we remark:

The structure equations of the canonical metrical $N$-connection $C\Gamma
(N)\,$ of the space $GH^{\left( k\right) n}$ are given by Theorem 7.8.1.

where $T_{\ jh}^i=0,$ $\underset{\left( \alpha \right) }{S_{\ jh}^i=0},$ $%
(\alpha =1,...,k-1),$ $S_i^{\ jh}=0.$

\section{Hamiltonian Relativistic Optics}

In the book: 'The Geometry of Higher Order Lagrange Spaces' -
Kluwer FTPH, Vol. 82 is given a generalized Lagrange metric,
formula (10.5.6) of the Relativistic Optics of order $k$. It is
rather complicated. In the dual spaces $T^{*k}M$ it can be
introduced much more simple.

Consider a semidefined Finsler space of order $k-1,$ $F^{(k-1)}M=$ \newline
$=(M,F(x,y^{(1)},...,y^{%
\left( k-1\right) }))$ and $\gamma _{ij}(x,y^{\left( 1\right)
},...,y^{\left( k-1\right) })$ its fundamental tensor field.

The projection $\pi _{k-1}^{*k}:(x,y^{(1)},...,y^{(k-1)},p)\in
T^{*k}M\rightarrow (x,y^{(1)},...,y^{(k-1)})\in T^{k-1}M$ allows to consider
$d$-tensor $\gamma _{ij}\circ \pi _{k-1}^{*k}$ on $T^{*k}M$. It will be
denoted by $\gamma _{ij}$, too. Its contravariant $\gamma ^{ij}$ will be
considered defined on the manifold $\widetilde{T^{*k}M}.$

Let us consider a differentiable function $n$ on $T^{*k}M$ with the property
$n>1$. It will be called a refractive index.

Some notations:

\begin{equation}
\left\{
\begin{array}{l}
\stackrel{\vee}{p^i}=\gamma ^{ij}p_j,\ p_i=\gamma _{ij}\stackrel{\vee}{p^j}
,\ p_i\stackrel{\vee}{p^i}=\gamma ^{ij}p_ip_j=\left\| p\right\| ^2 ,\\
\\
\displaystyle\frac 1{n(x,y^{\left( 1\right) },...,y^{\left( k-1\right)
},p)}=u(x,y^{\left( 1\right) },...,y^{\left( k-1\right) },p).
\end{array}
\right.  \tag{11.3.1}
\end{equation}

Now, we define on $\widetilde{T^{*k}M}$ the $d$-tensor field

\begin{equation}
\begin{array}{l}
g^{ij}(x,y^{\left( 1\right) },...,y^{\left( k-1\right) },p)=\\
\\
\qquad \qquad = \gamma
^{ij}(x,y^{(1)},...,y^{\left( k-1\right) })+\left( 1-\displaystyle\frac
1{n^2(x,y^{(1)},...,y^{\left( k-1\right) },p)}\right) \stackrel{\vee}{p^i}
\stackrel{\vee}{p^j}  \tag{11.3.2}
\end{array}
\end{equation}

\begin{teo}
1$^{\circ }$ $g^{ij}$ is a symmetric $d$-tensor field of type (2, 0) on the
manifold $\widetilde{T^{*k}M}$.

2$^{\circ }$ $rank\left\| g^{ij}\right\| =n.$
\end{teo}

Indeed, 1$^{\circ }$ $g^{ij}$ is a sum of two symmetric $d$-tensor of type
(2,0).

2$^{\circ }$ Consider the $d$-tensor

\begin{equation}
g_{ij}=\gamma _{ij}-\displaystyle\frac 1{a^2}\left( 1-\displaystyle\frac
1{n^2}\right) p_ip_j , \tag{11.3.3}
\end{equation}
where

\begin{equation}
a=1+\left( 1-\displaystyle\frac 1{n^2}\right) \left\| p\right\| ^2.
\tag{11.3.3a}
\end{equation}
It is easy to verify the following equality:

\begin{equation}
g^{ih}g_{hj}=\delta _j^i.  \tag{11.3.4}
\end{equation}

Consequently the pair $GH^{(k)n}=(M,g^{ij})$ is a generalized Hamilton space
of order $k.$

\begin{teo}
The space $GH^{\left( k\right) n}=(M,g^{ij})$ is not reducible to a
Hamilton space of order $k$.
\end{teo}

\textbf{Proof}. The tensor field $C^{ijk}$ from (11.1.3) is as follows
\[
-C^{ijh}=\stackrel{\cdot }{\partial ^h}\sigma \stackrel{\vee}{p^i}\stackrel{
\vee}{p^j}+\sigma (\gamma ^{ih}\stackrel{\vee}{p^j}+\gamma ^{jh}\stackrel{
\vee}{p^i}),\ \sigma =1-u^2
\]

If we assume that $C^{ijh}=C^{ihj},$ we obtain

\begin{equation}
\left( \stackrel{\cdot }{\partial ^h}\sigma \stackrel{\vee}{p^j}-\stackrel{
\cdot }{\partial ^j}\sigma \stackrel{\vee}{p^h}\right) \stackrel{\vee}{p^i}
+\sigma \left( \gamma ^{ih}\stackrel{\vee}{p^j}-\gamma ^{ij}\stackrel{\vee
}{p^h}\right) =0  \tag{*}
\end{equation}
Contracting by $p_i$ we have

\[
\stackrel{\cdot }{\partial ^h}\sigma \stackrel{\vee}{p^j}-\stackrel{\cdot }{
\partial ^j}\sigma \stackrel{\vee}{p^h}=0.
\]

Substituting in (*) we get $\gamma ^{ih}\stackrel{\vee}{p^j}-\gamma ^{ij}
\stackrel{\vee}{p^h}=0$. A new contraction with $\gamma _{ih}$ leads to $%
(n-1)\stackrel{\vee}{p^j}=0.$ Consequently $\stackrel{\vee}{p^j}=0,$ i.e $%
p_j=0.$ But this is impossible and the assumption we made is false. q.e.d.

The space $GH^{\left( k\right) n}$ will be called the generalized
Hamiltonian space of order $k$ of the Hamiltonian Relativistic Optics.

Let us consider a local $d$-vector field $V^i(x)$ and a local $d$-covector
field $\eta _i(x)$ on the manifold $M$. It is not difficult to see that the
mapping \newline  $S_{V,\eta }:M\rightarrow T^{*k}M$ defined locally by

\begin{equation}
\begin{array}{l}
x^i=x^{i^{\prime }}, \\
\\
y^{\left( 1\right) i}=V^i(x),...,y^{\left( k-1\right) i}=\displaystyle\frac
1{\left( k-2\right) !}\displaystyle\frac{d^{k-2}V^i}{dt^{k-2}}, \\
\\
p_i=\eta _i(x)
\end{array}
\tag{11.3.5}
\end{equation}
is a cross- section of the canonical projection $\pi
^{*k}:T^{*k}M\rightarrow M.$ It follows that $S_{V,\eta }(M)$ is a local
embedding of $M$ in the manifold $\widetilde{T^{*k}M}.$

The restriction of the fundamental tensor $g^{ij}$ of the space
$GH^{(k)n}$ to $S_{V,\eta }$ will be called the Synge metric of the \textit{%
\ Hamiltonian Relativistic Optics}. The restriction of function $%
n(x,y^{\left( 1\right) },...,y^{\left( k-1\right) },p)$ to $S_{V,\eta }(M)$
is the refractive index of the dispersive optic medium:

$(M,V(x),\eta (x),n(x,V(x),...,\displaystyle\frac 1{\left( k-2\right) !} %
\displaystyle\frac{d^{k-2}V}{dt^{k-2}},\eta (x)).$

Therefore, we say that the geometry of the generalized Hamilton space of
order $k,$ $GH^{\left( k\right) n}=(M,g^{ij}),$ with the fundamental tensor $%
g^{ij}$ in (11.3.2) is the geometrical theory of the previous optic medium,
endowed with Synge metric.

If the refractive index $n$ depend on $x\in M$ only, then the optic medium
is called nondispersive.

Let us consider the canonical nonlinear connection
$\stackrel{\circ }{N}$ with the dual coefficients
$\underset{\left( 1\right) }{M_j^i},..., \underset{\left(
k-1\right) }{M_j^{i\,}}$ of the Finsler space $F^{\left(
k-1\right) n}=(M,F).$ For simplicity we assume that
$\underset{\left( 1\right) }{M_j^i}=\underset{\left( 1\right)
}{M_j^i}(x,y^{\left( 1\right) }),$ (see ch. 8, \S 8). In this case
($\underset{\left( 1\right) }{M_j^i},..., \underset{\left(
k-1\right) }{M_j^{i\,}},N_{ij})$, with:

\begin{equation}
N_{ij}=\displaystyle\frac \delta {\delta y^{\left( 1\right) i}}(\underset{%
\left( 1\right) }{M_j^h})p_h  \tag{11.3.6}
\end{equation}
define a nonlinear connection $N\;$of the space $GH^{\left( k\right)
n}=(M,g^{ij})$ determined only by the fundamental tensor $g^{ij}.$

Let $(\displaystyle\frac \delta {\delta x^i},\displaystyle\frac \delta
{\delta y^{\left( 1\right) i}},...,\displaystyle\frac \delta {\delta
y^{\left( k-1\right) i}},\displaystyle\frac \delta {\delta p_i})$ and $%
(dx^i,\delta y^{\left( 1\right) i},...,\delta y^{\left( k-1\right) i},\delta
p_i)$ be the \newline adapted local basis and adapted local cobasis corresponding to
the nonlinear connection $N$.

We can determine the canonical metrical $N$-linear connection $C\Gamma (N)$
of the space $GH^{\left( k\right) n}$ starting from the expressions (11.3.2)
and (11.3.3) for the fundamental tensor $g^{ij}$ and its covariant tensor field
$g_{ij}$ and applying the usual tehniques of calculus.

The connections $N$ and $C\Gamma (N)$ allow to study the geometry of the
space $GH^{(k)n}.$ Taking the restriction of its main geometrical object
field to the section $S_{V,\eta }(M)$ we obtain the main notions and
properties of the Hamiltonian Relativistic Optics.

Evidently, there are some particular cases, as:

a) the nondispersive media

b) the Finsler space $F^{(k-1)n\text{ }}$is the prolongation of order $k-1$
of a Finsler space $F^n=(M,F(x,y^{\left( 1\right) })).$

Let us consider the absolute energy of the space $GH^{\left( k\right) n},$ [115]:

\begin{equation}
\mathcal{E}=g^{ij}p_ip_j  \tag{11.3.7}
\end{equation}

Taking into account the expression (11.3.2) of $g^{ij}$ and the function $a$
from (11.3.2a), we get:

\begin{equation}
\mathcal{E}=a\left\| p\right\| ^2  \tag{11.3.7a}
\end{equation}

Consequently, we can say that $\mathcal{E}(x,y^{\left( 1\right)
},...,y^{\left( k-1\right) },p)$ is a differentiable Hamiltonian uniquely
determined by the fundamental tensor $g^{ij}$ of the space $GH^{(k)n}$. It
allows to determine the Hamilton - Jacobi equations of the space. These equations are given by (5.1.17), (5.1.17'). The
energies of order $k-1,...,1$ are expressed in (5.3.1) with $H=\mathcal{%
\ E}$ and the low of conservation is mentioned in Theorem 5.3.2. Also,
Theorem 5.5.3 gives the N\"other symmetries for the Hamiltonian $H=\mathcal{E%
}$ from (11.3.7).

\section{The Metrical Almost Contact Structure of the Space $GH^{(k)n}$}

The generalized Hamilton space of order $k,$ $GH^{(k)n}=(M,g^{ij}),$ endowed
with an apriori given nonlinear connection $N$ determines a metrical almost
contact structure on the manifold $\widetilde{T^{*k}M}.$

The $N$-lift of the fundamental tensor field $g^{ij}$ is

\begin{equation}
\stackrel{\vee}{\Bbb{G}}=g_{ij}dx^i\otimes dx^j+\underset{\alpha =1}{%
\stackrel{k-1}{\sum }}g_{ij}\delta y^{\left( \alpha \right) i}\otimes \delta
y^{\left( \alpha \right) j}+g^{ij}\delta p_i\otimes \delta p_j  \tag{11.4.1}
\end{equation}

Evidently:

1$^{\circ }$ $\stackrel{\vee}{\Bbb{G}}$ is a pseudo-Riemannian structure on $T^{*k}M.$

2$^{\circ }$ The distributions
$\underset{0}{N},\underset{1}{N},..., \underset{k-1}{N},V_{k-1}$
and $W_k\,$are mutual orthogonal with respect to
$\stackrel{\vee}{\Bbb{G}}.$ \vspace{3mm}

3$^{\circ }$ $\Bbb{G}^H=g_{ij}dx^i\otimes dx^j,\ \Bbb{G}^{V_\alpha
}=g_{ij}\delta y^{\left( \alpha \right) i}\otimes \delta y^{\left( \alpha
\right) j},(\alpha =1,...,k-1),\ $\vspace{3mm}\\
$\Bbb{G}^{W_k}=g^{ij}\delta p_i\otimes \delta p_j$ are $d$-tensor fields.\vspace{3mm}

4$^{\circ }$ $\Bbb{G}^H,\Bbb{G}^{V_\alpha },(\alpha =1,...,k-1)$ and $\Bbb{G}
^{W_k}$ are covariant constant with respect to the canonical metrical $N$
-connection $C\Gamma (N).$

The geometrical object fields $N$ and $g^{ij}$ determine an almost contact $%
(k-1)n$- structure $\stackrel{\vee}{\Bbb{F}},$(given in (6.8.2)):

\begin{equation}
\stackrel{\vee}{\Bbb{F}}\left( \displaystyle\frac \delta {\delta
x^i}\right) =-g_{ij}\displaystyle\frac \delta {\delta p_j},\ \stackrel{\vee
}{\Bbb{F}}\left( \displaystyle\frac \delta {\delta y^{\left( \alpha \right)
i}}\right) =0,(\alpha =\overline{1,k-1}),\ \stackrel{\vee}{\Bbb{F}}\left( %
\displaystyle\frac \delta {\delta p_i}\right) =g^{ij}\displaystyle\frac
\delta {\delta x^j}  \tag{11.4.2}
\end{equation}

Theorem 6.8.1, from ch. 6 is valid:

\begin{teo}
1$^{\circ }$ The structure $\stackrel{\vee}{\Bbb{F}}$ is defined only by $N$
and $g^{ij}.$

2$^{\circ }$ $\stackrel{\vee}{\Bbb{F}}$ is the following tensor field of
type (11.1.1) on $T^{*k}M$:
\begin{equation}
\stackrel{\vee}{\Bbb{F}}=-g_{ij}\displaystyle\frac \delta {\delta
p_j}\otimes dx^i+g^{ij}\displaystyle\frac \delta {\delta x^i}\otimes \delta
p_j . \tag{11.4.3}
\end{equation}

3$^{\circ }$ $Ker\stackrel{\vee}{\Bbb{F}}=\underset{1}{N}\oplus
...\oplus
\underset{k-1}{N},$ $Im$ $\stackrel{\vee}{\Bbb{F}}=\underset{1}{N}%
\oplus W_k.$

4$^{\circ }$ $rank$ $\stackrel{\vee}{\Bbb{F}}=2n.$

5$^{\circ }$ $\stackrel{\vee}{\Bbb{F}^3}+\stackrel{\vee}{\Bbb{F}}=0.$
\end{teo}

Consequently $\stackrel{\vee}{\Bbb{F}}$ is an almost $(k-1)n$- contact
structure on the manifold $\widetilde{T^{*k}M}.$

The condition of normality of the structure $\stackrel{\vee}{\Bbb{F}}$ is
given by (see (6.6.5)):

\begin{equation}
\mathcal{N}_{\stackrel{\vee
}{\Bbb{F}}}(X,Y)+\underset{i=1}{\stackrel{n}{ \sum }}\left[
\underset{\alpha =1}{\stackrel{k-1}{\sum }}d(\delta y^{\left(
\alpha \right) i})(X,Y)+d(\delta p_i)(X,Y)\right] =0,\ \forall
X,Y\in \mathcal{X}(T^{*k}M)  \tag{11.4.4}
\end{equation}
where $\mathcal{N}_{\stackrel{\vee}{\Bbb{F}}}$ is the Nijenjuis tensor of $%
\stackrel{\vee}{\Bbb{F}}.$

Theorem 8.9.4 is valid for spaces $GH^{(k)n}.$

\begin{teo}
For any generalized Hamilton space of order $k$, $GH^{\left( k\right)
n}=(M,g^{ij})$ endowed with a nonlinear connection $N$ the following
properties hold.

1$^{\circ }$ The pair $(\stackrel{\vee}{\Bbb{G}},,\stackrel{\vee}{\Bbb{F}})
$ is a pseudo-Riemannian almost $(k-1)n$-contact structure determined only
by $N$ and $g^{ij}$.

2$^{\circ }$ The associated $2$-form is
\[
\theta =\delta p_i\wedge dx^i
\]

3$^{\circ }$ If the coefficients $N_{ij}$ of $N$ are symmetric, then
\[
\theta =dp_i\wedge dx^i
\]
is the canonical presymplectic structure on the manifold $\widetilde{T^{*k}M}.$

4$^{\circ }$ The conditions of normality of the structure $\stackrel{\vee}{%
\Bbb{F}}$ is expressed by (11.4.4).
\end{teo}
Finally, taking into account $\Bbb{G}$ from (11.4.1) and $\Bbb{F}$ from (11.4.3) it follows.

\begin{teo}
With respect to the canonical metrical connection $C\Gamma (N)\,$ of the
space $GH^{\left( k\right) n}\,$we have
\[
D_X\stackrel{\vee}{\Bbb{G}}=0,D_X\stackrel{\vee}{\Bbb{F}}=0
\]
\end{teo}

As such the geometry of the pseudo-Riemannian almost $(k-1)n$-contact
space $(T^{*k}M,\stackrel{\vee}{\Bbb{G}},,\stackrel{\vee}{\Bbb{F}})$ can
be studied by means of the canonical metrical $N$-linear connection $C\Gamma
(N)$ of the generalized Hamilton space of order $k,$ $GH^{(k)n}$.

\section{Generalized Cartan Space of Order $k$}

\begin{defi}
A generalized Cartan space of order $k,$ is a Generalized Hamilton space of
order $k$, $GH^{\left( k\right) n}=(M,g^{ij})$ in which the fundamental
tensor $g^{ij}$ satisfies the axioms:

1$^{\circ }$ $g^{ij}$ is positively defined on $\widetilde{T^{*k}M}.$

2$^{\circ }$ $g^{ij}$ is $0$-homogeneous on the fibres of the dual bundle $%
(T^{*k}M,\pi ^{*k},M).$
\end{defi}

We denote by $G\mathcal{C}^{(k)n}=(M,g^{ij})\,$a generalized Cartan space of
order $k.$

>From the axiom $2^{\circ }$ it follows

\begin{prop}
The following identities hold:

1$^{\circ }$ $g^{ij}$ being $0$-homogeneous, we have
\begin{equation}
\mathcal{L}_{\stackrel{k-1}{\Gamma }+kC^{*}}g^{ij}=0  \tag{11.5.1}
\end{equation}
or, developed:
\begin{equation}
y^{\left( 1\right) i}\displaystyle\frac{\partial g^{jh}}{\partial y^{\left(
1\right) i}}+...+(k-1)y^{\left( k-1\right) i}\displaystyle\frac{\partial
g^{jh}}{\partial y^{\left( k-1\right) i}}+kp_i\stackrel{\cdot }{\partial ^i}%
g^{jh}=0  \tag{11.5.1a}
\end{equation}

2$^{\circ }$ The absolute energy
\begin{equation}
\mathcal{E}=g^{ij}p_ip_j  \tag{11.5.2}
\end{equation}
is $2k$-homogeneous on the fibres of $T^{*k}M$.

3$^{\circ }$ $\mathcal{L}_{\stackrel{k-1}{\Gamma }+kC^{*}}\mathcal{E}=2k%
\mathcal{E}$.

4$^{\circ }$ $\mathcal{L}_{\stackrel{k-1}{\Gamma }+kC^{*}}C^{ijh}=-kC^{ijh}.$
\end{prop}

\textbf{Example 2.} Let $\stackrel{\circ }{g}^{ij}\,$be the
fundamental tensor of the Cartan
space $\mathcal{C}^{(k)n\text{ }}$and $\sigma \in \mathcal{F}(\widetilde{%
T^{*k}M})$ with the properties:

a) $\sigma $ is $0$-homogeneous;

b) $\stackrel{\cdot }{\partial ^i}\sigma $ nonvanishes.

The pair $G\mathcal{C}^{(k)n\text{ }}=e^{-2\sigma }\stackrel{\circ
}{g}^{ij}$ is a generalized Cartan space of order $k$.

In particular we can consider
\[
\sigma =\displaystyle\frac{p_iy^{\left( 1\right) i}}{\sqrt{\stackrel{\circ }{
g}^{ij}p_ip_j}\sqrt{\stackrel{\circ }{g}_{ij}y^{\left( 1\right) i}y^{\left(
1\right) j}}}.
\]

The previous example shows the existence of the spaces $G\mathcal{C}^{(k)n
\text{ }}$ are not reducible to a Cartan space of order $k$.

Let $N^{*}$ be a nonlinear connection on $T^{*k}M,$ having the
coefficients \newline $(\underset{\left( 1\right)
}{M_j^{*i}},...,\underset{\left( k-1\right)}{M_j^{*i\,}},N_{ij})$
homogeneous of degree $k-1,...,1$, $k\,$respectively.

\begin{teo}
There exists an unique canonical metrical $N^{*}$-connection of the space $%
\Bbb{G}\mathcal{C}^{\left( k\right) n}$. Its coefficients are given by
(9.5.5).
\end{teo}

Now, the geometry of $\Bbb{G}\mathcal{C}^{\left( k\right) n}\,$can be
studied like the geometry of $GH^{\left( k\right) n}$ spaces.

The $N^{*}$-lift to $\widetilde{T^{*k}M}$ of the fundamental tensor $g^{ij}$
of $G\mathcal{C}^{\left( k\right) n},\,$given by (11.3.2) is not homogeneous on
the fibres of the bundle $T^{*k}M$. Therefore, taking into account the
following Hamiltonians:

\begin{equation}
\mathcal{E}_1=g_{ij}\stackrel{\left( 1\right) i}{z}\stackrel{\left( 1\right)
j}{z},\mathcal{E}_2=g_{ij}\stackrel{\left( 2\right) i}{z}\stackrel{\left(
2\right) j}{z},...,\mathcal{E}_{k-1}=g_{ij}\stackrel{\left( k-1\right) i}{z}
\stackrel{\left( k-1\right) j}{z},\mathcal{E}=g^{ij}p_ip_j,  \tag{11.5.3}
\end{equation}
we can define the following tensor field:

\begin{equation}
\stackrel{\circ }{\Bbb{G}}=g_{ij}dx^i\otimes dx^j+\underset{\alpha =1}{%
\stackrel{k-1}{\sum }}\displaystyle\frac 1{\mathcal{E}_\alpha }g_{ij}\delta
y^{\left( \alpha \right) i}\otimes \delta y^{\left( \alpha \right) j}+ %
\displaystyle\frac 1{\mathcal{E}}g^{ij}\delta p_i\otimes \delta p_j.
\tag{11.5.4}
\end{equation}
Evidently, $\mathcal{E} >0.$
\begin{teo}
$\stackrel{\circ }{\Bbb{G}}$ is a Riemannian structure on $T^{*k}M$,
determined only by the fundamental tensor $g^{ij}$ of the space $G\mathcal{C}%
^{\left( k\right) n}$ and the nonlinear connection $N^{*}.$
\end{teo}

Consider the ${\cal F}(T^{*k}M)$-linear mapping $\stackrel{\circ }{\Bbb{F}}:
\mathcal{X}(\widetilde{T^{*k}M})\rightarrow \mathcal{X}(\widetilde{T^{*k}M})$
defined by:

\begin{equation}
\stackrel{\circ }{\Bbb{F}}\left( \displaystyle\frac \delta {\delta
x^i}\right) =-\mathcal{E}g_{ij}\displaystyle\frac \delta {\delta p_j},\
\stackrel{\circ }{\Bbb{F}}\left( \displaystyle\frac \delta {\delta y^{\left(
\alpha \right) i}}\right) =0,(\alpha =1,...,k-1),\ \stackrel{\circ }{\Bbb{F}}
\left( \displaystyle\frac \delta {\delta p_i}\right) =\displaystyle\frac{
g^{ij}}{\mathcal{E}}\displaystyle\frac \delta {\delta x^j}  \tag{11.5.5}
\end{equation}

We obtain:

\begin{teo}
We have:

1$^{\circ }$ $\stackrel{\circ }{\Bbb{F}}$ is a tensor field on $T^{*k}M$ of
type (1,1).

2$^{\circ }$ $\stackrel{\circ }{\Bbb{F}}$ is expressed in the adapted basis
by
\begin{equation}
\stackrel{\circ }{\Bbb{F}}=-\mathcal{E}g_{ij}\displaystyle\frac \delta
{\delta p_j}\otimes dx^i+\displaystyle\frac{g^{ij}}{\mathcal{E}}%
\displaystyle \frac \delta {\delta x^i}\otimes \delta p_j  \tag{11.5.6}
\end{equation}

3$^{\circ }$ $rank\stackrel{\circ }{\Bbb{F}}=2n$

4$^{\circ }$ $\stackrel{\circ }{\Bbb{F}^3}+\stackrel{\circ }{\Bbb{F}}=0$

5$^{\circ }$ $\stackrel{\circ }{\Bbb{F}}$ is determined only by $g^{ij}$ and
$N^{*}$.
\end{teo}

\begin{teo}
For a generalized Cartan spaces of order $k,$ $G\mathcal{C}^{(k)n}=(M,g^{ij})
$ endowed with a nonlinear connection $N^{*}$ the following properties hold:

1$^{\circ }$ The pair $(\stackrel{\circ }{\Bbb{G}},\stackrel{\circ }{\Bbb{F}}%
)$ is a Riemannian almost $(k-1)n$-contact structure determined only by $%
N^{*}$ and $g^{ij}$.

2$^{\circ }$ The associated $2$-form is
\[
\stackrel{\circ }{\theta }=\mathcal{E}g_{ij}\delta p_i\wedge dx^j
\]
\end{teo}

The proofs are made by usual methods.

Concluding, the space $(\widetilde{T^{*k}M},\stackrel{\circ }{\Bbb{G}},%
\stackrel{\circ }{\Bbb{F}})\,$is the geometrical model for the generalized
Cartan space of order $k$.

\begin{theindex}
\markboth{\it{THE GEOMETRY OF HIGHER-ORDER HAMILTON SPACES\ \ \ \ \ }}{\it{Index}}
\item Absolute energies 272
\item Accelerations of order $k$ 2
\item Adapted
\subitem basis 12, 121 \subitem cobasis 14, 123
\item Almost
\subitem complex structure \subitem contact structure 139, 204
\subitem Hermitian structure 70 \subitem K\"ahlerian structure 70
\subitem product structure 138  \subitem symplectic structure
\item Autoparallel curves 18, 131
\\
\item Berwald
\subitem connections 165 \subitem spaces 160
\item Bianchi identities 31, 166
\item Bundle
\subitem differentiable 2, 72 \subitem $k$-osculator 2, 72
\subitem of higher order accelerations 2 \subitem dual 72
\\
\item Canonical
\subitem spray 44, 60 \subitem semispray 43 \subitem $k$-semispray
43  \subitem metrical connections 46, 246 \subitem nonlinear
connections 44, 194 \subitem $N$-linear connections 46, 195
\item Cartan
\subitem nonlinear connection 61 \subitem metrical connection 67
\subitem spaces 233 \subitem spaces of higher order 233
\item Christoffel
\subitem symbols 19, 117 \subitem generalized symbols 46, 197
\item Coefficients
\subitem of a nonlinear connections 44 \subitem of a $N$-linear
connection 48, 155 \subitem dual of a nonlinear connection 44, 195
\subitem primal of a nonlinear connection 44, 195
\item Connections
\subitem nonlinear 11, 195 \subitem $N$-metrical 46, 197 \subitem
Levi-Civita \subitem Berwald 165
\item Covariant
\subitem $h$- and $v$- derivatives 25 \subitem $h$- and $v_\alpha
$- derivatives 25, 157 \subitem $h$- and $w_\alpha $- derivatives
157 \subitem differential 24
\item Craig-Synge equation 59
\item Curvature of an $N$-linear connection 29, 152
\\
\item $d$-vectors 22, 145
\subitem $d$-tensors 22, 145 \subitem distinguished tensors 22,
145
\item Deflection tensors 26
\item Direct decomposition 12, 120
\item Distributions 13
\subitem horizontal 11 \subitem vertical 6, 76
\item Dual
\subitem coefficients 15, 124 \subitem semisprays 113 \subitem of
$T^kM$ 72
\\
\item Electrodinamics 45, 201
\item Energy
\subitem of higher order 39, 99 \subitem absolute 272
\item Embeeding 212
\item Euler-Lagrange equations 36
\\
\item Finsler
\subitem metrics 56 \subitem spaces 55
\item Fundamental
\subitem functions 56, 179 \subitem tensors 56, 179
\\
\item Gauss-Weingarten formulas 226
\item Gauss-Codazzi equations 230
\item General Relativity 266
\item Generalized
\subitem Lagrange spaces 49 \subitem Hamilton spaces 259 \subitem
Cartan spaces 272
\item Geodesics 58
\\
\item Hamilton-Jacobi equations 42, 96
\item Hamilton-Jacobi-Ostrogradski equations 103
\item Hamiltonian system of higher order 186
\item Hamiltonian space of electrodinamics 201
\item Hamilton
\subitem space 179 \subitem vector field 77 \subitem 1-form 77
\item Higher order
\subitem Hamilton space 179 \subitem Lagrange space 42
\\
\item Integral of action 36, 91
\\
\item Jacobi-Ostrogradski momenta 42, 101
\item Jacobi method 42
\\
\item Lagrangian
\subitem differentiable 34 \subitem of higher order 34 \subitem
regular 35
\item Lagrange
\subitem space 42 \subitem space of higher order 42
\item Law of conservation 99
\item Legendre
\subitem mapping 187, 240 \subitem transformation 187, 240
\item Lie derivative 54
\item Lift
\subitem homogeneous 275 \subitem horizontal 11, 118 \subitem
$N$-lift 50, 118
\item Liouville
\subitem vector fields 6, 77 \subitem $d$-vector fields 127
\subitem 1-form 78
\\
\item Main invariants 35
\item $N$-linear connection 46, 147
\item Nonlinear connection 11, 118
\\
\item Poisson structure 83, 184
\item Presymplectic structure 78, 144
\item Projector 118
\subitem horizontal 118 \subitem vertical 188
\\
\item Relativistic Optics 266
\item Ricci identities 161
\\
\item Sections
\subitem in $T^kM$ 3 \subitem in $T^{*k}M$ 72
\item Structure equations 174, 252
\item Spaces
\subitem Finsler 55 \subitem Hamilton 179 \subitem Randers
\item Subspaces in Hamilton Spaces 212
\item Symplectic structure 184
\\
\item Torsions 26, 151
\\
\item Variational problem
\subitem in Lagrange spaces 36 \subitem in Hamilton spaces 91
\end{theindex}

\end{document}